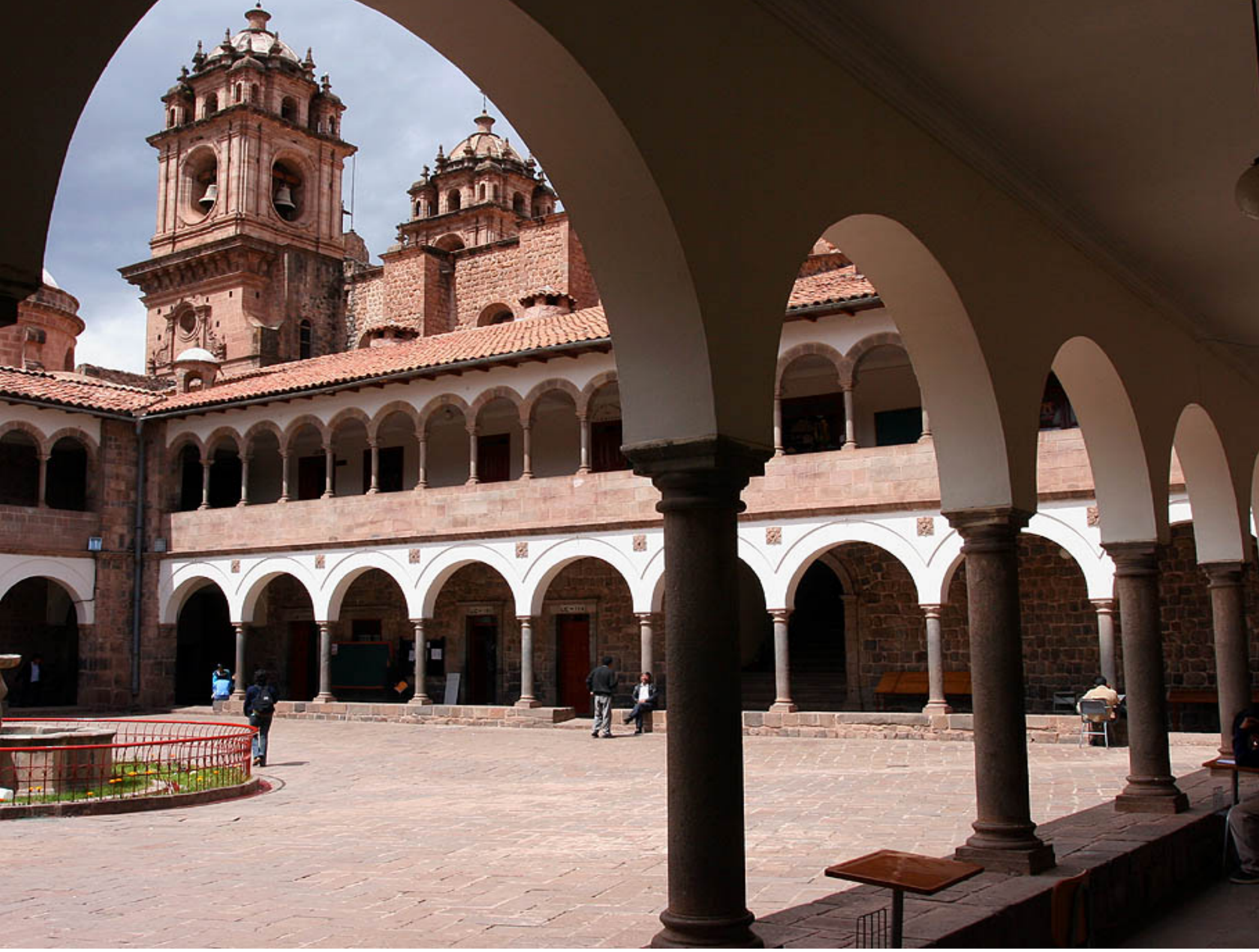

# Aritmética, Grupos y Análisis

AGRA II, Cusco, Perú

Mikhail Belolipetsky
Harald Andrés Helfgott
Carlos Gustavo Moreira (Editores)

# Aritmética, Grupos y Análisis
# AGRA II, Cusco, Perú

**Aritmética, Grupos y Análisis AGRA II, Cusco, Perú**

# Prefacio

He aquí las notas de la escuela AGRA II (Aritmética, Grupos y Análisis), realizada entre el 8 y el 22 de agosto del 2015, en la Universidad de San Antonio Abad del Cusco (UNSAAC), Perú. Una primera versión de estas notas circuló durante la escuela. La versión actual es el producto de un largo proceso de revisión que fue hecho posible gracias a la ayuda de árbitros anónimos.

La meta primaria de la serie de escuelas AGRA es la formación de estudiantes graduados y jóvenes investigadores en el proceso de especializarse en teoría de números, teoría de grupos, formas modulares y áreas afines.

El programa del AGRA II consistió de las siguientes unidades:

1. Formas modulares y curvas de Shimura
   En orden de intervención: Roberto Miatello (Universidad Nacional de Córdoba, Argentina), Gonzalo Tornaría (Universidad de la República, Uruguay), Ariel Pacetti (Universidad de Buenos Aires) y Michael Harris (Columbia University y Université de Paris VII).

2. Combinatoria aditiva
   Juanjo Rué (Freie Universität Berlin), Julia Wolf (University of Bristol, Reino Unido), Javier Cilleruelo (Universidad Autónoma de Madrid) y Pablo Candela (Instituto Rényi, Hungría), con la asistencia de Ana Zumalacárregui (University of New South Wales, Australia).

3. Introducción a la teoría de las curvas elípticas
   Marusia Rebolledo (Université Blaise Pascal, Clermont-Ferrand), Marc Hindry (Université de Paris VII).

4. Crecimiento en grupos y expansores
   Harald Andrés Helfgott (Georg-August-Universität Göttingen y Centre Na-

tional de la Recherche Scientifique/Université de Paris VI/VII) y Mikhail Belolipetsky (Instituto de Matemática Pura e Aplicada, Rio de Janeiro).

5. Análisis y geometría en grupos
   Andrzej Zuk (Université de Paris VII).

El trabajo dentro de cada unidad fue coordinado cuidadosamente. En lo que se refiere a la coordinación entre unidades, se trató de evitar las redundancias excesivas y de apreciar y utilizar las que eran menores.

Aparte de las unidades, hubo dos charlas de Ricardo Menares (Pontificia Universidad Católica de Valparaíso), cuyas notas se encuentran en el presente volúmen. Hubo también una charla por Fernando Rodríguez Villegas (ICTP – The Abdus Salam International Center for Theoretical Physics, Trieste), y una serie de dos charlas de Carlos Moreno (City University of New York), sobre valores de funciones zeta y funciones $L$.

Se agradece sinceramente el apoyo tanto financiero como logístico de ICTP, CIMPA (Centre International de Mathématiques Pures et Appliquées, Nice), Concytec/Fondecyt/ANC (Consejo Nacional de Ciencia, Tecnología e Innovación Tecnológica/Fondo Nacional de Desarrollo Científico, Tecnología e Innovación Tecnológica/Academia Nacional de Ciencias del Perú), IMPA (Instituto de Matemática Pura e Aplicada, Rio de Janeiro), y la institucion sede, UNSAAC (Universidad de San Antonio Abad del Cusco). Varias otras instituciones han contribuído a cubrir los costos de la asistencia de sus estudiantes y docentes.

El comité científico consistió de H. A. Helfgott, Carlos Gustavo Moreira (IMPA) y Gonzalo Tornaría (Universidad de la República). El comité nacional estuvo compuesto de Roxana Lopez Cruz (Universidad Nacional Mayor de San Marcos, Lima, Perú), Fernando Rodríguez Villegas (ICTP) y Oswaldo Velásquez (Instituto de Matemática y Ciencias Afines, Lima), mientras que el comité local fue formado por Guido Álvarez Jáuregui, Georgina Cruz Quin, Katia García Alfaro, Alejandro Ttito Ttica, todos de UNSAAC, asistidos por un gran número de profesores y estudiantes locales. Las notas se han visto suplementadas por videos de algunas charlas, tomados por personal local y editados por ICTP. Los videos están actualmente disponibles en YouTube.

Se agradece sincera y humildamente a los constructores del Cusco y del Valle Sagrado por haber proporcionado, con mucha antelación, un marco excepcional para la escuela.

Harald Andrés Helfgott
Cusco, agosto de 2015 – Gotinga, diciembre de 2018

# Índice

# 1 Formas modulares y curvas de Shimura

# FORMAS MODULARES Y OPERADORES DE HECKE

Roberto J. Miatello

**Resumen**

Las presentes notas contienen una introducción a las formas modulares desde el punto de vista clásico. En ellas se dan las principales propiedades enfatizando el rol de los coeficientes de Fourier. Como principales ejemplos se estudian las series de Eisenstein y las series theta. Se describe en detalle la estructura del álgebra de formas modulares en el caso del grupo modular $\Gamma(1) = SL_2(\mathbb{Z})$. Finalmente se estudian los operadores de Hecke, describiendo sin demostraciones la generalización a los grupos de Hecke $\Gamma_0(N)$. Se estudian las propiedades de las autoformas de Hecke y de las funciones $L$ asociadas, las cuales admiten productos de Euler. El objetivo de esta parte es servir de base a los aspectos más avanzados de la teoría a ser presentados en las exposiciones sobre curvas de Shimura.

## 1 Grupos Fuchsianos

Una superficie de Riemann es una variedad diferenciable compleja conexa de dimensión compleja 1. De acuerdo a un teorema clásico de Riemann–Koebe, hay sólo tres superficies de Riemann simplemente conexas: el plano complejo $\mathbb{C}$, el semiplano superior $H$ y la esfera de Riemann $\hat{\mathbb{C}} \simeq \mathbb{C}P^1$.

El semiplano superior $H$ es el cubrimiento universal de la 'gran mayoría' de las superficies de Riemann, en particular de todas las superficies compactas con género $g \geqslant 2$.

Los grupos $SL_2(\mathbb{R})$ y $GL_2(\mathbb{R})^+$ actúan transitivamente en $H$ por transformaciones de Moebius y por lo tanto se tienen los difeomorfismos

$$H \simeq SL_2(\mathbb{R})/SO(2,\mathbb{R}) \simeq GL_2(\mathbb{R})^+/\mathbb{R}^\times SO(2,\mathbb{R}),$$

dado que $SO(2,\mathbb{R})$ (respectivamente $\mathbb{R}^\times SO(2,\mathbb{R})$) es el subgrupo de isotropía de $i$ en $SL_2(\mathbb{R})$ (resp. $GL_2(\mathbb{R})^+$).

Además $H$ tiene la métrica Riemanniana definida por $g = \frac{dx^2+dy^2}{y^2}$ y el operador de Laplace dado por $y^2(\frac{\partial^2}{\partial x^2} + \frac{\partial^2}{\partial y^2})$. Las geodésicas en la métrica hiperbólica son semicírculos y rectas perpendiculares al eje real.

Un *grupo Fuchsiano de primera clase* es un subgrupo discreto $\Gamma \subset SL_2(\mathbb{R})$ que tiene covolumen finito o equivalentemente $\mathrm{vol}(\Gamma\backslash H) < \infty$. Estos grupos fueron introducidos y estudiados por Henri Poincaré (1882) y son sumamente interesantes. Poincaré los llamó grupos Fuchsianos a raíz de investigaciones de Fuchs (1880).

Una *región fundamental* $\mathfrak{F}$ de $\Gamma$ es un subconjunto abierto de $H$ tal que

(i) cada $z \in H$ es equivalente por algún $\gamma \in \Gamma$ a algún $z' \in \overline{\mathfrak{F}}$,

(ii) si dos elementos $z, z' \in \mathfrak{F}$ son $\Gamma$-equivalentes, entonces $z, z' \in \partial\mathfrak{F}$.

Se prueba que todo grupo Fuchsiano de primera clase tiene una región fundamental $\mathfrak{F}$ que es un polígono hiperbólico con finitos lados. Los ejemplos típicos de tales grupos son $SL_2(\mathbb{Z})$ y los llamados subgrupos de congruencia que definimos a continuación.

**Ejemplo 1.1.** Si $N \in \mathbb{N}$, sean

(1-1) $$\Gamma(N) = \{\gamma \in SL_2(\mathbb{Z}) : \gamma \equiv Id \mod N\},$$

(1-2) $$\Gamma_0(N) = \{\gamma = \left[\begin{smallmatrix} a & b \\ c & d \end{smallmatrix}\right] \in SL_2(\mathbb{Z}) : c \equiv 0 \mod N\},$$

el *subgrupo principal de congruencia de nivel* $N$ y el subgrupo de Hecke de nivel $N$, respectivamente.

Si $N = 1$, $\Gamma(N)$ y $\Gamma_0(N)$ coinciden con $\Gamma = SL_2(\mathbb{Z})$. Si $N = 2$, se tiene

$$\Gamma(2) = \{\gamma = \left[\begin{smallmatrix} a & b \\ c & d \end{smallmatrix}\right] \in SL_2(\mathbb{Z}) : b, c \equiv 0 \mod 2\},$$

un subgrupo de índice 6 de $SL_2(\mathbb{Z})$; $\Gamma_0(2)$ es un subgrupo de índice 2.

La región fundamental usual de $\Gamma$ es $\mathfrak{F} = \{z : |Re(z)| < \frac{1}{2}, |z| > 1\}$. Para obtener regiones fundamentales $\mathfrak{F}_1$ y $\mathfrak{F}_2$ para $\Gamma(2)$ y $\Gamma_0(2)$ respectivamente, conviene en primer lugar encontrar elementos que representen las coclases a izquierda (o a derecha) de tales grupos en $\Gamma$, teniendo en cuenta que $\Gamma$ está generado por los elementos $S = \left[\begin{smallmatrix} 0 & 1 \\ -1 & 0 \end{smallmatrix}\right]$ y $T = \left[\begin{smallmatrix} 1 & 1 \\ 0 & 1 \end{smallmatrix}\right]$. Luego, hay que trasladar $\mathfrak{F}$ por los elementos obtenidos (ejercicio), resultando regiones $\mathfrak{F}_1$ y $\mathfrak{F}_2$ que son polígonos hiperbólicos con 4 y 3 lados (es decir segmentos de geodésicas hiperbólicas), respectivamente. Los puntos de intersección con el eje real y el punto ideal $i_\infty$ son puntos fuera de $H$ que deben ser agregados a $\Gamma\backslash H$ para obtener una compactificación de este espacio. Se corresponden con las llamadas cúspides de $\Gamma$.

**Ejemplo 1.2.** Los ejemplos anteriores tienen regiones fundamentales no compactas. Un ejemplo de distinta naturaleza es el siguiente. Sean $n, p \in \mathbb{N}$, con $p$ primo tal que $n$ no es un cuadrado mod $p$. Sea

(1-3) $$\Gamma(n,p) = \left\{\left[\begin{smallmatrix} a+b\sqrt{n} & (c+d\sqrt{n})p \\ (c-d\sqrt{n})p & a-b\sqrt{n} \end{smallmatrix}\right] \in SL(2\mathbb{R}) : a^2 - b^2n - c^2p + d^2np = 1\right\}.$$

Se verifica que $|\mathrm{tr}\,\gamma| > 2$ para todo $\gamma \in \Gamma$ (ejercicio), lo que implica que $\Gamma(n,p)\backslash H$ es una superficie compacta. Estos grupos provienen de álgebras de cuaterniones y tendrán un rol central en las exposiciones sobre curvas de Shimura.

**1.1 Cúspides.** Cada $g \in SL_2(\mathbb{R})$ tiene uno o dos puntos fijos en $\mathbb{R} \cup i_\infty$; $g$ se dice *parabólica* si tiene un único punto fijo, *hiperbólica* si tiene dos puntos fijos y *elíptica* si tiene dos puntos fijos conjugados, $w, \bar{w}$, no reales. Dejamos como ejercicio verificar que $g = \left[\begin{smallmatrix} a & b \\ c & d \end{smallmatrix}\right]$ es parabólica sii $|a+d| = 2$, hiperbólica sii $|a+d| > 2$ y elíptica sii $|a+d| < 2$.

Alternativamente, puede darse una interpretación geométrica: $g$ es parabólica sii es conjugada a una traslación $\left[\begin{smallmatrix} 1 & b \\ 0 & 1 \end{smallmatrix}\right]$, $b \in \mathbb{R}$, hiperbólica sii es conjugada a una dilatación $\left[\begin{smallmatrix} a & 0 \\ 0 & a^{-1} \end{smallmatrix}\right]$, $a \in \mathbb{R}^\times$ y elíptica sii es conjugada a una rotación $\left[\begin{smallmatrix} \cos\theta & \operatorname{sen}\theta \\ -\operatorname{sen}\theta & \cos\theta \end{smallmatrix}\right]$, $\theta \in \mathbb{R}$. El único punto fijo de una traslación es el punto ideal $i_\infty$.

Si $\gamma \in \Gamma$ es parabólico, $r \in \mathbb{R} \cup i_\infty$ su punto fijo y $\beta \in \Gamma$, entonces $r' = \beta \cdot r$ es fijo por $\beta\gamma\beta^{-1} \in \Gamma$ el cual es también parabólico. En este caso se dice que $r$ y $\beta \cdot r$ son vértices equivalentes. Luego $\Gamma$ actúa en el subconjunto $\mathcal{C} \subset \mathbb{R} \cup i_\infty$ de vértices parabólicos de $G$ y cada clase de $\Gamma$-equivalencia de $\mathcal{C}$ se llama una *cúspide* de $\Gamma$. Frecuentemente identificaremos a cada vértice parabólico $v$ con su clase de equivalencia y diremos simplemente que $v$ es una cúspide de $\Gamma$.

Se prueba que el número de cúspides de un grupo Fuchsiano de primera clase es finito. El espacio $\overline{H} = H \cup \mathcal{C}$ puede ser munido de una topología Hausdorff y $\Gamma$ actúa continuamente en $\overline{H}$ por transformaciones de Möbius. Además, el cociente $\Gamma\backslash\overline{H}$ tiene estructura de superficie de Riemann compacta, obtenida por adición a $\Gamma\backslash H$ de las cúspides de $\Gamma$, ver Shimura (1971).

## 2 Formas modulares

Sea, para $g = \left[\begin{smallmatrix} a & b \\ c & d \end{smallmatrix}\right] \in GL_2(\mathbb{R})$, $z \in H$, $j(g,z) = cz + d$. La función $j(g,z)$ es un cociclo, esto es, satisface la relación $j(g_1 g_2, z) = j(g_1, g_2 z) j(g_2, z)$ para todo $g_1, g_2 \in G$ y $z$ en $H$.

Definamos para $f : H \to \mathbb{C}$, $\delta \in SL_2(\mathbb{R})^+$ y $k \in \mathbb{N}$,

$$f|\delta(z) = j(\delta, z)^{-k} f(\delta z). \tag{2-1}$$

La propiedad del cociclo de $j(\gamma, z)$ se traduce en la identidad $f|\delta_1\delta_2 = f|\delta_1|\delta_2$, esto es, la correspondencia $f \to f|\delta$ es una acción a derecha.

**Definición 2.1.** Dado un grupo Fuchsiano de primera clase $G$, una *forma modular de peso $k$* para $G$ es una función holomorfa $f : H \to \mathbb{C}$ tal que

$$f(\gamma z) = j(\gamma, z)^k f(z) \tag{2-2}$$

para todo $\gamma \in \Gamma$, $z \in H$. Equivalentemente, en la notación en (2-1), (2-2) se traduce en la identidad $f|\delta = f$ para todo $\delta \in \Gamma$.

En segundo lugar, para ser forma modular se requiere que $f$ sea holomorfa en las cúspides de $\Gamma$.

La segunda condición en la definición de forma modular significa lo siguiente. Si $\Gamma$ contiene una traslación es decir un elemento de la forma $g = \left[\begin{smallmatrix} 1 & h \\ 0 & 1 \end{smallmatrix}\right]$, equivalentemente, si $i_\infty$ es una cúspide, entonces (2-2) implica que $f(z+h) = f(z)$ para todo $z \in H$. Por lo tanto se tiene que

$$f(z) = \sum_{-\infty}^{+\infty} a_n(f)\, e^{2\pi i \left(\frac{nz}{h}\right)}. \tag{2-3}$$

La holomorfía de $f$ en $i_\infty$ requiere que $a_n(f) = 0$ para todo $n < 0$. Los coeficientes $a_n(f)$ se llaman los *coeficientes de Fourier* de $f$ en $i_\infty$ y la expansión (2-3) se llama la *expansión de Fourier* de $f(z)$ en $i_\infty$. Para $SL_2(\mathbb{Z})$ y $\Gamma_0(2)$ se tiene que $h = 1$ y para $\Gamma = \Gamma(2)$, $h = 2$. Usaremos con frecuencia la notación corriente $q = e^{2\pi i z}$.

Para analizar los desarrollos de Fourier en las restantes cúspides se procede por reducción al caso anterior. Si $r$ es una cúspide de $\Gamma$ fija por el elemento parabólico $\gamma \in \Gamma$, sea $g \in SL_2(\mathbb{R})$ tal que $g \cdot r = i_\infty$. Se tiene entonces que $i_\infty$ es una cúspide del grupo $g\Gamma g^{-1}$ y el requerimiento de holomorfía de $f$ en $r$ consiste en pedir que la forma modular $f|g$ (asociada al grupo $g\Gamma g^{-1}$) sea holomorfa en $i_\infty$.

Si $a_0 = 0$ para cada cúspide de $\Gamma$ se dice que $f$ es una *forma cuspidal* de peso $k$.

Denotaremos por $\mathfrak{M}_k(\Gamma)$ al espacio de formas $\Gamma$-modulares de peso $k$ y por $\mathcal{S}_k(\Gamma)$ al subespacio de $\mathfrak{M}_k(\Gamma)$ de formas $\Gamma$-cuspidales. Si $\Gamma = SL_2(\mathbb{Z})$, escribiremos simplemente $\mathfrak{M}_k$ y $\mathcal{S}_k$.

Las formas modulares son objetos naturales que aparecen muy frecuentemente en matemática; se corresponden con secciones holomorfas de fibrados lineales naturales sobre la superficie de Riemann $\Gamma\backslash\overline{H}$, ver Shimura (ibíd.).

**2.1 Series de Eisenstein y series theta.** Entre los ejemplos más comunes de formas modulares, están las series de Eisenstein y las series theta.

*Series de Eisenstein.* Sea $\Gamma = SL_2(\mathbb{Z})$. Es fácil verificar que si el peso es impar la única forma modular es $f = 0$, luego supondremos que el peso es par e igual a $2k$ con $k \geqslant 2$. Sea

(2-4) $$G_{2k}(z) = {\sum}'_{(m,n)\in\mathbb{Z}^2} (mz+n)^{-2k}$$

donde ' indica omitir $(0,0)$ en la suma. Se prueba que la serie converge absoluta y uniformemente sobre compactos si $k \geqslant 2$ y por lo tanto define una función holomorfa en $H$. Ahora bien, para $\gamma \in \Gamma = SL_2(\mathbb{Z})$,

(2-5) $$\begin{aligned} G_{2k}(\gamma z) &= {\sum}'_{(m,n)\in\mathbb{Z}^2} (m\left(\tfrac{az+b}{cz+d}\right)+n)^{-2k} \\ &= (cz+d)^{2k} {\sum}'_{(m,n)\in\mathbb{Z}^2} (m(az+b)+n(cz+d))^{-2k} \\ &= (cz+d)^{2k} {\sum}'_{(m,n)\in\mathbb{Z}^2} ((ma+nc)z+(mb+nd))^{-2k} \\ &= (cz+d)^{2k} G_{2k}(z)\,. \end{aligned}$$

Para probar que $G_{2k}(z)$ es una forma modular de peso $2k$, resta analizar el desarrollo de Fourier de $G_{2k}$ en $i_\infty$. Hallaremos este desarrollo a continuación, calculando los coeficientes de Fourier. En este caso $\Gamma$ tiene a $i_\infty$ como única cúspide.

Recordemos el desarrollo

(2-6) $$\pi \cot \pi z = \frac{1}{z} + \sum_1^\infty \left( \frac{1}{z+m} + \frac{1}{z-m} \right)$$

con convergencia uniforme sobre compactos en $\mathbb{C} \smallsetminus \mathbb{Z}$. Por otra parte

$$\pi \cot \pi z = \pi i \frac{e^{2\pi i z}+1}{e^{2\pi i z}-1} = \pi i - 2\pi i \sum_{n=0}^\infty e^{2\pi i n z}\,.$$

Luego, derivando $k-1$ veces obtenemos para todo $k \geqslant 2$

(2-7) $$\sum_{m\in\mathbb{Z}} \frac{1}{(z+m)^k} = \frac{(-2\pi i)^k}{(k-1)!} \sum_{n=1}^\infty n^{k-1}\, e^{2\pi i n z}.$$

**Proposición 2.2.** *Se tiene, para todo $k \geqslant 2$,*

$$G_{2k}(z) = 2\zeta(2k) + \frac{2(2\pi i)^{2k}}{(2k-1)!} \sum_{n \geqslant 1} \sigma_{2k-1}(n)\, e^{2\pi i n z} \tag{2-8}$$

*donde se denota $\sigma_h(n) = \sum_{d|n} d^h$ para $h \in \mathbb{N}_0$ y donde $\zeta(s)$ es la función zeta de Riemann.*

*Prueba.* Usando (2-7) se obtiene,

$$G_{2k}(z) = \sum_{(n,m) \neq (0,0)} \frac{1}{(nz+m)^{2k}} \tag{2-9}$$

$$= 2\zeta(2k) + 2\sum_{n=1}^{\infty} \sum_{m \in \mathbb{Z}} \frac{1}{(nz+m)^{2k}} \tag{2-10}$$

$$= 2\zeta(2k) + \frac{2(2\pi i)^{2k}}{(2k-1)!} \sum_{d=1}^{\infty} \sum_{h=1}^{\infty} d^{2k-1}\, e^{2\pi i h d z} \tag{2-11}$$

$$= 2\zeta(2k) + \frac{2(2\pi i)^{2k}}{(2k-1)!} \sum_{n=1}^{\infty} \sigma_{2k-1}(n)\, e^{2\pi i n z}\,. \tag{2-12}$$

□

**Nota 2.3.** Se suele normalizar $G_{2k}$, poniendo $E_{2k} = (2\zeta(2k))^{-1} G_{2k}(z)$. En esta notación resulta que

$$E_4(z) = 1 + 240 \sum_{i=1}^{\infty} \sigma_3(n) q^n$$

$$E_6(z) = 1 - 504 \sum_{i=1}^{\infty} \sigma_5(n) q^n$$

$$E_8(z) = 1 + 480 \sum_{i=1}^{\infty} \sigma_7(n) q^n$$

donde $q = e^{2\pi i z}$. Con estas normalizaciones, teniendo en cuenta la estructura del álgebra de formas modulares para $SL_2(\mathbb{Z})$, a ser descripta en la próxima subsección (en particular que $\mathcal{S}_{2k} = 0$ si $k < 6$ y $\dim \mathfrak{M}_8 = \dim \mathfrak{M}_{10} = 1$), se obtienen las identidades

$$E_4^2 = E_8, \qquad E_6 E_4 = E_{10}$$

las cuales implican

$$\sigma_7(n) = \sigma_3(n) + 120 \sum_{i=1}^{m-1} \sigma_3(m)\sigma_3(n-m) \tag{2-13}$$

$$11\sigma_9(n) = 21\sigma_5(n) - 10\sigma_3(n) + 5040 \sum_{m=1}^{n-1} \sigma_3(n)\sigma_5(n-m)\,. \tag{2-14}$$

Para pesos mayores las relaciones son más complicadas. En general, dado $k \geqslant 2$, $k$ admite una expresión, o más, $k = 2\alpha + 3\beta$ con $\alpha, \beta \geqslant 0$ y en este caso se tiene que

$$E_{2k} = E_4^{\alpha} E_6^{\beta} + f$$

con $f \in \mathcal{S}_{2k}$.

**Nota 2.4.** Si $\Gamma$ es cualquier subgrupo de índice finito de $SL(2, \mathbb{Z})$, una manera alternativa más general de definir la serie de Eisenstein es

$$\widetilde{G}_h(z) = \sum_{\gamma \in \Gamma_\infty \backslash \Gamma} j(\gamma, z)^{-h}, \tag{2-15}$$

donde $\Gamma_\infty = \{\gamma \in \Gamma : \gamma i_\infty = i_\infty\}$, $h \geqslant 4$ y $j(\gamma, z) = cz + d$ si $\gamma = \left[\begin{smallmatrix} a & b \\ c & d \end{smallmatrix}\right]$.

En el caso $\Gamma = SL_2(\mathbb{Z})$, $h = 2k$, se tiene $\Gamma_\infty = \{\gamma = \pm \left[\begin{smallmatrix} 1 & m \\ 0 & 1 \end{smallmatrix}\right], m \in \mathbb{Z}\}$ y se verifica que $\Gamma_\infty \backslash \Gamma$ está parametrizado por $\{(c, d) \in \mathbb{Z}^2 : (c, d) = 1\}/\pm 1$.

Luego,

$$\widetilde{G}_{2k}(z) = \frac{1}{2} \sum_{(c,d) \in \mathbb{Z}^2 : (c,d) = 1} (cz + d)^{-2k}$$

y por otra parte

$$\begin{aligned} G_{2k}(z) &= \sum\nolimits'_{(c,d) \in \mathbb{Z}^2} (cz + d)^{-2k} \\ &= \sum_{n=1}^{\infty} \sum_{(c,d) = n} (cz + d)^{-2k} \\ &= \sum_{n=1}^{\infty} n^{-2k} \sum_{(c,d) = 1} (cz + d)^{-2k} = 2\zeta(2k) \widetilde{G}_{2k}(z). \end{aligned} \tag{2-16}$$

En particular $\widetilde{G}_{2k}(z) = E_{2k}(z)$, la serie anteriormente definida.

*Series theta.* A continuación pasamos a considerar otro ejemplo típico de forma modular, el caso de las llamadas *series theta*. Dado $L$ un retículo en $\mathbb{R}^n$ sea

$$\theta_L(z) = \sum_{\beta \in L} e^{\pi i \|\beta\|^2 z} \tag{2-17}$$

Usando que $|e^{\pi i \|\beta\|^2 z}| = e^{-\pi \|\beta\|^2 y}$ se demuestra que la serie converge uniformemente sobre compactos en $H$ y por lo tanto $\theta_L(z)$ es una función holomorfa. Veremos que, bajo ciertas condiciones, $\theta_L(z)$ satisface una relación de modularidad para $SL_2(\mathbb{Z})$.

Dado $L$, sea $L' = \{v \in \mathbb{R}^n : \beta . v \in \mathbb{Z}, \ \forall \beta \in L\}$, el retículo dual de $L$. Un retículo $L$ se dice *autodual* si $L = L'$ y se dice *par* si $\|\beta\|^2 \in 2\mathbb{Z}$ para todo $\beta \in L$. Notar que si $L = L'$ entonces $\mathrm{vol}(L) = 1$.

**Teorema 2.5.** *Si $n = 8k$, y $L = L'$ es par entonces $\theta_L(z) \in \mathfrak{M}_{4k} \smallsetminus \mathcal{S}_{4k}$.*

*Prueba.* Para verificar la condición (2-2) es suficiente hacerlo para $S(z) = -\frac{1}{z}$ y $T(z) = z + 1$ pues $S$ y $T$ generan $SL_2(\mathbb{Z})$. Como $L$ es par $\theta_L(z)$ es claramente $T$-invariante. Queda por verificar que para todo $z \in H$

$$\theta_L(-\tfrac{1}{z}) = z^{4k}\theta_L(z)\,.$$

que por holomorfía basta verificar para todo $z = it$ con $t > 0$. Es decir, hay que probar que para todo $t > 0$

$$\sum_{\beta\in L} e^{-\pi\|\beta\|^2/t} = t^{4k}\sum_{\beta\in L} e^{-t\pi\|\beta\|^2}\,. \tag{2-18}$$

Esta identidad es consecuencia de la conocida fórmula de sumación de Poisson dada por la identidad

$$\sum_{\beta\in L} f(\beta) = \mathrm{vol}(L)^{-1}\sum_{\beta'\in L'}\hat{f}(\beta'), \tag{2-19}$$

donde $f \in \mathcal{S}(\mathbb{R}^n)$, el espacio de Schwartz, $L$ es un retículo en $\mathbb{R}^n$, $L'$ es el retículo dual de $L$ y $\hat{f}(y) = \int_{\mathbb{R}^n} e^{2\pi i x.y} f(x)\,dx$ es la transformada de Fourier de $f$.

Si se toma la función $f(x) = e^{-\pi\|x\|^2}$ para $x \in \mathbb{R}^n$, es un hecho conocido que $f(x) \in \mathcal{S}(\mathbb{R}^n)$ satisface $\hat{f}(x) = f(x)$.

Para $t \in \mathbb{R}^+$ sea $L_t = \sqrt{t}L$. Entonces $(L_t)' = t^{-1/2}L$ y $\mathrm{vol}(L_t) = t^{4k}\mathrm{vol}(L)$. Aplicando directamente (2-19) con $f(x) = e^{-\pi\|x\|^2}$ y $L = L_{t^{-1}}$ se obtiene (2-18) de la cual resulta la modularidad de $\theta_L(z)$.

Observamos además que de la propia definición resulta que

$$\theta_L(z) = 1 + \sum_{m=1}^{\infty} r_{2m}(L)\,e^{2m\pi i z}, \tag{2-20}$$

donde $r_{2m}(L) = \{\beta \in L : \|\beta\|^2 = 2m\}$, el número de representaciones de $2m$ por la forma cuadrática $\|x\|^2$ en $\mathbb{R}^n$, lo cual da el desarrollo de Fourier de $\theta_L(z)$ en $i_\infty$. □

Los desarrollos de Fourier de las series de Eisenstein y de las series theta muestran que ambas definen formas no cuspidales e ilustran el hecho general de que los coeficientes de Fourier de las formas modulares suelen ser funciones aritméticas importantes. La acotación de los coeficientes de Fourier es también un problema clásico importante.

**2.2 Estructura de $\mathfrak{M}_{2k}$ y $\mathcal{S}_{2k}$.** Los espacios de formas modulares para $\Gamma = SL_2(\mathbb{Z})$ tienen una estructura algebraica muy precisa que describiremos en esta subsección. Para comenzar, observamos las siguientes afirmaciones que son de rápida verificación:

(i) Si $f \in \mathfrak{M}_{2k}, g \in \mathfrak{M}_{2h}$, se tiene que $fg \in \mathfrak{M}_{2(k+h)}$. En particular, el espacio $\mathfrak{M} := \sum_{m=0}^{\infty}\mathfrak{M}_{2k}$ es una $\mathbb{C}$-álgebra graduada.

(ii) Si $f, g \in M_{2k}$, existen $\lambda, \mu \in \mathbb{C}$ tales que $\lambda f + \mu g \in \mathcal{S}_{2k}$. Es decir $S_{2k}$ tiene codimensión 1 en $\mathfrak{M}_{2k}$.

(iii) $\mathfrak{M}_{2k} = \mathcal{S}_{2k} \oplus \mathbb{C}G_{2k}$.

Para continuar, necesitamos hacer uso de un resultado que es consecuencia del teorema de los residuos esencialmente (ver Serre (1973), Thm 3, p.85). Si $p \in H \cup i_\infty$, sea $v_p(f)$ el orden de $f$ en $p$ y como anteriormente, denotamos $\rho = \frac{1}{2}(1 + i\sqrt{3})$.

**Teorema 2.6.** *Si $f$ es una forma modular de peso $2k$, $f \neq 0$ entonces*

$$v_{i_\infty}(f) + \frac{1}{3}v_\rho(f) + \frac{1}{2}v_i(f) + \sum_{p \neq \rho, i_\infty} v_p(f) = \frac{k}{6}. \tag{2-21}$$

Éste es el hecho principal que permite describir la estructura del álgebra de las formas modulares.

Como ilustración, examinemos la ecuación (2-21) en los casos de $G_4, G_6$ y $G_8$. Para $G_4$ el miembro de la derecha da $\frac{1}{3}$, luego necesariamente $v_p(f) = 0$ para $p \neq \rho$ y $v_\rho(f) = 1$, es decir $G_4$ tiene un cero simple en $\rho$ y no posee otros ceros. Similarmente $G_6$ posee un único cero simple en $p = i$ y $G_8$ posee un único cero doble en $p = \rho$ al igual que $G_4^2$. Esto implica necesariamente que existe una constante $c$ tal que $G_8 - cG_4^2 = 0$ pues de lo contrario tendría más ceros que los permitidos por la ecuación (2-21). Con argumentos similares se deduce que $\mathfrak{M}_2 = 0$ y $\mathfrak{M}_{2k} = \mathbb{C}G_{2k}$ si $k = 2, 3, 4, 5$. Si $k = 6$ observamos, usando el desarrollo de Fourier, que la forma $\Delta = 60G_4^3 - 140G_6^2$ tiene $a_0 = 0$, es decir es una forma cuspidal no nula de peso 12 con un único cero en $i_\infty$ (necesariamente simple, por (2-21)). Es fácil ver que $M_{12}$ está generado por $\Delta$ y $G_4^3$ (ó equivalentemente por $\Delta$ y $G_6^2$, ó por $\Delta$ y $G_{12}$).

Consideremos, para $f \in \mathfrak{M}_{2k}$, $k \geqslant 1$, la transformación lineal $\Phi : f \to \Delta f \in \mathcal{S}_{2k+12}$. Claramente $\Phi$ es inyectiva pues $\Delta \not\equiv 0$. Por otra parte, si $g \in \mathcal{S}_{2k+12}$, se tiene que $f := g/\Delta$ es analítica en $i_\infty$ (pues $\Delta$ tiene en $i_\infty$ un cero simple), tiene peso $2k$ y $\Phi(f) = g$. Luego $\Phi$ es un isomorfismo. Resumiendo

**Proposición 2.7.** *Si $k \geqslant 2$, entonces*

*(i) La aplicación $\Phi : \mathfrak{M}_{2k} \to \mathcal{S}_{2k+12}$, dada por $\Phi(f) = \Delta f$ es un isomorfismo.*

*(ii)* $\dim(\mathfrak{M}_{2k}) = \begin{cases} [k/6] & \text{si } k \equiv 1 \mod 6, \\ [k/6] + 1 & \text{si } k \not\equiv 1 \mod 6. \end{cases}$

*(iii) $\mathfrak{M} \simeq \mathbb{C}[G_4, G_6]$ como $\mathbb{C}$-álgebra graduada.*

*Prueba.* La primera afirmación ya fue probada. La segunda es consecuencia de observar que la igualdad vale si $k \leqslant 6$ y que ambos miembros aumentan 1 al cambiar $k$ por $k + 6$. Esto sigue de (i) y de la igualdad $\mathfrak{M}_{2k} = \mathbb{C}G_{2k} \oplus \mathcal{S}_{2k}$.

Pasamos a la prueba de (iii) que también es consecuencia del Teorema 2.6. En primer lugar, veamos que los elementos de la forma ${G_4}^\alpha {G_6}^\beta$ con $k = 2\alpha + 3\beta$ generan $\mathfrak{M}_{2k}$. Fijo $k > 0$, existen $\alpha \geqslant 0$ y $\beta \geqslant 0$ tales que $k = 2\alpha + 3\beta$ (ejercicio), por lo tanto, existe $c \neq 0$ tal que $G_{2k} = c{G_4}^\alpha {G_6}^\beta$. Como $\mathfrak{M}_{2k} = \mathbb{C}G_{2k} \oplus \mathcal{S}_{2k}$, restaría probar que los elementos de la forma ${G_4}^\alpha\, {G_6}^\beta$ con $k = 2\alpha + 3\beta$ generan $\mathcal{S}_{2k}$. Esto se prueba por inducción. Es válido para $0 < k \leqslant 6$ y siendo válido para $\mathcal{S}_{2k-12}$, lo es también para $\mathcal{S}_{2k} = \Delta\mathcal{S}_{2k-12}$ pues $\Delta = 60{G_4}^3 - 140{G_6}^2$. Para concluir la prueba faltaría ver que $G_4$ y $G_6$ son algebraicamente independientes. Supongamos que

$$\sum_{2\alpha+3\beta=k} c_{\alpha,\beta}\, {G_4}^\alpha\, {G_6}^\beta = 0.$$

Si $G_6$ aparece en todos los términos con coeficientes no nulos, se lo puede simplificar y queda

$$\sum_{2\alpha+3(\beta-1)=k-3} c_{\alpha,\beta} {G_4}^{\alpha} {G_6}^{\beta-1} = 0.$$

Por hipótesis inductiva, $c_{\alpha,\beta} = 0$ para todo $\alpha, \beta$, una contradicción. Si $G_6$ no aparece en algún término, o sea $c_{\alpha,0} \neq 0$, evaluando en $z = i$ se obtiene $c_{\alpha,0} = 0$, un absurdo. □

**Definición 2.8.** Sea la función $j(z) = 1728 G_4(z)^3/\Delta(z)$; $j(z)$ es una forma modular de peso 0 con un único polo simple en $i_\infty$.

Un hecho importante es que $j(z)$ define una biyección de $\Gamma\backslash\overline{H}$ en $\hat{\mathbb{C}}$. Notemos que $\Delta(z)$ tiene un único cero simple en $i_\infty$, luego $j(z)$ tiene un polo simple en $i_\infty$ y no tiene otros polos. Siendo cociente de formas de peso 12, $j(z)$ es una forma modular de peso 0, con un único cero (de orden 3) en $z = \rho$. Según (2-21) tenemos $-1 + 3.\frac{1}{3} = 0$ y no puede haber otros ceros.

Si planteamos para $\alpha \in \mathbb{C}$ la ecuación $j(z) - \alpha = 0$ (2-21) nos da

(2-22)
$$\begin{aligned} -1+1 &= 0 \text{ si } \alpha = j(u), u \neq i, \rho, \\ -1+2\frac{1}{2} &= 0 \text{ si } \alpha = j(i), \\ -1+3\frac{1}{3} &= 0 \text{ si } \alpha = 0, z = \rho. \end{aligned}$$

Luego $j : \Gamma\backslash H \to \mathbb{C}$ es una biyección, excepto que 0 se toma con multiplicidad 3 y $j(i)$ con multiplicidad 2. Esta situación es coincidente con la estructura de superficie de Riemann de $\Gamma\backslash H$ en entornos de $\rho$ e $i$ respectivamente, luego $j$ induce la biyección que se afirma.

**Proposición 2.9.** *Si $f$ es meromorfa en $H$, las siguientes condiciones son equivalentes.*

*(i)* *$f$ es una función $\Gamma$-modular de peso 0.*

*(ii)* *$f$ es cociente de dos formas $\Gamma$-modulares del mismo peso.*

*(iii)* *$f$ es una función racional de $j$.*

Dejamos la prueba como ejercicio para el lector.

**Nota 2.10.** El desarrollo de Fourier de $j(z)$ tiene la forma $j(z) = q^{-1} + 744 + \sum_{n=1}^{\infty} a_n(j) q^n$, donde los coeficientes $a_n(j)$ son números enteros conectados con las dimensiones de representaciones irreducibles de ciertos grupos finitos simples (ver moonshine (Conway–Norton, Borcherds)).

**Nota 2.11.** Las formas modulares $G_4$ y $G_6$ tienen un rol importante en la teoría de curvas elípticas.

En primer lugar, cada retículo $L$ en $\mathbb{C}$ define un toro complejo $T_L$ y dos toros complejos $T_L$, $T_{L'}$ asociados a $L$ y $L'$ son biholomorfos si y sólo si $L' = zL$ para algún $z \in \mathbb{C}$. Alternativamente, si se toma $L = \mathbb{Z} \oplus \mathbb{Z}\tau$, $L' = \mathbb{Z} \oplus \mathbb{Z}\tau'$ con $\tau, \tau' \in H$, es fácil ver que $T_L$ y $T_{L'}$ son biholomorfos si y sólo si $\tau' = \gamma\tau$ para algún $\gamma \in SL_2(\mathbb{Z})$.

En efecto, como $\tau' = a\tau + b$ y $1 = c\tau + d$, con $\left[\begin{smallmatrix} a & b \\ c & d \end{smallmatrix}\right] \in SL_2(\mathbb{Z})$ se tiene que $\tau' = \frac{a\tau+b}{c\tau+d}$.

En otras palabras, el espacio cociente $\Gamma\backslash H$ parametriza las estructuras complejas en el toro $T^2$ salvo biholomorfía; es el llamado *espacio de móduli* de superficies de Riemann de género 1 (i.e curvas elípticas).

Paralelamente, para cada $L$ existe una incrustación $\Psi : T_L \to \mathbb{C}P(2)$, definido por la correspondencia $z \to [(\wp(z), \wp'(z), 1)]$, si $z \notin L$ y $z \to [(0, 1, 0)]$ si $z \in L$, donde $\wp(z)$ es la función $\wp$ de Weierstrass asociada al retículo $L = \mathbb{Z} \oplus \mathbb{Z}\tau$. Como es bien conocido, $\wp(z)$ y $\wp'(z)$ satisfacen la ecuación $\wp(z)'^2 = 4\wp(z)^3 - g_2(\tau)\wp(z) - g_3(\tau)$, donde $g_2(\tau) = 60G_4(\tau)$ y $g_3(\tau) = 140G_6(\tau)$ con $\Delta(\tau) = g_2(\tau)^3 - 27g_3(\tau)^2 \neq 0$, para todo $\tau \in H$. La función $\Delta$ es llamada función discriminante, corresponde al discriminante de la curva elíptica.

Por consiguiente, la imagen del embedding $\Psi$ es la curva elíptica no singular de ecuación $y^2 = 4x^3 - g_2(\tau)x - g_3(\tau)$. Inversamente, se prueba que para toda curva elíptica no singular $E$ de ecuación $y^2 = 4x^3 - c_2x - c_3$ existe un retículo $L = \mathbb{Z} \oplus \mathbb{Z}\tau$ tal que $c_2 = g_2(\tau)$ y $c_3 = g_3(\tau)$. Precisamente en la prueba de este hecho se usa la suryectividad de la función $j(z)$ definida anteriormente.

**Nota 2.12.** Aprovechamos para hacer mención de la famosa conjetura de Taniyama: si $E$ es una curva elíptica sobre $\mathbb{Q}$, entonces existe una autoforma normalizada de peso 2 para $\Gamma_0(N)$, donde $N$ es el conductor de $E$, tal que $L(E, s) = L_f(s)$. Una descripción de la misma fue sugerida por Taniyama en los años 50, promovida por Shimura en los años 60 y por Weil (1967) quien dio fuertes evidencias de su validez. Como es posible dar una lista de autoformas normalizadas de peso 2 para $\Gamma_0(N)$ para cada $N$ fijo, la conjetura predice cuántas curvas elípticas hay con conductor N sobre $\mathbb{Q}$. Búsquedas con computadora han confirmado este número para pequeños valores de $N$.

Como se sabe, de la conjetura y se sigue el último teorema de Fermat por resultados previos de Ribet. La misma fue probada para la mayoría de las curvas elípticas por Wiles, Taylor y Diamond (1995), posteriormente extendida a toda curva elíptica no singular sobre $\mathbb{Q}$, por Breuil, Conrad, Diamond y Taylor en 2001.

Actualmente este resultado está comprendido en el programa de Langlands, que predice que todas las series de Dirichlet provenientes de variedades algebraicas o más generalmente, de motivos, se realizan como funciones $L$ de representaciones automorfas de grupos algebraicos reductivos.

# 3 Operadores de Hecke

El objetivo de esta sección es estudiar los operadores de Hecke actuando en $\mathfrak{M}_{2k}$, con $\Gamma = SL_2(\mathbb{Z})$. Probaremos que estos operadores generan un álgebra conmutativa $\mathcal{H}$ de operadores autoadjuntos en $\mathcal{S}_{2k}$. La familia se puede diagonalizar simultáneamente y las autofunciones comunes de tales operadores, llamadas *autoformas* de Hecke, son formas modulares con propiedades especiales, así como sus autovalores.

Como se dijo, por su mayor simplicidad, desarrollaremos en detalle la teoría en el caso de $SL_2(\mathbb{Z})$ y al final describiremos sin demostración los cambios necesarios para su generalización a los subgrupos de Hecke $\Gamma_0(N)$, resultados debidos a Atkin y Lehner (1970).

**Definición 3.1.** Si $n \in \mathbb{N}$ sea $\mathbb{M}_n = \{g \in M_2(\mathbb{Z}) : det(g) = n\}$. Si $f \in \mathfrak{M}_{2k}$, se define el operador de Hecke $T_n$ por

(3-1) $$T_n(f) = n^{2k-1} \sum_{\delta \in SL_2(\mathbb{Z})\backslash \mathbb{M}_n} f|\delta$$

donde $\delta$ recorre un sistema completo de representantes de $\mathbb{M}_n$ módulo $SL_2(\mathbb{Z})$.

Claramente $T_n(f)$ no depende del sistema de representantes pues $f|\gamma\delta = f|\delta$ si $\gamma \in SL_2(\mathbb{Z})$, es holomorfa en $H$ y además, si $\gamma \in SL_2(\mathbb{Z})$,

$$T_n(f)|\gamma = n^{2k-1} \sum_{\delta \in SL_2(\mathbb{Z})\backslash \mathbb{M}_n} f|\delta\gamma = T_n(f).$$

ya que si $\delta_1, \dots, \delta_r$ es un sistema de representantes de $SL_2(\mathbb{Z})\backslash \mathbb{M}_n$, $\delta_1\gamma, \dots, \delta_r\gamma$ es también un sistema de representantes de $SL_2(\mathbb{Z})\backslash \mathbb{M}_n$.

Ahora, dada $\delta \in \mathbb{M}_n$ existe $\gamma \in SL_2(\mathbb{Z})$ tal que $\gamma\delta = \left[\begin{smallmatrix} a' & b' \\ 0 & d' \end{smallmatrix}\right]$ y además, multiplicando a izquierda por $\pm \left[\begin{smallmatrix} 1 & r \\ 0 & 1 \end{smallmatrix}\right]$ se puede cambiar $b'$ por $b' + dr$. Luego, si tomamos las matrices triangulares tales que $ad = n$, $a, d > 0$ y $0 \leqslant b < d$, es fácil verificar que tales matrices constituyen un sistema completo de representantes de $SL_2(\mathbb{Z})\backslash \mathbb{M}_n$. Por lo tanto obtenemos la expresión

(3-2) $$T_n(f)(z) = n^{2k-1} \sum_{ad=n,\, ad>0} \sum_{b=0}^{d-1} d^{-2k} f\left(\frac{az+b}{d}\right)$$

Para concluir que $T_n(f) \in M_{2k}$ resta hallar la expansión de Fourier en $i_\infty$.

**Proposición 3.2.** *Sea $f \in \mathfrak{M}_{2k}$, $f(z) = \sum_{m=0}^{\infty} a_m(f) q^m$. Entonces*

(3-3) $$T_n(f)(z) = \sum_{m=0}^{\infty} \left( \sum_{d|(m,n)} d^{2k-1} a_{\frac{nm}{d^2}}(f) \right) q^m \text{ , y}$$

(3-4) $$a_m(T_n(f)) = \sum_{d|(m,n)} d^{2k-1} a_{\frac{nm}{d^2}}(f)$$

*para cada $m \geqslant 0$.*

*Prueba.* Por (3-2) se tiene

$$\begin{aligned} T_n(f)(z) &= n^{2k-1} \sum_{ad=n,\, ad>0} \sum_{b=0}^{d-1} d^{-2k} f\left(\frac{az+b}{d}\right) \\ &= n^{2k-1} \sum_{ad=n,\, ad>0} \sum_{b=0}^{d-1} d^{-2k} \left( \sum_{m \geqslant 0} a_m \, e^{2\pi i m \left(\frac{az+b}{d}\right)} \right) \\ &= n^{2k-1} \sum_{ad=n,\, ad>0} d^{-2k+1} \left( \sum_{m' \geqslant 0} a_{m'd} \, e^{2\pi i m' a z} \right) \\ &= \sum_{m'' \geqslant 0} \left( \sum_{a|(n,m'')} a^{2k-1} a_{\frac{m''n}{a^2}} \right) q^{m''} . \end{aligned}$$

En la tercera identidad se han cancelado los términos tales que $d \nmid m$, pues la suma de $e^{2\pi i m \frac{b}{d}}$ entre $b = 0$ y $b = d - 1$ es nula. A continuación se ha reemplazado $m'' = m'a$. $\square$

En particular,

$$a_0(T_n(f)) = \left(\sum_{a|n} a^{2k-1}\right) a_0(f) = \sigma_{2k-1}(n) a_0(f). \tag{3-5}$$

**Corolario 3.3.** *Sea $n \in \mathbb{N}$ y $2k \geqslant 0$. Entonces $T_n$ preserva $\mathfrak{M}_{2k}$ y $\mathcal{S}_{2k}$.*

**Corolario 3.4.** *Si $p$ es primo,*

$$a_m(T_p(f)) = \begin{cases} p^{2k-1} a_{\frac{m}{p}}(f) + a_{mp}(f), & si\ p \mid m\,, \\ a_{mp}(f), & si\ p \nmid m \end{cases}. \tag{3-6}$$

*Bajo la convención de que $a_{\frac{m}{p}}(f) = 0$ si $p \nmid m$, (3-6) se resume en*

$$a_m(T_p(f)) = p^{2k-1} a_{\frac{m}{p}}(f) + a_{mp}(f). \tag{3-7}$$

**Proposición 3.5.** *Si $(n, m) = 1$ se tiene que $T_n T_m = T_m T_n = T_{nm}$.*

*Prueba.* Dados $n, m$ con $(n, m) = 1$ probaremos que para todo $r$, es $a_r(T_n T_m(f)) = a_r(T_{nm}(f))$.

En primer lugar, según se probó, $a_r(T_m(f)) = \sum_{d|(r,m)} a_{rm/d^2}(f)$. Luego

$$\begin{aligned} a_r(T_n T_m(f)) &= \sum_{e|(r,n)} e^{2k-1} a_{rn/e^2}(T_m(f)) \\ &= \sum_{d|(m,rn/e^2)} \sum_{e|(r,n)} d^{2k-1} e^{2k-1} a_{rnm/e^2 d^2}(f) \\ &= \sum_{h|(nm,r)} h^{2k-1} a_{rnm/h^2}(f) = a_r(T_{nm}(f)). \end{aligned}$$

Dejamos al lector la verificación de la última igualdad. En el caso $n = p, m = q$ primos distintos, se tiene

$$\begin{aligned} a_m(T_p(T_q f)) &= a_{mp}(T_q(f)) + p^{2k-1} a_{m/p}(T_q(f)) \\ &= a_{mpq(f)} + q^{2k-1} a_{mp/q}(f) + p^{2k-1} a_{mq/p}(f) + p^{2k-1} q^{2k-1} a_{m/pq}(f) \end{aligned}$$

Notar que, en esta expresión, si $p$ y $q$ no dividen a $m$ queda simplemente $a_m(T_p(T_q f)) = a_m(T_{pq} f) = a_{mpq}(f)$. $\square$

**Proposición 3.6.** *Si $p$ es primo, entonces*

$$T_p T_{p^s} = T_{p^{s+1}} + p^{2k-1} T_{p^{s-1}}. \tag{3-8}$$

*Prueba.* Sea $f, g \in \mathfrak{M}_{2k}$. Entonces por (3-2)

$$T_{p^s} f(z) = p^{s(2k-1)} \sum_{0 \leqslant i \leqslant s} p^{-i2k} \sum_{0 \leqslant b < p^i} f\left(\frac{p^{s-i}z+b}{p^i}\right) \tag{3-9}$$

$$T_p g(z) = p^{(2k-1)} g(pz) + p^{-1} \sum_{0 \leqslant b' < p} g\left(\frac{z+b'}{p}\right)$$

Luego

$$\begin{aligned} T_p T_{p^s} f(z) &= p^{(s+1)(2k-1)} \sum_{0 \leqslant i \leqslant s} p^{-i2k} \sum_{0 \leqslant b < p^i} f\left(\frac{p^{s+1-i}z+b}{p^i}\right) \\ &+ p^{-1} p^{s(2k-1)} \sum_{0 \leqslant b' < p} \sum_{0 \leqslant i \leqslant s} p^{-i2k} \sum_{0 \leqslant b < p^i} f\left(\frac{p^{s-i}(z+b')+pb}{p^{i+1}}\right). \end{aligned} \tag{3-10}$$

Si hacemos $i = s$ en el segundo sumando obtenemos

$$p^{-1-s} \sum_{0 \leqslant b' < p} \sum_{0 \leqslant b < p^s} f\left(\frac{z+b'+pb}{p^{s+1}}\right) = p^{-1-s} \sum_{0 \leqslant b < p^{s+1}} f\left(\frac{z+b}{p^{s+1}}\right).$$

Si este término se agrega al primero se obtiene (3-9) con $s+1$ en lugar de $s$, es decir $T_{p^{s+1}} f(z)$. Los términos restantes son

$$p^{-1} p^{s(2k-1)} \sum_{0 \leqslant b' < p} \sum_{0 \leqslant i \leqslant s-1} p^{-i2k} \sum_{0 \leqslant b < p^i} f\left(\frac{p^{s-1-i}z+b+p^{s-1-i}b')}{p^i}\right). \tag{3-11}$$

Ahora, para cada $i$, los $p^{i+1}$ números $b + p^{s-1-i}b'$ con $0 \leqslant b < p^i$, $0 \leqslant b' < p$ recorren todas las clases módulo $p^i$, $p$ veces. Por ejemplo, si $i \leqslant \frac{p-1}{2}$ entonces $b + p^{s-1-i}b' \equiv b \mod p^i$ para todo $0 \leqslant b' < p$. Como $f(u+1) = f(u)$ se obtiene siempre el mismo valor.

Luego (3-11) es igual a

$$p^{s(2k-1)} \sum_{0 \leqslant i \leqslant s-1} p^{-i2k} \sum_{0 \leqslant b < p^i} f\left(\frac{p^{s-1-i}z+b)}{p^i}\right) = p^{2k-1} T_{p^{s-1}} f(z) \tag{3-12}$$

que es el segundo sumando de la ecuación (3-8). □

Se define el *álgebra de Hecke*, $\mathcal{H}$, como la $\mathbb{Q}$-álgebra generada por los operadores $T_n$, con $n \in \mathbb{N}_0$. De las dos proposiciones previas resulta que $\mathcal{H}$ es una $\mathbb{Q}$-álgebra conmutativa, generada por los operadores $T_p$ con $p$ primo. Más aun, cada $T_n$ es una combinación lineal entera de productos de operadores $T_p$ con $p$ primo.

A continuación veremos que los operadores de Hecke son autoadjuntos con respecto a un producto interno canónico en $\mathcal{S}_k$, Petersson (1939).

Recordemos que una medida riemanniana $SL(2,\mathbb{R})$-invariante en $H$ está dada por $\mu(z) = \frac{dxdy}{y^2} = -\frac{dz\wedge\overline{dz}}{\mathrm{Im}(z)^2}$. En efecto, si $g \in SL_2(\mathbb{R})$, $z \in H$ se tiene que $\mathrm{Im}(gz) = \mathrm{Im}(z)/|cz+d|^2$ y $g'(z) = (cz+d)^{-2}$. Luego

$$\frac{d(gz) \wedge \overline{d(gz)}}{\mathrm{Im}(gz)^2} = \frac{|g'(z)|^2 \, |cz+d|^4 \, dz \wedge \overline{dz}}{\mathrm{Im}(z)^2} = \frac{dz \wedge \overline{dz}}{\mathrm{Im}(z)^2}.$$

**Definición 3.7.** Se define el producto interno de Petersson en $L^2(\Gamma\backslash H)$ por

$$\langle f, g\rangle = \int_{\Gamma\backslash H} f(z)\overline{g(z)}\, y^{2k}\, \frac{dxdy}{y^2}. \tag{3-13}$$

Observamos que el integrando en (3-13) es una función $\Gamma$-invariante, ya que

$$f(\gamma z)\overline{g(\gamma z)}\,\mathrm{Im}(\gamma z)^{2k} = |cz+d|^{4k}\,\frac{\mathrm{Im}(z)^{2k}}{|cz+d|^{4k}}\, f(z)\overline{g(z)} = \mathrm{Im}(z)^{2k} f(z)\overline{g(z)}.$$

El producto interno $\langle f, g\rangle$ está definido para $f, g \in \mathfrak{M}_{2k}(\Gamma)$ si al menos una de ellas es cuspidal. La convergencia es consecuencia de la siguiente acotación.

**Lema 3.8.** *Si $f \in \mathcal{S}_{2k}$, existe $C > 0$ tal que $f$ satisface $|f(z)| \leqslant C|Im(z)|^{-k}$ para todo $z \in H$. Además*

$$|a_n(f)| \leqslant O(n^k)\,. \tag{3-14}$$

*Prueba.* La función $|f(z)y^k|$ es $\Gamma$-invariante en $H$. Por ser $f$ cuspidal, existe $M > 0$ tal que $|f(z)y^k| \leqslant M$ para todo $z \in H$. O sea $|f(z)| \leqslant My^{-k}$ para todo $z \in H$. Por lo tanto

$$|a_n(f)| = \left|\int_0^1 f(x+iy)e^{2\pi inz}\,dx\right| \leqslant My^{-k}e^{2\pi ny}$$

para todo $y > 0$. Tomando $y = 1/n$ resulta (3-14). □

**Teorema 3.9.** *$T_n$ es un operador autoadjunto en $(\mathcal{S}_{2k}, \langle\cdot, \cdot\rangle)$.*

*Prueba.* Se puede asumir que $n = p$ es primo.

Definamos $\alpha = \left[\begin{smallmatrix}1 & 0\\ 0 & p\end{smallmatrix}\right]$, $\beta = \left[\begin{smallmatrix}p & 0\\ 0 & 1\end{smallmatrix}\right]$, $\alpha_j = \left[\begin{smallmatrix}1 & j\\ 0 & 1\end{smallmatrix}\right]$ $(0 \leqslant j \leqslant p-1)$, $\alpha\alpha_j = \begin{bmatrix}1 & j\\ 0 & p\end{bmatrix}$, $\gamma = \left[\begin{smallmatrix}0 & -1\\ 1 & 0\end{smallmatrix}\right]$ y $\nu = \begin{bmatrix}\sqrt{p} & 0\\ 0 & \sqrt{p}\end{bmatrix}$. Se tienen las relaciones $\gamma\beta = \alpha\gamma = \nu\delta$, $\alpha\alpha_j = \nu\delta_j$, con $\delta_j = \begin{bmatrix}1/\sqrt{p} & j/\sqrt{p}\\ 0 & \sqrt{p}\end{bmatrix}$ y $\delta = \begin{bmatrix}0 & -1/\sqrt{p}\\ \sqrt{p} & 0\end{bmatrix} \in SL(2,\mathbb{R})$.

Sabemos que el conjunto de $p+1$ elementos $\{\alpha\alpha_j,\ 0 \leqslant j \leqslant p-1;\ \alpha\gamma\}$ forman un sistema de representantes de $\Gamma\backslash\mathbb{M}_p$ y similarmente el conjunto $\{\gamma\alpha\alpha_j = \beta\gamma\alpha_j,\ 0 \leqslant j \leqslant p-1;\ \beta\}$. Luego tenemos las siguientes expresiones

$$T_pf = p^{2k-1}\sum_{j=0}^{p-1} f|\alpha\alpha_j + f|\alpha\gamma\ , \qquad T_pf = p^{2k-1}\sum_{j=0}^{p-1} f|\beta\beta_j + f|\beta. \tag{3-15}$$

Usando la primera expresión de $T_p$ y la invariancia de la medida Riemanniana $d\mu(z)$ tenemos que la primera suma en $\langle T_p f, g\rangle$ es igual a

$$p^{2k-1}\sum_{j=0}^{p-1}\int_{\Gamma\backslash H} j(\nu\delta_j, z)^{-2k}\, f(\nu\delta_j z)\,\overline{g(z)}\,\mathrm{Im}(z)^{2k}\, d\mu(z) \tag{3-16}$$
$$= p^{2k-1}\sum_{j=0}^{p-1}\int_{\Gamma\backslash H} p^{-2k} f(z)\,\overline{g(\delta_j{}^{-1}z)}\,\mathrm{Im}(\delta_j{}^{-1}z)^{2k}\, d\mu(z)$$
$$= p^{2k}\int_{\Gamma\backslash H} f(z)\,\overline{g(pz)}\,\mathrm{Im}(z)^{2k}\, d\mu(z)$$

En la última igualdad hemos usado que $j(\delta_j^{-1}, z) = p^{-1/2}$ y $g(\delta_j^{-1}z) = g(pz-j) = g(pz)$. Similarmente vemos que $p^{2k-1}\langle f|\alpha, g\rangle$ está dado por

$$p^{2k-1}\int_{\Gamma\backslash H} j(\alpha\gamma, z)^{-2k}\, f(\alpha\gamma z)\,\overline{g(z)}\,\mathrm{Im}(z)^{2k}\, d\mu(z) \tag{3-17}$$
$$= p^{2k-1}\int_{\Gamma\backslash H} (\gamma^{-1}z)^{-2k} p^{-2k}\, f(\alpha z)\,\overline{g(\gamma^{-1}z)}\,\mathrm{Im}(\gamma^{-1}z)^{2k}\, d\mu(z)$$
$$= p^{2k-1}\int_{\Gamma\backslash H} z^{2k} p^{-2k}\, f(\delta z)\,\overline{z}^{2k}\overline{g(z)}\,\mathrm{Im}(z)^{2k}|z|^{-4k}\, d\mu(z)$$
$$= p^{2k-1}\int_{\Gamma\backslash H} p^{-2k}\, f(z)\,\overline{g(\delta^{-1}z)}\,\mathrm{Im}(\delta^{-1}z)^{2k}\, d\mu(z)$$
$$= p^{2k-1}\int_{\Gamma\backslash H} f(z)\,\overline{g(pz)}\,\mathrm{Im}(z)^{2k}\, d\mu(z)$$

el cual es igual a $p^{2k-1}\langle f, g|\beta\rangle$, el segundo término de $\langle f, T_p g\rangle$. Para el cálculo de la primera suma en $\langle f, T_p g\rangle$ usamos la segunda expresión en (3-15).

$$\langle f, T_p g\rangle = p^{2k-1}\sum_{j=0}^{p-1}\int_{\Gamma\backslash H} f(z)\,\overline{j(\beta\beta_j, z)^{-2k} g(\beta\beta_j z)}\,\mathrm{Im}(z)^{2k}\, d\mu(z)$$
$$= p^{2k-1}\sum_{j=0}^{p-1}\int_{\Gamma\backslash H} f(\beta_j{}^{-1}z)\,\overline{j(\beta_j, \beta_j{}^{-1}z)^{-2k} g(\beta z)}\,\mathrm{Im}(\beta_j{}^{-1}z)^{2k}\, d\mu(z)$$
$$= p^{2k}\int_{\Gamma\backslash H} f(z)\,\overline{g(pz)}\,\mathrm{Im}(z)^{2k}\, d\mu(z)$$

donde hemos usado que $j(\beta_j{}^{-1}, z)^{2k}\,\overline{j(\beta_j, \beta_j{}^{-1}z)^{-2k}} = |j(\beta_j{}^{-1}, z)|^{4k}$ y que $\beta z = pz$. Esto prueba el teorema, ya que se ha obtenido para $\langle f, T_p g\rangle$ la misma expresión que en (3-16). □

**Corolario 3.10.** *Para cada $k \geqslant 1$ existe una base ortonormal de $S_{2k}(\Gamma)$ de autofunciones de todos los operadores de Hecke, $T_n(f) = \lambda(n) f$ con $\lambda(n) \in \mathbb{R}$, $\forall\, n \in \mathbb{N}$.*

**Definición 3.11.** Una autofunción $f$ como en el corolario anterior se denomina una autoforma de Hecke. Si además verifica que $a_1(f) = 1$ se dice que $f$ es una autoforma normalizada.

**Proposición 3.12.** *Sea $f \in \mathfrak{M}_k(\Gamma)$ una autoforma de Hecke con $T_n(f) = \lambda(n)f$ para todo $n \in \mathbb{N}$, $f \neq 0$.*

*(i) Si $k > 0$, se tiene que $a_1(f) \neq 0$. Si $k = 0$ y $a_1(f) = 0$ entonces $f$ es constante.*

*(ii) Si $k > 0$ y $a_1(f) = 1$, entonces $\lambda(n) = a_n(f)$ para todo $n$ y además $a_n(f)a_m(f) = a_{nm}(f)$ si $(n,m) = 1$.*

*(iii) $a_p(f)a_{p^s}(f) = a_{p^{s+1}}(f) + p^{2k-1}a_{p^{s-1}}(f)$ para todo $p$ primo, $s \geqslant 0$.*

*Prueba.* Se tiene que $T_n(f) = \lambda(n)\sum_{m=0}^{\infty} a_m(f)q^m = \sum_{m=0}^{\infty} a_m(T_n f)q^m$ y $a_1(T_n f) = a_n(f)$ por (3-4).

Esto implica que $a_1(T_n f) = a_n(f) = \lambda(n)a_1(f)$. Por lo tanto, si $a_1(f) = 0$, es $a_n(f) = 0$ para todo $n \in \mathbb{N}$, luego $f = a_0(f)$ es constante y $k = 0$, luego vale (i).

Si $k > 0$, es $a_1(f) \neq 0$ y si se supone que $f$ es una forma cuspidal normalizada, o sea $a_1(f) = 1$, sigue que $a_n(f) = \lambda(n)$ para todo $n \in \mathbb{N}$. Como $T_nT_m = T_{nm}$ si $(n,m) = 1$, es también $\lambda(nm) = \lambda(n)\lambda(m)$ si $(n,m) = 1$ y lo mismo para $a_{nm}(f)$, lo cual implica (ii).

Similarmente, como $T_pT_{p^s} = T_{p^{s+1}} + p^{2k-1}T_{p^{s-1}}$ la misma identidad vale para $\lambda_{p^s} = a_{p^s}(f)$, luego se tiene (iii). □

**Corolario 3.13.** *Dos autoformas cuspidales normalizadas, o bien son ortogonales o son iguales.*

*Prueba.* Sean $f$, $g$, tales que $T_n f = \lambda(n)f$, $T_n g = \mu(n)g$ para todo $n$, se tiene $\lambda(n)\langle f, g\rangle = \langle T_n(f), g\rangle = \langle f, T_n(g)\rangle = \mu(n)\langle f, g\rangle$. Luego, o bien $\langle f, g\rangle = 0$, o $\lambda(n) = \mu(n)$ para todo $n$, o sea $a_n(f) = a_n(g)$ para todo $n$, luego $f = g$. □

Los resultados anteriores sobre multiplicatividad de coeficientes de Fourier de las autoformas de Hecke normalizadas se traducen en propiedades de las funciones $L$ asociadas.

Dada $f$ una forma modular de peso $2k$, Hecke le asoció la serie de Dirichlet $L$ asociada a $f$ por $L_f(s) = \sum_{n=0}^{\infty} a_n(f)n^{-s}$. Por (3-14) los coeficientes $a_n(f)$ admiten la cota $a_n(f) = O(n^k)$, luego la serie converge absoluta y uniformemente sobre compactos en el semiplano $\mathrm{Re}(s) > k+1$ y define una función holomorfa.

Una sucesión $a_n$ que satisface $a_{nm} = a_na_m$ (resp. $a_{nm} = a_na_m$ si $(n,m) = 1$) se dice *multiplicativa* (resp. *débilmente multiplicativa)*. Dejamos la verificación del siguiente hecho como ejercicio.

**Proposición 3.14.** *Sea la serie de Dirichlet $\sum_{n=0}^{\infty} a_nn^{-s}$ donde $a_n = O(n^N)$, con $N > 0$. Entonces, si $Re(s) > N+1$*

*(i) $\sum_{n=0}^{\infty} a_nn^{-s} = \prod_{p\in\wp}\left(1 - \frac{a_p}{p^{-s}}\right)^{-1}$ sii $a_n$ es multiplicativa.*

*(ii) $\sum_{n=0}^{\infty} a_nn^{-s} = \prod_{p\in\wp}\left(\sum_{j\geqslant 0} a_{p^j}p^{-js}\right)$ sii $a_n$ es débilmente multiplicativa.*

*(iii) Si $a_n$ es débilmente multiplicativa, entonces se tiene además*

$$\sum_{n=0}^{\infty} a_n n^{-s} = \prod_{p\, primo} \left(1 - a_p p^{-s} + p^{2k-1} p^{-2s}\right)^{-1}$$

*sii se satisface*

$$a_p a_{p^n} = a_{p^{n+1}} + p^{2k-1} a_{p^{n-1}}$$

*para todo $p$ primo.*

Como consecuencia de las proposiciones anteriores se obtiene la siguiente caracterización debida a Hecke.

**Teorema 3.15.** *Si $f \in \mathcal{S}_{2k}$ es una forma cuspidal normalizada, entonces son equivalentes*

*(i) $f$ es una autoforma de Hecke .*

*(ii) $L_f(s) = \prod_{p\, primo} \left(1 - a_p p^{-s} + p^{2k-1} p^{-2s}\right)^{-1}$ .*

*(iii) La sucesión de coeficientes de Fourier de $f$ es débilmente multiplicativa y satisface $a_p(f) a_{p^n}(f) = a_{p^{n+1}}(f) + p^{2k-1} a_{p^{n-1}}(f)$ para cada $p$ primo .*

*(iv) la sucesión $a_n(f)$ satisface para todo $m, n$*

$$a_n(f) a_m(f) = \sum_{d\,|\,(n,m)} d^{2k-1} a_{\frac{nm}{d^2}}.$$

*Prueba.* Si $T_n(f) = \lambda(n) f$ para todo $n$, entonces $a_n(f) = \lambda(n)$ para cada $n$ y se satisfacen (ii) y (iii). Además es fácil ver que para formas normalizadas (iii) equivale a (ii).

Para probar que (iii) implica (i) veamos primero que $f$ es autofunción de $T_p$ para todo $p$ primo. Es decir, hay que probar que

$$\sum_{m \geqslant 0} a_m(T_p f) q^m = \lambda(p) \sum_{m \geqslant 0} a_m(f) q^m,$$

o sea $a_m(T_p f) = \lambda(p) a_m(f)$ para todo $m$ (con $\lambda(p)$ a determinar).

Si $m = 1$, esto dice que $a_p(f) = \lambda(p) a_1(f) = \lambda(p)$.

Ahora bien, si $p \nmid m$ $a_m(T_p f) = a_{pm}(f) = a_p(f) a_m(f)$ por (iii).

En el caso $p \mid m$, $m = p^s m'$ con $(p, m') = 1$, se tiene que

$$\begin{aligned} a_m(T_p f) &= a_{mp}(f) + p^{2k-1} a_{m/p}(f) = a_{p^{s+1} m'}(f) + p^{2k-1} a_{p^{s-1} m'}(f) \\ &= (a_{p^{s+1}}(f) + p^{2k-1} a_{p^{s-1}}(f)) a_{m'}(f) \\ &= a_p(f) a_{p^s}(f) a_{m'}(f) = a_p(f) a_m(f). \end{aligned} \tag{3-18}$$

Luego $T_p(f) = a_p(f) f$ para todo $p$ primo. Afirmamos que de esto resulta que $T_n(f) = a_n(f) f$ para todo $n \in \mathbb{N}$.

En efecto, si $n = \prod p_j^{i_j}$, es $T_n(f) = T_{p_1{}^{i_1}} \dots T_{p_r{}^{i_r}}(f)$, luego basta ver que si $p$ es primo es $T_{p^n}(f) = a_{p^n}(f)$ para todo $n \in \mathbb{N}$. Para esto procedemos inductivamente, usando (iii)

$$\begin{aligned} T_{p^{n+1}}(f) &= T_p T_{p^n}(f) - p^{2k-1} T_{p^{n-1}}(f) \\ &= a_p(f) a_{p^n}(f) - p^{2k-1} a_{p^{n-1}}(f) = a_{p^{n+1}}(f). \end{aligned}$$

Finalmente, es claro que (iv) implica (iii). La prueba de que (iii) implica (iv) queda como ejercicio para el lector. □

**Ejemplo 3.16.** Observamos que si $\dim \mathcal{S}_{2k} = 1$, cualquier $f \in \mathcal{S}_{2k}$ es una autoforma. Este es el caso para todo $k$ tal que $6 \leqslant k < 12$; por ejemplo, $\Delta(z) \in \mathcal{S}_{12}$ es autofunción de $T_n$ para todo $n$.

Se tiene que

$$\Delta(z) = \sum_{n=1}^{\infty} \tau(n) q^n, \quad \tau(1) = 1.$$

La función $\tau(n)$ fue estudiada por Ramanujan quien conjeturó que

$$\tau(nm) = \tau(n)\tau(m) \text{ si } (n,m) = 1 \tag{3-19}$$

$$\tau(p)\tau(p^n) = \tau(p^{n+1}) + p^{11}\tau(p^{n-1}), \tag{3-20}$$

probado por Mordell en 1920.

Además Ramanujan conjeturó que $|\tau(p)| \leqslant 2p^{11/2}$ y más generalmente que

$$\tau(n) = O_\varepsilon(n^{(11/2)+\varepsilon})$$

para todo $\varepsilon > 0$. Estos hechos fueron probados por Deligne como consecuencia de las conjeturas de Weil (1969 y 1973). Deligne probó que para $f \in \mathcal{S}_k$ normalizada se tiene

$$|a_p(f)| \leqslant 2p^{k-1/2}, \quad |a_n(f)| \leqslant \sigma_0(n) n^{k-1/2}. \tag{3-21}$$

Como ya se vio no es difícil de probar la acotación $a_n(f) = O(n^k)$, sin embargo (3-21) sí es muy difícil.

Notamos que por los hechos anteriores se tiene

$$L_\Delta(s) = \sum_{n \geqslant 1} \tau(n) n^{-s} = \prod_p (1 - \tau(p)p^{-s} + p^{11}p^{-2s})^{-1} \text{ si } \operatorname{Re} s > 7.$$

Sea $f$ una autoforma normalizada para un subgrupo de congruencia. Factorizando $1 - a_p(f)x + p^{2k-1}x^2 = (1 - \alpha_p x)(1 - \alpha'_p x)$, se tiene $\alpha_p + \alpha'_p = a_p(f)$ y $\alpha_p \alpha'_p = p^{2k-1}$.

La llamada conjetura de Ramanujan–Petersson afirma que para todo $p$ se cumple que $\alpha'_p = \overline{\alpha_p}$, o sea $|\alpha_p| = |\alpha'_p| = p^{(2k-1)/2}$.

**Ejemplo 3.17.** A continuación probaremos que las series de Eisenstein $G_{2k}(z)$ son también autofunciones de todos los operadores de Hecke. La autofunción normalizada y la función $L$ asociada son respectivamente

$$G'_{2k}(z) = \frac{(-1)^k B_k}{4k} + \sum_{1}^{\infty} \sigma_{2k-1}(m) q^m,$$

$$L_{G'_{2k}}(s) = \prod_p \left(1 - \sigma_{2k-1}(p)p^{-s} + p^{2k-1}p^{-2s}\right)^{-1}.$$

Verificaremos que $\lambda(n) = \sigma_{2k-1}(n)$ para todo $n$ y $\lambda(p) = \sigma_{2k-1}(p) = 1 + p^{2k-1}$.

En primer lugar

$$a_0\left(T_n\left(G'_{2k}\right)\right) = \left(\sum_{d|n} d^{2k-1}\right) a_0(f) = \sigma_{2k-1}(n)a_0(f).$$

Resta verificar que $\lambda(p) = 1 + p^{2k-1}$ o equivalentemente que

$$a_p(T_m(G'_{2k})) = \sigma_{2k-1}(p)\sigma_{2k-1}(m) \quad \forall\, p, m.$$

Supondremos primero que $p \nmid m$. Entonces

$$\begin{aligned} a_p(T_m(G'_{2k}) &= \sigma_{2k-1}(pm) = \sum_{d|pm} d^{2k-1} = \sum_{d|m} d^{2k-1} + \sum_{d|m} p^{2k-1}d^{2k-1} \\ &= (1 + p^{2k-1})\sigma_{2k-1}(m) = \sigma_{2k-1}(p)\sigma_{2k-1}(m). \end{aligned}$$

Dejamos la verificación del caso $p \mid m$ como ejercicio, siendo similar al argumento en (3-18).

Los ejemplos anteriores y la descripción del álgebra $\mathfrak{M}$ en el caso de $\Gamma = SL_2(\mathbb{Z})$, generada como álgebra por las series de Eisenstein, dice en particular que para cada $k$, $\mathfrak{M}_{2k}$ posee una base de formas modulares cuyos coeficientes de Fourier son números enteros. Por otra parte, las autoformas de Hecke son formas especiales en cada dimensión. Por ejemplo $\Delta(z)^2$ es una forma cuspidal de peso 24 que no es una autoforma de Hecke. Se tiene el siguiente hecho.

**Proposición 3.18.** *Si $\Gamma = SL_2(\mathbb{Z})$, los coeficientes de Fourier de una autoforma de Hecke normalizada son enteros algebraicos totalmente reales.*

*Prueba.* Sea $f \in S_{2k}(\Gamma)$. Como ya se observó, $\mathfrak{M}_{2k}$ admite una base de formas con coeficientes de Fourier en $\mathbb{Z}$ ($\Delta(z)$ y las series de Eisenstein $G_m(z)$ tienen esta propiedad). Claramente, $T_n$ preserva el retículo $\mathfrak{L}$ de formas con coeficientes enteros y en una $\mathbb{Z}$-base de $\mathfrak{L}$ está representado por una matriz a coeficientes enteros. Sus autovalores $\lambda(n)$ son números reales que son raíces del polinomio característico de esta matriz que es mónico con coeficientes enteros, luego cada $\lambda(n)$ es un entero algebraico totalmente real. □

**3.1 Subgrupos de nivel $N$.** La teoría de formas modulares se simplifica en el caso en que $\Gamma = SL_2(\mathbb{Z})$. Para subgrupos de congruencia es bastante más complicada y hay importantes preguntas sin respuesta. En esta subsección describiremos los cambios necesarios en la teoría en el caso de los subgrupos de congruencia de Hecke $\Gamma_0(N)$, resultados debidos a Atkin y Lehner (1970).

En primer lugar se define el operador $T_n$ para $f \in S_k(\Gamma_0(N))$ por

$$T_n(f)(z) = n^{k-1} \sum_{ad=n} \sum_{0\leqslant b<d} d^{k-1} f\left(\tfrac{az+b}{d}\right)$$

donde $\left[\begin{smallmatrix} a & b \\ c & d \end{smallmatrix}\right] \in \Gamma_0(N)$.

Interesan especialmente los $T_n$ con $(n, N) = 1$. Bajo esta hipótesis, la Proposición 3.2 es válida con igual prueba, lo que implica que $T_n(f)$ es holomorfa (y cuspidal si $f$ lo es) en $i_\infty$. Se sigue que $T_n$ preserva $\mathfrak{M}_k(\Gamma_0(N))$ y $\mathcal{S}_k(\Gamma_0(N))$, previo análisis del comportamiento de $T_n(f)$ en las cúspides de $\Gamma$ distintas de $i_\infty$. Para todo $n, m$ con $(nm, N) = 1$, con igual prueba que para $N = 1$ se obtiene que $T_n T_m = T_m T_n = T_{nm}$ así como la expresión (3-8) para todo primo $p \nmid N$. En consecuencia, los operadores $T_n$ con $(n, N) = 1$ generan un álgebra conmutativa $\mathcal{H}_N$.

El producto de Petersson se define por (3-13) para $f, g$ en $\mathcal{S}_k(\Gamma_0(N))$ (basta que una de ellas, $f$ o $g$ sea cuspidal) y se prueba que los operadores $T_n$ con $(n, N) = 1$ son autoadjuntos, luego existe una base ortonormal de $S_k(\Gamma_0(N))$ de autofunciones comunes a todos los $T_n$ con $(n, N) = 1$.

Los operadores de Hecke en el caso en que $p \nmid N$ tienen propiedades diferentes que los $T_p$ en que $p \mid N$. Estos $p$ suelen llamarse los primos ramificados. Para estos primos $T_p$ no es autoadjunto (ver Shimura (1971)).

Para una autoforma $f$ que es autofunción de $T_n$ para todo $(n, N) = 1$ se tiene como anteriormente, que $\lambda(n) a_1(f) = a_n(f)$. Sin embargo, no se puede concluir en este caso que si $a_1(f) = 0$ entonces $f = 0$, ya que esto implica $a_n(f) = 0$ sólo para los $n$ tales que $(n, N) = 1$.

Para el grupo $\Gamma_0(N)$ de nivel es conveniente distinguir el espacio de las formas que provienen de un nivel inferior a $N$. Como primera observación sea la secuencia de grupos

$$\Gamma_0(N) \subset \Gamma_0(M) \subset \Gamma(1) = SL(2, \mathbb{Z})$$

donde $M \mid N$. Entonces, si $f \in S_k(\Gamma_0(M))$, para todo $h \mid \frac{N}{M}$ la forma $f|\left[\begin{smallmatrix} h & 0 \\ 0 & 1 \end{smallmatrix}\right](z) = f(hz)$ está en $S_k(\Gamma_0(N))$. Esto es consecuencia de la identidad

$$\left[\begin{smallmatrix} h & 0 \\ 0 & 1 \end{smallmatrix}\right] \left[\begin{smallmatrix} a & b \\ cN & d \end{smallmatrix}\right] \left[\begin{smallmatrix} h^{-1} & 0 \\ 0 & 1 \end{smallmatrix}\right] = \begin{bmatrix} a & bh \\ c\frac{N}{h} & d \end{bmatrix}$$

que pertenece a $\Gamma_0(M)$ pues $h \mid \frac{N}{M}$, luego $f|\left[\begin{smallmatrix} h & 0 \\ 0 & 1 \end{smallmatrix}\right]|\gamma = f|\left[\begin{smallmatrix} h & 0 \\ 0 & 1 \end{smallmatrix}\right]$ para todo $\gamma \in \Gamma_0(N)$.

Se denota por $S_k(\Gamma_0(N))_o$ ($M_k(\Gamma_0(N))_o$) al subespacio de $S_k(\Gamma_0(N))$ generado por las formas construidas por el procedimiento anterior, o sea $g(z) = f(hz)$ donde $f \in S_k(\Gamma_0(M))$ y $h \mid \frac{N}{M}$. Este es el espacio de las llamadas *formas viejas*, o sea las formas que provienen de divisores propios de $N$. Es decir

$$\mathcal{S}_k(\Gamma_0(N))_o = \oplus_{n|N} \oplus_{d|\frac{N}{n}} S_k(\Gamma_0(\tfrac{N}{n}))\,.$$

Se tiene además

$$S_k(\Gamma_0(N)) = S_k(\Gamma_0(N))_o \oplus S_k(\Gamma_0(N))_p$$

donde $S_k(\Gamma_0(N))_p := S_k(\Gamma_0(N))_o^\perp$ es el subespacio ortogonal a $S_k(\Gamma_0(N))_o$, el espacio de las llamadas *formas nuevas* o *primitivas*. Observamos que como $f|d(z) = \sum_{m\geqslant 1} a_m(f) e^{2\pi i dmz}$, el coeficiente $a_1(f|d) = 0$ si $d > 1$, esto es, las formas viejas nunca son normalizadas. Claramente, si $(n, N) = 1$, $T_n$ conmuta con cada operador $f \to f|d$, con $d \mid N$, luego $T_n$ preserva $S_k(\Gamma_0(N))_o$ y por ser autoadjunto, $T_n$ también preserva $S_k(\Gamma_0(N))_n$.

La situación fue clarificada por Atkin y Lehner (1970). A continuación enunciamos algunos de sus principales resultados, donde recuperan resultados válidos en el caso $N = 1$, en el contexto de autoformas primitivas.

**Teorema 3.19.** *(i) Si $f, g \in S_k(\Gamma_0(N))_p$ son autofunciones de $T_n$ con el mismo autovalor $\lambda(n)$ para todo $(n, N) = 1$ entonces $f = cg$ para algún $c \in \mathbb{C}$.*

*(ii) Sea $f \in S_k(\Gamma_0(N))_p$ autofunción de $T_n$ para todo $(n, N) = 1$. Entonces $f$ es una autoforma, vale $a_n(f) = \lambda(n)a_1(f)$ para todo $n$ y si $f \neq 0$ entonces $a_1(f) \neq 0$.*

*(iii) Dos autoformas primitivas normalizadas son perpendiculares o son iguales y el espacio $S_k(\Gamma_0(N))_p$ tiene una base de autoformas de Hecke normalizadas que es única salvo permutaciones.*

## Referencias

Roberto J. Miatello
miatellorj@gmail.com
miatello@famaf.unc.edu.ar
CIEM–FaMAF
Universidad Nacional de Córdoba
5000-Córdoba
Argentina

# FORMAS MODULARES CUATERNIÓNICAS

Gonzalo Tornaría

## Índice general

El propósito general de estas notas es dar una introducción muy breve a las formas modulares cuaterniónicas. Están basadas en un minicurso dictado en Cusco, Perú, en la escuela AGRA de 2015.

En la primera sección se presentan las álgebras de cuaterniones sobre $\mathbb{Q}$, su clasificación local y global, y los aspectos aritméticos básicos que son necesarios para desarrollar la teoría de formas modulares cuaterniónicas (ordenes, ideales, clases de ideales).

En la segunda sección se presentan las formas modulares cuaterniónicas para el caso particular de álgebras de cuaterniones definidas. Esta situación se corresponde, en el sentido de las notas de Harris (2020), al caso de variedades de Shimura de dimensión 0, que tiene la particularidad de no involucrar análisis.

La sección culmina explorando la correspondencia de Eichler, que explica la relación entre las formas modulares cuaterniónicas así construidas y las formas modulares clásicas, y que es el precursor de la correspondencia de Jacquet–Langlands.

Para el lector interesado sugerimos algunas referencias adicionales. Para álgebras de cuaterniones la referencia clásica es Vignéras (1980). Para la demostración de algunos resultados nos referiremos a Voight (s.f.) que es una referencia más completa y actualizada, y que se encuentra disponible para descargar en la página de su autor. Para la correspondencia de Eichler, sugerimos consultar los artículos de Eichler (1973) y de Gross (1987).

El lector encontrará provechoso complementar estas notas con otras en este volumen, en particular las notas de Miatello (2020) introducen las formas modulares clásicas, las notas de Pacetti (2020) introducen las formas de Hilbert que generalizan las formas modulares clásicas a cuerpos de números totalmente reales, y las notas de Harris (2020) presentan una versión más general de formas modulares cuaterniónicas.

> Ocasionalmente haremos un interludio indicado por los márgenes como aquí, en el que agregaremos alguna información adicional que podría ser omitida en una primera lectura.

# 1 Álgebras de cuaterniones

Sea $F$ un cuerpo de característica $\neq 2$. Dados $a, b \in F^\times = F - \{0\}$ construimos un álgebra sobre $F$ con base $\{1, i, j, k\}$, donde la multiplicación está dada por

$$i^2 = a, \quad j^2 = b, \quad \text{y} \quad k := ij = -ji \; .$$

La tabla completa de la multiplicación puede obtenerse a partir de éstas relaciones (ver ejercicio 1). Por ejemplo $k^2 = ijij = i(-ij)j = -i^2 j^2 = -ab$. A esta álgebra la denotaremos $(a, b)_F$ y diremos que es un *álgebra de cuaterniones*.

**Ejercicio 1.** Escribir la tabla completa de la multiplicación en $(a, b)_F$. Concluir que que es un álgebra de dimensión 4, y determinar su centro.

Esta construcción es una generalización del álgebra de *cuaterniones de Hamilton*, que es $\mathbb{H} = (-1, -1)_{\mathbb{R}}$. Para $x = x_0 + x_1 i + x_2 j + x_3 k \in \mathbb{H}$ definimos $x^\star := x_0 - x_1 i - x_2 j - x_3 k$ y la norma $\nu(x) := x x^\star = x_0^2 + x_1^2 + x_2^2 + x_3^2$. Como $\nu(x) \in \mathbb{R}^\times$ para $x \neq 0$ deducimos que $\mathbb{H}$ es un álgebra de división (todo $x \neq 0$ es invertible).

En general, un álgebra de cuaterniones no necesariamente es un álgebra de división:

**Proposición 1.1.** *Si $a$ es un cuadrado en $F^\times$, entonces $(a, b)_F \cong M_2(F)$.*

*Prueba.* Si $a = \alpha^2$ con $\alpha \in F^\times$, podemos definir un homomorfismo de álgebras $\varphi : (a, b)_F \to M_2(F)$ por

$$\varphi(i) := \begin{pmatrix} \alpha & 0 \\ 0 & -\alpha \end{pmatrix}, \qquad \varphi(j) := \begin{pmatrix} 0 & 1 \\ b & 0 \end{pmatrix}, \qquad \varphi(k) := \begin{pmatrix} 0 & \alpha \\ -b\alpha & 0 \end{pmatrix} .$$

Es fácil verificar que cumple las relaciones, y como $\{\varphi(1), \varphi(i), \varphi(j), \varphi(k)\}$ es linealmente independiente sobre $F$ (pues $\alpha \neq -\alpha$) se sigue que $\varphi$ es un isomorfismo. ☐

Cuando $a$ no es un cuadrado en $F$, podemos considerar el cuerpo $K = F(\sqrt{a})$ y resulta que $(a, b)_F$ es una subálgebra de $(a, b)_K \cong M_2(K)$. Denotemos $\sigma$ al automorfismo no trivial de $K/F$, entonces podemos describir explícitamente $(a, b)_F$ como

$$(a, b)_F \cong \{\begin{pmatrix} u & v \\ b v^\sigma & u^\sigma \end{pmatrix} : u, v \in K\} .$$

$$i \mapsto \begin{pmatrix} \sqrt{a} & 0 \\ 0 & -\sqrt{a} \end{pmatrix}$$

$$j \mapsto \begin{pmatrix} 0 & 1 \\ b & 0 \end{pmatrix}$$

El álgebra $(a, b)_K$ se obtiene por *extensión de escalares* a partir de $(a, b)_F$, es decir que $K \otimes (a, b)_F = (a, b)_K \cong M_2(K)$.

**Ejemplo 1.2.** En el caso de los cuaterniones de Hamilton, extendiendo escalares a $\mathbb{C} = \mathbb{R}(\sqrt{-1})$ tenemos que $\mathbb{H} \cong \{\begin{pmatrix} u & v \\ -\bar{v} & \bar{u} \end{pmatrix} : u, v \in \mathbb{C}\} \subseteq M_2(\mathbb{C})$.

**Ejemplo 1.3.** Si $a, b \in \mathbb{Q}$ con $a > 0$, entonces $F(\sqrt{a}) \subseteq \mathbb{R}$. Entonces tenemos una representación de $(a, b)_{\mathbb{Q}}$ en $M_2(\mathbb{R})$ dada por

$$(a, b)_{\mathbb{Q}} \cong \left\{ \begin{pmatrix} x_0 + x_1\sqrt{a} & x_2 + x_3\sqrt{a} \\ (x_2 - x_3\sqrt{a})\, b & x_0 - x_1\sqrt{a} \end{pmatrix} : x_0, x_1, x_2, x_3 \in \mathbb{Q} \right\} \subseteq M_2(\mathbb{R}) .$$

Recordar que un álgebra sobre $F$ se dice *central* si su centro es $F$ y *simple* si no contiene ideales biláteros no triviales.

**Proposición 1.4.** *El álgebra $(a,b)_F$ es un álgebra central simple sobre $F$.*

*Prueba.* Como $(a,b)_F$ es de dimensión 4 y no conmutativa es central (ejercicio 2). Si $a$ es un cuadrado, entonces $(a,b)_F \cong M_2(F)$ que es simple (ejercicio 3). Si $a$ no es un cuadrado, podemos considerar el cuerpo $K = F(\sqrt{a})$ como arriba. Como $(a,b)_F \otimes K \cong M_2(K)$ es un álgebra simple sobre $K$, se sigue fácilmente que $(a,b)_F$ es un álgebra simple sobre $F$ (ejercicio 4). □

**Ejercicio 2.** Sea $F$ un cuerpo y sea $D$ un álgebra sobre $F$.

1. Probar que si $\dim_F D = 2$ entonces $D$ es conmutativa.
2. Probar que si $\dim_F D = 4$ entonces $D$ es conmutativa o central.
   *Sugerencia: ¿Qué puede decirse sobre la dimensión del centro de $D$ sobre $F$?*

Concluir que $(a,b)_F$ es central.

**Ejercicio 3.** Probar que $M_2(F)$ es simple.
*Sugerencia: determinar todos los ideales a izquierda de $M_2(F)$.*

**Ejercicio 4.** Probar que $(a,b)_F$ es simple.
*Sugerencia: usar que $(a,b)_F \subseteq M_2(K)$ con $K = F(\sqrt{a})$.*

**Corolario 1.5.** *O bien $(a,b)_F$ es un álgebra de división, o bien $(a,b)_F \cong M_2(F)$.*

*Prueba.* Se sigue de la clasificación de las álgebras simples, o del ejercicio 5. □

**Ejercicio 5.** Sea $D$ un álgebra simple de dimensión 4 sobre un cuerpo $F$. Probar que si $D$ no es de división, entonces $D \cong M_2(F) \cong (1,1)_F$.
*Sugerencia: si $\mathfrak{a} \subset D$ es un ideal propio a izquierda, de dimensión $n$ sobre $F$, la acción de $D$ en $\mathfrak{a}$ dada por multiplicación a izquierda induce un homomorfismo $D \to M_n(F)$ inyectivo.*

**Ejercicio 6.** Sea $D$ un álgebra de dimensión 4 sobre un cuerpo $F$ de característica $\neq 2$. Supongamos que $D$ es de división y no conmutativa.

1. Si $x \in D \setminus F$, entonces $F[x]$ es un álgebra conmutativa de dimensión 2 sobre $F$.
2. Existe $i \in D \setminus F$ tal que $i^2 = a \in F^\times$.
3. Considerar $\varphi : D \to D$ dada por $\varphi(y) = iyi^{-1}$. Probar que $\varphi^2 = \mathrm{id}_D$.
4. Concluir que existe $j \in D$ tal que $ij = -ji$.
5. Probar que $j^2 = b \in F^\times$ y concluir que $D \cong (a,b)_F$.

Los últimos dos ejercicios muestran que si $D$ es un álgebra central simple de dimensión 4 sobre un cuerpo $F$ de característica $\neq 2$ entonces $D \cong (a,b)_F$ para algunos $a, b \in F^\times$. Por eso habitualmente se define un álgebra de cuaterniones como un álgebra central simple de dimensión 4. Esta última definición resulta ser adecuada incluso en característica 2. Consultar Voight (s.f., Capítulo 6) para la construcción de álgebras de cuaterniones en característica 2.

Veremos ahora cómo definir la norma, lo que nos permitirá caracterizar las álgebras de cuaterniones de división. Un álgebra de cuaterniones $D = (a, b)_F$ posee una *involución canónica* $\star : D \to D$ que definimos, para $x = x_0 + x_1\, i + x_2\, j + x_3\, k \in D$, como

$$x^\star := x_0 - x_1\, i - x_2\, j - x_3\, k\,.$$

Notemos que $\star$ es una involución $F$-lineal que deja fijo $F$, y que $(xy)^\star = y^\star\, x^\star$. Definimos la *traza* y la *norma* (reducidas) de $x$ como

$$\tau(x) := x + x^\star = 2\,x_0\,, \qquad \nu(x) := x\,x^\star = x^\star\,x = x_0^2 - a\,x_1^2 - b\,x_2^2 + ab\,x_3^2\,.$$

Como $x^2 - (x + x^\star)\,x + x^\star\,x = 0$, tenemos que $x \in D$ es raíz de

$$T^2 - \tau(x)\,T + \nu(x) \in F[T]\,,$$

que llamamos *polinomio característico* (reducido) de $x$.

*Observación.* Por lo anterior, si $x \notin F$ entonces $F[x]$ es un álgebra de dimensión 2 con base $\{1, x\}$. Se sigue que la traza y la norma quedan determinadas por la relación $x^2 = \tau(x)\,x - \nu(x)$ y la involución canónica por $x^\star = \tau(x) - x$.

Es claro de la definición que $x \in F$ si y sólo si $x = x^\star$. Por otro lado el discriminante del polinomio característico de $x$ es $\Delta(x) = \tau(x)^2 - 4\nu(x) = (x - x^\star)^2$, por lo cual $x = x^\star$ si y sólo si su polinomio característico es un cuadrado.

**Ejemplo 1.6.** En el álgebra de matrices $M_2(F)$ la involución canónica está dada por

$$\left(\begin{smallmatrix} u & v \\ w & t \end{smallmatrix}\right)^\star = \left(\begin{smallmatrix} t & -v \\ -w & u \end{smallmatrix}\right),$$

mientras que la traza y norma reducidas estarán dadas por la traza y el determinante de la matriz, respectivamente:

$$\tau\left(\begin{smallmatrix} u & v \\ w & t \end{smallmatrix}\right) = u + t, \qquad \nu\left(\begin{smallmatrix} u & v \\ w & t \end{smallmatrix}\right) = ut - vw.$$

La norma nos permite caracterizar los elementos invertibles en un álgebra de cuaterniones $D$. En efecto, $x \in D$ es invertible si y solamente si $\nu(x) \in F^\times$, y en tal caso el inverso de $x$ está dado por $x^{-1} = \nu(x)^{-1}\,x^\star$.

**Proposición 1.7.** *Sea $D = (a, b)_F$ un álgebra de cuaterniones sobre $F$. Son equivalentes*

1. *$D \cong M_2(F)$.*
2. *$D$ no es un álgebra de división.*
3. *La forma cuadrática $\nu$ es isotrópica en $D$ (tiene ceros no triviales).*
4. *La forma cuadrática $a\,x^2 + b\,y^2$ representa* 1.

*Prueba.* El ejercicio 5 prueba (1) $\Leftrightarrow$ (2). La equivalencia (2) $\Leftrightarrow$ (3) sigue de la caracterización de los elementos invertibles en $D$. La implicación (4) $\Rightarrow$ (3) es clara. Falta probar (3) $\Rightarrow$ (4). Sea $x \in D$, $x \neq 0$ tal que $\nu(x) = 0$. El ideal $Dx$, como espacio vectorial sobre $F$, tiene dimensión 2, entonces existe un $y \in Dx$, $y \neq 0$ tal que $y = y_0 + y_1\, i + y_2\, j$, es decir $\nu(y) = y_0^2 - a\,y_1^2 - b\,y_2^2 = 0$. Para concluir, si $y_0 \neq 0$ se sigue que $a\,x^2 + b\,y^2$ representa 1, mientras que si $y_0 = 0$ se sigue que $a\,x^2 + b\,y^2$ es isotrópica, pero una forma cuadrática isotrópica es universal y por lo tanto representa 1. $\square$

**1.1 Clasificación de álgebras de cuaterniones.** En esta sección vamos a mencionar resultados de clasificación de álgebras de cuaterniones sobre algunos cuerpos, particularmente sobre $\mathbb{Q}$.

**Proposición 1.8.** *Toda álgebra de cuaterniones sobre $\mathbb{C}$ es isomorfa a $M_2(\mathbb{C})$.*

*Prueba.* Todo $a \in \mathbb{C}^\times$ es un cuadrado, entonces $(a,b)_\mathbb{C} \cong M_2(\mathbb{C})$. □

**Proposición 1.9.** *Un álgebra de cuaterniones sobre $\mathbb{R}$ es isomorfa a $M_2(\mathbb{R})$ o a $\mathbb{H}$.*

*Prueba.* Si $a > 0$ o si $b > 0$ entonces $(a,b)_\mathbb{R} \cong M_2(\mathbb{R})$. Por otra parte si $a = -x^2$ y $b = -y^2$ entonces es claro que $(-x^2,-y^2)_\mathbb{R} \cong (-1,-1)_\mathbb{R} \cong \mathbb{H}$. □

**Proposición 1.10.** *Toda álgebra de cuaterniones sobre un cuerpo finito $\mathbb{F}_q$ es isomorfa a $M_2(\mathbb{F}_q)$.*

*Prueba.* Vamos a probar que $a\,x^2 + b\,y^2 = 1$ tiene solución en $\mathbb{F}_q$. En efecto, la imagen de $a\,x^2$ tiene $(q+1)/2$ elementos y la imagen de $1 - b\,y^2$ también tiene $(q+1)/2$ elementos. Por el principio del palomar, tienen un valor en común, es decir una solución de $a\,x^2 = 1 - b\,y^2$ como afirmamos. □

Sea $D$ un álgebra de cuaterniones sobre $\mathbb{Q}$. Consideramos $D_\infty := D \otimes \mathbb{R}$, que es un álgebra de cuaterniones sobre $\mathbb{R}$. Decimos que $D$ es *definida* cuando $D_\infty \cong \mathbb{H}$ y que $D$ es *indefinida* cuando $D_\infty \cong M_2(\mathbb{R})$. Esto corresponde a que la forma norma $\nu$ sea definida o indefinida, respectivamente. Es claro que la clase de isomorfismo de $D_\infty$ es un invariante de $D$; por ejemplo, $(-1,-1)_\mathbb{Q} \ncong (-1,3)_\mathbb{Q}$.

Sin embargo, esto no alcanza para tener una clasificación sobre $\mathbb{Q}$. Por ejemplo las álgebras $D = (-1,-1)_\mathbb{Q}$ y $D' = (-1,-3)_\mathbb{Q}$ son ambas definidas pero no son isomorfas. Consideremos las respectivas formas normas:

$$\nu(x) = x_0^2 + x_1^2 + x_2^2 + x_3^2\,, \qquad \nu'(y) = y_0^2 + y_1^2 + 3\,y_2^2 + 3\,y_3^2\,.$$

Es fácil ver que la segunda no tiene ceros no triviales módulo 9, mientras que la primera los tiene ($1^2+2^2+2^2=9$). ¿Es posible usar esto para probar que $(-1,-1)_\mathbb{Q} \ncong (-1,-3)_\mathbb{Q}$? Esto presenta algunas dificultades. En primer lugar los enteros módulo 9 no son un cuerpo. Si trabajamos módulo 3, para tener un cuerpo, entonces ambas formas tienen ceros no triviales módulo 3, y de todas maneras ya vimos que hay una única clase de isomorfismo de álgebras de cuaterniones sobre $\mathbb{F}_3$. Por otra parte $\mathbb{Q}$ no está contenido en los enteros módulo 9 o módulo 3, por lo que no es tan clara la extensión de escalares.

Para resolver estas dificultades, debemos emplear el cuerpo de los números 3-ádicos $\mathbb{Q}_3$. En primer lugar $\mathbb{Q} \subset \mathbb{Q}_3$ por lo cual podemos considerar el cambio de base $D \otimes \mathbb{Q}_3$. Por otra parte, las afirmaciones hechas arriba acerca de los ceros módulo 9 de las formas normas se traducen en: $\nu$ es isotrópica en $D \otimes \mathbb{Q}_3$ (tiene ceros no triviales), y $\nu'$ es anisotrópica en $D' \otimes \mathbb{Q}_3$ (no tiene ceros excepto el trivial $\nu'(0) = 0$). Entonces $D \otimes \mathbb{Q}_3 \cong M_2(\mathbb{Q}_3)$ mientras que $D' \otimes \mathbb{Q}_3$ es un álgebra de división.

**Interludio: los números $p$-ádicos.** Si consideramos el valor absoluto usual en $\mathbb{Q}$ obtenemos una métrica en $\mathbb{Q}$ cuya completación son los números reales. De la misma manera podemos construir, para cada primo $p$, los *números $p$-ádicos* como la completación de $\mathbb{Q}$ con respecto al valor absoluto $p$-ádico.

Al cuerpo de los números $p$-ádicos lo denotaremos $\mathbb{Q}_p$, y en ocasiones denotaremos $\mathbb{Q}_\infty$ al cuerpo de los números reales; estos son ejemplos de *cuerpos locales*. Usaremos la letra $v$ para referirnos a un primo $p$ o al *lugar arquimediano* $\infty$.

Para construir el valor absoluto $p$-ádico definimos la *valuación $p$-ádica* $v_p : \mathbb{Q}^\times \to \mathbb{Z}$ dada por $v_p(p^r \frac{a}{b}) := r$ siempre que $p \nmid a$, $p \nmid b$. Finalmente definimos $|x|_p := p^{-v_p(x)}$ para $x \in \mathbb{Q}^\times$, y $|0|_p := 0$.

**Ejercicio 7.** Verificar que $|x|_p$ es un valor absoluto que satisface la *desigualdad ultramétrica* $|x+y|_p \leqslant \max(|x|_p, |y|_p)$.

Como $\mathbb{Q}_p$ es la completación de $\mathbb{Q}$, todo $x \in \mathbb{Q}_p$ puede ser escrito como límite de una sucesión de números racionales. Observar que $|p^r|_p \to 0$ cuando $r \to \infty$, y podemos escribir $x$ como una serie $\sum_{n \geqslant N_0} a_n\, p^n$. Las sumas parciales corresponden a considerar $x$ módulo $p^r$.

Para un primo $p$ podemos definir $\mathbb{Z}_p$, el anillo de los *enteros $p$-ádicos*, como la clausura de $\mathbb{Z}$ en $\mathbb{Q}_p$.

**Ejercicio 8.** Probar las siguientes propiedades de $\mathbb{Z}_p$.

1. $\mathbb{Z}_p$ es una completación de $\mathbb{Z}$ con respecto al valor absoluto $p$-ádico.
2. $\mathbb{Z}_p = \{x \in \mathbb{Q}_p : |x| \leqslant 1\} = \{x \in \mathbb{Q}_p : |x| < 2\}$.
3. $\mathbb{Z}_p$ es un anillo, es compacto y es abierto en $\mathbb{Q}_p$.
4. Todo elemento de $\mathbb{Z}_p$ puede escribirse como una serie $\sum_{n \geqslant 0} a_n\, p^n$.
5. Para todo $r \geqslant 1$ existe un homomorfismo canónico $\mathbb{Z}_p \to \mathbb{Z}/p^r\mathbb{Z}$.

En la sección 2.1 de las notas de Pacetti (2020) puede encontrarse una introducción algo más detallada a los números $p$-ádicos.

En nuestro ejemplo, la forma cuadrática $y_0^2 + y_1^2 + 3\,y_2^2 + 3\,y_3^2$ no tiene ceros no triviales en $\mathbb{Q}_3$, pues no los tiene módulo 9. Por otra parte, la forma cuadrática $x_0^2 + x_1^2 + x_2^2 + x_3^2$ tiene ceros módulo 9, e.g. $1^2 + 2^2 + 2^2 \equiv 0 \pmod 9$. A partir de este cero es posible construir ceros en $\mathbb{Q}_3$ mediante una construcción conocida como *Lema de Hensel*. Una solución es $1^2 + 2^2 + x^2 = 0$ con $x = 2 + 2 \cdot 3^2 + 3^4 + 3^6 + 3^9 + \cdots \in \mathbb{Q}_3$.

**Proposición 1.11** (Clasificación local)**.** *Existe un álgebra de cuaterniones de división sobre $\mathbb{Q}_v$, única a menos de isomorfismo. Cuando $v = p \neq 2$ está dada por $(r, p)_{\mathbb{Q}_p}$, donde $r \in \mathbb{Z}$ es un no residuo cuadrático módulo $p$.*

En otras palabras, existen dos álgebras de cuaterniones sobre $\mathbb{Q}_v$ a menos de isomorfismo: el álgebra de matrices y un álgebra de cuaterniones de división.

Para $v = \infty$ el resultado se demuestra en la Proposición 1.9. El siguiente ejercicio lo demuestra para $v \neq 2, \infty$.

**Ejercicio 9.** Sea $p \neq 2$ primo, y sea $r \in \mathbb{Z}$ un no residuo cuadrático módulo $p$.

1. Probar que $(r, p)_{\mathbb{Q}_p}$ es un álgebra de división.
2. Para cada par $a, b \in \{1, r, p, pr\}$, determinar si $(a, b)_{\mathbb{Q}_p}$ es un álgebra de matrices o un álgebra de división.

3. Probar que las 6 álgebras de división del paso anterior son isomorfas.

4. Usar que $\mathbb{Q}_p^\times/(\mathbb{Q}_p^\times)^2 = \{1, r, p, pr\}$ para concluir la Proposición 1.11 para $v = p \neq 2$.

Cuando $v = 2$ el álgebra $(-1, -1)_{\mathbb{Q}_2}$ es un álgebra de división. La misma idea del ejercicio anterior sirve para probar que es única a menos de isomorfismo, aunque resulta más engorroso ya que $\mathbb{Q}_2^\times/(\mathbb{Q}_2^\times)^2 = \{1, 3, 5, 7, 2, 6, 10, 14\}$ tiene más elementos.

Cuando $D$ es un álgebra de cuaterniones sobre $\mathbb{Q}$, podemos considerar $D_v := D \otimes \mathbb{Q}_v$ para todo $v$ (primo o $\infty$). Decimos que $D$ ramifica en $v$ cuando $D_v$ es un álgebra de división. La ramificación de $D$, definida como el conjunto de los $v$ en los que $D$ ramifica, es un invariante de $D$.

**Ejercicio 10.** Probar que si $p \nmid 2ab$ entonces $(a, b)_{\mathbb{Q}_p} \cong M_2(\mathbb{Q}_p)$, y concluir que la ramificación de $(a, b)_{\mathbb{Q}}$ es finita.

**Teorema 1.12** (Clasificación global)**.** *Sean $D$ y $D'$ álgebras de cuaterniones sobre $\mathbb{Q}$.*

1. $D \cong D' \iff D_v \cong D'_v$ *para todo* $v \iff D$ *y* $D'$ *tienen igual ramificación.*

2. *La ramificación de $D$ es un conjunto finito de cardinal par.*

3. *Dado un conjunto de lugares que sea finito y de cardinal par, existe un álgebra de cuaterniones sobre $\mathbb{Q}$ con esa ramificación.*

*Prueba.* Ver Serre (1973, Capítulos III-IV) o Voight (s.f., Capítulo 14). □

De acuerdo al teorema, dado un primo $N$ (finito) existe un álgebra de cuaterniones ramificada en es $\{N, \infty\}$.

**Ejercicio 11.** Sea $N$ un primo (finito). Probar que:

1. El álgebra $(-1, -1)_{\mathbb{Q}}$ ramifica en $\{2, \infty\}$.

2. Si $N \equiv 3 \pmod 4$, el álgebra $(-1, -N)_{\mathbb{Q}}$ ramifica en $\{N, \infty\}$.

3. Si $N \equiv 5 \pmod 8$, el álgebra $(-2, -N)_{\mathbb{Q}}$ ramifica en $\{N, \infty\}$.

4. Si $N \equiv 1 \pmod 8$, existe un primo $q \equiv 3 \pmod 4$ tal que $-q$ no es un cuadrado módulo $N$, y el álgebra $(-q, -N)_{\mathbb{Q}}$ ramifica en $\{N, \infty\}$.

Cuando $F$ es un cuerpo de números (extensión finita de $\mathbb{Q}$), los lugares de $F$ corresponden a las inmersiones de $F$ en $\mathbb{R}$ o en $\mathbb{C}$ (lugares arquimedianos) y a los ideales primos en su anillo de enteros (lugares no arquimedianos). Las completaciones de $F$ con respecto a los lugares arquimedianos son $\mathbb{R}$ o $\mathbb{C}$, y las completaciones con respecto a los lugares no arquimedianos son extensiones finitas de los $\mathbb{Q}_p$.

Los resultados mencionados se generalizan para $F$ de la siguiente manera: si $D$ es un álgebra de cuaterniones sobre $F$, la ramificación de $D$ (el conjunto de lugares de $F$ donde $D$ es de división) es un conjunto finito de cardinal par que no contiene ningún lugar complejo. Además, dado un conjunto de lugares de $F$ que sea finito y de cardinal par, y que no contenga ningún lugar complejo, existe un álgebra de cuaterniones sobre $F$ con esa ramificación, única a menos de isomorfismo. Para la demostración consultar Voight (ibíd., Sección 14.6).

**1.2 Órdenes de cuaterniones.** Sea $D$ un álgebra de cuaterniones sobre $\mathbb{Q}$. Un *retículo* en $D$ es un subgrupo $M \subset D$ tal que $M$ es finitamente generado y tal que $\mathbb{Q}M = D$. Equivalentemente $M$ es un $\mathbb{Z}$-módulo libre de rango 4; en otras palabras existe una base $\{d_1, d_2, d_3, d_4\}$ de $D$ que genera $M$ como $\mathbb{Z}$-módulo. Recíprocamente, cualquier base de $D$ genera un retículo en $D$. Usaremos la notación $[d_1, d_2, d_3, d_4]$ para referirnos al retículo generado por esos elementos. La *norma de un retículo* $M$ se define como

$$\nu(M) := \operatorname{mcd}\{\nu(d) \ : \ d \in M\},$$

que existe pues $M$ es finitamente generado. En efecto si $\{d_1, d_2, d_3, d_4\}$ es una base de $M$, entonces $\nu(M) = \operatorname{mcd}\{\nu(d_1), \nu(d_2), \nu(d_3), \nu(d_4)\}$.

Un elemento $x \in D$ es *integral* si $\tau(x), \nu(x) \in \mathbb{Z}$. Equivalentemente $\mathbb{Z}[x]$ es finitamente generado como $\mathbb{Z}$-módulo. Un *orden* en $D$ es un retículo $R \subset D$ que es un subanillo, esto es tal que $1 \in R$ y $R$ es cerrado por el producto. Los elementos de $R$ deben ser integrales, pues si $x \in R$ entonces $\mathbb{Z}[x] \subset R$ debe ser finitamente generado; sin embargo, *¡el conjunto de los elementos integrales de D no forma un anillo!*

Por ejemplo, tanto $x = \begin{pmatrix} 1 & 1/2 \\ 0 & 1 \end{pmatrix}$ como $y = \begin{pmatrix} 1 & 0 \\ 1/2 & 1 \end{pmatrix}$ son elementos integrales en $M_2(\mathbb{Q})$, pero $x + y$ no lo es. En efecto, tanto $\mathbb{Z}[x]$ como $\mathbb{Z}[y]$ son finitamente generados como $\mathbb{Z}$-módulos, pero no así $\mathbb{Z}[x + y]$.

Este fenómeno no es exclusivo de las álgebras de matrices:

**Ejercicio 12.** Sea $D = (-1, -1)_{\mathbb{Q}}$, que es un álgebra de división. Mostrar que existen elementos $x, y$ integrales tales que $x + y$ no es integral.
*Sugerencia:* $x = x_0 + x_1 i + x_2 j + x_3 k$ *puede ser integral con* $x_i \notin \mathbb{Z}$, *por ejemplo* $x = \frac{3i+4k}{5}$ *es integral ya que* $3^2 + 4^2 = 25$.

Si $R$ es un orden en $D$ se define su *determinante* $\det R := \det(\tau(d_i^\star\, d_j)) \in \mathbb{Z}$ donde $\{d_1, d_2, d_3, d_4\}$ es una base de $R$. El siguiente ejercicio muestra que $\det R$ no depende de la elección de la base y que $\det R$ es un cuadrado.

**Ejercicio 13.** Sea $R$ un orden en $D = (a, b)_{\mathbb{Q}}$.

1. Probar que $\det R$ no depende de la base.
2. Si $R' \subseteq R$ entonces $\det R' = [R : R']^2 \det R$.
3. Si $a, b \in \mathbb{Z}$, entonces $[1, i, j, k]$ es un orden con determinante $(4ab)^2$.
4. El determinante de cualquier orden en $D$ es un cuadrado.

**Ejemplo 1.13.** Sean $a, b \in \mathbb{Z}$ con $a > 0$, consideramos el orden $R = [1, i, j, k]$ en $(a, b)_{\mathbb{Q}}$. En el Ejemplo 1.3 vimos que podemos representar $(a, b)_{\mathbb{Q}} \hookrightarrow M_2(\mathbb{R})$. Esta inmersión preserva la norma (que en matrices es el determinante) y entonces la imagen de $\{x \in R \ : \ \nu(x) = 1\}$ en $M_2(\mathbb{R})$ es un subgrupo discreto $\Gamma(a, b) \subseteq \mathrm{SL}_2(\mathbb{R})$. Concretamente

$$\Gamma(a, b) = \left\{ \begin{pmatrix} x_0 + x_1\sqrt{a} & x_2 + x_3\sqrt{a} \\ (x_2 - x_3\sqrt{a})\, b & x_0 - x_1\sqrt{a} \end{pmatrix} \ : \ x_0^2 - a\, x_1^2 - b\, x_2^2 + ab\, x_3^2 = 1, x_i \in \mathbb{Z} \right\}.$$

Cuando $(a,b)_{\mathbb{Q}}$ es un álgebra de división, puede probarse que $|\tau(x)| > 2$ para todo $x \in \Gamma(a,b)$, $x \neq 1$; esto implica que $\Gamma(a,b)\backslash H$ es compacto, donde $H$ es el semiplano superior. Comparar con el Ejemplo 1.2 de las notas de Miatello (2020). Ver también la Proposición 1.2 en las notas de Harris (2020).

Definimos el *discriminante* (reducido) de $R$ como $\operatorname{disc} R := \sqrt{\det R} > 0$. Cuando $R' \subset R$ tenemos que $\operatorname{disc} R' = [R' : R] \operatorname{disc} R$ y como el discriminante es un entero positivo deducimos que existen órdenes maximales (aquellos con discriminante minimal).

**Proposición 1.14.** *Si $R$ es un orden maximal en $D$ entonces* $\operatorname{disc} R$ *es el producto de todos los primos donde $D$ ramifica. En particular todos los órdenes maximales tienen igual discriminante, que es libre de cuadrados; a este número lo denotamos* $\operatorname{disc} D$.

*Prueba.* El resultado puede probarse estudiando los órdenes maximales en álgebras de cuaterniones sobre $\mathbb{Q}_p$, mostrando que la valuación de su discriminante es 0 para álgebras de matrices y 1 para álgebras de división. □

**Ejemplo 1.15.** Sea $D = (-1,-1)_{\mathbb{Q}}$. Tenemos un orden $R' = [1, i, j, k]$ (cuaterniones de Lipschitz) pero su discriminante es 4, entonces $R$ no es maximal. Sea $\rho = \frac{1+i+j+k}{2}$, que es integral, entonces $R = R' + \mathbb{Z}\,\rho$ es un orden y es maximal (cuaterniones de Hurwitz).

**Ejercicio 14.** Probar las afirmaciones del ejemplo 1.15.

**Ejemplo 1.16.** Sea $D = (-1,-N)_{\mathbb{Q}}$ donde $N \equiv 3 \pmod 4$ es primo. Entonces el retículo $[1, i, \frac{1+j}{2}, \frac{i+k}{2}]$ es un anillo y tiene discriminante $N$, por lo tanto es maximal.

**Ejercicio 15.** Probar las afirmaciones del ejemplo 1.16.

**Ejercicio 16.** Encontrar órdenes maximales para las otras álgebras del ejercicio 11.

Fijemos un orden $R$ en un álgebra de cuaterniones $D$. Un *ideal a derecha para* $R$ es un retículo $I \subset D$ tal que $IR = I$.

En general resulta necesario trabajar con ideales para $R$ que sean localmente principales. Supongamos que $D$ es un álgebra de cuaterniones sobre $\mathbb{Q}$. Para cada primo $p$ definimos $R_p := R \otimes \mathbb{Z}_p$, e $I_p := I \otimes \mathbb{Z}_p$.

**Ejercicio 17.** Probar que

1. $R_p$ es un orden en $D_p = D \otimes \mathbb{Q}_p$.
2. $I_p$ es un ideal a derecha para $R_p$.
3. Si $R$ es maximal entonces $R_p$ es maximal.
4. Si $R$ es maximal entonces existe $x_p \in D_p$ tal que $I_p = x_p\, R_p$.
   *Sugerencia: considerar $x_p \in I_p$ con norma de valuación mínima.*

Un ideal $I$ a derecha para $R$ es *localmente principal* si para todo $p$ existe $x_p \in D_p$ tal que $I_p = x_p\, R_p$. Por el ejercicio, cuando $R$ es maximal todos los ideales son localmente principales. Consultar Voight (s.f., Capítulo 16) por más detalles, particularmente el Teorema 16.1.3 y la Sección 16.2.

Denotamos $\mathcal{I}(R)$ al conjunto de ideales a derecha propios para $R$. El grupo multiplicativo $D^\times$ actúa en $\mathcal{I}(R)$ por multiplicación a la izquierda. Las órbitas de $\mathcal{I}(R)$ por esta acción se llaman *clases de ideales*. Llamaremos $h(R)$ al número de clases de ideales de $R$.

**Teorema 1.17.** *El número de clases de ideales $h(R)$ es finito.*

*Prueba.* Ver Voight (ibíd., Sección 17.5), en particular el Teorema 17.5.6. □

En general hay fórmulas para el número de clases de ideales de un orden de cuaterniones. A modo de ejemplo:

**Proposición 1.18.** *Si $R$ es un orden de discriminante $N$ primo, entonces*

$$h(R) = \begin{cases} [N/12] & si\ N \equiv 1 \pmod{12} \\ [N/12]+1 & si\ N \not\equiv 1, 11 \pmod{12} \\ [N/12]+2 & si\ N \equiv 11 \pmod{12} \end{cases}$$

*Prueba.* Es una consecuencia de la fórmula de la masa de Eichler. Ver Voight (ibíd., Capítulos 25 y 30), en particular los teoremas 25.3.18 y 30.1.5. □

## 2 Formas modulares cuaterniónicas

En esta sección haremos una introducción a las formas modulares cuaterniónicas en un caso particular. En las notas de Harris (2020) se verá una definición más general.

Sea $D$ un álgebra de cuaterniones sobre $\mathbb{Q}$; suponemos que $D$ es definida (es decir, $D_\infty = D \otimes \mathbb{R}$ es un álgebra de división).

**Definición 2.1.** Dado un orden $R \subset D$, una *forma modular cuaterniónica* para $R$ es una función $f : \mathcal{I}(R) \to \mathbb{C}$ tal que $f(d\,I) = f(I)$ para todo $d \in D^\times$.

Denotaremos $\mathcal{M}(R)$ al espacio de formas modulares cuaterniónicas para $R$. Dado un ideal $I$, la función característica de la clase de $I$ es una forma modular cuaterniónica que denotaremos $[I]$. A efectos computacionales es conveniente fijar un conjunto de representantes de las clases de ideales $\{I_1, \dots, I_h\}$, de modo que $\{[I_1], \dots, [I_h]\}$ es una base de $\mathcal{M}(R)$.

Para cada ideal $I \in \mathcal{I}(R)$ el grupo $\Gamma_I = \{d \in D^\times : d\,I = I\}/\mathbb{Z}^\times$ es finito, ya que es discreto dentro de $D_\infty^\times/\mathbb{R}^\times \cong SO_3(\mathbb{R})$ que es compacto. Denotamos $w_I = \#\Gamma_I$, que depende solamente de la clase de $I$.

Definimos un producto interno en $\mathcal{M}(R)$ que está dado en la base por

$$\langle [I], [J] \rangle := \tfrac{1}{2}\#\{d \in D^\times : I = d\,J\} = \begin{cases} w_I & \text{si } [I] = [J], \\ 0 & \text{si } [I] \neq [J]. \end{cases}$$

Notar que $\{[I_1], \dots, [I_h]\}$ es una base ortogonal de $\mathcal{M}(R)$.

El *grado* se define como el funcional lineal en $\mathcal{M}(R)$ tal que $\mathrm{gr}([I]) = 1$. Sea

$$e_0 := \sum_{i=1}^{h} \frac{1}{w_{I_i}} [I_i]\,.$$

Entonces $\mathrm{gr}(f) = \langle f, e_0 \rangle$.

Decimos que $f$ es *cuspidal* si es ortogonal a $e_0$, es decir si $\mathrm{gr}(f) = 0$.

**2.1 Operadores de Hecke.** A continuación vamos a definir una familia de operadores en $\mathcal{M}(R)$ que llamaremos *operadores de Hecke*. Dado un ideal $I \in \mathcal{I}(R)$ y $n \geqslant 1$ un entero sea

$$\mathcal{T}_n(I) := \left\{ I' \in \mathcal{I}(R) \ : \ I' \subset I, \quad \nu(I') = n\,\nu(I) \right\}.$$

El operador de Hecke $t_n : \mathcal{M}(R) \to \mathcal{M}(R)$ se define, en la base, como

$$t_n[I] := \sum_{I' \in \mathcal{T}_n(I)} [I'].$$

Observar que $t_1 = \mathrm{id}_{\mathcal{M}(\mathrm{R})}$.

**Lema 2.2.** *Sean $I, J \in \mathcal{I}(R)$ y $n \geqslant 1$. Entonces $J \in \mathcal{T}_n(I) \Leftrightarrow n\,I \in \mathcal{T}_n(J)$.*

*Prueba.* Supongamos que $J = d\,I$ con $d \in D^\times$. Si $d\,I \in \mathcal{T}_n(I)$ entonces $\nu(d) = n$ y $d\,I \subset I$. Entonces también $d^\star I \subset I$ y concluimos que $n\,I = d\,d^\star I \subset d\,I$, y comparando normas concluimos que $n\,I \in \mathcal{T}_n(d\,I)$. Cuando $I$ y $J$ no son equivalentes se prueba del mismo modo pero localmente, usando que $I$ y $J$ son localmente principales y por lo tanto localmente equivalentes. □

**Proposición 2.3.** *Los operadores de Hecke en $\mathcal{M}(R)$ satisfacen*

1. *$t_n$ es autoadjunto.*
2. *Si $(n, n') = 1$ entonces $t_n\,t_{n'} = t_{n'}\,t_n = t_{n\,n'}$.*
3. *Si $p \nmid \operatorname{disc} R$ es primo entonces $t_{p^{k+2}} = t_{p^{k+1}}\,t_p - p\,t_{p^k}$.*
4. *Si $(n\,n', \operatorname{disc} R) = 1$ entonces $t_n\,t_{n'} = \sum_{d \mid (n,n')} d\,t_{n\,n'/d^2}$*

*En particular los operadores $t_p$ con $p \nmid \operatorname{disc} R$ primo generan un álgebra conmutativa $\mathbb{T}_R$ que contiene todos los operadores $t_n$ con $(n, \operatorname{disc} R) = 1$. Además éstos últimos generan $\mathbb{T}_R$ como $\mathbb{Z}$-módulo.*

*Prueba.* Calculamos

$$\begin{aligned}
\langle t_n[I], [J] \rangle &= \sum_{I' \in \mathcal{T}_n(I)} \big\langle [I'], [J] \big\rangle \\
&= \sum_{I' \in \mathcal{T}_n(I)} \tfrac{1}{2}\,\#\big\{ d \in D^\times \ : \ I' = d\,J \big\} \\
&= \tfrac{1}{2}\,\#\{ d \in D^\times \ : \ d\,J \in \mathcal{T}_n(I) \}
\end{aligned}$$

El Lema 2.2 implica que $d\,J \in \mathcal{T}_n(I) \Leftrightarrow n\,I \in \mathcal{T}_n(d\,J)$. Entonces la última expresión es

$$\begin{aligned}
&= \tfrac{1}{2}\,\#\{ d \in D^\times \ : \ n\,I \in \mathcal{T}_n(d\,J) \} \\
&= \tfrac{1}{2}\,\#\big\{ d \in D^\times \ : \ n\,d^{-1}\,I \in \mathcal{T}_n(J) \big\} \\
&= \tfrac{1}{2}\,\#\big\{ d' \in D^\times \ : \ d'\,I \in \mathcal{T}_n(J) \big\} \\
&= \langle [I], t_n[J] \rangle
\end{aligned}$$

Las afirmaciones (2) y (3) se prueban localmente. La afirmación (4) es un ejercicio a partir de (2) y (3). □

**Corolario 2.4.** *El espacio $\mathcal{M}(R)$ tiene una base ortogonal de vectores propios para $\mathbb{T}_R$.*

*Prueba.* Como $\mathbb{T}_R$ es un álgebra conmutativa de operadores autoadjuntos, este resultado es una consecuencia del Teorema Espectral. □

**Ejercicio 18.** Utilizando el Teorema Espectral para un operador autoadjunto, completar los detalles de la demostración del Corolario.
*Sugerencia: probar primero que si tenemos dos operadores autoadjuntos que conmutan entre sí, entonces el espacio se puede escribir como suma directa ortogonal de subespacios propios comunes.*

**Proposición 2.5.** *Si $p \nmid \operatorname{disc} R$ primo, entonces $\#\mathcal{T}_p(I) = p + 1$.*

*Prueba.* Se prueba localmente. Como $I$ es localmente principal, tenemos que $I_p = x_p\, R_p$ para algún $x_p \in D_p$. Entonces $I' \in \mathcal{T}_p(I)$ si y sólo si

$$\begin{cases} I'_q = I_q & \text{para todo } q \neq p \\ I'_p = x_p\, \alpha\, R_p & \text{con } \alpha \in R_p,\ \nu(\alpha) = p. \end{cases}$$

Además $x_p\, \alpha'\, R_p = x_p\, \alpha\, R_p$ si y sólo si $\alpha' \in \alpha R_p^\times$. Por lo tanto

$$\#\mathcal{T}_p(I) = \#\big\{\alpha \in R_p \ :\ \nu(\alpha) = p\big\}/R_p^\times .$$

Cuando $p \nmid \operatorname{disc} R$ tenemos que $R_p \cong M_2(\mathbb{Z}_p)$ y el resultado se deduce del siguiente ejercicio. □

**Ejercicio 19.** Sea $p$ primo. Probar que $\big\{\alpha \in M_2(\mathbb{Z}_p) \ :\ \det(\alpha) = p\big\}/\mathrm{GL}_2(\mathbb{Z}_p)$ tiene $p + 1$ elementos.
*Sugerencia:* $\big\{\left(\begin{smallmatrix} p & 0 \\ 0 & 1 \end{smallmatrix}\right), \left(\begin{smallmatrix} 1 & 0 \\ 0 & p \end{smallmatrix}\right), \left(\begin{smallmatrix} 1 & 0 \\ 1 & p \end{smallmatrix}\right), \dots, \left(\begin{smallmatrix} 1 & 0 \\ p-1 & p \end{smallmatrix}\right)\big\}$ *es un conjunto de representantes.*

La proposición muestra que $t_p$ es homogéneo de grado $p + 1$ en el sentido de que $\operatorname{gr}(t_p f) = (p + 1)\operatorname{gr}(f)$. Entonces $e_0$ es un vector propio para $\mathbb{T}_R$ con $t_p\, e_0 = (p + 1)\, e_0$. Cualquier otro vector propio será ortogonal a $e_0$, es decir cuspidal.

**Ejemplo 2.6.** Consideramos el álgebra $D = (-1, -11)_{\mathbb{Q}}$ ramificada en $\{11, \infty\}$ y el orden maximal $R = [1, i, \frac{1+j}{2}, \frac{i+k}{2}]$ de discriminante 11. El número de clases es 2 (ver Proposición 1.18), y un conjunto de representantes es

$$I_1 = [1, i, \tfrac{1+j}{2}, \tfrac{i+k}{2}]$$
$$I_2 = [2, 2i, i + \tfrac{1+j}{2}, 1 + i + \tfrac{i+k}{2}]$$

Los estabilizadores son $\Gamma_{I_1} = \{\pm 1, \pm i\}$ y $\Gamma_{I_2} = \left\{\pm 1, \pm\frac{2-i+k}{4}, \pm\frac{2+i-k}{4}\right\}$, de modo que $w_{I_1} = 2$ y $w_{I_2} = 3$. Se puede calcular

$$\mathcal{T}_2(I_1) = \left\{(1 + i)\, I_1,\ I_2,\ \tfrac{2-i-k}{4}\, I_2\right\},$$

$$\mathcal{T}_2(I_2) = \left\{2\, I_1,\ \tfrac{2+i-k}{2}\, I_1,\ \tfrac{2-i+k}{2}\, I_1\right\}.$$

Entonces el operador de Hecke $t_2$ está dado por

$$t_2[I_1] = [I_1] + 2[I_2]\,, \qquad t_2[I_2] = 3[I_1]\,,$$

y tiene vectores propios $e_0 = \frac{1}{2}[I_1] + \frac{1}{3}[I_2]$ con valor propio 3, y $f_1 = [I_1] - [I_2]$ con valor propio $-2$. Podemos calcular otros operadores de Hecke, que en la base $\{[I_1], [I_2]\}$ resultan:

$$t_2 = \left(\begin{smallmatrix}1 & 3\\ 2 & 0\end{smallmatrix}\right), \quad t_3 = \left(\begin{smallmatrix}2 & 3\\ 2 & 1\end{smallmatrix}\right), \quad t_5 = \left(\begin{smallmatrix}4 & 3\\ 2 & 3\end{smallmatrix}\right), \quad t_7 = \left(\begin{smallmatrix}4 & 6\\ 4 & 2\end{smallmatrix}\right), \quad \ldots \quad .$$

El vector propio $e_0$ tiene valores propios 3, 4, 6, 8, …, como ya habíamos observado, mientras que $f_1$ tiene valores propios -2, -1, 1, -2, ….

**2.2 Correspondencia de Eichler.** Para evitar dificultades técnicas nos limitaremos al caso de un orden $R$ de discriminante $N$ primo en un álgebra definida. Necesariamente $R$ es un orden maximal en un álgebra de cuaterniones $D$ ramificada en $\{N, \infty\}$.

Hemos visto que el espacio $\mathcal{M}(R)$ tiene una acción por operadores autoadjuntos del álgebra conmutativa $\mathbb{T} = \mathbb{T}_R$, que es similar a la acción de Hecke en espacios de formas modulares clásicas. Para una introducción a las formas modulares clásicas y a sus operadores de Hecke, referimos al lector a las notas de Miatello (2020).

Observemos que las relaciones de la Proposición 2.3 son las mismas que para los operadores de Hecke en formas modulares de peso 2 para $\Gamma_0(N)$. Por ejemplo, ver las proposiciones 3.5 y 3.6 en Miatello (ibíd.).

Eichler calcula la traza de $t_n$ actuando en $\mathcal{M}(R)$ y por otra parte calcula la traza de $T_n$ actuando en $\mathfrak{M}_2(\Gamma_0(N))$. Comparando ambas obtiene el siguiente resultado.

**Teorema 2.7** (Eichler (1955))**.** *Para todo $n \geqslant 1$ tenemos*

$$\operatorname{Tr}(t_n \leftrightarrow \mathcal{M}(R)) = \operatorname{Tr}(T_n \leftrightarrow \mathfrak{M}_2(\Gamma_0(N)))$$

*Prueba.* Eichler tiene una presentación autocontenida de este resultado y algunas aplicaciones en Eichler (1973) □

Como consecuencia de este resultado, y puesto que ambas álgebras de Hecke son semisimples y están generadas como $\mathbb{Z}$-módulo por los operadores de Hecke, se deduce que son isomorfas y que hay una correspondencia

$$\{\text{vectores propios en } \mathcal{M}(R)\}/_{\mathbb{C}^\times} \longleftrightarrow \{\text{vectores propios en } \mathfrak{M}_2(\Gamma_0(N))\}/_{\mathbb{C}^\times}$$

donde formas correspondientes tienen los mismos valores propios. Otra forma de enunciar lo mismo es decir que existe un isomorfismo $\mathcal{M}(R) \cong \mathfrak{M}_2(\Gamma_0(N))$ que preserva la acción de Hecke. Sin embargo este isomorfismo no es canónico.

Sabemos que en $\mathfrak{M}_2(\Gamma_0(N))$ los espacios propios tienen dimensión 1 (no hay formas viejas porque $N$ es primo y $\mathfrak{M}_2(\Gamma_0(1)) = \{0\}$), y lo mismo vale entonces para $\mathcal{M}(R)$.

Veremos ahora una construcción de formas modulares de peso 2 para $\Gamma_0(N)$ a partir de formas modulares cuaterniónicas. Esto permite una realización explícita de la correspondencia de Eichler como una función $\mathbb{T}$-bilineal.

**Definición 2.8.** Sean $f, g \in \mathcal{M}(R)$. Definimos

$$\phi(f,g) := \frac{\operatorname{gr}(f)\cdot\operatorname{gr}(g)}{2} + \sum_{n\geqslant 1} \langle t_n f, g\rangle\, q^n \qquad (q = e^{2\pi i z}).$$

**Proposición 2.9.** *$\phi(f,g)$ es una forma modular de peso 2 para $\Gamma_0(N)$ y*

$$\phi(t_n\, f, g) = \phi(f, t_n\, g) = T_n\, \phi(f,g) \qquad \textit{si } (n, \operatorname{disc} R) = 1.$$

*En otras palabras*

$$\phi : \mathcal{M}(R) \underset{\mathbb{T}}{\otimes} \mathcal{M}(R) \to \mathfrak{M}_2(\Gamma_0(N))$$

*es $\mathbb{T}$-equivariante, y como $\mathcal{M}(R)$ es un $\mathbb{T}$-módulo libre de rango 1 concluye que $\phi$ es un isomorfismo de $\mathbb{T}$-módulos.*

*Prueba.* Sean $I, J \in \mathcal{I}(R)$. Consideramos el retículo $M = \{d \in D \;:\; d\,J \subset I\}$. En la demostración de la Proposición 2.3 calculamos

$$\langle t_n[I], [J]\rangle = \tfrac{1}{2}\#\{d \in D^\times \;:\; d\,J \in \mathcal{T}_n(I)\}\,.$$

Ahora $d\,J \in \mathcal{T}_n(I)$ es equivalente a $d \in M$ con $\nu(d) = n\,\nu(M)$, es decir que

$$\langle t_n[I], [J]\rangle = \tfrac{1}{2}\#\{d \in M \;:\; \nu(d) = n\,\nu(M)\}\,.$$

Luego

$$\phi([I],[J]) = \tfrac{1}{2}\sum_{d\in M} q^{Q(d)}$$

es la serie theta asociada a la forma cuadrática $Q(d) = \nu(d)/\nu(M)$, que es una forma modular en $\mathfrak{M}_2(\Gamma_0(N))$. En general $\phi(f,g)$ es combinación lineal de estas series theta.

La igualdad $\phi(t_n\, f, g) = \phi(f, t_n\, g)$ se deduce fácilmente pues $t_n$ es autoadjunto y es homogéneo respecto al grado. La última igualdad basta probarla con $n = p$ primo. Usando la fórmula de $T_p$ en términos de coeficientes de Fourier calculamos

$$T_p\,\phi(f,g) = (p+1)\,\frac{\operatorname{gr}(f)\cdot\operatorname{gr}(g)}{2} + \sum_{n\geqslant 1}\big(\langle t_{np}\, f, g\rangle + p\,\langle t_{n/p}\, f, g\rangle\big)\, q^n$$

bajo la convención que $t_{n/p} = 0$ si $p \nmid n$. Por otra parte

$$\phi(t_p\, f, g) = \frac{\operatorname{gr}(t_p\, f)\cdot\operatorname{gr}(g)}{2} + \sum_{n\geqslant 1}\langle t_n\, t_p\, f, g\rangle\, q^n\,.$$

El resultado se sigue pues por la Proposición 2.3 tenemos $t_n\, t_p = t_{np} + p\, t_{n/p}$. □

*Observación.* Para cada $f \in \mathcal{M}(R)$ tenemos una transformación lineal $\mathbb{T}$-equivariante $\phi(\cdot, f) : \mathcal{M}(R) \to \mathfrak{M}_2(\Gamma_0(N))$, pero no necesariamente es un isomorfismo. En todo caso existe $f$ tal que $\phi(\cdot, f)$ es un isomorfismo $\mathcal{M}(R) \cong \mathfrak{M}_2(\Gamma_0(N))$, pero no será canónico.

**Ejemplo 2.10.** Continuando con el ejemplo 2.6, las series theta asociadas a la base de $\mathcal{M}(R)$ son $\theta_{ij} = \phi([I_i], [I_j])$:

$$\theta_{11} = \tfrac{1}{2} + 2q + 2q^2 + 4q^3 + 10q^4 + 8q^5 + 16q^6 + 8q^7 + 18q^8 + 14q^9 + O(q^{10})$$
$$\theta_{12} = \tfrac{1}{2} + 6q^2 + 6q^3 + 6q^4 + 6q^5 + 12q^6 + 12q^7 + 18q^8 + 18q^9 + O(q^{10})$$
$$\theta_{22} = \tfrac{1}{2} + 3q + 3q^3 + 12q^4 + 9q^5 + 18q^6 + 6q^7 + 18q^8 + 12q^9 + O(q^{10})$$

Con estas series theta podemos calcular

$$\begin{aligned} E_0 &= \phi(e_0, [I_1]) = \tfrac{1}{2}\,\theta_{11} + \tfrac{1}{3}\,\theta_{12} = \phi(e_0, [I_2]) = \tfrac{1}{2}\,\theta_{12} + \tfrac{1}{3}\,\theta_{22} \\ &= \tfrac{5}{12} + q + 3q^2 + 4q^3 + 7q^4 + 6q^5 + 12q^6 + 8q^7 + 15q^8 + 13q^9 + O(q^{10}) \end{aligned}$$

$$\begin{aligned} F_1 &= \phi(f_1, \tfrac{1}{2}[I_1]) = \tfrac{1}{2}\,(\theta_{11} - \theta_{12}) = \phi(f_1, -\tfrac{1}{3}[I_2]) = \tfrac{1}{3}\,(\theta_{22} - \theta_{12}) \\ &= q - 2q^2 - q^3 + 2q^4 + q^5 + 2q^6 - 2q^7 - 2q^9 + O(q^{10}) \end{aligned}$$

$$0 = \phi(f_1, e_0) = \tfrac{1}{2}\,\theta_{11} - \tfrac{1}{6}\,\theta_{12} - \tfrac{1}{3}\,\theta_{22}$$

Las dos primeras son las dos autoformas modulares normalizadas de peso 2 para $\Gamma_0(11)$. La forma $E_0$ es una serie de Eisenstein, y $F_1$ es una autoforma cuspidal que corresponde a la curva elíptica $y^2 + y = x^3 - x^2 - 10x - 20$.

*Observación.* Notamos que $e_0 = \frac{1}{2}\,[I_1] + \frac{1}{3}\,[I_2] \equiv \frac{1}{2}\,[I_1] - \frac{1}{2}\,[I_2] = \frac{1}{2} f_1 \pmod 5$. Se deduce que $E_0 = \phi(e_0, [I_1]) \equiv \phi(\frac{1}{2}\, f_1, [I_1]) = F_1 \pmod 5$. En efecto, $E_0 - F_1 = \frac{5}{6}\,\theta_{12}$. Esta congruencia está relacionada con el hecho que la curva elíptica $y^2 + y = x^3 - x^2 - 10x - 20$ tiene torsión de orden 5.

**2.3 Formas modulares con valores vectoriales.** Para generalizar esta formulación de formas modulares cuaterniónicas a otros pesos podemos emplear formas con valores vectoriales.

Sea $(V, \rho)$ una representación irreducible unitaria de $D_\infty^\times/\mathbb{R}^\times$. Entonces definimos

$$\mathcal{M}_\rho(R) := \{ f : \mathcal{I}(R) \to V \ : \ f(dI) = \rho(d) \cdot f(I) \ \forall d \in D^\times \}$$

Observemos que $f(I) \in V^{\Gamma_I}$ para todo $I \in \mathcal{I}(R)$.

El producto interno se define en este espacio como

$$\langle f, g \rangle := \sum_{I \in D^\times \backslash \mathcal{I}(R)} [f(I), g(I)] \cdot w_I$$

donde $[,]$ es el producto interno invariante en $V$. Es fácil ver que $[f(I), g(I)] \cdot w_I$ depende solamente de la clase de $I$.

El espacio $\mathcal{M}_\rho(R)$ tiene una descomposición ortogonal $\mathcal{M}_\rho(R) = \oplus_{i=1}^h \mathcal{M}_\rho^{[I_i]}(R)$ donde $\mathcal{M}_\rho^I(R)$ es el subespacio de formas modulares con soporte en la clase de $I$. Dado $I \in \mathcal{I}(R)$ y $v \in V^{\Gamma_I}$ existe una única $f = f_{I,v} \in \mathcal{M}_\rho^{[I]}(R)$ tal que $f(I) = v$, y no es difícil probar que esto induce un isomorfismo $\mathcal{M}_\rho^I(R) \cong V^{\Gamma_I}$. Concluimos que

$$\mathcal{M}_\rho(R) \cong \oplus_{i=1}^h V^{\Gamma_{I_i}} \ .$$

Los operadores de Hecke se definen en $\mathcal{M}^I_\rho(R)$ como

$$t_n f_{I,v} := \sum_{I' \in \mathcal{T}_n(I)} f_{I',v'}$$

donde en cada término $v'$ es la proyección ortogonal de $v$ en $V^{\Gamma_{I'}}$.

## Referencias

GONZALO TORNARÍA
tornaria@cmat.edu.uy
tornaria@gmail.com
CMAT
URUGUAY

AGRA

# FORMAS DE HILBERT

Ariel Pacetti

## Índice general

## 1 Reinterpretación de las formas modulares clásicas

**1.1 Formas modulares como formas en $\mathrm{SL}_2(\mathbb{R})$.** La principal referencia para esta parte es el libro Gelbart (1975). Recordar que las formas modulares eran formas holomorfas en el semiplano de Poincaré $\mathfrak{h} = \{z \in \mathbb{C} \, : \, \Im(z) > 0\}$. Para entender la relación de estas funciones con el cuerpo de números racionales, debemos entender como es que este cuerpo aparece en la teoría. Recordar que el grupo $\mathrm{SL}_2(\mathbb{R})$ de matrices reales de $2 \times 2$ con determinante 1 actúa en $\mathfrak{h}$ de manera transitiva, dado que si $z = x + iy$ es un elemento de $\mathfrak{h}$,

$$x + iy = \begin{pmatrix} y^{1/2} & xy^{-1/2} \\ 0 & y^{-1/2} \end{pmatrix} \cdot i$$

Luego podemos identificar $\mathfrak{h}$ con el cociente de $\mathrm{SL}_2(\mathbb{R})$ módulo el estabilizador de $i$, o sea

$$\mathfrak{h} \simeq \mathrm{SL}_2(\mathbb{R}) / \mathrm{SO}(2), \tag{1}$$

donde $\mathrm{SO}(2)$ es el grupo ortogonal, dado por matrices de la forma

$$r(\theta) = \begin{pmatrix} \cos(\theta) & \sin(\theta) \\ -\sin(\theta) & \cos(\theta) \end{pmatrix}.$$

**Ejercicio 1.** Probar que si $M = \left(\begin{smallmatrix} a & b \\ c & d \end{smallmatrix}\right)$ es una matriz de $2 \times 2$ de determinante 1 que estabiliza el punto $i$, entonces existe un $\theta \in [0, 2\pi)$ tal que $M = r(\theta)$.

Por $\Gamma$ vamos a denotar un subgrupo de congruencia, como $\Gamma_0(N)$. Para una primera lectura de estas notas, basta suponer que $\Gamma$ es todo el grupo $\mathrm{SL}_2(\mathbb{Z})$, y vamos a asumirlo en varias ocasiones (diciendo que sucede en el caso general). A una forma modular $f \in M_k(\Gamma)$, le podemos asociar una función $\phi_f : \mathrm{SL}_2(\mathbb{R}) \to \mathbb{C}$, de la siguiente manera:

$$\phi_f(g) = f(g \cdot i) j(g,i)^{-k}. \tag{2}$$

Es claro que $\phi_f$ es cero si y sólo si $f$ lo es. Notar que $\mathrm{SL}_2(\mathbb{R})$ posee una acción de $\Gamma$ a izquierda (por el que estudiamos invariancia de las formas modulares), y una acción de $\mathrm{SO}(2)$ a derecha también. Ambas acciones son de vital importancia para entender las formas modulares como formas automorfas.

**Proposición 1.** La función $\phi_f$ satisface las siguientes propiedades:

1. $\phi_f(\gamma g) = \phi_f(g)$ para todo $\gamma \in \Gamma$.
2. $\phi_f(gr(\theta)) = \exp(ik\theta)\phi_f(g)$ para todo $k \in \mathbb{Z}$.
3. Si $f$ es cuspidal, entonces $\phi_f$ es una función acotada y
$$\int_{\Gamma \backslash \mathrm{SL}_2(\mathbb{R})} |\phi_f(g)|^2 dg < \infty.$$
4. Si $f$ es cuspidal, entonces $\phi_f$ es cuspidal, lo que quiere decir que para todo $g \in \mathrm{SL}_2(\mathbb{R})$ y todo $\gamma \in \mathrm{SL}_2(\mathbb{Z})$,
$$\int_0^1 \phi_f\left(\gamma \left(\begin{smallmatrix} 1 & x \\ 0 & 1 \end{smallmatrix}\right) g\right) dx = 0.$$

*Prueba.* 1. Por definición,

$$\begin{aligned}\phi_f(\gamma g) &= f(\gamma g \cdot i) j(\gamma g, i)^{-k} \\ &= j(\gamma, g \cdot i)^k f(g \cdot i) j(\gamma, g \cdot i)^{-k} j(g,i)^{-k} = \phi_f(g).\end{aligned}$$

2. Recordar que $\mathrm{SO}(2)$ es justamente el normalizador de $i$, con lo cual

$$\begin{aligned}\phi_f(gr(\theta)) &= f(gr(\theta) \cdot i) j(gr(\theta), i)^{-k} \\ &= f(g \cdot i) j(g, r(\theta) \cdot i)^{-k} j(r(\theta), i)^{-k} = \phi_f(g) \exp(-i\theta)^{-k}.\end{aligned}$$

3. Como se mencionó en las notas de R. Miatello, la compactificación de $\Gamma \backslash \mathfrak{h}$ se obtiene agregando la cúspide de infinito al semiplano complejo superior (o finitos puntos para subgrupos de congruencia) y definiendo un abierto base de él. Si $f(z) \in M_k(\Gamma)$, entonces la función $F(z) = \Im(z)^{k/2}|f(z)|$ es invariante por la acción de $\Gamma$. Ahora si $f(z)$ es cuspidal, decrece de manera exponencial al acercarse a la cúspide de infinito, con lo cual $F(z)$ tiende a cero al acercarnos a tal cúspide. Luego $F(z)$ se extiende de manera continua a un compacto, y por lo tanto está acotada.

Ahora, si tomamos $z = g \cdot i = \frac{ai+b}{ci+d}$, donde $g = \left(\begin{smallmatrix} a & b \\ c & d \end{smallmatrix}\right) \in \mathrm{SL}_2(\mathbb{R})$, entonces $\Im(z) = \frac{1}{|ci+d|^2}$, con lo cual

$$\left|F\left(\frac{ai+b}{ci+d}\right)\right| = \left|f\left(\frac{ai+b}{ci+d}\right)\right| |ci+d|^{-k} = |\phi_f(g)|.$$

Para poder siquiera enunciar la finitud de la integral (que implica que la función está en un espacio $L^2$), debemos decir con respecto a que medida estamos integrando. $\mathrm{SL}_2(\mathbb{R})$ lo podemos parametrizar por $(x, y, \theta)$, donde $\theta \in [0, 2\pi)$ (corresponde al compacto), $x \in (-\infty, \infty)$ e $y \in (0, \infty)$. Luego la medida de Haar es la medida invariante del semiplano de Poincaré $\frac{dxdy}{y^2}$ multiplicada por la medida $\frac{1}{2\pi}d\theta$ (esta normalización hace con que el compacto tenga medida 1). Así,

$$\begin{aligned}\int_{\Gamma\backslash \mathrm{SL}_2(\mathbb{R})} |\phi_f(g)|^2 dg &= \frac{1}{2\pi}\int_{\Gamma\backslash\mathfrak{H}}\int_0^{2\pi} \left|\phi_f\left(\begin{pmatrix} y^{1/2} & xy^{-1/2} \\ 0 & y^{1/2}\end{pmatrix} r(\theta)\right)\right|^2 \frac{dxdy}{y^2} d\theta \\ &= \int_{\Gamma\backslash\mathfrak{H}} |f(x+iy)|^2\, y^k\, \frac{dxdy}{y^2} = ||f||,\end{aligned}$$

donde $||f||$ denota la norma de Petersson introducida en las notas de R. Miatello.

4. Como estamos asumiendo que $\Gamma = \mathrm{SL}_2(\mathbb{Z})$, la acción a derecha de $\gamma$ es trivial. Luego, si llamamos $z = g \cdot i$, tenemos que

$$\int_0^1 \phi_f\left(\gamma \left(\begin{smallmatrix} 1 & x \\ 0 & 1\end{smallmatrix}\right) g\right) dx = j(g,i)^{-k}\int_0^1 f(z+x)dx = a_0(f) = 0.$$

En el caso general, nos queda la integral no de $f$, sino de $f|[\gamma]_k$, que corresponde a mirar la expansión de Fourier en otras cúspides. Pero el hecho de ser cuspidal implica que el primer coeficiente de Fourier en todas las cúspides es cero. □

Luego a cada forma modular $f \in S_k(\Gamma)$ le hemos asociado una función $\phi_f$ en $L^2(\Gamma\backslash \mathrm{SL}_2(\mathbb{R}))$. Como caracterizar la imagen?

Dentro del conjunto de funciones $C^\infty$ de $\Gamma\backslash \mathrm{SL}_2(\mathbb{R})$ (que son densas en $L^2(\Gamma\backslash \mathrm{SL}_2(\mathbb{R}))$), hay un operador que juega un rol fundamental, el llamado *operador de Casimir*. Este operador, en las coordenadas $(x, y, \theta)$, se puede definir como

$$\Delta = -y^2\left(\frac{\partial^2}{\partial x^2} + \frac{\partial^2}{\partial y^2}\right) - y\, \frac{\partial^2}{\partial x \partial\theta}. \tag{3}$$

**Proposición 2.** La función $f \to \phi_f$ da una biyección entre $S_k(\Gamma)$ y las funciones $\phi$ en $\mathrm{SL}_2(\mathbb{R})$ que satisfacen:

1. $\phi(\gamma g) = \phi_f(g)$ para todo $\gamma \in \Gamma$.
2. $\phi(gr(\theta)) = \exp(ik\theta)\phi_f(g)$.
3. $\Delta\phi = -\frac{k}{2}(\frac{k}{2} - 1)\phi$.
4. $\phi$ es una función acotada y cuspidal.

*Prueba.* Corroboremos que $\phi_f$ satisface todas las propiedades. Vimos que si $f \in S_k(\Gamma)$, entonces $\phi_f$ satisface todas las propiedades salvo la tercera. Para ver esta, recordar que

$$\phi_f\left(\left(\begin{smallmatrix} y^{1/2} & xy^{-1/2} \\ 0 & y^{-1/2} \end{smallmatrix}\right) r(\theta)\right) = f(x+iy)y^{k/2}e^{ik\theta}$$

Claramente:

- $\frac{\partial^2}{\partial x \partial \theta}\phi_f(z,\theta) = (ik)\frac{\partial f(z)}{\partial x} y^{k/2}e^{ik\theta}$.
- $\frac{\partial^2}{\partial x^2}\phi_f(z,\theta) = \frac{\partial^2 f(z)}{\partial x^2} y^{k/2}e^{ik\theta}$.
- $\frac{\partial^2}{\partial y^2}\phi_f(z,\theta) = (\frac{\partial^2 f(z)}{\partial y^2} y^{k/2} + k\frac{\partial f(z)}{\partial y} y^{k/2-1} + f(z)\frac{k}{2}(\frac{k}{2}-1)y^{k/2-2})e^{ik\theta}$.

Luego,

$$\Delta\phi_f = \left(-y^{k/2+2}\left(\frac{\partial^2}{\partial x^2}+\frac{\partial^2}{\partial y^2}\right) f(z) - y^{k/2+1} ik \left(\frac{\partial}{\partial x} - i\frac{\partial}{\partial y}\right) f(z) \right.$$
$$\left. -\frac{k}{2}\left(\frac{k}{2}-1\right) y^{k/2} f(z)\right) e^{ik\theta}.$$

Como $f(z)$ es holomorfa, los dos primeros términos se anulan. Para probar la recíproca, necesitamos demostrar que si $\phi$ satisface las condiciones, y definimos $f(x+iy) = \phi(g)j(g,i)^k$, donde $g = \left(\begin{smallmatrix} y^{1/2} & xy^{-1/2} \\ 0 & y^{-1/2} \end{smallmatrix}\right)$, obtenemos una forma modular. La condición de invariancia se deduce fácilmente de las fórmulas de la demostración de la Proposición 1. Lo que resulta mas complicado con esta definición es verificar justamente que la función $f$ es holomorfa. Daremos una idea de esta demostración en la siguiente sección. □

**1.2 Interpretación del operador de Casimir.** Notar que $\mathrm{SL}_2(\mathbb{R})$ actúa en $L^2(\Gamma\backslash \mathrm{SL}_2(\mathbb{R}))$ vía la acción regular a derecha, esto es $g\cdot\phi(h) = \phi(hg)$, con lo cual podemos considerar las representaciones irreducibles [1] de $\mathrm{SL}_2(\mathbb{R})$ en tal espacio. Denotemos por $d : \mathrm{SL}_2(\mathbb{R}) \times L^2(\Gamma\backslash \mathrm{SL}_2(\mathbb{R})) \to L^2(\Gamma\backslash \mathrm{SL}_2(\mathbb{R}))$ a dicha representación.

**Definición.** Un *grupo de Lie* es un grupo con una estructura de variedad diferenciable, para el cual la operación de grupo y la función inversa son $C^\infty$.

El grupo $\mathrm{SL}_2(\mathbb{R})$ es un grupo de Lie real, y lo que obtuvimos fue una representación de él. Recordar que en general, si $G$ es un grupo topológico localmente compacto, una representación de $G$ es un par $(\pi, H)$, donde $H$ es un espacio de Hilbert, y $\pi : G \times H \to H$ es tal que $(g, f) \to \pi(g)f$ es continua. Si $(\pi, H)$ es una representación de un grupo de Lie $G$, resulta conveniente estudiar las representaciones del álgebra de Lie de $G$.

**Definición.** Un *álgebra de Lie* es un espacio vectorial $\mathfrak{g}$ junto con una operación binaria $[\cdot,\cdot] : \mathfrak{g}\times\mathfrak{g} \to \mathfrak{g}$ llamada corchete de Lie, que satisface las siguientes propiedades:

- (Bilineal) $[aX + bY, Z] = a[X, Z] + b[Y, Z]$ y $[X, aY + bZ] = a[X, Y] + b[X, Z]$.

[1] Recordar que una representación irreducible es aquella que no tiene subespacios invariantes propios.

- (Alternada) $[X, X] = 0$.
- (Jacobi) $[X, [Y, Z]] + [Z, [X, Y]] + [Y, [Z, X]] = 0$.

Todo grupo de Lie tiene asociada un álgebra de Lie, que está dada como espacio vectorial por el espacio tangente en la identidad, y el corchete está dado por la función exponencial. Recordar que la función exponencial es una función

$$\exp : \mathfrak{g} \to G, \tag{4}$$

definida como $\exp(X) = \gamma(1)$, donde $\gamma : (-\epsilon, \epsilon) \to G$ es tal que $\gamma(0) = 1$, $\gamma'(0) = X$ y $F_X(\gamma(t)) = \gamma'(t)$, donde $F_X$ es un campo invariante a derecha asociado a $X$ (lo que implica que $\gamma(t)$ es un morfismo de grupos). La exponencial es un difeomorfismo local (ya que su diferencial es la identidad). Luego el corchete satisface

$$\exp(X)\exp(Y)\exp(X)^{-1}\exp(Y)^{-1} = \exp(\frac{1}{2}[X, Y] + \cdots), \tag{5}$$

donde los puntos suspensivos quieren decir que son términos de orden al menos 3.

**Ejercicio 2.** Probar que si $G = \mathrm{SL}_2(\mathbb{R})$, entonces su álgebra de Lie esta dada por las matrices de traza cero, la exponencial es la exponencial usual de matrices, y el corchete está dado por $[X, Y] = XY - YX$.

Si $(\pi, H)$ es una representación de un grupo localmente compacto $G$, decimos que $f \in H$ es $C^1$ si vale que para todo $X \in \mathfrak{g}$ existe el límite

$$\pi(X)f = \frac{d}{dt}\Big|_{t=0} \pi(\exp(tX))f = \lim_{t\to 0} \frac{1}{t}(\pi(\exp(tX))f - f). \tag{6}$$

Recursivamente, decimos que $f$ es $C^k$ si es $C^1$ y $\pi(X)f$ es $C^{k-1}$ para todo $X \in \mathfrak{g}$. Definimos el espacio de vectores suaves $H^\infty$ como el espacio de vectores $f$ que están en $C^k$ para todo $k \in \mathbb{N}$. Por definición, el álgebra de Lie $\mathfrak{g}$ de $G$ actúa en $H^\infty$. No es difícil ver que $\pi([X, Y]) = \pi(X)\pi(Y) - \pi(Y)\pi(X)$, con lo cual $(\pi, H)$ es una representación de $\mathfrak{g}$.

**Proposición 3.** El subespacio $H^\infty$ satisface:

- $H^\infty$ es $G$ invariante.
- $H^\infty$ es denso en $H$.

*Prueba.* Ver Bump (1997, Proposición 2.4.2). □

Supongamos que $(\pi, H)$ es una representación unitaria[2] e irreducible de $G$, y sea $K$ un subgrupo compacto de $G$ (nuestro interés es en $G = \mathrm{SL}_2(\mathbb{R})$ y $K = \mathrm{SO}(2)$). Luego podemos estudiar la restricción a $K$ de nuestra representación. Si $\tau$ es una representación irreducible de $K$, denotamos

[2] Recordar que una representación unitaria en un espacio de Hilbert es aquella que preserva el producto interno, o sea $\langle \pi(v), \pi(w) \rangle = \langle v, w \rangle$

por $H_\tau$ el subespacio de $H$ donde $K$ actúa por $\tau$. Por el teorema de Peter–Weyl (ver por ejemplo Bump (ibíd., Teorema 2.4.1)),

$$H = \bigoplus_{\tau}^{\perp} H_\tau. \tag{7}$$

Una representación se dice *admisible* si los subespacios $H_\tau$ son de dimensión finita. Vale que toda representación unitaria irreducible de $\mathrm{SL}_2(\mathbb{R})$ es admisible, con lo cual en la descomposición anterior los espacios $H_\tau$ tiene dimensión finita. Notar que en dicha descomposición, los vectores cuya proyección a los subespacios $H_\tau$ son cero para casi toda representación (salvo finitas), son densos y tienen la particularidad de que si miramos el subespacio

$$\pi(K)v = \langle \pi(k)v \ : \ k \in K \rangle,$$

es de dimensión finita. Los vectores $K$ *finitos* son los que satisfacen esta propiedad, y denotamos $H_{fin}$ es subespacio de vectores finitos. El compacto $K$ actúa en $H_{fin}$. Sea $H_0 = H_{fin} \cap H^\infty$ (claramente $H_0$ es $K$-invariante).

**Proposición 4.** El subespacio $H_0$ es denso en $H$, y tiene una estructura de $(\mathfrak{g}, K)$-módulo, esto es un $\mathfrak{g}$-módulo y un $K$-módulo donde ambas representaciones satisfacen la compatibilidad

$$\pi(g)\pi(X)\pi(g^{-1})f = \pi(\mathrm{Ad}(g)X)f, \qquad \forall\ X \in \mathfrak{g}, g \in K. \tag{8}$$

*Prueba.* Mirando la descomposición (7), es claro que $H_\tau^\infty = H^\infty \cap H_\tau$ es denso en $H_\tau$, con lo cual la primera afirmación se deduce del hecho de que sumas finitas de vectores suaves es suave. La compatibilidad es inmediata. □

Recordar que al álgebra de Lie $\mathfrak{g}$ podemos asociarle su álgebra envolvente, que se define como

$$U(\mathfrak{g}) = \bigoplus_{k=0}^{\infty} \otimes_{\mathbb{R}}^{k} \mathfrak{g}/I, \qquad I = \langle X \otimes Y - Y \otimes X - [X, Y] \ : \ X, Y \in \mathfrak{g} \rangle. \tag{9}$$

Donde $\otimes_{\mathbb{R}}^{k}\mathfrak{g}$ es el producto tensorial usual, y el producto esta dado por el producto tensorial $\otimes^i \mathfrak{g} \times \otimes^j \mathfrak{g} \to \otimes^{i+j}\mathfrak{g}$. Es casi inmediato que cualquier representación de $\mathfrak{g}$ se extiende a una representación de $U(\mathfrak{g})$ (por preservar el corchete, ver Bump (ibíd., Proposición 2.3)).

**Proposición 5.** Sea $V$ un $(\mathfrak{g}, K)$ modulo admisible irreducible, y sea $\mathfrak{z}$ un elemento del centro de $U(\mathfrak{g})$. Entonces $\mathfrak{z}$ actúa como un escalar en $V$.

*Prueba.* Como $X$ es un elemento del centro, $\mathrm{Ad}(g)\mathfrak{z} = \mathfrak{z}$, con lo cual la acción de $K$ conmuta con la de $\mathfrak{z}$. Luego $X$ deja fijo los espacios $V_\tau$ de (7). Como estos son de dimensión finita e irreducibles, por el Lema de Schur, $\mathfrak{z}$ actúa como un escalar en cada uno de ellos. Por otro lado, si $\lambda$ es un autovalor en algún espacio, y miramos $V_\lambda = \{v \in V \ : \ \pi(\mathfrak{z})v = \lambda v\}$, esto es un subespacio no vacío, e invariante por la acción de $\mathfrak{g}$ (porque $\mathfrak{z}$ conmuta con $\mathfrak{g}$) y de $K$. Como $V$ es irreducible, $V_\lambda = V$. □

En nuestro caso, $K = \mathrm{SO}(2)$, que es un grupo abeliano compacto, con lo cual sus representaciones están dadas por caracteres. Las mismas están parametrizadas por los enteros, y están dadas por $\chi_k(\theta) = e^{ik\theta}$. Luego cualquier representación admisible de $\mathrm{SL}_2(R)$, se descompone como

$$V = \bigoplus_{k\in\mathbb{Z}} V(k).$$

El conjunto $\sigma = \{k \in \mathbb{Z} \ : \ V(k) \neq \varnothing\}$ se llama el conjunto de *$K$-tipos*.

Tomemos los siguientes generadores de $U(\mathfrak{g})$:

$$\hat{R} = \left(\begin{smallmatrix} 0 & 1 \\ 0 & 0 \end{smallmatrix}\right), \qquad \hat{L} = \left(\begin{smallmatrix} 0 & 0 \\ 1 & 0 \end{smallmatrix}\right) \quad \text{y} \quad \hat{H} = \left(\begin{smallmatrix} 1 & 0 \\ 0 & -1 \end{smallmatrix}\right). \tag{10}$$

Es fácil ver que satisfacen las relaciones: $[\hat{H}, \hat{R}] = 2\hat{R}$, $[\hat{H}, \hat{L}] = -2\hat{L}$ y $[\hat{R}, \hat{L}] = \hat{H}$. Si definimos

$$-4\Delta = \hat{H}^2 + 2\hat{R}\hat{L} + 2\hat{L}\hat{R},$$

(donde multiplicación es en $U(\mathfrak{g})$), entonces vale que $\Delta$ es un elemento del centro de $U(\mathfrak{g})$ (ver Bump (1997, Teorema 2.2.1)) que se llama *Casimir*. Para relacionar el operador de Casimir, con el introducido en la sección anterior, es preciso extender la acción de $\mathfrak{g}$ a su complexificación (para poder hacer un cambio de coordenadas acorde). Si $V$ es un $(\mathfrak{g}, K)$-módulo, podemos mirar la complexificación de $V$, dada por $V_{\mathbb{C}} = V \otimes_{\mathbb{R}} \mathbb{C}$. A la vez, podemos mirar la complexificación de $\mathfrak{g}$, donde el espacio vectorial es el complexificado, y el corchete lo extendemos de manera que sea $\mathbb{C}$-bilineal. Definamos los elementos

$$R = \frac{1}{2}\begin{pmatrix} 1 & i \\ i & -1 \end{pmatrix}, \qquad L = \frac{1}{2}\begin{pmatrix} 1 & -i \\ -i & -1 \end{pmatrix} \quad \text{y} \quad H = \begin{pmatrix} 0 & -i \\ i & 0 \end{pmatrix}. \tag{11}$$

Estos son conjugados a los anteriores por la matriz $-\frac{1+i}{2}\left(\begin{smallmatrix} i & 1 \\ i & -1 \end{smallmatrix}\right)$, con lo cual satisfacen las mismas relaciones (y dan el mismo $\Delta$ por estar en el centro).

Si $V$ es un $(\mathfrak{g}, K)$-modulo irreducible, sabemos que

$$V = \sum_k V(k),$$

donde cada $V(k)$ es de dimensión finita. Denotemos $W = \left(\begin{smallmatrix} 0 & -1 \\ 1 & 0 \end{smallmatrix}\right)$.

**Proposición 6.** En la notación anterior valen.

1. $W\phi = ik\phi$, con lo cual $H\phi = k\phi$.
2. $V(k) = \{v \in V \ : \ Hv = kv\}$.
3. $RV(k) \subset V(k+2)$.
4. $LV(k) \subset V(k-2)$.
5. Cada $V(k)$ es de dimensión a lo sumo 1, y si $V(j)$, $V(k)$ son no nulos, entonces $j \equiv k \pmod 2$.

*Prueba.* (1) Por definición, si $\phi \in V(k)$, $W\phi(g) = \frac{d}{dt}\Big|_{t=0} \phi(g\exp(tW))$. Como $\exp(t\left(\begin{smallmatrix} 0 & -1 \\ 1 & 0 \end{smallmatrix}\right)) = \begin{pmatrix} \cos(t) & -\sin(t) \\ \sin(t) & \cos(t) \end{pmatrix}$,

$$dW\phi(g) = \frac{d}{dt}\Big|_{t=0} \exp(ikt)\phi(g) = ik\phi(g).$$

(2) Es claro de (1).
(3) Recordar que $[H, R] = 2R$, con lo cual

$$HR\phi = [H, R]\phi + RH\phi = 2R\phi + Rk\phi = (k+2)R\phi.$$

Luego de (2) se obtiene el enunciado. (4) se obtiene de una cuenta similar.
(5) Sea $V(k)$ no nulo, y tomemos $\phi$ un elemento no nulo en él. Luego podemos mirar

$$W = \mathbb{C}\phi \oplus \bigoplus_{n>0} \mathbb{C}R^n\phi \oplus \bigoplus_{n>0} \mathbb{C}L^n\phi. \tag{12}$$

Si vemos que esto es un $(\mathfrak{g}, K)$-submódulo no nulo, por ser $V$ irreducible deben coincidir, lo que demuestra la primera parte. Basta calcular $RL\phi$ y $LR\phi$. Sabemos que $\Delta$ actúa por un escalar en $V$, denotémoslo por $\lambda$. Utilizando que $-4\Delta = H^2 + 2RL + 2LR = H^2 - 2H + 4RL$, y $\phi \in V(k)$, tenemos que

$$-4\lambda\phi = \Delta\phi = (k^2 + 2k)\phi + 4LR\phi = (k^2 - 2k)\phi + 4RL\phi.$$

Luego tanto $LR\phi$ como $RL\phi$ están en $\mathbb{C}\phi$. La segunda parte es automática ahora, dado que sabemos que los espacios indexados por pares son invariantes, y los impares también, luego al ser el espacio irreducible, uno sólo de ellos puede ser no nulo. □

Recordar que denotamos por $d : \mathfrak{g} \times C^\infty(\mathrm{SL}_2(\mathbb{R})) \to C^\infty(\mathrm{SL}_2(\mathbb{R}))$ la acción regular a derecha. Dada $f \in S_k(\Gamma)$, su representación asociada está dada por el módulo que $\phi_f$ genera vía dicha acción, o sea $V_f = \mathrm{SL}_2(\mathbb{R}) \cdot \phi_f$. Mirándolo como $(\mathfrak{g}, K)$-módulo, sabemos que al restringir la representación a $K$ se descompone como

$$V_f = \sum_{j \in \mathbb{Z}} V_f(j).$$

**Proposición 7.** La representación $V_f$ es irreducible. Mas aún, vale:

1. $\Delta\phi_f = \frac{k}{2}(1 - \frac{k}{2})\phi_f$, y actúa de igual manera en $V_f$.

2. $V_f(j) = 0$ si $j < k$.

3. $V_f(k) = \mathbb{C}\phi_f$.

*Prueba.* Como vimos anteriormente, la acción de $dH$ es diagonal en cada sumando, con lo cual basta calcular la acción de $dL$ y $dR$. Como $L = \frac{1}{2}(\hat{H} - 2i\hat{R} - H)$, vamos a calcular la acción de

$d\hat{H}$ y $d\hat{R}$ en $\phi_f$.

$$d\hat{R}\phi_f \begin{pmatrix} y^{1/2} & xy^{-1/2} \\ 0 & y^{-1/2} \end{pmatrix} = \frac{d}{dt}\Big|_{t=0} \phi_f \left( \begin{pmatrix} y^{1/2} & xy^{-1/2} \\ 0 & y^{-1/2} \end{pmatrix} \left(\begin{smallmatrix} 1 & t \\ 0 & 1 \end{smallmatrix}\right) \right) =$$
$$\frac{d}{dt}\Big|_{t=0} \phi_f \begin{pmatrix} y^{1/2} & xy^{-1/2}+y^{1/2}t \\ 0 & y^{-1/2} \end{pmatrix} = \frac{d}{dt}\Big|_{t=0} y^{k/2} f(iy+x+ty) = y^{k/2+1}\frac{\partial f}{\partial x}(x+iy).$$

$$d\hat{H}\phi_f \begin{pmatrix} y^{1/2} & xy^{-1/2} \\ 0 & y^{-1/2} \end{pmatrix} = \frac{d}{dt}\Big|_{t=0} \phi_f \left( \begin{pmatrix} y^{1/2} & xy^{-1/2} \\ 0 & y^{-1/2} \end{pmatrix} \begin{pmatrix} e^t & 0 \\ 0 & e^{-t} \end{pmatrix} \right) =$$
$$\frac{d}{dt}\Big|_{t=0} e^{tk} y^{k/2} f(iye^{2t}+x) = y^{k/2}\left( kf(x+iy) - 2y\frac{\partial f}{\partial y}(x+iy) \right).$$

Juntando todo, tenemos que

$$dL(\phi_f) = \frac{1}{2}y^{k/2+1}(-2i)\left(\frac{\partial f}{\partial x}(z) - i\frac{\partial f}{\partial y}(z)\right) = 0, \tag{13}$$

por ser la función $f$ holomorfa.

Calculemos la acción de $-4\Delta$ en $V_f$. Como $V_f$ es el generado por $\phi_f$ y el operador de Casimir está en el centro, su acción es igual en todo $V_f$, con lo cual basta calcularla en $\phi_f$. Recordar que

$$-4\Delta = H^2 + 2RL + 2LR. \tag{14}$$

Sabemos que $L\phi_f = 0$, con lo cual

$$-4\Delta\phi_f = (H^2 + 2LR)\phi_f = (H^2 - 2[R,L])\phi_f = (H^2 - 2H)\phi_f. \tag{15}$$

Sabemos que $H\phi_f = -k\phi_f$, con lo cual $-2\Delta\phi_f = (k^2-2k)\phi_f$. Dividiendo por $-4$ obtenemos (1).

Para ver (2), sabemos que $L\phi_f = 0$, resta ver que pasa al aplicar $LR\phi_f$. Como tanto el operador de Casimir como $H$ actúan de manera diagonal en $V_f(k)$, de (15) se ve que $LR\phi_f$ es un múltiplo de $\phi_f$. En general, si calculamos $LR^n\phi_f$ (con $n \geqslant 2$), utilizando (14), vemos que $LR^n\phi_f$ es un múltiplo de $R^{n-1}\phi_f$ mas $RL(R^{n-1}\phi_f)$. La hipótesis inductiva dice que $LR^{n-1}\phi_f$ es un múltiplo de $R^{n-2}\phi_f$, con lo cual $LR^n\phi_f$ es un múltiplo de $R^{n-1}\phi_f$. Esto termina la demostración. □

*Observación.* El punto crucial de esta construcción, es que en la representación $V_f$ asociada a $f$, el elemento $\phi_f$ es un vector de peso mínimo, en el sentido que genera el primer subespacio no nulo.

Terminemos ahora la demostración de la Proposición 4: al tomar una función $\phi_f$ como en el enunciado, la representación que obtenemos tiene autovalor para el Casimir $\frac{k}{2}(1-\frac{k}{2})$. Haciendo una cuenta parecida a las hechas anteriormente, es fácil verificar que la representación como $(\mathfrak{g}, K)$-módulo que induce es exactamente como la descripta en la Proposición 7. Luego $L\phi_f = 0$, lo que por (13) implica que la función es holomorfa.

*Observación.* Lo fundamental del operador de Casimir, es que $\mathbb{C}[\Delta]$ es el centro de $U(\mathfrak{g})$, con lo cual es (salvo escalares) el único operador que conmuta con la acción de $K$ (y por lo tanto se diagonaliza).

**Definición.** La representación de $\mathrm{GL}_2(\mathbb{R})$ descripta anteriormente se llama *serie discreta* de peso $k$.

## 2 Formas automorfas como formas en los adèles

Al pensar a las formas modulares como formas de $\mathrm{SL}_2(\mathbb{R})$, nos surge el problema natural de tratar de definir a los operadores de Hecke en ellas. Dichos operadores no tienen una definición elemental en este mundo automorfo. Es por esto que necesitamos extender el dominio de definición para incluir los adèles y así sí poder definirlos.

**2.1 Adèles e idèles.** El anillo de adèles fue introducido por Claude Chevalley (1909-1984). En estas notas (y durante el curso), utilizaremos esta noción solamente para el caso de cuerpos de números. Son de gran utilidad en la teoría algebraica de números. Algunas referencias son Neukirch (1992), Stein (2012), Weil (1995) y Milne (2014).

Sea $\mathbb{Q}$ el cuerpo de números racionales. Una valuación en $\mathbb{Q}$ es una función $v : \mathbb{Q}\backslash\{0\} \to \mathbb{Z}$ que satisface:

- $v(nm) = v(n) + v(m)$.
- $v(n+m) \geqslant \min\{v(n), v(m)\}$.

A veces resulta útil extender la valuación a todo el cuerpo, decretando que $v(0) = \infty$. Los ejemplos que vamos a considerar de valuaciones son los siguientes: dado $p$ un primo (que podemos suponer positivo), todo número racional no nulo lo podemos escribir de manera única como $p^r \frac{a}{b}$, donde $r \in \mathbb{Z}$ y $p \nmid ab$. Esto nos permite definir una función

$$v_p : \mathbb{Q}\backslash\{0\} \to \mathbb{Z}, \qquad v_p\left(p^r \frac{a}{b}\right) = r. \tag{16}$$

Es fácil ver que esta función es una valuación, y se la denomina la *valuación p-ádica*.

**Ejercicio 3.** Probar que si $p$ es un número primo, entonces $v_p$ es efectivamente una valuación.

**Ejercicio 4.** Probar que si miramos al anillo de polinomios $\mathbb{Q}[x]$, y definimos la valuación dada por el orden de anulación en cero de un polinomio no nulo, satisface las propiedades de una valuación.

A una valuación $v$ uno puede asociarle un valor absoluto sobre el cuerpo. Si elegimos $\alpha \in \mathbb{R}$ de manera que $0 < \alpha < 1$, entonces definimos el valor absoluto *noarquimediano*

$$|\ .\ |_v : \mathbb{Q} \to \mathbb{R}_{\geqslant 0}, \text{ por}$$

$$|m|_v = \begin{cases} \alpha^{v(m)} & \text{si } m \neq 0. \\ 0 & \text{si } m = 0. \end{cases}$$

Distintas elecciones de $\alpha$ dan valores absolutos equivalentes, en el sentido que dan los mismos los abiertos. El valor absoluto $p$-ádico es el asociado a la valuación $v_p$, tomando la normalización usual $\alpha = \frac{1}{p}$. El mismo se denota por $|m|_p$. Así, $|3|_3 = \frac{1}{3}$ y $|\frac{4}{9}|_3 = 9$.

Claramente el valor absoluto $| \ . \ |_v$ asociado a una valuación satisface:

1. $|m|_v = 0$ si y sólo si $m = 0$.
2. $|nm|_v = |n|_v |m|_v$.
3. $|n + m|_v \leqslant \max\{|n|_v, |m|_v\}$.

La última desigualdad implica la desigualdad triangular clásica, pero es mucho mas fuerte que esta (y se suele llamar *ultramétrica*).

**Ejercicio 5.** Probar que con el valor absoluto $p$-ádico, el conjunto de números enteros está acotado. Mas aún, $|n|_p \leqslant 1$ para todo $n \in \mathbb{Z}$ (o sea todos los números enteros están a distancia a lo sumo 1 del cero).

No es difícil ver que un valor absoluto sobre un cuerpo satisface la propiedad ultramétrica si y sólo si los enteros están acotados (ver Milne (2014, Proposition 7.2)).

**Teorema 8** (Ostrowski)**.** Todo valor absoluto no trivial de $\mathbb{Q}$ es equivalente a uno de los siguientes:

- El valor absoluto clásico $|x|_\infty$.
- El valor absoluto $p$-ádico $|x|_p$ para algún primo $p$.

Recordar que la condición de que dos valores absolutos sean equivalentes es que den los mismos abiertos en el conjunto. Para la demostración, ver por ejemplo Milne (ibíd., Theorem 7.12).

**Ejercicio 6.** Probar que los valores absolutos que aparecen en el Teorema de Ostrowski son todos no equivalentes.

Para tener una manera unificada de hablar de los distintos valores absolutos, se suele introducir la noción de lugares, donde los lugares están compuestos por los valores absolutos noarquimedianos (correspondientes a los números primos), y uno extra asociado al valor absoluto usual, denominado lugar del infinito.

Dado un valor absoluto en un cuerpo, podemos completarlo con respecto a dicho valor absoluto.

**Ejercicio 7.** Probar que al completar un cuerpo con respecto a un valor absoluto obtenemos nuevamente un cuerpo, o sea las operaciones se pueden extender a la completación, y siguen satisfaciendo los axiomas de un cuerpo. Probar a la vez que el valor absoluto se extiende a la completación.

Si completamos a los números racionales con el valor absoluto clásico, obtenemos el cuerpo $\mathbb{R}$ de números reales. No obstante, si completamos a $\mathbb{Q}$ con el valor absoluto $p$-ádico, obtenemos otro cuerpo, que notamos por $\mathbb{Q}_p$. Por ejemplo, la sucesión $\{p^n\}_{n \in \mathbb{N}}$ tiende a cero en $\mathbb{Q}_p$. Esto permite caracterizar al conjunto de números $p$-ádicos como sumas

$$\mathbb{Q}_p = \left\{ \sum_{n \geqslant N_0} a_n p^n \ : \ a_n \in \{0, \dots, p-1\} \right\}.$$

**Ejercicio 8.** Probar esta caracterización de la siguiente manera:

- Probar que todo entero tiene una tal descripción (como suma finita), que corresponde a mirar al número en base $p$.
- Probar que si $\{a_n\}$ es una sucesión de Cauchy con el valor absoluto $p$-ádico, entonces los primeros términos de tal expresión se vuelven constantes a partir de un cierto $n$.
- Probar que si $n \in \mathbb{Z}$ es coprimo con $p$, entonces $\frac{1}{n}$ tiene una tal expresión. (Sugerencia: invertir la expresión de $n$ formalmente). Deducir que todo racional se escribe de esta manera.
- Terminar el argumento.

De manera análoga, podemos definir $\mathbb{Z}_p$ como el completado de $\mathbb{Z}$ con respecto al valor absoluto $p$-ádico. Notar que todos los elementos así construidos dentro de $\mathbb{Q}_p$ tienen valor absoluto a lo sumo 1. Es fácil ver que

$$\mathbb{Z}_p = \{x \in \mathbb{Q}_p \ : \ |x|_p \leqslant 1\},$$

o sea corresponde a la bola unidad. Notar que $\mathbb{Z}_p$ es cerrado (por definición), y a la vez es abierto, porque el valor absoluto $p$-ádico toma un conjunto discreto de valores (así no alcanza ningún valor entre 1 y $p$, luego podemos definir $\mathbb{Z}_p$ como la bola abierta de radio $3/2$). Esto hace con que $\mathbb{Q}_p$ sea totalmente disconexo como espacio topológico. Notar que $\mathbb{Z}_p$ es compacto (por ser cerrado y acotado). Este hecho será de vital importancia.

**Ejercicio 9.** El propósito de este ejercicio es probar que $\sqrt{-1}$ es un elemento en $\mathbb{Z}_5$. Como vimos antes, todo elemento de $\mathbb{Z}_5$ se puede escribir como

$$a_0 + a_1 \cdot 5 + a_2 \cdot 5^2 + a_3 \cdot 5^3 + \cdots$$

Luego debemos determinar los coeficientes $a_n$.

- Probar que existe un $a_0$ que hace que $a_0^2 \equiv -1 \pmod 5$.
- Construido $a_0$, probar que existe $a_1$ tal que $(a_0 + 5a_1)^2 \equiv -1 \pmod{25}$.
- Demostrar el enunciado con algún argumento inductivo.

Luego no sólo $\mathbb{Q}_5$ es un cuerpo totalmente distinto a $\mathbb{R}$, sino que hay polinomios que tienen raíces en él, pero no son reales.

**Teorema 9** (Fórmula del producto)**.** Si $n \in \mathbb{Q}$, entonces

$$|n| \prod_{p \text{ primos}} |n|_p = 1.$$

La demostración es totalmente elemental, y la dejamos como ejercicio para el lector.

Sea $S$ un conjunto finito de lugares que contenga al lugar del infinito. Definimos al conjunto de $S$-adèles, y lo denotamos $\mathbb{A}_{\mathbb{Q}}^S$, al conjunto

$$\mathbb{A}_{\mathbb{Q}}^S = \prod_{p \notin S} \mathbb{Z}_p \times \prod_{p \in S} \mathbb{Q}_p \times \mathbb{R}. \tag{17}$$

Este anillo es un anillo topológico con la topología producto. El anillo de *adèles* es la unión sobre los conjuntos $S$ que son finitos y contienen el lugar del infinito, de los conjuntos $\mathbb{A}_{\mathbb{Q}}^S$. Equivalentemente, podemos definir los adèles como el producto directo restringido con respecto al maximal $\mathbb{Z}_p$ para los primos finitos, donde la topología es la producto. Así, $\mathbb{A}_{\mathbb{Q}}$ resulta un anillo topológico localmente compacto.

Notar que claramente $\mathbb{Q}$ es un subanillo de $\mathbb{A}_{\mathbb{Q}}$, donde podemos pensar los números racionales como vectores infinitos del mismo elemento, y es claro que cualquier número racional $r$ está en $\mathbb{Z}_p$ para todos los primos $p$ salvo finitos (que dividen al denominador), con lo cual es un adèle.

**2.2 Formas modulares como funciones de los adèles.** Queremos extender las definiciones de los adèles a otros grupos, como el grupo de matrices de $2\times 2$. Dado un primo $p$, definimos $\mathrm{GL}_2(\mathbb{Q}_p)$ al grupo de matrices inversibles de $2 \times 2$ con coeficientes en el cuerpo de números $p$-ádicos $\mathbb{Q}_p$.

**Ejercicio 10.** Vimos anteriormente que el anillo de enteros $p$-ádicos $\mathbb{Z}_p$ es compacto. Probar que el grupo $\mathrm{GL}_2(\mathbb{Z}_p)$ es compacto. (Sugerencia: ver que es cerrado dentro del compacto $\mathbb{Z}_p^4$).

Análogamente, al mirar el valor absoluto usual, y su completación, denotamos por $\mathrm{GL}_2(\mathbb{R})$ al grupo de matrices de $2 \times 2$ con coeficientes en el cuerpo de números reales $\mathbb{R}$. Con todos estos grupos, queremos formar al grupo $\mathrm{GL}_2(\mathbb{A}_{\mathbb{Q}})$, el grupo de matrices de $2 \times 2$ con coeficientes en el anillo de adèles.

Dicho grupo se puede definir como el producto directo restringido de los grupos $\mathrm{GL}_2(\mathbb{Q}_p)$ con respecto a los compactos $\mathrm{GL}_2(\mathbb{Z}_p)$. Esto es, los elementos de $\mathrm{GL}_2(\mathbb{A}_{\mathbb{Q}})$ son tiras de elementos $(g_2, g_3, \ldots, g_\infty)$ donde $g_p \in \mathrm{GL}_2(\mathbb{Q}_p)$ y $g_\infty \in \mathrm{GL}_2(\mathbb{R})$ que satisfacen que para todos los primos $p$ salvo finitos, $g_p \in \mathrm{GL}_2(\mathbb{Z}_p)$. Sea $\{K_p\}_p$ un conjunto de subgrupos compactos donde cada $K_p \in \mathrm{GL}_2(\mathbb{Q}_p)$ (por ejemplo $K_p = \mathrm{GL}_2(\mathbb{Z}_p)$) de forma tal que:

- $K_p = \mathrm{GL}_2(\mathbb{Z}_p)$ para casi todo primo $p$ (o sea todos salvo finitos).
- El determinante $\det : K_p \to \mathbb{Z}_p^\times$ es suryectivo.

Nos interesa particularmente al tratar de estudiar las formas modulares en $S_2(\Gamma_0(N))$, los compactos $K_p(N)$ dados por

$$K_p(N) = \left\{ \begin{pmatrix} a & b \\ c & d \end{pmatrix} \in \mathrm{SL}_2(\mathbb{Z}_p) \ : \ c \equiv 0 \pmod{N} \right\}. \tag{18}$$

Entonces, vale el siguiente resultado.

**Teorema 10** (Aproximación fuerte)**.** Con las notaciones anteriores, la siguiente igualdad se verifica:

$$\mathrm{GL}_2(\mathbb{A}_{\mathbb{Q}}) = \mathrm{GL}_2(\mathbb{Q})\,\mathrm{GL}_2(\mathbb{R})^+ \prod_p K_p(N).$$

*Prueba.* Ver por ejemplo Hida (2000, Theorem 3.2). En el caso de matrices de $1 \times 1$ (o sea $\mathrm{GL}_2$), este resultado se deduce del Teorema Chino del Resto. □

Luego para definir una forma automorfa en $\mathrm{GL}_2(\mathbb{A}_{\mathbb{Q}})$, tenemos que entender como se comporta en $\mathrm{GL}_2(\mathbb{R})^+$, como actúa $\mathrm{GL}_2(\mathbb{Q})$, y como actúan los subgrupos compactos $K_p$. Si $g \in \mathrm{GL}_2(\mathbb{A}_{\mathbb{Q}})$, lo descomponemos como $g = \gamma g_\infty k_0$, donde $\gamma \in \mathrm{GL}_2(\mathbb{Q})$, $g_\infty \in \mathrm{GL}_2(\mathbb{R})^+$ y $k_0 \in \prod_p K_p$. Dada $f \in S_k(\Gamma_0(N))$, definimos $\phi_f : \mathrm{GL}_2(\mathbb{A}_{\mathbb{Q}}) \to \mathbb{C}$ por

$$\phi_f(g) = f(g_\infty \cdot i) j(g_\infty, i)^{-k}, \text{ donde } g = \gamma g_\infty k_0.$$

El Teorema 10 nos asegura que esto da una función bien definida, dado que es fácil ver que $\mathrm{GL}_2(\mathbb{Q}) \cap \mathrm{GL}_2(\mathbb{R})^+ \prod_p K_p = \Gamma_0(N)$. Juntándolo con la Proposición 2, tenemos:

**Proposición 11.** La función $f(z) \to \phi_f(g)$ da un isomorfismo entre $S_k(\Gamma_0(N))$ y el espacio de funciones $\phi$ en $\mathrm{GL}_2(\mathbb{A}_{\mathbb{Q}})$ que satisfacen las siguientes condiciones:

1. $\phi(\gamma g) = \phi(g)$ para todo $\gamma \in \mathrm{GL}_2(\mathbb{Q})$.

2. $\phi(g k_0) = \phi(g)$ para todo $k_0 \in \prod_p K_p$.

3. $\phi(g r(\theta)) = \exp(ik\theta)\phi(g)$ para todo $\theta \in [0, 2\pi]$.

4. La función $\phi$, vista como función de $\mathrm{GL}_2(\mathbb{R})^+$ es $C^\infty$ y satisface la ecuación diferencial
$$\Delta\phi = -\frac{k}{2}\left(\frac{k}{2} - 1\right)\phi.$$

5. $\phi(zg) = \phi(g)$ para todo $z \in Z_{\mathbb{A}} = \left\{\left(\begin{smallmatrix} t & 0 \\ 0 & t \end{smallmatrix}\right) \;:\; t \in \mathbb{A}_{\mathbb{Q}}^\times\right\}$.

6. (Crecimiento moderado) Para todo $c > 0$ y todo conjunto compacto $\Omega$ de $\mathrm{GL}_2(\mathbb{A}_{\mathbb{Q}})$, existen constantes $C, N$ tales que
$$\left|\phi\left(\left(\begin{smallmatrix} a & 0 \\ 0 & 1 \end{smallmatrix}\right) g\right)\right| \leqslant C|a|^N.$$

7. $\phi$ es cuspidal, o sea
$$\int_{\mathbb{Q}\backslash\mathbb{A}_{\mathbb{Q}}} \phi\left(\left(\begin{smallmatrix} 1 & x \\ 0 & 1 \end{smallmatrix}\right) g\right) dx = 0 \text{ para casi todo } g.$$

Queremos hacer dos observaciones de esta Proposición.

*Observación.* Si en lugar de trabajar con $S_k(\Gamma_0(N))$, permitimos que las formas tengan Nebentypus, o sea miramos formas en $S_k(\Gamma_0(N),\epsilon)$, entonces usando la relación entre caracteres de $\mathbb{Z}/N\mathbb{Z}$ y los caracteres de Hecke, podemos asociarle a $\epsilon$ un carácter $\psi : \mathbb{A}_{\mathbb{Q}}^\times \to \mathbb{C}^\times$ trivial en $\mathbb{Q}^\times$ (esto es el caso unidimensional de la correspondencia anterior), y obtenemos una función $\phi_f$ que satisface las relaciones anteriores, pero agregando a la acción del compacto y del centro la acción de $\psi$, o sea reemplazamos 2. y 5. por

2′. $\phi(gk_0) = \phi(g)\psi(k_0)$.
5′. $\phi(zg) = \phi(g)\psi(z)$.

*Observación.* Sin entrar demasiado en detalles técnicos, el grupo $\mathrm{GL}_2(\mathbb{A}_{\mathbb{Q}})$ por ser localmente compacto posee una medida de Haar invariante (consultar por ejemplo el libro Bump (1997) para ver una descripción explícita; en realidad al ser $\mathrm{GL}_2$ un grupo unimodular dicha medida es invariante tanto a izquierda como a derecha), lo que permite integrar funciones. Además, la medida la suponemos normalizada de manera que el compacto maximal tenga medida 1. La definición de la medida es compatible con la identificación

$$Z_{\mathbb{A}_{\mathbb{Q}}} \mathrm{GL}_2(\mathbb{Q})\backslash \mathrm{GL}_2(\mathbb{A}_{\mathbb{Q}})/K_0K_\infty \simeq \Gamma_0(N)\backslash \mathrm{SL}_2(\mathbb{R})/\mathrm{SO}(2),$$

donde $K_0 = \prod_p K_p(N)$ como en (18). Luego si queremos integrar $\phi_f$, tenemos que

$$\int_{Z_{\mathbb{A}_{\mathbb{Q}}} \mathrm{GL}_2(\mathbb{Q})\backslash \mathrm{GL}_2(\mathbb{A}_{\mathbb{Q}})/K_0K_\infty} |\phi_f(g)|^2 dg = \int_{\Gamma_0(N)\backslash\mathfrak{h}} |f(z)|^2 y^k \frac{dx\,dy}{y^2} < \infty. \tag{19}$$

Decimos que una función $\phi$ en $\mathrm{GL}_2(\mathbb{A}_{\mathbb{Q}})$ es de cuadrado integrable si vale que

$$\int_{Z_{\mathbb{A}_{\mathbb{Q}}} \mathrm{GL}_2(\mathbb{Q})\backslash \mathrm{GL}_2(\mathbb{A}_{\mathbb{Q}})/K_0K_\infty} |\phi(g)|^2 dg < \infty$$

y denotamos a este espacio por $L^2(\mathrm{GL}_2(\mathbb{Q})\backslash \mathrm{GL}_2(\mathbb{A}_{\mathbb{Q}}))$. Lo que acabamos de ver es que $\phi_f$ es de cuadrado integrable. El espacio de funciones de cuadrado integrable y que además son holomorfas, en el sentido de la Proposición 11, 7. lo denotamos por $L_0^2(\mathrm{GL}_2(\mathbb{Q})\backslash \mathrm{GL}_2(\mathbb{A}_{\mathbb{Q}}))$.

**2.3 Operadores de Hecke.** Ahora que podemos pensar a las formas modulares como funciones en $\mathrm{GL}_2(\mathbb{A}_{\mathbb{Q}})$ de cuadrado integrable, podemos definir a los operadores de Hecke como ciertos operadores de convolución en este espacio. Esta noción se generaliza fácilmente a otros cuerpos de números.

Sea $p$ un primo tal que $p \nmid N$, con lo cual $K_p(N) = \mathrm{GL}_2(\mathbb{Z}_p)$. Definimos $H_p$ como el conjunto

$$H_p = \mathrm{GL}_2(\mathbb{Z}_p)\begin{pmatrix} p & 0 \\ 0 & 1 \end{pmatrix}\mathrm{GL}_2(\mathbb{Z}_p), \tag{20}$$

esto es el conjunto de todas las matrices de $2\times2$ con coeficientes en $\mathbb{Z}_p$ de determinante $p\mathbb{Z}_p^\times$. Notar que este es un conjunto que es invariante bajo el producto tanto a derecha como izquierda por el compacto (comparar con la sección 2.3. de las notas de M. Harris). Luego definimos el operador de Hecke $\widetilde{T}(p)$ como convolucionar a derecha con la característica de $H_p$ (en la componente $p$-ésima simplemente), esto es

$$\widetilde{T}(p)\phi(g) = \int_{H_p} \phi(gh_p)dh_p.$$

Aquí estamos haciendo un poco de abuso de notación, porque el elemento $h_p$ representa al elemento de $\mathrm{GL}_2(\mathbb{A}_\mathbb{Q})$ que tiene a la identidad en todas las coordenadas, salvo la $p$-ésima, donde tiene a $h_p$.

Es bastante claro ver que si $\phi \in L^2(\mathrm{GL}_2(\mathbb{Q})\backslash \mathrm{GL}_2(\mathbb{A}_\mathbb{Q}))$, entonces $\widetilde{T}(p)\phi$ también es de cuadrado integrable, y además satisface las mismas propiedades de la Proposición 11. Para verificar la segunda propiedad, notar que la invariancia en todas las coordenadas salvo la $p$-ésima es clara. En la coordenada $p$-ésima, notar que $\mathrm{GL}_2(\mathbb{Z}_p)H_p = H_p$, y la medida de Haar es invariante por traslación, así

$$\int_{H_p} \phi(gk_0h_p)dh_p = \int_{H_p} \phi(gh_p)dh_p.$$

Luego $\widetilde{T}(p)$ define un operador tanto en $L^2(\mathrm{GL}_2(\mathbb{Q})\backslash \mathrm{GL}_2(\mathbb{A}_\mathbb{Q}))$ como en $L^2_0(\mathrm{GL}_2(\mathbb{Q})\backslash \mathrm{GL}_2(\mathbb{A}_\mathbb{Q}))$.

**Lema 12.** Si $f \in S_k(\Gamma_0(N))$ entonces $p^{k/2-1}\widetilde{T}(p)\phi_f = \phi_{T(p)f}$.

*Prueba.* La idea de la demostración es mirar el método de eliminación Gaussiana para matrices de $2 \times 2$, que nos permite dar la descomposición disjunta

$$H_p = \bigcup_{i=0}^{p-1} \begin{pmatrix} p & i \\ 0 & 1 \end{pmatrix} \mathrm{GL}_2(\mathbb{Z}_p) \cup \begin{pmatrix} 1 & 0 \\ 0 & p \end{pmatrix} \mathrm{GL}_2(\mathbb{Z}_p).$$

**Ejercicio 11.** Si $s_i$ es cualquiera de los $p+1$ representantes de las descomposición anterior de $H_p$, y $M = \left(\begin{smallmatrix} a & b \\ c & d \end{smallmatrix}\right) \in \mathrm{GL}_2(\mathbb{Z}_p)$, probar que existe un único representante $s_j$ tal que $s_j^{-1}Ms_i \in \mathrm{GL}_2(\mathbb{Z}_p)$.

Como $\phi_f$ es invariante por la acción de $\mathrm{GL}_2(\mathbb{Z}_p)$, y como la medida de Haar está normalizada de manera que $\mu(\mathrm{GL}_2(\mathbb{Z}_p)) = 1$, nos queda que

$$p^{k/2-1}\widetilde{T}(p)\phi_f(g) = p^{k/2-1}\sum_{i=0}^{p-1} \phi_f\left(g\begin{pmatrix} p & i \\ 0 & 1 \end{pmatrix}\right) + p^{k/2-1}\phi_f\left(g\begin{pmatrix} 1 & 0 \\ 0 & p \end{pmatrix}\right) \tag{21}$$

Recordar el abuso de notación de la definición del operador de Hecke. Aquí la matriz $\left(\begin{smallmatrix} p & i \\ 0 & 1 \end{smallmatrix}\right)$ es el elemento en $\mathrm{GL}_2(\mathbb{A}_\mathbb{Q})$ que tiene a la matriz identidad en todas las coordenadas, salvo la $p$-ésima, donde tiene a dicha matriz.

Escribamos $g = \gamma g_\infty k_0$ (usando el Teorema 10). Denotemos por $k_p$ a la componente $p$-ésima del elemento $k_0 \in K$. Dado $s_p$ cualquiera de los $p+1$ elementos de la descomposición de $H_p$, denotemos por $\tilde{s}_p$ al representante (que es una matriz en $\mathrm{GL}_2(\mathbb{Q})$) que hace que $\tilde{s}_p^{-1}k_ps_p \in \mathrm{GL}_2(\mathbb{Z}_p)$ (que existe por el Ejercicio 11). Luego definimos el siguiente elemento de $\mathrm{GL}_2(\mathbb{A}_\mathbb{Q})$:

$$(t_p)_q = \begin{cases} 1 & \text{si } q \neq p \\ \tilde{s}_p & \text{si } q = p. \end{cases}$$

Notar que este elemento no está en el compacto (en $p$ no tiene determinante una unidad), pero el elemento $\tilde{s}_p/t_p$ (o sea poner en todas las coordenadas el mismo elemento, salvo en $p$, donde ponemos un 1) sí lo está. Entonces

$$g t_p = \gamma \tilde{s}_p (\tilde{s}_p^{-1} g_\infty)(\tilde{s}_p^{-1} k_0 t_p).$$

El último factor está en el compacto (claramente lo está en todas las coordenadas salvo la $p$-ésima y acabamos de ver anteriormente que esta también lo esta). Luego $\phi_f(g t_p) = f(\tilde{s}_p^{-1} g_\infty \cdot i) j(\tilde{s}_p^{-1} g_\infty, i)^{-k} = f\left(\frac{z+j}{p}\right) p^{-k/2}$. Reemplazando en (21), tenemos la bien conocida fórmula

$$p^{k/2-1} \tilde{T}(p) \phi_f(g) = \frac{1}{p} \sum_{i=0}^{p-1} f\left(\frac{z+i}{p}\right) + p^{k-1} f(pz).$$

□

Notar que es claro que los operadores de Hecke conmutan entre ellos (por actuar en distintas coordenadas), y no es difícil ver que son autoadjuntos para el producto interno de $L^2(\mathrm{GL}_2(\mathbb{Q}) \backslash \mathrm{GL}_2(\mathbb{A}_\mathbb{Q}))$, con lo cual nuestra correspondencia preserva autofunciones.

**2.4 Correspondencia entre formas modulares y representaciones adélicas.** Como vimos anteriormente, en $L^2(\mathrm{GL}_2(\mathbb{Q}) \backslash \mathrm{GL}_2(\mathbb{A}_\mathbb{Q}))$ tenemos una acción natural de $\mathrm{GL}_2(\mathbb{A}_\mathbb{Q})$ a derecha. Esta representación se parte como una parte continua (que proviene de las series de Eisenstein) y una parte discreta, que corresponde a las formas en $L_0^2(\mathrm{GL}_2(\mathbb{Q}) \backslash \mathrm{GL}_2(\mathbb{A}_\mathbb{Q}))$ (i.e. las cuspidales). Si restringimos la acción regular al espacio de formas cuspidales, se puede ver que dicho operador es un operador compacto, con lo cual el espacio se parte como suma (infinita) de representaciones irreducibles, cada una de ellas con multiplicidad finita. Mas aún, Jacquet–Langlands demostraron que la multiplicidad de las representaciones irreducibles es exactamente uno.

Toda representación irreducible en $L_0^2(\mathrm{GL}_2(\mathbb{Q}) \backslash \mathrm{GL}_2(\mathbb{A}_\mathbb{Q}))$ de $\mathrm{GL}_2(\mathbb{A}_\mathbb{Q})$ se parte como producto tensorial restringido de representaciones "locales", o sea de representaciones de $\mathrm{GL}_2(\mathbb{Q}_p)$ y de $\mathrm{GL}_2(\mathbb{R})$, luego es muy importante entender las representaciones de dichos grupos. Las mismas se conocen completamente, y por ejemplo se pueden ver en Gelbart (1975), Capítulo 4 (ver también las notas de M. Harris, Sección 2.4.1 en adelante). La relación entre las formas automorfas y las representaciones de grupo $\mathrm{GL}_2(\mathbb{A}_\mathbb{Q})$ está dada por la siguiente equivalencia:

- Si $f \in S_k(\Gamma_0(N))$ es autofunción para los operadores de Hecke, el subespacio generado por $\phi_f$ dentro de $L_0^2(\mathrm{GL}_2(\mathbb{Q}) \backslash \mathrm{GL}_2(\mathbb{A}_\mathbb{Q}))$ bajo la acción de $\mathrm{GL}_2(\mathbb{A}_\mathbb{Q})$ es una representación irreducible, o sea $H_f = \pi_{\mathrm{GL}_2(\mathbb{A}_\mathbb{Q})} \phi_f$ es un subespacio irreducible.

- (Casselman) Sea $\pi$ una representación unitaria irreducible de $\mathrm{GL}_2(\mathbb{A}_\mathbb{Q})$. Si $\pi$ es tal que aparece dentro de la representación regular a derecha de $L_0^2(\mathrm{GL}_2(\mathbb{Q}) \backslash \mathrm{GL}_2(\mathbb{A}_\mathbb{Q}))$ y es tal que la representación $\pi_\infty$ en la componente del infinito es una serie discreta de peso $k$, entonces existe un nivel $N$ (que se puede definir a partir de las representaciones locales $\pi_p$), y una forma modular $f \in S_k(\Gamma_0(N))$ autofunción para los operadores de Hecke, tal que $\pi \simeq H_f$.

Para ver los detalles de esta correspondencia, ver el capítulo 5 de Gelbart (ibíd.) y las referencias que allí se mencionan.

# 3 Formas de Hilbert

**3.1 Definiciones y propiedades básicas.** El estudio de las formas fue introducido por David Hilbert y por Otto Blumenthal a finales del siglo XIX. Referencias estándares para su estudio son los libros van der Geer (1988) y Freitag (1990), aunque las notas Bruinier (2008) son muy amenas también. La idea es obtener análogos de dimensión mayor de cocientes del plano complejo superior. Para facilitar la exposición, nos concentraremos en el caso de formas de Hilbert para cuerpos cuadráticos reales (que fueron las primeras estudiadas), que ya poseen suficiente grado de generalidad, y así evitar el uso de índices. Al final diremos como funciona el caso general (que es totalmente análogo).

Por $K$ denotaremos un cuerpo cuadrático real, esto es $K = \mathbb{Q}(\sqrt{d})$, donde $d > 0$ es un no cuadrado en $\mathbb{Q}$. Cualquiera de estos cuerpos tiene exactamente dos inmersiones reales, o sea dos maneras de meter al cuerpo dentro del cuerpo de números reales. Ellas están dadas por:

$$\tau_1(a + b\sqrt{d}) = a + b\sqrt{d}, \qquad \tau_2(a + b\sqrt{d}) = a - b\sqrt{d}.$$

A partir de las inmersiones $\tau_1$ ó $\tau_2$ podemos mirar al grupo $\mathrm{GL}_2(K)$ dentro de $\mathrm{GL}_2(\mathbb{R})$. Como en el semiplano complejo superior sólo actúan las matrices de determinante positivo, si $\Gamma$ es un subgrupo cualquiera de $\mathrm{GL}_2(K)$, denotamos por $\Gamma^+ = \{\gamma \in \Gamma \;:\; \det(\gamma) \gg 0\}$, donde $\nu \gg 0$ significa que tanto $\tau_1(\nu) > 0$ como $\tau_2(\nu) > 0$.

Notar que no podemos mirar el cociente por la acción de $\mathrm{GL}_2(K)^+$ dado que para que un cociente tenga estructura de variedad, precisamos que el grupo actuando sea discreto. La manera de solucionarlo es tomar no el semiplano complejo superior, sino dos copias del mismo, y mirar la acción de $\mathrm{GL}_2(K)^+$ a través de ambas. Luego, si denotamos por $\mathcal{O}_K$ el anillo de enteros de $K$, tiene sentido mirar la acción de $\mathrm{GL}_2(\mathcal{O}_K)^+$ en $\mathfrak{h} \times \mathfrak{h}$ dada por

$$\begin{pmatrix} a & b \\ c & d \end{pmatrix} \cdot (z_1, z_2) = \left( \frac{\tau_1(a)z_1 + \tau_1(b)}{\tau_1(c)z_1 + \tau_1(d)}, \frac{\tau_2(a)z_2 + \tau_2(b)}{\tau_2(c)z_2 + \tau_2(d)} \right).$$

Como sucede con el caso clásico, este cociente tiene una estructura de variedad diferenciable, pero no puede ser proyectiva por no ser compacta. Para obtener una compactificación, uno puede agregar las cúspides usuales, y considerar no $\mathfrak{h} \times \mathfrak{h}$, sino $\mathfrak{h}^2 \cup \mathbb{P}^1(K)$, y la acción de $\mathrm{GL}_2(\mathcal{O}_K)^+$ en dicho espacio, donde la acción de $\mathrm{GL}_2(\mathcal{O}_K)^+$ en $\mathbb{P}^1(K)$ está dada por el producto matricial, o sea $\left(\begin{smallmatrix} a & b \\ c & d \end{smallmatrix}\right) [\alpha : \beta] = [a\alpha + b\beta : c\alpha + d\beta]$.

**Ejercicio 12.** Recordar que todo ideal $\mathfrak{a}$ de $K$ está generado por dos elementos de $K$, o sea $\mathfrak{a} = \langle \alpha, \beta \rangle$. Probar que la asignación entre $\mathbb{P}^1(K)$ y clases de ideales de $K$, que a $[\alpha : \beta]$ le asocia $\langle \alpha, \beta \rangle$ da una biyección entre las cúspides de $\mathrm{GL}_2(\mathcal{O}_K)$ y el grupo de clases de ideales. Luego, el número de cúspides necesarias para compactificar el cociente depende de un invariante del cuerpo, a saber el número de clases.

La superficie que obtenemos al compactificar agregando las cúspides (llamada la compactificación de Baily–Borel), resulta una variedad algebraica proyectiva (o sea está dada por ceros de polinomios dentro de algún espacio proyectivo), pero es una variedad singular. El problema es que los puntos que pusimos para hacerla compacta, son todos singulares (porque deberíamos haber puesto cosas de codimensión 1, o sea de dimensión 1, que son curvas, no puntos). Uno puede obtener una superficie no singular mediante un proceso de "blow up" de puntos singulares, pero no nos vamos a preocupar por estos tecnicismos (ver van der Geer (1988) para más detalles).

El problema de que pueden existir ideales no principales, no solamente afecta al número de cúspides, sino a toda la teoría de las formas de Hilbert, en esta formulación clásica. Dado un ideal $\mathfrak{b}$, definimos el grupo

$$\mathrm{GL}_2^+(\mathcal{O}_K,\mathfrak{b}) = \left\{\begin{pmatrix} a & b \\ c & d\end{pmatrix} \;:\; a,d\in\mathcal{O}_K, b\in\mathfrak{b}^{-1}, c\in\mathfrak{b}\right\}.$$

Estos grupos son todos maximales (¡y no conjugados si $\mathfrak{b},\mathfrak{b}_2$ pertenecen a distintas clases de ideales!), con lo cual para estudiar las formas modulares, debemos mirarlos a todos ellos. Igual que para las formas clásicas, también tenemos los subgrupos de congruencias. Dado un ideal $\mathfrak{n}$, definimos

$$\Gamma_0(\mathfrak{n},\mathfrak{b}) = \left\{\begin{pmatrix} a & b \\ c & d\end{pmatrix} \in \mathrm{GL}_2^+(\mathcal{O}_K,\mathfrak{b}) \;:\; c\in\mathfrak{n}\mathfrak{b}\right\}.$$

Vamos a denotar por $X_0(\mathfrak{n},\mathfrak{b})$ al cociente $\Gamma_0(\mathfrak{n},\mathfrak{b})\backslash\mathfrak{H}^2\cup\mathbb{P}^1(K)$. Esto tiene estructura de superficie, salvo las cúspides, y los puntos elípticos (los puntos donde hay matrices no diagonales que los tienen como puntos fijos) donde uno debería desingularizar la superficie.

**Definición.** Una función holomorfa $f:\mathfrak{H}^2\to\mathbb{C}$ es una forma modular de peso $\mathbf{k}=(k_1,k_2)$, donde $k_1,k_2\in\mathbb{N}$ para $\Gamma_0(\mathfrak{n},\mathfrak{b})$ si para todo $\gamma=\left(\begin{smallmatrix} a & b \\ c & d\end{smallmatrix}\right)\in\Gamma_0(\mathfrak{n},\mathfrak{b})$ vale

$$f(\gamma\cdot(z_1,z_2)) = \tau_1(\det(\gamma))^{-k_1/2}\tau_2(\det(\gamma))^{-k_2/2} (\tau_1(c)z_1+\tau_1(d))^{k_1}(\tau_2(c)z_2+\tau_2(d))^{k_2} f(z_1,z_2). \tag{22}$$

Denotamos por $M_{\mathbf{k}}(\Gamma_0(\mathfrak{n},\mathfrak{b}))$ al espacio de dichas funciones.

Notar la similitud con la definición clásica, salvo que ahora al tener la función mas variables, aparecen productos de los términos $j(g,i)$ del caso unidimensional. Hay un tema no menor que llama la atención en la definición, que es que en ningún momento pedimos que la función sea holomorfa en las cúspides (que sí lo pedimos en el caso clásico). Para poder entender esto, debemos entender quién es el grupo de isotropía de la cúspide de infinito (las matrices que la fijan), y poder entender en particular el desarrollo de Fourier de las formas de Hilbert. Es fácil ver, que dicho grupo es justamente

$$\left\{\begin{pmatrix} \varepsilon & \mu \\ 0 & 1\end{pmatrix} \;:\; \varepsilon\in\mathcal{O}_K^{\times,+}, \mu\in\mathfrak{b}^{-1}\right\}.$$

Si nos olvidamos de las unidades por un momento (ellas jugarán un rol fundamental mas adelante), tenemos invariancia por el retículo $\mathfrak{b}^{-1}$, vía la identificación $\mu \to \left(\begin{smallmatrix} 1 & \mu \\ 0 & 1 \end{smallmatrix}\right)$. Pero si nuestra función es invariante por un retículo $M$, al considerar su expansión de Fourier, debemos mirar los caracteres de dicho retículo. Ellos están dados por las funciones $\exp(2\pi i(\nu_1 z_1 + \nu_2 z_2))$, donde $(\nu_1, \nu_2) \in M^\perp$, donde

$$M^\perp = \{\nu \in K \ : \ \mathrm{Tr}(\nu\mu) \in \mathbb{Z} \ \forall \mu \in M\},$$

esto es el dual con respecto a la forma cuadrática dada por la traza (como se vió en las notas de R. Miatello). En resumen, si $f(z_1, z_2) \in M_{\mathbf{k}}(\Gamma_0(\mathfrak{n}, \mathfrak{b}))$, su expansión de Fourier es de la forma

$$f(z_1, z_2) = \sum_{\nu \in \mathfrak{b}\delta^{-1}} a_\nu \exp(2\pi i z_1 \tau_1(\nu)) \exp(2\pi i z_2 \tau_2(\nu)), \tag{23}$$

donde $\delta$ es el diferente de $\mathcal{O}_K$ (su dual con respecto a la forma traza anterior). Como la expansión está definida en términos de pares, la condición de holomorfía en las cúspides debería decir que si $\nu$ tiene alguna de sus dos inmersiones negativas, entonces el coeficiente $a_\nu$ es cero.

**Teorema 13** (Principio de Koecher)**.** Si $f(z_1, z_2)$ es una forma de Hilbert, entonces todos los coeficientes de Fourier $a_\nu$ donde $\tau_1(\nu) < 0$ o $\tau_2(\nu) < 0$ son nulos.

El principio de Koecher dice justamente que la condición de holomorfía es automática para cualquier forma de Hilbert, y por lo tanto no es necesario agregarla. Decimos que una forma de Hilbert es *cuspidal* si el término constante de la expansión de Fourier (notar que el 0 siempre está en el dual de un retículo) es nulo en todas las cúspides. Dicho espacio lo denotamos $S_{\mathbf{k}}(\Gamma_0(\mathfrak{n}, \mathfrak{b}))$.

*Observación.* Si el peso $\mathbf{k} = (k_1, k_2)$ tiene las dos coordenadas iguales, se llama un *peso paralelo*, y se dice que la forma es de peso paralelo $k$ (estas son las formas mas interesantes desde un punto de vista geométrico de la superficie). Si el peso no es paralelo, entonces el Principio de Koecher implica que toda forma modular es automáticamente cuspidal también.

Si $f, g \in M_{\mathbf{k}}(\Gamma)$, y una de ellas es cuspidal, podemos definir un producto interno entre ellas (el producto de *Petersson*) de manera similar al caso clásico, esto es

$$\langle f, g \rangle = \iint_{\Gamma \backslash \mathfrak{H}^2} f(z_1, z_2) \overline{g(z_1, z_2)} y_1^{k_1} y_2^{k_2} \frac{dx_1 dy_1}{y_1^2} \frac{dx_2 dy_2}{y_2^2}. \tag{24}$$

**Ejercicio 13.** Ver que para la acción de $\mathrm{SL}_2(\mathcal{O}_K)$, $\frac{dx_1 dy_1}{y_1^2} \frac{dx_2 dy_2}{y_2^2}$ es una medida invariante en $\mathfrak{H}^2$. Verificar que la integral recién definida converge si alguna de las dos funciones es cuspidal.

**Teorema 14.** El espacio vectorial $M_{\mathbf{k}}(\Gamma_0(\mathfrak{n}, \mathfrak{b}))$ tiene dimensión finita.

Existen fórmulas para calcular la dimensión del espacio de formas cuspidales (que se obtienen de aplicar el Teorema de Riemann Roch en dimensiones mayores), pero son bastante difíciles de expresar e involucran cálculos de desingularización de puntos elípticos y cúspides. Ver van der Geer (ibíd.) Proposición 4.1 por ejemplo para el caso en que no hay puntos elípticos.

**3.2 Serie L de formas de Hilbert.** Queremos poner toda la información de los coeficientes de Fourier de una forma modular en un objeto analítico que codifique su información. El problema es que hay demasiados coeficientes de Fourier, y los coeficientes de Fourier de una forma de Hilbert no son todos independientes (como en el caso clásico), sino que hay muchas relaciones entre ellos, debido a las unidades (que juegan un rol fundamental como mencionamos anteriormente). Si $\varepsilon \in \mathcal{O}_K^\times$ es una unidad, la matriz $\left(\begin{smallmatrix} \varepsilon & 0 \\ 0 & 1 \end{smallmatrix}\right) \in \Gamma_0(\mathfrak{n}, \mathfrak{b})$, y su acción en $(z_1, z_2)$ está dada por $(\tau_1(\varepsilon)z_1, \tau_2(\varepsilon)z_2)$. Luego, si $f \in M_{\mathbf{k}}(\Gamma)$,

$$f(z_1, z_2) = \tau_1(\varepsilon)^{k_1/2}\tau_2(\varepsilon)^{k_2/2} f(\tau_1(\varepsilon)z_1, \tau_2(\varepsilon)z_2).$$

Si miramos la expansión de Fourier de ambos lados de la igualdad llegamos a la fórmula

$$a_{\varepsilon\nu} = a_\nu \tau_1(\varepsilon)^{k_1/2}\tau_2(\varepsilon)^{k_2/2}. \tag{25}$$

Denotemos por $U_K^+$ al conjunto de unidades de $\mathcal{O}_K$ que son totalmente positivas. Luego si $\mathbf{k}$ es paralelo y $\varepsilon \in U^+$, $\tau_1(\varepsilon)^{k/2}\tau_2(\varepsilon)^{k/2} = \mathcal{N}(\varepsilon)^{k/2} = 1$, con lo cual el coeficiente de Fourier es invariante por esta acción (que da infinitos valores de $\nu$ distintos).

Para un valor cualquiera $\mathbf{k}$, la ecuación (25) implica que la expresión

$$a_\nu \tau_1(\nu)^{-k_1/2}\tau_2(\nu)^{-k_2/2}$$

es invariante por la acción de $U_K^+$, o sea que es un número asociado a un generador totalmente positivo del ideal $\langle \nu \rangle$. Notar que el exponente elegido no es único, si cambiamos ambos exponentes por un mismo número, también queda invariante (porque los elementos de $U_K^+$ tienen norma 1). Para evitar problemas que dependen de la normalización, vamos a restringirnos al caso de peso paralelo.

Definimos la L-serie asociada a una forma de Hilbert $f$ de peso paralelo como

$$L(f,s) = \sum_{\substack{\nu \in \mathfrak{b}\delta^{-1}/U^+ \\ \nu \gg 0}} \frac{a_\nu}{\mathcal{N}(\nu\delta\mathfrak{b}^{-1})^s} = \sum_{[\mathfrak{a}]=[\mathfrak{b}^{-1}\delta]} \frac{a_\mathfrak{a}}{\mathcal{N}(\mathfrak{a})^s}, \tag{26}$$

donde la suma es sobre los ideales de $\mathcal{O}_K$ en la misma clase estricta que $\mathfrak{b}^{-1}\delta$ del grupo de clases, y el coeficiente $a_\mathfrak{a}$ es cualquiera asociado a un generador totalmente positivo de dicho ideal.

Notar que la suma no involucra a todos los ideales enteros, sólo una clase del grupo de clases estricto. Para formas de peso paralelo y nivel 1, si miramos la función $L$ completa $\Lambda(f,s) = \mathcal{D}(K)^s(2\pi)^{-2s}\Gamma(s)^2L(f,s)$, vale

**Teorema 15.** La función $\Lambda(f,s)$ se extiende de manera meromorfa a todo el plano complejo y satisface la ecuación funcional

$$\Lambda(f,s) = (-1)^k\Lambda(f,k-s).$$

*Prueba.* Ver Bruinier (2008), Teorema 1.44. En dichas notas, sólo se consideran matrices de determinante 1 para el ideal $\mathfrak{b} = \mathcal{O}_K$, pero mirando cualquier otro maximal, y utilizando la ecuación funcional para la matriz $\begin{pmatrix} 0 & \frac{1}{B} \\ B & 0 \end{pmatrix}$, donde $B = \mathcal{N}(\mathfrak{b})$, la misma cuenta funciona. □

Como veremos en la próxima sección, para poder estudiar los operadores de Hecke, es preciso mirar las formas modulares para un conjunto de representantes del grupo de clases estricto de $K$ juntas. Sea $\{\mathfrak{b}_1, \ldots, \mathfrak{b}_h\}$ un tal conjunto de representantes. Definimos

$$M_{\mathbf{k}}(\mathfrak{n}) = \bigoplus_{i=1}^{h} M_{\mathbf{k}}(\Gamma_0(\mathfrak{n}, \mathfrak{b}_i)). \tag{27}$$

Con lo cual una forma de nivel $\mathfrak{n}$ es una tupla de funciones invariantes cada una de ellas por un subgrupo de congruencias diferente. Es en este conjunto donde podremos definir la acción de los operadores de Hecke. Si $\mathbf{f} = (f_1, \ldots, f_h) \in M_{\mathbf{k}}(\mathfrak{n})$, y $\mathfrak{m} \subset \mathcal{O}_K$ es un ideal cualquiera, definimos el coeficiente de Fourier $c(\mathfrak{m}, \mathbf{f}) = a_\nu(f_i)$, donde $\mathfrak{m} = \nu\mathfrak{b}_i$, con $\nu \gg 0$.

Luego la L-serie de una forma modular de Hilbert $\mathbf{f} \in M_{\mathbf{k}}(\mathfrak{n})$ es

$$L(\mathbf{f}, s) = \sum_{\mathfrak{m}} \frac{c(\mathfrak{m}, \mathbf{f})}{\mathcal{N}(\mathfrak{m})^s}.$$

Es claro que dicha L-serie coincide con la suma de las L-series de cada una de ellas.

**3.3 Operadores de Hecke.** La teoría de operadores de Hecke esta bien explicada en el artículo Shimura (1978) en término de cocientes dobles. La noción de cocientes dobles es el análogo discreto de la definición del conjunto $H_p$ de (20), y la demostración del Lema 12 muestra como pasar de clases dobles a coclases a derecha y obtener la formulación clásica de los operadores de Hecke. A pesar de la gran utilidad de los cocientes dobles y sus propiedades, resultan un poco avanzados para un primer estudio de las formas modulares (al lector interesado le recomendamos mirar la referencia clásica Shimura (1994)).

Lo que vamos a hacer en cambio es dar la definición del operador de Hecke en términos de coeficientes de Fourier de funciones en $M_{\mathbf{k}}(\mathfrak{n})$ como en (27).

**Definición.** Si $\mathbf{f} \in M_{\mathbf{k}}(\mathfrak{n})$ y $\mathfrak{p}$ es un ideal primo que no divide a $\mathfrak{n}$, el operador de Hecke $T(\mathfrak{p})$ actuando en $\mathbf{f}$ corresponde a la forma modular cuya expansión de Fourier está dada por

$$c(\mathfrak{m}, T(\mathfrak{p})\mathbf{f}) = \mathcal{N}(\mathfrak{p})c(\mathfrak{p}\mathfrak{m}, \mathbf{f}) + c(\mathfrak{m}\mathfrak{p}^{-1}, \mathbf{f}),$$

donde $c(\mathfrak{m}\mathfrak{p}^{-1}, \mathbf{f}) = 0$ si $\mathfrak{p} \nmid \mathfrak{m}$.

Notar que para calcular $T(\mathfrak{p})\mathbf{f}$, es necesario conocer no una sola componente de la forma $\mathbf{f}$, sino varias de ellas. Solamente los operadores de Hecke correspondientes a ideales principales están definidos en cada componente por separado (y definen una acción en cada $M_{\mathbf{k}}(\mathfrak{n}, \mathfrak{b})$).

Como en el caso clásico, el álgebra de operadores de Hecke con índice coprimo con $\mathfrak{n}$ es un álgebra conmutativa, y cada operador de Hecke es autoadjunto con respecto al producto de Petersson. Existe además una familia de operadores (de Atkin–Lehner) $W_{\mathfrak{p}}$ para cada primo $\mathfrak{p} \mid \mathfrak{n}$, y las autofunciones para todos estos operadores son las que satisfacen una ecuación funcional. En resumen, con las definiciones correctas, las formas modulares de Hilbert se comportan de manera análoga a las formas modulares clásicas.

**3.4 Interpretación automorfa.** Queremos hacer el análogo al Capítulo 2 para formas de Hilbert. Un ideal primo $\mathfrak{p}$ de $\mathcal{O}_K$ define una valuación, y podemos completar $K$ con respecto a dicha valuación. El cuerpo (completo) resultante lo denotamos $K_{\mathfrak{p}}$. Para definir los adèles necesitamos considerar todos los valores absolutos de $K$.

**Teorema 16** (Ostrowski)**.** Si $K$ es un cuerpo cuadrático real, todos los valores absolutos de $K$ son:

- Los arquimedianos dados por las inmersiones $\tau_1, \tau_2$, o sea $|x|_{\tau_i} = |\tau_i(x)|$.
- El valor absoluto noarquimediano $|x|_{\mathfrak{p}}$ asociado a un ideal primo $\mathfrak{p}$.

Si $S$ es un conjunto finito de primos, los $S$-adèles son

$$\mathbb{A}_K^S = \prod_{\mathfrak{p} \notin S} \mathcal{O}_{\mathfrak{p}} \times \prod_{\mathfrak{p} \in S} K_{\mathfrak{p}} \times \mathbb{R}_{\tau_1} \times \mathbb{R}_{\tau_2}.$$

Luego definimos el anillo de adèles como la unión de los $S$-adèles sobre todos los conjuntos finitos $S$. De manera análoga, definimos $\mathrm{GL}_2(\mathbb{A}_K)$. Definamos los compactos $K_{\mathfrak{p}}(\mathfrak{n})$ como en (18), y denotemos $K(\mathfrak{n}) = \prod_{\mathfrak{p}} K_{\mathfrak{p}}(\mathfrak{n})$.

**Teorema 17** (Aproximación fuerte)**.** Con las notaciones anteriores, la siguiente igualdad se verifica:

$$\mathrm{GL}_2(\mathbb{A}_K) = \bigsqcup_{i=1}^{h} \mathrm{GL}_2(K) \begin{pmatrix} 1 & 0 \\ 0 & t_i \end{pmatrix} \left( \mathrm{GL}_2(\mathbb{R})_{\tau_1}^{+} \times \mathrm{GL}_2(\mathbb{R})_{\tau_2}^{+} \times K(\mathfrak{n}) \right),$$

donde $t_i \in \mathbb{A}_K^{\times}$ corresponde al ideal $\mathfrak{b}_i^{-1}$, esto es $\mathfrak{b}_i = \prod_{\mathfrak{p}} \mathfrak{p}^{v_{\mathfrak{p}}(t_i(\mathfrak{p}))}$.

En el caso en que el número de clases estricto sea 1, entonces aproximación fuerte dice exactamente lo mismo que en $\mathbb{Q}$, pero en general hacen falta todas las superficies $X_0(\mathfrak{n}, \mathfrak{b})$.

**Proposición 18.** Hay una identificación entre

$$\mathrm{GL}_2(K)\backslash \mathrm{GL}_2(\mathbb{A}_K)/ \left(\mathrm{SO}(2)_{\tau_1} \times \mathrm{SO}(2)_{\tau_2} \times K(\mathfrak{n})\right) \leftrightarrow \bigcup_{i=1}^{n} \Gamma(\mathfrak{n}, \mathfrak{b}_{\mathfrak{i}})\backslash \mathfrak{h}^2.$$

*Prueba.* Denotemos por $K_\infty$ al grupo $\mathrm{SO}(2)_{\tau_1} \times \mathrm{SO}(2)_{\tau_2}$ y por $M_i = \begin{pmatrix} 1 & 0 \\ 0 & t_i \end{pmatrix}$. Por el Teorema de aproximación fuerte, tenemos

$$\begin{aligned}
&\mathrm{GL}_2(K)\backslash \mathrm{GL}_2(\mathbb{A}_K)/K_\infty K(\mathfrak{n}) = \\
&\qquad = \bigsqcup_{i=1}^{h} \mathrm{GL}_2(K)\backslash \left(\mathrm{GL}_2(K) M_i \,\mathrm{GL}_2(\mathbb{R})_\infty^+ \times K(\mathfrak{n})\right)/K_\infty K(\mathfrak{n}) \\
&\qquad\qquad = \bigsqcup_{i=1}^{h} (M_i^{-1}\,\mathrm{GL}_2(\mathbb{R})_\infty^+ K(\mathfrak{n}) M_i \cap \mathrm{GL}_2(K))\backslash \mathrm{GL}_2(\mathbb{R})_\infty^+/K_\infty.
\end{aligned}$$

Pero $\mathrm{GL}_2(\mathbb{R})_\infty^+$ lo podemos identificar con $\mathfrak{h}^2$ mirando la acción en el punto $(i,i)$, y justamente el estabilizador es $K_\infty$. Es elemental verificar que $M_i^{-1}\,\mathrm{GL}_2(\mathbb{R})_\infty^+ K(\mathfrak{n}) M_i \cap \mathrm{GL}_2(K) = \Gamma_0(\mathfrak{n}, \mathfrak{b}_i)$. □

Luego podemos definir las formas automorfas de Hilbert copiando la definición clásica, y sabemos que coincide con las formas $\mathbf{f} \in M_{\mathbf{k}}(\mathfrak{n})$. Esto justifica el por qué debemos mirar tuplas de funciones en lugar de una sola de ellas.

**Definición.** Una *forma automorfa de Hilbert* de peso $\mathbf{k}$ y nivel $K(\mathfrak{n})$ para $K$ es una función $\phi$ : $\mathrm{GL}_2(\mathbb{A}_K) \to \mathbb{C}$ que satisface:

1. $\phi(\gamma g) = \phi(g)$ para todo $\gamma \in \mathrm{GL}_2(K)$.
2. $\phi(g k_0) = \phi(g)$ para todo $k_0 \in K(\mathfrak{n})$.
3. $\phi(g(r(\theta_1)_{\tau_1}, r(\theta_2)_{\tau_2})) = \exp(i(k_1\theta_1 + k_2\theta_2))\phi(g)$ para todo $\theta 1, \theta_2 \in [0, 2\pi]$.
4. La función $\phi$, vista como función de $\mathrm{GL}_2(\mathbb{R})_\infty^+$ es $C^\infty$ y satisface las ecuaciones diferenciales
$$\Delta\phi_{\tau_j} = -\frac{k_j}{2}\left(\frac{k_j}{2} - 1\right)\phi, \text{ para } j = 1, 2.$$
5. $\phi(zg) = \phi(g)$ para todo $z \in Z_{\mathbb{A}} = \left\{\left(\begin{smallmatrix} t & 0 \\ 0 & t \end{smallmatrix}\right) : t \in \mathbb{A}_K^\times\right\}$.
6. (Crecimiento moderado) Para todo $c > 0$ y todo conjunto compacto $\Omega$ de $\mathrm{GL}_2(\mathbb{A}_K)$, existen constantes $C, N$ tales que
$$\left|\phi\left(\left(\begin{smallmatrix} a & 0 \\ 0 & 1 \end{smallmatrix}\right) g\right)\right| \leqslant C|a|^N.$$
7. Además, $\phi$ se dice cuspidal si satisface
$$\int_{K\backslash \mathbb{A}_K} \phi\left(\left(\begin{smallmatrix} 1 & x \\ 0 & 1 \end{smallmatrix}\right) g\right) dx = 0 \text{ para casi todo } g.$$

Existe una biyección entre las formas automorfas de Hilbert y las formas de Hilbert, que sigue la correspondencia clásica. Además, usando las formas automorfas de Hilbert es claro como definir los operadores de Hecke, dados como un operador de convolución, y una demostración similar a la dada en el Lema 12 muestra la buena definición de dichos operadores, y sus propiedades.

**3.5 Digresión de formas de Hilbert.** Como ya mencionamos en reiteradas ocasiones, las formas de Hilbert son el análogo de las formas clásicas, pero para cuerpos de números totalmente reales. El trato que le hemos dado hasta acá sigue la línea original de estudio, y resulta muy fructífero para muchas aplicaciones. No obstante, el hecho de pasar de una curva (dimensión 1) a variedades de dimensión mayor, tiene sus claras desventajas. La primera de ellas es de índole computacional, dado que para calcular el espacio de formas de Hilbert es necesario trabajar con espacios de dimensión grande para los cuales es difícil calcular la cohomología (en la práctica, los métodos modernos para calcularlas se separan en dos: calcular las formas usando funciones theta, como explicado en el curso de G. Tornaría, o calcular explícitamente la cohomología de las curvas de Shimura que se introducirán en el curso de M. Harris).

Un gran problema teórico es que las formas de Hilbert también deberían tener representaciones de Galois asociadas, pero esto no es nada inmediato a partir de su definición. En el caso en que $[K : \mathbb{Q}]$ sea impar, o en que la forma de Hilbert $\mathbf{f}$ posea un primo Steinberg o supercuspidal, esto se puede hacer gracias a las curvas de Shimura (y será explicado en el curso de M. Harris). En general, la manera es mediante un proceso de aproximación $p$-ádico, y fue hecho por Taylor.

A la vez, si nos restringimos a formas de peso paralelo 2, estas están asociadas a formas diferenciales, pero ahora la forma diferencial $f(z_1, z_2)dz_1dz_2$ es una 2-forma diferencial de la superficie, o sea que vive en el $H^2(X_0(\mathfrak{n}, \mathfrak{b}), \mathbb{C})$. No hay ninguna construcción natural para asociarle un retículo de rango 2 a una forma modular que sea autofunción para los operadores de Hecke con autovalores racionales (como sucede en el caso clásico). Nuevamente, el uso de Curvas de Shimura nos permitirá en muchos casos poder resolver este problema también.

**3.6 Caso general.** A pesar de que sólo miramos formas de Hilbert para cuerpos cuadráticos reales, la misma se generaliza de manera directa a cuerpo de números totalmente real $K$. Estos son cuerpos de la forma $K = \mathbb{Q}[\alpha]$, donde $\alpha$ es raíz de un polinomio racional irreducible cuyas raíces son todas reales. A tal cuerpo podemos asociarle las llamadas *inmersiones reales*, estas son todas las maneras distintas de meter al cuerpo dentro del cuerpo de números reales, y justamente están dadas por enviar $\alpha$ a cualquier otra raíz de su polinomio minimal (hay justamente $n$ de ellas, donde $n$ es el grado del minimal de $\alpha$ que coincide con el grado de la extensión $[K : \mathbb{Q}]$). Ahora todas las cuentas están indexadas no por dos, sino por $n$ parámetros, y dichos parámetros se pueden enumerar a partir de las inmersiones del cuerpo $K$. Notar que el peso de una forma de Hilbert está indexado por inmersiones, y no hay manera de "ordenarlas". Todo lo hecho en las secciones anteriores sigue valiendo con exactamente la misma demostración, pero hemos obviado las cuentas para evitar el uso de índices que simplemente dificultan la presentación.

## Referencias

Ariel Pacetti
apacetti@famaf.unc.edu.ar
arielpacetti@gmail.com
Universidad de Buenos Aires

# CURVAS DE SHIMURA

Michael Harris

## Índice general

## 1 Grupos cuaterniónicos como grupos fuchsianos

**1.1 Grupos discretos de cuaterniones y plano superior de Poincaré.** Sea $D$ un álgebra de cuaterniones sobre un cuerpo totalmente real $F$. Sea $\Sigma$ el conjunto de primos arquimedianos de $F$. Sea $G$ el grupo multiplicativo de $D$. Para cada primo $v$ de $F$ (arquimediano o no) tenemos la completación $D_v$ sobre el cuerpo topológico $F_v$, y ponemos $G_v = D_v^\times$. El álgebra de adèles de $D$ es el producto restringido $D(\mathbf{A}) = \prod'_v D_v$ sobre todos los primos de $F$, definido exactamente como en el caso de un cuerpo de números. Es decir, en cada primo no arquimediano $v$, $D_v$ contiene una subálgebra compacta maximal $\mathcal{O}_{D_v} \subset D_v$; si $D_v \xrightarrow{\sim} M(2, F_v)$ entonces $\mathcal{O}_{D_v} = M(2, \mathcal{O}_v)$. El álgebra $D(\mathbf{A})$ es el subconjunto de $(x_v) \in \prod_v D_v$ donde $x_v \in \mathcal{O}_{D_v}$ salvo en un número finito de primos. Además hay un homomorfismo inyectivo $F_{\mathbf{A}} := \mathbf{A}_F \to D(\mathbf{A})$ cuya imagen es igual al centro de $D(\mathbf{A})$.

Del mismo modo, definimos el grupo de $D$ como el producto restringido $D^\times(\mathbf{A}) = \prod'_v D_v^\times = \prod'_v G_v$. Escribimos $D^\times(\mathbf{A}) = D_\infty^\times \times D^\times(\mathbf{A}_f)$, donde $D_\infty^\times = \prod_{v\in\Sigma} D_v^\times$ y $D^\times(\mathbf{A}_f) = \prod'_{v \text{ finito}} D_v^\times$. La norma reducida $\nu : D \to F$ es un mapeo multiplicativo y así define homomorfismos locales y globales:

$$\nu_v : G_v \to F_v^\times, \quad \nu : G \to F^\times, \quad \nu_{\mathbf{A}} : D^\times(\mathbf{A}) \to F_{\mathbf{A}}^\times.$$

La norma local $\nu_v$ es suryectiva si $v$ es un primo finito o si $D_v \simeq M(2,\mathbb{R})$; si $D_v \simeq \mathbb{H}$, el álgebra de cuaterniones de Hamilton, la imagen de la norma local $\nu_v$ es el grupo de números reales *positivos*.

Escribimos $D_F^\times$ para los elementos globales, es decir el grupo multiplicativo del álgebra $D$. Podemos considerar $D_F^\times$ como subgrupo de $D^\times(\mathbf{A})$. Dado un subgrupo compacto abierto $K \subset$ $D^\times(\mathbf{A}_f)$, el *grupo de congruencia de nivel $K$*, $\Gamma_K \subset D_F^\times$, es la intersección en $D^\times(\mathbf{A})$ de $D_F^\times$ con $K$.

Sea $\Sigma_D \subset \Sigma$ (resp. $\Sigma'_D \subset \Sigma$) el subconjunto de los $v$ no ramificados (resp. ramificados), para los cuales $D_v \simeq M(2,\mathbb{R})$ (resp. $D_v \simeq \mathbb{H}$). Sea $\mathfrak{H}^\pm = \mathbb{C} \smallsetminus \mathbb{R}$, la unión de los semiplanos superior e inferior. Hay una acción del grupo $D_\infty^\times = \prod_{v\in\Sigma} D_v^\times$ sobre $\mathfrak{H}^{\pm,\Sigma_D} = \prod_{v\in\Sigma_D} \mathfrak{H}_v^\pm$: si $v \in \Sigma_D$ el factor $D_v^\times \xrightarrow{\sim} GL(2,\mathbb{R})$ actúa sobre el factor $\mathfrak{H}_v^\pm$ y los factores $D_v^\times$ con $v \in \Sigma'_D$ actúan trivialmente. Sea $\Gamma_K \subset D_F^\times$ el grupo de congruencia de nivel $K$. Vía la inclusión $\Gamma_K \hookrightarrow D_\infty^\times$ tenemos una acción de $\Gamma_K$ sobre $\mathfrak{H}^{\pm,\Sigma_D}$.

Fijamos un punto $h \in \mathfrak{H}^{\pm,\Sigma_D}$ y definimos $\tilde{K}_\infty = \tilde{K}_h \subset D_\infty^\times$ [con virgulilla] como el estabilizador de $h$:

$$\tilde{K}_\infty = \{g \in D_\infty^\times \mid g(h) = h\}. \tag{1-1}$$

El grupo $\tilde{K}_\infty$ contiene al centro $Z_\infty = \prod_{v\in\Sigma} F_v^\times$ de $D_\infty^\times$, y el cociente $\tilde{K}_\infty/Z_\infty$ es compacto. Sea $K_\infty = K_h$ el subgrupo compacto maximal de $\tilde{K}_\infty$ (hay solo uno), que además es compacto maximal conexo en $D_\infty^\times$. La definición de formas modulares sobre $D^\times(\mathbf{A})$ (ver más abajo) depende de la elección de un subgrupo compacto maximal $K_\infty \subset D_\infty^\times$, pero distintos $K_\infty$ dan teorías equivalentes (espacios isomorfos) de formas modulares.

**Proposición 1-2.** *La acción de $\Gamma_K$ sobre $\mathfrak{H}^{\pm,\Sigma_D}$ es propiamente discontinua. Si $D$ es un álgebra de división, entonces el cociente $\Gamma_K\backslash\mathfrak{H}^{\pm,\Sigma_D}$ es compacto.*

La discontinuidad propia se muestra como en el caso de la acción de $GL(2,\mathbb{R})^+$ sobre el semiplano superior.

**Ejercicio 1.1.** Demostrar la discontinuidad propia.

Para mostrar la compacidad, utilizamos el teorema siguiente del libro *Basic Number Theory* de Weil (Theorem 4 del capítulo IV, sección 3) :

**Teorema 1-3.** *Sea $D$ un álgebra de división de dimensión finita sobre $F$. Para cada número real $\mu \geqslant 1$, sea*

$$D_\mu = \{d \in D^\times(\mathbf{A}) \ : \ ||d||_{\mathbf{A}} \leqslant \mu, ||d||^{-1} \geqslant \mu^{-1}\}$$

*Entonces $D_\mu$ es un conjunto cerrado en $D^\times(\mathbf{A})$ y la imagen de $D_\mu$ en $D^\times\backslash D^\times(\mathbf{A})$ es compacta.*

Si $D = F$, el enunciado sigue del argumento de la geometría de números utilizado para demostrar la finitud del número de clases y el teorema de Dirichlet. El caso de un álgebra no conmutativa es exactamente lo mismo. En particular, si $D_1 = \ker || \bullet || : D^\times(\mathbf{A}) \to \mathbb{R}^\times$, el teorema de Weil implica que $D^\times\backslash D_1/K \cdot K_\infty$ es compacto.

El cociente $\Gamma_K\backslash\mathfrak{H}^{\pm,\Sigma_D}$ es un ejemplo de una *variedad de Shimura conexa* de dimensión $|\Sigma_D|$. Cuando $|\Sigma_D| = 1$, es una *curva de Shimura conexa*. (Es un abuso escribir eso; esa curva no es necesariamente conexa, porque $\mathfrak{H}^\pm$ tiene 2 componentes conexas. Pero es más sencillo no separar las componentes conexas de $\mathfrak{H}^\pm$.) Resulta de la proposición que el subgrupo de $\Gamma_K$ que estabiliza la componente conexa $\mathfrak{H}^+$ de $\mathfrak{H}^\pm$ es un grupo fuchsiano, de modo que $\Gamma_K\backslash\mathfrak{H}^{\pm,\Sigma_D}$ es isomorfa a la unión de uno o dos cocientes de $\mathfrak{H}^+$ por un grupo fuchsiano. Todavía no es ese el objeto que nos interesa. Cuando el grupo de clases $h_F$ es de orden mayor a 1 no se pueden definir operadores de Hecke de modo natural sobre esas variedades conexas. Para eso tenemos que trabajar con las variedades de Shimura adélicas. Además, para las aplicaciones aritméticas, hay que construir modelos de las variedades de Shimura sobre cuerpos de números.

En el resto del curso, siempre vamos a suponer que $|\Sigma_D| = 1$; entonces la variedad es una curva de Shimura. Para evitar dificultades técnicas es mejor en este curso trabajar con el grupo $PD = D^\times/F^\times$, y con los grupos locales $PD_\infty = D_\infty^\times/F_\infty^\times$ y adélicos $PD(\mathbf{A}) = D^\times(\mathbf{A})/F_\mathbf{A}^\times$, $PD(\mathbf{A}_f) = D^\times(\mathbf{A}_f)/F^\times(\mathbf{A}_f)$. Eso tiene sentido porque el centro $F_\infty^\times$ de $D_\infty^\times$ actúa trivialmente sobre $\mathfrak{H}^{\pm,\Sigma_D} = \mathfrak{H}^\pm$. Una *curva de Shimura adélica* de nivel $K \subset PD(\mathbf{A}_f) = D^\times(\mathbf{A}_f)/F^\times(\mathbf{A}_f)$ es el cociente

$$_KS(D) = PD\backslash[\mathfrak{H}^{\pm,\Sigma_D} \times PD(\mathbf{A}_f)/K].$$

Exactamente como en el caso de la variedad modular de Hilbert, hay un conjunto finito $U = \{u_i, i \in I\}$ de elementos de $PD(\mathbf{A}_f)$ y una descomposición

$$PD(\mathbf{A}) = D^\times(\mathbf{A})/F_\mathbf{A}^\times = \coprod_{i\in I} PDu_i[PD_\infty \times K].$$

Pero esa descomposición es más sencilla que en el caso de $GL(2, F)$ cuando $D$ es un álgebra de división; como Weil muestra en su libro *Basic Number Theory*, la demostración de la finitud del número de clases se aplica sin cambios a un tal cociente. Entonces, si escribimos

$$_KS(D) = PD\backslash[\mathfrak{H}^{\pm,\Sigma_D} \times (PD(\mathbf{A}_f)/K)] = PD\backslash[\coprod_{i\in I} \mathfrak{H}^{\pm,\Sigma_D} \times PDu_iK/K]$$

se puede escribir

$$\coprod_{i\in I}[PD \cap u_iKu_i^{-1}\backslash\mathfrak{H}^{\pm,\Sigma_D}] = \coprod_{i\in I} \Gamma_i\backslash\mathfrak{H}^{\pm,\Sigma_D}$$

donde hemos escrito

$$\Gamma_i = PD \cap u_iKu_i^{-1} = \Gamma_{u_iKu_i^{-1}}$$

en nuestra notación anterior.

En efecto, si $x, x' \in \mathfrak{H}^{\pm,\Sigma_D}$, entonces las imagenes de $xu_i$ y $x'u_j$ en el cociente

$$PD\backslash[\coprod_{i\in I} \mathfrak{H}^{\pm,\Sigma_D} \times PDu_i \times K]/K$$

coinciden si, y solamente si, $u_i = u_j$ y existen $d \in PD, k \in K$, con

$$dxu_ik = x'u_i \Leftrightarrow x' = dx[u_iku_i^{-1}].$$

Pero como $x', x \in \mathfrak{H}^{\pm,\Sigma_D}$ y $u_iku_i^{-1} \in D^\times(\mathbf{A}_f)/F^\times(\mathbf{A}_f)$, eso quiere decir que $d[u_iku_i^{-1}] = 1 \in D^\times(\mathbf{A}_f)/F^\times(\mathbf{A}_f)$; de modo que $d \in \Gamma_i$. Así la imagen de $xu_i$ en $\mathfrak{H}^{\pm,\Sigma_D}$ está bien definida en $\Gamma_i\backslash\mathfrak{H}^{\pm,\Sigma_D}$.

Ya hemos demostrado la siguiente proposición:

**Proposición 1-4.** *El cociente adélico*

$$_KS(D) = (D^\times/F^\times)\backslash[\mathfrak{H}^{\pm,\Sigma_D} \times (D^\times(\mathbf{A}_f)/F^\times(\mathbf{A}_f))/K]$$

*es una unión finita de curvas de Shimura conexas.*

**1.2 Curvas de Shimura: teoría analítica de formas modulares.** Todas las consideraciones en esta sección son válidas para cualquier álgebra cuaterniónica de división $D$. Pero como ya hemos dicho, siempre vamos a suponer que $|\Sigma_D| = 1$, y escribiremos $\mathfrak{H}^{\pm}$ en vez de $\mathfrak{H}^{\pm,\Sigma_D}$.

Sea $v$ el único primo arquimediano en $\Sigma_D$, y sea $\tilde{K}_v = \tilde{K}_\infty \cap D_v^\times$. Entonces $D_v \simeq M(2,\mathbb{R})$, y tenemos un isomorfismo $D_v^\times/\tilde{K}_v = \mathfrak{H}^{\pm}$: entonces $\tilde{K}_v$ es el estabilizador de $h$ y el mapeo

$$p_h : D_v^\times/\tilde{K}_v \ \to \ \mathfrak{H}^{\pm}, d \mapsto d(h)$$

es un isomorfismo $C^\infty$. El estabilizador $\tilde{K}_\infty$ de $h$ es entonces igual a $\tilde{K}_v \times \prod_{w\neq v} D_w^\times$, y tenemos también un isomorfismo $C^\infty$

$$p_{h,\infty} : \prod_w D_w^\times/(\tilde{K}_v \times \prod_{w\neq v} D_w^\times) \ \xrightarrow{\sim} \mathfrak{H}^{\pm}. \tag{1-5}$$

El grupo $\tilde{K}_v$ estabiliza el punto $h$ y por esta razón actúa sobre el espacio tangente complexificado $T_{\mathfrak{H}^{\pm},h,\mathbb{C}} = T_{\mathfrak{H}^{\pm},h} \otimes \mathbb{C}$ de $\mathfrak{H}^{\pm}$ en $h$. Sabemos además que la acción de $D_v^\times$ sobre $\mathfrak{H}^{\pm}$ conserva la estructura analítica del espacio. Así la acción de $\tilde{K}_v$ sobre $T_{\mathfrak{H}^{\pm},h,\mathbb{C}}$ conserva la descomposición en subespacios holomorfo y anti-holomorfo, cada uno de dimensión uno:

$$T_{\mathfrak{H}^{\pm},h,\mathbb{C}} \ \xrightarrow{\sim} T_h^{hol} \oplus \bar{T}_h^{hol}.$$

El diferencial del mapeo $p_h$ es un mapeo suryectivo

$$dp_h : T_{D_v^\times,1,\mathbb{C}} \ \to \ T_{\mathfrak{H}^{\pm},h,\mathbb{C}}.$$

Denotamos $\mathfrak{g}_v$ el álgebra de Lie del grupo de Lie $D_v^\times \xrightarrow{\sim} GL(2,\mathbb{R})$, e identificamos $\mathfrak{g}_v$ con $M(2,\mathbb{R})$. Como $T_{D_v^\times,1} = \mathfrak{g}_v$, el diferencial es un mapeo suryectivo

$$M(2,\mathbb{C}) \ \to \ T_{\mathfrak{H}^{\pm},h,\mathbb{C}}$$

que además es un homomorfismo de representaciones del grupo $\tilde{K}_v$. Aquí la acción de $\tilde{K}_v$ sobre $M(2,\mathbb{C})$ está dada por conjugación (representación adjunta):

$$ad(k)(X) = kXk^{-1}, \ \ k \in \tilde{K}_v, X \in M(2,\mathbb{C}).$$

De aquí en adelante, $h$ designa el punto $i \in \mathfrak{H}^+$. Con esta convención, el estabilizador $\tilde{K}_v$ en $GL(2,\mathbb{R})$ es nada más que el subgrupo de matrices $\begin{pmatrix} a & b \\ -b & a \end{pmatrix}$ con $a, b \in \mathbb{R}$ y $a^2 + b^2 \neq 0$. Y el mapeo $a + bi \mapsto \begin{pmatrix} a & b \\ -b & a \end{pmatrix}$ define un isomorfismo $\mathbb{C}^\times \xrightarrow{\sim} \tilde{K}_v$. Es fácil descomponer $M(2,\mathbb{C})$ en espacios propios para la acción de $\tilde{K}_v$: si

$$\tilde{\mathfrak{k}}_v = \mathrm{Lie}(\tilde{K}_v) \otimes \mathbb{C} = \left\{ \begin{pmatrix} a & b \\ -b & a \end{pmatrix}, \ \ a, b \in \mathbb{C} \right\},$$

entonces

$$M(2,\mathbb{C}) = \tilde{\mathfrak{k}}_v \oplus \mathfrak{p}^+ \oplus \mathfrak{p}^-$$

con
$$\mathfrak{p}^+ = \mathbb{C}X^+, \mathfrak{p}^- = \mathbb{C}X^-, X^+ = \begin{pmatrix} 1 & i \\ i & -1 \end{pmatrix}, X^- = \begin{pmatrix} 1 & -i \\ -i & -1 \end{pmatrix}.$$
Resulta de un cálculo fácil que

**Lema 1-6.** *El diferencial $dp_h$ identifica $\mathfrak{p}^+ \xrightarrow{\sim} T_h^{hol}$, $\mathfrak{p}^- \xrightarrow{\sim} \bar{T}_h^{hol}$. Una función $\phi : \mathfrak{H}^\pm \to \mathbb{C}$ es holomorfa si y solamente si la función compuesta $\Phi = \phi \circ p_h : D_v^\times \to \mathbb{C}$ es solución de la ecuación diferencial $X^-\Phi = 0$.*

**Ejercicio 1.2.** Demostrar el Lema 1-6.

Más generalmente, ponemos $\tilde{\mathfrak{k}} = \mathrm{Lie}(\tilde{K}_\infty) = \tilde{\mathfrak{k}}_v \oplus (\oplus_{w \neq v} \mathrm{Lie}(D_w^\times))$.

**1.2.1 Formas modulares y funciones adélicas.** Una *forma modular clásica* sobre la curva (no conexa) ${}_KS(D)$ es una función holomorfa sobre $\mathfrak{H}^{\pm,\Sigma_D} \times (D^\times(\mathbf{A}_f)/K)$ que satisface una ecuación funcional que corresponde a sus pesos, que en nuestra situación son enteros pares (porque todas nuestras formas son supuestas invariantes por la acción de $F^\times(\mathbf{A})$). Pero como la componente $D_w^\times$, con $w \neq v$, es compacta módulo $F_w^\times$, no podemos definir un factor como $c_w z_w + d_w$ en la coordenada $w$. En vez de eso, tenemos que utilizar las representaciones irreducibles del grupo $D_w^\times$, que generalmente son de dimensión superior a 1; es decir, tenemos que introducir *formas modulares con valores vectoriales*.

Podemos tratar todos los primos arquimedianos de manera homogénea. Sea
$$\rho : \tilde{K}_\infty / (\prod_{w|\infty} F_w^\times) \to GL(W)$$
una representación irreducible, con $W = W_\rho$ un espacio vectorial complejo de dimensión finita. Una tal representación admite una factorización
$$\rho = \rho_v \otimes (\otimes_{w \neq v} \rho_w),$$
$$\rho_v : \tilde{K}_v / F_v^\times \to \mathbb{C}^\times, \ \ \rho_w : \tilde{K}_w / F_w^\times \xrightarrow{\sim} \mathbb{H}^\times / \mathbb{R}^\times \to GL(W_w)$$
donde todas las $\rho_w$ son irreducibles. Definimos el *factor de automorfía*
$$j_\rho : D_v^\times \times \prod_{w \neq v} D_w^\times \times \mathfrak{H}^\pm \to GL(W) = GL(\otimes_{w \neq v} W_w);$$

$$j_\rho\left(\begin{pmatrix} a & b \\ c & d \end{pmatrix}, (g_w), z\right) = \rho_v(cz + d) \cdot \otimes_{w \neq v} \rho_w(g_w). \tag{1-7}$$

**Ejercicio 1.3.** Demostrar que $j_\rho$ satisface la ecuación funcional de factores de automorfía:
$$j_\rho(g_v g_v', (g_w)(g_w'), z) = j_\rho(g_v, (g_w), g_v'(z)) j_\rho(g_v', (g_w'), z).$$

Vamos a escribir $j_\rho(g, z)$ con $g = (g_v, (g_w)) \in D_\infty^\times$.

Así podemos definir una *forma modular clásica de peso $\rho$* por la fórmula habitual:

**Definición 1-8.** Sea $K \subset D^\times(\mathbf{A}_f)$ un subgrupo abierto compacto. Una forma modular de peso $\rho$ y de nivel $K$ para $D^\times$ es una función holomorfa

$$f : \mathfrak{H}^\pm \times D^\times(\mathbf{A}_f)/K \to W$$

que satisface la ecuación funcional

$$f(d(z), \iota_f(d) \cdot g_f) = j_\rho(\iota_\infty(d)) f(z), z \in \mathfrak{H}^\pm, d \in D^\times.$$

Aquí hemos designado por $\iota_\infty$ (respectivamente $\iota_f$) la inclusión

$$\iota_\infty : D^\times \hookrightarrow D_\infty^\times \text{ (respectivamente } \iota_f : D^\times \hookrightarrow D^\times(\mathbf{A}_f)).$$

Vamos a designar por $M_\rho(D^\times, K)$ al espacio de formas modulares clásicas de peso $\rho$ y de nivel $K$. Si $\rho_v$ es el homomorfismo $x \mapsto x^2$ y si $\rho_w$ es la representación trivial para $w \neq v$, escribimos $M_{(2,2,\ldots,2)}(D^\times, K)$ en vez de $M_\rho(D^\times, K)$. Como en el caso de formas modulares de Hilbert, tenemos la descripción siguiente de $M_{(2,2,\ldots,2)}(D^\times, K)$:

**Proposición 1-9.** *Hay un isomorfismo natural*

$$M_{(2,2,\ldots,2)}(D^\times, K) \xrightarrow{\sim} \Omega^1({}_K S(D))$$

*donde* $\Omega^1({}_K S(D))$ *designa el espacio de formas diferenciales holomorfas sobre* ${}_K S(D)$.

## 2 Formas modulares cuaterniónicas y curvas de Shimura

Antes de presentar las definiciones formales, quisiera explicar que una forma modular sobre el grupo $G = D^\times$ admite una descomposición como combinación lineal de autoformas para operadores de Hecke, exactamente como las formas modulares de Hilbert, con algunas modificaciones naturales. El aspecto más importante de la teoría es que una autoforma sobre el grupo multiplicativo de $D$ contiene la misma información que una autoforma modular de Hilbert. Hay una *correspondencia* entre autoformas modulares sobre $D^\times$ y autoformas modulares de Hilbert. La correspondencia es inyectiva: todas las autoformas (más precisamente, sistemas de autovalores) sobre $D^\times$ tienen realizaciones como autoformas modulares de Hilbert, pero no todas las autoformas de Hilbert surgen de esta manera de $D^\times$. La correspondencia se llama *correspondencia de Jacquet–Langlands* y es una de las motivaciones principales para el estudio de formas modulares sobre $D^\times$. La propiedad principal de esta correspondencia es la preservación de funciones $L$.

Hay pocas referencias sobre las formas modulares cuaterniónicas, pero el libro Godement y Jacquet (1972) es una buena introducción a los principales temas de esta sección.

**2.1 Definiciones.** En esta sección, $D$ es un álgebra de cuaterniones *de división* sobre el cuerpo totalmente real $F$. Escribimos $D^\times(\mathbf{A}) = D_\infty^\times \times D^\times(\mathbf{A}_f)$ como en la sección precedente.

**Definición 2-1.** Una forma modular (forma automorfa) sobre $D^\times$ es una función $f : D^\times(\mathbf{A}) \to \mathbb{C}$ con las propiedades siguientes:

1. Para cualquier $g_f \in D^\times(\mathbf{A}_f)$, $g_\infty \mapsto f(g_\infty, g_f)$ es una función $C^\infty$ sobre $D_\infty^\times$.

2. Para cualquier $g_\infty \in D_\infty^\times$, $g_f \mapsto f(g_\infty, g_f)$ es una función localmente constante sobre $D^\times(\mathbf{A}_f)$.

3. Para cualquier $\gamma \in D^\times$ y $g \in D^\times(\mathbf{A})$, $f(\gamma \cdot g) = f(g)$.

4. Existe un subgrupo $K_f \subset D^\times(\mathbf{A}_f)$, abierto y compacto, tal que $f(gk) = f(g)$ para cualquier $g \in D^\times(\mathbf{A})$ y $k \in K_f$.

5. Sea $C(D^\times(\mathbf{A}))$ el espacio de funciones continuas complejas sobre $D^\times(\mathbf{A})$. El subespacio de $C(D^\times(\mathbf{A}))$ generado por las funciones $g \mapsto f(gk_\infty z)$, con $k_\infty \in K_\infty$, y con $z \in F^\times(\mathbf{A})$ (el centro de $D^\times(\mathbf{A})$) es de dimensión finita.

6. Sea $v \in \Sigma_D$ un primo arquimediano de $F$ con $D_v \simeq M(2,\mathbb{R})$. Sea $C_v \in U(\mathrm{Lie}(D_v))$ el operador de Casimir (ver el curso de A. Pacetti). Entonces el subespacio de $C(D^\times(\mathbf{A}))$ generado por las funciones $\prod_{v\in\Sigma_D} C_v^{k_v} f$, $k_v = 0, 1, 2, 3, \ldots$, es de dimensión finita.

A diferencia de lo que pasa en el caso de formas modulares de Hilbert, el cociente $(D_\infty^\times \times D^\times(\mathbf{A}_f))/F^\times(\mathbf{A})$ es compacto (como ya hemos visto), y no es necesario imponer una condición de crecimiento moderado. Vamos a designar por $\mathfrak{A}(D^\times)$ al espacio de formas modulares sobre $D^\times$; esta notación sigue siendo válida cuando $D = M(2, F)$ es el álgebra de matrices. Si $W$ es un espacio vectorial complejo de dimensión finita, podemos definir una forma modular sobre $D^\times$ con valores en $W$ como una función $f : D^\times(\mathbf{A}) \to W$ tal que, para toda forma lineal $\lambda : W \to \mathbb{C}$, $\lambda \circ f$ es una forma modular en el sentido de la Definición 2-1.

Como en el caso de formas modulares de Hilbert, las formas modulares clásicas para $D^\times$ se identifican con formas modulares adélicas. Si $d = (d_\infty, d_f) \in D^\times(\mathbf{A})$, $f \in M_\rho(D^\times, K_f)$, definimos

$$\begin{aligned} \Phi = \Phi(f) :\ & D^\times\backslash D^\times(\mathbf{A})/K_f \to W; \\ \Phi(d_\infty, d_f) &= j_\rho(d_\infty, i)^{-1} f(d_\infty(i), d_f). \end{aligned} \tag{2-2}$$

**Proposición 2-3.** *El mapeo $f \mapsto \Phi(f)$ define un isomorfismo entre $M_\rho(D^\times, K_f)$ y el espacio $\mathfrak{A}^{hol}(D^\times, K_f, \rho)$ de formas modulares $\Phi : D^\times\backslash D^\times(\mathbf{A})/K_f \to W$ tales que*

*(i) $dr(X^-)\Phi = 0$, donde $r : D_\infty^\times \to Aut(\mathfrak{A}(D^\times))$ es la representación regular a derecha (de multiplicación a derecha) y $dr : \mathrm{Lie}(D_\infty^\times) \to \mathrm{End}(\mathfrak{A}(D^\times))$ es su diferencial.*

*(ii) Para todo $k \in \tilde{K}_\infty$ y todo $d = (d_\infty, d_f) \in D^\times(\mathbf{A})$, $\Phi(d_\infty k, d_f) = \rho(k)^{-1}\Phi(d)$.*

La demostración es idéntica a la del caso de formas modulares de Hilbert.

Para tratar todas las formas automorfas de manera uniforme, es más natural reemplazar el espacio $\mathfrak{A}^{hol}(D^\times, K_f, \rho)$ de formas modulares vectoriales por la imagen de $\mathfrak{A}^{hol}(D^\times, K_f, \rho) \otimes \mathrm{Hom}(W_\rho, \mathbb{C})$ en $\mathfrak{A}(D^\times)$ bajo el mapeo natural

$$\begin{aligned} \mathfrak{A}^{hol}(D^\times, K_f, \rho) \otimes \mathrm{Hom}(W_\rho, \mathbb{C}) \hookrightarrow \mathfrak{A}(D^\times(\mathbf{A}), W_\rho) \otimes \mathrm{Hom}(W_\rho, \mathbb{C}) \\ \to \mathfrak{A}(D^\times), \end{aligned}$$

donde $\mathfrak{A}(D^\times, W_\rho)$ designa el espacio de funciones de $D^\times(\mathbf{A})$ con valores en $W_\rho$ que satisfacen las condiciones de la Definición 2-1, y la última flecha es inducida por contracción del producto tensorial

$$W_\rho \otimes \mathrm{Hom}(W_\rho, \mathbb{C}) \to \mathbb{C}.$$

Así vamos a trabajar únicamente con formas modulares adélicas con valores complejos.

Si $f \in \mathcal{A}(D^\times)$ y si $g \in D^\times(\mathbf{A})$ definimos $r(g)(f) \in \mathcal{A}(D^\times)$ por la fórmula

$$r(g)(f)(d) = f(dg)$$

(traslación a derecha). La acción $g \mapsto r(g) \in Aut(\mathcal{A}(D^\times))$ define una representación de $D^\times(\mathbf{A})$ en el espacio $\mathcal{A}(D^\times)$.

**Definición 2-4.** Una representación automorfa de $D^\times(\mathbf{A})$ es una subrepresentación irreducible de $\mathcal{A}(D^\times)$.

Esta definición sigue siendo válida cuando $D = M(2, F)$. En este caso, $\mathcal{A}(GL(2, F))$ contiene el subespacio $\mathcal{A}^0(GL(2, F))$ de *formas cuspidales* como subrepresentación. Una *representación automorfa cuspidal* de $GL(2, F_\mathbf{A})$ es una subrepresentación irreducible de $\mathcal{A}^0(D^\times)$.

**Comentario 2-5.** Si trabajamos con el grupo $PD = D^\times/F^\times$ y consideramos el subespacio

$$\mathcal{A}(PD) = \mathcal{A}(D^\times) \cap C(D^\times(\mathbf{A})/F^\times(\mathbf{A}))$$

entonces $\mathcal{A}(PD)$ es también una representación de $D^\times(\mathbf{A})$ que además es isomorfa a una suma *numerable* de representaciones irreducibles, si $D$ es un álgebra de división. (En el caso contrario, hay que considerar también los espacios provenientes de series de Eisenstein.)

**Observación 2-6.** *Sea $\Pi$ una representación irreducible de $\mathcal{A}^0(D^\times)$. Entonces hay una factorización $\Pi \xrightarrow{\sim} \otimes'_v \Pi_v$ como producto tensorial restringido sobre los primos de $F$. Aquí $\Pi_v$ es una representación irreducible de $D_v^\times$. Ver por ejemplo Godement y Jacquet (1972, §10).*

**2.2 Operadores de Hecke.** Fijamos un grupo de nivel abierto compacto $K \subset D^\times(\mathbf{A}_f)$. El álgebra de Hecke $\mathcal{H}(K)$ (respectivamente $\mathcal{H}(K)_\mathbb{Z}$) de nivel $K$ es el álgebra $C_c(D^\times(\mathbf{A}_f)//K)$ (respectivamente $C_c(D^\times(\mathbf{A}_f)//K, \mathbb{Z})$) de funciones continuas sobre $D^\times(\mathbf{A}_f)$ con soporte compacto con valores en $\mathbb{C}$ (respectivamente en $\mathbb{Z}$), invariantes bajo multiplicación por ambos lados por elementos de $K$. Es un álgebra para la convolución

$$\phi_1 \star \phi_2(g) = \int_{D^\times(\mathbf{A}_f)} \phi_1(h)\phi_2(gh^{-1})d^\times h,$$

donde tomamos la medida de Haar $d^\times h$ normalizada de tal forma que $\int_K d^\times h = 1$. Entonces la función característica $1_K$ de $K$ es el elemento neutro del álgebra $\mathcal{H}(K)$. Suponemos que $K = \prod_v K_v$ con $K_v$ abierto compacto en $D_v^\times$. El álgebra local $\mathcal{H}(K_v)$ (respectivamente $H_\mathbb{Z}(D_v^\times, K_v)$) es el álgebra $C_c(D_v^\times//K_v)$ (respectivamente $C_c(D_v^\times//K_v, \mathbb{Z})$), definida de la misma manera.

**2.2.1 Álgebra de Hecke no ramificada (esférica).** Si $v$ es un primo no ramificado para $D$ y si $K_v = GL(2, \mathcal{O}_v)$, entonces $\mathcal{H}(K_v)$ es el álgebra de Hecke clásica. Tomamos la medida de Haar $dh_v$ con $\int_{K_v} dh_v = 1$. El álgebra $\mathcal{H}(K_v)$ tiene como generadores las funciones características de las dobles coclases $R_v = K_v \cdot \begin{pmatrix} \varpi_v & 0 \\ 0 & \varpi_v \end{pmatrix} \cdot K_v$, $R_v^{-1}$, y $T_v = K_v \cdot \begin{pmatrix} \varpi_v & 0 \\ 0 & 1 \end{pmatrix} \cdot K_v$, con $\varpi_v$ un uniformizador de $F_v$.

**2.2.2 Acción de operadores de Hecke sobre formas modulares.** Sea $f \in M_\rho(D^\times, K_f)$ una forma modular clásica de peso $\rho$ y de nivel $K_f$. Suponemos que $v$ es un primo *no ramificado* para $K_f$: $K_f \supset K_v = GL(2, \mathcal{O}_v)$. Entonces el álgebra de Hecke esférica $\mathcal{H}(K_v)$ actúa sobre el espacio $\mathcal{A}^{hol}(D^\times, K_f, \rho)$. Si $T \in \mathcal{H}(K_v)$, definimos $T^{\text{class}}(f)$ por la fórmula (cf. Prop. 2-3)

$$\Phi(T^{\text{class}}(f)) = T(\Phi(f)) \tag{2-7}$$

En particular, si la forma automorfa $\Phi(f) \in \mathcal{A}^{hol}(D^\times, K_f, \rho)$ es un vector propio de el álgebra $\mathcal{H}(K_v)$, entonces $f$ es también un vector propio de los operadores de $\mathcal{H}(K_v)$. Los operadores $T^{\text{class}}$ coinciden con los operadores de Hecke clásicos, salvo multiplicación por factores escalares de normalización.

**2.3 Funciones $L$.** Sea $f$ una autoforma modular sobre $D^\times$. Éste no es un grupo conmutativo, pero el método de la tesis de Tate se aplica casi sin cambios en esta situación para permitir la definición de una función $L$ de $f$. Seguimos la versión de este método en el libro *Zeta Functions of Simple Algebras* de Godement y Jacquet (1972).

Más precisamente, el grupo adélico multiplicativo $D^\times(\mathbf{A})$ está contenido en el grupo aditivo $D(\mathbf{A})$. Sea $\varphi : D^\times(\mathbf{A}) \to \mathbb{C}$ una forma modular y sea $\Phi : D(\mathbf{A}) \to \mathbb{C}$ una función con soporte compacto. (Suponemos que $\varphi$ es invariante bajo $F^\times_{\mathbf{A}}$ para simplificar las fórmulas.) Podemos definir una integral zeta:

$$Z(\varphi, \Phi, s) = \int_{D^\times(\mathbf{A})} \varphi(g)\Phi(g)||\nu(g)||^s d^\times g \tag{2-8}$$

que es absolutamente convergente para $Re(s) >> 0$. Suponemos que $\Phi$ tiene una factorización $\Phi = \otimes'_v \Phi_v$ donde $\Phi_v = 1_{\mathcal{O}_{D_v}}$ en casi todo primo no arquimediano $v$. Suponemos también que $\varphi$ tiene una factorización análoga $\varphi = \otimes'_v \varphi_v$. Para interpretar esa condición es necesario utilizar la teoría de representaciones; tenemos que suponer que $\varphi$ es un vector en una subrepresentación irreducible $\Pi$ de $L^2(F^\times_{\mathbf{A}} D^\times \backslash D^\times(\mathbf{A}))$. Como en 2-6 hay una factorización $\Pi \xrightarrow{\sim} \otimes'_v \Pi_v$, donde $\Pi_v$ es una representación irreducible de $D^\times_v$. Entonces para $Re(s) >> 0$ la integral zeta tiene un producto de Euler convergente:

$$Z(\varphi, \Phi, s) = \prod_v Z_v(\varphi_v, \Phi_v, s),$$

donde

$$Z_v(\varphi_v, \Phi_v, s) = \int_{D^\times_v} \varphi_v(g_v)\Phi_v(g_v)||\nu_v(g_v)||^s d^\times g_v.$$

Además, podemos definir la transformada de Fourier $\hat{\Phi}$ de $\Phi$. Sea

$$\psi = \otimes_v \psi_v : F\backslash F_{\mathbf{A}} \to \mathbb{C}^\times$$

un carácter aditivo no trivial, con $\psi_v : F_v \to \mathbb{C}^\times$. Ponemos

$$\hat{\Phi}(x) = \int_{D(\mathbf{A})} \Phi(y)\psi \circ Tr_D(xy)dy; \;\; \hat{\Phi}_v(x) = \int_{D_v} \Phi_v(y)\psi_v \circ Tr_D(xy)dy,$$

donde $Tr_D : D \to F$ es la traza reducida y $dy$ es una medida de Haar autodual. Si ponemos $\varphi^\vee(g) = \varphi(g^{-1})$, $\varphi^\vee_v(g_v) = \varphi(g_v^{-1})$, entonces

**Teorema 2-9.** *Para cada primo $v$ se pueden definir factores locales $L_v(\Pi_v, s)$ y $\varepsilon(\Pi_v, \psi_v, s)$ como en la tesis de Tate con las propiedades siguientes:*

1. *Para toda $\Phi_v$ y toda $\varphi_v \in \Pi_v$, el cociente*

$$\Xi(\varphi_v, \Phi_v, s) = \frac{Z_v(\varphi_v, \Phi_v, s + \frac{1}{2})}{L_v(\Pi_v, s)}$$

   *es una función entera de $s$.*

2. *Si $v$ es un primo no arquimediano, $q = Nv$, entonces o bien $L_v(\Pi_v, s) = 1$, o bien $L_v(\Pi_v, s)$ es producto de uno o dos factores de Euler de la tesis de Tate, y la función $\Xi(\varphi_v, \Phi_v, s)$ es un polinomio en $q^s$ y $q^{-s}$.*

3. *Si $v$ es un primo arquimediano entonces $L_v(\Pi_v, s)$ es un producto de potencias de $\pi$ y de factores $\Gamma$.*

4. *El factor $\varepsilon(\Pi_v, \psi_v, s)$ es una función entera de $s$, y hay una ecuación funcional local:*

$$\Xi(\varphi_v^\vee, \hat{\Phi}_v, 1 - s) = (-1)^{e_v} \varepsilon(\Pi_v, \psi_v, s) Z_v(\varphi_v, \Phi_v, s)$$

   *donde $e_v = 0$ si $D_v = M(2, F_v)$ y $e_v = 1$ si $D_v$ es un álgebra de división.*

Exactamente como en la tesis de Tate, hay también una ecuación funcional global:

**Teorema 2-10.** *Sea $\varphi$ y $\Phi$, y $Z(\varphi, \Phi, s)$ como en (2-8). Entonces*

$$Z(\varphi, \Phi, s) = Z(\varphi^\vee, \hat{\Phi}, 2 - s).$$

Finalmente, definimos la función $L$ de $\Pi$ como el producto de Euler

(2-11) $$L(\Pi, s) = \prod_v L_v(\Pi_v, s).$$

El producto es convergente para $Re(s) >> 1$ (de hecho, para $Re(s) > 1$ con nuestras hipótesis). Definimos el factor épsilon del mismo modo:

(2-12) $$\varepsilon(\Pi, s) = \prod_v \varepsilon_v(\Pi_v, \psi_v, s).$$

El producto es independiente del carácter $\psi$ escogido. Como en la tesis de Tate, deducimos formalmente la ecuación funcional de la función $L$.

**Corolario 2-13.**

$$L(\Pi, s) = \varepsilon(\Pi, s) L(\Pi^\vee, 1 - s),$$

*donde $\Pi^\vee$ es la representación automorfa contragradiente de $\Pi$.*

**2.4 Teoría local de representaciones, caso no arquimediano.** En este párrafo $k$ es un cuerpo $p$-ádico, con anillo de enteros $\mathcal{O}$ y valor absoluto $|\bullet|_k$. Los dos lemas siguientes son básicos.

**Lema 2-14** (lema de Schur)**.** *Sea $\pi$ una representación lisa e irreducible de $GL(2,k)$. Sea $Z = k^\times \subset GL(2,k)$ el centro de $GL(2,k)$. Entonces existe un homomorfismo $\xi_\pi : Z \to \mathbb{C}^\times$, el* carácter central *de $\pi$, tal que, para todo $v \in \pi$ y todo $z \in Z$,*

$$\pi(z)v = \xi_\pi(z)v.$$

En el lema de Schur, la hipótesis que $\pi$ sea lisa es esencial.

**Ejercicio 2.1.** Demostrar el lema de Schur: Sea $U \subset GL(2,k)$ un subgrupo abierto tal que el subespacio

$$\pi^U = \{v \in \pi \mid \pi(u)v = v \ \forall\, u \in U\}$$

es de dimensión positiva. Sea $z \in Z$. Mostrar que el grupo $Z$ estabiliza el subespacio $\pi^U$ y que existe un vector $v \in \pi^U$, $v \neq 0$, que es vector propio de todos los elementos de $Z$. Entonces utilizar la irreducibilidad de $\pi$ para terminar la demostración.

**Lema 2-15.** *Sea $\pi$ una representación lisa e irreducible de dimensión finita de $GL(2,k)$. Entonces* $\dim \pi = 1$ *y existe un carácter continuo (localmente constante) $\chi : k^\times \to \mathbb{C}^\times$ tal que $\pi = \chi \circ$*det *(el carácter se factoriza por medio del determinante).*

**Ejercicio 2.2.** (a) Demostrar que, si $\pi$ es lisa y de dimensión finita, entonces hay un subgrupo abierto $U \subset GL(2,k)$ tal que $\pi^U = \pi$. En particular, existe $\varepsilon > 0$ tal que, si $a \in \mathcal{O}$, $|a|_k < \varepsilon$, entonces, para todo $v \in \pi$,

$$\pi\left(\begin{pmatrix} 1 & a \\ 0 & 1 \end{pmatrix}\right) v = \pi\left(\begin{pmatrix} 1 & 0 \\ a & 1 \end{pmatrix}\right) v = v.$$

(b) Demostrar el Lema 2-15.

Ahora pasamos a las tres clases de representaciones lisas de dimensión infinita.[1]

**2.4.1 Serie principal.** Sea $(\chi_1, \chi_2)$ un par ordenado de caracteres de $k^\times$. Sean $G = GL(2,k)$ y $B \subset G$ el subgrupo de Borel triangular superior, que se puede descomponer como $B = A \cdot N$, con

$$A = \left\{\begin{pmatrix} a_1 & 0 \\ 0 & a_2 \end{pmatrix}\right\}, \quad N = \left\{\begin{pmatrix} 1 & x \\ 0 & 1 \end{pmatrix}\right\}.$$

Definimos

$$I(\chi_1,\chi_2) = \\ \{f : G \to \mathbb{C} \mid f(ang) = \chi_1(a_1)|a_1|^{\frac{1}{2}}\chi_2(a_2)|a_2|^{-\frac{1}{2}} \cdot f(g),\ a \in A,\ n \in N\}$$

[1] Las diapositivas de un seminario de Pilar Bayer, http://www.icmat.es/seminarios/langlands/14.01.10/bayer.pdf, contienen una buena introducción a este material.

(todas las funciones son supuestas localmente constantes y por lo tanto continuas). Esta es una representación inducida normalizada, y $G$ actúa sobre $I(\chi_1, \chi_2)$ por traslación a derecha:

$$r(g)f(h) = f(hg).$$

Las potencias $\frac{1}{2}$ de la norma garantizan que casi siempre tenemos un isomorfismo

$$I(\chi_1, \chi_2) \xrightarrow{\sim} I(\chi_2, \chi_1)$$

Más precisamente, tenemos la siguiente proposición:

**Proposición 2-16.** *(a) $I(\chi_1, \chi_2)$ es irreducible a menos que $\chi_1/\chi_2 = |\bullet|^{\pm 1}$.*

*(b) Si $\chi_1/\chi_2 \neq |\bullet|^{\pm 1}$ entonces*

$$I(\chi_1, \chi_2) \xrightarrow{\sim} I(\chi_2, \chi_1)$$

*como representaciones irreducibles admisibles de $GL(2, F)$.*

*(c) Si $\chi_1/\chi_2 \neq |\bullet|^{\pm 1}$ entonces la representación contragradiente de $I(\chi_1, \chi_2)^\vee$ de $I(\chi_1, \chi_2)$ es isomorfa a la representación $I(\chi_1^{-1}, \chi_2^{-1})$.*

‘Las representaciones $I(\chi_1, \chi_2)$ se llaman *representaciones de la serie principal.*

Sea $K = GL(2, \mathcal{O}) \subset G$, donde $\mathcal{O}$ es el anillo de enteros de $k$. Sea $\varpi \in \mathcal{O}$ un uniformizador, $q = |\mathcal{O}/\varpi\mathcal{O}|$. Suponemos $\chi_1$ y $\chi_2$ *no ramificados*. Entonces $I(\chi_1, \chi_2)$ contiene un vector $K$-invariante canónico $f_0$ definido por

$$f_0(k) = 1 \ \forall\, k \in K.$$

La descomposición de Iwasawa $G = B{\cdot}K$ implica que todos los valores de $f_0$ quedan determinados por este propiedad. Se puede definir un álgebra local de Hecke $H_{\mathbb{Z}}(G, K)$, la subálgebra del álgebra compleja $\mathcal{H}(K)$ definida antes, de funciones con valores en $\mathbb{Z}$; entonces para cualquier anillo $A$ definimos $H_A(G, K) = H_{\mathbb{Z}}(G, K) \otimes_{\mathbb{Z}} A$. El álgebra $H_A(G, K)$ es conmutativa, y tiene como en el caso complejo dos generadores

$$T = K \cdot \begin{pmatrix} \varpi & 0 \\ 0 & 1 \end{pmatrix} \cdot K, \quad R = K \cdot \begin{pmatrix} \varpi & 0 \\ 0 & \varpi \end{pmatrix} \cdot K;$$

$$H_A(G, K) = A[T, R, R^{-1}].$$

El álgebra de convolución $H(G)$ de todas la funciones localmente constantes y con soporte compacto opera sobre cualquier representación lisa de $G$, y la función $f_0$ es un autovector para su subálgebra $H_{\mathbb{C}}(G, K)$, con

$$Tf_0 = q^{\frac{1}{2}}(\chi_1(\varpi) + \chi_2(\varpi))f_0, \quad Rf_0 = \chi_1(\varpi)\chi_2(\varpi)f_0. \tag{2-17}$$

Una representación irreducible con un vector $v_0$ fijo por $K$ está determinada, salvo isomorfismo, por sus autovalores $t_0$ y $r_0$ con $T(v_0) = t_0 v_0$, $R(v_0) = r_0 v_0$. Una tal representación se llama *esférica.*

Si $k = \mathbb{Q}_p$ entonces $T$ es el operador clásico $T(p)$, y $R = T(p, p)$ (salvo normalización). La teoría clásica de los operadores de Hecke queda completamente sustituida por la teoría de representaciones esféricas.

### 2.4.2 Representaciones de Steinberg.

Sea $\chi : k^\times \to \mathbb{C}^\times$ un carácter liso. Sean

$$\chi_1 = \chi \cdot |\bullet|^{-\frac{1}{2}}, \chi_2 = \chi \cdot |\bullet|^{\frac{1}{2}}.$$

En ese caso es fácil ver que la función $f_\chi(g) = \chi(det(g))$ pertenece a $I(\chi_1, \chi_2)$. Hay una sucesión exacta corta de representaciones admisibles

$$0 \to \mathbb{C} f_\chi \to I(\chi_1, \chi_2) \to St(\chi) \to 0$$

donde $G$ actúa por el carácter $\chi \circ \det$ sobre $\mathbb{C} f_\chi$ y $St(\chi)$ es irreducible; las $St(\chi)$ son las representaciones de Steinberg. Cuando $\chi$ es el carácter trivial, llamamos a $St(1)$ *la* representación de Steinberg.

La representación $I(\chi_2, \chi_1)$ es también reducible y tiene a $St(\chi)$ como subrepresentación, y a $\chi \circ \det$ como cociente.

### 2.4.3 Representaciones supercuspidales.

**Definición 2-18.** Sea $\pi$ una representación irreducible de $GL(2,k)$. Decimos que $\pi$ es *supercuspidal* si $Hom_{GL(2,F)}(\pi, I(\chi_1, \chi_2)) = 0$ para cualquier par $(\chi_1, \chi_2)$ de caracteres.

Las representaciones de Steinberg y las representaciones de la serie principal no son supercuspidales, todos las demás lo son.

No hay una construcción elemental de representaciones supercuspidales. La clasificación de estas representaciones más importante es la *correspondencia de Langlands*. Para explicar eso, necesitamos algunas definiciones suplementarias. Sea $\mathbb{F} = \mathcal{O}/\varpi\mathcal{O}$ el cuerpo residual de $k$, y $Frob : \bar{\mathbb{F}} \to \bar{\mathbb{F}}$ el automorfismo de Frobenius: $Frob(x) = x^q$. Sea $W_k$ el grupo de Weil de $F$; lo podemos definir como el subgrupo de $\gamma \in Gal(\bar{k}/k)$ que actúan sobre $\bar{\mathbb{F}}$ como potencias enteras de $Frob$.

**Teorema 2-19** (Tunnell, Kutzko). *Hay una biyección entre las (clases de equivalencia de) representaciones supercuspidales de $GL(2,k)$ y (clases de equivalencia de) representaciones irreducibles de dimensión 2 de $W_k$.*

La correspondencia de Langlands vale también para las representaciones supercuspidales de $GL(n,k)$ (Teorema de Harris—Taylor y Henniart), pero la clasificación de representaciones no supercuspidales es técnicamente más complicada.

### 2.4.4 Representaciones discretas y temperadas.

**Definición 2-20.** Sea $\pi$ una representación irreducible de $GL(2,k)$. Decimos que $\pi$ es *discreta* si $\pi$ es isomorfa, o bien a una representación de Steinberg, o bien a una representación supercuspidal. Decimos que $\pi$ es *esencialmente temperada* si $\pi$ es isomorfa, sea a una representación discreta, sea a una representación $I(\chi_1, \chi_2)$ de la serie principal con $|\chi_1(x)| = |\chi_2(x)|$ para $x \in k^\times$. Decimos que $\pi$ es *temperada* si $\pi$ es esencialmente temperada y con carácter central unitario.

**Comentario 2-21.** Las representaciones temperadas son las que contribuyen a la fórmula de Plancherel en análisis armónico sobre el grupo (descomposición del espacio $L^2(GL(2,k))$ como integral directa de representaciones irreducibles). Las representaciones discretas con carácter central unitario son las que contribuyen discretamente a la ' de Plancherel (con medida puntual positiva).

**2.4.5 Representaciones del grupo multiplicativo de un álgebra de división.** Ahora suponemos $G = D^\times$ con $D$ un álgebra cuaterniónica de división sobre $k$. Sea $\pi$ una representación lisa irreducible de $G$. Exactamente como en el caso de $GL(2,k)$, el centro $Z_G \simeq k^\times$ actúa sobre $\pi$ por un carácter $\xi_\pi$. Como el cociente $G/Z_G$ es un grupo compacto, eso implica que $\pi$ es necesariamente de dimensión *finita*.

**2.5 Teoría local de representaciones, caso arquimediano.** En este párrafo vamos a estudiar algunas representaciones irreducibles de los grupos $GL(2,\mathbb{R})$ y $\mathbb{H}^\times$. Comenzamos por el grupo $\mathbb{H}^\times$. El lema de Schur sigue siendo válido, y como el cociente $\mathbb{H}^\times/\mathbb{R}^\times$ es compacto y se identifica con el grupo $SO(3)$ de rotaciones de $\mathbb{R}^3$, toda representación irreducible y continua de $\mathbb{H}^\times$ es de dimensión finita.

**Proposición 2-22.** *Para todo entero impar $t = 1, 3, 5, \ldots$, existe una única clase de equivalencia de representaciones irreducibles $\pi_t^{\mathbb{H}}$ de $\mathbb{H}^\times/\mathbb{R}^\times$ de dimensión $t$.*

*Prueba.* Sea $G$ el grupo recubridor universal de $\mathbb{H}^\times/\mathbb{R}^\times \simeq SO(3)$; entonces $G$ es isomorfo al grupo $SU(2)$: hay una sucesión exacta corta

$$1 \to \{\pm 1\} \to SU(2) \to SO(3) \to 1.$$

Las representaciones irreducibles de un grupo de Lie $G$ compacto y simplemente conexo son clasificadas por las representaciones irreducibles de su álgebra de Lie (complexificada) $\mathfrak{g}$. Cuando $G = SU(2)$, el álgebra de Lie $\mathfrak{g}$ es isomorfa al álgebra de Lie

$$\mathfrak{sl}(2) = \{X \in M(2,\mathbb{C}) \mid Tr(X) = 0\}.$$

Un teorema básico de la teoría de álgebras de Lie afirma que las representaciones irreducibles de $\mathfrak{sl}(2)$ son clasificadas por enteros positivos $1, 2, 3, \ldots$; la representación trivial es la única de dimensión 1, la representación natural $\mathfrak{sl}(2) \to \mathrm{End}(\mathbb{C}^2)$ es la única de dimensión 2, y para cualquier entero $t \geqslant 2$, la representación

$$\mathfrak{sl}(2) \to \mathrm{End}(\mathrm{Sym}^{t-1}(\mathbb{C}^2))$$

es de dimensión $t$ (es la representación natural sobre los polinomios homogéneos de grado $t-1$ en dos variables).

Si ponemos $\mathbb{C} = \mathrm{Sym}^0(\mathbb{C}^2)$, con la representación trivial, entonces el grupo $SU(2)$ actúa también sobre $\mathrm{Sym}^{t-1}(\mathbb{C}^2)$ para cualquier $t > 0$. El centro $C$ de $SU(2)$ es el grupo $\left\{\pm\begin{pmatrix}1 & 0\\ 0 & 1\end{pmatrix}\right\}$; la acción de $C$ sobre $\mathrm{Sym}^{t-1}(\mathbb{C}^2)$ es trivial si y solamente si $t$ es impar. En consecuencia, las representaciones irreducibles de $\mathbb{H}^\times/\mathbb{R}^\times \simeq SU(2)/C$ son las $\mathrm{Sym}^{t-1}(\mathbb{C}^2)$ con $t$ impar. □

La teoría de representaciones de $G = GL(2,\mathbb{R})$ es más complicada. En realidad, en la teoría de formas automorfas, es más natural hablar de representaciones del álgebra de Lie (complexificada) $\mathfrak{g} = M(2,\mathbb{C})$ que tienen una estructura de $(\mathfrak{g}, K)$-*módulo* en el sentido de Harish-Chandra. Sea $K \subset G$ un grupo compacto maximal; cuando $G = GL(2,\mathbb{R})$, ponemos $K = O(2)$, el grupo ortogonal con respecto al producto escalar euclidiano.

**Definición 2-23.** Un $(\mathfrak{g}, K)$-módulo es un espacio vectorial $V$ (sobre $\mathbb{C}$) que admite una acción lineal $\pi$ de $\mathfrak{g}$ y una acción continua $\rho$ de $K$ que son compatibles:

1. La acción de $\mathfrak{g}$ sobre $V$, restringida al álgebra de Lie de $K$, coincide con el diferencial de $\rho$;
2. Si $k \in K$, $X \in \mathfrak{g}$, entonces $\rho(k)\pi(X)\rho(k)^{-1} = \pi(ad(k)X)$.

De un $(\mathfrak{g}, K)$-módulo obtenemos automáticamente una representación del *álgebra universal envolvente* de $\mathfrak{g}$. Sin embargo, si

$$\pi : G \rightarrow Aut(V)$$

es una representación continua sobre un espacio vectorial topológico razonable, se puede definir el subespacio $V^\infty \subset V$ de vectores diferenciables; entonces $V^\infty$ es denso en $V$ (teorema de Gårding) y admite una estructura canónica de $(\mathfrak{g}, K)$-módulo. Eso se aplica, por ejemplo, a las representaciones de la serie principal de $GL(2, \mathbb{R})$, que son construidas exactamente como en el caso no arquimediano.

Para las aplicaciones a las formas modulares holomorfas, necesitamos solo la clase de *representaciones holomorfas* de $GL(2, \mathbb{R})$, o más bien los $(\mathfrak{g}, K)$-módulos holomorfos. Nos interesan solo las representaciones con carácter central trivial (es decir, las representaciones de $PGL(2, \mathbb{R}) = GL(2, \mathbb{R})/\mathbb{R}^\times$). Como en el caso de $\mathbb{H}^\times/\mathbb{R}^\times$, la clasificación de $(\mathfrak{g}, K)$-módulos con carácter central trivial se reduce a la clasificación de módulos sobre $\mathfrak{sl}(2)$. Como grupo compacto maximal de $SL(2, \mathbb{R})$ tomamos el grupo $SO(2) = O(2) \cap SL(2, \mathbb{R})$.

**Ejercicio 2.3.** Determinar el normalizador de $SO(2)$ en $SL(2, \mathbb{R})$.

**Proposición 2-24.** *Sea $t \geqslant 0$ un entero. Existe un único $(\mathfrak{sl}(2), SO(2))$-módulo $\pi'_t$ irreducible generado por un vector $v_t$ caracterizado por las dos propiedades siguientes:*

*(i) $X^- v_t = 0$, con $X^- = \begin{pmatrix} 1 & -i \\ -i & -1 \end{pmatrix}$;*

*(ii) $v_t$ es un vector propio de $SO(2)$ con carácter*

$$\begin{pmatrix} a & b \\ -b & a \end{pmatrix} v_t = \alpha_{t+1}(a + ib) v_t,$$

*donde $\alpha_t(e^{i\theta}) = e^{it\theta}$.*

**Ejercicio 2.4.** Sea $v_t \in \pi'_t$ el vector de la proposición. Sea $X^+ = \begin{pmatrix} 1 & i \\ i & -1 \end{pmatrix} \in \mathfrak{sl}(2)$. Demostrar que, para todo entero $a \geqslant 0$, el vector $(X^+)^a v_t$ es un vector propio de $SO(2)$, y calcular su valor propio.

**Definición 2-25.** Una representación irreducible $\pi$ de $GL(2, \mathbb{R})$ (es decir, un $(\mathfrak{gl}(2), O(2))$-módulo)) es de la *serie discreta* si contiene al $(\mathfrak{sl}(2), SO(2))$-módulo $\pi'_t$ con $t > 0$. La representación irreducible $\pi$ es un *límite de la serie discreta* si contiene al $(\mathfrak{sl}(2), SO(2))$-módulo $\pi'_0$.

La proposición siguiente es una aplicación fácil de la teoría de representaciones inducidas:

**Proposición 2-26.** *Sea $\tau_1$ y $\tau_2$ dos representaciones irreducibles de $GL(2, \mathbb{R})$ que contienen la misma representación $\pi'_t$ de $SL(2, \mathbb{R})$. Entonces $\tau_1$ y $\tau_2$ son equivalentes.*

Así tenemos el derecho de hablar de *la* representación irreducible $\pi_t$ de $GL(2, \mathbb{R})$ que contiene $\pi'_t$.

**Comentario 2-27.** Sea $F = \mathbb{Q}$ y sea $f$ una forma modular clásica de peso $k > 0$ para un subgrupo de $SL(2, \mathbb{Z})$ de índice finito. Sea $\pi$ la representación automorfa de $GL(2, \mathbf{A})$ que corresponde a $f$. Entonces la componente $\pi_\infty$ de $\pi$ correspondiente al primo arquimediano es isomorfa a la representación $\pi_{k-1}$.

Más generalmente, si $f$ es una forma modular de Hilbert de peso $(2, 2, \ldots, 2)$, con representación automorfa asociada $\pi$, entonces la componente $\pi_v$ en cualquier primo arquimediano es isomorfa a la representación $\pi_1$.

**2.6 Correspondencia de Shimizu y de Jacquet–Langlands.** Para describir la correspondencia de Jacquet–Langlands entre formas modulares (automorfas) de Hilbert y formas modulares (automorfas) sobre $D^\times(\mathbf{A})$, lo más sencillo es de comenzar con la *correspondencia local de Jacquet–Langlands* $JL$ entre representaciones de los grupos locales $D_v^\times$ y de $GL(2, F_v)$. Con las definiciones ya introducidas eso se hace bastante rápidamente y muy conceptualmente. Sin embargo, las demostraciones de la correspondencia local están basadas en la correspondencia global.

**2.6.1 Correspondencias locales.** Sea $v$ un primo de $F$. Ponemos $G = GL(2, F_v)$, $J = H_v^\times$, donde $H_v$ es un álgebra cuaterniónica de división de dimensión $4$ sobre $F_v$. Las representaciones irreducibles de $J$ son todas de dimensión finita. Designamos por $Rep(J)$ (respectivamente $Rep(G)$) el conjunto de (clases de equivalencia de) representaciones irreducibles de $J$ (respectivamente $G$).

**Teorema 2-28** (Correspondencia (local) de Jacquet–Langlands)**.** *Existe una biyección canónica*

$$JL : Rep(J) \simeq Rep_{disc}(G) \subset Rep(G)$$

*entre el conjunto $Rep(J)$ de representaciones irreducibles de $J$ y el conjunto $Rep_{disc}(G)$ de representaciones (irreducibles) discretas de $G$. Esta biyección conserva los factores locales de la ecuación funcional: si $\Pi \in Rep(J)$, y $\psi : F_v \to \mathbb{C}^\times$ es un carácter aditivo no trivial, entonces*

$$L(\Pi, s) = L(JL(\Pi), s) \quad \varepsilon(\Pi, \psi, s) = \varepsilon(JL(\Pi), \psi, s)$$

*donde $L(\Pi, s)$ y $\varepsilon(\Pi, \psi, s)$ son los factores locales de Godement–Jacquet introducidos en el Teorema 2-9.*

*Las representaciones de $J$ de dimensión $1$ son las que se corresponden con las representaciones de Steinberg de $G$.*

Las representaciones no discretas de $GL(2, F_v)$, en particular las representaciones de la serie principal, no tienen correspondientes en $Rep(J)$.

Le demostración depende de la fórmula de trazas y necesita la introducción de la teoría de caracteres de representaciones de grupos reductivos sobre cuerpos locales. El carácter de una representación irreducible $\pi$ de $J$ es nada más que la traza habitual, que tiene un sentido porque $\pi$ es de dimensión finita. En cambio, el carácter de una representación de dimensión infinita es una distribución, y su existencia y propiedades hacen parte de la teoría de Harish-Chandra. Con las buenas definiciones, se puede caracterizar la correspondencia $JL$ por un relación explícita de caracteres. Sin embargo, la conservación de factores locales es suficiente para caracterizar la correspondencia.

En general, la determinación explícita de la correspondencia local de Jacquet–Langlands es difícil, porque no hay una descripción elemental de las representaciones supercuspidales. En ciertos casos hay una descripción directa de la correspondencia. Nos limitamos a las representaciones con carácter central trivial.

**Proposición 2-29.** *Sea $v$ un primo arquimediano. Entonces para todo entero impar $t > 0$, $JL(\pi_t^{\mathbb{H}}) = \pi_t$.*

Por supuesto, no se puede separar la demostración de esta proposición de la demostración de la correspondencia en general.

**Proposición 2-30.** *Sea $v$ un primo no arquimediano, y sea $\chi : F_v^\times \to \mathbb{C}^\times$ un carácter liso. Sea $\pi(\chi) = \chi \circ \nu : J \to \mathbb{C}^\times$ una representación de dimensión 1 de $J$. Entonces $JL(\pi(\chi)) = St(\chi)$.*

**2.6.2 La correspondencia global.** Ahora sea $D$ un álgebra de división global sobre $F$. Sea $S'_D$ el conjunto de primos de $F$ (arquimedianos o no) de ramificación para $D$: $D_v$ es un álgebra de división sobre $F_v$ si y solamente si $v \in S'_D$. Entonces $\Sigma'_D$ es el subconjunto de primos arquimedianos en $S'_D$. La cardinalidad $|S'_D|$ de $S'_D$ es un número *par*; como hemos supuesto que $D$ es un álgebra de división, $|S'_D| \geqslant 2$.

La correspondencia global de Jacquet–Langlands es una biyección del conjunto $Aut(D^\times)$ de representaciones automorfas de $D^\times$ y un subconjunto del conjunto $Aut(GL(2, F))$ de representaciones automorfas de $GL(2, F_{\mathbf{A}})$.

**Teorema 2-31.** *Sea $\pi$ una representación automorfa de $D^\times(\mathbf{A})/F_{\mathbf{A}}^\times$. Sea $\Pi = JL(\pi)$ la representación admisible de $GL(2, F_{\mathbf{A}})$ definida como producto tensorial restringido $\otimes'_v \Pi_v$, donde*

1. *Si $v \notin S'_D$, $D_v^\times \simeq GL(2, F_v)$ y $\Pi_v$ es equivalente a $\pi_v$;*
2. *Si $v \in S'_D$, $\Pi_v$ es equivalente a $JL(\pi_v)$.*

*Entonces $\Pi$ es una representación automorfa de $GL(2, F_{\mathbf{A}})$; además, $\Pi$ esta contenida en el espacio de formas cuspidales.*

*El mapeo $\pi \mapsto JL(\pi)$ define una biyección entre $Aut(D^\times)$ y el subconjunto de $Aut(GL(2, F))$ de representaciones automorfas $\Pi$ de $GL(2, F_{\mathbf{A}})$ que satisfacen la condición de que $\Pi_v$ esta en la serie discreta para todo $v \in S'_D$.*

La demostración de este teorema por Jacquet y Langlands está basada en la fórmula de trazas. Se puede leer un bosquejo de la demostración en la sección 2.2 de Harris (2011). Ese bosquejo presupone la existencia de la correspondencia local, pero en la práctica, las correspondencias local y global son demostradas simultáneamente.

**2.7 Formas modulares de Hilbert de peso $(2, 2, \ldots, 2)$.** El Teorema 2-31 tiene una interpretación más elemental en términos de formas modulares holomorfas. Sea $_K S(D)$ la curva de Shimura asociada al álgebra de división $D$ con $\Sigma_D = \{v\}$. Sea $(\rho, W_\rho)$ una representación irreducible de $\tilde{K}_\infty$ trivial en el centro de $D_\infty^\times$. Si $w \neq v$, designamos por $t_w$ la dimensión de $W_{\rho_w}$, donde $\rho_w : D_w^\times/F_w^\times = \mathbb{H}^\times/\mathbb{R}^\times \to GL(W_{\rho_w})$ es la componente local de $\rho$ en $w$. Sea $\rho_v$ el carácter $\alpha^{-t_v}$ de $SO(2) \subset \tilde{K}_v$.

Si $u$ es un primo no arquimediano, $u \notin S'_D$, entonces $D_u^\times \xrightarrow{\sim} GL(2, F_u)$. Así podemos identificar las álgebras de Hecke esféricas de $D_u^\times$ y de $GL(2, F_u)$, relativamente a sus subgrupos compactos maximales respectivos.

**Corolario 2-32.** *Sea $S$ un conjunto finito de primos no arquimedianos que contiene todos los primos no arquimedianos en $S'_D$. Sea $K_f \subset D^\times(\mathbf{A}_f)$ un subgrupo compacto abierto que contiene $GL(2, \mathcal{O}_u)$ para todo primo no arquimediano $u \notin S$. Sea $f : \mathfrak{H}^\pm \times D^\times(\mathbf{A}_f)/K_f \to W$ una forma modular clásica que es autoforma para los operadores de Hecke en los primos fuera de $S$. Entonces existe una forma modular de Hilbert $f^{Hilb}$ de peso $t_v$ en $v$ y de peso $t_w + 1$ en el primo arquimediano $w \neq v$, que es autoforma para los operadores de Hecke en los primos fuera de $S$ con los mismos valores propios que $f$.*

Para aplicaciones a las curvas elípticas, nos bastará enfocar nuestra atención en las formas modulares sobre curvas de Shimura que corresponden a las formas modulares de Hilbert de peso (2,2,...,2), es decir, con $t_v = 2$ y con $t_w = 0$ para todo $w \neq v$.

**2.7.1 Formas nuevas.** La teoría de *formas nuevas*, introducida por Atkin y Lehner y generalizada por Miyake y otros, permite definir un espacio natural de dimensión uno en una representación automorfa de $GL(2, F_\mathbf{A})$. Esta teoría se generaliza fácilmente al caso de $D^\times(\mathbf{A})$ si nos limitamos a las representaciones automorfas asociadas a formas de peso (2,2,...,2) cuyas componentes locales en los primos de $S'_D$ son de dimensión uno.

Introduciremos los grupos de nivel más importantes. Si $v \in S'_D$, sea $\mathcal{O}_{D_v}$ el orden maximal de $D_v$.

(2-33) $$K(D) = \prod_{v \notin S'_D} GL(2, \mathcal{O}_v) \times \prod_{v \in S'_D} \mathcal{O}_{D_v}^\times \subset D^\times(\mathbf{A}_f)$$

(producto sobre primos no arquimedianos) es un grupo de nivel maximal. Sea $\mathfrak{n} \subset \mathcal{O}_F$ un ideal relativamente primo al conjunto $S'_D$, y definimos

(2-34) $$K_0(\mathfrak{n}, D) = \{k = (k_v) \in K(D) \mid k_v \equiv \begin{pmatrix} * & * \\ 0 & * \end{pmatrix} \pmod{\mathfrak{n}}\ \forall v \notin S'_D\};$$

(2-35) $$K_1(\mathfrak{n}, D) = \{k = (k_v) \in K(D) \mid k_v \equiv \begin{pmatrix} * & * \\ 0 & 1 \end{pmatrix} \pmod{\mathfrak{n}}\ \forall v \notin S'_D\}.$$

**Teorema 2-36.** *Sea $S$ un conjunto finito de primos no arquimedianos que contiene todos los primos no arquimedianos en $S'_D$. Sea $K_f \subset D^\times(\mathbf{A}_f)$ un subgrupo compacto abierto que contiene $GL(2, \mathcal{O}_u)$ para todo primo no arquimediano $u \notin S$. Sea $f : \mathfrak{H}^\pm \times D^\times(\mathbf{A}_f)/K_f \to \mathbb{C}$ una forma modular clásica de peso (2,2,...,2) que es autoforma para los operadores de Hecke en los primos fuera de $S$. Sea $\Pi$ la representación automorfa de $D^\times(\mathbf{A})$ asociada a la forma modular clásica $f$. Suponemos que, en todo primo no arquimediano $v \in S$, la componente $\Pi_v$ de $\Pi$ es de dimensión uno.*

*Entonces existe un ideal* $\mathfrak{n} \subset \mathcal{O}_F$ *relativamente primo al conjunto* $S'_D$*, el* conductor *de* $\Pi$ *(fuera de* $S'_D$*), tal que*

1. $\dim \Pi^{K_1(\mathfrak{n},D)} = 1$*;*
2. *Si* $\mathfrak{n}' \supsetneq \mathfrak{n}$ *es un ideal de* $\mathcal{O}_F$*, entonces* $\Pi^{K_1(\mathfrak{n}',D)} = \{0\}$.

Cuando $\Pi$ es una representación automorfa de $D^\times(\mathbf{A})$ como en el Teorema 2-36 con conductor $\mathfrak{n}$ (fuera de $S'_D$), llamaremos a los elementos de $\Pi^{K_1(\mathfrak{n},D)}$ *vectores nuevos*. Las formas modulares clásicas que corresponden a los vectores nuevos en $\Pi$ se llaman *formas nuevas* (de nivel $K_1(\mathfrak{n}, D)$). Una forma nueva de nivel $K_1(\mathfrak{n}, D)$ que es invariante por la acción del grupo conmutativo $K_0(\mathfrak{n}, D)/K_1(\mathfrak{n}, D)$ se llama una *forma nueva* de nivel $K_0(\mathfrak{n}, D)$. Sea

$$M_\rho(D^\times, K_1(\mathfrak{n}, D))^{\text{nuevo}} \subset M_\rho(D^\times, K_1(\mathfrak{n}, D))$$

(respectivamente

$$\mathfrak{A}^{hol}(D^\times, K_1(\mathfrak{n}, D), \rho)^{\text{nuevo}} \subset \mathfrak{A}^{hol}(D^\times, K_1(\mathfrak{n}, D), \rho))$$

los subespacios de formas nuevas (respectivamente de vectores nuevos), y definimos $M_\rho(D^\times, K_0(\mathfrak{n}, D))^{\text{nuevo}}$ y $\mathfrak{A}^{hol}(D^\times, K_0(\mathfrak{n}, D), \rho)^{\text{nuevo}}$ de la misma manera.

El subespacio de $M_\rho(D^\times, K)$ ortogonal al subespacio $M_\rho(D^\times, K)^{\text{nuevo}}$, con respecto al producto escalar de Petersson, se llama el subespacio de *formas viejas*:

$$M_\rho(D^\times, K) = M_\rho(D^\times, K)^{\text{nuevo}} \oplus M_\rho(D^\times, K)^{\text{viejo}}. \tag{2-37}$$

En el caso de formas de peso (2,2,...,2), identificamos $M_{(2,2,\ldots,2)}(D^\times, K) = \Omega^1({}_KS(D))$; entonces (2-37) se escribe

$$\Omega^1({}_KS(D)) = \Omega^1({}_KS(D))^{\text{nuevo}} \oplus \Omega^1({}_KS(D))^{\text{viejo}}. \tag{2-38}$$

La principal propiedad del espacio de formas nuevas está expresada por el siguiente teorema (teorema de multiplicidad uno).

**Teorema 2-39.** *Sea* $K = K_0(\mathfrak{n}, D)$ *o* $K_1(\mathfrak{n}, D)$*. El espacio de formas nuevas de nivel* $K$ *es un módulo semisimple sobre las álgebras de Hecke*

$$\mathbb{T}(K) \subset \text{End}(M_\rho(D^\times, K)), \mathbb{T}^{\text{nuevo}}(K) \subset \text{End}(M_\rho(D^\times, K)^{\text{nuevo}})$$

*generada por los operadores* $T^{class}$ *con* $T \in \mathcal{H}(K_v)$*,* $v$ *relativamente primo al conjunto* $S'_D$ *y al ideal* $\mathfrak{n}$.

*Además, si* $\lambda : \mathbb{T}(K) \to \mathbb{C}$ *es un carácter entonces el subespacio* $M_\rho(D^\times, K)^{\text{nuevo}}[\lambda]$ *de vectores propios de* $\mathbb{T}(K)$ *para el carácter* $\lambda$ *es de dimensión* $\leqslant 1$.

*Si* $\dim M_\rho(D^\times, K)^{\text{nuevo}}[\lambda] = 1$ *entonces el carácter* $\lambda$ *aparece con multiplicidad* $1$ *en* $M_\rho(D^\times, K)$ *(es decir, la multiplicidad en espacio de formas viejas es igual a* $0$*).*

Definimos $\mathbb{T}(K)_\mathbb{Q} \subset \mathbb{T}(K)$ como la $\mathbb{Q}$-subálgebra de $\mathbb{T}(K)$ generada sobre $\mathbb{Q}$ por los operadores $T^{class}$ con $T \in H_\mathbb{Z}(D_v^\times, K)$, con $v$ relativamente primo al conjunto $S'_D$ y al ideal $\mathfrak{n}$.

**Proposición 2-40.** *La descomposición* (2-38) *está definida sobre* $\mathbb{Q}$*. Más precisamente,*

1. *Los caracteres de* $\mathbb{T}(K)$ *están todos definidos sobre la clausura algebraica* $\overline{\mathbb{Q}}$ *de* $\mathbb{Q}$*, y si* $\lambda : \mathbb{T}(K)_{\mathbb{Q}} \to \overline{\mathbb{Q}}$ *es un carácter no trivial entonces, para* $\tau \in Gal(\overline{\mathbb{Q}}/\mathbb{Q})$*,* $\tau \circ \lambda$ *es también un carácter no trivial de* $\mathbb{T}(K)_{\mathbb{Q}}$*:*

$$M_\rho(D^\times, K)[\lambda] \neq 0 \Rightarrow M_\rho(D^\times, K)[\tau \circ \lambda] \neq 0, \ \forall \tau \in Gal(\overline{\mathbb{Q}}/\mathbb{Q}).$$

2. *Si* $\dim M_\rho(D^\times, K)^{\text{nuevo}}[\lambda] = 1$ *entonces,*

$$\forall \tau \in Gal(\overline{\mathbb{Q}}/\mathbb{Q}), M_\rho(D^\times, K)^{\text{viejo}}[\tau \circ \lambda] = 0.$$

*Prueba.* La primera parte de la proposición es un caso particular de un teorema de Shimura sobre la racionalidad de formas modulares. Ahora, si $\lambda'$ es un carácter de $\mathbb{T}(K)$ tal que $M_\rho(D^\times, K)^{\text{viejo}}[\lambda'] \neq 0$, entonces hay un ideal $\mathfrak{n}' \supsetneq \mathfrak{n}$ con $M_\rho(D^\times, K_1(\mathfrak{n}', D))[\lambda'] \neq 0$. La segunda parte se obtiene aplicando la primera parte a los subespacios $M_\rho(D^\times, K_1(\mathfrak{n}', D))$ de $M_\rho(D^\times, K_1(\mathfrak{n}, D))$ con $\mathfrak{n}' \supsetneq \mathfrak{n}$. □

**Ejercicio 2.5.** Sea $\mathfrak{p}$ un ideal primo de $\mathcal{O}_F$. Suponemos que

$$D_{\mathfrak{p}}^\times \xrightarrow{\sim} GL(2, F_{\mathfrak{p}}).$$

Sea $a \geqslant 1$ un entero y sea $\mathfrak{n} = \mathfrak{p}^a$. Sea $\varpi \in \mathcal{O}_{\mathfrak{p}}$ un uniformizador en $F_{\mathfrak{p}}$ y sea $\gamma_b = \begin{pmatrix} \varpi^b & 0 \\ 0 & 1 \end{pmatrix}$, $b = 1, \ldots, a$. Sea $\mathfrak{m}$ un ideal de $\mathcal{O}_F$ relativamente primo a $\mathfrak{p}$. Demostrar que la traslación por $\gamma_b^{-1} \in D_{\mathfrak{p}}^\times \xrightarrow{\sim} GL(2, F_{\mathfrak{p}})$:

$$U_b(f)(g) = f(g\gamma_b^{-1})$$

define un mapeo inyectivo

$$U_b : M_\rho(D^\times, K_0(\mathfrak{m} \cdot \mathfrak{p}^{a-b}, D)) \to M_\rho(D^\times, K_0(\mathfrak{m} \cdot \mathfrak{p}^a, D)).$$

El espacio de formas viejas es la suma de las imagenes de los $U_b$.

**Ejercicio 2.6.** Sea $\mathfrak{p}$ un ideal primo de $\mathcal{O}_F$. Suponemos que

$$D_{\mathfrak{p}}^\times \xrightarrow{\sim} GL(2, F_{\mathfrak{p}}).$$

Sea $w_{\mathfrak{p}} = \begin{pmatrix} 0 & 1 \\ \varpi & 1 \end{pmatrix}$. Sea $\mathfrak{m}$ como en el Ejercicio 2.5. Demostrar que la traslación por $w_{\mathfrak{p}}$ (definida como en el Ejercicio 2.5) define una involución de $M_\rho(D^\times, K_0(\mathfrak{m} \cdot \mathfrak{p}, D))$.

## 3 Formas modulares cuaterniónicas y curvas elípticas

**3.1 Curvas elípticas sobre un cuerpo totalmente real.** Una *curva elíptica* $E$ sobre un cuerpo $F$ es una curva no singular proyectiva de género 1 definida por ecuaciones con coeficientes en $F$, con un punto racional $\infty$ sobre $F$. Toda curva elíptica $E$ es isomorfa a una curva plana definida por una ecuación cúbica en tres variables, o en forma inhomogénea de Weierstrass

$$y^2 + a_1xy + a_3y = x^3 + a_2x^2 + a_4x + a_6, \;\; a_i \in F, i = 1, 2, 3, 4, 6.$$

Decir que una curva elíptica esta definida sobre $F$ es equivalente a decir que su invariante $j$ pertenece a $F$.

Como $E$ es de género 1, su jacobiano $J(E)$ es de dimensión 1, y la aplicación

$$E \rightarrow J(E), x \mapsto (x) - (\infty)$$

es un isomorfismo. En particular, los puntos de $E$ sobre cualquier cuerpo forman un grupo conmutativo, y designamos por 0 su elemento neutro. El teorema siguiente[2] es fundamental:

**Teorema 3-1.** *Sea $\bar{F}$ una clausura algebraica de $F$. Si $N \in \mathbb{N}$ es coprimo con la característica de $F$, el grupo $E[N] = \{x \in X(\bar{F}), Nx = 0\}$ es isomorfo a $(\mathbb{Z}/N\mathbb{Z})^2$.*

El grupo de Galois $Gal(\bar{F}/F)$ actúa sobre $E[N]$ y si fijamos una base obtenemos una representación

$$Gal(\bar{F}/F) \rightarrow GL(2, \mathbb{Z}/N\mathbb{Z}).$$

Sea $p$ un número primo invertible en $F$. Definimos el *módulo de Tate*

$$T_p(E) = \varprojlim_n E[p^n]$$

donde el mapeo $E[p^n] \rightarrow E[p^{n-1}]$ está dado por la multiplicación por $p$. Entonces $T_p(E) \xrightarrow{\sim} \mathbb{Z}_p^2$; más aún, la acción de $Gal(\bar{F}/F)$ sobre $T_p(E)$ es una representación continua

$$\rho_{E,p} : Gal(\bar{F}/F) \rightarrow GL(2, \mathbb{Z}_p).$$

Ahora suponemos que $F$ es un cuerpo de números totalmente real. Entonces es un teorema de Serre (en el caso $F = \mathbb{Q}$) y otros que la representación $\rho_{E,p}$ es absolutamente irreducible. Sea $v$ un primo de $F$ no arquimediano de característica residual distinta de $p$. Existe un conjunto finito $S = S(E)$ de primos de $F$ tal que la curva $E$ tiene buena reducción en el cuerpo residual $k(v)$ de $v$ para $v \notin S$. Sea $v \notin S$, $q_v = |k(v)|$. Entonces $E \pmod v$ es una curva elíptica sobre el cuerpo finito $k(v)$, y el grupo $E(k(v))$ es un grupo finito. Sea $N_v(E) = |E(k(v))|$, y escribimos

$$N_v(E) = 1 + q_v - a_v(E).$$

Entonces el número $a_v(E) \in \mathbb{Z}$ determina la restricción $\rho_v$ de la representación $\rho_{E,p}$ a un grupo de descomposición $\Gamma_v \subset Gal(\bar{F}/F)$. De hecho, esa restricción es no ramificada, porque $E$ tiene buena reducción, y $\rho_v$ está determinada por la traza del Frobenius $\rho_{E,p}(Frob_v)$; y

$$Tr(\rho_{E,p}(Frob_v)) = a_v(E).$$

Como el teorema de densidad de Chebotarev implica que la unión de las clases de conjugación de los $\rho_{E,p}(Frob_v)$, con $v \notin S$, es densa en la imagen de $\rho_{E,p}$, y como $\rho_{E,p}$ es absolutamente irreducible, su clase de isomorfismo está completamente determinada por los números $N_v(E)$ con $v \notin S$.

[2] Ver, por ejemplo, J. Silverman, *The Arithmetic of Elliptic Curves*, Corollary 6.4.

Si $E$ y $E'$ son dos curvas elípticas isógenas sobre $F$, entonces $S(E) = S(E')$ y $N_v(E) = N_v(E')$ para todo $v \notin S$. En consecuencia, tenemos la primera parte del teorema siguiente:

**Teorema 3-2.** *(a) Si $E$ y $E'$ son dos curvas elípticas isógenas sobre $F$, entonces $\rho_{E,p}$ y $\rho_{E',p}$ son equivalentes para todo $p$.*

*(b) [Faltings] Si $E$ y $E'$ son dos curvas elípticas sobre $F$, y si $\rho_{E,p}$ y $\rho_{E',p}$ son equivalentes como representaciones de $Gal(\bar{F}/F)$ para un número primo $p$, entonces $E$ y $E'$ son isógenas.*

**Comentario 3-3.** En particular, el teorema de Faltings implica que, si $\rho_{E,p}$ y $\rho_{E',p}$ son equivalentes como representaciones de $Gal(\bar{F}/F)$ para un número primo $p$, entonces $\rho_{E,q}$ y $\rho_{E',q}$ son equivalentes para todo primo $q$. Pero éste es un resultado mucho más elemental que el teorema de Faltings. Lo mencionamos aquí porque es importante en la demostración de la automorfía de curvas elípticas sobre cuerpos reales cuadráticos.

**3.2 El jacobiano de una curva de Shimura.** Si $X = \coprod X_i$ es una curva proyectiva suave con componentes conexas $X_i$, definimos $Jac(X) = \prod_i Jac(X_i)$. Sean $K = K_f \subset D^\times(\mathbf{A}_f)$ y $X$ la curva de Shimura ${}_KS(D)$, considerada como superficie de Riemann. El jacobiano se construye del siguiente modo. La cohomología $H^1({}_KS(D), \mathbb{C})$ admite una descomposición de Hodge

$$H^1({}_KS(D), \mathbb{C}) = \Omega^1({}_KS(D)) \oplus \bar{\Omega}^1({}_KS(D)) \tag{3-4}$$

en suma directa de espacios de diferenciales holomorfas y anti-holomorfas. Por otro lado, como tenemos el isomorfismo de la Proposición 1-9, podemos identificar

$$H^1({}_KS(D), \mathbb{C}) \xrightarrow{\sim} M_{(2,2,\ldots,2)}(D^\times, K) \oplus \bar{M}(2,2,\ldots,2)(D^\times, K) \tag{3-5}$$

Por dualidad, la descomposición (3-4) corresponde a una descomposición de la homología:

$$H_1({}_KS(D), \mathbb{C}) = \Omega^1({}_KS(D))^\perp \oplus \bar{\Omega}^1({}_KS(D))^\perp. \tag{3-6}$$

El teorema siguiente es válido para toda curva proyectiva suav e compleja:

**Teorema 3-7.** *Hay isomorfismos naturales (funtoriales) de grupos*

$$H_1({}_KS(D), \mathbb{Z})\backslash Hom_{\mathbb{C}}(\Omega^1({}_KS(D), \mathbb{C})) \xrightarrow{\sim}$$
$$H_1({}_KS(D), \mathbb{Z})\backslash H_1({}_KS(D), \mathbb{C})/\Omega^1({}_KS(D))^\perp \xrightarrow{\sim} Jac({}_KS(D)).$$

**Comentario 3-8.** Si $X$ es una curva suave proyectiva compleja, la inclusión $H_1(X, \mathbb{R}) \hookrightarrow H_1(X, \mathbb{C})$ define un isomorfismo

$$H_1(X, \mathbb{R}) \xrightarrow{\sim} H_1(X, \mathbb{C})/\Omega^1(X)^\perp$$

de espacios vectoriales reales. Así podemos escribir

$$Jac(X) \xrightarrow{\sim} H_1(X, \mathbb{Z})\backslash H_1(X, \mathbb{C})/\Omega^1(X)^\perp \xrightarrow{\sim} H_1(X, \mathbb{R})/H_1(X, \mathbb{Z})).$$

**3.3 El álgebra de Hecke como álgebra de endomorfismos del jacobiano.** Sean $X$ y $Y$ dos curvas proyectivas suaves conexas sobre el cuerpo $k$. Una *correspondencia* entre $X$ e $Y$ es una subvariedad cerrada $C \subset X \times Y$ (no necesariamente conexa) tal que las dos proyecciones $p_1 : C \to X$ y $p_2 : C \to Y$ son morfismos finitos. Sean $n_1$ y $n_2$ los grados de los morfismos finitos $p_1$ y $p_2$. Una correspondencia define un homomorfismo

$$[C] : Jac(X) \to Jac(Y)$$

del siguiente modo. Si $x \in X$, la fibra $C_x = p_1^{-1}(x) = \sum a_i c_i(x)$ es un divisor efectivo en $C$ de grado $n_1 = \sum a_i$; aquí $a_i \in \mathbb{N}$ y $c_i(x) \in C_x$. Entonces $[C](x) = \sum a_i p_2(c_i(x)) \in Jac(Y)$. Si $C$ es una curva suave (esa será el caso en nuestros ejemplos), la definición es más sencilla. Por funtorialidad del jacobiano tenemos morfismos $p_{1,*} : Jac(C) \to Jac(X)$ y $p_{2,*} : Jac(C) \to Jac(Y)$. Por dualidad (autodualidad del jacobiano), tenemos también un morfismo

$$p_1^* : Jac(X) = \widehat{Jac(X)} \to \widehat{Jac(C)} = Jac(C),$$

y $[C] = p_{2,*} \circ p_1^*$.

Si $k \subset \mathbb{C}$ es un subcuerpo de $\mathbb{C}$, y si $C$ es una curva suave, tenemos también morfismos de homología

$$p_1^* : H_1(X, \mathbb{Z}) \to H_1(C, \mathbb{Z}); p_{2,*} : H_1(C, \mathbb{Z}) \to H_1(Y, \mathbb{Z})$$

y los morfismos análogos para homología con coeficientes reales. El primer homomorfismo es definido por dualidad de Poincaré. El homomorfismo $[C]$ se define explícitamente por

$$p_{2,*} \circ p_1^* : Jac(X) = H_1(X, \mathbb{R})/H_1(X, \mathbb{Z}) \to H_1(C, \mathbb{R})/H_1(C, \mathbb{Z}) \to H_1(Y, \mathbb{R})/H_1(Y, \mathbb{Z}) = Jac(Y).$$

La definición de una correspondencia sigue siendo válida con modificaciones naturales cuando $X$ y $Y$ son curvas suaves y proyectivas pero no necesariamente conexas. Ahora sea $K = \prod_v K_v \subset D^\times(\mathbf{A}_f)$ y consideremos el caso cuando $X = Y$ es la curva de Shimura ${}_K S(D)$ (no necesariamente conexa) de nivel $K$. Esta curva es proyectiva y suave si $K$ es un subgrupo suficientemente pequeño, lo que vamos a suponer. Sea $v$ un primo no ramificado para $D$ y con $K_v = GL(2, \mathcal{O}_v)$ (decimos que $v$ es *relativamente primo al nivel de* $K$). Sea

$$I_v = \{k = \begin{pmatrix} a & b \\ c & d \end{pmatrix} \in K_v \mid c \in \varpi_v \mathcal{O}_v\}$$

con $\varpi_v$ un uniformizador de $F_v$. Sea $K_0(v) = \prod_{w \neq v} K_w \times I_v$ (hemos reemplazado $K_v$ por $I_v$). La inclusión $\iota_v : K_0(v) \hookrightarrow K_v$ define un morfismo natural $p_1 : {}_{K_0(v)}S(D) \to {}_K S(D)$. Pero hay una segunda inclusión $\iota'_v : K_0(v) \hookrightarrow K_v$, definida por

$$\iota'_v(k) = n_v k n_v^{-1}, \text{ con } n_v = \begin{pmatrix} 0 & \varpi_0^{-1} \\ 1 & 0 \end{pmatrix} \in D_v^\times$$

La inclusión $\iota'_v$ define un segundo morfismo $p_2 : {}_{K_0(v)}S(D) \to {}_K S(D)$. La imagen del morfismo

$$(p_1, p_2) : {}_{K_0(v)}S(D) \to {}_K S(D) \times {}_K S(D)$$

es una correspondencia $T(v)$.

Escribimos $Jac_K(D) = Jac({}_K S(D))$. El homomorfismo

$$[T(v)] = p_{2,*} \circ p_1^* : Jac_K(D) \;\to\; Jac_K(D)$$

es una *correspondencia de Hecke*. Por otro lado, hemos visto que $Jac_K(D)$ es isomorfo al grupo

$$H_1({}_K S(D), \mathbb{Z}) \backslash Hom_{\mathbb{C}}(\Omega^1({}_K S(D)), \mathbb{C})$$

que es cociente de $Hom_{\mathbb{C}}(M_{(2,2,\ldots,2)}(D^\times, K), \mathbb{C})$. La proposición siguiente se demuestra por un cálculo explícito elemental.

**Proposición 3-9.** *Sea $v$ un primo no arquimediano que es relativamente primo al nivel de $K$. La correspondencia de Hecke $[T(v)] \in \mathrm{End}(Jac_K(D))$ es inducida por el operador de Hecke $T_v$ actuando en $M_{(2,2,\ldots,2)}(D^\times, K)$ via la suryección*

$$Hom_{\mathbb{C}}(M_{(2,2,\ldots,2)}(D^\times, K), \mathbb{C}) \;\to\; Jac_K(D).$$

*Además, hay una correspondencia $R(v) \subset {}_K S(D) \times {}_K S(D)$ tal que $[R(v)] \in \mathrm{End}(Jac_K(D))$ es inducida por la acción del operador de Hecke $R_v$ sobre $M_{(2,2,\ldots,2)}(D^\times, K)$.*

**Corolario 3-10.** *La subálgebra $\mathbb{T}_K \subset \mathrm{End}(Jac_K(D)) \otimes \mathbb{Q}$ generada por las correspondencias $[T(v)]$ y $[R(v)]$, para todos los primos $v$ que son primos al nivel de $K$, es conmutativa.*

**3.3.1 El módulo de Tate y la cohomología étale de ${}_K S(D)$.** Recordamos un teorema fundamental de Shimura:

**Teorema 3-11** (Shimura)**.** *Sea $\sigma : F \hookrightarrow \mathbb{R}$ el único primo (real) en $\Sigma_D$. El conjunto de curvas de Shimura $\{{}_K S(D)\}$, con $K \subset D^\times(\mathbf{A}_f)$ abierto compacto, tiene un* modelo canónico ${}_K S(D)_\sigma$ *sobre el cuerpo $\sigma(F)$. Más precisamente, para todo $K$ hay un isomorfismo canónico*

$${}_K S(D) \xrightarrow{\sim} {}_K S(D)_\sigma \times_{F,\sigma} \mathbb{C}$$

*y una acción del grupo $D^\times(\mathbf{A}_f)$ sobre el conjunto de ${}_K S(D)_\sigma$ que conserva estos isomorfismos.*

En consecuencia, el sistema de jacobianos $Jac_K(D)$ tiene también un modelo canónico sobre $\sigma(F)$ (que identificamos con $F$), y podemos definir el módulo de Tate $p$-ádico $T_p(Jac_K(D)) = \varprojlim_n Jac_K(D)[p^n]$ como representación del grupo de Galois $Gal(\bar{F}/F)$. Por otro lado, hay isomorfismos canónicos entre cohomología $p$-ádica étale y espacios duales de módulos de Tate:

$$H^1({}_K S(D), \mathbb{Q}_p) \xrightarrow{\sim} Hom(T_p(Jac_K(D)), \mathbb{Q}_p) \tag{3-12}$$

como representaciones de $Gal(\bar{F}/F)$. Además, los isomorfismos (3-12) son isomorfismos de módulos sobre el álgebra de Hecke $\mathbb{T}_K$.

Para simplificar la notación, fijamos $\Sigma_D = \{\sigma\}$ y escribimos

$$X(D) = {}_{K(D)}S(D)_\sigma; \quad X_0(\mathfrak{n}, D) = {}_{K_0(\mathfrak{n},D)}S(D)_\sigma; X_1(\mathfrak{n}, D) = {}_{K_1(\mathfrak{n},D)}S(D)_\sigma;$$

y definimos los jacobianos $Jac(D) = Jac(X(D))$, $Jac_0(\mathfrak{n}, D)$, y $Jac_1(\mathfrak{n}, D)$ y las álgebras de Hecke $\mathbb{T}(D)$, $\mathbb{T}_0(\mathfrak{n}, D)$, y $\mathbb{T}_1(\mathfrak{n}, D)$ del modo evidente.

**Comentario 3-13.** En general, si una de las curvas $X(D)$, $X_0(\mathfrak{n}, D)$, y $X_1(\mathfrak{n}, D)$ es singular, la reemplazamos por su modelo suave para definir su jacobiano. Eso no cambia nada en lo esencial.

Ya hemos visto en la sección 2.7.1 que la teoría de *formas nuevas* es válida para curvas de Shimura $X_1(\mathfrak{n}, D)$ y $X_0(\mathfrak{n}, D)$. Como por un lado tenemos

$$H^1(X_i(\mathfrak{n}, D), \mathbb{Q}_p) \xrightarrow{\sim} H^1(X_i(\mathfrak{n}, D), \mathbb{Q}) \otimes \mathbb{Q}_p, i = 0, 1;$$

$$H^1(X_i(\mathfrak{n}, D), \mathbb{C}) \xrightarrow{\sim} H^1(X_i(\mathfrak{n}, D), \mathbb{Q}) \otimes \mathbb{C}, i = 0, 1$$

y como por otro lado tenemos la descomposición de cada uno de los dos factores del miembro derecho de (3-4) en formas nuevas y formas viejas, siguiendo (2-38), eso implica que el miembro izquierdo de (3-12) admite una descomposición análoga a la de (2-38):

$$\begin{aligned} &H^1(X_i(\mathfrak{n}, D), \mathbb{Q}_p) = \\ &H^1(X_i(\mathfrak{n}, D), \mathbb{Q}_p)^{\text{nuevo}} \oplus H^1(X_i(\mathfrak{n}, D), \mathbb{Q}_p)^{\text{viejo}},\ i = 0, 1. \end{aligned} \tag{3-14}$$

Por dualidad, (3-14) y (3-12) implican una descomposición del módulo de Tate de $Jac_i(\mathfrak{n}, D)$:

$$\begin{aligned} &T_p(Jac_i(\mathfrak{n}, D)) \otimes \mathbb{Q}_p = \\ &[T_p(Jac_i(\mathfrak{n}, D)) \otimes \mathbb{Q}_p]^{\text{nuevo}} \oplus [T_p(Jac_i(\mathfrak{n}, D)) \otimes \mathbb{Q}_p]^{\text{viejo}},\ i = 0, 1. \end{aligned} \tag{3-15}$$

Es una consecuencia de la Proposición 2-40 que la descomposición (3-15) proviene de una descomposición de la homología

$$\begin{aligned} &H_1(Jac_i(\mathfrak{n}, D), \mathbb{Q}) = H_1(X_i(\mathfrak{n}, D), \mathbb{Q}) \\ &= H_1(Jac_i(\mathfrak{n}, D), \mathbb{Q})^{\text{nuevo}} \oplus H_1(Jac_i(\mathfrak{n}, D), \mathbb{Q})^{\text{viejo}}. \end{aligned} \tag{3-16}$$

Definimos el *cociente nuevo* $Jac_i(\mathfrak{n}, D)^{\text{nuevo}}$ de $Jac_i(\mathfrak{n}, D)$, $i = 0, 1$, como la variedad abeliana cociente maximal $p : Jac_i(\mathfrak{n}, D) \to A$, con la propiedad de que el núcleo del mapeo $p_* : H_1(Jac_i(\mathfrak{n}, D), \mathbb{Q}) \to H_1(A, \mathbb{Q})$ sea igual a $H_1(Jac_i(\mathfrak{n}, D), \mathbb{Q})^{\text{viejo}}$.

**Proposición 3-17.** *1. La inclusión* $\mathbb{T}_i(\mathfrak{n}, D) \hookrightarrow \text{End}(Jac_i(\mathfrak{n}, D)) \otimes \mathbb{Q}$ *induce un homomorfismo*

$$\mathbb{T}_i(\mathfrak{n}, D) \to \text{End}(Jac_i(\mathfrak{n}, D)^{\text{nuevo}}) \otimes \mathbb{Q},\ i = 0, 1.$$

2. *La homología* $H_1(Jac_i(\mathfrak{n}, D)^{\text{nuevo}}, \mathbb{Q})$ *es un* $\mathbb{T}_i(\mathfrak{n}, D)$*-módulo semisimple.*

3. *Sean* $A$ *y* $B$ *dos variedades abelianas de dimensión* $\geqslant 1$, *y sean*

$$p : Jac_i(\mathfrak{n}, D)^{\text{nuevo}} \to A; \quad q : Jac_i(\mathfrak{n}, D)^{\text{nuevo}} \to B$$

*homomorfismos de variedades abelianas. Suponemos que los núcleos de* $p$ *y* $q$ *son conexos. Si* $A$ *y* $B$ *son isomorfas, entonces* $\ker p = \ker q$.

*Prueba.* Es consecuencia inmediata de la Proposición 2-40. □

**3.4 Curvas elípticas como cocientes de una curva de Shimura.** Sea $E$ una curva elíptica sobre el cuerpo totalmente real $F$. Es un celebrado teorema de Breuil, Conrad, Diamond y Taylor (2001), basado en los métodos introducidos por Andrew Wiles, que si $F = \mathbb{Q}$, $E$ es isomorfa a un cociente del jacobiano de la curva modular clásica $X_0(N)$ donde $N$ es igual al conductor de $E$. (El caso de curvas elípticas semiestables había sido demostrado un poco antes por Wiles y Taylor–Wiles, y había permitido a Wiles demostrar el Último Teorema de Fermat.)

Un teorema análogo, válido para cualquier cuerpo real cuadrático $F$ ($[F : \mathbb{Q}] = 2$) ha sido demostrado muy recientemente por Freitas, Le Hung y Siksek (2015). Pero el enunciado del teorema es un poco diferente en el caso general. En el siguiente teorema consideramos $F$ como subcuerpo del cuerpo $\mathbb{R} \subset \mathbb{C}$ por una de las inyecciones $\sigma : F \hookrightarrow \mathbb{R}$.

**Teorema 3-18** (Freitas, Le Hung, y Siksek)**.** *Sea $E$ una curva elíptica sobre el cuerpo real cuadrático $F$. Entonces $E$ es modular: existe una forma modular de Hilbert $f$ de peso $(2, 2)$, que es vector propio de los operadores de Hecke de nivel relativamente primo al nivel $K$, tal que hay una identidad de funciones $L$:*

$$L(s, f) = L(s, E).$$

*Si además existe un primo no arquimediano $v$ donde $E$ tiene reducción multiplicativa o supercuspidal [ver más abajo], entonces $E$ es isomorfa a un cociente del jacobiano de una curva de Shimura de la forma $X_0(\mathfrak{n}, D)$, donde $D$ es un álgebra de cuaterniones sobre $F$ con $|\Sigma_D| = 1$.*

En particular, la condición "$E$ es modular" no implica necesariamente que $E$ es isomorfa a un cociente del jacobiano de una curva de Shimura si $d = [F : \mathbb{Q}]$ es par. Existen ejemplos de curvas elípticas sobre cuerpos reales cuadráticos con buena reducción en todos los primos no arquimedianos; una tal curva no puede ser cociente del jacobiano de una curva de Shimura.

Tenemos que explicar el significado de "reducción multiplicativa o supercuspidal." Esas nociones tienen significados puramente geométricos, pero para nosotros es más fácil explicarlas en términos del sistema de representaciones $\{\rho_{E,p}\}$ de $Gal(\overline{\mathbb{Q}}/F)$. Sea $v$ un primo no arquimediano de $F$, y sea $\Gamma_v \subset Gal(\overline{\mathbb{Q}}/F)$ un grupo de descomposición del primo $v$ – $\Gamma_v$ es el estabilizador de una extensión $\bar{v}$ de la valuación sobre $F$ asociada al primo $v$ a la extensión $\overline{\mathbb{Q}}$. Fijamos un primo racional $p$ y sea $\rho_{p,E,v}$ la restricción de $\rho_{E,p}$ al subgrupo $\Gamma_v$ de $Gal(\overline{\mathbb{Q}}/F)$. La clase de isomorfismo de $\rho_{p,E,v}$ no depende de la elección de la extensión $\bar{v}$ de $v$. Para un $v$ fijo podemos siempre suponer que $v$ no divide a $p$. Si $\rho_{p,E,v}$ es una representación irreducible, entonces el Teorema 2-19 implica la existencia de una representación supercuspidal $\pi_v(E) = \pi(\rho_{p,E,v})$ irreducible y lisa de $GL(2, F_v)$. En este caso decimos que $E$ *tiene reducción supercuspidal*.

Más generalmente, el Teorema 2-19 es la parte difícil de la *correspondencia local de Langlands para* $GL(2)$ (sobre un cuerpo no arquimediano). Sea $k$ un cuerpo $p$-ádico, con $v$ relativamente primo a $p$, y sea $\rho : \Gamma_v \to GL(2, k)$ una representación continua. La correspondencia local de Langlands da una representación $\pi(\rho)$, irreducible y lisa, de $GL(2, F_v)$. Si $\rho$ es irreducible, es la representación del Teorema 2-19. Si $\rho$ es reducible, tenemos que distinguir entre dos casos. Sea $I_v \subset \Gamma_v$ el subgrupo de inercia.

1. Si $\rho$ es reducible y la imagen de la restricción de $\rho$ al subgrupo $I_v$ es infinita, entonces $\pi(\rho)$ es una representación de Steinberg.

2. Si $\rho$ es reducible y la imagen de la restricción de $\rho$ al subgrupo $I_v$ es finita, entonces reemplazamos $\rho$ por su semisimplificación $\rho^{ss}$. Es la suma de dos caracteres de $\Gamma_v$ con valores en $k$. Es decir, $\rho^{ss} = \sigma_1 \oplus \sigma_2$,

$$\sigma_i : \Gamma_v \to k^\times, \ i = 1, 2.$$

Sea $\chi_i : F_v^\times \to L^\times$ el carácter de $F_v^\times$ asociado a $\sigma_i$ por la teoría de cuerpo de clases para $F_v$. Entonces $\pi(\rho) = I(\chi_1, \chi_2)$ cuando esta representación inducida es irreducible (lo que el único caso que nos interesa).

3. En particular, $\pi(\rho)$ es una representación esférica si y solamente si $\rho$ es una representación no ramificada.

**Comentario 3-19.** La definición de $\pi(\rho)$ parece depender de la identificación de los caracteres $\chi_i : F_v^\times \to k^\times$ con caracteres con valores *complejos*. Cuando $\rho = \rho_{p,E,v}$, o más generalmente cuando $\rho$ proviene de una representación sobre la cohomología $p$-ádica de una curva de Shimura, es un teorema que los caracteres $\chi_i$ tienen valores algebraicos, y existe un método para reemplazar los $\chi_i$ $p$-ádicos por $\chi_i$ complejos.

Decimos que $E$ *tiene reducción multiplicativa* si $\pi(\rho_{p,E,v})$ es una representación de Steinberg.

Hemos definido para todo primo $v$ no arquimediano una representación lisa e irreducible $\pi_v(E) = \pi(\rho_{p,E,v})$ asociada a la curva elíptica $E$, o más bien a su módulo de Tate. El siguiente teorema implica que $\pi_v(E)$ es esférica para todos los primos salvo un número finito:

**Teorema 3-20.** *La curva elíptica $E$ tiene buena reducción en el primo $v$ si y solamente si, para todo primo racional $p$ relativamente primo a $v$, $\rho_{p,E,v}$ es una representación no ramificada (es decir, si y solamente si $\pi(\rho)$ es una representación esférica).*

Siguiendo el Comentario 2-27, para un primo $\sigma$ arquimediano definimos $\pi_v(E) = \pi_1$, la representación de $GL(2, \mathbb{R})/\mathbb{R}^\times$ que corresponde a las formas modulares de peso 2. Así podemos definir la representación irreducible $\pi(E) = \otimes'_v \pi_v(E)$ de $GL(2, F_{\mathbf{A}})$, teniendo en cuenta el hecho (Teorema 3-20) que casi todas las $\pi_v(E)$ son no esféricas. El Teorema 3-18 admite la siguiente reformulación:

**Teorema 3-21** (Freitas, Le Hung, y Siksek). *Sea $E$ una curva elíptica sobre el cuerpo real cuadrático $F$. Entonces la representación $\pi(E)$ de $GL(2, F_{\mathbf{A}})$ es isomorfa a una representación automorfa cuspidal.*

*Si además existe un primo no arquimediano $v$ donde $E$ tiene reducción multiplicativa o supercuspidal , entonces $\pi(E)$ es de la forma $\pi(E) = JL(\pi^D(E))$, donde $\pi^D(E)$ es una representación automorfa de $D^\times(\mathbf{A})$ con $D$ un álgebra cuaterniónica de división sobre $F$ con $|\Sigma_D| = 1$.*

Más generalmente, tenemos la siguiente conjetura.

**Conjetura 3-22.** *Sea $E$ una curva elíptica sobre el cuerpo totalmente real $F$. Entonces la representación $\pi(E)$ de $GL(2, F_{\mathbf{A}})$ es isomorfa a una representación automorfa cuspidal.*

*Si además el grado $d = [F : \mathbb{Q}]$ es impar, o si existe un primo no arquimediano $v$ donde $E$ tiene reducción multiplicativa o supercuspidal, entonces $\pi(E)$ es de la forma $\pi(E) = JL(\pi^D(E))$, donde $\pi^D(E)$ es una representación automorfa de $D^\times(\mathbf{A})$ con $D$ un álgebra cuaterniónica de división sobre $F$ con $|\Sigma_D| = 1$.*

## 3.5 Algunas ideas de la demostración de Freitas et al.

**3.5.1 La parte automorfa.** Como $d = 1$ es un número impar, la Conjetura 3-22 se aplica al caso $d = 1$, es decir, al caso $F = \mathbb{Q}$. Cuando $F = \mathbb{Q}$, la Conjetura 3-22 es el teorema de Breuil, Conrad, Diamond y Taylor (2001). La demostración en el caso de un cuerpo real cuadrático tiene la misma estructura que la demostración de Breuil, Conrad, Diamond y Taylor (ibíd.), que generaliza la demostración del caso ya tratado por Wiles (1995). Wiles trabajaba no con curvas de Shimura construidas a partir de álgebras de división sino con las curvas modulares $X_0(N)$ asociadas al grupo $GL(2) = M(2)^\times$ sobre $\mathbb{Q}$:

$$X_0(N) = [GL(2,\mathbb{Q})\backslash\mathfrak{H}^\pm \times GL(2,\mathbf{A}_\mathbb{Q})/K_0(N)\cdot\mathbf{A}_\mathbb{Q}^\times]^*$$

donde el supraíndice * designa la compactificación. Aquí hemos reemplazado el ideal $\mathfrak{n}$ por el entero $N > 1$. Como $D = M(2,\mathbb{Q})$, escribimos $\mathbb{T}_N$ en vez de $\mathbb{T}_0(N, D)$ para el álgebra de Hecke correspondiente.

Sea $E$ una curva elíptica sobre $\mathbb{Q}$. Para demostrar que $E$ es isógena a un cociente de $J_0(N) = Jac(X_0(N))$, es suficiente (Teorema 3-2) demostrar que la representación $\rho_{E,p}$ de $Gal(\overline{\mathbb{Q}}/\mathbb{Q})$ sobre el módulo de Tate de $E$ es isomorfa a una subrepresentación del módulo de Tate de $Jac(X_0(N))$. Una consecuencia de Breuil, Conrad, Diamond y Taylor (2001) es el siguiente teorema más preciso:

**Teorema 3-23.** *De hecho, $\rho_{E,p}$ es isomorfa a la representación de $Gal(\overline{\mathbb{Q}}/\mathbb{Q})$ sobre un $\mathbb{T}_N$-submódulo de*

$$T_p(Jac(X_0(N))^{\text{nuevo}} = \varprojlim_n Jac(X_0(N))^{\text{nuevo}}[p^n]$$

*de $\mathbb{Z}_p$-rango 2.*

Como $T_p(Jac(X_0(N))^{\text{nuevo}}) \xrightarrow{\sim} H_1(Jac(X_0(N))^{\text{nuevo}}, \mathbb{Z}) \otimes \mathbb{Z}_p$, un $\mathbb{T}_N$-submódulo $M \subset H_1(Jac(X_0(N))^{\text{nuevo}}, \mathbb{Z})$ define un $\mathbb{T}_N$-submódulo de $T_p(Jac(X_0(N))^{\text{nuevo}})$. Por otra parte, el isomorfismo de Hodge (3-4) identifica $M \otimes \mathbb{C}$ a un $\mathbb{T}_N$ submódulo de

$$Hom(\Omega^1(X_0(N)), \mathbb{C}) \oplus Hom(\bar{\Omega}^1(X_0(N)), \mathbb{C}).$$

Como el álgebra de Hecke $\mathbb{T}_N$ es autoadjunta para el producto escalar de Petersson, se puedo demostrar que

**Proposición 3-24.** *(a) Los $\mathbb{T}_N$-módulos $\Omega^1(X_0(N)), \mathbb{C})$ y $\bar{\Omega}^1(X_0(N)), \mathbb{C})$ son isomorfos.*

*(b) Sea $h \in H^1(Jac(X_0(N))^{\text{nuevo}}, \mathbb{C})$ un vector propio de $\mathbb{T}_N$:*

$$T(h) = \lambda(T)h \;\; \forall T \in \mathbb{T}_N.$$

*Entonces hay una forma modular nueva*

$$\omega \in \Omega^1(X_0(N)), \mathbb{C}) \xrightarrow{\sim} M_2(\Gamma_0(N))$$

*que es vector propio de $\mathbb{T}_N$ con carácter $\lambda$.*

*(c) Sea $M \subset H_1(Jac(X_0(N))^{\text{nuevo}}, \mathbb{Z})$ un $\mathbb{T}_N$-submódulo de $\mathbb{Z}$-rango 2. Entonces la acción de $\mathbb{T}_N$ sobre $M$ se factoriza por un carácter $\mathbb{T}_N \to \mathbb{Z}$.*

Sea $M \subset H_1(Jac(X_0(N))^{\text{nuevo}}, \mathbb{Z})$ un $\mathbb{T}_N$-submódulo de $\mathbb{Z}$-rango 2. La Proposición 3-24 implica que existe una forma modular nueva $\omega \in M_2(\Gamma_0(N))$ con los mismos valores propios que $M$ para el álgebra de Hecke $\mathbb{T}_N$. Ya hemos visto que $\omega$ corresponde a una representación automorfa $\Pi$ de $GL(2, \mathbf{A}_\mathbb{Q})$. Y es fácil demostrar que $\Pi$ es una representación *cuspidal*.

Por otra parte, si $\Pi$ es una representación automorfa cuspidal de $GL(2, \mathbf{A}_\mathbb{Q})$ tal que $\dim \Pi^{K_0(N)} = 1$, y si el carácter de $\mathbb{T}_N$ sobre $\Pi^{K_0(N)}$ tiene valores en $\mathbb{Z}$, entonces existe un cociente $E$ de $Jac(X_0(N))^{\text{nuevo}}$ de dimensión 1, invariante por las correspondencias de Hecke, y la acción de $\mathbb{T}_N$ sobre $T_p(E) \subset T_p(Jac(X_0(N))^{\text{nuevo}})$ (o sobre el cociente $T_p(Jac(X_0(N))^{\text{nuevo}}) \to T_p(E)$) es idéntica a la acción sobre $\Pi^{K_0(N)}$.

La misma construcción es válida si reemplazamos $\mathbb{Q}$ por un cuerpo totalmente real $F$ y $X_0(N)$ por la curva de Shimura $X_0(\mathfrak{n}, D)$. El siguiente teorema de Carayol es una generalización de un resultado celebrado de Eichler y Shimura.

**Teorema 3-25** (Carayol)**.** *Sea $E$ una curva elíptica sobre el cuerpo totalmente real $F$. Suponemos que $E$ es isomorfa al cociente de $Jac(X_0(\mathfrak{n}, D))^{\text{nuevo}}$ asociado a la representación automorfa $\Pi^D$ de $D^\times(\mathbf{A})$. Entonces $\pi(E) \xrightarrow{\sim} JL(\Pi^D)$. En particular, tenemos*

$$L(s, E) \xrightarrow{\sim} L(s, \Pi^D).$$

La demostración del teorema de Carayol (y del teorema de Eichler–Shimura) está basada en una comparación de las correspondencias $T(q)$ (para un número primo $q$) con la correspondencia de Frobenius, actuando sobre un modelo de la curva de Shimura en característica $q$. Esta comparación se llama la *Fórmula de congruencia de Eichler–Shimura*.

**3.5.2 La parte aritmética.** Hasta ahora hemos explicado (¡muy rápidamente!) cómo obtener una curva elíptica sobre $F$ a partir de una representación automorfa con ciertas propiedades. Pero el teorema de Freitas et al. va en la dirección opuesta. ¿Cómo obtener una representación automorfa $\pi$ cuando tenemos solamente la curva elíptica $E$, o más bien su representación $p$-ádica $\rho_{E,p}$?

Designamos por $\rho_E[p]$ la representación de $Gal(\overline{\mathbb{Q}}/F)$ sobre $E[p]$, que es un espacio vectorial de dimensión 2 sobre $\mathbb{F}_p$. La idea de Wiles es comenzar por suponer que $p = 3$. ¿Por qué 3? Porque un teorema de Langlands y un teorema de Tunnell implican el siguiente teorema.

**Teorema 3-26.** *Sea $F$ un cuerpo totalmente real. Sea*

$$\rho : Gal(\overline{\mathbb{Q}}/F) \to GL(2, \mathbb{F}_3)$$

*una representación continua. Entonces hay un ideal $\mathfrak{n} \subset \mathcal{O}_F$ y una forma modular nueva de Hilbert $f$ de peso $(2, 2, \ldots, 2)$ para $K_0(\mathfrak{n})$ tal que $\rho_{f,3} \equiv \rho \pmod 3$.*

No hemos definido las representaciones $\rho_{f,p} : Gal(\overline{\mathbb{Q}}/F) \to GL(2, \bar{\mathbb{Q}}_p)$ asociadas a las formas nuevas de Hilbert. Si $f$ es de peso $(2, 2, \ldots, 2)$ y de nivel $K_0(\mathfrak{n})$, y si $d = [F : \mathbb{Q}]$ es un número impar, entonces $\rho_{f,p}$ es la representación sobre un submódulo del módulo de Tate del jacobiano de una curva de Shimura, y la reducción módulo $p$ de $\rho_{f,p}$ es la acción sobre un subgrupo del grupo de elementos de $p$-torsión del jacobiano. Más generalmente la primera construcción de estas representaciones fue encontrada por Taylor en su tesis.

Como consecuencia del Teorema 3-26, sabemos que $\rho_E[3]$ es isomorfa a una representación de la forma $\rho_{f,3} \pmod 3$ para una forma modular nueva de Hilbert $f$ de peso $(2,2,\dots,2)$. El siguiente teorema del artículo Freitas, Le Hung y Siksek (2015) está basado en el método de Wiles y Taylor–Wiles, y más directamente es una consecuencia de un resultado de Breuil y Diamond (y de resultados anteriores de Kisin, Gee, y Barnet–Lamb–Gee–Geraghty):

**Teorema 3-27.** *Sea $E$ una curva elíptica sobre un cuerpo totalmente real $F$. Sea $p$ un número primo impar. Suponemos que*

1. *$\rho_E[p]$ es isomorfa a una representación de la forma $\rho_{f,p} \pmod p$ donde $f$ es una forma modular nueva de Hilbert para $F$ de peso $(2,2,\dots,2)$;*
2. *La restricción de $\rho_E[p]$ al grupo de Galois $Gal(\overline{\mathbb{Q}}/F(\zeta_p))$ es absolutamente irreducible.*

*Entonces $E$ es modular.*

Ya sabemos que $\rho_E[3]$ satisface la primera condición. Pero ¿qué hacer si la restricción de $\rho_E[3]$ a $Gal(\overline{\mathbb{Q}}/F(\zeta_3))$ es reducible? Cuando $F = \mathbb{Q}$, Wiles utiliza una propiedad geométrica de la curva $X_0(15)$ – es una curva de género uno, con un número finito de puntos racionales sobre $\mathbb{Q}$ – para mostrar que, si la imagen de $\rho_E[3]$ es demasiado pequeña, entonces se tiene (i) se puede mostrar que $\rho_E[5]$ es isomorfa a una representación de la forma $\rho_{f,5} \pmod 5$ y (ii) las imagenes de $\rho_E[3]$ y $\rho_E[5]$ no pueden ser simultáneamente pequeñas. Así el Teorema 3-27 implica que $E$ es modular. (Ver el artículo Buzzard (2012) para una explicación más detallada de la demostración de Wiles.)

Cuando $F$ es un cuerpo real cuadrático, la curva $X_0(15)$ puede tener un número infinito de puntos racionales sobre $F$, y el argumento de Wiles no funciona. Los autores de Freitas, Le Hung y Siksek (2015) se vieron obligados a generalizar este argumento. En vez de la curva $X_0(15)$ y solo las representaciones $\rho_E[3]$ y $\rho_E[5]$ utilizan una familia de 27 curvas modulares y las tres representaciones $\rho_E[3]$, $\rho_E[5]$ y $\rho_E[7]$. El análisis de puntos racionales de estas 27 curvas, sobre cuerpos reales, es muy sutil, y es poco probable que permita demostrar la Conjetura 3-22 en toda su generalidad.

## Referencias

Michael Harris
harris@math.columbia.edu
Columbia University

# 2 | Combinatoria aditiva

# LEMAS DE ELIMINACIÓN Y APLICACIONES

Juanjo Rué

**Resumen**

En este capítulo demostramos el teorema de Roth, que afirma que si un conjunto $A \subseteq \{1, \dots, n\}$ no contiene progresiones aritméticas de longitud 3, entonces $|A| = o(n)$. Para demostrar dicho resultado utilizaremos técnicas procedentes de la combinatoria. Más concretamente, enunciamos el lema de regularidad de Szemerédi (resultado fundamental en teoría extremal de grafos) y lo aplicamos para obtener el denominado lema de eliminación de triángulos. Dicho lema es el paso previo a la demostración del teorema de Roth. Finalmente, se discuten extensiones del método a configuraciones distintas de las progresiones aritméticas, así como en un contexto más general de los grupos finitos.

## 1 Introducción

En esta sección del curso de combinatoria aditiva vamos a mostrar aplicaciones de la teoría de grafos en el contexto de la teoría de números. Para ello, empezaremos hablando de una técnica fundamental en combinatoria extremal denominada *Lema de regularidad de Szemerédi*. Dicho resultado dice (de manera informal) que todo grafo denso suficientemente grande puede descomponerse en subgrafos que tienen un comportamiento muy parecido a los grafos bipartitos aleatorios. En particular mostraremos como aplicar esta técnica en el denominado *Lema de eliminación de triángulos*, y como consecuencia hallaremos una prueba del teorema de Roth puramente combinatoria, así como extensiones de ésta.

**Notación:** si no se dice lo contrario, todos los grafos que se considerarán serán simples. Usando la notación habitual, para un grafo $G = (V, E)$, $V$ denotará el conjunto de vértices y $E$ el correspondiente conjunto de aristas. Dado un vértice $x \in V$, el conjunto de vecinos de $x$ (que denotaremos por $N(x)$) es el conjunto de vértices incidentes con $x$. El grado de $x$ (que denotaremos por $d(x)$) es el cardinal de $N(x)$. Finalmente, denotaremos el conjunto de enteros $\{1, \dots, n\}$ por $[n]$.

J. R. ha sido parcialmente financiado por los proyectos FP7-PEOPLE-2013-CIG CountGraph (referencia 630749), por el proyecto colateral PROCOPE-DAAD *RanConGraph* (referencia 57134837) el DFG *Research Training Group Methods for Discrete Structures* (GRK1408), y la *Berlin Mathematical School*.

## 2 El lema de regularidad de Szemerédi.

En el curso de su demostración de que existen progresiones aritméticas arbitrariamente largas en conjuntos de densidad positiva, Szemerédi introdujo un método en teoría de grafos que ha resultado ser fundamental en combinatoria, y especialmente útil en el estudio de problemas de tipo extremal. Es el que se denomina *Lema de regularidad de Szemerédi*. *Grosso modo* dicho resultado dice lo siguiente: fijada una tolerancia $\varepsilon > 0$, existe siempre un valor $N_0 := N_0(\varepsilon)$ tal que todo grafo con un número de vértices $n > N_0(\varepsilon)$ admite una partición de sus vértices tal que:

1. El número de bloques de la partición únicamente depende de $\varepsilon$, y no de $n$.
2. La estructura de las aristas entre la mayoría de las parejas de bloques puede modelarse como en un grafo aleatorio bipartito.

Formalicemos estas ideas. Sea $G = (V, E)$ un grafo. Sean $X$ e $Y$ subconjuntos disjuntos de $V$. Denotemos por $e(X, Y)$ el número de aristas entre $X$ e $Y$, y definamos

$$d(X, Y) = \frac{e(X, Y)}{|X||Y|},$$

es decir, la correspondiente densidad de aristas definida por el par de conjuntos de vértices. La definición fundamental es ahora la siguiente:

**Definición 2.1** (Par $\varepsilon$-regular). *Sea $G = (V, E)$ un grafo y $\{X, Y\}$ una pareja de subconjuntos disjuntos de $V$. Sea $\varepsilon > 0$. Decimos que el par $\{X, Y\}$ es $\varepsilon$-regular si para todo $A \subseteq X$, $B \subseteq Y$, $|A| \geqslant \varepsilon|X|$ y $|B| \geqslant \varepsilon|Y|$ se cumple que*

$$|d(X, Y) - d(A, B)| \leqslant \varepsilon. \qquad (1)$$

Esta definición nos dice que la densidad de aristas es en cierto modo "uniforme" bajo la elección de subconjuntos suficientemente grandes. Esta propiedad es una versión más débil de lo que ocurre en el siguiente modelo aleatorio de tipo binomial: tomemos dos conjuntos disjuntos de vértices $V_1$ y $V_2$. Construyamos un grafo bipartito incluyendo cada una de las posibles $|V_1||V_2|$ aristas independientemente con probabilidad $p$. Si denotamos por $d(V_1, V_2)$ la variable aleatoria que mide la densidad de aristas, entonces su esperanza es igual a $\mathbb{E}[d(V_1, V_2)] = p$. Si realizamos el mismo cálculo para subconjuntos $W_1 \subseteq V_1$ y $W_2 \subseteq V_2$ obtenemos que $\mathbb{E}[d(W_1, W_2)] = p$, y por lo tanto $\mathbb{E}[d(W_1, W_2) - d(V_1, V_2)] = 0$.

El siguiente paso es ahora descomponer un grafo dado en término de pares $\varepsilon$-regulares. Este hecho viene descrito en la siguiente definición:

**Definición 2.2** (Partición $\varepsilon$-regular). *Sea $G = (V, E)$ un grafo y $\varepsilon > 0$ un valor dado. Sea $V_0 \cup V_1 \cup V_2 \cdots \cup V_k$ una partición de $V$. Decimos que esta partición es $\varepsilon$-regular si:*

*(a)* $|V_0| \leqslant \varepsilon|V|$,

*(b)* $|V_1| = \cdots = |V_k|$,

*(c) Todos excepto $\varepsilon k^2$ de los pares $\{V_i, V_j\}$ son $\varepsilon$-regulares.*

La existencia de una tal partición no es del todo clara. Afortunadamente lo que ocurre es que, para grafos suficientemente grandes, dichas particiones siempre existen. Este es el resultado que se conoce como el *lema de regularidad de Szemerédi:*

**Teorema 2.3** (Lema de regularidad de Szeméredi)**.** *Sea* $\varepsilon > 0$ *y* $m > 0$. *Existe entonces un entero* $M := M(\varepsilon, m)$ *tal que todo grafo con más de* $m$ *vértices admite una partición* $\varepsilon$*-regular* $\{V_0, \ldots, V_k\}$ *con* $m \leqslant k \leqslant M$.

Notar que en este resultado el número de bloques depende *únicamente* de $\varepsilon$ y del parámetro $m$, pero no del tamaño del grafo ($n$). La idea de prueba del teorema (ver por ejemplo Diestel (2005)) es la siguiente: empecemos con una partición arbitraria de $V = V_0 \cup V_1 \cup V_2 \cdots \cup V_r$ del mismo tamaño, y procedamos repetidamente a refinar dicha partición (es decir, subdividir cada uno de los $V_i$ en subconjuntos disjuntos). Entonces, se puede definir una *función de energía* sobre dicha partición que se incrementa un valor constante en cada refinamiento realizado. Después de un número finito de iteraciones, necesariamente la partición debe ser $\varepsilon$-regular, ya que de otro modo se llegaría a una contradicción al mirar la energía final obtenida.

Actualmente este lema es una de las piedras angulares de la teoría extremal de grafos. Desafortunadamente la dependencia de los parámetros en el teorema es muy mala: el algoritmo que se he utiliza para refinar la partición utiliza $\frac{1}{\varepsilon^4}$ iteraciones para llegar a una partición $\varepsilon$-regular.

## 3 Una aplicación: contando triángulos

Veamos una aplicación de este resultado que será muy útil en la teoría de números. El lema de regularidad de Szemerédi nos permite contar triángulos en grafos de manera sencilla:

**Teorema 3.1** (Lema contador de triángulos)**.** *Sea* $G = (V, E)$ *un grafo y* $X \cup Y \cup Z$ *una partición* $V$. *Supongamos que* $d(X, Y) = \alpha$, $d(X, Z) = \beta$ *y* $d(Y, Z) = \gamma$. *Sea* $\varepsilon > 0$ *tal que* $\min\{\alpha, \beta, \gamma\} \geqslant 2\varepsilon$. *Asumamos también que los pares* $\{X, Y\}$, $\{Y, Z\}$ *y* $\{X, Z\}$ *son* $\varepsilon$*-regulares.*

*Entonces, el número de triángulos* $\triangle xyz$, *con* $x \in X$, $y \in Y$ *y* $z \in Z$ *es mayor o igual que*

$$(1 - 2\varepsilon)(\alpha - \varepsilon)(\beta - \varepsilon)(\gamma - \varepsilon)|X||Y||Z|.$$

*Prueba.* Para un vértice $v \in V$ denotemos por $d_X(v)$ $d_Y(v)$ y $d_Z(v)$ el número de vecinos de $v$ en $X$, $Y$ y $Z$, respectivamente.

Empecemos la prueba demostrando que controlamos el número de vértices en $X$ con pocos vecinos en $Y$ y en $Z$: debemos demostrar este hecho para poder aplicar más tarde la condición de $\varepsilon$-regularidad. Veamos de hecho que se cumple lo siguiente:

$$|\{x \in X : d_Y(x) < (\alpha - \varepsilon)|Y|\}| < \varepsilon|X|.$$

Supongamos lo contrario: $|\{x \in X : d_Y(x) < (\alpha - \varepsilon)|Y|\}| \geqslant \varepsilon|X|$. Para simplificar la notación, escribamos $\{x \in X : d_Y(x) < (\alpha - \varepsilon)|Y|\} = X'$. Entonces, $|X'| \geqslant \varepsilon|X|$. Usando la definición de pareja $\varepsilon$-regular, $|d(X', Y) - d(X, Y)| < \varepsilon$. Como $d(X, Y) = \alpha$ y $d(X', Y) < \alpha - \varepsilon$ (observar que cada vértice $x$ en $X'$ contribuye en como mucho $(\alpha - \varepsilon)|Y|$ aristas) tenemos que

$$d(X', Y) - d(X, Y) < \alpha - \varepsilon - \alpha < -\varepsilon,$$

pero esto contradice el hecho que $|d(X', Y) - d(X, Y)| < \varepsilon$. Permutando letras, el mismo argumento es válido cuando estudiamos vértices con pocos vecinos en $Y$ y en $Z$. Esto nos lleva a concluir que

$$|\{x \in X : d_Y(x) \geqslant (\alpha - \varepsilon)|Y| \text{ y } d_Z(x) \geqslant (\beta - \varepsilon)|Z|\}| \geqslant (1 - 2\varepsilon)|X|.$$

Tomemos ahora $x \in \{x \in X : d_Y(x) \geqslant (\alpha - \varepsilon)|Y| \text{ y } d_Z(x) \geqslant (\beta - \varepsilon)|Z|\}$ y estudiemos en cuantos triángulos está involucrado. En particular, $|N(x) \cap Y| = d_Y(x) \geqslant (\alpha - \varepsilon)|Y| \geqslant \varepsilon|Y|$ y $|N(x) \cap Z| = d_Z(x) \geqslant (\beta - \varepsilon)|Z| \geqslant \varepsilon|Z|$ (recordar que teníamos que $\min\{\alpha, \beta, \gamma\} \geqslant 2\varepsilon$). Ahora utilizando la condición de $\varepsilon$-regularidad,

$$|d(Y, Z) - d(N(x) \cap Y, N(x) \cap Z)| < \varepsilon \Rightarrow \gamma - \frac{e(N(x) \cap Y, N(x) \cap Z)}{|N(x) \cap Y||N(x) \cap Z|} < \varepsilon.$$

Concluimos entonces que $e(N(x) \cap Y, N(x) \cap Z) > (\alpha - \varepsilon)(\beta - \varepsilon)(\gamma - \varepsilon)|Y||Z|$. Cada arista en este conjunto define un triángulo usando el vértice $x$, por lo tanto el número de triángulos es mayor o igual que el número de triángulos donde el vértice $x$ satisface la anterior propiedad. Dicho valor es $(1 - 2\varepsilon)|X|e(N(x) \cap Y, N(x) \cap Z)$, como queríamos demostrar. □

Con este lema preliminar ya podemos pasar a demostrar la primera aplicación importante de regularidad, el denominado lema de eliminación de triángulos:

**Teorema 3.2** (Lema de eliminación de triángulos)**.** *Por cada $\varepsilon > 0$ existe un valor $\delta := \delta(\varepsilon) > 0$ (tal que $\delta \to 0$ cuando $\varepsilon \to 0$) y un $n_0 := n_0(\varepsilon)$ tal que para cada grafo $G$ con $n \geqslant n_0$ vértices y como mucho $\delta n^3$ triángulos, se puede hacer libre de triángulos eliminando como mucho $\varepsilon n^2$ aristas.*

Antes de demostrar el teorema, veamos que dice: tomemos un grafo $G$ con $n$ vértices ($n$ suficientemente grande). Si éste tiene $o(n^2)$ aristas, el lema es trivial: podemos eliminar todos los triángulos eliminando todas las aristas. El resultado es altamente no trivial cuando el grafo es denso: en este caso, el lema nos dice que se pueden eliminar todos los triángulos eliminando únicamente una proporción pequeña del total de las aristas.

Veamos ahora su prueba usando el Lema de regularidad de Szemerédi:

*Prueba.* Veamos la implicación opuesta: si necesitamos eliminar como mínimo $\varepsilon n^2$ aristas para hacer que el grafo resultante sea libre de triángulos, entonces veremos que el grafo inicial tenía más de $\delta n^3$ triángulos.

Tomemos $\varepsilon > 0$, y elijamos $m = \lfloor \frac{4}{\varepsilon} \rfloor$. Ahora, consideremos una partición $\frac{\varepsilon}{4}$-regular de $G$ con partición $V_0 \cup V_1 \cup \cdots \cup V_k$. Dicha partición existe por el Lema de regularidad de Szemerédi. Escribamos $c = |V_1| = \cdots = |V_k|$. Observar que $\lfloor \frac{4}{\varepsilon} \rfloor < k$, y que $kc < n$ (en la última desigualdad no estamos añadiendo la contribución del número de vértices en $V_0$, que satisface que $|V_0| \leqslant \frac{\varepsilon}{4}n$).

Empecemos eliminando algunas pocas aristas de $G$ para los siguientes conjuntos de vértices:

1. Todas las aristas que son incidentes con $V_0$ (en $V_0$ o exteriores): tenemos como mucho $|V_0|n = \frac{\varepsilon}{4}n^2$ aristas de este tipo.

2. Todas las aristas interiores en $V_1, \ldots, V_k$: tenemos como mucho $k\binom{c}{2} < kc^2 < \frac{n^2}{k} < \frac{\varepsilon}{4}n^2$ aristas de este tipo.

3. Todas las aristas definidas por pares que no son $\frac{\varepsilon}{4}$-regulares: recordar que hay como mucho $\frac{\varepsilon}{4}k^2$ de estas parejas, con lo que resulta que el número total de aristas en esta situación es menor que $\frac{\varepsilon}{4}k^2c^2 < \frac{\varepsilon}{4}n^2$.
4. Todas las aristas entre pares $\{V_i, V_j\}$ $\frac{\varepsilon}{4}$-regulares, con $d = d(V_i, V_j) < \frac{\varepsilon}{2}$. En este caso se tiene como mucho $\binom{k}{2}$ de estas parejas (cota trivial), y el número de aristas está acotado entonces por $\binom{k}{2}d(V_i, V_j)|V_i||V_j| < \frac{k^2}{2}\frac{\varepsilon}{2}c^2 < \frac{\varepsilon}{4}n^2$.

Resumiendo, añadiendo las contribuciones anteriores, hemos eliminado como mucho $\varepsilon n^2$ aristas en total. Si en este momento el grafo resultante es libre de triángulos, entonces hemos acabado. De no ser así, el grafo todavía contiene triángulos, y necesitamos eliminar más aristas con el fin de obtener un grafo sin triángulos. Observar que las aristas que han sobrevivido están definidas entre pares $\frac{\varepsilon}{4}$-regulares cuya densidad es mayor que $\frac{\varepsilon}{2}$. Basta tomar entonces 3 de estos estos conjuntos (a los que llamamos $V_i$, $V_j$ y $V_k$) y observar que las condiciones necesarias en el Lema contador de triángulos se satisfacen (tomando $\frac{\varepsilon}{4}$ en lugar de $\varepsilon$). Dicho de otro modo:

- $d(V_i, V_j) = \alpha$, $d(V_j, V_k) = \beta$, $d(V_i, V_k) = \gamma$, y $\min\{\alpha, \beta, \gamma\} \geqslant \frac{\varepsilon}{2}$.
- Cada par es $\frac{\varepsilon}{4}$-regular (por hipótesis).

Es por ello que por el Lema contador de triángulos, estos tres conjuntos definen como mínimo

$$\left(1 - \frac{\varepsilon}{2}\right)\left(\frac{\varepsilon}{4}\right)^3 c^3$$

triángulos. Esta cota puede escribirse en términos de $n$: como $n = |V_0| + kc$, $|V_0| \leqslant \frac{\varepsilon}{4}n$ y $k \leqslant M(m, \frac{\varepsilon}{4}) := M(\varepsilon)$ tenemos que

$$n = |V_0| + ck \Rightarrow c > \frac{1}{k}\left(1 - \frac{\varepsilon}{4}\right)n > \frac{1}{M(\varepsilon)}\left(1 - \frac{\varepsilon}{4}\right)n,$$

y por lo tanto, el número de triángulos (definidos por esta tripleta, y por lo tanto, en el grafo total) es como mínimo

$$\left(1 - \frac{\varepsilon}{2}\right)\left(\frac{\varepsilon}{4}\right)^3 \frac{1}{M(\varepsilon)^3}\left(1 - \frac{\varepsilon}{4}\right)^3 n^3.$$

Eligiendo ahora $\delta = (1 - \frac{\varepsilon}{2})\left(\frac{\varepsilon}{4}\right)^3 \frac{1}{M(\varepsilon)^3}\left(1 - \frac{\varepsilon}{4}\right)^3$ tenemos demostrado el resultado: si no hubiéramos eliminado todos los triángulos al eliminar $\varepsilon n^2$ aristas, entonces el grafo hubiera tenido más de $\delta n^3$ triángulo desde un buen principio). En particular, $\delta \to 0$ cuando $\varepsilon \to 0$. □

## 4 Una prueba del teorema de Roth.

Veamos finalmente una aplicación de toda esta tecnología en el contexto de la teoría de números. Dicha prueba es debida a Szemerédi y Ruzsa con el propósito de demostrar el ya clásico teorema de Roth.

**Teorema 4.1** (Teorema de Roth). *Sea $A \subseteq [n]$ sin progresiones aritméticas de longitud 3. Entonces* $|A| = O\left(\frac{n}{\log\log(n)}\right)$.

Este resultado fue demostrado por primera vez por Klaus Roth en el año 1954 mediante el uso de técnicas de análisis de Fourier. Actualmente, la mejor cota superior para el tamaño de un conjunto de enteros libre de progresiones aritméticas de longitud tres es debida a Bloom (2016), quién demuestra que si $A \subseteq [n]$ no contiene progresiones aritméticas de longitud 3 entonces

$$|A| = O\left(\frac{n(\log\log(n))^4}{\log(n)}\right)$$

Nuestro objetivo es más modesto, y consiste en demostrar una versión un poco más débil de este resultado ($o(n)$ en lugar de $O\left(\frac{n}{\log(\log(n))}\right)$, pero que puede demostrarse únicamente con técnicas procedentes de la teoría de grafos.

**Teorema 4.2** (Teorema de Roth, versión 2). *Sea $A \subseteq [n]$. Si $A$ no contiene una progresión aritmética de longitud 3, entonces $|A| = o(n)$.*

*Prueba.* Escribiremos 3-AP para denotar *progresión aritmética de longitud 3*. Mostremos que para todo $\varepsilon > 0$ y $A \subseteq [n]$ cumpliendo que $|A| > \varepsilon n$ (donde $n$ es suficientemente grande, se determinará más tarde), entonces $A$ debe contener una 3-AP. Con esta idea en mente, vamos a construir un grafo adecuado y aplicar el lema de eliminación de triángulos. Para un conjunto $S \subseteq [n]$, definimos el grafo $H(S) = (V, E)$ cuyo conjunto de vértices es igual a $\{(i,1) : i \in [n]\} \cup \{(j,2) : j \in [2n]\} \cup \{(k,3) : k \in [3n]\}$ (por lo tanto, $|V| = 6n$) y su conjunto de aristas $E$ viene definido por:

- $(i,1)$ y $(j,2)$ son adyacentes si y sólo si $j - i \in S$.
- $(j,2)$ y $(k,3)$ son adyacentes si y sólo si $k - j \in S$.
- $(i,1)$ y $(k,3)$ son adyacentes si y sólo si $k - i \in 2 \cdot S = \{2s : s \in S\}$.

Observar ahora que si $(i,1), (j,2)$ y $(k,3)$ definen un triángulo en $H(S)$, escribiendo $j - i = a_1$, $k - j = a_2$ y $k - i = 2a_3$ ($a_i \in S$), entonces $\{a_1, a_2, a_3\}$ define una 3-AP (de hecho $x, y, z$ define una 3-AP si y sólo si $x + y = 2z$, que en este caso se satisface: $a_1 + a_2 = (j - i) + (k - j) = k - i = 2a_3$). Adicionalmente, la ternas $(i,1)$, $(i + s, 2)$, $(i + 2s, 3)$ con $s \in S$ son triángulos asociados a las 3-AP's triviales (diferencia 0) $s$, $s + 0$, $s + 2 \cdot 0$. De este tipo de triángulos hay entonces $|S| \cdot n$. Obviamente esos triángulos son disjuntos, es por ello que debemos eliminar como mínimo $|S|n$ aristas para hacer que el grafo resultante sea libre de triángulos.

Supongamos ahora que $|A| > \varepsilon n$ (y que, obviamente $|A| \leqslant n$). Entonces, el número de aristas necesarias para eliminar todos los triángulos en $H(A)$ es como mínimo $\varepsilon n^2 = \frac{\varepsilon}{36}(6n)^2 = \frac{\varepsilon}{36}|V|^2$ (necesitamos estas aristas como mínimo para eliminar los triángulos triviales). Invocando ahora el Lema de eliminación de triángulos, existe un $\delta := \delta(\frac{\varepsilon}{36})$ tal que el grafo $H(A)$ tiene como mínimo $\delta|V|^3 = \delta 6^3 n^3$ triángulos.

Por lo tanto, el número de triángulos no triviales en $H(A)$ es como mínimo $\delta 6^3 n^3 - n^2$. De este modo, tomando $n$ tal que

$$0 < \delta 6^3 n^3 - \varepsilon n^2 \Rightarrow n > \frac{\varepsilon}{6^3 \delta},$$

aseguramos la existencia de algún triángulo no trivial en $H(A)$, y por lo tanto $A$ debe contener alguna 3-AP. □

Volviendo a la prueba, observamos que el número de triángulos es cúbico. Este hecho es un caso particular de un fenómeno muy general denominado *supersaturación*. La propiedad de supersaturación nos da de hecho algo más fuerte que lo que hemos afirmado: el número de 3-AP es de hecho cuadrático. A este resultado se le conoce como *teorema de Varnavides*.

Actualmente existe una gran variedad de resultados inspirados en esta técnica. Green (2005), utilizando técnicas del análisis de Fourier, demostró que para grupos abelianos existen lemas de eliminación en el siguiente sentido:

**Teorema 4.3** (Green (ibíd.)). *Sea $G$ un grupo abeliano finito de tamaño $N$, y sea $A$ un subconjunto de $G$. Entonces, si el número de soluciones de la ecuación $x_1 + \cdots + x_r = 0$, $x_i \in A$ es $o(N^{r-1})$, entonces se pueden eliminar $o(N)$ elementos de $A$ (obteniendo $A'$) tal que la ecuación $x_1 + \cdots + x_r = 0$ es libre de soluciones en $A'$.*

Las resultados obtenidos por Green se basan en técnicas analíticas, que dejan de ser eficaces en el estudio de ecuaciones en grupos no abelianos. Posteriormente, por medio de argumentos combinatorios, Králʹ, Serra y Vena consiguieron estudiar el mismo problema en grupos no abelianos, ver Král, Serra y Vena (2009). Estas ideas se han aplicado posteriormente con gran eficacia en contextos más generales, como sistemas de ecuaciones sobre grupos y cuerpos finitos (ver Serra y Vena (2014), Král, Serra y Vena (2012) y Shapira (2010)).

## 5 Ejercicios

1. **Problema 1:** Demostrar que una partición $\varepsilon$-regular de un grafo $G = (V, E)$ define también una partición $\varepsilon$-regular del grafo complementario (es decir, del grafo obtenido tomando el mismo conjunto de vértices $V$ pero tomando todas las aristas *no* contenidas en $E$).

2. **Problema 2:** Si $A$ y $B$ definen un par $\varepsilon$-regular con $d(A, B) = \delta$, y $A'$, $B'$ verifican que $|A'| \geqslant \gamma|A|$ and $|B'| \geqslant \gamma|B|$ para algún $\gamma \geqslant \varepsilon$, entonces el par $\{A', B'\}$ es $\varepsilon \cdot \max\left\{\frac{1}{\gamma}, 2\right\}$ – regular con densidad $[\delta - \varepsilon, \delta + \varepsilon]$.

3. **Problema 3:** (Ajtai–Szemerédi) Un *triángulo rectángulo isósceles* en $[2n]^2$ es un conjunto de puntos de la forma $\{(x, y), (x + h, y), (x, y + h) : x, y, h \in [n]\}$. Demostrar que para cada $\varepsilon > 0$ existe $n_0 := n_0(\varepsilon)$ tal que si $n > n_0$ entonces si $A \subseteq [2n]^2$ no contiene ningún triángulo rectángulo isósceles, entonces $A = o(n^2)$.

4. **Problema 4:** (Frankl, Graham y Rödl) Demostrar la siguiente generalización del teorema de Roth: para cada $\varepsilon > 0$ existe un valor $n_0 = n_0(\varepsilon)$ con la siguiente propiedad: suponer que $G$ es un grupo abeliano de orden impar, $|G| > n_0$. Entonces cada subconjunto $B \subset G$ con $|B| > \varepsilon|G|$ contiene tres elementos distintos $x, y, z$ que satisfacen $x + y = 2z$.

## Referencias

JUANJO RUÉ
juan.jose.rue@upc.edu
FREIE UNIVERSITÄT BERLIN
INSTITUT FÜR MATHEMATIK UND INFORMATIK
ARNIMALLEE 3
14195 BERLIN
GERMANY

and

U. POLITÈCNICA DE CATALUNYA
DEPARTMENT OF MATHEMATICS
C3 BUILDING (EETAC), OFFICE 003
C. ESTEVE TERRADAS, 5
08860 CASTELLDEFELS
SPAIN

# COMBINATORIA ADITIVA, PRIMERA PARTE: CONJUNTOS SUMA, EL FENÓMENO SUMA-PRODUCTO Y EL MÉTODO POLINOMIAL

Julia Wolf

## Índice general

**Resumen**

En esta parte del curso se estudiarán los primeros resultados en combinatoria aditiva. Dicha área se inició en el contexto del estudio de la estructura de conjuntos de enteros en relación a las operaciones de suma y producto. En esta introducción a la materia, mostraremos que muchas de las cuestiones planteadas en los enteros se pueden considerar también en grupos abelianos o en cuerpos finitos.

El desarrollo espectacular en los últimos años de esta área de investigación ha demostrado que no sólo se ocupa del estudio de problemas aislados, sino que ha jugado un papel importante como punto de interacción de diversas áreas de las matemáticas. En esta primera parte usaremos especialmente técnicas combinatorias, geométricas y algebraicas, y ya más adelante en el curso se mostrará el uso de técnicas probabilísticas, analíticas y procedentes de la teoría de grafos. Finalmente destacar las conexiones íntimas de esta área con otros temas que se tratarán en esta escuela, como son los expansores o el análisis armónico mediante sus conexiones con el denominado problema de Kakeya.

## 1 La estructura de conjuntos suma

En combinatoria aditiva investigamos la estructura de conjuntos bajo la operación "+".[1] Dados conjuntos $A$ y $B$ en un grupo conmutativo, podemos definir el *conjunto suma* $A + B := \{a + b : a \in A, b \in B\}$. Es fácil ver que este conjunto suma tiene tamaño mayor o igual que $\max\{|A|, |B|\}$, y que dicho cardinal puede llegar a ser $|A||B|$ cuando todas la sumas son distintas.

En los enteros, grupo donde se inició el estudio de estas cuestiones en los años 60, observamos que el cardinal de $A + A$ es mínimo, por ejemplo, cuando $A$ es una progresión aritmética (en este caso se puede calcular $|A + A| = 2|A| - 1$). De manera similar, se observa que una progresión geométrica, por ejemplo, tiene un conjunto suma muy grande. [2]

[1] No podemos dar la historia completa de este campo de investigación, pero sugerimos al lector el libro Nathanson (1996) para un tratamiento extenso de la teoría clásica.

En el lema siguiente mostramos que todo conjunto $A$ cuyo conjunto suma tiene tamaño $2|A|-1$ es una progresión aritmética.

**Lema 1.1.** *Sea $A \subseteq \mathbb{Z}$ un conjunto finito. Entonces $|A + A| \geqslant 2|A| - 1$, con igualdad si y sólo si $A$ es una una progresión aritmética.*

DEMOSTRACIÓN: Escribamos el conjunto $A$ como una sucesión finita $a_1, \ldots, a_n$, escrita en orden creciente. Entonces

$$a_1 + a_1 < a_1 + a_2 < a_1 + a_3 < a_1 + a_4 < \ldots$$
$$< a_1 + a_n < a_2 + a_n < a_3 + a_n < \cdots < a_n + a_n.$$

Por lo tanto, se deduce que hay como mínimo $2|A| - 1$ elementos distintos en el conjunto suma $A + A$. Observar que hubiera sido igualmente lícito el uso de la cadena de desigualdades

$$a_1 + a_1 < a_1 + a_2 < a_2 + a_2 < a_2 + a_3 < \ldots$$
$$< a_2 + a_n < a_3 + a_n < a_4 + a_n < \cdots < a_n + a_n.$$

En el caso extremo, es decir cuando $|A + A| = 2|A| - 1$, las dos cadenas presentadas deben ser idénticas. Concluimos pues que $a_2 + a_i = a_1 + a_{i+1}$ para $i = 2, \ldots, n - 1$, lo que obliga al conjunto $A$ a ser una progresión aritmética. □

**Ejercicio 1.2.** *Sean $A, B \subset \mathbb{Z}$ conjuntos finitos verificando $|A + B| = |A| + |B| - 1$. Demostrar que $A$ y $B$ son progresiones aritméticas con la misma diferencia.*

Indicación: *El caso $|A| = |B|$ es más fácil.*

Como segundo ejemplo estudiaremos el mismo problema en un grupo cíclico módulo un primo $p$, que denotamos por $\mathbb{Z}_p$. Como ejercicio de calentamiento vamos a mostrar que cuando el tamaño de $A \subseteq \mathbb{Z}_p$ es mayor que la mitad del cardinal del grupo $\mathbb{Z}_p$, entonces el conjunto suma $A + A$ cubre todo el grupo.

**Ejercicio 1.3.** *Sea $p$ un número primo, y sean $A, B \subseteq \mathbb{Z}_p$ conjuntos tales que $|A| + |B| > p$. Demostrar que $A + B = \mathbb{Z}_p$.*

El siguiente es un teorema clásico que modela la versión modular del Lema 1.1.

**Teorema 1.4** (Teorema de Cauchy–Davenport)**.** *Sean $A, B \subseteq \mathbb{Z}_p$. Entonces*

$$|A + B| \geqslant \min(p, |A| + |B| - 1).$$

DEMOSTRACIÓN: Sin pérdida de generalidad suponemos que $0 \in B$. Nótese que el teorema se cumple cuando $|A| + |B| > p$ (ver el Ejercicio 1.3), pero también trivialmente cuando $|A| = 1$ o $|B| = 1$. Supongamos entonces que los conjuntos $A$ y $B$ son tales que $|A| \geqslant 2, |B| \geqslant 2$,

[2] Es un primer indicio del hecho que las operaciones de adición y multiplicación no son compatibles. Veremos con más detalle este fenómeno en el Capítulo 2.

$|A|+|B| < p-1$, pero $|A+B| < |A|+|B|-1$. Elegimos conjuntos $A$ y $B$ con estas propiedades tales que el tamaño de $B$ es mínimo.

Nuestro objetivo es construir conjuntos $A', B'$ que satisfagan las mismas desigualdades, pero tal que $|B'| < |B|$. Primero observamos que, dado que $|B| \geqslant 2$, existe un elemento $b \in B, b \neq 0$. Si $a+b \in A$ para todo $a \in A$, entonces $a+jb \in A$ para todo $j = 0,1,2,\ldots$. Pero en este caso $A \supseteq \{a+jb : j = 0,1,2\ldots\} = \mathbb{Z}_p$, lo que es imposible dado que $\mathbb{Z}_p \supseteq A \neq \mathbb{Z}_p$. En consecuencia existe $a^* \in A$ tal que $a^*+b \notin A$, y entonces $b \notin A - a^*$. Consideramos los conjuntos

$$A(a^*) = A \cup (B+a^*) \quad \text{y} \quad B(a^*) = B \cap (A-a^*)$$

(estos conjuntos $A(a^*)$ y $B(a^*)$ se conocen como la *e-transformada* de *Dyson* de los conjuntos $A$ y $B$). Observamos que $b \notin B(a^*)$, y en consecuencia $|B(a^*)| < |B|$. Veamos ahora una propiedad muy útil para nuestra prueba:

**Afirmación 1.5.** *Se cumple la igualdad* $|A|+|B| = |A(a^*)|+|B(a^*)|$.

Demostración la Afirmación: La demostración es la siguiente:

$$|A(a^*)|-|A| = |A(a^*) \setminus A| = |a^* + (B \setminus B(a^*))| = |B \setminus B(a^*)| = |B|-|B(a^*)|.$$

□

Esto implica que

$$|A(a^*)+B(a^*)| \leqslant |A+B| < |A|+|B|-1 = |A(a^*)|+|B(a^*)|-1,$$

donde la segunda desigualdad se debe a la hipótesis, y la última igualdad se deduce de la Afirmación 1.5. Hemos llegado pues a la contradicción deseada. □

Para concluir esta sección preliminar estudiaremos ahora el caso de un espacio vectorial de dimensión finita $n$ (pero arbitrariamente grande) sobre un cuerpo finito de característica $p$. Como veremos varias veces en este curso, el grupo $\mathbb{F}_p^n$ nos sirve a menudo como modelo que nos permitirá luego atacar problemas sobre la estructura de conjuntos de enteros (que pueden ser más complejos desde el punto de vista técnico).

Un ejercicio simple es demostrar que un conjunto $A \subseteq \mathbb{F}_p^n$ verifica $A+A = A$, y entonces $|A+A| = |A|$, si y sólo si $A$ es un subespacio (posiblemente afín) de $\mathbb{F}_p^n$. Un ejercicio un poco más avanzado muestra que incluso si relajamos la condición sobre el cardinal del conjunto suma, el conjunto $A$ sigue pareciéndose a un subespacio en el sentido siguiente:

**Ejercicio 1.6.** *Sea $p$ un número primo, y sea $A \subseteq \mathbb{F}_p^n$ un conjunto tal que $|A+A| \leqslant K|A|$ con $K < 3/2$. Demostrar que existe un subespacio (posiblemente afín) $V$ de $\mathbb{F}_p^n$ que contiene $A$ cumpliendo que $|A| \geqslant 2|V|/3$.*

En general, para cualquier grupo abeliano $G$ y subconjunto finito $A$ de $G$, denotamos por $K$ la menor constante tal que $|A+A| \leqslant K|A|$. Esta constante la llamamos la *constante de doblamiento* de $A$. En la siguiente sección consideramos el régimen en el que la constante de doblamiento es grande pero permanece acotada (en comparación con el tamaño del conjunto $A$ y el del grupo $G$, que tomamos tendiendo a infinito).

**1.1 Crecimiento de conjuntos suma.** El tema principal de esta capítulo es el análisis de conjuntos suma bajo la condición $|A + A| \leqslant K|A|$, donde $K$ es una constante. Veremos que en este caso podemos concluir que $A$ es un conjunto estructurado en un cierto sentido que describiremos con precisión más adelante. Concretamente, en la Sección 1.2 demostraremos que cuando $A \subseteq \mathbb{F}_p^n$ es un conjunto que cumple que $|A + A| \leqslant K|A|$, entonces tiene que estar contenido de manera "eficiente" en un subespacio de $\mathbb{F}_p^n$ cuyo cardinal está acotado en función de $K$.

Naturalmente, cuando un conjunto está contenido en un subespacio, su crecimiento bajo la operación "+" está limitado: el conjunto suma iterado llenará todo el subespacio y no puede crecer más. Antes de demostrar el resultado principal sobre los conjuntos cuyo conjunto suma está acotado, demostramos que el crecimiento de un tal conjunto bajo la operación "+" es bastante lento. Este fenómeno fue cuantificado por Plünnecke (1969) en el teorema siguiente, conocido como la *desigualdad de Plünnecke*. Recientemente, Petridis (2012) logró demostrar este teorema de una manera muy elemental, Y nosotros seguiremos sus pasos.[3]

**Teorema 1.7** (Desigualdad de Plünnecke)**.** *Sean $A, B \subseteq G$ conjuntos finitos tal que $|A + B| \leqslant K|A|$. Entonces para toda pareja de enteros $k, m \geqslant 1$*

$$|kB - mB| \leqslant K^{k+m}|A|.$$

Demostración: Sin pérdida de generalidad supongamos que $|A + B| = K|A|$. Elegimos un subconjunto no vacío $A' \subseteq A$ que minimice el cociente $|A' + B|/|A'|$, y escribimos $K' := |A' + B|/|A'|$. Observamos que $K' \leqslant K$, y que $|A' + B| = K'|A'|$ además de $|A'' + B| \geqslant K'|A''|$ para todo $A'' \subseteq A$.

Para continuar, supongamos que la siguiente afirmación es cierta (más tarde procederemos a demostrarla):

**Afirmación 1.8.** *Sean $A', B, K'$ como descritos anteriormente. Entonces para todo conjunto $C$, $|A' + B + C| \leqslant K'|A' + C|$.*

Terminamos primero la demostración del teorema suponiendo que la afirmación es cierta. En primer lugar comprobamos que para todo $m \in \mathbb{N}$,

$$|A' + mB| \leqslant K'^m|A'|.$$

Efectivamente, para $m = 1$ la desigualdad se cumple por la hipótesis. Para $m > 1$, se supone que la desigualdad se cumple para $m - 1$ y se sustituye $C = (m - 1)B$ en la afirmación 1.8 para obtener $|A' + mB| \leqslant K'|A' + (m - 1)B|$. Por la hipótesis inductiva, el término a la derecha es menor que $K'^m|A'|$.

El enunciado completo se deduce ahora a partir de la denominada *desigualdad de Ruzsa*: para todos los conjuntos $U$, $V$ y $W$, tenemos que $|U||V - W| \leqslant |U + V||U + W|$. (La validez de esta desigualdad se prueba fácilmente al definir una aplicación $\phi : U \times (V - W) \to (U + V) \times (U + W)$, que asigna a $(u, x = v - w)$ el elemento $(u + v, u + w)$, donde para todo $x \in V - W$ fijamos una representación $v - w$, de manera inyectiva). Esta desigualdad demuestra que

$$|A'||kB - mB| \leqslant |A' + kB||A' + mB| \leqslant K'^k|A'| \cdot K'^m|A'| \leqslant K^{k+m}|A'|^2,$$

[3] El argumento vale para todo grupo abeliano (y puede adaptarse a grupos no conmutativos), pero continuamos teniendo en mente al caso $G = \mathbb{F}_p^n$.

lo que implica $|kB - mB| \leqslant K^{k+m}|A'| \leqslant K^{k+m}|A|$. □

Demostración la Afirmación 1.8: Lo demostraremos por inducción sobre el cardinal del conjunto $C$. Cuando $|C| = 1$, la afirmación es trivial. Supongamos que el resultado es válido para $C$, y consideramos el conjunto $C' = C \cup \{x\}$. Observamos que

$$A' + B + C' = (A' + B + C) \cup [(A' + B + x) \setminus (D + B + x)],$$

donde $D$ es el conjunto $D := \{a \in A' : a + B + x \subseteq A' + B + C\}$. Pero la propiedad que definía la constante $K'$ implica que $|D + B| \geqslant K'|D|$, así que

$$|A' + B + C'| \leqslant |A' + B + C| + |A' + B| - |D + B| \leqslant K'(|A' + C| + |A'| - |D|). \tag{1-1}$$

Aplicamos el mismo argumento otra vez, escribiendo

$$A' + C' = (A' + C) \cup [(A' + x) \setminus (E + x)],$$

donde $E$ es el conjunto $E := \{a \in A' : a + x \in A' + C\}$ y la unión es disjunta. Concluimos que, como $E \subseteq D$,

$$|A' + C'| = |A' + C| + |A'| - |E| \geqslant |A' + C| + |A'| - |D|,$$

lo que junto a (1-1) implica el enunciado de la afirmación. □

**1.2 El teorema de Freiman–Ruzsa.** Una vez demostrada la desigualdad de Plünnecke, nos podemos embarcar a obtener el resultado estructural prometido. El teorema 1.9 es un resultado espectacular de Ruzsa (1999), que adaptó un teorema anterior de Freiman en el grupo de enteros Freĭman (1973) al grupo $\mathbb{F}_p^n$. Es uno de los casos mas convincentes para el uso del grupo modelo.

**Teorema 1.9** (Teorema de Freiman–Ruzsa). *Sea $A \subseteq G = \mathbb{F}_p^n$ un conjunto tal que $|A+A| \leqslant K|A|$. Entonces $A$ está contenido en un subespacio $H \leqslant \mathbb{F}_p^n$ de cardinal menor que $K^2 p^{K^4}|A|$.*

Demostración: La idea crucial (y ingeniosa) es elegir un subconjunto $X \subseteq 2A - A$ que sea máximo con la propiedad que las traslaciones $x + A$ para $x \in X$ sean disjuntas. En primer lugar mostramos que un conjunto $X$ que cumpla esta propiedad no puede ser demasiado grande. Efectivamente, tenemos $X + A \subseteq 3A - A$, y por la desigualdad de Plünnecke (el Teorema 1.7) sabemos que $|3A - A| \leqslant K^4|A|$. Ya que los conjuntos $x + A$ son disjuntos y de cardinal $|A|$ cada uno, obtenemos

$$K^4|A| \geqslant |3A - A| \geqslant |X + A| = \sum_{x \in X} |x + A| = |X||A|,$$

y en consecuencia $|X| \leqslant K^4$.

Enseguida comprobamos que

$$2A - A \subseteq X + (A - A).$$

Para ver esta afirmación, observamos que si $y \in 2A - A$, entonces $y + A \cap x + A \neq \varnothing$ para algún $x \in X$: si $y \in X$ el enunciado es trivial, y si $y \notin X$ se deduce de haber supuesto que con $X$ elegimos un conjunto máximo. En ambos casos concluimos que $y \in X + (A - A)$.

Sumando $A$ repetidamente a ambos lados de la inclusión antedicha, obtenemos

$$kA - A \subseteq (k-1)X + (A - A) \tag{1-2}$$

para todo $k \geqslant 2$. Así pues, hemos conseguido codificar cada vez más sumandos de $A$ dentro de pocas traslaciones de $A - A$ (el conjunto $X$ es de tamaño constante).

Denotemos por $H$ el subgrupo de $\mathbb{F}_p^n$ generado por $A$ y $Y$ para el subgrupo generado por $X$. Deducimos de (1-2) que

$$H = \bigcup_{k \geqslant 1} (kA - A) \subseteq Y + (A - A).$$

Pero todo elemento de $Y$ se puede escribir como suma de menos de $|X|$ elementos con coeficientes entre 1 y $p$, así que $|Y| \leqslant p^{|X|} \leqslant p^{K^4}$. La observación que

$$|H| \leqslant |Y||A - A| \leqslant K^2 p^{K^4} |A|$$

concluye la demostración. □

La dependencia del cardinal de $H$ sobre la constante de doblamiento $K$ se puede mejorar con más argumentos (en particular, utilizando las denominadas compresiones combinatorias, véase Green y Tao (2009), Even-Zohar (2012) y Even-Zohar y Lovett (2014)). Sin embargo, el ejemplo siguiente muestra que la dependencia debe ser de naturaleza exponencial.

**Ejemplo 1.10.** *Consideremos el conjunto $A \subseteq \mathbb{F}_p^n$, que consiste en una unión de un subespacio $H$ (muy grande) y $K - 1$ vectores elegidos al azar (no contenidos en $H$). Entonces la constante de doblamiento de $A$ es aproximadamente igual a $K$, pero todo subespacio $H'$ conteniendo $A$ debe ser de cardinal mayor que $p^{K-2}|A|$.*

Sin embargo, el conjunto en este ejemplo se considera también muy estructurado: aparte de un número constante de elementos, el conjunto $A$ está contenido en un subespacio de cardinal menor o igual que $|A|$ (a saber, el subespacio $H$ mismo). Esta observación nos invita a dar la reformulación siguiente del teorema de Freiman–Ruzsa.

**Teorema 1.11** (Teorema de Freiman–Ruzsa, reformulado)**.** *Sea $A \subseteq G = \mathbb{F}_p^n$ un conjunto tal que $|A + A| \leqslant K|A|$. Entonces existe un subespacio $H \leqslant \mathbb{F}_p^n$ de cardinal como máximo $C_1(K)|A|$ tal que para algún $x \in G$,*

$$|A \cap (x + H)| \geqslant \frac{|A|}{C_2(K)},$$

*donde $C_1(K)$ y $C_2(K)$ son constantes que dependen únicamente de $K$.*

Esta versión del teorema de Freiman–Ruzsa nos conduce a proponer la denominada *Conjetura Polinomial de Freiman–Ruzsa* (abreviada 'PFR' en inglés), que sigue siendo uno de los problemas abiertos fundamentales en combinatoria aditiva (para mas detalles véase la Sección 10 de Webb (2005)).

**Conjetura 1.12** (Conjetura polinomial de Freiman–Ruzsa)**.** *Las constantes $C_1(K)$ y $C_2(K)$ en el Teorema 1.11 dependen polinomialmente de $K$.*

El mejor resultado hasta la fecha se debe a Sanders (2012), quién consigue cotas casi óptimas para este problema.

El equivalente del teorema de Freiman–Ruzsa (el Teorema 1.11) en los enteros, demostrado por Freĭman (1973), dice que un conjunto finito de enteros cuyo conjunto suma es pequeño está contenido de manera eficiente en una progresión aritmética multidimensional.

**Teorema 1.13.** *Sea $A \subseteq \mathbb{Z}$ un conjunto finito de enteros tal que $|A+A| \leqslant C|A|$ para alguna constante $C$. Entonces $A$ está contenido en una traslación de una progresión aritmética multidimensional se la forma*

$$\{x \in \mathbb{Z} : x = \sum_{i=1}^{C_1} m_i x_i : m_i \in \mathbb{Z}, |m_i| \leqslant l_i\}$$

*para algunos $x_1, \ldots, x_{C_1} \in \mathbb{Z}$, que es de dimensión $C_1$ y de cardinal menor o igual que $C_2|A|$, donde $C_1$ y $C_2$ constantes que dependen únicamente de $C$.*

En ese contexto se puede formular una conjetura análoga a la Conjetura 1.12.

En la pasada década hemos visto también varias generalizaciones del Teorema 1.9 a otros grupos, incluso no conmutativos. Hasta ahora el resultado más general en esta dirección es Breuillard, Green y Tao (2012), y existen algunos artículos de revisión excelentes sobre este tema (ver Green (2014), Sanders (2013) y Helfgott (2015)).

## 2 Geometría de incidencia y el fenómeno suma-producto

**2.1 El teorema de Szemerédi–Trotter.** Para empezar trataremos de responder a una pregunta aparentemente inocente en el marco de la geometría euclidiana: dado un conjunto $P$ de puntos y un conjunto $L$ de rectas en el plano, ¿cuántas incidencias puede haber entre puntos de $P$ y rectas de $L$? Dicho de otro modo, buscamos una cota superior para la cantidad

$$I(P,L) := |\{(p,\ell) \in P \times L : p \in \ell\}|.$$

Trivialmente $I(P,L) \leqslant |P||L|$, pero esta desigualdad está muy lejos de ser óptima. El ejercicio siguiente mejora la cota.

**Ejercicio 2.1.** *Usar la desigualdad de Cauchy–Schwarz para demostrar que todo conjunto finito $P$ de puntos y $L$ de rectas en el plano cumplen*

$$I(P,L) \leqslant \min(|P|^{1/2}|L| + |P|, |L|^{1/2}|P| + |L|).$$

De hecho esta cota se puede mejorar. Para demostrarlo, vamos a utilizar argumentos procedentes de la teoría de grafos. Recuérdese que un *grafo* $G = (V,E)$ es un conjunto $V$ de *vértices* (muchas veces representados por puntos en el plano) junto con un conjunto $E \subseteq V \times V$ de *aristas* (representadas por líneas conectando dos vértices). El *número de cruce*, también llamado número de cruzamiento, de un grafo $G$ es el menor número de cruces de aristas en un diagrama plano del grafo $G$. Si el número de cruce de $G$ es 0, se dice que el grafo es *plano* (es decir, puede representarse sin cortes de aristas). En un grafo plano conexo [4] se cumple la *fórmula de Euler*, es decir

$$|F| - |E| + |V| = 2,$$

donde $F$ es el conjunto de "caras" del diagrama plano del grafo $G$ (es decir, cada una de las regiones 2-dimensionales en la que el dibujo del grafo divide el plano).

**Lema 2.2.** *Sea $G$ un grafo con $n$ vértices y $e$ aristas. Entonces el número de cruce de $G$ es mayor o igual a $e-(3n-6)$.*

DEMOSTRACIÓN: Suponemos que $H$ es un grafo plano, y sus conjuntos de vértices, aristas y caras se denotan por $V, E, F$, respectivamente. Podemos suponer también que $H$ es conexo, porque si no lo fuera podríamos añadir aristas hasta que lo sea, y ello no afecta a la demostración. Por la fórmula de Euler sabemos que

$$|F| - |E| + |V| = 2.$$

Entonces si $t$ es el número de pares $(f, e)$ tal que $f \in F, e \in E$ y la cara $f$ está limitada por la arista $e$, entonces $t \leqslant 2|E|$ and $t \geqslant 3|F|$, so $|E| \leqslant 3|V| - 6$. Concluimos que si $e > 3n - 6$, entonces el grafo $G$ no puede ser plano, y en consecuencia todo diagrama del grafo $G$ debe contener como mínimo un cruce. Quitamos una de las aristas de este cruce para obtener un grafo $G'$ con $n$ vértices y $e-1$ aristas. Podemos repetir este proceso $e-(3n-6)$ veces, así que $G$ tenía como mínimo $e-(3n-6)$ cruces. □

**Lema 2.3.** *Sea $G$ un grafo con $n$ vértices y $e \geqslant 4n$ aristas. Entonces el número de cruce de $G$ es mayor o igual a $e^3/(64n^2)$.*

DEMOSTRACIÓN: Supongamos que $G$ tiene un diagrama con $s$ cruces. Vamos a elegir un subgrafo $H$ de $G$ al azar, y contar el número de aristas y cruces inducidos en $H$.

Para obtener $H$, elegimos cada vértice de $G$ independientemente con probabilidad $p$ (un parámetro que fijaremos al final). Sea $H$ el subgrafo de $G$ inducido en este conjunto de vértices. Está claro que el número esperado de vértices de $H$ es $pn$, mientras el número esperado de aristas de $H$ es $p^2 e$ (ambos vértices conectados por la arista deben haber sido elegidos). De manera similar, el número esperado de cruces en $H$ es $p^4 s$, lo que implica, por el Lema 2.2, que

$$p^4 s \geqslant p^2 e - (3pn - 6) \geqslant p^2 e - 3n.$$

En consecuencia $s \geqslant e/p^2 - 3n/p^3$, y cuando escribimos $p := 4n/e < 1$ (par la hipótesis sobre $G$), obtenemos $s \geqslant e/(4n/e)^2 - 3n/(4n/e)^3 = e^3/(64n^2)$, que es precisamente la desigualdad deseada. □

La demostración precedente se debe a Ajtai, Chvátal, Newborn y Szemerédi, y independientemente Leighton. Ahora ya tenemos todos los ingredientes para demostrar una cota superior óptima para $I(P, L)$, conocida bajo el nombre *Teorema de Szemerédi–Trotter* Szemerédi y Trotter (1983), cuya demostración elegante se debe a Szekely.

[4]Un grafo $G$ se dice *conexo* si, para cualquier par de vértices $u$ y $v$ en $G$, existe al menos un camino (una sucesión de vértices adyacentes, conectados por aristas) de $u$ a $v$.

**Teorema 2.4** (Teorema de Szemerédi–Trotter). *Sea $P$ un conjunto finito de puntos en $\mathbb{R}^2$, y $L$ un conjunto finito de rectas. Entonces el número de incidencias entre $P$ y $L$ satisface*

$$I(P, L) \leqslant 4|L|^{2/3}|P|^{2/3} + 4|P| + |L|.$$

DEMOSTRACIÓN: Denotamos $|P| = n$ y $|L| = m$. Definimos un grafo $G$ con conjunto de vértices $P$ tal que $(p, p') \in P \times P$ es una arista si y sólo si los puntos $p$ y $p'$ definen una recta $\ell \in L$. Para todo $\ell \in L$, sea $s(\ell) := |\{p \in P : p \in \ell\}|$. Ahora

$$|E(G)| \geqslant \sum_{\ell \in L}(s(\ell) - 1) = I(P, L) - m,$$

y entonces si $I(P, L) - m \geqslant 4n$, el número de cruce de $G$ está acotado superiormente por $m^2$, y, por el Lema 2.3, también acotado inferiormente por $(I(P, L) - m)^3/(64n^2)$. Deducimos que

$$I(P, L) - m \leqslant (64n^2m^2)^{1/3} = 4n^{2/3}m^{2/3},$$

y así $I(P, L) \leqslant 4n^{2/3}m^{2/3} + m$. Para obtener el resultado del teorema, nos queda por observar que en el caso contrario al que acabamos de estudiar, el número de incidencias cumple que $I(P, L) < 4n + m$. □

**Ejercicio 2.5.** *Se denota por $[N]$ el conjunto de naturales $\{0, 1, 2, \ldots, N\}$. Considerando los elementos del retículo $[N] \times [2N^2]$, y el conjunto de todas las rectas con pendientes entre $1$ y $N$ que pasen por uno de los puntos del retículo, demostrar que la cota superior en el enunciado del Teorema 2.4 es óptima (a excepción de la constante).*

**Ejercicio 2.6.** *Generalizar el Teorema 2.4 en el sentido siguiente. Sea $P$ un conjunto de $n$ puntos en $\mathbb{R}^2$, y sea $L$ una familia de curvas simples de tamaño $m$ tal que cada par de curvas se intersecan en un máximo de $t$ puntos, y tal que cada par de puntos pertenece a un máximo de $s$ curvas. Entonces el número de incidencias entre $P$ y $L$ satisface*

$$I(P, L) \ll |L|^{2/3}|P|^{2/3}t^{1/3}s^{1/3} + s|P| + |L|.$$

Véase Pach y Sharir (1998), y Wang, Yang y Zhang (2013) para generalizaciones a curvas con $d$ grados de libertad.

**Ejercicio 2.7.** *Sea $P$ un conjunto finito de elementos en $\mathbb{R}^2$, y sea $k \geqslant 2$ un entero. Demostrar que el número de rectas conteniendo como mínimo $K$ puntos de $P$ está acotado superiormente por $O(|P|^2/k^3 + |P|/k)$.*

Como ejercicio suplementario, deducimos del Teorema 2.4 el siguiente resultado de Beck.

**Ejercicio 2.8.** *Demostrar que existen constantes $C_1$ y $C_2$ tal que, cuando $P$ es un conjunto finito de puntos en el plano, y $L_P$ es el conjunto de rectas generadas por $P$, una de las siguientes afirmaciones es cierta:*

- *Existe una recta $\ell \in L_P$ conteniendo como mínimo $C_1|P|$ puntos de $P$;*
- $|L_P| \geqslant C_2|P|^2$.

**2.2 El fenómeno suma-producto en los números reales.** En esta sección estudiamos una aplicación importante del Teorema de Szemerédi–Trotter en la combinatoria aritmética, que mostramos en primer lugar en el contexto de los números reales. Dado un conjunto finito $A \subseteq \mathbb{R}$, definimos como ya hemos hecho anteriormente el conjunto suma $A + A := \{a + a' : a, a' \in A\}$. De manera similar, podemos definir el *conjunto producto* $A \cdot A := \{aa' : a, a' \in A\}$, cuyo cardinal también varia entre $|A|$ y $|A|^2$. Nótese, sin embargo, que los casos de mínimo cardinal son muy distintos para el conjunto suma y el conjunto producto: el tamaño del primer conjunto se minimiza cuando $A$ es una progresión aritmética, pero en este caso el conjunto producto es muy grande. A la inversa, cuando $A$ es una progresión geométrica, el conjunto producto está minimizado, pero el conjunto suma alcanza el cardinal máximo. La siguiente pregunta es por lo tanto natural: ¿Es posible que el conjunto suma y el conjunto producto sean pequeños simultáneamente?

Efectivamente, la estructura aditiva y la estructura multiplicativa no parecen coexistir fácilmente. Con este propósito Erdős y Szemerédi (1983) conjeturaron lo siguiente.

**Conjetura 2.9** (Conjetura de Erdős–Szemerédi). *Sea $A \subseteq \mathbb{R}$ un conjunto finito. Entonces para todo $\epsilon > 0$,*

$$\max(|A + A|, |A \cdot A|) \gg |A|^{2-\epsilon}.$$

Dicho de otra manera, por lo menos uno de los dos conjuntos –suma o producto– debe alcanzar su cardinal máximo. En Erdős y Szemerédi (ibíd.) Erdős y Szemerédi demostraron un exponente de $1 + \epsilon$, y este exponente ha sido mejorado varias veces desde entonces. Aún así, la Conjetura 2.9 se mantiene como uno de los grandes problemas abiertos en combinatoria aritmética.

Comenzamos presentando un argumento simple pero muy elegante debido a Elekes (1997).

**Teorema 2.10.** *Sea $A \subseteq \mathbb{R}$ un conjunto finito. Entonces*

$$|A + A|^2 |A \cdot A|^2 \gg |A|^5,$$

*y en particular,*

$$\max(|A + A|, |A \cdot A|) \gg |A|^{5/4}.$$

DEMOSTRACIÓN: Consideramos el conjunto de puntos $P = (A + A) \times (A \cdot A)$, junto con el conjunto de rectas $L = \{y = a(x - b) : a, b \in A\}$. Inmediatamente observamos que $|P| = |A + A||A \cdot A|$ y que $|L| = |A|^2$. Cada recta de la forma $y = a(x - b)$ contiene como mínimo $|A|$ puntos de $P$, es decir, los de la forma $(b + a', aa')$ para $a' \in A$. Por lo tanto,

$$I(P, L) \geqslant |L||A| = |A|^3.$$

Pero por el Teorema de Szemerédi–Trotter (el Teorema 2.4) también tenemos la desigualdad

$$I(P, L) \ll |P|^{2/3}|L|^{2/3} + |P| + |L| = |A + A|^{2/3}|A \cdot A|^{2/3}|A|^{4/3} + |A + A||A \cdot A| + |A|^2.$$

Combinando las cotas se obtiene la desigualdad deseada. □

**Ejercicio 2.11.** *Modificando el argumento de Elekes, demostrar que para todos conjuntos finitos $A, B, C \subseteq \mathbb{R}$ tales que $|A| = |B| = |C|$,*

$$|AB||A + C| \gg \sqrt{|A|^3|B||C|}.$$

**Ejercicio 2.12.** *Considerando el conjunto de puntos $P = A(A+1) \times A(A+1)$ y el conjunto de rectas $L = \{y = ax/b + a : a \in A, b \in A+1\}$, demostrar que*

$$|A(A+1)| \gg |A|^{5/4}.$$

*(El exponente ha sido mejorado a $24/19$ por Jones y Roche-Newton (2013).)*

Algunos años más tarde la estimación suma-producto de Elekes fue mejorada por Solymosi (2005), logrando

$$|A+A|^8 |A \cdot A|^3 \gg \frac{|A|^{14}}{\log^3 |A|},$$

y en particular

$$\max(|A+A|, |A \cdot A|) \gg \frac{|A|^{14/11}}{\log^{3/11} |A|}.$$

**Ejercicio 2.13.** *Deducir de este resultado de Solymosi que cuando el cardinal del conjunto suma es mínimo (es decir, $|A+A| \ll |A|$), el cardinal del conjunto producto debe alcanzar casi el máximo (es decir, $|A \cdot A| \gg |A|^{2-o(1)}$). Nótese que este resultado no se puede deducir del Teorema 2.10, ni con la palabras "suma" y "producto" intercambiadas.*

**Ejercicio 2.14.** *Considerando el conjunto $A = [N]$, demostrar que el término $o(1)$ en el exponente del ejercicio precedente es necesario.*

En Solymosi (2009) mejoró aún más el exponente utilizando un argumento distinto, pero igualmente simple y elegante.

**Teorema 2.15.** *Sea $A \subseteq \mathbb{R}$ un conjunto finito. Entonces*

$$|A+A|^2 |A \cdot A| \gg \frac{|A|^4}{\log |A|},$$

*y en particular*

$$\max(|A+A|, |A \cdot A|) \gg \frac{|A|^{4/3}}{\log^{1/3} |A|}.$$

Demostración: La demostración usa la noción de la *energía multiplicativa*, definida por

$$E^{\times}(A) := |\{(a,b,c,d) \in A^4 : ad = bc\}|.$$

Nótese que trivialmente $E^{\times}(A) \leqslant |A|^3$. Comprobamos también la cota inferior siguiente.

**Afirmación 2.16.**

$$E^{\times}(A) \geqslant \frac{|A|^4}{|A \cdot A|}. \tag{2-1}$$

Demostración de la Afirmación 2.16: Para todo $s \in \mathbb{R}$, denotamos por $r_A(s)$ el número de representaciones de $s$ como $s = a \cdot a'$ con $a, a' \in A$. Entonces, por la desigualdad de Cauchy–Schwarz,

$$E^{\times}(A) = \sum_{s \in A \cdot A} r_A(s)^2 \geqslant \frac{(\sum_{s \in A \cdot A} r_A(s))^2}{|A \cdot A|},$$

y el simple hecho de que $|A|^2 = \sum_{s \in A \cdot A} r_A(s)$ completa la demostración. □

En vista de (2-1), esperaríamos poder acotar superiormente la energía multiplicativa por una cantidad que dependiera del cardinal del conjunto suma de $A$, lo que nos proporcionaría una estimación suma-producto en el sentido de la Conjetura 2.9.

Para ello, definamos primero, para cada $\lambda \in \mathbb{R}$, la función $r_A^*(\lambda)$ cuenta el número de representaciones de $\lambda$ como $\lambda = a/a'$, como siempre con $a, a' \in A$. Esto nos permite expresar la energía multiplicativa de $A$ como

$$E^\times(A) = \sum_{\lambda \in A/A} r_A^*(\lambda)^2 = \sum_{\lambda \in A/A} |\lambda \cdot A \cap A|^2,$$

donde $A/A := \{a/a' : a, a' \in A\}$ y $\lambda \cdot A := \{\lambda a : a \in A\}$. El logaritmo en la cota final tiene su origen en la división en intervalos diádicos que aplicamos al escribir

$$E^\times(A) = \sum_{\lambda \in A/A} |\lambda \cdot A \cap A|^2 = \sum_{i=0}^{\log|A|} \sum_{\lambda \in \Lambda_i} |\lambda \cdot A \cap A|^2,$$

donde $\Lambda_i := \{\lambda \in A/A : 2^i \leqslant |\lambda \cdot A \cap A| < 2^{i+1}\}$. Deducimos que existe un entero $i_0$ para el cual

$$\frac{E^\times(A)}{\log|A|} \leqslant \sum_{\lambda \in \Lambda_{i_0}} |\lambda \cdot A \cap A|^2. \tag{2-2}$$

Para facilitar la notación, enumeramos los elementos de $\Lambda_{i_0}$ en orden creciente, es decir $\lambda_1, \lambda_2, \ldots, \lambda_k$ para algún entero $k$, así que $E^\times(A)/\log|A| \leqslant k2^{2(i_0+1)}$. Llega el momento de considerar la geometría del problema. Por construcción, cada recta de la forma $y = \lambda_j x$, $j = 1, \ldots, k$, contiene entre $2^{i_0}$ y $2^{i_0+1}$ puntos de $A \times A$. Para todo $j = 1, \ldots, k$, denotamos el conjunto de estos puntos $P_j$. Dado que todo par de rectas distintas forma una base de $\mathbb{R}^2$, llegamos a la conclusión que el conjunto $P_j + P_{j'}$ es de cardinal $|P_j||P_{j'}| \geqslant 2^{2i_0}$ cuando $j \neq j'$. Además, para todos $j \neq j'$, los conjuntos suma $P_j + P_{j+1}$ y $P_{j'} + P_{j'+1}$ son disjuntos. Entonces

$$k2^{2i_0} \leqslant |\bigcup_{j=1}^{k} (P_j + P_{j+1})| \leqslant |(A \times A) + (A \times A)| = |A + A|^2,$$

pero de (2-2) y (2-1) tenemos también

$$k2^{2i_0+2} \geqslant \frac{E^\times(A)}{\log|A|} \geqslant \frac{|A|^4}{|A \cdot A| \log|A|},$$

lo que implica el resultado final. □

**2.3 El fenómeno suma-producto en cuerpos finitos.** Es natural plantear el mismo problema suma-producto en cuerpos que no sean $\mathbb{R}$. Nos concentramos en esta sección en el caso de cuerpos de característica positiva, en particular un primo $p$.

Para entender el porqué este problema es más complicado en $\mathbb{F}_p$, tomemos un conjunto $A \subseteq \mathbb{F}_p$ que consista en todo $\mathbb{F}_p$ excepto un conjunto pequeño de elementos (digamos 10). Observar entonces que un tal conjunto $A$ no puede crecer mucho ni bajo adición ni multiplicación, ya que está restringido por el cardinal del cuerpo ambiente. Para poder demostrar cualquier fenómeno suma-producto en un cuerpo finito, hace falta suponer que el cardinal de $A$ no sea del mismo orden que $p$.

Un primer paso importante en esta dirección, en parte debido a sus implicaciones para el problema de Kakeya que veremos en el Capítulo 3, fue el siguiente resultado de Bourgain, Katz y Tao (2004).

**Teorema 2.17.** *Sea $p$ un primo. Entonces para todo $\delta > 0$, existe $\epsilon = \epsilon(\delta) > 0$ con la propiedad siguiente. Sea $A \subseteq \mathbb{F}_p$ un conjunto de cardinal $p^\delta < |A| < p^{1-\delta}$. Entonces*

$$\max(|A + A|, |A \cdot A|) \gg |A|^{1+\epsilon}.$$

Poco después, Bourgain, Glibichuk y Konyagin (2006) quitaron la condición $|A| > p^\delta$. No daremos los detalles de la prueba (que es larga y técnica). Queda por remarcar que aunque la demostración del Teorema 2.17 no hacía referencia a la geometría de incidencia, Bourgain, Katz y Tao dedujeron un teorema de incidencia (el Corolario 2.20 más abajo) de su teorema suma-producto, mejorando la estimación primitiva en este caso.

**Ejercicio 2.18.** *Demostrar que la cota del Ejercicio 2.1 vale también en el plano finito $\mathbb{F}_p^2$.*

**Ejercicio 2.19.** *Demostrar que se puede cumplir la igualdad en la desigualdad del Ejercicio 2.18 cuando $|P| = p^2$.*

**Corolario 2.20.** *Sea $p$ un primo impar. Entonces para todo $\delta > 0$, existe $\epsilon = \epsilon(\delta) > 0$ con la propiedad siguiente. Para todo conjunto $P$ de puntos y $L$ de rectas en $\mathbb{F}_p^2$ de cardinal $|P|, |L| \leqslant N = p^\delta$ para algún $0 < \delta < 2$,*

$$I(P, L) \ll N^{3/2-\epsilon}.$$

La deducción del corolario a partir del Teorema 2.17 es tediosa, y le referimos al lector al artículo original Bourgain, Katz y Tao (2004).

Como ya lo hemos observado, subconjuntos grandes de cuerpos finitos no pueden crecer en el sentido de Erdős–Szemerédi porque están contenidos en un ambiente demasiado restrictivo. El ejemplo siguiente, que se debe a Bourgain (2005) y Chang (2008), proporciona una versión mas precisa de este enunciado.

**Ejemplo 2.21.** *Sea $g$ un generador de $\mathbb{F}_p^\times$ y sea $M$ un entero a determinar. Consideramos para todo entero $L \leqslant p$ el conjunto*

$$A_L := \{g^x : 1 \leqslant x \leqslant M\} \cap \{L+1, \ldots, L+M\} \subseteq \mathbb{F}_p.$$

*Por el principio del palomar, existe un entero $L$ tal que el cardinal de $A_L$ esta acotado inferiormente por $M/(p/M) = M^2/p$. Escribamos $A := A_L$. Entonces $A$ satisface la relación $|A+A| \ll M$ además de $|A \cdot A| \ll M$, así que denotando $M := \sqrt{p|A|}$ obtenemos un conjunto cuyo conjuntos* *suma y producto son de cardinal menor o igual que* $\sqrt{p|A|}$. *Entonces si* $|A| \gg p^{1/3}$, *no podemos esperar obtener una estimación del tipo*

$$\max(|A+A|, |A\cdot A|) \gg |A|^{2-\epsilon},$$

*como lo requiere la conjetura de Erdős–Szemerédi.*

El análogo ingenuo de la Conjetura 2.9 tiene que ser sustituido por la versión siguiente (explicitada por Garaev (2008), por ejemplo).

**Conjetura 2.22.** *Sea* $A \subseteq \mathbb{F}_p$. *Entonces*

$$\max(|A+A|, |A\cdot A|) \gg \min(|A|^{2-\epsilon}, |A|^{1/2} p^{1/2-\epsilon}).$$

En primer lugar nos concentramos en el caso en el que el cardinal de $A$ es grande, que está completamente resuelto. De hecho, es posible demostrar un teorema de incidencia en este contexto, de lo cual se deduce fácilmente el correspondiente teorema suma-producto.

**Teorema 2.23.** *Sea* $p$ *un primo impar. Sea* $P$ *un conjunto de puntos y* $L$ *un conjunto de rectas en* $\mathbb{F}_p^2$. *Entonces*

$$I(P,L) \leqslant \frac{|P||L|}{p} + \sqrt{p|P||L|}.$$

Demostración: Esta demostración se debe a Vinh (2011), y utiliza de nuevo la teoría de grafos. Para simplificar el argumento, trabajaremos en el espacio proyectivo $\mathbb{P}\mathbb{F}_p^2$, en el cual consideramos una inmersión de una copia de $\mathbb{F}_p^2$ de manera habitual (identificando un punto $x = (x_1, x_2)$ con la clase de equivalencia de $(x_1, x_2, 1)$, la que denotamos $[x]$). Definimos un grafo $G = (V, E)$ poniendo $V := \mathbb{P}\mathbb{F}_p^2$ y conectando $[x]$ y $[y]$ par una arista si y sólo si $x\cdot y := x_1y_1 + x_2y_2 + x_3y_3 = 0$. Dicho de otra manera, $[x]$ está conectado a $[y]$ si el punto representado por $[x]$ se encuentra en la recta representada por $[y]$, y viceversa. Nótese que $G$ tiene $n := p^2 + p + 1$ vértices y es regular de grado $d := p + 1$, con $d$ vértices de $G$ presentando bucles (porque el número de soluciones no nulas de $x_1^2 + x_2^2 + x_3^2 = 0$ sobre $\mathbb{F}_p$ es igual a $p^2 - 1$).

Calculamos los valores propios la matriz de adyacencia de $G$, es decir, la matriz $A = (a_{ij})_{i,j\in V}$ tal que $a_{ij} = 1$ si y sólo si hay una arista entre los vértices $i$ y $j$, y cero si no. Esta matriz es real y simétrica y por lo tanto admite una base ortonormal de vectores propios. Es bien conocido que el mayor valor propio de la matriz de adyacencia de un grafo regular es igual al grado de regularidad, y que está asociado con el vector propio $\mathbf{1} = (1, 1, \ldots, 1)$ normalizado. En otras palabras, el mayor valor propio de la matriz de adyacencia de nuestro grafo $G$ es $\lambda_0 := d$, correspondiente al vector propio $v_0 := \mathbf{1}/\sqrt{n}$.

Hay varias maneras de acotar la magnitud de los valores propios restantes. Una opción para comprobar que son todos de magnitud como máximo $\sqrt{d}$ se basa en estudiar la descomposición $A^2 = A^T A = J + (d-1)I$, donde $J$ es la matriz $n \times n$ que consiste sólo en unos, y $I$ es la matriz identidad de dimensión $n$. (Para ver eso, observamos que para distintos $[x]$ y $[y]$, existe precisamente un trayecto de longitud 2 entre ellos, lo que corresponde al hecho de que cada par de rectas se corta en un sólo punto). Dado que el rango de $J$ es 1 y que $J$ tiene sólo un valor propio no trivial (que corresponde al vector $v_0$), concluimos que los valores propios restantes son todos iguales a $d - 1$.

Es bien conocido también, y de hecho es el núcleo de la prueba, que todo grafo $d$-regular con $n$ vértices cuyos valores propios no triviales son pequeños se comporta aproximadamente como un grafo aleatorio del modelo $G_{n,d/n}$. Para no dejar nada en el tintero, enunciamos y demostramos el lema siguiente –muy estándar– con tal fin.

**Lema 2.24** (Lema de expansores)**.** *Sea $G = (V, E)$ un grafo $d$-regular con $n$ vértices. Sean $B, C$ subconjuntos de $V$. Entonces el número de aristas $e(B, C)$ entre $B$ y $C$ satisface*

$$\left| e(B,C) - \frac{d}{n}|B||C| \right| \leqslant \max_{\lambda \neq \lambda_0} |\lambda| \sqrt{|B||C|}.$$

Demostración del Lema 2.24: Escribimos $v_0, \ldots, v_n$ para la sucesión de vectores propios de la matriz de adyacencia $A$ de $G$, asociados a la sucesión (en orden no creciente) de valores propios $\lambda_0, \ldots, \lambda_n$. Denotamos $1_B$ el vector característico del conjunto $B$, es decir la coordenada $i$ de $1_B$ vale 1 si $i \in B$, y 0 si no. Desarrollamos $1_B = \sum_{k=0}^n b_k v_k$ y $1_C = \sum_{k=0}^n c_k v_k$ respecto a la base ortonormal $v_0, \ldots, v_n$. Observamos que $v_0 \cdot 1_B = b_0 = |B|/\sqrt{n}$, y de manera similar $v_0 \cdot 1_C = c_0 = |C|/\sqrt{n}$. Ahora

$$e(B,C) = 1_B^T A 1_C = \left( \sum_{k=0}^n b_k v_k \right)^T A \left( \sum_{j=0}^n c_j v_j \right) = \sum_{k=0}^n \lambda_k b_k c_k$$

por ortonormalidad de los vectores $v_i$. La discusión precediendo al lema implica que $\lambda_0 b_0 c_0 = d \cdot |B||C|/n$, y entonces

$$|e(B,C) - \frac{d}{n}|B||C|| \leqslant \max_{\lambda \neq \lambda_0} |\lambda| \sum_{k=1}^n b_k c_k \leqslant \max_{\lambda \neq \lambda_0} |\lambda| \left( \sum_{k=1}^n b_k^2 \right)^{1/2} \left( \sum_{k=1}^n c_k^2 \right)^{1/2}$$

como queríamos demostrar. □

El lema se aplica inmediatamente al conjunto $B := P$ de puntos y el conjunto $C := L$ de rectas, proporcionando

$$I(P, L) \leqslant \frac{p+1}{p^2+p+1}|P||L| + \sqrt{d|P||L|} \leqslant \frac{|P||L|}{p} + \sqrt{p|P||L|},$$

lo que completa la demostración del teorema. □

Por el mismo argumento que en el cuerpo de números reales, obtenemos una variante del teorema suma-producto a partir del teorema de incidencia (el Teorema 2.23).

**Corolario 2.25.** *Sea $p$ un primo impar, y sea $A \subseteq \mathbb{F}_p$. Entonces*

$$\max(|A+A|, |A \cdot A|) \gg \frac{2|A|^2}{p^{1/2} + (p + 4|A|^3/p)^{1/2}}.$$

**Ejercicio 2.26.** *Deducir el Corolario 2.25 a partir del Teorema 2.23, imitando la demostración del Teorema 2.10.*

**Ejercicio 2.27.** *Demostrar que el Corolario 2.25 implica los enunciados siguientes.*

- *Si* $p^{1/2} \ll |A| \leqslant p^{2/3}$*, entonces*

$$\max(|A+A|, |A\cdot A|) \gg \frac{|A|^2}{p^{1/2}}.$$

- *Si* $p^{2/3} \leqslant |A| \ll p$*, entonces*

$$\max(|A+A|, |A\cdot A|) \gg p^{1/2}|A|^{1/2}.$$

- *Si* $|A| \leqslant \sqrt{p}$*, no se puede decir nada que no sea trivial.*

**Ejercicio 2.28.** *Demostrar que el segundo enunciado del Ejercicio 2.27 es (esencialmente) óptimo.*

**Ejercicio 2.29.** *(Para los lectores que conocen la teoría básica de las sumas exponenciales.) Establecer una cota superior para la expresión*

$$\frac{1}{p}\sum_{n=0}^{p-1}\sum_{x\in A\cdot A}\sum_{a_1\in A}\sum_{a_2\in A}\sum_{y\in A+A} e(n(xa_1^{-1}+a_2-y)),$$

*donde* $e(x) := e^{2\pi i x/p}$*. Dar una demostración alternativa a los enunciados del Ejercicio 2.27.*

**Ejercicio 2.30.** *Generalizar el Corolario 2.25 al resultado siguiente: para todos conjuntos* $A, B, C \subseteq \mathbb{F}_p$*,*

$$|A+B||A\cdot C| \gg \min\left(p|A|, \frac{|A|^2|B||C|}{p}\right).$$

En el régimen en el que $A$ tiene tamaño pequeño, los métodos analíticos/espectrales fallan y se requieren herramientas combinatorias más sutiles. Los años recientes han visto una sucesión de avances en cuanto al exponente de crecimiento, empezando por Garaev (2007), que demostró que para todo $|A| < p^{7/13}(\log p)^{-4/13}$,

$$\max(|A+A|, |A\cdot A|) \gg |A|^{15/14-o(1)}.$$

En trabajos posteriores, Katz y Shen (2008) mejoraron el exponente a $14/13-o(1)$ para $|A| \leqslant p$, superado poco después por Bourgain y Garaev (2009) que consiguieron $13/12-o(1)$. Unos años más tarde, Rudnev (2012) obtuvo un exponente de $12/11-o(1)$ para $|A| \leqslant \sqrt{p}$. Por último, hacia un año, Roche-Newton, Rudnev y Shkredov (2016) obtuvieron el siguiente resultado notable:

**Teorema 2.31.** *Sea* $p$ *un primo y sea* $A \subseteq \mathbb{F}_p$*. Suponemos que* $|A| < p^{5/8}$*. Entonces*

$$\max(|A+A|, |A\cdot A|) \gg |A|^{6/5}.$$

Nótese que cuando $|A| \sim p^{5/8}$, la cota coincide con la del Corolario 2.25. El Teorema 2.31 se deduce del teorema de incidencia siguiente que se debe a Rudnev (2018).

**Teorema 2.32.** *Sea* $p$ *un primo,* $P$ *un conjunto de puntos y* $\Pi$ *un conjunto de planos en* $\mathbb{PF}_p^3$*. Sea* $\sigma$ *el máximo número de planos incidentes en una sola recta, y supongamos que* $|P| \geqslant |\Pi|$ *y que* $|\Pi| = O(p^2)$*. Entonces*

$$I(P,\Pi) \ll |P||\Pi|^{1/2} + \sigma|P|.$$

La demostración de este teorema de incidencia está basada en los métodos extraordinariamente poderosos de Guth y Katz, que no podemos abarcar en este curso. Aquí sólo mostramos como deducir un teorema de tipo suma-producto.

Demostración del Teorema 2.31 asumiendo el Teorema 2.32: Sea $A \subseteq \mathbb{F}_p$, y para facilitar la notación escribimos $B := A \cdot A$ y $C := A^{-1}$. El número de soluciones de la ecuación

$$a + bc = a + b'c' \tag{2-3}$$

con $a, a' \in A, b, b \in B, c, c' \in C$ se acota inferiormente por $E(A)|A|^2$, donde $E(A)$ es la *energía aditiva* definida por

$$E(A) := |\{(a_1, a_2, a_3, a_4) \in A^4 : a_1 + a_2 = a_3 + a_4\}|.$$

Pero toda solución de (2-3) corresponde a una incidencia entre un punto de $P := \{(a, c, b') : a \in A, b' \in B, c \in C\}$ y un plano de $\Pi := \{\pi : x + by - c'z = a', a' \in A, b \in B, c' \in C\}$. Cada uno de estos conjuntos es de cardinal $|\Pi| = |P| = |A||B||C| = |A|^2|A \cdot A|$, y el número máximo de planos colineales esta acotado por $\max(|A|, |B|, |C|) = \max(|A|, |A \cdot A|) = |A \cdot A|$. Observamos que la condición $|A| < p^{5/8}$ garantiza que $|\Pi| = |A|^2|A \cdot A| \ll p^{5/4}|A|^{6/5} < p^2$ siempre que $|A \cdot A| \ll |A|^{6/5}$, lo que evidentemente podemos suponer sin pérdida de generalidad. El teorema de incidencia de Rudnev (el Teorema 2.32) por lo tanto implica que el número de soluciones de (2-3) es

$$\ll (|A|^2|A \cdot A|)^{3/2} + |A|^2|A \cdot A|^2.$$

Es fácil de ver que el primer término en esta expresión es dominante. Además, la energía aditiva $E(A)$ se acota inferiormente por $|A|^4/|A+A|$ (imitando el argumento simple por Cauchy–Schwarz que hemos usado en la demostración del Teorema 2.15 para la energía multiplicativa), lo que resulta en la desigualdad

$$|A + A|^2|A \cdot A|^3 \gg |A|^6.$$

□

## 3 El método polinomial y el fenómeno Kakeya

**3.1 El método polinomial.** En esta sección introducimos el método polinomial, que se ha utilizado con gran éxito en el ámbito de la combinatoria aditiva, pero también en informática teórica (para un artículo de revisión reciente y completo, véase Tao (2014)). Empezamos por una generalización simple del teorema fundamental del álgebra sobre un cuerpo finito $\mathbb{F}$, que dice que todo polinomio en una variable de grado $d$ posee como máximo $d$ raíces. Dicho de otra manera, si $S \subseteq \mathbb{F}$ es un conjunto tal que $|S| > \mathrm{grad}(P)$, entonces existe $x \in S$ tal que $P(x) \neq 0$.

**Ejercicio 3.1** (Lema de Schwartz–Zippel). *Sea $f \in \mathbb{F}_p[x_1, \ldots, x_n]$ un polinomio no nulo de grado $d$. Entonces $f$ posee como máximo $dp^{n-1}$ raíces.*

Observamos que este lema implica inmediatamente que cuando el grado $d$ de un polinomio no nulo $f$ es estrictamente menor $p$, entonces $f$ posee $< p \cdot p^{n-1} = p^n$ raíces. En otras palabras, un polinomio no nulo de grado estrictamente menor que $p$ tiene que tomar un valor no cero sobre $\mathbb{F}_p^n$.

Una versión más elaborada de este principio es el famoso *Nullstellensatz combinatorio* de Alon (1999).

**Teorema 3.2** (Nullstellensatz combinatorio). *Sea $f \in \mathbb{F}_p[t_1, \ldots, t_n]$ un polinomio de grado $d$ con un coeficiente no nulo en $t_1^{d_1} \cdots t_n^{d_n}$, con $d_1 + \cdots + d_n = d$, y sean $S_1, \ldots, S_n \subseteq \mathbb{F}_p$ conjuntos tales que $|S_i| > d_i$ para todo $i = 1, 2, \ldots, n$. Entonces existen $x_1 \in S_1, \ldots, x_n \in S_n$ tales que $f(x_1, \ldots, x_n) \neq 0$.*

Demostración: La demostración se realiza por inducción sobre $n$. El caso $n = 1$ corresponde al teorema fundamental del álgebra. Supongamos que hemos demostrado el resultado para $n - 1 \geqslant 1$ variables.

Sea $g$ el polinomio

$$g(t_n) := \prod_{s_n \in S_n} (t_n - s_n),$$

en una variable de grado $|S_n|$. Podemos escribir $P(t_1, \ldots, t_n)$ como

$$P(t_1, \ldots, t_n) = q(t_1, \ldots, t_{n-1}) g(t_n) + r(t_1, \ldots, t_n),$$

donde $q$ es un polinomio de grado menor o igual a $d - |S_n|$, y $r$ es un polinomio de grado menor o igual a $d$. Escribimos

$$r(t_1, \ldots, t_n) = \sum_{j=0}^{|S_n|} r_j(t_1, \ldots, t_{n-1}) t_n^j,$$

y desarrollamos $qg$ como la suma de $q t_n^{|S_n|}$ y términos de orden inferior, cada uno de grado como máximo

$$(d - |S_n|) + (|S_n| - 1) < d = d_1 + \cdots + d_n.$$

Concluimos que los términos de orden inferior tienen un coeficiente nulo en $t_1^{d_1} \cdots t_n^{d_n}$. Dado que $|S_n| > d_n$, el polinomio $q t_n^{|S_n|}$ también tiene un coeficiente nulo en $t_1^{d_1} \cdots t_n^{d_n}$. Esto significa que el resto $r$ tiene que tener un coeficiente no nulo en $t_1^{d_1} \cdots t_n^{d_n}$, lo que implica en particular que $r_{d_n}$ tiene un coeficiente no nulo en $t_1^{d_1} \cdots t_{n-1}^{d_{n-1}}$.

Por la hipótesis inductiva, existen $x_1 \in S_1, \ldots, x_{n-1} \in S_{n-1}$ tales que $r_{d_n}(x_1, \ldots, x_{n-1}) \neq 0$. Pero el caso $n = 1$ aplicado a

$$r(t_n) := r(x_1, \ldots, x_{n-1}, t_n),$$

lo cual es un polinomio de grado $d_n$, nos permite encontrar también $x_n \in S_n$ tal que $r(x_n) \neq 0$. Ya que $g(x_n) = 0$, tenemos

$$P(x_1, \ldots, x_n) = q(x_1, \ldots, x_{n-1}) g(x_n) + r(x_1, \ldots, x_n) = r(x_1, \ldots, x_n) \neq 0,$$

la conclusión deseada. □

Como primera aplicación, damos una demostración alternativa del teorema de Cauchy–Davenport del capítulo 1.

Demostración alternativa del Teorema 1.4: El caso $|A| + |B| > p$ es trivial, supongamos entonces que $|A| + |B| \leqslant p$. Nuestro objetivo es demostrar que $|A + B| \geqslant |A| + |B| - 1$.

Supongamos, a fin de obtener una contradicción, que $|A + B| < |A| + |B| - 1$, es decir $K := |A + B| \leqslant |A| + |B| - 2$. Definimos el polinomio $P(t_1, t_2) := \prod_{m \in A+B} (t_1 + t_2 - m)$ de grado $K$.

Se observa que el coeficiente de $t_1^{|A|-1} t_2^{(K-|A|-1)}$ en $P$ es igual a $\binom{K}{|A|-1}$, lo cual es no nulo modulo $p$ puesto que $K \leqslant |A| + |B| - 2 < |A| + |B| \leqslant p$. Aplicamos el Teorema 3.2 con $d = d_1 + d_2$ y $d_1 = |A| - 1$, $d_2 = |B| - 1$, así que $A$ y $B$ son conjuntos tales que $|A| > d_1$ y $|B| > d_2$. Se deduce que existen $a \in A$ y $b \in B$ tales que $P(a, b) \neq 0$. Pero por construcción, $P(a', b') = 0$ para todo $a' \in A$ y todo $b' \in B$, lo que proporciona la contradicción que buscábamos. □

Existe una variante inocua del teorema de Cauchy–Davenport que resistía durante muchos años a toda tentativa de resolución; en particular, el método de la $e$-transformada del Capitulo 1 falla al intentarla aplicar en este caso que mostraremos. El Nullstellensatz combinatorio nos permite dar una solución a este problema (aunque fue resuelta por métodos distintos por da Silva y Hamidoune en 1994).

**Ejercicio 3.3.** *Dado un conjunto $A \subseteq \mathbb{F}_p$, definimos el* conjunto suma restringido $A \widehat{+} B := \{a + b : a \in A, b \in B, a \neq b\}$. *Demostrar que*

$$|A \widehat{+} B| \geqslant \min\{|A| + |B| - 3, p\}.$$

*Además, si $|A| \neq |B|$, entonces*

$$|A \widehat{+} B| \geqslant \min\{|A| + |B| - 2, p\}.$$

**3.2 El problema de Kakeya euclidiano.** El problema de Kakeya es uno de los grandes problemas abiertos de la matemáticas moderna, en parte porque tiene, a pesar de la simplicidad de su enunciado, importantes y profundas implicaciones en otros áreas, como por ejemplo en el análisis de ecuaciones en derivadas parciales, y en teoría de números. El objeto de estudio son los *conjuntos de Kakeya,* es decir conjuntos en un espacio euclídeo de dimensión $n$ que contengan un segmento unitario de línea en todas las direcciones.

**Definición 3.4** (Conjunto de Kakeya). *Un* conjunto de Kakeya *de dimensión $n$ es un subconjunto de $\mathbb{R}^n$ que contiene un segmento unitario de línea en todas las direcciones.*

La cuestión central es si se pueden construir conjuntos de Kakeya "pequeños", en un sentido que desarrollaremos en esta sección. El primer resultado en esta dirección es el famoso teorema de Besicovitch (1928), que afirma que existen conjuntos de Kakeya cuya medida de Lebesgue es nula. Referimos al lector a las notas de curso de Green (2003) para los detalles de la prueba.

**Teorema 3.5.** *Existe un conjunto de Kakeya de dimensión 2, cerrado y acotado, cuya medida de Lebesgue es nula.*

Mientras los conjuntos de Kakeya pueden ser muy pequeños en el sentido de la medida de Lebesgue, resulta que tienen un fuerte carácter bidimensional. Para dar una caracterización precisa, necesitamos la definición siguiente.

**Definición 3.6** (Dimensión de Minkowski). *Sea $B \subseteq \mathbb{R}^n$, y denotamos por $N_\delta(B) := \{x \in \mathbb{R}^n : \inf_{b \in B} |x - b| < \delta\}$ el $\delta$-entorno de $B$. Definimos la* dimensión de Minkowski inferior $\underline{d}(B)$ *por*

$$\underline{d}(B) := \inf\left\{d : \liminf_{\delta \to 0} \frac{|N_\delta(B)|}{\delta^{n-d}} = 0\right\},$$

*donde* $|\cdot|$ *es la medida de Lebesgue, y la* dimensión de Minkowski superior $\overline{d}(B)$ *por*

$$\overline{d}(B) := \inf\left\{d : \limsup_{\delta\to 0} \frac{|N_\delta(B)|}{\delta^{n-d}} = 0\right\}.$$

Hay otras nociones de dimensión que se han utilizado en este contexto, como por ejemplo la dimensión de Hausdorff. Nosotros nos centraremos en la primera.

**Ejemplo 3.7.** *El cuadrado* $S := [0,1]\times[0,1]\times\{0\}$*, como subconjunto de* $\mathbb{R}^3$*, tiene dimensión de Minkowski superior e inferior igual a* $\overline{d}(S) = \underline{d}(S) = 2$*. Para ver eso, se observa que* $2\delta \leqslant |N_\delta(S)| \leqslant 4\delta$*, así que*

$$\liminf_{\delta\to 0}\left(n - \frac{\log|N_\delta(S)|}{\log\delta}\right) \geqslant 2$$

*y*

$$\limsup_{\delta\to 0}\left(n - \frac{\log|N_\delta(S)|}{\log\delta}\right) \leqslant 2.$$

**Ejercicio 3.8.** *Calcular la dimensión de Minkowski superior e inferior del conjunto de Cantor*[5] *como subconjunto de* $\mathbb{R}$.

La célebre conjetura de Kakeya afirma que todo conjunto de Kakeya en $\mathbb{R}^n$ tiene dimensión de Minkowski $n$.

**Conjetura 3.9** (Conjetura de Kakeya)**.** *Sea* $d(n) := \inf_{B\subseteq\mathbb{R}^n} \overline{d}(B)$*, donde el ínfimo está definido sobre todos los conjuntos de Kakeya* $B \subseteq \mathbb{R}^n$*. Entonces*

$$d(n) = n.$$

Esta conjetura permanece abierta salvo el caso $n = 2$, en el cual fue resuelta por Davies.

**Teorema 3.10.** *Todo conjunto de Kakeya en el plano euclidiano es de dimensión de Minkowski superior igual a 2, es decir* $d(2) = 2$.

DEMOSTRACIÓN: Sea $B \subseteq \mathbb{R}^2$ un conjunto de Kakeya. Consideremos su $\delta$-entorno $N_\delta(B)$. Tenemos que demostrar que $|N_\delta(B)| \gg \delta^\epsilon$ para todo $\epsilon$. Dado que $B$ contiene un segmento unitario en cada dirección, $N_\delta(B)$ contiene un rectángulo de tamaño $\delta\times 1$ en cada dirección, y en particular en cada dirección definiendo un ángulo de $\pi j/2k$ con el eje $X$ positivo, donde $j = 1,\dots,k := \lfloor 1/\delta\rfloor$. Denotamos estos rectángulos $R_1,\dots,R_k$. Escribimos $1_{R_i}$ para la función indicatriz del rectángulo $i$, y $A := \cup_j R_j$. Por la desigualdad de Cauchy–Schwarz, obtenemos

$$k^2\delta^2 = \left(\int (1_{R_1} + \cdots + 1_{R_k})(x)dx\right)^2 \leqslant |A| \int (1_{R_1} + \cdots + 1_{R_k})(x)^2 dx = |A|\sum_{j,l} |R_j \cap R_l|. \tag{3-1}$$

[5] El conjunto de Cantor es el conjunto de todos los puntos del intervalo real $[0,1]$ que admiten una expresión en base 3 que no utilice el dígito 1.

Pero el área de intersección de cada pareja de rectángulos se calcula fácilmente: es igual a $\delta^2/\sin\theta$, donde $\theta := |j-l|\pi/2k$ es el ángulo entro los dos rectángulos $R_j$ y $R_l$. Dado que $k \leqslant 1/\delta$, tenemos $\delta^2/\sin\theta \leqslant 2\delta/|j-l|$, y pues para cada $j$,

$$\sum_l |R_j \cap R_l| \leqslant \delta + 2\sum_{s=1}^{k} \frac{2\delta}{s} \ll \delta \log(1/\delta).$$

Sumando todos los valores de $j$ y sustituyendo en (3-1) proporciona $|A| \gg 1/\log(1/\delta)$ y completa la demostración. □

Considerando de manera similar la intersección de tubos bien separados en dimensión $n \geqslant 3$, no es demasiado difícil demostrar que $d(n) \geqslant (n+1)/2$ (para los detalles, véase otra vez Green (2003)). Sin embargo, para $n$ grande esta cota está lejos del enunciado de la conjetura 3.9. Wolff (1999) la mejoró ligeramente a $d(n) \geqslant (n+2)/2$: veremos una versión simplificada de su argumento, sobre un cuerpo finito, en la sección siguiente. Utilizando métodos de la combinatoria aditiva, Bourgain (1999) obtuvo $d(n) \geqslant 13n/25$, y Katz y Tao (1999, 2002), refinando su método, obtuvieron el récord actual de $d(n) \geqslant (n-1)/\alpha + 1$, done $\alpha \approx 1.675$ es la mayor raíz del polinomio $x^3 - 4x + 2$. (Para una versión simplificada de este argumento, proporcionando la cota $d(n) \geqslant 4n/7$, véase el artículo de revisión de Dvir (2010)).

**3.3 El problema de Kakeya sobre un cuerpo finito.** Como ya vimos en el capítulo 1, muchas veces conviene transferir un problema difícil a un cuerpo finito. En el caso de Kakeya la cuestión sobre un cuerpo finito fue planteada por primera vez por Wolff (1999). Para empezar, aclaremos lo que significa un conjunto de Kakeya sobre un cuerpo finito.

**Definición 3.11** (Conjunto de Kakeya en $\mathbb{F}_p^n$)**.** *Un* conjunto de Kakeya *en* $\mathbb{F}_p^n$ *es un conjunto que contiene una recta en cada dirección.*

Una recta es simplemente un conjunto de la forma $\{x_0 + tx : t = 0, 1, \ldots, p-1\}$, así que la dirección $x$ de la recta está definida únicamente (hasta equivalencia proyectiva).

Como antes, nuestro objetivo es establecer la mínima dimensión de un conjunto de Kakeya.

**Definición 3.12** (Dimensión de Besicovitch)**.** *Definimos la* dimensión de Besicovitch $d_{\mathbb{F}_p}(n)$ *como el ínfimo de todo $d$ tal que existe una constante $C = C(d)$ y un conjunto de Kakeya en $\mathbb{F}_p^n$ de cardinal menor o igual a $Cp^d$.*

**Ejercicio 3.13.** *Considerando el número de incidencias entre puntos y rectas en $\mathbb{F}_p^n$, demostrar que*

$$d_{\mathbb{F}_p}(n) \geqslant \frac{1}{2}(n+1).$$

Después de alguna reflexión, llegamos a la conjetura siguiente.

**Conjetura 3.14** (Conjetura Kakeya en $\mathbb{F}_p^n$)**.** *Se cumple la igualdad*

$$d_{\mathbb{F}_p}(n) = n.$$

Resulta que en dimensión 2, es fácil de demostrar.

**Teorema 3.15.** *Todo conjunto de Kakeya en $\mathbb{F}_p^2$ es de cardinal mayor o igual a $p(p+1)/2$, así que*

$$d_{\mathbb{F}_p}(2) = 2.$$

Demostración: Ya hemos utilizado este tipo de argumentos en el caso euclidiano, pero aquí es todavía más fácil. Dado que el conjunto de Kakeya $A \subseteq \mathbb{F}_p^2$ contiene una recta en cada una de las $p+1$ direcciones, se puede suponer sin pérdida de generalidad que $A$ es la unión de estas rectas, denotadas $\ell_1, \dots, \ell_{p+1}$. Obsérvese que

$$p^2(p+1)^2 = \left( \sum_{x \in \mathbb{F}_p^2} (1_{\ell_1} + \cdots + 1_{\ell_{p+1}})(x) \right)^2 \leqslant$$

$$\leqslant |A| \sum_{x \in \mathbb{F}_p^2} (1_{\ell_1} + \cdots + 1_{\ell_{p+1}})(x)^2 = |A| \sum_{i,j} |\ell_i \cap \ell_j|.$$

Pero cada pareja de rectas se interseca en precisamente un punto, excepto cuando son iguales, en cual caso el cardinal de la intersección es igual a $p$. Concluimos que

$$p^2(p+1)^2 \leqslant 2|A|p(p+1),$$

lo que completa la prueba. □

De hecho, el ejercicio siguiente muestra que la constante $1/2$ en el Teorema 3.15 es óptima.

**Ejercicio 3.16.** *Sea $p > 2$. Considerando el conjunto*

$$S := \{(x,t) \in \mathbb{F}_p^2 : x + t^2 \text{ es un cuadrado en } \mathbb{F}_p\},$$

*demostrar que existe un conjunto de Kakeya en $\mathbb{F}_p^2$ de cardinal menor o igual a $p(p+3)/2$.*

Presentamos otro argumento que mejora la cota básica del Ejercicio 3.13. Este argumento, conocido bajo el nombre *argumento de cepillo de Wolff*, fue importante en el contexto euclidiano (véase la sección 3.2).

**Teorema 3.17.** *Se cumple la desigualdad*

$$d_{\mathbb{F}_p}(n) \geqslant \frac{1}{2}(n+2).$$

*Más precisamente, todo conjunto de Kakeya en $\mathbb{F}_p^n$ es de cardinal mayor o igual a*

$$\frac{1}{8} p^{(n+2)/2}.$$

Antes de empezar, planteamos un ejercicio en dimensión 2 cuyo resultado utilizamos más tarde.

**Ejercicio 3.18.** *Sea $A \subseteq \mathbb{F}_p^2$ una unión de $k$ rectas en direcciones distintas. Demostrar que $|A| \geqslant pk/2$.*

Demostración del Teorema 3.17: Supongamos que $A \subseteq \mathbb{F}_p^n$ es un conjunto de Kakeya de tamaño $|A| \leqslant p^{(n+2)/2}$, o sea , supongamos que $A$ consiste en una unión de $k$ rectas, con $p^{n-1} \leqslant k \leqslant 2p^{n-1}$. Afirmamos que existe alguna recta que se corta con un mínimo de $p^{n/2}/4$ rectas. Para ver eso, denotamos las $k$ rectas $\ell_1, \ldots, \ell_k$, y observamos que

$$p^2k^2 = \left( \sum_{x \in \mathbb{F}_p^n} (1_{\ell_1} + \cdots + 1_{\ell_k})(x) \right)^2 \leqslant |A| \sum_{x \in \mathbb{F}_p^n} (1_{\ell_1} + \cdots + 1_{\ell_k})(x)^2,$$

lo que es igual a

$$|A| \sum_{i,j=1}^{k} |\ell_i \cap \ell_j| \leqslant p^{(n+2)/2} \sum_{i,j=1}^{k} |\ell_i \cap \ell_j|.$$

Deducimos que existe $1 \leqslant i \leqslant k$ tal que

$$\sum_{j \neq i} |\ell_i \cap \ell_j| \geqslant \frac{1}{2} p^{n/2} - p \geqslant \frac{1}{4} p^{n/2}.$$

A partir de ahora nos concentramos en $\ell := \ell_i$ y el conjunto de rectas que se cortan con $\ell$, las que enumeramos $\ell'_1, \ldots, \ell'_m$, donde $m \geqslant p^{n/2}/4$. Nos referiremos a esta colección de rectas como el *cepillo*, y lo denotamos $H$.

Cada una de las rectas $\ell'_1, \ldots, \ell'_m$ está ubicada en un único plano de dimensión 2 que contiene también $\ell$. Trabajaremos con la colección completa de estos planos $(\Pi_j)_{j=1,\ldots,t}$, y supondremos que cada plano $\Pi_j$ contiene $\beta_j \geqslant 1$ de las líneas $\ell'_1, \ldots, \ell'_m$. Del Ejercicio 3.18 sabemos que para todo $j = 1, \ldots, t$,

$$|\Pi_j \cap H| \geqslant \frac{1}{2}(\beta_j + 1)p \geqslant (\frac{1}{4}\beta_j + 1)p.$$

En particular, se ve que todo plano $\Pi_j$ contiene como mínimo $\beta_j p/4$ puntos de $H$ que no están en $\ell$, y estas colecciones de puntos son disjuntas cuando $j$ varia entre 1 y $m$. Concluimos que

$$|A| \geqslant |H| \geqslant \frac{1}{4} p \sum_{j=1}^{t} \beta_j = \frac{1}{4} pm \geqslant \frac{1}{8} p^{(n+2)/2},$$

lo que completa la demostración. □

En dimensión 3, la cota inferior $5/2$ para la dimensión de un conjunto Kakeya permanecía el récord durante mucho tiempo. El trabajo de Bourgain, Katz y Tao (2004), que ya mencionamos en la sección 2.3, la mejoró por un épsilon (a $5/2 + \epsilon$!) Otros artículos, que no tenemos tiempo para detallar aquí pero que estuvieron basados en el método de Katz y Tao (1999), mejoraron la cota a

$$d_{\mathbb{F}_p}(n) \geqslant \frac{1}{7}(4n + 3).$$

Fue en este punto que Dvir (2009) anunció una resolución completa de la conjetura 3.14, utilizando un método completamente distinto, a saber el método polinomial. Antes de dar la prueba (sorprendentemente corta), introducimos la noción del *conjunto de Nikodým*.

**Definición 3.19** (Conjunto de Nikodým)**.** *Un conjunto $N \subseteq \mathbb{F}_p^n$ se dice* conjunto de Nikodým *si para todo $x \notin N$, existe una recta por $x$ que interseca $N$ en todo punto salvo uno. En otras palabras, para todo $x \notin N$, existe $y \in \mathbb{F}_p^n$ tal que $\{x + ty : t = 1, \dots, p-1\} \subseteq N$.*

Está claro que los conjuntos de Nikodým están relacionados con los conjuntos de Kakeya.

**Ejercicio 3.20.** *Sea $B \subseteq \mathbb{F}_p^n$ un conjunto de Kakeya. Demostrar que existe un conjunto de Nikodým $N \subseteq \mathbb{F}_p^n$ de cardinal menor o igual a $p|B|$.*

¿Cuál es el propósito de definir los conjuntos de Nikodým? Supongamos que existiera un conjunto de Kakeya $B$ (y pues un conjunto de Nikodým $N$) pequeño en $\mathbb{F}_p^n$. Si pudiésemos encontrar un polinomio $f$ de grado $d < p-1$ cuyos valores en $N$ se conocieran, entonces podríamos recuperar los valores de $f$ en todo punto de $\mathbb{F}_p^n$. En particular, si $f$ fuese no trivial pero nulo en $N$, entonces podríamos concluir que $f$ es nulo en todo el espacio. Pero esta conclusión estaría en contradicción con el Lema de Schwartz–Zippel, según el cual un tal polinomio tiene que tomar al menos un valor no cero.

**Teorema 3.21.** *Sea $B \subseteq \mathbb{F}_p^n$ un conjunto de Kakeya. Entonces $|B| \gg p^{n-1}/(n!)$.*

DEMOSTRACIÓN: Sea $d := p-2$, y supongamos, a fin de obtener una contradicción, que

$$|B| < \binom{n-1+d}{n-1}.$$

Entonces el número de monomios en $\mathbb{F}_p[x_1, \dots, x_n]$ de grado $d$, igual al número de asignaciones no negativas $j_1, \dots, j_n$ tales que $j_1 + \dots + j_n = d$, es estrictamente mayor que el cardinal de $B$. Se deduce que existe un polinomio homogéneo no trivial $g \in \mathbb{F}_p[x_1, \dots, x_n]$ de grado $d$ que se anula en todo punto de $B$. Dado que $g$ es homogéneo, $g$ se anula también en el conjunto de Nikodým $N := \{tx : x \in B\}$ asociado con $B$.

Consideramos ahora un punto $x \notin N$ y la recta perforada $\ell^* := \{x + ty : t = 1, \dots, p-1\}$ que para algún $y \in \mathbb{F}_p^n$ está contenida en $N$. Restringido a $\ell$, el polinomio $g$ es un polinomio de grado $d$ que se anula en $p-1 > d$ puntos, así que $g$ es nulo en todo punto de $\ell := \{x + ty : t \in \mathbb{F}_p\}$; en particular, $g(x) = 0$. Pero $x$ fue arbitrario, lo que implica que $g$ se anula en todo $\mathbb{F}_p^n$, contradiciendo el Lema de Schwartz–Zippel (el Lema 3.1) según el cual un polinomio no trivial $g$ de grado $d$ puede poseer como máximo $dp^{n-1}$ raíces. □

**Ejercicio 3.22.** *Sea $B \subseteq \mathbb{F}_p^n$ un conjunto Kakeya. Formando un producto cartesiano iterado, demostrar que $|B| \gg_{\epsilon,n} p^{n-\epsilon}$.*

Poco después de la primera publicación del trabajo de Dvir, Alon y Tao remarcaron que el argumento de Dvir se podría adaptar para resolver la conjetura completa, con un exponente de $n$ en vez de $n-1$. Lo dejamos como ejercicio (véase el teorema 3 del artículo original de Dvir (ibíd.)).

**Ejercicio 3.23.** *Sea $B \subseteq \mathbb{F}_p^n$ un conjunto Kakeya. Entonces $|B| \geqslant p^n/(n!)$.*

Esta cota ha sido mejorada a $p^n/2^n$ Dvir, Kopparty, Saraf y Sudan (2009). La construcción siguiente, que se debe a Mockenhaupt y Tao (2004) para $n = 2$, y Saraf y Sudan (2008) para $n > 2$ en el caso de característica impar y par, respectivamente, muestra que esta cota es óptima (salvo el factor de 2).

**Ejercicio 3.24.** *Sea $p$ un primo impar. Considerando el conjunto*

$$B' := \{(x_1^2/4 + x_1 t, \dots x_{n-1}^2/4 + x_{n-1}t, t) : x_1, \dots, x_{n-1}, t \in \mathbb{F}_p\},$$

*demostrar que existe un conjunto de Kakeya de cardinal menor o igual a*

$$\frac{p^n}{2^{n-1}} + O(p^{n-1}).$$

**Agradecimientos.** Estoy muy agradecida a Juanjo Rué por haber revisado este texto.

## Referencias

Julia Wolf
julia.wolf@dpmms.cam.ac.uk
School of Mathematics
University of Bristol
Bristol BS8 1TW
United Kingdom

# CONJUNTOS DE SIDON

Javier Cilleruelo†

## Índice general

## 1 Prefacio

Uno de los temas favoritos de Paul Erdős y que mejor describe su gusto por los problemas aritméticos con sabor combinatorio, ha sido el de los conjuntos de Sidon. Corría el año 1932 cuando Simon Sidon, analista húngaro, le preguntó a Erdős sobre el crecimiento de sucesiones de enteros positivos con la propiedad de que todas las sumas de dos elementos de la sucesión son distintas.

Estos conjuntos, que Erdős llamaría más tarde conjuntos de Sidon, son el objeto de este curso. Aunque el interés de Sidon por estos conjuntos tenía que ver con cuestiones del análisis de Fourier, el problema cautivó a un joven Erdős por su vertiente aritmética y combinatoria y se convertiría en un tema recurrente en su investigación hasta que nos abandonara en 1998 en busca de "El Libro", ese libro virtual donde Erdős afirmaba que se encuentran las demostraciones más elegantes e ingeniosas que jamás hayan sido escritas.

Son muchos los problemas que nos podemos plantear acerca de los conjuntos de Sidon. Casi todos ellos tienen que ver con el tamaño máximo que pueden llegar a tener estos conjuntos en un intervalo o un grupo finito dado, y en el caso de las sucesiones infinitas, con la construcción de sucesiones infinitas de Sidon $A$ cuya función contadora $A(x) = |A \cap [1, x]|$ crezca tanto como sea posible. Estos dos problemas, tratados en el primero y en el segundo capítulo respectivamente, serán los objetivos principales del curso.

Hemos dedicado un último capítulo a problemas sin resolver sobre los conjuntos de Sidon donde hemos aprovechado para comentar otros problemas sobre conjuntos de Sidon no tratados en los dos primeros capítulos. Es posible que el lector más ambicioso no resista la tentación de empezar por este capítulo y utilizar los capítulos anteriores y las referencias para ir completando información.

Los dos primeros capítulos están acompañados de ejercicios que ayudarán al lector a afianzar los contenidos del curso. Algunos son sencillos o consisten en completar algunos detalles de alguna demostración. Otros tienen mayor dificultad y permiten al lector explorar nuevos territorios.

Las notas de este curso son una parte del libro “Conjuntos de Sidon” Cilleruelo (2014a), que fue escrito para el curso que con el mismo título, fue impartido en Mérida (Venezuela) en septiembre de 2014. El capítulo II del libro “Sequences” de Halberstam y Roth (1983) es una referencia clásica sobre los conjuntos de Sidon hasta 1966. Erdős y Freud (1991) escribieron un survey muy completo hasta 1991, pero en húngaro. “A complete annotated bibliography of work related to Sidon sequences”, de O’Bryant (2004), es también una fuente de información valiosa sobre conjuntos de Sidon.

## 2 Conjuntos de Sidon finitos

**2.1 Los orígenes.** Erdős solía contar la siguiente anécdota referida a Simon Sidon. Una de las tardes que él y su amigo Paul Turán fueron a visitar al analista húngaro, Sidon abrió bruscamente la puerta y les gritó: *vengan en otro momento y busquen a otra persona*. No sería esa tarde sino otra cuando Simon Sidon despertó el interés de Erdős al preguntarle por los conjuntos de enteros positivos con la propiedad de que todas las sumas de dos elementos de la sucesión son distintas.

Aunque el interés de Sidon por estos conjuntos tenía que ver con cuestiones del análisis de Fourier, el problema cautivó a un joven Erdős por su vertiente aritmética y combinatoria y se convertiría en un tema recurrente en su investigación hasta que nos abandonara en 1996 en busca de “El Libro”, ese libro virtual donde Erdős afirmaba que se encuentran las demostraciones más ingeniosas y elegantes que jamás hayan sido escritas. Fue el propio Erdős quien bautizó con el nombre de conjuntos de Sidon a estos conjuntos. Una definición más formal y más general de un conjunto de Sidon es la siguiente.

**Definición 2.1.** *Sea $G$ un grupo conmutativo. Un conjunto $A \subset G$ es un conjunto de Sidon si*

$$a + b = c + d \implies \{a, b\} = \{c, d\}$$

*para todo $a, b, c, d \in A$.*

En otras palabras, $A$ es un conjunto de Sidon si todas las sumas de dos elementos de $A$ son distintas salvo por el orden de presentación de los sumandos.

Como $a+b = c+d$ si y sólo si $a-c = d-b$, los conjuntos de Sidon también se definen, indistintamente, como aquellos con la propiedad de que todas las diferencias no nulas de sus elementos son distintas.

Habitualmente utilizaremos el término conjunto de Sidon cuando nos refiramos a un conjunto finito y sucesión de Sidon cuando éste sea infinito. En este capítulo empezaremos estudiando los conjuntos de Sidon finitos y dejaremos para el siguiente las sucesiones de Sidon infinitas.

El conjunto siguiente es un conjunto de Sidon:

$$A = \{1,\ 2,\ 5,\ 10,\ 16,\ 23,\ 33,\ 35\}.$$

Una manera de organizar el cálculo de todas las diferencias positivas de dos elementos del conjunto para comprobar que efectivamente es un conjunto de Sidon, consiste en construir el triángulo de diferencias como se muestra a continuación. En la primera fila y en negrita hemos dispuesto los elementos del conjunto y en las filas inferiores todas las diferencias de dos elementos del conjunto que provienen de los elementos en negrita de las dos diagonales correspondientes. Nótese que todas las diferencias son distintas, así que el conjunto $A$ es, efectivamente, un conjunto de Sidon. Más adelante veremos que es de tamaño maximal; es decir, el intervalo $[1, 35]$ no contiene ningún conjunto de Sidon de 9 elementos.

| | | | | | | | | | | | | | | |
|---|---|---|---|---|---|---|---|---|---|---|---|---|---|---|
| **1** | | **2** | | **5** | | **10** | | **16** | | **23** | | **33** | | **35** |
| | 1 | | 3 | | 5 | | 6 | | 7 | | 10 | | 2 | |
| | | 4 | | 8 | | 11 | | 13 | | 17 | | 12 | | |
| | | | 9 | | 14 | | 18 | | 23 | | 19 | | | |
| | | | | 15 | | 21 | | 28 | | 25 | | | | |
| | | | | | 22 | | 31 | | 30 | | | | | |
| | | | | | | 32 | | 33 | | | | | | |
| | | | | | | | 34 | | | | | | | |

Parece totalmente improbable que exista una fórmula sencilla que nos proporcione, de manera exacta, el mayor tamaño de un conjunto de Sidon en el intervalo $[1, n]$. Sin embargo sí que sabemos dar una respuesta asintótica a este problema. Se tratará en la siguiente sección pero antes dedicaremos unas pocas líneas a la notación que utilizaremos a lo largo del curso.

Se recuerda que $f(x) = O(g(x))$ significa que existe una constante positiva $C$ tal que $f(x) < Cg(x)$ para $x$ suficientemente grande. Se utiliza principalmente cuando en una estimación hay un término principal y un término de error del que sólo nos interesa el orden de magnitud. Algunas veces se utiliza también la notación $f(x) \ll g(x)$ para indicar lo mismo. Ésta última se utiliza sobre todo cuando del término principal sólo nos interesa el orden de magnitud.

La notación $f(x) = o(g(x))$ indica, por el contrario, que $\lim_{x\to\infty} \frac{f(x)}{g(x)} = 0$. Así que, por ejemplo, los términos $o(1)$ se refieren siempre a cantidades que tienden a cero cuando $x$ tiende a infinito.

**2.2 Conjuntos de Sidon en intervalos.** ¿Cuál es el mayor número de elementos que puede tener un conjunto de Sidon en el intervalo $\{1, \ldots, n\}$?

Esta es una pregunta básica sobre los conjuntos de Sidon a la que se sabe dar una respuesta asintótica. Utilizaremos la notación clásica (ver Halberstam y Roth (1983), capítulo II)

$$F_2(n) = \max\{|A| : \ A \subset \{1, \ldots, n\}, \ A \text{ es Sidon}\}.$$

Un sencillo argumento de conteo proporciona una primera cota superior para esta cantidad. Como todas las diferencias positivas $a - a', \ a, a' \in A$ son distintas y menores que $n$ y hay exactamente $\binom{|A|}{2}$ de esas diferencias, se tiene la desigualdad $\binom{|A|}{2} \leqslant n - 1$, de la que se sigue la cota trivial

$$F_2(n) < \sqrt{2n} + 1/2. \tag{1}$$

Esta cota superior ya permite ver que, como decíamos en la sección anterior, un conjunto de Sidon en el intervalo $[1, 35]$ no puede tener más de 8 elementos:

$$F_2(35) < \sqrt{70} + 1/2 = 8.866\ldots$$

La cota superior (1) está lejos de ser una cota óptima cuando $n$ es grande. Pero se puede mejorar si en lugar de tener en cuenta todas las diferencias $a-a'$, $a, a' \in A$, se consideran sólo las diferencias pequeñas. De esta manera Erdős y Turán (1941) demostraron la desigualdad

$$F_2(n) < \sqrt{n} + O(n^{1/4}).$$

Años más tarde Lindström (1969a) precisó más el término de error con una demostración muy ingeniosa que reproducimos aquí.

**Teorema 2.1** (Lindström)**.**

$$F_2(n) < n^{1/2} + n^{1/4} + 1.$$

*Prueba.* Sean $1 \leqslant a_1 < \cdots < a_k \leqslant n$ los elementos de un conjunto de Sidon en $[1, n]$. Dado un $r$, que elegiremos al final, las siguientes desigualdades son claras:

$$\begin{aligned} (a_2 - a_1) + (a_3 - a_2) + \cdots + (a_k - a_{k-1}) &= a_k - a_1 < n \\ (a_3 - a_1) + (a_4 - a_2) + \cdots + (a_k - a_{k-2}) &= a_k + a_{k-1} - (a_1 + a_2) < 2n \\ &\cdots \\ (a_{r+1} - a_1) + (a_{r+2} - a_2) + \cdots + (a_k - a_{k-r}) &= (a_k + \cdots + a_{k-(r+1)}) - (a_1 + \cdots + a_r) < rn \end{aligned} \tag{2}$$

En la parte izquierda aparecen exactamente

$$(k-1) + (k-2) + \cdots + (k-r) = rk - r(r+1)/2$$

diferencias positivas de la forma $a_j - a_i$ y todas ellas son distintas por ser $\{a_1, \ldots, a_k\}$ un conjunto de Sidon. Así que si $l = rk - r(r+1)/2$ y llamamos $S_r$ a la suma de todas esas diferencias tenemos que

$$S_r = \sum_{\substack{i,j \\ 1 \leqslant i < j \leqslant i+r}} (a_j - a_i) \geqslant \sum_{n=1}^{l} n > \frac{l^2}{2} = (rk - r(r+1)/2)^2/2.$$

Por otra parte, las desigualdades de la parte derecha en (2) implican que

$$S_r < n + (2n) + \cdots + (rn) = nr(r+1)/2.$$

De las dos desigualdades sobre $S_r$ obtenemos que

$$\begin{aligned} k &< \sqrt{n(1 + 1/r)} + (r+1)/2 \\ &< \sqrt{n} + \frac{\sqrt{n}}{2r} + \frac{r+1}{2}. \end{aligned}$$

La elección de $r = \lceil n^{1/4} \rceil$ finaliza la demostración. □

Vamos a dar otra demostración distinta, más moderna e inspiradora, que es esencialmente la que dio Ruzsa (1993). Antes de proceder vamos a introducir algunas definiciones y notaciones que son habituales en la teoría combinatoria de números.

Sea $G$ un grupo conmutativo. Dados dos subconjuntos $A, B \subset G$, definimos el conjunto suma

$$A + B = \{a + b,\ a \in A,\ b \in B\}$$

y la función

$$r_{A+B}(x) = |\{(a,b) \in A \times B,\ a + b = x\}|,$$

que cuenta el número de representaciones de $x$ como suma de un elemento de $A$ y otro de $B$. A lo largo del curso haremos uso de las identidades triviales

$$r_{A-A}(0) = |A| \tag{3}$$

y

$$\sum_{x \in G} r_{A+B}(x) = |A||B|. \tag{4}$$

La cantidad

$$\sum_x r^2_{A+B}(x)$$

se denomina *energía aditiva* entre $A$ y $B$ y cuenta el número de soluciones de la ecuación $a + b = a' + b'$ con $a, a' \in A$ y $b, b' \in B$, que coincide con el número de soluciones de la ecuación $a - a' = b' - b$ con $a, a' \in A$ y $b, b' \in B$. Esta observación da lugar a la identidad

$$\sum_x r^2_{A+B}(x) = \sum_x r_{A-A}(x) r_{B-B}(x), \tag{5}$$

que aparecerá insistentemente a lo largo de este curso.

El siguiente lema se debe a Ruzsa (1993).

**Lema 2.1** (Ruzsa (ibíd.))**.** *Sea $A$ un conjunto de Sidon en un grupo conmutativo $G$ y sea $B$ cualquier subconjunto de $G$. Entonces se tiene que*

$$|A|^2 \leqslant |A + B| \left(1 + \frac{|A| - 1}{|B|}\right). \tag{6}$$

*Prueba.* Haciendo uso de la identidad (4), de la desigualdad de Cauchy–Schwarz y de (5) obtenemos

$$\begin{aligned} (|A||B|)^2 &= \left(\sum_{x \in A+B} r_{A+B}(x)\right)^2 \\ &\leqslant |A + B| \sum_x r^2_{A+B}(x) \\ &= |A + B| \sum_x r_{A-A}(x) r_{B-B}(x). \end{aligned} \tag{7}$$

Como $r_{A-A}(x) \leqslant 1$ para $x \neq 0$, tenemos que

$$\begin{aligned}\sum_x r_{A-A}(x) r_{B-B}(x) &= r_{A-A}(0) r_{B-B}(0) + \sum_{x\neq 0} r_{A-A}(x) r_{B-B}(x)\\ &\leqslant r_{A-A}(0) r_{B-B}(0) + \sum_{x\neq 0} r_{B-B}(x)\\ &= |A||B| + |B|^2 - |B|. \qquad (8)\end{aligned}$$

De (7) y (8) se sigue la desigualdad

$$(|A||B|)^2 \leqslant |A+B| \big(|A||B| + |B|^2 - |B|\big)$$

y la demostración del lema. □

El Lema 2.1 fue utilizado por Ruzsa para dar una demostración alternativa de la desigualdad de Lindström en el Teorema 2.1. Una pequeña modificación técnica permite una ligera mejora en la desigualdad.

**Teorema 2.2** (Cilleruelo (2010b) y Ruzsa (1993))**.** *Si $A \subset [1, n]$ es un conjunto de Sidon, entonces*

$$|A| < n^{1/2} + n^{1/4} + 1/2.$$

*Prueba.* Consideremos el conjunto $B = [0, l] \cap \mathbb{Z}$ donde

$$l = \lfloor \sqrt{n(|A|-1)} \rfloor.$$

Entonces $|A+B| \leqslant n + l$ y $|B| = l + 1$. Así que el Lema 2.1 implica la desigualdad

$$\begin{aligned}|A|^2 &\leqslant (n+l)\left(1 + \frac{|A|-1}{l+1}\right)\\ &< n + l + \frac{n(|A|-1)}{l+1} + |A| - 1\\ &\leqslant n + 2\sqrt{n(|A|-1)} + |A| - 1\\ &= (\sqrt{n} + \sqrt{|A|-1})^2,\end{aligned}$$

que una sencilla manipulación conduce a la desigualdad

$$(|A| - \sqrt{n})^2 < |A| - 1. \qquad (9)$$

Escribiendo $|A| = \sqrt{n} + cn^{1/4} + 1/2$ y sustituyendo esta expresión en (9) obtenemos

$$c^2 n^{1/2} + cn^{1/4} + 1/4 < n^{1/2} + cn^{1/4} - 1/2,$$

que da lugar a una contradicción cuando $c \geqslant 1$. □

**Ejercicio 2.1** (Cilleruelo (2010b))**.** *Demostrar que si $A \subset \{1, \ldots, n\}$ es tal que $r_{A-A}(x) \leqslant g$ para todo $x \neq 0$, entonces*

$$|A| \leqslant (gn)^{1/2} + (gn)^{1/4} + 1/2.$$

Es interesante observar que todas las demostraciones que se conocen de la estimación $F_2(n) < \sqrt{n} + O(n^{3/4})$ sólo utilizan el hecho de que todas las diferencias menores que $n^{3/4}$ son distintas. El siguiente ejercicio ilustra todavía mejor ese hecho.

**Ejercicio 2.2.** *Sea $\omega(x)$ una función que tiende a infinito y sea $A \subset [1,n]$ un conjunto de enteros positivos tal que $r_{A-A}(x) \leqslant 1$ para todo $x$, $1 \leqslant x \leqslant \omega(n)\sqrt{n}$. Demostrar que entonces $|A| \leqslant (1+o(1))\sqrt{n}$.*

La cota superior

$$F_2(n) < n^{1/2} + n^{1/4} + 1/2$$

parece ser el límite del método consistente en contar diferencias pequeñas. Aunque durante mucho tiempo Erdős conjeturó que $F_2(n) < \sqrt{n} + O(1)$, acabó afirmando Erdős (1989) que la conjetura era demasiado optimista y que la conjetura correcta debería ser $F_2(n) < \sqrt{n} + O(n^{\epsilon})$ para todo $\epsilon > 0$. El Ejercicio 2.15 apoya también la creencia de que

$$\limsup_{n\to\infty}(F_2(n) - \sqrt{n}) = \infty, \tag{10}$$

pero todavía no se ha encontrado una demostración de este hecho.

El Teorema 2.2 también puede deducirse a partir de la siguiente desigualdad que tiene su propio interés.

**Teorema 2.3.** *Si $A \subset [1,n]$ entonces*

$$|A| < \sqrt{n} + \sqrt{|A|^2 - |A-A|}.$$

Dejamos al lector la demostración de este teorema en el siguiente ejercicio.

**Ejercicio 2.3.** *Sean $A$ y $B$ dos subconjuntos de un grupo conmutativo $G$. Demostrar que*

$$|A|^2 \leqslant |A+B|\left(1 + \frac{|A|^2 - |A-A|}{|B|}\right).$$

*y utilizar esta desigualdad para demostrar el Teorema 2.3.*

**Ejercicio 2.4.** *Dar una demostración alternativa del Teorema 2.2 a partir del Teorema 2.3.*

**Ejercicio 2.5.** *Sea $A \subset [1,n]$ un conjunto de enteros positivos tal que $|A-A| = (1+o(1))|A|^2$. Demostrar que $|A| \leqslant (1+o(1))\sqrt{n}$.*

El Ejercicio 2.5 muestra que la condición de ser de Sidon no es estrictamente necesaria para obtener la cota superior $|A| \leqslant (1+o(1))\sqrt{n}$. Es suficiente con que $|A-A| \sim |A|^2$.

En contra de lo que se pudiera sospechar no se llega a la misma conclusión si asumimos que $|A+A| \sim |A|^2/2$. Erdős y Freud (1991) construyeron, para cada $n$, un conjunto $A \subset [1,n]$ con $|A| \sim \frac{2}{\sqrt{3}}\sqrt{n} = (1.154\ldots)\sqrt{n}$ elementos y tal que $|A+A| \sim |A|^2/2$.

Consideraron un conjunto de Sidon $B$ de tamaño máximo en $[1, n/3)$ que, como veremos en la próxima sección, tiene $\sim \sqrt{n/3}$ elementos. El conjunto $A = B \cup (n-B)$ tiene entonces $\sim \frac{2}{\sqrt{3}}\sqrt{n}$ elementos y es fácil comprobar que todas las sumas de dos elementos de $A$ con distintas excepto las $|B|$ sumas de la forma $b + (n-b)$. Dejamos los detalles como un ejercicio.

**Ejercicio 2.6.** *Demostrar que el conjunto* $A \subset [1, n]$ *construido por Erdős y Freud satisface que* $|A + A| \sim |A|^2/2$ *y tiene* $|A| \sim \frac{2}{\sqrt{3}}\sqrt{n}$ *elementos.*

El ejercicio siguiente da una estimación trivial en la otra dirección.

**Ejercicio 2.7.** *Sea* $A \subset [1, n]$ *un conjunto de enteros positivos tal que* $|A+A| = (1+o(1))|A|^2/2$*. Demostrar que* $|A| \leqslant (1 + o(1))2\sqrt{n}$*.*

O. Pikhurko (2006) ha demostrado que si $A \subset [1, n]$ es tal que $|A + A| = (1 + o(1))|A|^2/2$ entonces $|A| \leqslant (1.863\dots)\sqrt{n}$.

**2.3 Conjuntos de Sidon en grupos conmutativos finitos.** Para obtener buenas cotas inferiores para $F_2(n)$ necesitamos construir conjuntos de Sidon en $\{1, \dots, n\}$ tan grandes como sea posible. Las mejores construcciones conocidas provienen de construcciones algebraicas en grupos cíclicos. Es claro que si $A \subset \mathbb{Z}_n$ es un conjunto de Sidon en $\mathbb{Z}_n$ también lo es en el intervalo $[1, n]$. Sin embargo el recíproco no es cierto. Es sencillo comprobar que el conjunto $\{1, 2, 5, 10, 16, 23, 33, 35\}$, que era un conjunto de Sidon en el intervalo $[1, 35]$, no lo es, sin embargo, en $\mathbb{Z}_{35}$.

Sea $A$ un conjunto de Sidon en un grupo conmutativo finito $G$. Como todas las diferencias $a - a'$, $a \neq a' \in A$ son no nulas y distintas, se tiene la desigualdad trivial $|A|(|A| - 1) \leqslant |G| - 1$, que implica la cota superior

$$|A| \leqslant \sqrt{|G| - 3/4} + 1/2. \tag{11}$$

Es decir, si llamamos $F_2(G)$ al máximo cardinal de un conjunto de Sidon en $G$, siempre se tiene que

$$F_2(G) \leqslant \left\lfloor \sqrt{|G| - 3/4} + 1/2 \right\rfloor. \tag{12}$$

Esta cota superior es óptima para algunas familias infinitas de grupos finitos. El ejemplo más sencillo para el que se alcanza esta cota, y que es un caso particular del Ejemplo 1 que aparece más adelante, es la parábola en $\mathbb{Z}_p \times \mathbb{Z}_p$ cuando $p$ es un primo impar:

$$A = \{(x, x^2) : x \in \mathbb{Z}_p\} \subset G = \mathbb{Z}_p \times \mathbb{Z}_p.$$

**Ejercicio 2.8.** *Sea* $q$ *una potencia de un primo impar. Demostrar que el conjunto* $A = \{(x, x^2) : x \in \mathbb{F}_q\}$ *es un conjunto de Sidon en* $\mathbb{F}_q \times \mathbb{F}_q$*. Comprobar que para esta familia de grupos se alcanza la cota superior* (12)*.*

Probablemente Erdős y Turán (1941) se inspiraron en este conjunto para construir el primer ejemplo de un conjunto de Sidon $A$ en $\{1, \dots, n\}$ con $|A| \gg \sqrt{n}$.

**Ejercicio 2.9** (Erdős y Turán). *Demostrar que para todo* $p$ *primo, el conjunto*

$$A = \{(x^2)_p + 2xp : 0 \leqslant x \leqslant p - 1\}$$

*es un conjunto de Sidon en* $\{1, \dots, 2p^2\}$ *con* $p$ *elementos. Deducir de esto que* $F_2(n) \geqslant \sqrt{n/2}(1 + o(1))$*.*

A continuación describimos otras familias de conjuntos de Sidon de tamaño máximo en sus grupos ambientes correspondientes. En lo que sigue $q$ indicará un primo o la potencia de un primo.

**Ejemplo 1.** *Sea $q$ una potencia de un primo impar y sean $r(x), s(x) \in \mathbb{F}_q[X]$ polinomios de grado menor o igual que $2$, linealmente independientes en $\mathbb{F}_q$. Entonces, el conjunto*

$$A = \{(r(x), s(x)) : \ x \in \mathbb{F}_q\}$$

*es un conjunto de Sidon en $\mathbb{F}_q \times \mathbb{F}_q$ con $q$ elementos.*

**Ejercicio 2.10.** *Demostrar que los conjuntos del Ejemplo 1 son conjuntos de Sidon en $\mathbb{F}_q \times \mathbb{F}_q$.*

**Ejercicio 2.11.** *Sea $q = 2^n$. Demostrar que el conjunto*

$$A = \{(x, x^3) : \ x \in \mathbb{F}_q\}$$

*es un conjunto de Sidon en $\mathbb{F}_q \times \mathbb{F}_q$.*

**Ejemplo 2.** *Para todo generador $g$ de $\mathbb{F}_q^*$, el conjunto*

$$A = \{(x, g^x) : \ x \in \mathbb{Z}_{q-1}\} \tag{13}$$

*es un conjunto de Sidon en $\mathbb{Z}_{q-1} \times \mathbb{F}_q$ con $q-1$ elementos. Este conjunto también se puede describir de la forma*

$$A = \{(\log x, x) : \ x \in \mathbb{F}_q^*\}$$

*donde* $\log x = \log_g x$ *es el logaritmo discreto en base $g$.*

Para probar que $A$ es un conjunto de Sidon, tenemos que ver que dado un elemento $(a, b) \in \mathbb{Z}_{q-1} \times \mathbb{F}_q$, $(a, b) \neq (0, 0)$, la igualdad

$$(x, g^x) - (y, g^y) = (a, b) \tag{14}$$

determina los valores $x, y$. La igualdad (14) se puede escribir de la forma

$$\begin{aligned} x - y &\equiv a \pmod{q-1} \\ g^x - g^y &\equiv b \pmod{q}. \end{aligned}$$

De la primera ecuación se tiene que $g^x = g^{y+a}$, que sustituido en la segunda da lugar a la ecuación $g^y(g^a - 1) = b$. Si $a = 0$ entonces $b = 0$, en contra de lo supuesto. Si $a \neq 0$, entonces $g^a - 1 \neq 0$ y el valor de $y$ queda determinado y a su vez el de $x$.

**Ejemplo 3.** *Dados dos generadores $g_1, g_2$ de $\mathbb{F}_q^*$, el conjunto*

$$A = \{(x, y) \in \mathbb{Z}_{q-1} \times \mathbb{Z}_{q-1} : \ g_1^x + g_2^y = 1\} \tag{15}$$

*es un conjunto de Sidon en $\mathbb{Z}_{q-1} \times \mathbb{Z}_{q-1}$ con $q-2$ elementos.*

**Ejercicio 2.12.** *Demostrar que el conjunto del Ejemplo 3 es un conjunto de Sidon.*

Cuando $q = p$ es un primo, podemos identificar $\mathbb{F}_p$ con $\mathbb{Z}_p$. Las figuras de abajo corresponden a los conjuntos de Sidon descritos para el caso $p = 17$ .

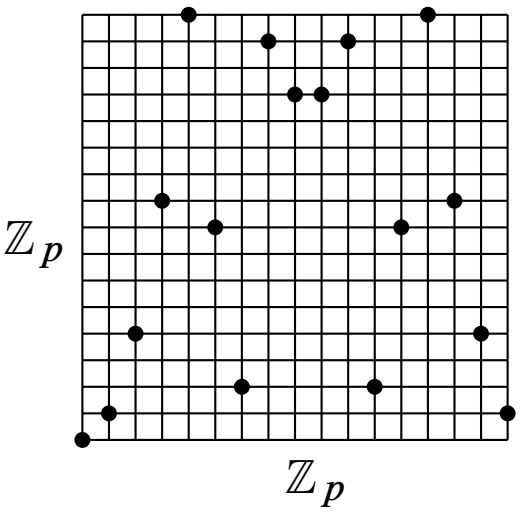

Ejemplo 1

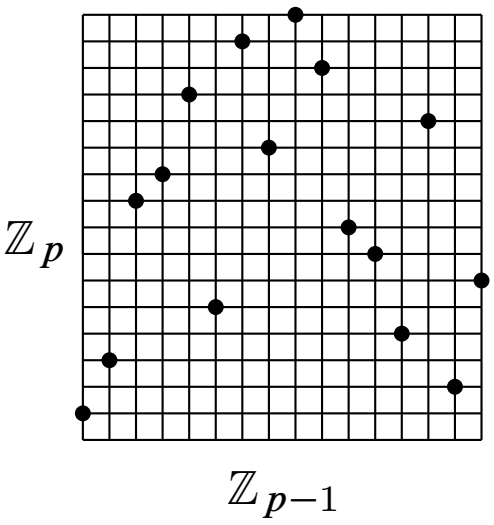

Ejemplo 2

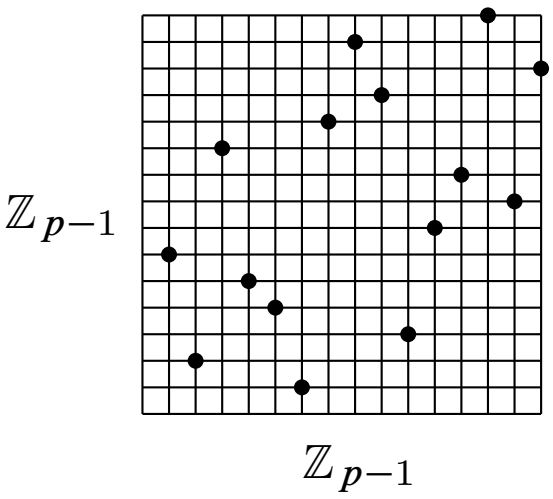

Ejemplo 3

Observar que los conjuntos de Sidon de estos ejemplos, con $q$, $q-1$ y $q-2$ elementos respectivamente, tienen cardinal maximal (los dos primeros) o casi (el último).

**Ejercicio 2.13.** Sea $\phi : G \to G'$ un isomorfismo entre los grupos $G$ y $G'$. Demostrar que si $A$ es un conjunto de Sidon en $G$, entonces el conjunto $\phi(A) = \{\phi(a) : a \in A\}$ es un conjunto de Sidon en $G'$.

Como muestra el Ejercicio 2.13, los isomorfismos entre grupos preservan la propiedad de ser de Sidon. Así que la imagen del conjunto de Sidon descrito en el Ejemplo 2 por el isomorfismo natural entre $\mathbb{Z}_{p-1} \times \mathbb{Z}_p$ y $\mathbb{Z}_{p(p-1)}$ es un conjunto de Sidon en $\mathbb{Z}_{p(p-1)}$ de $p-1$ elementos. Esta observación se debe a Ruzsa (1993) y proporciona la construcción más sencilla de conjuntos de Sidon de tamaño maximal en grupos cíclicos

**Proposición 2.1** (Ruzsa). *Sea $g$ un generador de $\mathbb{F}_p^*$. El conjunto*

$$A = \{px - (p-1)(g^x)_p : 0 \leqslant x \leqslant p-2\}$$

*es un conjunto de Sidon en $\mathbb{Z}_{p(p-1)}$.*

Se deja al lector la demostración de la Proposición 2.1 en el siguiente ejercicio.

**Ejercicio 2.14.** *Demostrar que el conjunto $A$ de la Proposición 2.1 es la imagen del conjunto del Ejemplo 2 bajo el isomorfismo natural entre $\mathbb{Z}_{p-1} \times \mathbb{Z}_p$ y $\mathbb{Z}_{(p-1)p}$, y que por lo tanto el conjunto $A$ es un conjunto de Sidon en $\mathbb{Z}_{(p-1)p}$.*

La siguiente reflexión y el Ejercicio 2.15 invitan a pensar que la conjetura (10) debería ser cierta. Si se eligen al azar $p-1$ elementos de $\mathbb{F}_p$ con repetición (un mismo elemento puede ser elegido varias veces), con alta probabilidad ocurriría que algunos elementos serían elegidos más de una vez, otros ninguna e incluso habría un intervalo de longitud aproximadamente $\log p$ cuyos elementos no habrían sido elegidos. Es decir, habría lagunas de longitud $\sim \log p$. Se piensa que la sucesión $g^x - x \pmod p$, $x = 1, \dots, p-1$ se comporta como una sucesión aleatoria de este tipo, pero ni siquiera se ha podido demostrar la siguiente conjetura.

**Conjetura 2.1.** *Para todo $M$ existe un primo $p$ y un generador $g$ de $\mathbb{F}_p^*$ tal que la sucesión $g^x - x \pmod p$, $x = 0, \dots, p-1$ no contiene ningún elemento en algún intervalo $I$ de longitud $M$.*

**Ejercicio 2.15.** *Utilizar el Teorema 2.1 para demostrar que la Conjetura 2.1 implicaría que*

$$\limsup_{n\to\infty}(F_2(n)-\sqrt{n})=\infty.$$

*Es decir, que la conjetura de Erdős, $F_2(n)<\sqrt{n}+O(1)$, no sería cierta.*

Los conjuntos de Sidon construidos en la Proposición 2.1 permiten obtener una buena cota inferior para $F_2(n)$.

**Teorema 2.4.** *Sea $\theta$ con la propiedad de que para todo $m$ suficientemente grande, el intervalo $(m, m+m^\theta)$ contiene algún primo. Entonces*

$$F_2(n)\geqslant n^{1/2}+O(n^{\theta/2}).$$

*Prueba.* Dado $n$, sea $p$ el mayor primo tal que $p(p-1)\leqslant n$. Si $n$ es suficientemente grande podemos encontrar un tal primo $p$ con $p>n^{1/2}-2n^{\theta/2}$. Es claro que el conjunto construido en la Proposición 2.1 es, en particular, un conjunto de Sidon en $\{1,\dots,p(p-1)\}$ y por lo tanto en $\{1,\dots n\}$ □

Se conjetura que en el Teorema 2.4 es posible tomar cualquier $\theta>0$ pero sólo se sabe cierto Baker, Harman y Pintz (2001) para $\theta\geqslant 0.525$.

De los Teoremas 2.2 y 2.4 se deduce la estimación asintótica para $F_2(n)$.

**Corolario 2.1.** $F_2(n)\sim\sqrt{n}$.

Las conjeturas sobre las cotas superior y la cota inferior se pueden resumir en la siguiente.

**Conjetura 2.2.** $F_2(n)=\sqrt{n}+O(n^\epsilon)$ *para todo* $\epsilon>0$.

Existen otras familias de grupos cíclicos que contienen conjuntos de Sidon que alcanzan la cota superior (12). La primera de ellas fue encontrada por Singer (1938). Utilizando geometría proyectiva finita construyó un conjunto de Sidon $A$ de $q+1$ elementos en $\mathbb{Z}_{q^2+q+1}$. Nótese que el tamaño de su conjunto diferencia es $|A-A|=|A|^2-|A|+1=q^2+q+1$. Es decir todo elemento no nulo del grupo se escribe, de manera única, como diferencia de dos elementos de $A$. Los conjuntos con esta propiedad se denominan conjuntos de diferencias perfectas y fue en este contexto donde este conjunto fue encontrado. Pasaron algunos años hasta que esta construcción fuera conocida por Erdős y popularizada en el mundo de los conjuntos de Sidon. Posteriormente Bose (1942) construyeron un conjunto de Sidon de $q$ elementos en $\mathbb{Z}_{q^2-1}$. Estas construcciones, los Ejemplos 1, 2, 3 y la cota superior en (11) prueban el valor exacto de $F_2(G)$ para algunas familias de grupos:

$$\begin{aligned}
F_2(\mathbb{F}_q\times\mathbb{F}_q)&=q\\
F_2(\mathbb{F}_q\times\mathbb{Z}_{q-1})&=q-1\\
F_2(\mathbb{Z}_{q^2+q+1})&=q+1\\
F_2(\mathbb{Z}_{q^2-1})&=q\\
F_2(\mathbb{Z}_{q-1}\times\mathbb{Z}_{q-1})&\in\{q-2,q-1\}.
\end{aligned}$$

Observar que si $q = p^n$, entonces el grupo aditivo en $\mathbb{F}_q$ es isomorfo al grupo $\mathbb{Z}_p \times \overbrace{\cdots}^{n} \times \mathbb{Z}_p$.

Se desconoce si $F_2(\mathbb{Z}_n)$ tiene comportamiento asintótico cuando $n \to \infty$. Lo único que se sabe se resume en el ejercicio siguiente.

**Ejercicio 2.16.** *Demostrar las siguientes desigualdades.*

$$\frac{1}{\sqrt{2}} \leqslant \liminf_{n\to\infty} \frac{F_2(\mathbb{Z}_n)}{\sqrt{n}} \leqslant \limsup_{n\to\infty} \frac{F_2(\mathbb{Z}_n)}{\sqrt{n}} = 1. \tag{16}$$

**2.4 Conjuntos $B_h$.** Los conjuntos $B_h$ son una generalización natural de los conjuntos de Sidon en los que todas las sumas de $h$ elementos del conjunto son distintas. De hecho, a los conjuntos de Sidon también se los denomina conjuntos $B_2$.

**Definición 2.2.** *Sea $G$ un grupo conmutativo. Decimos que $A \subset G$ es un conjunto $B_h$ si todas las sumas*

$$a_1 + \cdots + a_h, \quad a_1, \ldots, a_h \in A$$

*son distintas salvo por el orden de presentación de los sumandos. En general decimos que $A$ es un conjunto $B_h[g]$ si para todo $x \in G$ la ecuación*

$$x = x_1 + \cdots + x_h,\ x_i \in A$$

*tiene a lo más $g$ soluciones distintas salvo por permutaciones de los sumandos.*

Se define $F_h(n)$ como el mayor número de elementos que puede llegar a tener un conjunto $B_h$ en $\{1, \ldots, n\}$. De manera análoga se define $F_h(G)$ como el mayor tamaño de un conjunto $B_h$ en $G$.

Los conjuntos $B_h$ son bastante más esquivos que los conjuntos de Sidon debido a que no tienen una caracterización natural en términos de diferencias. Eso hace que algunos resultados, análogos a los que son conocidos para los conjuntos de Sidon, no se hayan logrado demostrar para los conjuntos $B_h$. En particular se desconoce el comportamiento asintótico de $F_h(n)$, incluso se desconoce si dicho comportamiento asintótico existe. En cualquier caso no es difícil obtener una cota superior utilizando el resultado del problema siguiente.

**Ejercicio 2.17.** *Demostrar que el número de sumas de $h$ elementos (no necesariamente distintos) de un conjunto $A$ es a lo más $\binom{|A|+h-1}{h}$.*

Sea $A \subset [1, n]$ un conjunto $B_h$ de enteros positivos. Como todas las sumas de $h$ elementos de $A$ son distintas y menores o iguales que $hn$, el resultado del ejercicio anterior implica que

$$\frac{|A|^h}{h!} < \binom{|A| + h - 1}{h} \leqslant hn \implies |A| < (h \cdot h!\, n)^{1/h}.$$

De manera análoga, en el caso de un grupo finito $G$ con $n$ elementos tenemos que

$$\frac{|A|^h}{h!} < \binom{|A| + h - 1}{h} \leqslant n \implies |A| < (h!\, n)^{1/h}.$$

Es decir,

$$F_h(n) < (h \cdot h!\, n)^{1/h}$$

y

$$F_h(G) < (h!\, |G|)^{1/h}.$$

Por otro lado se conocen tres construcciones de conjuntos $B_h \subset \{1, \dots, n\}$ de tamaño $\sim n^{1/h}$. Las dos primeras Bose y Chowla (1962/1963) son una generalización de los conjuntos de Sidon de Singer y Bose y la tercera es una generalización, obtenida por Gómez Ruiz y Trujillo Solarte (2011), que combina las ideas de la construcción de Ruzsa para $h = 2$ con la de Bose–Chowla para $h \geqslant 2$. Estas construcciones prueban las siguientes cotas inferiores:

$$F_h(\mathbb{Z}_{q^h-1}) \geqslant q$$
$$F_h(\mathbb{Z}_{q^h+\cdots+q+1}) \geqslant q+1$$
$$F_h(\mathbb{F}_{q^{h-1}} \times \mathbb{Z}_{q-1}) \geqslant q-1.$$

En Gómez Ruiz y Trujillo Solarte (ibíd.) se da una presentación unificada de estas tres construcciones. Aquí seguimos Gómez Ruiz y Trujillo Solarte (ibíd.) para exponer la construcción de Bose–Chowla, que es la más sencilla de todas ellas.

**Teorema 2.5** (Bose–Chowla)**.** *Sea $\mathbb{F}_{q^h}$ un cuerpo de $q^h$ elementos y sea $\theta$ un generador del grupo multiplicativo $\mathbb{F}_{q^h}^*$. Entonces el conjunto*

$$A = \{\log_\theta(\theta + a) : a \in \mathbb{F}_q\}$$

*es un conjunto $B_h$ en $\mathbb{Z}_{q^h-1}$.*

*Prueba.* Por ser $\theta$ un generador de $\mathbb{F}_{q^h}$, su polinomio minimal tiene grado $h$. Como consecuencia de esto vamos a ver que el conjunto

$$\theta + \mathbb{F}_q = \{\theta + a : a \in \mathbb{F}_q\}$$

es un conjunto $B_h$ multiplicativo en $\mathbb{F}_{q^h}^*$. Supongamos que no es así y que se da la igualdad siguiente entre dos productos de $h$ elementos de $\theta + \mathbb{F}_q$:

$$\prod_{i=1}^{h} (\theta + a_i) = \prod_{i=1}^{h} (\theta + a_i')$$

con $\{a_1, \dots, a_h\} \neq \{a_1', \dots, a_h'\}$. En ese caso el polinomio

$$\prod_{i=1}^{h} (X + a_i) - \prod_{i=1}^{h} (X + a_i'),$$

que es de grado menor que $h$, se anula en $\theta$, lo que contradice que el polinomio mínimo de $\theta$ es de grado $h$. Una vez visto que $\theta + \mathbb{F}_q$ es un conjunto de Sidon multiplicativo, es claro que el conjunto

$$A = \{\log_\theta(\theta + a) : a \in \mathbb{F}_q\}$$

es un conjunto $B_h$ en $\mathbb{Z}_{q^h-1}$. □

Razonando de la misma manera que lo hicimos para obtener la cota inferior sobre $F_2(n)$ a partir de construcciones de conjuntos de Sidon en grupos cíclicos, tenemos que

$$n^{1/h}(1+o(1)) \leqslant F_h(n) \leqslant (h \cdot h!)^{1/h}.$$

Con argumentos puramente combinatorios se puede mejorar la cota superior. Lindström (1969b) lo hizo por primera vez para sucesiones $B_4$ demostrando que $F_4(n) \leqslant (8n)^{1/4}(1+o(1))$, y fue generalizado por otros autores Jia (1993) y Chen (1994) para todo $h \geqslant 3$.

**Teorema 2.6.**

$$\begin{aligned} F_{2h-1}(n) &\leqslant (h!^2 n)^{1/(2h-1)}(1+o(1)) \\ F_{2h}(n) &\leqslant (h \cdot h!^2 n)^{1/(2h)}(1+o(1)). \end{aligned}$$

Una demostración más sencilla que las originales se puede obtener a partir de los dos ejercicios siguientes.

**Ejercicio 2.18.** *Demostrar que si $A \subset [1,n]$ es un conjunto $B_{2h}$ entonces el conjunto $C = A + \overset{h}{\cdots} + A$ satisface que $|C - C| \sim |C|^2$.*

**Ejercicio 2.19.** *Combinar el ejercicio anterior y el Teorema 2.3 para demostrar que si $A \subset [1,n]$ es un conjunto $B_{2h}$ entonces*

$$|A| \leqslant (1+o(1)) \left(h!^2 \cdot hn\right)^{1/(2h)}.$$

Todavía existen ligeras mejoras sobre estas cotas. Cuando $h$ es grande, se puede sacar partido de que las sumas $a_1 + \cdots + a_h$ tenderán a estar concentradas sobre su media. Resultados en esta dirección aparecen en D'yachkov y Rykov (1984) y mejoras posteriores aparecen en Shparlinskiĭ (1986), Cilleruelo (2001) y Green (2001). Para $h$ pequeño se puede utilizar análisis de Fourier para sacar partido del hecho de que los elementos de un conjunto suma $A + A$ no pueden estar bien distribuido en intervalos. Esta estrategia, basada en parte en las ideas de Cilleruelo, Ruzsa y Trujillo (2002), aparecen en Cilleruelo (2001) y Green (2001). En particular Ben Green (ibíd.) ha demostrado que $F_4(n) \leqslant (7n)^{1/4}(1+o(1))$.

## 3 Sucesiones de Sidon infinitas

La manera de cuantificar el tamaño de una sucesión infinita $A$ de enteros positivos es a través de su función contadora

$$A(x) = |A \cap [1,x]|,$$

que cuenta el número de términos de la sucesión menores o iguales que $x$. En lo que se refiere a sucesiones de Sidon infinitas, el principal objetivo consiste en construir sucesiones de Sidon infinitas que tengan una función contadora tan grande como sea posible.

Es fácil construir sucesiones infinitas de Sidon. Por ejemplo la sucesión de las potencias de 2 es una sucesión de Sidon porque todas las sumas de dos potencias de 2 son distintas. Pero esta sucesión no es muy interesante. Es muy poco densa, crece demasiado deprisa. Su función contadora es $A(x) \sim \log_2 x$.

La construcción más ingenua de un conjunto de Sidon con una función contadora decente es la generada por el algoritmo voraz. Consiste en empezar con $a_1 = 1,\ a_2 = 2$, y una vez construidos $a_1, ..., a_{n-1}$, añadir el menor entero positivo $a_n$ que no viole la condición de ser de Sidon; es decir, el siguiente que no sea de la forma $a_i - a_j + a_k,\ 1 \leqslant i, j, k \leqslant n-1$. Los primeros términos de esta sucesión, introducida por Erdős pero conocida como sucesión de Mian–Chowla, son los siguientes:

$$1, 2, 4, 8, 13, 21, 31, 45, 66, 81, 97, 123, 148, 182, 204, 252, 290, ...$$

Aunque se desconoce cómo crece realmente esta sucesión, como a lo más hay $(n-1)^3$ enteros prohibidos de la forma

$$a_i - a_j + a_k,\ 1 \leqslant i, j, k \leqslant n-1,$$

siempre es cierto que $a_n \leqslant (n-1)^3 + 1$, lo que nos permite seleccionar un conjunto de Sidon en $\{1, ..., m\}$ con $m^{1/3}$ elementos por lo menos. Es decir la función contadora de la sucesión voraz de Sidon satisface $A(x) \gg x^{1/3}$.

Se desconoce cuál es el verdadero comportamiento de la función contadora de esta sucesión. Los datos computacionales sugieren que $A(x)/x^{1/3} \to \infty$ y de hecho existe un modelo heurístico bastante sólido Cilleruelo (s.f.) y avalado por los datos computacionales, que permite conjeturar que

$$A(x) \sim c(x \log x)^{1/3},$$

donde $c$ es una constante explícita cuyo valor aproximado es $c = 1.7107\ldots.$

**3.1 Crecimiento de las sucesiones de Sidon infinitas.** Es claro que si $A$ es una sucesión de Sidon infinita entonces el conjunto formado por elementos hasta $x$ es un conjunto de Sidon finito en el intervalo $[1, x]$. De la cota superior trivial (1) para el máximo tamaño de un conjunto de Sidon en $[1, x]$ se obtiene:

$$A(x) \leqslant \sqrt{2x} + 1/2. \tag{17}$$

En principio podríamos pensar que pudiera existir una sucesión de Sidon infinita con $A(x) \gg x^{1/2}$, de manera análoga a lo que ocurre en el caso finito. Sin embargo Erdős demostró que tal sucesión no puede existir. En Stöhr (1955) aparece el siguiente resultado.

**Teorema 3.1** (Erdős). *Si $A$ es una sucesión de Sidon entonces*

$$\liminf_{x\to\infty} A(x) \left(\frac{\log x}{x}\right)^{1/2} \ll 1.$$

*Prueba.* Dividimos el intervalo $[1, N^2]$ en $N$ intervalos de longitud $N$:

$$I_l = ((l-1)N + 1, lN], \qquad l = 1, \ldots, N$$

y llamamos

$$D_l = |A \cap I_l| = A(lN) - A((l-1)N)$$

al número de elementos de $A$ en $I_l$. Contando las diferencias positivas de todas las parejas de dos elementos pertenecientes a un mismo intervalo tenemos que

$$\sum_{l=1}^{N}\binom{D_l}{2}\leqslant N$$

debido a que todas las diferencias que estamos contando son distintas y menores que $N$. Manipulando esta desigualdad y utilizando la estimación

$$A(N^2)\leqslant 2N+1/2$$

vista en (17), llegamos a la desigualdad

$$\sum_{l=1}^{N} D_l^2 \leqslant 2N+2\sum_{l=1}^{N} D_l = 2N+2A(N^2)\leqslant 7N. \tag{18}$$

Aplicando la desigualdad de Cauchy–Schwarz y utilizando (18) y la parte derecha de la desigualdad

$$\log N \leqslant \sum_{l=1}^{N}\frac{1}{l}\leqslant 2\log N \tag{19}$$

para $N\geqslant 3$, obtenemos

$$\sum_{l=1}^{N}\frac{D_l}{\sqrt{l}}\leqslant\left(\sum_{l=1}^{N} D_l^2\right)^{1/2}\left(\sum_{l=1}^{N}\frac{1}{l}\right)^{1/2}\leqslant\sqrt{14N\log N}. \tag{20}$$

Por otra parte, sumando por partes y utilizando la desigualdad

$$\frac{1}{\sqrt{l}}-\frac{1}{\sqrt{l+1}}\geqslant\frac{1}{4l^{3/2}}$$

para $l\geqslant 1$ obtenemos

$$\begin{aligned}\sum_{l=1}^{N}\frac{D_l}{\sqrt{l}} &= \sum_{l=1}^{N}\frac{A(lN)-A((l-1)N)}{\sqrt{l}}\\ &= \frac{A(N^2)}{\sqrt{N+1}}+\sum_{l=1}^{N}A(lN)\left(\frac{1}{\sqrt{l}}-\frac{1}{\sqrt{l+1}}\right)\\ &\geqslant \frac{1}{4}\sum_{l=1}^{N}\frac{A(lN)}{l^{3/2}}.\end{aligned}$$

Si definimos

$$\tau_N=\inf_{t\geqslant N}A(t)\left(\frac{\log t}{t}\right)^{1/2},$$

es claro que $A(lN) \geqslant \tau_N \left(\frac{lN}{\log(lN)}\right)^{1/2}$ para todo $l \geqslant 1$. Así que

$$\begin{aligned}\sum_{l=1}^{N} \frac{D_l}{\sqrt{l}} &\geqslant \tau_N \left(\frac{N}{\log(N^2)}\right)^{1/2} \frac{1}{4} \sum_{l=1}^{N} \frac{1}{l} \\ &\geqslant \tau_N \left(\frac{N}{\log N}\right)^{1/2} \frac{1}{4\sqrt{2}} \sum_{l=1}^{N} \frac{1}{l} \\ &\geqslant \tau_N \frac{(N \log N)^{1/2}}{4\sqrt{2}},\end{aligned}$$

donde en el último pase se ha utilizado la parte izquierda de la desigualdad (19). Esta desigualdad y (20) prueban que $\tau_N \leqslant 8\sqrt{7}$, y por lo tanto que

$$\liminf_{x\to\infty} A(x) \left(\frac{\log x}{x}\right)^{1/2} \leqslant \lim_{N\to\infty} \tau_N \leqslant 8\sqrt{7}.$$

□

**Ejercicio 3.1.** *Refinar la demostración del Teorema 3.1 para demostrar que si $A$ es una sucesión de Sidon infinita entonces*

$$\liminf_{x\to\infty} A(x) \left(\frac{\log x}{x}\right)^{1/2} \leqslant 4.$$

El siguiente ejercicio es interesante porque muestra que se puede llegar a una conclusión similar a la del Teorema 3.1 asumiendo solo que $|A_x - A_x| \sim |A_x|^2$. Curiosamente no sabemos si la conclusión también es cierta asumiendo que $|A_x + A_x| \sim |A_x|^2/2$.

**Ejercicio 3.2.** *Sea $A$ una sucesión de enteros positivos y para cada $x$ consideremos el conjunto $A_x = A \cap [1, x]$. Demostrar que si*

$$|A_x - A_x| \sim |A_x|^2$$

*entonces*

$$\liminf_{x\to\infty} \frac{|A_x|}{\sqrt{x}} = 0.$$

El Teorema 3.1 ha sido generalizado Chen (1993) para sucesiones $B_{2h}$ combinando la estrategia del Teorema 3.1 con el hecho de que el conjunto $hA = A + \overset{h}{\cdots} + A$ es "casi" un conjunto de Sidon.

**Teorema 3.2** (Chen). *Si $A$ es una sucesión $B_{2h}$ entonces*

$$\liminf_{n\to\infty} A(n) \left(\frac{\log n}{n}\right)^{\frac{1}{2h}} \ll 1.$$

Se desconoce si $\liminf_{x\to\infty} A(x)/x^{1/h} = 0$ cuando $A$ es una sucesión $B_h$ y $h$ es impar. En relación con este problema, Helm (1996) ha demostrado que no existe ninguna sucesión $B_3$ infinita con comportamiento asintótico de la forma $A(x) \sim cx^{1/3}$.

Como contrapunto al Teorema 3.1 Erdős demostró que existe una sucesión infinita de Sidon con $\limsup_{x\to\infty} A(x)/\sqrt{x} = 1/2$. Este resultado fue mejorado por Krückeberg (1961).

**Teorema 3.3** (Krukenberger). *Existe una sucesión infinita de Sidon $A$ con*

$$\limsup_{x\to\infty} \frac{A(x)}{\sqrt{x}} \geqslant \frac{1}{\sqrt{2}}.$$

*Prueba.* Sea $m_j = 2^{4^j}$. Para cada intervalo

$$I_j = (m_j + 2m_{j-1}, 2m_j]$$

consideremos un conjunto de Sidon $A_j \subset I_j$ de tamaño maximal

$$|A_j| \sim (m_j - 2m_{j-1})^{1/2} \sim m_j^{1/2}.$$

La última estimación asintótica se debe a que $m_{j-1} = o(m_j)$. Al conjunto $A_j$ le quitamos ahora todos los elementos $a$ para los que exista un $a' \in A_j$ con $0 < a - a' \leqslant 2m_{j-1}$. Como $A_j$ es de Sidon, existen a lo más $2m_{j-1}$ de estos elementos $a$. Así que el conjunto $A_j^*$ resultante tendrá todavía

$$|A_j^*| \geqslant |A_j| - 2m_{j-1} \sim m_j^{1/2}$$

elementos porque de nuevo $m_{j-1} = o\left(m_j^{1/2}\right)$. Veamos que la sucesión infinita

$$A = \bigcup_j A_j^*$$

satisface las condiciones del teorema.

Veamos primero que $A$ es un conjunto de Sidon. Supongamos que $a_1 + a_2 = a_1' + a_2'$ con $a_1 > a_1' \geqslant a_2' \geqslant a_2$ y todos ellos pertenecientes al conjunto $A$ . Supongamos que $a_1 \in A_j^*$. Veamos que necesariamente $a_1' \in A_j^*$. Si no fuera así entonces

$$a_1 = a_1' + a_2' - a_2 \leqslant a_1' + a_2' \leqslant 2 \cdot 2m_{j-1} < m_j,$$

lo que contradice el hecho de que $a_1 \in A_j^*$.

Veamos que también $a_2' \in A_j^*$. Es claro que $a_2' = a_1 - a_1' + a_2 > a_1 - a_1' > 2m_{j-1}$, debido a que, por construcción, es imposible que $a_1 - a_1' \leqslant 2m_{j-1}$. Eso implica que $a_2' \in A_j^*$.

Por último veamos que $a_2 \in A_j^*$. Escribimos

$$a_2 = a_1' + a_2' - a_1 > 2(m_j + 2m_{j-1}) - 2m_j = 4m_{j-1},$$

que en particular implica que $a_2 \in A_j^*$. Pero como $A_j^*$ es un conjunto de Sidon, no es posible que los cuatro elementos de la identidad $a_1 + a_2 = a_1' + a_2'$ pertenezcan a $A_j^*$.

Una vez visto que $A$ es de Sidon vayamos con el límite superior.

$$\frac{A(2m_j)}{(2m_j)^{1/2}} \geqslant \frac{|A_j^*|}{(2m_j)^{1/2}} \sim \frac{m_j^{1/2}}{(2m_j)^{1/2}} = \frac{1}{\sqrt{2}},$$

y por lo tanto,

$$\limsup_{x\to\infty} \frac{A(x)}{x^{1/2}} \geqslant \limsup_{j\to\infty} \frac{A(2m_j)}{(2m_j)^{1/2}} \geqslant \frac{1}{\sqrt{2}}.$$

□

Erdős se preguntaba si podría haber una sucesión infinita de Sidon tal que

$$\limsup_{x\to\infty} A(x)/\sqrt{x} = 1.$$

La constante 1 claramente no puede ser sustituida por una más grande porque, como hemos visto en el primer capítulo, siempre se tiene la cota superior $A(x) \leqslant \sqrt{x}(1+o(1))$. Una respuesta afirmativa al siguiente problema de Erdős implicaría la existencia de una sucesión infinita de Sidon con $\limsup_{x\to\infty} A(x)/\sqrt{x} = 1$.
**Problema** (Erdős): Sean $b_1, \dots, b_k$ enteros positivos que forman un conjunto de Sidon. ¿Existirán infinitos conjuntos de Sidon $A_n \subset [1,n]$ de tamaño $|A_n| \sim \sqrt{n}$ y que contengan a $b_1, \dots, b_k$?

**Ejercicio 3.3.** *Demostrar que una respuesta afirmativa a la pregunta de Erdős implicaría la existencia de una sucesión infinita de Sidon $A$ con* $\limsup_{x\to\infty} A(x)/\sqrt{x} = 1$.

**3.2 Construcción de sucesiones de Sidon infinitas.** La sucesión $(a_n)$ generada por el algoritmo avaricioso (la sucesión de Mian–Chowla) es la construcción más sencilla de una sucesión de Sidon infinita. Ya vimos en el capítulo anterior que $a_n \leqslant (n-1)^3+1$, lo que implica que $A(x) \gg x^{1/3}$. Esta construcción fue durante 50 años la mejor de la que se disponía hasta que en Ajtai, Komlós y Szemerédi (1981) demostraron la existencia de una sucesión infinita de Sidon con $A(x) \gg (x\log x)^{1/3}$.

Este resultado fue dramáticamente mejorado por Ruzsa (1998b) al demostrar la existencia de una sucesión infinita de Sidon con $A(x) = x^{\sqrt{2}-1+o(1)}$. La demostración de Ruzsa es muy ingeniosa. Ruzsa observó que los primos forman una sucesión de Sidon multiplicativa y por lo tanto la sucesión $\{\log p\}$ es una sucesión de Sidon de números reales.

A grandes rasgos la demostración de Ruzsa es como sigue. Considera un parámetro $\alpha \in [1,2]$ y para cada $\alpha$ construye una sucesión $B_\alpha = \{b_p\}$ indexada en los primos donde cada $b_p$ se construye a partir de los dígitos del desarrollo binario de $\alpha \log p$. Lo que Ruzsa demuestra es que para casi todo $\alpha \in [1,2]$ la sucesión $B_\alpha$ es "casi" de Sidon en el sentido de que podemos eliminar unos pocos términos de $B_\alpha$ para conseguir que la sucesión resultante sea de Sidon.

El siguiente ejercicio se puede considerar como la versión finita de la construcción de Ruzsa.

**Ejercicio 3.4.** Demostrar que el conjunto

$$A = \{\lfloor n\log(4p/\sqrt{n})\rfloor : \ \sqrt{n}/4 < p \leqslant \sqrt{n}/2\}$$

es un conjunto de Sidon en el intervalo $[1,n]$ de tamaño $|A| \gg \frac{\sqrt{n}}{\log n}$.

Una construcción similar a la de Ruzsa se puede hacer también utilizando los argumentos de los primos de Gauss en lugar de los logaritmos de los primos racionales.

**Ejercicio 3.5.** Sea $\phi(\mathbf{p})$ el argumento del primo Gaussiano $\mathbf{p} = |\mathbf{p}|e^{i\phi(\mathbf{p})}$. Demostrar que el conjunto

$$A = \{\lfloor 4n\phi(\mathbf{p})\rfloor : \ |\mathbf{p}| < \sqrt{n}, |\phi(\mathbf{p})| < \pi/4\}$$

es un conjunto de Sidon en el intervalo $[1,4n]$ de tamaño $|A| \gg \frac{\sqrt{n}}{\log n}$.

Volviendo a las sucesiones de Sidon infinitas, Erdős ofreció 1000 dólares por la resolución de la siguiente conjetura, que sigue sin estar resuelta.

**Conjetura 3.1.** *Para todo* $\epsilon > 0$ *existe una sucesión infinita de Sidon con* $A(x) \gg x^{1/2-\epsilon}$.

El Teorema 3.1 muestra que esta conjetura no es cierta para $\epsilon = 0$. Tanto las construcciones de Ajtai, Komlos y Szemerédi como la de Ruzsa son construcciones probabilísticas. No son explícitas. Recientemente Cilleruelo (2014b) ha encontrado una construcción explícita de una sucesión infinita de Sidon con función contadora similar a la de Ruzsa.

**Teorema 3.4** (Cilleruelo (ibíd.)). *Existe una sucesión de Sidon infinita A, que puede ser descrita explícitamente, con función contadora*

$$A(x) = x^{\sqrt{2}-1+o(1)}. \tag{21}$$

En esta sección nos dedicaremos a demostrar el Teorema 3.4 construyendo, de manera explícita, la sucesión de Sidon infinita a la que hace referencia el teorema.

**3.2.1 El método del logaritmo discreto.** La principal dificultad para construir sucesiones de Sidon infinitas densas reside en que las construcciones finitas que se han visto en el primer capítulo, provienen todas ellas de construcciones algebraicas en grupos finitos y no se sabe cómo extenderlas a sucesiones infinitas. La siguiente construcción de un conjunto finito de Sidon (en general de un conjunto $B_h$) es una construcción que podemos denominar semialgebraica en el sentido de que, aunque el grupo ambiente es $\mathbb{Z}_{q-1}$ (la parte algebraica), los elementos se describen a partir de los primos racionales. El tamaño de estos conjuntos de Sidon es menor (por un factor logarítmico) que el de los conjuntos de Sidon de procedencia algebraica, pero éstos ofrecen una mayor flexibilidad a la hora de extenderlos a una sucesión infinita.

**Teorema 3.5.** *Sea q un primo y g un generador de* $\mathbb{F}_q^*$. *El conjunto*

$$A = \{x : \ g^x \equiv p \text{ para algún primo } p \leqslant q^{1/h}\}$$

*es un conjunto* $B_h$ *en* $\mathbb{Z}_{q-1}$ *con* $\pi(q^{1/h})$ *elementos.*

*Prueba.* Supongamos que

$$x_1 + \cdots + x_h = y_1 + \cdots + y_h \pmod{q-1}$$

con $x_1, \ldots, x_h, y_1, \ldots, y_h \in A$. En ese caso tenemos que

$$g^{x_1} \cdots g^{x_h} \equiv g^{y_1} \cdots g^{y_h} \pmod{q}$$

y por construcción

$$p_1 \cdots p_h \equiv p'_1 \cdots p'_h \pmod{q}.$$

Como tanto el lado derecho como el izquierdo son menores que $q$, la congruencia es en realidad una igualdad en enteros

$$p_1 \cdots p_h = p'_1 \cdots p'_h$$

y el teorema fundamental de la aritmética implica que $\{p_1, \ldots, p_h\} = \{p'_1, \ldots, p'_h\}$, que a su vez implica que $\{x_1, \ldots, x_h\} = \{y_1, \ldots, y_h\}$. □

El conjunto $A$ también podía haber sido descrito de la forma

$$A = \{\log_g p : \ p \text{ primo} \ , \ p \leqslant q^{1/h}\},$$

donde $\log_g p$ es el logaritmo discreto de $x$ en base $g$.

Esta construcción fue la que inspiró la construcción de la sucesión de Sidon infinita que pasamos a describir.

**3.2.2 Bases generalizadas.** La manera de representar los números en una base dada (normalmente la base 10) es algo bien conocido por todos. Este hecho se puede generalizar de la manera siguiente.

Dada una sucesión de enteros positivos $2 \leqslant b_1 \leqslant \cdots \leqslant b_j \leqslant \cdots$ (la base), todo entero no negativo se puede escribir de manera única de la forma

$$a = x_1 + x_2 b_1 + \cdots + x_j b_1 \cdots b_{j-1} + \cdots$$

donde los $x_i$ (los dígitos) son enteros tales que $0 \leqslant x_i < b_i$.

**Ejercicio 3.6.** *Utilizar el algoritmo de Euclides para demostrar la afirmación anterior sobre las bases generalizadas.*

La base en la que expresaremos los elementos de nuestra sucesión será de la forma

$$\overline{q} := 4q_1, \ldots, 4q_i, \ldots$$

donde los $q_i$ son primos que satisfacen la desigualdad

$$2^{2i-1} < q_i \leqslant 2^{2i+1}.$$

Es claro que todo entero positivo $a$ se puede expresar de manera única de la forma

$$a = x_1 + x_2(4q_1) + \cdots + x_i(4q_1 \cdots 4q_{i-1}) + \cdots$$

donde los $x_i$ (los dígitos) son enteros tales que $0 \leqslant x_i < 4q_i$. El factor 4 aparece en los elementos de la base por razones técnicas que se verán más adelante.

Fijada la base, representamos cada entero $a$ mediante sus dígitos:

$$a := \cdots x_k \cdots x_1.$$

La ventaja de representar los enteros mediante dígitos es que podemos ver los enteros como si fueran vectores. Este hecho había sido utilizado anteriormente por otros autores.

Una observación importante es la siguiente. Supongamos que los dígitos de los enteros de nuestra sucesión satisfacen la desigualdad

$$q_i < x_i < 2q_i \tag{22}$$

y sean $a, a'$ dos elementos de dicha sucesión, con dígitos

$$\begin{aligned} a &= x_k \ldots\ldots x_1 \\ a' &= \quad x'_j \ldots x'_1. \end{aligned}$$

Como $x_i + x_i' < 4q_i$ para todo $i$, los dígitos de $a + a'$ en la base $\overline{q} := 4q_1, \dots, 4q_i, \dots$ se pueden calcular sumando los dos dígitos de cada posición:

(23) $$a + a' = (x_k + 0)\dots(x_{j+1} + 0)(x_j + x_j')\dots(x_1 + x_1').$$

Observar además que los dígitos en las posiciones $1, \dots, j$ son todos mayores que $2q_i$ y que el resto son menores que $2q_i$. Es decir, los dígitos de $a + a'$ determinan el número de dígitos de $a$ y $a'$ cuando éstos satisfacen (22).

Vamos a describir los elementos de nuestra sucesión en una base como la descrita y de tal manera que los dígitos de los elementos van a satisfacer la desigualdad $q_i < x_i < 2q_i$. De esta manera podremos sacar ventaja de la observación que acabamos de hacer.

**3.2.3 La distribución de los números primos.** La sucesión que vamos a construir va a estar indexada con la sucesión de los números primos. Por esa razón es conveniente recordar cómo se distribuyen los números primos en la sucesión de los enteros. La función $\pi(x)$ es la que cuenta el número de primos menores o iguales que $x$.

Uno de los teoremas fundamentales de la teoría de números es el teorema de los números primos que nos habla del comportamiento asintótico de la función $\pi(x)$.

**Teorema 3.6** (Teorema de los números primos). *Cuando $x \to \infty$ se tiene que*

$$\pi(x) \sim \frac{x}{\log x}.$$

El teorema de los números primos, pronosticado por matemáticos como Gauss y Riemann, fue demostrado independientemente por Jacques Hadamard y Charles-Jean de la Vallée Poussin en 1896.

**3.2.4 Una sucesión de Sidon infinita explícita.** Empezaremos la construcción de la sucesión $A$ del Teorema 3.4 indicando en qué base vamos a expresar sus elementos:
**La base.** Fijamos una sucesión de primos $(q_i)$ con

(24) $$2^{2i-1} < q_i < 2^{2i+1}$$

y consideramos la base generalizada

$$\overline{q} := 4q_1, \dots, 4q_i, \dots$$

Es esta base la que utilizaremos para describir, mediante sus dígitos, los elementos de la sucesión infinita de Sidon $A$ del Teorema 3.4.

Por comodidad utilizaremos la notación

$$Q_k = \prod_{i=1}^{k} q_i$$

para indicar el producto de los primeros $k$ primos de esa sucesión. Por (24) es claro que

(25) $$2^{k^2} < Q_k < 2^{(k+1)^2}.$$

**El conjunto de índices**: Vamos a enumerar los elementos de nuestra sucesión utilizando el conjunto de los primos $\mathcal{P}$ como conjunto de índices.

$$A = (a_p)_{p\in\mathcal{P}}$$

y representaremos cada elemento mediante sus dígitos en la base $\overline{q}$:

$$a_p = \cdots x_k(p)\cdots x_1(p).$$

**La función contadora**: Para determinar el número de dígitos de cada elemento, que será lo que a su vez determine el crecimiento de la sucesión, fijamos un número real $c,\ 0 < c < 1/2$ (que acabará siendo el exponente en la función contadora de $A$), y consideramos la siguiente partición de los números primos:

$$\mathcal{P} = \bigcup_{k\geqslant 2} \mathcal{P}_k,$$

donde

$$\mathcal{P}_k = \left\{p \text{ primo} :\ \frac{Q_{k-1}^c}{k-1} < p \leqslant \frac{Q_k^c}{k}\right\}.$$

En la sucesión que construiremos los elementos $a_p$ con $p \in \mathcal{P}_k$ tendrán exactamente $k$ dígitos. El cálculo del número de elementos de $\mathcal{P}_k$ lo dejamos como ejercicio.

**Ejercicio 3.7.** *Demostrar que*

$$|\mathcal{P}_k| \asymp \frac{Q_k^c}{k^3}. \tag{26}$$

Al final de la demostración tomaremos $c = \sqrt{2}-1$, pero ahora preferimos escribir simplemente $c$ para que se aprecie en la demostración por qué no es posible tomar otro valor mayor.

**Lema 3.1.** *Sea $A = (a_p)$ una sucesión indexada con los primos. Supongamos que todos los elementos $a_p$ con $p \in \mathcal{P}_k$ tienen exactamente $k$ dígitos. Entonces*

$$A(x) = x^{c+o(1)}.$$

*Prueba.* Consideremos, para cada $x$, el entero $k$ tal que

$$4^k Q_k < x \leqslant 4^{k+1} Q_{k+1}. \tag{27}$$

De (25) se sigue fácilmente que

$$x = 2^{k^2(1+o(1))}. \tag{28}$$

Observemos que si $p \leqslant \frac{Q_k^c}{k}$ entonces $p \in \mathcal{P}_j$ para algún $j \leqslant k$. Eso quiere decir que

$$a_p = x_1 + x_2(4q_1) + \cdots + x_j(4q_1\cdots 4q_{j-1})$$

para algunos $0 \leqslant x_i \leqslant q_i - 1,\ 1 \leqslant i \leqslant j$. En particular

$$a_p \leqslant (4q_1)\cdots(4q_j) \leqslant (4q_1)\cdots(4q_k) = 4^k Q_k < x$$

y por lo tanto $A(x) \geqslant \pi(Q_k^c/k)$. Finalmente el teorema de los números primos y (28) implican la desigualdad

$$\pi(Q_k^c/k) \geqslant \pi(2^{ck^2}/k) \gg 2^{ck^2}/k^3 = 2^{ck^2(1+o(1))} = x^{c+o(1)}.$$

Para la cota superior observemos que si $p > \frac{Q_{k+1}^c}{k+1}$ entonces $p \in \mathcal{P}_j$ para algún $j \geqslant k+2$ y entonces

$$a_p > (4q_1)\cdots(4q_{j-1}) \geqslant (4q_1)\cdots(4q_{k+1}) = 4^{k+1}Q_{k+1} \geqslant x.$$

Es decir, $A(x) \leqslant \pi\left(Q_{k+2}^c/(k+2)\right)$. De nuevo tenemos

$$\pi(Q_{k+2}^c/(k+2)) \leqslant Q_{k+2}^c = 2^{k^2(1+o(1))} = x^{c+o(1)}.$$

□

**Los dígitos**: Para terminar de describir nuestra sucesión

$$A = (a_p)_{p\in\mathcal{P}}$$

tenemos que decir quiénes son los dígitos de cada elemento $a_p$ en nuestra base $\overline{q}$.

Cada entero $a_p$ con $p \in P_k$ va a ser un entero $a_p = x_k \dots x_1$ con exactamente $k$ dígitos en nuestra base, lo que garantiza, gracias al Lema 3.1 que la función contadora de la sucesión satisface $A(x) = x^{c+o(1)}$.

El dígito $x_i(p)$, para $i \leqslant k$, es la solución de la congruencia

$$g_i^{x_i(p)} \equiv p \pmod{q_i}, \qquad q_i + 1 \leqslant x_i(p) \leqslant 2q_i - 1 \tag{29}$$

donde $g_i$ es un generador del grupo multiplicativo $\mathbb{F}_{q_i}^*$, que habremos fijado previamente para cada $q_i$. Definimos $x_i(p) = 0$ for $i > k$.

Es decir, si $p \in \mathcal{P}_k$, los dígitos de $a_p$ en la base $\overline{q} := 4q_1, \dots, 4q_i, \dots$ son

$$a_p = x_k x_{k-1} \cdots x_2 x_1,$$

donde los $x_i = x_i(p),\ i = 1, \dots, k$ han sido definidos en (29).

Utilizaremos la notación $A_c$ para designar a nuestra sucesión y enfatizar la dependencia de $c$. La siguiente proposición concierne a las propiedades de Sidon de la sucesión $A_c$.

**Proposición 3.1.** *Supongamos que existen $a_{p_1}, a_{p_2}, a_{p_1'}, a_{p_2'} \in A_c$ con $a_{p_1} > a_{p_1'} \geqslant a_{p_2'} > a_{p_2}$ y tales que*

$$a_{p_1} + a_{p_2} = a_{p_1'} + a_{p_2'}.$$

*Entonces tenemos que:*

*i) existen $j, k,\ j \leqslant k$ tales que $p_1, p_1' \in \mathcal{P}_k,\ p_2, p_2' \in \mathcal{P}_j$.*

*ii) $p_1 p_2 \equiv p_1' p_2' \pmod{Q_j}$*

*iii) $p_1 \equiv p_1' \pmod{Q_k/Q_j}$.*

*iv)* $Q_k^{1-c} < Q_j < Q_k^{\frac{c}{1-c}}$.

*Prueba.* Como $0 \leqslant x_j(p_1) + x_j(p_2) < 4q_j$ para todo $j$, la igualdad $a_{p_1} + a_{p_2} = a_{p_1'} + a_{p_2'}$ implica que los dígitos de ambas sumas son iguales:

$$x_j(p_1) + x_j(p_2) = x_j(p_1') + x_j(p_2') \tag{30}$$

para todo $j$. Por construcción y la observación (23) podemos ver que $p_1 \in \mathcal{P}_k$ y $p_2 \in \mathcal{P}_j$ donde $k$ es el mayor entero para el que

$$x_k(p_1) + x_k(p_2) \geqslant q_k + 1$$

y $j$ es el mayor para el que

$$x_j(p_1) + x_j(p_2) \geqslant 2q_j + 2.$$

Esta observación prueba la parte i). Para probar ii) y iii) observemos que (30) implica que para todo $i$ tenemos

$$g_i^{x_i(p_1)+x_i(p_2)} \equiv g_i^{x_i(p_1')+x_i(p_2')} \pmod{q_i}.$$

También sabemos que si $p \in \mathcal{P}_k$, entonces

$$g_i^{x_i(p)} \equiv p \pmod{q_i},\ i \leqslant k \tag{31}$$

$$g_i^{x_i(p)} \equiv 1 \pmod{q_i},\ i > k. \tag{32}$$

Así que $p_1p_2 \equiv p_1'p_2' \pmod{q_i}$ para todo $i \leqslant j$. Esta observación y el teorema chino del resto implica la congruencia en ii)

$$p_1p_2 \equiv p_1'p_2' \pmod{Q_j}.$$

Como $p_1p_2 \neq p_1'p_2'$, entonces $|p_1p_2 - p_1'p_2'| \geqslant Q_j$. Por otra parte, como $p_1, p_1' \in \mathcal{P}_k$ y $p_2, p_2' \in \mathcal{P}_j$ tenemos que $p_1p_2 \leqslant Q_k^cQ_j^c$ y tenemos la desigualdad

$$Q_k^cQ_j^c > |p_1p_2 - p_1'p_2'| \geqslant Q_j \Longrightarrow Q_j < Q_k^{\frac{c}{1-c}}. \tag{33}$$

En particular se tiene que $j < k$ porque $c < 1/2$. Las congruencias (31) implican que $p_1 \equiv p_1' \pmod{q_i}$ para todo $i$ con $j+1 \leqslant i \leqslant k$. El teorema chino del resto implica la congruencia

$$p_1 \equiv p_1' \pmod{Q_k/Q_j}.$$

Procediendo como antes obtenemos la desigualdad

$$Q_k^c > |p_1 - p_1'| \geqslant Q_k/Q_j \Longrightarrow Q_j > Q_k^{1-c}. \tag{34}$$

Las desigualdades (33) y (34) implican iv). □

Observemos que si hubiera sumas repetidas, entonces la Proposición 3.1, iv) implicaría que $1 - c < \frac{c}{1-c}$, lo cual no es cierto para $c = \frac{3-\sqrt{5}}{2} = 0.38...$ Así que la sucesión $A_c$ es una sucesión de Sidon para este valor de $c$, que es mayor que $1/3$. Esta observación nos proporciona el siguiente corolario.

**Corolario 3.1.** *La sucesión* $A = A_c$ *con* $c = \frac{3-\sqrt{5}}{2}$ *es una sucesión infinita de Sidon con función contadora* $A(x) = x^{\frac{3-\sqrt{5}}{2}+o(1)}$.

Si $c > \frac{3-\sqrt{5}}{2}$ ya no es cierto que $A_c$ vaya a ser una sucesión de Sidon. Aparecerán infinitas sumas que se repiten. Pero si aparecen con poca frecuencia podemos eliminar los elementos de $A_c$ involucrados en esas sumas para así obtener una verdadera sucesión de Sidon. Por supuesto hay que controlar que los elementos que vamos a descartar no sean demasiados para que eso no afecte demasiado al orden de la función contadora de la nueva sucesión. Eso lo vamos a poder hacer para todo $c \leqslant \sqrt{2}-1$ y esa es la estrategia que seguiremos para demostrar el Teorema 3.4.

Consideremos la sucesión

$$A_c^* = (a_p)_{p \in \mathcal{P}^*}$$

donde los números $a_p$ se definen como antes pero ahora $\mathcal{P}^*$ es el conjunto de los primos que quedan después de eliminar de cada $\mathcal{P}_k$ un subconjunto $\mathcal{R}_k$ con el propósito de evitar la presencia de algunas sumas repetidas que pudieran aparecer. Para que no parezca extraña la definición de los conjuntos $\mathcal{R}_k$, esperaremos a que la propia demostración nos diga quiénes tienen que ser estos conjuntos. En cualquier caso sea

$$\mathcal{P}^* = \bigcup_k (\mathcal{P}_k \setminus \mathcal{R}_k)$$

donde los conjuntos de primos eliminados $\mathcal{R}_k$ los definiremos más adelante.

Vamos a demostrar que para $c = \sqrt{2}-1$, la sucesión $A_c^* = \{a_p\}_{p \in \mathcal{P}^*}$ es una sucesión infinita de Sidon con $A_c^*(x) = x^{\sqrt{2}-1+o(1)}$.

Para ver que $A_c^*$ es una sucesión de Sidon, supongamos que

$$a_{p_1} + a_{p_2} = a_{p_1'} + a_{p_2'}$$

con $a_{p_1} > a_{p_1'} \geqslant a_{p_2'} > a_{p_2}$ y $p_1, p_1', p_2, p_2' \in \mathcal{P}^*$. La proposición 3.1 implica que

$$p_1, p_1' \in \mathcal{P}_k \setminus \mathcal{R}_k \quad \text{ y } \quad p_2, p_2' \in \mathcal{P}_j \setminus \mathcal{R}_j$$

para algún par de índices $j, k$ que satisface $Q_j < Q_k^{\frac{c}{1-c}}$. Esta última restricción es consecuencia de la Proposición 3.1, iv).

Seguidamente observemos que gracias a las partes ii) and iii) de la Proposición 3.1 podemos escribir

$$p_1(p_2 - p_2') = s_1 Q_j + s_2 Q_k / Q_j$$

para los enteros no nulos

$$s_1 = \frac{p_1 p_2 - p_1' p_2'}{Q_j}, \qquad s_2 = \frac{(p_1' - p_1) p_2'}{Q_k / Q_j},$$

los cuales satisfacen las desigualdades

$$1 \leqslant |s_1| \leqslant \frac{Q_j^c Q_k^c}{jkQ_j}, \quad 1 \leqslant |s_2| \leqslant \frac{Q_j^c Q_k^c Q_j}{jkQ_k}.$$

Esto implica que $p_1$ es un primo de $\mathcal{P}_k$ que divide a algún elemento $s$ de alguno de los conjuntos

$$S_{j,k} = \left\{ s = s_1 Q_j + s_2 Q_k / Q_j \, : \, 1 \leqslant |s_1| \leqslant \frac{Q_j^c Q_k^c}{jkQ_j}, \, 1 \leqslant |s_2| \leqslant \frac{Q_j^c Q_k^c Q_j}{jkQ_k} \right\}$$

para algún $j$ tal que $Q_j < Q_k^{\frac{c}{1-c}}$.

Ahora parece más claro cómo deberíamos definir el conjunto $\mathcal{R}_k$ al que hemos aludido al principio de la demostración. El conjunto $\mathcal{R}_k$ es el conjunto de los primos $p_1$ en $\mathcal{P}_k$ que dividen a algún elemento de algún $S_{j,k}$ para algún $j$ tal que $Q_j < Q_k^{\frac{c}{1-c}}$.

De esta manera es claro que la sucesión $A_c^*$ es de Sidon. Si hubiera una suma repetida $a_{p_1} + a_{p_2} = a_{p_1'} + a_{p_2'}$ con $p_1, p_1' \in \mathcal{P}_k$, $p_2, p_2' \in \mathcal{P}_j$ entonces tendríamos que $p_1 \in \mathcal{R}_k$ y por tanto $a_{p_1} \notin A_c^*$.

Para demostrar que $A_{\bar{q},c}^*(x) = x^{c+o(1)}$ sólo necesitamos probar que $|\mathcal{R}_k| = o(|\mathcal{P}_k|)$.

Primero veamos que para cada $s \in S_{j,k}$ y $k$ suficientemente grande, existe a lo más un $p \in \mathcal{P}_k$ dividiendo a $s$. Si $p, p' \mid s$ tendríamos que

$$\frac{Q_{k-1}^{2c}}{(k-1)^2} < pp' \leqslant |s| \leqslant 2 \cdot \frac{Q_j^c Q_k^c}{jk} < \frac{Q_k^{\frac{c}{1-c}}}{k},$$

lo cual no puede ser cierto para $k$ grande porque $2c > \frac{c}{1-c}$ para $c < 1/2$. Por lo tanto,

$$|S_{j,k}| \leqslant \left(2 \cdot \frac{Q_j^c Q_k^c}{jkQ_j}\right) \left(\frac{Q_j^c Q_k^c Q_j}{jkQ_k}\right) \leqslant 2 \frac{Q_j^{2c} Q_k^{2c-1}}{j^2 k^2}. \tag{35}$$

Utilizando (35), la identidad

$$\frac{2c^2}{1-c} + 2c - 1 = c$$

para $c = \sqrt{2} - 1$ y la estimación (26) en el último paso, tenemos, para $k$ suficientemente grande, la estimación requerida,

$$\begin{aligned} |\mathcal{R}_k| \leqslant \sum_{Q_j < Q_k^{\frac{c}{1-c}}} |S_{j,k}| &\ll \sum_{Q_j < Q_k^{\frac{c}{1-c}}} \frac{Q_j^{2c} Q_k^{2c}}{j^2 k^2 Q_k} \\ &\ll \frac{Q_k^{\frac{2c^2}{1-c} + 2c - 1}}{k^4} = \frac{Q_k^c}{k^4} = o\left(|\mathcal{P}_k|\right). \end{aligned}$$

El último paso es consecuencia de la estimación (26). Precisamente el valor $c = \sqrt{2} - 1$ sale de la igualdad $\frac{2c^2}{1-c} + 2c - 1 = c$.

### 3.3 Sucesiones $B_h$ infinitas.

**Ejercicio 3.8.** Sea $A_h$ la sucesión $B_h$ construida con el algoritmo avaricioso. Demostrar que $A_h(x) \gg x^{\frac{1}{2h-1}}$.

Como ya comentamos en el capítulo anterior, el exponente $1/(2h-1)$ que se obtiene con el algoritmo voraz fue mejorado por Ruzsa para el caso $h=2$. Recientemente Cilleruelo y Tesoro (2015) se ha utilizado una variante del método de Ruzsa que permite mejorar dicho exponente en los casos $h=3$ y $h=4$. Sin embargo el método de Ruzsa se extiende mal para valores más grandes de $h$.

Para valores mayores de $h$ adaptaremos el método utilizado en la construcción de la sucesión que aparece en el Teorema 3.4 y lo combinaremos con un argumento probabilístico para demostrar la existencia de sucesiones $B_h$ que mejoran el exponente $\frac{1}{2h-1}$ para todo $h$.

**Teorema 3.7** (Cilleruelo (2014b)). *Para todo $h \geqslant 3$ existe una sucesión $B_h$, $\mathcal{B}$, con*

$$\mathcal{B}(x) = x^{\sqrt{(h-1)^2+1}-(h-1)+o(1)}.$$

*Prueba.* Fijemos

$$c = \sqrt{(h-1)^2+1} - (h-1)$$

y consideremos la función $f(t) = ct^2 - t^2/\sqrt{\log t}$. Sea

$$\mathcal{P} = \bigcup_{k \geqslant 3} \mathcal{P}_k$$

donde

$$\mathcal{P}_k = \left\{ p \text{ primo } : \ e^{f(k-1)} < p \leqslant e^{f(k)} \right\}.$$

Sea $\overline{q} := q_1 < q_2 < \cdots$ una sucesión de primos con

$$e^{2j-1} < q_j \leqslant e^{2j+1}$$

y sea $g_j$ un generador de $\mathbb{F}_{q_j}^*$. Para cada $p \in \mathcal{P}_k$ definimos el entero

$$b_p = x_1(p) + \sum_{2 \leqslant j \leqslant k} x_j(p)(h^2 q_1)\cdots(h^2 q_{j-1}),$$

donde $x_j(p)$ es la solución de la congruencia

$$g_j^{x_j(p)} \equiv p \pmod{q_j}, \qquad (h-1)q_j + 1 \leqslant x_j(p) \leqslant hq_j - 1.$$

Definimos $x_j(p) = 0$ para $j > k$.

Es claro que la sucesión $\mathcal{B}_{\overline{q},c} = \{b_p\}$ será una sucesión $B_h$ si y sólo si para todo $l$, $2 \leqslant l \leqslant h$ no existe una suma repetida de la forma

(36)
$$\begin{gathered} b_{p_1} + \cdots + b_{p_l} = b_{p'_1} + \cdots + b_{p'_l} \\ \{b_{p_1}, \ldots, b_{p_l}\} \cap \{b_{p'_1}, \ldots, b_{p'_l}\} = \varnothing \\ b_{p_1} \geqslant \cdots \geqslant b_{p_l} \\ b_{p'_1} \geqslant \cdots \geqslant b_{p'_l}. \end{gathered}$$

La siguiente proposición es una generalización de la Proposición 3.1.

**Proposición 3.2.** *Supongamos que existen* $p_1, \dots, p_l, p'_1, \dots, p'_l \in \mathcal{B}_{\overline{q},c}$ *satisfaciendo* (36). *Entonces se tiene que:*

*i)* $p_i, p'_i \in \mathcal{P}_{k_i},\ i = 1, \dots, l$ *for some* $k_l \leqslant \cdots \leqslant k_1$.

$$p_1 \cdots p_l \equiv p'_1 \cdots p'_l \pmod{Q_{k_l}}$$

*ii)*
$$\begin{aligned} p_1 \cdots p_{l-1} &\equiv p'_1 \cdots p'_{l-1} \pmod{Q_{k_{l-1}}/Q_{k_l}} \quad \text{if } k_l < k_{l-1} \\ &\cdots \\ p_1 &\equiv p'_1 \pmod{Q_{k_1}/Q_{k_2}} \quad \text{if } k_2 < k_1. \end{aligned}$$

*iii)* $k_l^2 < \frac{c}{1-c}\left(k_1^2 + \cdots + k_{l-1}^2\right)$.

*iv)* $q_1 \cdots q_{k_1} \mid \prod_{i=1}^{l}\left(p_1 \cdots p_i - p'_1 \cdots p'_i\right)$.

*Prueba.* La demostración es similar a la de la Proposición 3.1: En este caso $k_i$ es el $j$ más grande tal que

$$x_j(p_1) + \cdots + x_j(p_l) \geqslant i((h-1)q_j + 1).$$

La parte iii) es consecuencia de la primera congruencia de ii). La parte iv) es también una consecuencia obvia de la parte ii). □

La sucesión $\mathcal{B}_{\overline{q},c}$ definida al principio de esta sección puede no ser una sucesión $B_h$ para el valor de $c$ que hemos fijado. El plan de la demostración es quitar de $\mathcal{B}_{\overline{q},c} = (b_p)_{p\in\mathcal{P}}$ el mayor elemento que aparezca en cada repetición para obtener una verdadera sucesión $B_h$.

Más precisamente, definimos $\mathcal{P}^* = \mathcal{P}^*(\overline{q})$ como el conjunto

$$\mathcal{P}^* = \bigcup_k (\mathcal{P}_k \setminus \mathcal{R}_k(\overline{q}))$$

donde $\mathcal{R}_k(\overline{q}) = \{p \in \mathcal{P}_k :\ b_p \text{ es el mayor elemento en alguna ecuación (36)}\}$.

Al haber eliminado todas las posibles sumas repetidas es claro que la sucesión

$$\mathcal{B}^*_{\overline{q},c} = (b_p)_{p\in\mathcal{P}^*}$$

es una sucesión $B_h$.

Recordemos que en el Lema 3.1 demostrábamos que $\mathcal{B}_{\overline{q},c}(x) = x^{c+o(1)}$. Y si $|\mathcal{R}_k(\overline{q})| = o(|\mathcal{P}_k|)$, tenemos que

$$\mathcal{B}^*_{\overline{q},c}(x) \sim \mathcal{B}_{\overline{q},c}(x) = x^{c+o(1)}.$$

Así que la demostración del Teorema 3.7 se completará si probamos que existe una sucesión $\overline{q}$ tal que $|\mathcal{R}_k(\overline{q})| = o(|\mathcal{P}_k|)$ cuando $k \to \infty$.

Para $2 \leqslant l \leqslant h$ escribimos

$$\text{Bad}_l(\overline{q}, k_l, \dots, k_1) = \{(p_1, \dots, p'_l) :\ p_i, p'_i \in \mathcal{P}_{k_i},\ i = 1, \dots, l \text{ satisfacen (36)}\}.$$

Observemos que cada $p \in \mathcal{R}_k(\overline{q})$ proviene de alguna $2l$-tupla

$$(p_1, \dots, p'_l) \in \text{Bad}_l(\overline{q}, k_l, \dots, k_1),$$

con $2 \leqslant l \leqslant h,\ k_l \leqslant \cdots \leqslant k_1 = k$. Entonces,

$$(37)\qquad \begin{aligned} |\Re_k(\overline{q})| &\leqslant \sum_{l=2}^{h} \sum_{k_l \leqslant \cdots \leqslant k_1 = k} |\mathrm{Bad}_l(\overline{q}, k_l, \ldots, k_1)| \\ &\leqslant hk^{h-1} \max_{\substack{2 \leqslant l \leqslant h \\ k_l \leqslant \cdots \leqslant k_1 = k}} |\mathrm{Bad}_l(\overline{q}, k_l, \ldots, k_1)|. \end{aligned}$$

No sabemos dar una buena cota superior para $|\mathrm{Bad}_l(\overline{q}, k_l, \ldots, k_1)|$ para una sucesión concreta de primos $\overline{q} := q_1 < q_2 < \cdots$, pero lo sabemos hacer en media y es aquí donde introducimos el argumento probabilístico. Si el lector está familiarizado con el trabajo de Ruzsa, verá que la sucesión $\overline{q}$ jugará el mismo papel que el parámetro $\alpha$ en la construcción de Ruzsa.

Consideremos el espacio probabilístico de las sucesiones $\overline{q} := q_1 < q_2 < \cdots$ donde cada $q_j$ se elige uniformemente entre todos los primos en el intervalo $(e^{2j-1}, e^{2^j+1}]$. Usaremos que

$$\pi(e^{2k+1}) - \pi(e^{2k-1}) \gg e^{2k}/k = e^{2k+O(\log k)}$$

para deducir que dados $q_1 < \cdots < q_{k_1}$ satisfaciendo

$$e^{2j-1} < q_j \leqslant e^{2j+1}$$

tenemos que

$$\begin{aligned} \mathbb{P}(q_1, \ldots, q_{k_1} \in \overline{q}) &= \prod_{k=1}^{k_1} \frac{1}{\pi(e^{2k+1}) - \pi(e^{2k-1})} \\ &\leqslant e^{-k_1^2 + O(k_1 \log k_1)}. \end{aligned}$$

Entonces, para una $2l-$tupla dada $(p_1, \ldots, p_l')$, usamos la Proposición 3.2, iv) y la estimación $\tau(n) = n^{O(1/\log\log n)}$ para la función divisor para deducir que

$$\begin{aligned} &\mathbb{P}((p_1, \ldots, p_l') \in \mathrm{Bad}_l(\overline{q}, k_l, \ldots, k_1)) \\ &\qquad \leqslant \sum_{\substack{q_1, \ldots, q_{k_1} \\ q_1 \cdots q_{k_1} \mid \prod_{i=1}^{l} (p_1 \cdots p_i - p_1' \cdots p_i')}} \mathbb{P}(q_1, \ldots, q_{k_1} \in \overline{q}) \\ &\qquad \leqslant \tau\left(\prod_{i=1}^{l} (p_1 \cdots p_i - p_1' \cdots p_i')\right) e^{-k_1^2 + O(k_1 \log k_1)} \\ &\qquad \leqslant e^{-k_1^2 + O(k_1^2/\log k_1)}. \end{aligned}$$

Usando la Proposition 3.2 iii) en la última desigualdad tenemos que

$$\begin{aligned}
&\mathbb{E}(|\{(p_1,\dots,p_l'):\ p_i,p_i'\in \wp_{k_i},\ i=1,\dots,l \text{ satisfying } (36)\}|)\\
&\qquad\leqslant e^{-k_1^2+O(k_1^2/\log k_1)}|\{(p_1,\dots,p_l'):\ p_i,p_i'\in \wp_{k_i}\}|\\
&\qquad\leqslant e^{-k_1^2+O(k_1^2/\log k_1)}|\wp_{k_1}|^2\cdots|\wp_{k_l}|^2\\
&\qquad\leqslant e^{-k_1^2+O(k_1^2/\log k_1)}\cdot e^{2f(k_1)+\cdots+2f(k_l)}\\
&\qquad\leqslant e^{-k_1^2+\frac{2c}{1-c}\left(k_1^2+\cdots+k_{l-1}^2\right)-(2c+o(1))k_1^2/\sqrt{\log k_1}}\\
&\qquad\leqslant e^{\left(-1+\frac{2c(l-1)}{1-c}\right)k_1^2-(2c+o(1))k_1^2/\sqrt{\log k_1}}.
\end{aligned}$$

Usando (37), obtenemos

$$\mathbb{E}(|\Re_k(\overline{q})|)\leqslant e^{\left(-1+\frac{2c(h-1)}{1-c}\right)k^2-(2c+o(1))k^2/\sqrt{\log k}}.$$

Finalmente usamos la identidad $-1+\frac{2c(h-1)}{1-c}-c=0$ para $c=\sqrt{(h-1)^2+1}-(h-1)$ para obtener

$$\begin{aligned}
\mathbb{E}\left(\sum_k\frac{|\Re_k(\overline{q})|}{|\wp_k|}\right)&\leqslant\sum_k k^2e^{\left(-1+\frac{2c(h-1)}{1-c}-c\right)k^2-(c+o(1))k^2/\sqrt{\log k}}\\
&\leqslant\sum_k k^2e^{-(c+o(1))k^2/\sqrt{\log k}}.
\end{aligned}$$

Como la serie es convergente tenemos que para casi toda sucesión $\overline{q}$ la serie

$$\sum_k\frac{|\Re_k(\overline{q})|}{|\wp_k|}$$

es convergente. Así que, para cualquiera de estas sucesiones $\overline{q}$ tenemos $|\Re_k(\overline{q})|=o(|\wp_K|)$, que es lo que pretendíamos probar. □

## 4 Problemas sin resolver sobre conjuntos de Sidon

Este capítulo está dedicado a aquellos problemas sobre los conjuntos de Sidon a los que todavía no se sabe dar respuesta. Muchos de ellos ya han aparecido en los dos capítulos anteriores. Gran parte de ellos fueron propuestos por Erdős hace muchos años y son, presumiblemente, muy difíciles. Pero hay otros más recientes que quizás no se hayan pensado lo suficiente y sean más asequibles. Nuestra intención es que esta colección de problemas pueda servir de fuente de inspiración para jóvenes investigadores.

**4.1 Conjuntos de Sidon en intervalos.** Erdős y Turán (1941) demostraron que si $A\subset[1,n]$ es un conjunto de Sidon entonces $|A|<n^{1/2}+O(n^{1/4})$. Lindström (1969a) y Ruzsa (1998b) obtuvieron demostraciones más limpias que dan lugar a la cota $|A|<n^{1/2}+n^{1/4}+1$. Esta

estimación lleva 54 años sin mejorarse salvo por la insignificante mejora Cilleruelo (2010b) (que ya estaba implícita en Ruzsa (1998b)) y que supone la cota $|A| < n^{1/2} + n^{1/4} + 1/2$. Mejorar esta cota superior (incluso eliminar el $1/2$ para $n$ grande) es un problema que se resiste.

Durante mucho tiempo Erdős estuvo convencido de que las cotas anteriores se podían sustituir por $|A| \leqslant \sqrt{n} + O(1)$. Pero, según afirma el propio Erdős (1989), Ruzsa y H. Taylor le convencieron de que había razones que sugerían que podía no ser cierto. La conjetura más plausible es la siguiente.

**Conjetura 4.1** (Erdős). *Para todo $\epsilon > 0$, si $A \subset \{1, \dots, n\}$ es un conjunto de Sidon, entonces $|A| \leqslant \sqrt{n} + O(n^{\epsilon})$.*

Como ya hemos comentado, se cree que la conjetura anterior no es cierta para $\epsilon = 0$. De otra manera:

**Conjetura 4.2.** *Para todo $M$ existe un $n$ y un conjunto de Sidon $A \subset [1, n]$ con $|A| > \sqrt{n} + M$.*

Los resultados computacionales avalan esta última conjetura. Pero también existen argumentos heurísticos sólidos a su favor. Uno de ellos está relatado en el Ejercicio 2.15. Otro resultado que apoya esta conjetura es que la conjetura análoga en dimensión 2 sí se sabe cierta Cilleruelo (2010b).

Si seguimos con detalle cualquiera de los argumentos que proporcionan la cota superior $|A| < \sqrt{n} + O(n^{1/4})$ cuando $A$ es un conjunto de Sidon en $[1, n]$ observamos que sólo se utiliza el hecho de que todas las diferencias positivas $a - a'$ menores que $n^{3/4}$ son distintas. Sería interesante estudiar si bajo estas condiciones, menos restrictivas que las de ser conjunto de Sidon, la cota superior está más cerca de lo mejor posible. Un primer paso sería demostrar la siguiente conjetura, previsiblemente más sencilla.

**Conjetura 4.3.** *Para todo $M$ existe un $n$ conjunto $A \subset [1, n]$ con $|A| > \sqrt{n} + M$ y con la propiedad de que todas las diferencias $a - a'$, $a, a' \in A$ con $0 < a - a' < n^{1/2}$ son distintas.*

### 4.2 Conjuntos de Sidon en dimensiones superiores.

**Conjetura 4.4.** *Para todo $\epsilon > 0$, si $A \subset [1, n]^d$ es un conjunto de Sidon, entonces $|A| \leqslant n^{d/2} + O(n^{\epsilon})$.*

Si esta conjetura es cierta es sin duda muy difícil de demostrar, al menos en dimensión $d = 2$. La razón es que en ese caso daría una respuesta positiva a una antigua conjetura de Vinogradov sobre el menor residuo no cuadrático módulo $p$: para todo $\epsilon > 0$ y $p$ primo suficientemente grande existe un residuo no cuadrático menor que $p^{\epsilon}$.

**Conjetura 4.5.** *Para todo $d \geqslant 1$ y para todo $M$ existe un conjunto de Sidon en $[1, n]^d$ con $|A| > n^{d/2} + M$.*

Como hemos comentado anteriormente, se sabe que es cierto para $d = 2$. Quizás los casos $d \geqslant 3$ no se hayan intentado lo suficiente.

**4.3 Conjuntos de Sidon en grupos finitos.** ¿Cuál es el mayor tamaño de un conjunto de Sidon en $\mathbb{Z}_n$? Si $n$ es de la forma $n = p(p-1)$, $n = q^2 + q + 1$ o $n = q^2 - 1$, donde $q$ es la potencia de un primo, entonces $\mathbb{Z}_n$ contiene un conjunto de Sidon de tamaño $\sqrt{n}+O(1)$. Cualquiera de estas familias implica que

$$\limsup_{n\to\infty} \frac{F_2(\mathbb{Z}_n)}{\sqrt{n}} = 1.$$

Estas construcciones obedecen a ciertos milagros algebraicos que no tienen que ocurrir para todo $n$. Por otra parte es claro que cualquier conjunto de Sidon en el intervalo $[1, n/2]$ es en particular un conjunto de Sidon en $\mathbb{Z}_n$. Este argumento implica que $\liminf_{n\to\infty} \frac{F_2(\mathbb{Z}_n)}{\sqrt{n}} \geqslant 1/\sqrt{2}$. Aunque creemos que el valor de este límite es $1/\sqrt{2}$, ni siquiera sabemos demostrar la siguiente conjetura.

**Conjetura 4.6.**

$$\liminf_{n\to\infty} \frac{F_2(\mathbb{Z}_n)}{\sqrt{n}} < 1.$$

Los conjuntos de la forma

$$A = \{(f(x), g(x)) : x \in \mathbb{Z}_p\} \subset \mathbb{Z}_p \times \mathbb{Z}_p \tag{38}$$

donde $f$ y $g$ son polinomios independientes de grados $1 \leqslant \deg f, \deg g \leqslant 2$ son conjuntos de Sidon con $p$ elementos. La demostración de la siguiente conjetura sería el primer resultado inverso cobre conjuntos de Sidon.

**Conjetura 4.7.** *Todos los conjuntos de Sidon en $\mathbb{Z}_p\times\mathbb{Z}_p$ con $p$ elementos son de la forma descrita en* (38).

Se sabe que si $A = \{(x, g(x)) : x \in \mathbb{Z}_p\}$ es un conjunto de Sidon en $\mathbb{Z}_p \times \mathbb{Z}_p$ entonces $g$ es un polinomio cuadrático.

**4.4 Sucesiones infinitas de Sidon.** La sucesión avariciosa de Sidon (la sucesión de Mian–Chowla) es aquella que, empezando en $a_1 = 1$, tiene como término $a_n$ el menor entero positivo que se puede elegir de tal manera que $\{a_1, \ldots, a_n\}$ sea un conjunto de Sidon. Mientras que un argumento sencillo muestra que $A(x) \geqslant x^{1/3}$, se desconoce cuál es el verdadero orden de magnitud de $A(x)$.

**Conjetura 4.8.** *La función contadora de la sucesión de Mian–Chowla satisface que*

$$\frac{A(x)}{x^{1/3}} \to \infty.$$

De hecho, argumentos heurísticos y computacionales sugieren que $A(x)$ debería tener un comportamiento asintótico de la forma

$$A(x) \sim c(x \log x)^{1/3}.$$

Sin embargo, no se sabe ni siquiera si $A(x) \ll x^{1/2-\epsilon}$ para algún $\epsilon > 0$.

Los $n$ primeros términos de la sucesión de Mian–Chowla forman una sucesión maximal en el sentido de que no es posible añadir un elemento entre 1 y $a_n$ sin que se destruya la propiedad de ser de Sidon. Sea $M(n)$ el menor tamaño de un conjunto de Sidon en $[1, n]$ con la propiedad de ser maximal. Se desconoce cuál es el orden de magnitud de $M(n)$. Es claro que $M(n) \geqslant n^{1/3}$ y por otra parte Ruzsa (1998a) ha demostrado que $M(n) \ll (n \log n)^{1/3}$.

**Problema 4.1.** *¿Es cierto que $M(n)/n^{1/3} \to \infty$?*

Una respuesta afirmativa a este problema implicaría inmediatamente la Conjetura 4.8.

La siguiente conjetura aparece insistentemente en varios artículos de Erdős.

**Conjetura 4.9** (Erdős). *Para todo $\alpha < 1/2$ existe una sucesión infinita de Sidon $A$ con $A(x) \gg x^{\alpha}$.*

Ruzsa demostró la existencia de una sucesión de Sidon $A$ con $A(x) = x^{\sqrt{2}-1+o(1)}$. Una construcción explícita con la misma función contadora fue dada por Cilleruelo. Ver el capítulo 2 para más detalles.

El propio Erdős demostró que la conjetura anterior no es cierta para $\alpha = 1/2$ pero demostró la existencia de una sucesión infinita de Sidon con

$$\limsup_{x\to\infty} \frac{A(x)}{\sqrt{x}} \geqslant c$$

con $c = 1/2$. Posteriormente Krückeberg (1961) lo demostró para $c = 1/\sqrt{2}$.

**Conjetura 4.10** (Erdős). *Existe una sucesión infinita de Sidon con*

$$\limsup_{x\to\infty} \frac{A(x)}{\sqrt{x}} = 1.$$

Erdős también observó que la conjetura 4.10 seguiría de la siguiente.

**Conjetura 4.11** (Erdős). *Dados elementos $a_1, \dots, a_k$ de un conjunto de Sidon, existe para todo $n$ un conjunto de Sidon $A \subset [1, n]$ con $|A| \sim \sqrt{n}$ que contiene $a_1, \dots, a_k$.*

Es posible que ninguna de estas dos conjeturas sea cierta.

Sea $A_x = A \cap [1, x]$. En el Ejercicio 3.2 se pide demostrar que si $A$ es una sucesión infinita con $|A_x - A_x| \sim |A_x|^2$ entonces $\liminf_{x\to\infty} \frac{A(x)}{\sqrt{x}} = 0$. Sería interesante saber si se puede conseguir la misma conclusión asumiendo una hipótesis análoga para $A_x + A_x$.

**Problema 4.2.** *Sea $A$ una sucesión infinita con $|A_x + A_x| \sim |A_x|^2/2$. ¿Es cierto que $\liminf_{x\to\infty} \frac{A(x)}{\sqrt{x}} = 0$?*

### 4.5 Sucesiones $B_h$.

La diferencia esencial entre las sucesiones de Sidon (sucesiones $B_2$) y las sucesiones $B_h$ con $h \geqslant 3$, radica en que, a diferencia de las primeras, las sucesiones $B_h$ con $h \geqslant 3$ no se pueden caracterizar en términos de sus diferencias. Eso hace que muchos resultados, que son conocidos para $h = 2$, se desconozcan para $h \geqslant 3$. El primero de ellos se refiere al máximo tamaño de un conjunto $B_h$ en $[1, n]$. Mientras que es bien conocido que $F_2(n) \sim \sqrt{n}$, se desconoce el comportamiento asintótico de $F_h(n)$. Ni siquiera se sabe si tiene.

**Problema 4.3.** *Hallar el valor asintótico de $F_h(n)$ para $h \geqslant 3$.*

Hay construcciones que demuestran que $F_h(n) \geqslant n^{1/h}(1+o(1))$, pero las cotas superiores son bastante peores. Por ejemplo para $h = 4$, la mejor cota superior se debe a Ben Green: $F_4(n) \leqslant (7n)^{1/4}(1 + o(1))$. Probablemente, $F_h(n) \sim n^{1/h}$.

Respecto a las sucesiones $B_h$ infinitas sucede algo parecido. Si $A$ es una sucesión $B_h$ con $h$ par, entonces $B = A + \cdots + A$ es casi una sucesión de Sidon en el sentido de que $|B - B| \sim |B|^2$ y se puede demostrar que en ese caso $\liminf_{x\to\infty} \frac{A(x)}{x^{1/h}} = 0$. Sin embargo si $h$ es impar ese argumento no funciona y el resultado análogo se desconoce.

**Conjetura 4.12.** *Si $A$ es una sucesión $B_h$ infinita entonces*

$$\liminf_{x\to\infty} A(x)/x^{1/h} = 0.$$

La conjetura se ha demostrado para $h$ par.

**4.6 Conjuntos de Sidon con condiciones adicionales.** Erdős consideró conjuntos $A$ que no son de Sidon pero que satisfacen $|A + A| \sim |A|^2/2$. A estos conjuntos los llamó conjuntos quasi-Sidon.

**Problema 4.4.** *Dar estimaciones no triviales de*

$$Q(n) = \max |A| : \; A \subset [1, n], \; A \textit{ es quasi-Sidon}.$$

Se sabe que

$$1.154\cdots = \frac{2}{\sqrt{3}}(1 + o(1)) \leqslant \frac{Q(n)}{\sqrt{n}} \leqslant \left(\frac{1}{4} + \frac{1}{(\pi + 2)^2}\right)(1 + o(1)) = 1.863\ldots.$$

La cota inferior se debe a una construcción de Erdős y Freud (1991) y la cota superior a Pikhurko (2006).

**Conjetura 4.13.** *Demostrar que si $A \subset \{1, \ldots, n\}$ es de Sidon y convexo, entonces $|A| = o(\sqrt{n})$.*

Esta interesante conjetura la escuché (creo que a Ruzsa) en el workshop que Ruzsa organizó en Budapest en el año 2000. Se dice que una sucesión es convexa si las diferencias entre dos términos consecutivos de la sucesión son crecientes. Por ejemplo, la sucesión de los cuadrados es una sucesión convexa.

**Problema 4.5.** *Construir un conjunto, lo más grande posible, $A \subset \{1, \ldots, n\}$ que sea de Sidon y convexo.*

No es difícil demostrar que si $A \subset [1, n]$ es una sucesión de Sidon formada por cuadrados entonces $|A| = o(\sqrt{n})$.

**Problema 4.6.** *¿Es cierto que para todo $\epsilon > 0$ existe un conjunto de Sidon en $[1, n]$ formado por cuadrados y de tamaño $|A| \gg n^{1/2-\epsilon}$?*

No es difícil demostrar que existe uno de tamaño $|A| \geqslant n^{1/3-o(1)}$. Bastante más complicado, aunque se sabe cierto, es la demostración de que existe un conjunto de Sidon de cuadrados $A \subset [1, n]$ tal que $|A| \gg n^{1/3}$. Probablemente no se pueda mejorar el exponente $1/3$. La razón para sospechar esto es un trabajo reciente de Saxton y Thomason (2015). Uno de sus corolarios es que si elegimos un conjunto de $\sqrt{n}$ elementos al azar en $[1, n]$, con probabilidad tendiendo a 1, el mayor conjunto de Sidon que contiene tiene tamaño $n^{1/3+o(1)}$. Si los cuadrados en $[1, n]$ se comportan como un conjunto aleatorio para este problema entonces no se debería esperar que contuvieran un conjunto de Sidon de tamaño $n^{1/3+\epsilon}$.

Komlós, Sulyok y Szemeredi (1975) demostraron que cualquier conjunto de $n$ elementos contiene un conjunto de Sidon de tamaño $|A| \gg \sqrt{n}$. Erdős hizo la siguiente conjetura con respecto a este problema.

**Conjetura 4.14** (Erdős)**.** *Todo conjunto de enteros de $n$ elementos contiene un conjunto de Sidon con $\sim \sqrt{n}$ elementos.*

Mi opinión es que esta conjetura es falsa.

**4.7 Bases y sucesiones de Sidon.** Una de las conjeturas más importantes de la teoría combinatoria de números es la que se conoce como Conjetura de Erdős–Turán.

**Conjetura 4.15.** *Si $A$ es una base asintótica de orden 2 entonces su función de representación no está acotada.*

La siguiente conjetura, conocida como conjetura fuerte de Erdős–Turán, implica la anterior porque si $A$ es una base de orden 2 entonces $A(x) \gg x^{1/2}$.

**Conjetura 4.16.** *Si $A(x) \gg x^{1/2}$ entonces su función de representación no está acotada.*

En la otra dirección, Erdős conjeturó lo siguiente.

**Conjetura 4.17** (Erdős)**.** *Existe una sucesión de Sidon que es base asintótica de orden 3.*

Esta conjetura parece difícil pero recientemente se ha demostrado Cilleruelo (2015) que para todo $\epsilon > 0$ existe una sucesión de Sidon $A$ tal que todo $n$ suficientemente grande se puede escribir de la forma $n = a_1 + a_2 + a_3 + a_4$ con $a_1, a_2, a_3, a_4 \in A$ y $a_4 < n^{\epsilon}$.

## Referencias

Javier Cilleruelo†
Universidad Autónoma de Madrid
Instituto de Ciencias Matemáticas (ICMAT)

# INTRODUCCIÓN A MÉTODOS DE ANÁLISIS DE FOURIER EN COMBINATORIA ARITMÉTICA

Pablo Candela

**Resumen**

Este capítulo ofrece una introducción al análisis de Fourier sobre grupos abelianos finitos, orientada hacia el estudio de ciertas aplicaciones de esta área en teoría de números. Las aplicaciones en cuestión, que incluyen los teoremas de Roth y de Meshulam, y el lema de Bogolyubov, constituyen resultados importantes del área de teoría de números conocida como la combinatoria aritmética.

## Índice general

## Introducción

Desde sus comienzos hace aproximadamente dos siglos, el análisis de Fourier ha adquirido una importancia fundamental y creciente en diversas áreas matemáticas, tanto puras como aplicadas (Chamizo s.f.). En teoría de números, varias de estas aplicaciones constituyen resultados clásicos; entre estas se encuentran, por ejemplo, las aplicaciones del método del círculo como el teorema de Vinográdov, que nos dice que todo entero impar suficientemente grande es la suma de tres números primos (resultado cuya reciente mejora por Helfgott resultó en una prueba de la conjetura débil de Goldbach (Helfgott 2013)), o el problema de Waring, que consiste en determinar el número de representaciones de un entero como suma de un número dado de potencias $k$-ésimas (Chamizo, Cristóbal y Ubis 2006).

Alrededor de la mitad del siglo pasado, empezaron a darse aplicaciones combinatorias del análisis de Fourier en teoría de números que formaron algunos de los primeros resultados principales del área hoy conocida como *combinatoria aritmética*. Hoy en día esta área conoce un desarrollo muy activo, dentro del cual el análisis de Fourier ha generado un amplio abanico de métodos y aplicaciones. En este curso estudiaremos algunos ejemplos centrales de estos métodos y algunas de sus aplicaciones principales. Para ello nos concentraremos en el análisis de Fourier sobre grupos abelianos finitos (también llamado análisis de Fourier *discreto*).

Empezamos en la sección siguiente introduciendo los conceptos básicos. En la sección 2 estudiaremos algunas de las ideas principales que dan a esta teoría su utilidad en combinatoria aritmética. En la última sección daremos ejemplos concretos del uso de estas ideas, demostrando algunos resultados centrales en el área, como el teorema de Roth y el lema de Bogolyubov.

## 1 Análisis de Fourier en grupos abelianos finitos – teoría básica

Para todo conjunto finito $X$ y toda función $f : X \to \mathbb{C}$, denotamos por $|X|$ la cardinalidad de $X$ y por $\mathbb{E}_{x\in X} f(x)$ el promedio de $f$ sobre $X$, es decir

$$\mathbb{E}_{x\in X} f(x) = \frac{1}{|X|}\sum_{x\in X} f(x).$$

Las funciones de valores complejos $f : X \to \mathbb{C}$ forman un espacio vectorial, que denotamos por $\mathbb{C}^X$, de dimensión $N = |X|$. Se puede definir la operación de producto interno (o producto escalar) siguiente:

$$\langle\cdot,\cdot\rangle \;:\; \mathbb{C}^X \times \mathbb{C}^X \to \mathbb{C}, \quad (f,g) \;\mapsto\; \langle f,g\rangle = \mathbb{E}_{x\in X} f(x)\overline{g(x)},$$

obteniendo así el *espacio vectorial con producto interno*[1] $(\mathbb{C}^X, \langle\cdot,\cdot\rangle)$. Una noción muy útil relacionada con estos espacios es la de *base ortonormal*. Un conjunto de elementos $\{v_x : x \in X\}$ es una base ortonormal en $\mathbb{C}^X$ si satisface la condición siguiente:

$$\text{(1-1)} \qquad \text{para todo } x,y\in X,\ \text{tenemos } \langle v_x, v_y\rangle = \delta_{xy} = \begin{cases} 1, & x=y \\ 0, & x\neq y \end{cases}.$$

De esta condición se deduce fácilmente que los $N$ elementos $v_x$ son linealmente independientes, luego forman en efecto una base para $\mathbb{C}^X$. Toda función en $\mathbb{C}^X$ se puede por lo tanto descomponer de manera única como combinación lineal de estos elementos. La relación (1-1) permite además expresar los coeficientes de esta descomposición simplemente.

**Lema 1.1.** *Sea $\{v_x : x \in X\}$ una base ortonormal de $\mathbb{C}^X$. Entonces para toda función $f : X \to \mathbb{C}$ tenemos $f = \sum_{x\in X} \langle f, v_x\rangle\, v_x$.*

*Prueba.* Existen coeficientes $c_x \in \mathbb{C}$ tales que $f = \sum_{x\in X} c_x v_x$. Dado cualquier $x_0 \in X$, tomando el producto interno con $v_{x_0}$ de ambos lados de esta ecuación, y usando (1-1), deducimos que $c_{x_0} = \langle f, v_{x_0}\rangle$. □

Sea $G$ un grupo abeliano finito. En $\mathbb{C}^G$ (y en general en $\mathbb{C}^X$) existen muchas bases ortonormales. Entre estas, hay una base específica que tiene propiedades adicionales muy útiles, y que constituye una de las nociones fundamentales del análisis de Fourier discreto. Esta es la base de los *caracteres*

[1] Más precisamente, es un espacio hilbertiano, pero no usaremos este hecho.

sobre $G$.

Denotemos por $\mathbb{C}^\times$ el grupo de números complejos nonulos con multiplicación.

**Definición 1.2.** Un *carácter* sobre $G$ es un homomorfismo $G \to \mathbb{C}^\times$.

Esta definición implica que todo carácter debe tomar todos sus valores en $\{z \in \mathbb{C} : |z| = 1\}$. En efecto, si hubiera un carácter $\chi : G \to \mathbb{C}^\times$ y un elemento $x \in G$ con $|z| = |\chi(x)| \neq 1$, entonces tendríamos $1 = \chi(0) = \chi(|G|\,x) = z^{|G|} \neq 1$, una contradicción.

Podemos equipar el conjunto de caracteres sobre $G$ con la operación de multiplicación punto a punto: si $\chi_1, \chi_2$ son caracteres, el carácter producto $\chi_1\chi_2$ se define para todo $x \in G$ por $\chi_1\chi_2(x) = \chi_1(x)\chi_2(x)$. Obtenemos así un grupo abeliano, llamado el *dual* de $G$ y denotado por $\widehat{G}$. El elemento neutro de $\widehat{G}$ es el carácter que manda todo $x \in G$ a 1, llamado el *carácter principal*. Tenemos también que el inverso de un carácter $\chi$ es el carácter $x \mapsto \overline{\chi(x)} = 1/\chi(x)$ (donde $\overline{z}$ denota el conjugado de $z \in \mathbb{C}$).

Vamos a estudiar el grupo dual en más detalle, en particular para demostrar que sus elementos forman en efecto una base ortonormal en $\mathbb{C}^G$. Empezamos con las observaciones siguientes.

**Lema 1.3.** *Sea $G$ un grupo abeliano finito, y sea $n$ el exponente*[2] *de $G$. Entonces*

1. *Todo carácter $\chi \in \widehat{G}$ toma sus valores en $\{z \in \mathbb{C} : z^n = 1\}$.*
2. *Si $G$ es un grupo cíclico $\mathbb{Z}_N = \mathbb{Z}/N\mathbb{Z}$, entonces*

$$\widehat{G} = \{x \mapsto \exp(2\pi i\, rx/N) : r \in \mathbb{Z}_N\}.$$

Dejamos la prueba como ejercicio. La parte (2) de este lema nos da una descripción explícita de los caracteres sobre $\mathbb{Z}_N$: a todo $\chi \in \widehat{\mathbb{Z}_N}$ le corresponde un único elemento $r \in \mathbb{Z}_N$ tal que $\chi(x) = \exp(2\pi i\, rx/N)$. Llamaremos este elemento $r$ la *frecuencia* de $\chi$. Denotando el carácter de frecuencia $r$ por $\chi_r$, tenemos que la función $r \mapsto \chi_r$ es un isomorfismo de grupos $\mathbb{Z}_N \to \widehat{\mathbb{Z}_N}$.

Resulta que esta descripción de los caracteres sobre $\mathbb{Z}_N$ tiene una generalización muy satisfactoria a todo grupo abeliano finito. Para dar esta descripción usaremos la notación $\mathbb{T}$ para el grupo circular $\mathbb{R}/\mathbb{Z}$, y la notación $e$ para la función $\mathbb{T} \to \mathbb{C}$, $\theta \mapsto e(\theta) = \exp(2\pi i\, \theta)$. Podemos ahora reescribir todo carácter $\chi$ sobre $G = \mathbb{Z}_N$ de la forma $\chi(x) = e(rx/N)$.

Para nuestra descripción general de los caracteres, conviene primero generalizar la función $(r, x) \mapsto rx/N \mod 1$.

**Definición 1.4.** Sea $G$ un grupo abeliano. Una *forma bilineal* de $G^2$ a $\mathbb{T}$ es una función $G \times G \to \mathbb{T}$, $(r, x) \mapsto r \cdot x$ con la propiedad siguiente: para todo $r \in G$, la función $x \mapsto r \cdot x$ es un homomorfismo, y para todo $x \in G$ la función $r \mapsto r \cdot x$ es un homomorfismo. Decimos que la forma bilineal es *no degenerada* si para todo $r \in G \setminus \{0\}$ existe $x \in G$ tal que $r \cdot x \neq 0$, y para todo $x \in G \setminus \{0\}$ existe $r$ tal que $r \cdot x \neq 0$. Decimos que la forma es *simétrica* si $r \cdot x = x \cdot r$ para todo $r, x \in G$.

El teorema fundamental de los grupos abelianos finitos nos dice que $G$ es isomorfo a una suma directa de grupos cíclicos, $G \cong \mathbb{Z}_{N_1} \oplus \cdots \oplus \mathbb{Z}_{N_t}$. Usando este teorema, podemos completar la descripción general del grupo dual.

[2] El exponente de $G$ es el menor entero positivo $m$ tal que para todo $g \in G$ se tiene $mg = 0$.

**Proposición 1.5.** *Sea $G$ un grupo abeliano finito. Entonces $\widehat{G} \cong G$. Existe una forma bilineal de $G$ a $\mathbb{T}$ simétrica y no degenerada $(r, x) \mapsto r \cdot x$. Para todo carácter $\chi \in \widehat{G}$ existe un único elemento $r \in G$ tal que $\chi(x) = e(r \cdot x)$.*

En particular, denotando el carácter $x \mapsto e(r \cdot x)$ de frecuencia $r$ por $e_r$, tenemos que la función $r \mapsto e_r$ es un isomorfismo explícito de $G$ a $\widehat{G}$.

*Prueba.* Empezamos observando los dos hechos siguientes, fácilmente demostrados (ejercicio 2): primero, si $G \cong H$ entonces $\widehat{G} \cong \widehat{H}$; segundo, para grupos abelianos finitos $G_1, G_2$ cualesquiera, el dual de $G_1 \oplus G_2$ es isomorfo a $\widehat{G_1} \oplus \widehat{G_2}$. Combinando esto con el teorema fundamental de los grupos abelianos finitos y el hecho que $\widehat{\mathbb{Z}_N} \cong \mathbb{Z}_N$, deducimos que $\widehat{G} \cong G$.

Para encontrar la forma bilineal deseada, podemos suponer que $G = \mathbb{Z}_{N_1} \oplus \cdots \oplus \mathbb{Z}_{N_t}$. Sobre cada grupo $\mathbb{Z}_{N_j}$ tenemos ya una forma bilineal con las propiedades deseadas, a saber la forma $(r, x) \mapsto rx/N_j \mod 1$. Por otro lado, se verifica fácilmente que si $G_1, G_2$ tienen cada uno una forma $\cdot$ con estas propiedades, entonces la función $G_1 \oplus G_2 \to \mathbb{R}/\mathbb{Z}$, $((r_1, r_2), (x_1, x_2)) \mapsto r_1 \cdot x_1 + r_2 \cdot x_2$ también es una forma bilineal adecuada. Razonando por inducción sobre $j$, encontramos la forma sobre $G$.

Para ver la última parte, nótese que, por un lado, cada función $x \mapsto e(r \cdot x)$, $r \in G$ es un carácter, y por otro lado no pueden haber más caracteres sobre $G$. □

De ahora en adelante, supondremos siempre que un grupo abeliano finito $G$ viene equipado con una tal forma bilineal. [3]

Verificamos ahora que los caracteres forman efectivamente una base ortonormal en $\mathbb{C}^G$.

**Proposición 1.6.** *Sea $G$ un grupo abeliano finito. Entonces tenemos $\langle e_r, e_s \rangle = \delta_{rs}$ para todo $r, s \in G$.*

*Prueba.* Tenemos $\langle e_r, e_s \rangle = \mathbb{E}_{x \in G}\, e((r - s) \cdot x)$. Este promedio tiene valor 1 cuando $r = s$, con lo cual solo nos queda por probar que tiene valor 0 cuando $r \neq s$. Sea $t = r - s \neq 0$. Siendo no degenerada la forma bilineal $\cdot$, existe $x_0 \in G$ tal que $t \cdot x_0 \neq 0$. Por otra parte, tenemos que $e(t \cdot x_0)\, \mathbb{E}_{x \in G}\, e(t \cdot x) = \mathbb{E}_{x \in G}\, e(t \cdot (x + x_0)) = \mathbb{E}_{x \in G}\, e(t \cdot x)$, luego $(e(t \cdot x_0) - 1)\, \mathbb{E}_{x \in G}\, e(t \cdot x) = 0$. Como $e(t \cdot x_0) - 1 \neq 0$, deducimos que, efectivamente, $\mathbb{E}_{x \in G}\, e(t \cdot x) = 0$. □

Ha llegado el momento de definir la transformada de Fourier y describir sus primeras propiedades básicas.

**Definición 1.7** (Transformada de Fourier discreta)**.** Para toda función $f \in \mathbb{C}^G$ y todo $\chi \in \widehat{G}$, la *transformada de Fourier* de $f$ es la función $\widehat{f} : \widehat{G} \to \mathbb{C}$ definida por

$$\widehat{f}(\chi) = \langle f, \chi \rangle = \mathbb{E}_{x \in G} f(x) \overline{\chi(x)}.$$

Los valores $\widehat{f}(\chi)$ de la transformada se llaman los *coeficientes* de Fourier de $f$.

**Proposición 1.8.** *Sean $f, g \in \mathbb{C}^G$. Entonces tenemos*

[3] Se puede verificar que para nuestro uso de esta teoría no tiene importancia qué forma particular escogemos. Si se quiere ser más explícito, se puede elegir la forma $r \cdot x = \sum_{r=1}^{t} \phi(r)_j\, \phi(x)_j / N_j \mod 1$, donde $\phi$ es un isomorfismo $G \to \mathbb{Z}_{N_1} \oplus \cdots \oplus \mathbb{Z}_{N_t}$.

- *Fórmula de inversión:*

$$f(x) = \sum_{\chi \in \widehat{G}} \widehat{f}(\chi)\, \chi(x). \tag{1-2}$$

- *Teorema de Plancherel:*

$$\langle f, g \rangle = \sum_{\chi \in \widehat{G}} \widehat{f}(\chi)\, \overline{\widehat{g}(\chi)}. \tag{1-3}$$

- *Identidad de Parseval:*

$$\mathbb{E}_{x \in G} |f(x)|^2 = \sum_{\chi \in \widehat{G}} |\widehat{f}(\chi)|^2. \tag{1-4}$$

*Prueba.* La fórmula de inversión es un caso especial del lema 1.1. El teorema de Plancherel se obtiene substituyendo las fórmulas de inversión para $f$ y $g$ dentro de $\langle f, g \rangle$ y usando la proposición 1.6. El caso especial $f = g$ nos da la identidad de Parseval. □

**Nota 1.9.** Formular los coeficientes de Fourier como lo hacemos, usando el promedio $\mathbb{E}_G$, supone que estamos usando la *medida de probabilidad uniforme* sobre $G$, es decir la medida que a cada conjunto $A \subset G$ asigna el valor $|A|/|G|$. La fórmula de inversión (1-2), compatible con esta definición de los coeficientes, supone en cambio que la medida que usamos sobre el grupo dual $\widehat{G}$ es la *medida de conteo*, que da a $A$ el valor $|A|$. Siguiendo esta convención, es habitual definir la norma euclidiana sobre $\mathbb{C}^G$ usando la medida de probabilidad uniforme, escribiendo $\|f\|_{L^2(G)} = (\mathbb{E}_{x \in G} |f(x)|^2)^{1/2}$, y sobre $\mathbb{C}^{\widehat{G}}$ usando la medida de conteo, escribiendo $\|g\|_{\ell^2(\widehat{G})} = \left( \sum_\chi |g(\chi)|^2 \right)^{1/2}$. La identidad de Parseval se puede entonces escribir $\|f\|_{L^2(G)} = \|\widehat{f}\|_{\ell^2(\widehat{G})}$.

La transformada de Fourier de $f$ es una función $\widehat{f}$ definida sobre $\widehat{G}$. Sin embargo, en la práctica la trataremos a menudo como una función sobre $G$, usando el isomorfismo de la proposición 1.5, escribiendo $\widehat{f}(r)$ en vez de $\widehat{f}(e_r)$.

Usaremos también la famosa desigualdad siguiente.

**Proposición 1.10** (Desigualdad de Cauchy–Schwarz). *Para funciones $f, g \in \mathbb{C}^G$ cualesquiera,*

$$|\langle f, g \rangle| \leqslant \|f\|_{L^2(G)} \, \|g\|_{L^2(G)}. \tag{1-5}$$

(Nótese la versión equivalente: $\left| \sum_{x \in G} f(x) \overline{g(x)} \right| \leqslant \|f\|_{\ell^2(G)} \, \|g\|_{\ell^2(G)}$.)

A continuación empezamos a ilustrar la utilidad combinatoria de la transformada de Fourier.

## 2 Primeros usos en combinatoria aritmética

Uno de los temas centrales de la combinatoria aritmética, que se remonta a los orígenes de esta área, consiste en estudiar bajo qué condiciones (las más naturales y débiles posibles) un subconjunto $A$ de un grupo abeliano debe contener configuraciones aritméticas de un tipo dado. Estas configuraciones consisten generalmente en tuplas de elementos de $A$ que satisfacen ecuaciones lineales prescritas. En esta dirección, un resultado central en el área es el famoso teorema de Szemerédi, conjeturado por (Erdős y Turán 1936). Denotemos por $[N]$ el conjunto $\{1, 2, \ldots, N\}$.

**Teorema 2.1** (Szemerédi)**.** *Sea $\delta > 0$ y sea $k$ un entero positivo. Entonces existe $N_0 > 0$ tal que para todo entero $N > N_0$ y todo conjunto de enteros $A \subset [N]$ de cardinalidad $|A| \geqslant \delta N$, existe una progresión aritmética de longitud $k$ incluida en $A$.*

Nótese que una tal progresión es una $k$-tupla de elementos de $A$ que satisfacen cierto sistema de $k-2$ ecuaciones lineales. Por ejemplo para $k = 3$ el sistema consiste en la única ecuación $x_1 - 2x_2 + x_3 = 0$.

Szemerédi dio la primera prueba de este teorema en 1975 usando argumentos complicados de teoría de grafos (Szemerédi 1975). El caso notrivial más simple del teorema, a saber el caso $k = 3$, había sido demostrado ya por (Roth 1953). La prueba de Roth usa el análisis de Fourier sobre grupos $\mathbb{Z}_N$, y las ideas principales que introdujo han pasado a formar parte importante del conjunto de herramientas en el área. Por ello, en la sección siguiente estudiaremos el teorema de Roth como primera aplicación central. Con vistas a este objetivo, en esta sección presentaremos algunas ideas, en su mayoría relacionadas con la prueba del teorema de Roth, como principios generales del uso combinatorio de la transformada de Fourier. [4]

Dado un conjunto $X$ y un subconjunto $A \subset X$, denotaremos por $1_A$ la *función indicadora* de $A$ sobre $X$, definida como sigue: $1_A(x) = 1$ si $x \in A$, y $1_A(x) = 0$ si $x \in X \setminus A$.

Dada una ecuación lineal, y dado $A \subset G$, nos interesa contar cuantas soluciones de esta ecuación hay en $A$. Por ejemplo, en el caso del teorema de Roth, esto consiste esencialmente en analizar la cantidad siguiente:

$$\mathbb{E}_{\substack{x_1,x_2,x_3 \in G:\\ x_1-2x_2+x_3=0}} 1_A(x_1)\, 1_A(x_2)\, 1_A(x_3). \tag{2-1}$$

Más generalmente, podemos considerar funciones $f_1, \ldots, f_t : G \to \mathbb{C}$ arbitrarias (no necesariamente funciones indicadoras de conjuntos) y considerar promedios de tales funciones tomados sobre el conjunto de soluciones a una ecuación lineal $c_1x_1 + \cdots + c_tx_t = 0$ con coeficientes $c_i$ enteros:

$$\mathbb{E}_{\substack{x_1,\ldots,x_t \in G:\\ c_1x_1+\cdots+c_tx_t=0}} f_1(x_1) \cdots f_t(x_t). \tag{2-2}$$

El primer principio básico que vamos a estudiar es que la transformada de Fourier permite simplificar promedios de este tipo.

## 2.1 Simplificación de promedios relativos a ecuaciones lineales.

Empezamos con el promedio (2-1). Siempre y cuando el orden de $G$ sea impar, usando la fórmula de inversión (1-2) para $1_A$ obtenemos la ecuación siguiente:

[4]Nótese, no obstante, que los métodos de Roth no permiten demostrar el teorema de Szemerédi para $k > 3$. Más generalmente, los métodos Fourier-analíticos conciernen sobre todo a las configuraciones dadas por sistemas de ecuaciones lineales llamados *sistemas de complejidad 1*. (Véase el ejercicio 10.)

$$\mathbb{E}_{\substack{x_1,x_2,x_3\in G:\\ x_1-2x_2+x_3=0}} 1_A(x_1)\,1_A(x_2)\,1_A(x_3) = \sum_r \widehat{1_A}(r)^2\,\widehat{1_A}(-2r). \tag{2-3}$$

La suma a la derecha en (2-3) simplifica el promedio original, en el sentido que este promedio involucra más de una variable mientras que la suma es sobre la única variable $r$. [5] Este tipo de simplificación es muy útil, como veremos más adelante. Además, es bastante general.

**Proposición 2.2.** *Sea $c_1x_1+\cdots+c_tx_t=0$ una ecuación lineal con coeficientes enteros $c_i$. Sea $G$ un grupo abeliano finito tal que la multiplicación $x\mapsto c_ix$ es sobreyectiva para todo $i\in[t]$, y sean $f_1,\dots,f_t$ funciones $G\to\mathbb{C}$ arbitrarias. Entonces tenemos la fórmula siguiente:*

$$\mathbb{E}_{\substack{x_1,\dots,x_t\in G:\\ c_1x_1+\cdots+c_tx_t=0}} f_1(x_1)\cdots f_t(x_t) \;=\; \sum_{r\in G}\widehat{f_1}(c_1r)\cdots\widehat{f_t}(c_tr). \tag{2-4}$$

Este es el resultado principal de esta subsección. Para demostrarlo, utilizaremos la operación de convolución, una operación que juega un papel crucial en esta área, sobre todo por su relación con la transformada de Fourier.

**Definición 2.3** (Convolución)**.** Sean $f,g\in\mathbb{C}^G$. La *convolución* de $f$ y $g$ es la función $G\to\mathbb{C}$ denotada por $f*g$ y definida como sigue:

$$f*g(x) \;=\; \mathbb{E}_{y\in G}\, f(y)\,g(x-y).$$

La importancia de esta operación va mucho más allá de su utilidad en la prueba de la proposición 2.2. Por ejemplo, más adelante veremos que es muy útil para estudiar la estructura de conjuntos suma $A+B=\{a+b : a\in A, b\in B\}$, ya que tenemos la relación siguiente:

$$A+B=\operatorname{supp}(1_A*1_B),$$

donde $\operatorname{supp}(f)=\{x\in G : f(x)\neq 0\}$ es el *soporte* de la función $f\in\mathbb{C}^G$.

Por ahora nos concentramos en demostrar la proposición 2.2. El hecho principal que usaremos para ello es que la transformada de Fourier convierte una convolución de dos funciones en simple producto de sus coeficientes:

$$\text{para todo } r\in G \text{ y } f,g\in\mathbb{C}^G,\quad \text{tenemos}\quad \widehat{f*g}(r) \;=\; \widehat{f}(r)\,\widehat{g}(r). \tag{2-5}$$

Esta *fórmula de producto* se verifica aplicando la fórmula de inversión (1-2) para $f$ y $g$, y la ortonormalidad de los caracteres (véase el ejercicio 3).

*Prueba de la proposición 2.2.* La idea es reescribir el promedio

$$\mathbb{E}_{\substack{x_1,\dots,x_t\in G:\\ c_1x_1+\cdots+c_tx_t=0}} f_1(x_1)\cdots f_t(x_t)$$

[5] Esta ecuación es un caso especial de la fórmula (2-4) probada a continuación, pero como primer ejemplo puede ser instructivo intentar demostrar la ecuación directamente.

de una forma que revele que se trata de una convolución iterada. Hacemos primero el cambio de variables siguiente: $y_1 = c_1 x_1, \dots, y_t = c_t x_t$. Nótese que, siendo cada función $x \mapsto c_i x$ sobreyectiva (luego biyectiva) de $G$ a $G$, tiene una función inversa, que denotaremos por $c_i^{-1}$. El promedio original queda escrito de la forma siguiente:

$$\mathbb{E}_{\substack{y_1,\dots,y_t \in G:\\ y_1+\dots+y_t=0}} \; f_1(c_1^{-1} y_1) \;\cdots\; f_t(c_t^{-1} y_t).$$

Definiendo las funciones $f_1'(y) = f_1(c_1^{-1} y), \dots, f_t'(y) = f_t(c_t^{-1} y)$, podemos reescribir el último promedio simplemente y ver que es efectivamente una convolución iterada:

$$\mathbb{E}_{\substack{y_1,\dots,y_t \in G:\\ y_1+\dots+y_t=0}} \; f_1'(y_1) \;\cdots\; f_t'(y_t) = f_1' * \cdots * f_t'(0).$$

La fórmula de inversión nos dice que esto es igual a $\sum_r \widehat{f_1' * \cdots * f_t'}(r)$, y la fórmula de producto (aplicada $t-1$ veces) nos dice que esto es igual a $\sum_r \widehat{f_1'}(r)\widehat{f_2'}(r)\cdots\widehat{f_t'}(r)$. Finalmente, verificamos mediante un simple cálculo[6] que $\widehat{f_i'}(r) = \widehat{f_i}(c_i r)$. □

**Nota 2.4.** Esta proposición se puede demostrar de otra manera usando la fórmula de Poisson (véase el ejercicio 4).

¿Cómo se puede utilizar la simplificación de promedios dada en (2-4)? Para responder debemos estudiar, primero, cómo la estructura combinatoria de un conjunto $A$ se puede ver reflejada en los coeficientes $\widehat{1_A}(r)$, y segundo, cómo una información sobre estos coeficientes se puede transformar a su vez en nueva información combinatoria provechosa.

## 2.2 Transformar datos combinatorios en datos Fourier-analíticos, y viceversa.

Empezamos con un conjunto $A \subset G$ dado. Denotamos por $\alpha$ la *densidad* de $A$ en $G$, a saber $\alpha = |A|/|G|$. Una primera observación es que el llamado *coeficiente principal* $\widehat{1_A}(0)$ (correspondiente al carácter principal, a saber el carácter con valor constante igual a 1) es igual a $\alpha$, independientemente de la estructura combinatoria de $A$:

$$\widehat{1_A}(0) = \mathbb{E}_{x \in G} 1_A(x) = \alpha.$$

Por lo tanto, para estudiar esta estructura hemos de concentrarnos en los coeficientes noprincipales $\widehat{1_A}(r)$, $r \neq 0$. Para ello trabajaremos con la llamada *función nivelada* de $A$, denotada por $f_A$ y definida como sigue:

$$f_A(x) = 1_A(x) - \alpha.$$

Nótese que tenemos $\widehat{f_A}(r) = \widehat{1_A}(r)$ para todo $r \neq 0$, y $\widehat{f_A}(0) = 0$.

Como mencionamos previamente, nuestro problema combinatorio general es el de analizar el promedio de $1_A$ relativo a una ecuación lineal dada. Para ser concretos tomemos de nuevo el promedio

[6] En el cálculo se usa el hecho que para todo entero $n$ tenemos $r \cdot (nx) = (nr) \cdot x$.

de progresiones aritméticas de longitud 3, que denotaremos de ahora en adelante por $T_3(A)$:

$$T_3(A) = \mathbb{E}_{\substack{x_1,x_2,x_3\in G:\\ x_1-2x_2+x_3=0}} 1_A(x_1)\, 1_A(x_2)\, 1_A(x_3).$$

La fórmula (2-3) nos dice que esto es igual a $\sum_r \widehat{1_A}(r)^2\, \widehat{1_A}(-2r)$. Aquí también es muy útil separar el carácter principal de los demás, como sigue:

$$T_3(A) = \alpha^3 + \sum_{r\neq 0} \widehat{1_A}(r)^2\, \widehat{1_A}(-2r) = \alpha^3 + \sum_r \widehat{f_A}(r)^2\, \widehat{f_A}(-2r). \tag{2-6}$$

Con esta ecuación se precisa nuestro problema: hemos de analizar la diferencia entre un promedio (aquí $T_3(A)$) y lo que podemos llamar su *valor principal* (aquí $\alpha^3$, y para una ecuación más general en $t$ variables, $\alpha^t$). La ecuación también deja claro que en este análisis la *magnitud* de los coeficientes $\widehat{f_A}(r)$ va a jugar un papel importante.

**Definición 2.5** (Uniformidad de Fourier)**.** Definimos la *norma de uniformidad de Fourier* para toda función $f : G \to \mathbb{C}$ como sigue:

$$\|f\|_u = \sup_r \big|\widehat{f}(r)\big|.$$

Como primer ejemplo de la utilidad de esta norma, tenemos el resultado siguiente.

**Proposición 2.6.** *Para todo conjunto $A \subset G$ de densidad $\alpha$ tenemos*

$$\big|\, T_3(A) - \alpha^3 \,\big| \leqslant \|f_A\|_u\, \alpha. \tag{2-7}$$

Este es un caso especial del resultado principal de esta sección, que nos dice esencialmente que la norma $\|\cdot\|_u$ controla la diferencia entre un promedio de tipo (2-2) y su valor principal:

**Proposición 2.7.** *Sea $c_1x_1 + \cdots + c_tx_t = 0$ una ecuación lineal con coeficientes enteros $c_i$, con $t \geqslant 3$. Sea $G$ un grupo abeliano finito tal que la multiplicación $x \mapsto c_ix$ es sobreyectiva para todo $i \in [t]$, y sean $A_1,\dots,A_t$ subconjuntos de $G$ arbitrarios, de densidades $\alpha_1,\dots,\alpha_t$ respectivamente. Entonces tenemos*

$$\Big|\alpha_1\cdots\alpha_t - \mathbb{E}_{\substack{x_1,\dots,x_t\in G:\\ c_1x_1+\cdots+c_tx_t=0}} 1_{A_1}(x_1)\cdots 1_{A_t}(x_t)\Big| \leqslant \|f_{A_1}\|_u\cdots\|f_{A_{t-2}}\|_u\, \alpha_{t-1}^{1/2}\,\alpha_t^{1/2}. \tag{2-8}$$

*En particular tenemos* $\big|\,\alpha^t - \mathbb{E}_{\substack{x_1,\dots,x_t\in G:\\ c_1x_1+\cdots+c_tx_t=0}} 1_A(x_1)\cdots 1_A(x_t)\,\big| \leqslant \|f_A\|_u^{t-2}\,\alpha.$

*Prueba.* La fórmula (2-4) nos dice que el promedio

$$\mathbb{E}_{\substack{x_1,\dots,x_t\in G:\\ c_1x_1+\cdots+c_tx_t=0}} 1_{A_1}(x_1)\cdots 1_{A_t}(x_t)$$

es igual a

$$\alpha_1\cdots\alpha_t + \sum_{r\neq 0} \widehat{1_{A_1}}(c_1r)\cdots\widehat{1_{A_t}}(c_tr) = \alpha_1\cdots\alpha_t + \sum_r \widehat{f_{A_1}}(c_1r)\cdots\widehat{f_{A_t}}(c_tr).$$

Por lo tanto, el lado izquierdo de (2-8) vale como mucho

$$\sum_r \big|\widehat{f_{A_1}}(c_1 r)\big| \cdots \big|\widehat{f_{A_t}}(c_t r)\big| \leqslant$$
$$\leqslant \sup_r \big|\widehat{f_{A_1}}(c_1 r)\big| \cdots \sup_r \big|\widehat{f_{A_{t-2}}}(c_{t-2} r)\big| \sum_r \big|\widehat{f_{A_{t-1}}}(c_{t-1} r)\big|\, \big|\widehat{f_{A_t}}(c_t r)\big|$$
$$\leqslant \|f_{A_1}\|_u \cdots \|f_{A_{t-2}}\|_u \sum_r \big|\widehat{1_{A_{t-1}}}(c_{t-1} r)\big|\, \big|\widehat{1_{A_t}}(c_t r)\big|.$$

Por la desigualdad de Cauchy–Schwarz y la identidad de Parseval, tenemos (usando también la biyectividad de las funciones $r \mapsto c_i r$)

$$\sum_r \big|\widehat{1_{A_{t-1}}}(c_{t-1} r)\big|\, \big|\widehat{1_{A_t}}(c_t r)\big| \leqslant \Big(\sum_r \big|\widehat{1_{A_{t-1}}}(r)\big|^2\Big)^{1/2} \Big(\sum_r \big|\widehat{1_{A_t}}(r)\big|^2\Big)^{1/2}$$
$$= \big(\mathbb{E}_{x\in G}\, 1_{A_{t-1}}(x)^2\big)^{1/2} \big(\mathbb{E}_{x\in G}\, 1_{A_t}(x)^2\big)^{1/2}$$
$$= \alpha_{t-1}^{1/2}\ \alpha_t^{1/2}.$$

□

**Nota 2.8.** Una norma estrechamente relacionada con $\|\cdot\|_u$ y de gran utilidad es la llamada *norma $U^2$ de Gowers*, definida para toda función $f : G \to \mathbb{C}$ por

$$\|f\|_{U^2} = \big(\mathbb{E}_{x,h,k\in G}\, f(x)\ \overline{f(x+h)}\ \overline{f(x+k)}\ f(x+h+k)\big)^{1/4}$$

(Véase el ejercicio 5.)

Este resultado nos permite reducir gran parte de nuestro problema combinatorio inicial al de entender la relación entre la estructura de $A$ y el valor de $\|f_A\|_u$. Empezamos a elucidar esta relación con la pregunta siguiente, un poco imprecisa:

> *¿Qué tipo de estructura del conjunto $A$ garantiza que $\|f_A\|_u$ sea pequeña comparada con $\alpha$?*

Podemos pensar en un coeficiente de Fourier noprincipal $\widehat{1_A}(r)$ geométricamente como la suma en el plano complejo de los números $e_r(x)/|G|$, $x \in A$. Esta suma tendrá menor magnitud cuanto más haya cancelación entre estos números. Para que $\|f_A\|_u$ sea pequeña, tiene que haber una tal cancelación para todo $r \neq 0$, y un modo claro de garantizar esto es que para todo $r \neq 0$ los números $e_r(x)$, $x \in A$ se repartan de manera más o menos uniforme en el círculo unidad. Un momento de reflexión nos conduce a un tipo de estructura que garantiza esto de manera muy natural: un conjunto *$A$ aleatorio de probabilidad $\alpha$*, es decir un conjunto aleatorio tal que para cada $x \in G$ el evento $x \in A$ tiene probabilidad $\alpha$, y estos eventos son independientes.

**Proposición 2.9** (Uniformidad de un conjunto aleatorio)**.** *Sea $A \subset G$ un conjunto aleatorio de probabilidad $\alpha \leqslant 1/2$ y sea $t > 1/4|G|$. Si $\alpha \geqslant 8\log(4t|G|)/|G|$, entonces tenemos con probabilidad al menos $1 - 1/t$ que*

$$\big|\,|A|/|G| - \alpha\,\big| \leqslant 4\sqrt{\log(4t|G|)/|G|}, \quad y \quad \|f_A\|_u \leqslant 4\sqrt{\log(4t|G|)/|G|}.$$

*Prueba.* Usamos la desigualdad de Chernoff en su versión para variables aleatorias complejas. Este resultado dice lo siguiente (véase Tao y Vu (2006, teorema 1.8 y ejercicio 1.3.4)): sean $X_1, \dots, X_n$ variables aleatorias complejas independientes tales que $|X_i - \mathbb{E}(X_i)| \leqslant 1$ para todo $i \in [n]$. Sea $X = X_1 + \cdots + X_n$ y sea $\sigma = \sqrt{\mathrm{Var}(X)}$ la desviación estándar de $X$. Entonces para todo $\lambda > 0$ tenemos

$$\mathbb{P}(|X - \mathbb{E}(X)| \geqslant \lambda\sigma) \leqslant 4\max(e^{-\lambda^2/8}, e^{-\lambda\sigma/4}). \tag{2-9}$$

Dado un elemento $r \in G$ cualquiera, aplicamos (2-9) a las variables $X_x = e(-r\cdot x)1_A(x)$, $x \in G$. Nótese que para la variable $X = \sum_{x\in G} X_x$ tenemos que $X = |G|\,\widehat{1_A}(r)$, que $\mathbb{E}(X) = \alpha\,|G|$ para $r = 0$, y que $\mathbb{E}(X) = 0$ para $r \neq 0$. Tenemos por otro lado que

$$\begin{aligned}\sigma &= \Big(\sum_x \mathrm{Var}(X_x)\Big)^{1/2} = \Big(\sum_x \mathbb{E}|X_x - \mathbb{E}(X_x)|^2\Big)^{1/2}\\ &\in \Big[\Big(\sum_x \mathbb{E}(1_A(x)-\alpha)^2\Big)^{1/2}, \Big(\sum_x \mathbb{E}(1_A(x)+\alpha)^2\Big)^{1/2}\Big]\\ &= \Big[\big(\alpha(1-\alpha)|G|\big)^{1/2}, \big((\alpha+3\alpha^2)|G|\big)^{1/2}\Big] \subset \big[\sqrt{\alpha|G|/2}, \sqrt{3\alpha|G|}\,\big],\end{aligned}$$

donde para la última inclusión usamos que $\alpha \leqslant 1/2$. En particular, si aplicamos (2-9) con $\lambda^2 \leqslant \alpha|G|$ entonces, como $\lambda < 2\sigma$, luego $\lambda^2/8 < \lambda\sigma/4$, tenemos que la cota superior en (2-9) es $4e^{-\lambda^2/8}$. Aplicando pues (2-9) para todo $r$, obtenemos (usando que $\sigma < \sqrt{2|G|}$) que las probabilidades

$$\mathbb{P}\big(\big|\,|A| - \alpha|G|\,\big| \geqslant 2\lambda|G|^{1/2}\big),\quad \mathbb{P}\big(|G|\,|\widehat{1_A}(r)| \geqslant 2\lambda|G|^{1/2}\big),\quad r \neq 0$$

son todas inferiores o iguales a $4e^{-\lambda^2/8}$. Ahora, observemos que la probabilidad del evento $\{\big|\,|A| - \alpha|G|\,\big| < 2\lambda|G|^{1/2}\ \ \&\ \forall\, r \neq 0,\ |G|\,|\widehat{1_A}(r)| < 2\lambda|G|^{1/2}\}$ es 1 menos la probabilidad del evento $\{\big|\,|A| - \alpha|G|\,\big| \geqslant 2\lambda|G|^{1/2}$ o bien $\exists\, r \neq 0,\ |G|\,|\widehat{1_A}(r)| \geqslant 2\lambda|G|^{1/2}\}$. Deducimos que con probabilidad al menos $1 - 4|G|e^{-\lambda^2/8}$ se da que $\big|\,|A|/|G| - \alpha\,\big| \leqslant 2\lambda|G|^{-1/2}$ y que $|\widehat{1_A}(r)| \leqslant 2\lambda|G|^{-1/2}$ para todo $r \neq 0$. Eligiendo $\lambda = \sqrt{8\log(4t|G|)}$, completamos la prueba. □

Nótese que, si denotamos de nuevo por $\alpha$ la densidad de $A$ en $G$, entonces por la identidad de Parseval tenemos $\sqrt{\alpha(1-\alpha)} = \big(\sum_r |\widehat{f_A}|^2\big)^{1/2} \leqslant \|f_A\|_u\,|G|^{1/2}$. De modo que la proposición 2.9 nos dice que los conjuntos aleatorios son extremadamente uniformes (sus normas $\|f_A\|_u$ son muy pequeñas, casi tan pequeñas como el mínimo posible). Junto con la proposición 2.7, esto indica que la norma $\|f_A\|_u$ se puede usar para medir cuanto se parece el conjunto $A$, en sus propiedades combinatorias relativas a ecuaciones lineales, a un conjunto aleatorio de la misma densidad. Se dice que, cuanto más pequeño es el valor de $\|f_A\|_u$, más *quasi-aleatorio* es el conjunto. Nótese, no obstante, que los conjuntos aleatorios no son los únicos ejemplos de conjuntos con norma de uniformidad muy pequeña (este matiz conduce a temas muy interesantes; véase el ejercicio 8).

De todos modos, con la proposición 2.7 nuestro problema inicial queda resuelto en el caso de conjuntos suficientemente quasi-aleatorios. En efecto, si $\|f_A\|_u$ es suficientemente pequeña como función de $\alpha$, entonces cualquier promedio del tipo de (2-8) tiene que ser cercano a su valor principal (que es también aproximadamente el valor esperado de este promedio para un conjunto aleatorio de la misma densidad).

A continuación exploramos el caso no quasi-aleatorio, a partir de la pregunta siguiente:

**Pregunta 2.10.** ¿Suponiendo que $\|f_A\|_u$ es mayor que una función dada de $\alpha$, qué información combinatoria podemos deducir acerca de $A$?

Empecemos buscando algunos ejemplos de conjuntos $A$ que tengan $\|f_A\|_u$ grande (al menos una fracción fija de $\alpha$). Dado que los conjuntos aleatorios fueron tan buenos ejemplos de la propiedad opuesta, es natural buscar ahora entre conjuntos con una estructura "ordenada". Para precisar esto, volvamos a la intuición geométrica que nos condujo al ejemplo del conjunto aleatorio. Vemos pronto que una manera natural de que el promedio $\mathbb{E}_{x\in G} 1_A(x)e(r\cdot x)$ tenga la mayor magnitud posible es que para $x \in A$ los números complejos $e_r(x)$ apunten lo más posible en la misma dirección en el plano (minimizando así su cancelación mutua). La manera óptima de garantizar esto es que $x \mapsto r\cdot x$ sea constante sobre $A$. Esto se da, por ejemplo, si $A$ es una clase lateral del llamado *conjunto anulador* de $r$.

**Definición 2.11.** Dado un conjunto $S \subset G$, el *anulador* de $S$ es el subgrupo de $G$ denotado $S^\perp$ y definido como sigue

$$S^\perp := \{x \in G : r\cdot x = 0 \text{ para todo } r \in S\}.$$

(Cuando $S$ es un singleton $\{r\}$ escribiremos $r^\perp$ en vez de $\{r\}^\perp$.)

Tenemos que $r\cdot x$ es constante para todo $x \in A$ si y sólo si $A$ está contenido en una clase lateral de $r^\perp$. (En ámbitos del análisis de Fourier más generales, donde no se da el isomorfismo $G \cong \widehat{G}$, el conjunto anulador se define como subgrupo del grupo dual: $S^\perp := \{\chi \in \widehat{G} : \chi(r) = 1 \text{ para todo } r \in S\}$.)

Trabajar con estos conjuntos $r^\perp$ nos será de muy gran utilidad siempre y cuando se pueda garantizar que el tamaño de un tal subgrupo es comparable al orden de $G$. Ahora bien, esto no es el caso en cualquier grupo $G$, y por esta razón a partir de ahora es conveniente distinguir algunos tipos diferentes de grupos abelianos finitos.

Una familia de grupos $G$ especialmente ricos en subgrupos anuladores de gran tamaño es la familia de espacios vectoriales finitos $G = \mathbb{F}_p^n$, donde $p$ es un primo fijo y $n$ es la dimensión (que típicamente se toma mucho más grande que $p$). Como veremos en la sección siguiente, esta familia es muy útil, por la riqueza de subgrupos (subespacios) de densidades variadas que contiene, y porque en estos grupos muchos cálculos de análisis de Fourier se reducen a simples argumentos de álgebra lineal. En estos grupos, los subespacios constituyen ejemplos óptimos de conjuntos densos con grandes coeficientes de Fourier noprincipales.

**Ejemplo 2.12** (Transformada de Fourier de un subespacio)**.** *Sea $V$ un subespacio de $\mathbb{F}_p^n$, y sea $N = |\mathbb{F}_p^n| = p^n$. Entonces tenemos*[7]*, para toda frecuencia $r \in \mathbb{F}_p^n$,*

(2-10) $$\widehat{1_V}(r) = \frac{|V|}{N} 1_{V^\perp}(r).$$

[7]Véase el ejercicio 4.

*En particular, para todo* $r \in V^{\perp}$ *la transformada alcanza su valor máximo posible* $\frac{|V|}{N}$.

En el otro extremo tenemos los grupos cíclicos $\mathbb{Z}_N$ con $N$ primo, que no contienen ningún subgrupo notrivial. En particular, en estos grupos un conjunto anulador $r^{\perp}$ con $r \neq 0$ no es más que el singleton $\{0\}$. ¿Qué noción podemos usar en este tipo de grupos para reemplazar los conjuntos anuladores? Como veremos, resulta sumamente útil trabajar con una familia de conjuntos que constituyen *anuladores aproximativos* de caracteres. Estos conjuntos son los llamados *conjuntos de Bohr*.

**Definición 2.13.** Sea $S \subset G$ y sea $\delta \geqslant 0$. El *conjunto de Bohr* con *frecuencias* en $S$ y *radio* $\delta$ es el conjunto $B(S,\delta) = \{x \in G : \|r \cdot x\|_{\mathbb{T}} \leqslant \delta \text{ para todo } r \in S\}$.

Aquí, dado un elemento $\theta \in \mathbb{T}$, representado por un punto del intervalo $[-1/2, 1/2)$, denotamos por $\|\theta\|_{\mathbb{T}}$ la distancia entre $\theta$ y el conjunto $\mathbb{Z}$, es decir $\|\theta\|_{\mathbb{T}} = |\theta|$. Recordamos también las desigualdades básicas

$$4\|\theta\|_{\mathbb{T}} \leqslant |1 - e(\theta)| \leqslant 2\pi\|\theta\|_{\mathbb{T}}. \tag{2-11}$$

Estas desigualdades se pueden comprobar considerando el cociente de longitudes $\frac{|1-e(\theta)|}{2\pi\|\theta\|_{\mathbb{T}}}$ en el plano complejo. En efecto, supongamos que $\theta \in (0, 1/2]$, y notemos que el denominador es la longitud del arco de círculo unitario que va del punto 1 en el plano complejo hasta el punto $e(\theta)$ (estando este punto en la mitad del plano $\mathrm{Im}(z) > 0$). Vemos entonces claramente que el cociente en cuestión es como mucho 1, lo cual nos da la segunda desigualdad en (2-11). Vemos también que el cociente alcanza su valor mínimo para $\|\theta\|_{\mathbb{T}} = 1/2$, a saber el valor $2/\pi$, lo cual nos da la primera desigualdad.

Los conjuntos de Bohr tienen varias propiedades que los hacen semejantes a subgrupos; veremos ejemplos de tales propiedades en la sección siguiente. Por ahora, confirmemos que los conjuntos de Bohr pueden ser de gran densidad y que tienen grandes coeficientes de Fourier.

**Proposición 2.14.** *Sea* $\delta \leqslant 1/(4\pi)$, *sea* $S \subset G$, *y denotemos por* $\beta$ *la densidad* $|B(S,\delta)|\,/\,|G|$. *Entonces para todo* $r \in S$ *tenemos* $|\widehat{1_{B(S,\delta)}}(r)| \geqslant \beta/2$. *Tenemos también*

$$\delta^{|S|} \leqslant \beta \leqslant 4/|S|. \tag{2-12}$$

*Prueba.* Denotando $B(S,\delta)$ por $B$, tenemos

$$|G|\,|\widehat{1_B}(r)| = \Big|\sum_{x\in B} e_r(x)\Big| = \Big||B| - \sum_{x\in B}(1 - e_r(x))\Big| \geqslant |B| - \sum_{x\in B}\big|1 - e_r(x)\big|.$$

Si $r \in S$ entonces $\sum_{x\in B}\big|1 - e_r(x)\big| \leqslant |B|\,2\pi\,\sup_{x\in B}\|r\cdot x\|_{\mathbb{T}} \leqslant |B|\,/2$, de donde sigue nuestra primera afirmación.

Combinando esta estimación con la identidad de Parseval obtenemos

$$\beta = \sum_r |\widehat{1_B}(r)|^2 \geqslant |S|\beta^2/4$$

de donde sigue la segunda desigualdad de (2-12).

Para todo elemento fijo $z = (z_r)_{r\in S} \in \mathbb{T}^S$, tenemos

$$\sum_{x\in G}\prod_{r\in S} 1_{\|r\cdot x - z_r\|<\delta/2}(x,z) = |\{x\in G : \|r\cdot x - z_r\| < \delta/2,\ \forall\, r\in S\}|.$$

Si fijamos cualquier elemento $x_0$ que verifica $\|r\cdot x_0 - z_r\| < \delta/2$ para todo $r\in S$, entonces para todo otro elemento $x$ con esta propiedad tenemos que $y = x - x_0$ verifica $\|r\cdot y\| \leqslant \|r\cdot x - z_r\| + \|r\cdot x_0 - z_r\| < \delta$. Deducimos que

$$\sum_{x\in G}\prod_{r\in S} 1_{\|r\cdot x - z_r\|<\delta/2}(x,z) \leqslant |\{y\in G : \|r\cdot y\| < \delta,\ \forall\, r\in S\}| = |B(S,\delta)|.$$

Integrando esto sobre $z\in\mathbb{T}^S$ deducimos que $\sum_{x\in G}\delta^{|S|} \leqslant |B(S,\delta)|$, de donde sigue la primera desigualdad en (2-12). □

Nótese que en grupos $\mathbb{F}_p^n$ un conjunto de Bohr $B(S,\delta)$ siempre contiene el subespacio anulador $S^{\perp}$ (de hecho es igual a este subespacio si $\delta < 1/p$). Esta propiedad tiene una versión análoga en $\mathbb{Z}_N$, si reemplazamos los subespacios por progresiones aritméticas simétricas con respecto a 0.

**Ejemplo 2.15.** *Un conjunto de Bohr $B(S,\delta)$ en $\mathbb{Z}_N$ contiene una progresión aritmética centrada en 0 y de cardinalidad al menos $\delta N^{1/|S|}$.*

Esto se demuestra usando el teorema de aproximación simultánea de Dirichlet (véase el ejercicio 6). Se dice que la progresión en este ejemplo es "larga" porque su cardinalidad es una función de $N$ que tiende al infinito con $N$, si $|S|$ y $\delta$ están fijados.[8] (Si se quiere ser más preciso, se dice que la progresión tiene *tamaño polinomial* en $N$.) Como veremos, a menudo en $\mathbb{Z}_N$ es técnicamente más conveniente trabajar con tales progresiones aritméticas que con conjuntos de Bohr enteros.

Las progresiones aritméticas son nuestro tercer y último ejemplo de conjuntos con grandes coeficientes de Fourier.

**Lema 2.16.** *Sea $P$ una progresión aritmética en $\mathbb{Z}_N$ ($N$ primo), de cardinalidad $\ell\in(0,N/2)$ y diferencia común $d$. Entonces $|\widehat{1_P}(d^{-1})| \geqslant 2\ell/(\pi N)$. Tenemos también para todo $r\neq 0$ $|\widehat{1_P}(r)| \leqslant \min\{\frac{\ell}{N}, \frac{1}{2N\|rd/N\|_{\mathbb{T}}}\}$.*

*Prueba.* Nótese primero que la magnitud de los coeficientes de Fourier de un conjunto no cambia si se traslada el conjunto, de modo que podemos suponer que $P = \{0, d, 2d, \dots, (\ell-1)d\}$. Usando (2-11), tenemos entonces

$$N\,|\widehat{1_P}(d^{-1})| = \Big|\sum_{t=0}^{\ell-1} e(t/N)\Big| = \frac{|1-e(\ell/N)|}{|1-e(1/N)|} \geqslant 2\ell/\pi$$

[8]Se conoce un resultado más fuerte, a saber que un conjunto de Bohr $B(S,\delta)$ en $\mathbb{Z}_N$ siempre contiene una progresión aritmética generalizada (o multidimensional) de *densidad positiva* dependiente sólo de $|S|, \delta$ (véase Tao y Vu (2006, Proposición 4.22)).

Por otro lado, para todo $r \neq 0$ tenemos

$$N\left|\widehat{1_P}(r)\right| = \left|\sum_{t=0}^{\ell-1} e(trd/N)\right| = \frac{|1-e(\ell rd/N)|}{|1-e(rd/N)|} \leqslant \frac{1}{2\|rd/N\|_{\mathbb{T}}}.$$

□

Hemos visto algunos ejemplos principales de conjuntos $A$ con gran norma $\|f_A\|_u$. Sin embargo, se ve fácilmente que estos ejemplos no son exhaustivos: si $A$ es un tal conjunto y $B$ es un conjunto suficientemente quasi-aleatorio disjunto de $A$, entonces (por la linealidad de la transformada de Fourier) el conjunto $C = A \cup B$ también tiene $\|f_C\|_u$ grande. Por lo tanto, aun no hemos dado una respuesta general satisfactoria a la pregunta 2.10.

El ejemplo del conjunto $C$ que acabamos de ver sugiere una posible respuesta: si $\|f_A\|_u$ es grande, puede que $A$ no sea uno de los conjuntos estructurados que hemos visto, pero a lo mejor es necesario que $A$ tenga una intersección de gran tamaño con un tal conjunto.

Roth confirmó una idea de este tipo en (Roth 1953) trabajando en grupos $\mathbb{Z}_N$, y la usó para demostrar su teorema sobre progresiones de longitud 3 (teorema que veremos en la sección siguiente). Más precisamente, la idea de Roth es que si $\|f_A\|_u$ es grande entonces existe una progresión aritmética larga dentro de la cual el conjunto $A$ tiene *mayor densidad* que en el grupo original. Terminaremos esta sección con este resultado, pero primero presentamos un resultado análogo en los grupos $\mathbb{F}_p^n$, donde el argumento se puede expresar más limpiamente usando subespacios de codimensión 1 en vez de progresiones aritméticas.

**Lema 2.17** (Incremento de densidad en $\mathbb{F}_p^n$). *Sea $A \subset \mathbb{F}_p^n$ de densidad $\alpha$, y sean $r \in \mathbb{F}_p^n \setminus \{0\}$ y $c > 0$ tales que $|\widehat{1_A}(r)| \geqslant c\alpha$. Entonces existe un subespacio $V \leqslant \mathbb{F}_p^n$ de codimensión 1 tal que $|A \cap (x+V)| / |V| \geqslant \alpha(1+c/2)$ para algún $x \in \mathbb{F}_p^n$.*

*Prueba.* Como es de esperar, el subespacio que tomaremos es el conjunto anulador $V = r^{\perp}$. Las clases laterales de este subgrupo son las preimagenes de los elementos de $\{0, \frac{1}{p}, \dots, \frac{p-1}{p}\}$ por la función $\mathbb{F}_p^n \to \mathbb{T}, x \mapsto r \cdot x$. Fijando un elemento $x_j$ en cada una de estas preimagenes, tenemos $\mathbb{F}_p^n = \bigsqcup_{j=0}^{p-1} x_j + V$. Denotemos por $\alpha_j$ la densidad que tiene $A$ dentro de $x_j + V$, es decir $\alpha_j = |A \cap (x_j + V)|/|V|$. Podemos multiplicar $\widehat{1_A}(r)$ por algún número $e(\theta)$ de modo que tengamos $\widehat{1_A}(r)e(\theta) > 0$. Entonces, usando el hecho que $\mathbb{E}_{x\in\mathbb{F}_p^n} f_A(x) = 0$, tenemos

$$c\alpha \leqslant \widehat{1_A}(r)e(\theta) = \widehat{f_A}(r)e(\theta) = \mathbb{E}_{x\in\mathbb{F}_p^n}\, f_A(x)\,(e(r\cdot x+\theta)+1).$$

Tomando la parte real, tenemos $c\alpha \leqslant \mathbb{E}_{x\in\mathbb{F}_p^n}\, f_A(x)\,\mathrm{Re}(e(r\cdot x+\theta)+1)$. La función $x \mapsto \mathrm{Re}(e(r\cdot x+\theta)+1)$ toma un valor constante en el intervalo $[0,2]$ sobre cada clase $x_j + V$, valor que denotaremos por $\lambda_j$. Tenemos entonces

$$\begin{aligned} c\alpha \leqslant \frac{1}{p^n}\sum_{j=0}^{p-1}\sum_{x\in x_j+V}(1_A(x)-\alpha)\,\lambda_j &= \frac{|V|}{p^n}\sum_{j=0}^{p-1}\lambda_j\,\mathbb{E}_{x\in x_j+V}(1_A(x)-\alpha) \\ &= \mathbb{E}_{j\in[0,p-1]}\,\lambda_j(\alpha_j-\alpha). \end{aligned}$$

Tiene por lo tanto que existir $j \in [0, p-1]$ tal que $c\alpha \leqslant \lambda_j(\alpha_j - \alpha)$, de lo cual se deduce que $\alpha_j \geqslant \alpha + c\alpha/2$. □

En este argumento hemos obtenido el incremento de densidad sobre una traslación del subespacio anulador de un carácter dominante de $A$ (un carácter $e_r$ para el cual $|\widehat{f_A}(r)| = \|f_A\|_u$).

Pasamos ahora a los grupos cíclicos de orden primo. Aquí, como ya mencionamos, no hay subgrupos anuladores notriviales, pero podemos aun obtener un incremento de densidad sobre una traslación de una progresión aritmética larga que es *anuladora aproximativa* de un carácter dominante. Este aspecto aproximativo hace que este argumento sea un poco más técnico que el de $\mathbb{F}_p^n$, por dos razones principales.

Primero, en $\mathbb{Z}_N$ la partición del grupo en traslaciones de una misma progresión es solo aproximativa (es decir que no se da exactamente el análogo de la partición de $\mathbb{F}_p^n$ en clases laterales de $V$). Segundo, para las aplicaciones, a menudo querremos garantizar que, si empezamos con $A \subset \mathbb{Z}_N$ visto como un conjunto de enteros en $[N]$, entonces el incremento de densidad tiene lugar en una progresión aritmética en $\mathbb{Z}_N$ que es también una progresión en $\mathbb{Z}$ (por ejemplo, el conjunto $\{4, 5, 1\}$ es una progresión en $\mathbb{Z}_5$ pero no lo es en $\mathbb{Z}$).

Estas dificultades se pueden resolver usando el resultado siguiente.

**Lema 2.18.** *Sea $I$ un intervalo en $\mathbb{Z}_N$ de cardinalidad $m$, y sea $r$ un elemento nonulo de $\mathbb{Z}_N$. Para todo $\epsilon > 0$, existe una partición de $I$ en progresiones aritméticas $P_i$ en $\mathbb{Z}_N$ de cardinalidad al menos $\epsilon\sqrt{m}/2$, y tales que para todo $x, y \in P_i$ tenemos $\|(rx - ry)/N\|_{\mathbb{T}} \leqslant \epsilon$, para todo $i$.*

Aquí $rx/N$ es nuestra forma bilineal favorita $r \cdot x$ sobre $\mathbb{Z}_N$.

*Prueba.* El resultado es invariante por traslación en $\mathbb{Z}_N$, con lo cual podemos suponer que $I$ es un intervalo de enteros en $[N]$. Sea $u = \sqrt{m}$, y consideremos los números $0, r, 2r, \ldots, ur$. Por el principio del palomar, existen $v < w$ en $[0, u]$ tales que $\|r \cdot w - r \cdot v\|_{\mathbb{T}} = \|r \cdot (w - v)\|_{\mathbb{T}} \leqslant 1/u$. Denotando $w - v$ por $s$, dividimos $I$ en clases módulo $s$, y observamos que cada clase tiene cardinalidad al menos $\lfloor m/s \rfloor \geqslant \lfloor m/u \rfloor$. Cada una de estas clases puede ser dividida en progresiones aritméticas de la forma $a, a + s, \ldots, a + \ell s$, con $\epsilon u/2 < \ell \leqslant \epsilon u$. Para cualquier tal progresión $P_i$, tenemos $\|r \cdot x - r \cdot y\|_{\mathbb{T}} \leqslant \ell \|r \cdot s\|_{\mathbb{T}} \leqslant \epsilon$ para todo $x, y \in P_i$. □

Dado un conjunto $X$ y una función $\phi : X \to \mathbb{C}$, llamaremos la cantidad $\sup_{x,y \in X} |\phi(x) - \phi(y)|$ el *diámetro* de $\phi$ sobre $X$. Del último lema deducimos, usando (2-11), que el carácter $e_r$ tiene diámetro como mucho $2\pi\epsilon$ sobre cada progresión $P_i$. En otras palabras, cada progresión $P_i$ es un *conjunto de nivel aproximado* para $e_r$ (un conjunto sobre el cual $e_r$ es casi constante).

Podemos ahora demostrar el incremento de densidad en grupos $\mathbb{Z}_N$, con la garantía adicional de que si identificamos el conjunto subyacente con $[N]$, entonces la progresión donde ocurre el incremento es una progresión en $\mathbb{Z}$ (y no solo en $\mathbb{Z}_N$).

**Lema 2.19** (Incremento de densidad en $[N]$). *Sea $B \subset [N]$ de densidad $|B|/N = \beta$ y sean $r \neq 0$ en $\mathbb{Z}_N$ y $c > 0$ tales que en $\mathbb{Z}_N$ tengamos $|\widehat{1_B}(r)| \geqslant c\beta$. Entonces existe una progresión aritmética $P \subset [N]$ de cardinalidad al menos $c_1\sqrt{N}$ tal que*

$$|B \cap P| / |P| \geqslant \beta(1 + c_2).$$

*Podemos tomar $c_1 = c/60$, $c_2 = c/4$.*

*Prueba.* Fijando $\epsilon = c/(8\pi)$, usando el lema 2.18 dividimos $[N]$ en progresiones aritméticas de longitud al menos $\epsilon\sqrt{N}/2$, sobre cada una de las cuales el carácter $e_r$ tiene diámetro como mucho $2\pi\epsilon$. Fijando un elemento $x_i$ arbitrario en cada progresión $P_i$, deducimos lo siguiente:

$$\begin{aligned} |\widehat{f_B}(r)|\,N &\leqslant \sum_i \Big| \sum_{x\in P_i} f_B(x) e(r\cdot x)\Big| \\ &= \sum_i \Big| \sum_{x\in P_i} f_B(x) e_r(x_i) + \sum_{x\in P_i} f_B(x)\big(e_r(x) - e_r(x_i)\big)\Big| \\ &\leqslant \sum_i \left[ \Big| \sum_{x\in P_i} f_B(x)\Big| + \sum_{x\in P_i} |f_B(x)| 2\pi\epsilon \right] \\ &\leqslant 4\beta\pi\epsilon N + \sum_i \Big| \sum_{x\in P_i} f_B(x)\Big|. \end{aligned}$$

Usando el valor elegido para $\epsilon$, el hecho que $\mathbb{E}_{\mathbb{Z}_N} f_B = 0$, y la suposición $|\widehat{f_B}(r)| \geqslant c\beta$, deducimos que $\sum_i \big[\big|\sum_{x\in P_i} f_B(x)\big| + \sum_{x\in P_i} f_B(x)\big] \geqslant c\beta N/2 = \sum_i c\beta|P_i|/2$.

Existe por lo tanto un valor de $i$ tal que

$$\Big| \sum_{x\in P_i} f_B(x)\Big| + \sum_{x\in P_i} f_B(x) \geqslant c\beta|P_i|/2,$$

de donde sigue que $\sum_{x\in P_i} f_B(x) \geqslant c\beta|P_i|/4$, lo cual implica el resultado. □

Existen también argumentos de incremento de densidad con respecto a conjuntos de Bohr. Estos argumentos pueden ser mucho más eficientes que los lemas 2.17 y 2.19 de un punto de vista cuantitativo (dando mayores incrementos de densidad); son también bastante más técnicos (ver por ejemplo (Sanders 2011)).

# 3 Aplicaciones

En esta sección combinamos las ideas vistas en la sección 2 para demostrar algunos resultados clásicos de combinatoria aritmética.

## 3.1 Teoremas de Roth y Meshulam.

¿Cuán grande puede ser un subconjunto de un grupo abeliano finito si el conjunto no contiene progresiones aritméticas de longitud 3? Los teoremas de Roth y de Meshulam dan respuestas notriviales a esta pregunta, en los grupos $\mathbb{Z}_N$ y $\mathbb{F}_p^n$ respectivamente.

Históricamente el teorema de (Roth 1953) es anterior al de (Meshulam 1995), pero la demostración de este último es más sencilla (gracias a la ya mencionada riqueza algebraica de $\mathbb{F}_p^n$), de modo que empezaremos con este resultado. Para esto, no se pierde generalidad con trabajar en $\mathbb{F}_3^n$ en vez de $\mathbb{F}_p^n$.

**Teorema 3.1** (Meshulam). *Sea $A \subset G = \mathbb{F}_3^n$, sea $\alpha = |A|/|G|$, y supongamos que $A$ no contiene ninguna progresión aritmética notrivial de longitud 3. Entonces*[9]

$$\alpha = O(1/\log|G|). \tag{3-1}$$

*Prueba.* De la suposición inicial deducimos que $T_3(A) = \alpha/|G|$ (contando las progresiones $(x, x+d, x+2d)$ con $d = 0$). La Proposición 2.6 implica entonces que $\|f_A\|_u \geqslant \alpha^{-1}(\alpha^3 - \alpha/|G|)$. Si $|G| \geqslant 2\alpha^{-2}$, obtenemos que $\|f_A\|_u \geqslant \alpha^2/2$.

Podemos entonces obtener un incremento de densidad: el Lema 2.17 nos da un subespacio $V_1$ y un elemento $x_1$ tal que $|A \cap (x_1 + V_1)|/|V_1| \geqslant \alpha(1 + \alpha/4)$. Denotemos el conjunto $A \cap (x_1 + V_1)$ por $A_1$, y su densidad en $V_1$ por $\alpha_1$. Observemos que este $A_1$, siendo un subconjunto de $A$, tampoco contiene progresiones. Como una traslación de una progresión sigue siendo una progresión, podemos trasladar $A$ por $-x_1$ para suponer que $A_1$ es un subconjunto de $V_1 \cong \mathbb{F}_3^{n-1}$. La idea ahora es repetir este argumento en $A_1$. Obtenemos así otro subconjunto $A_2$, dentro de un nuevo subespacio $V_2 \cong \mathbb{F}_3^{n-2}$, con aun mayor densidad $\alpha_2 \geqslant \alpha_1(1 + \alpha_1/4) \geqslant \alpha + 2\alpha^2/4$.

Cada vez que repetimos el argumento, perdemos una dimensión. Por otro lado, incrementamos la densidad de $\alpha^2/4$. En particular, pasamos de densidad $\alpha$ a $2\alpha$ en $\alpha/(\alpha^2/4) = 4\alpha^{-1}$ repeticiones, pasamos de $2\alpha$ a $4\alpha$ en $2\alpha/((2\alpha)^2/4) = \frac{1}{2}4\alpha^{-1}$ repeticiones, y en general de $2^j\alpha$ a $2^{j+1}\alpha$ en $\frac{1}{2^j}4\alpha^{-1}$ repeticiones. Como la densidad de un conjunto es como mucho 1, el número de veces que podemos repetir este argumento tiene un máximo $t \leqslant (1 + \frac{1}{2} + \frac{1}{4} + \ldots)4\alpha^{-1} = 8\alpha^{-1}$. En particular, después de $t$ repeticiones es necesario que $|V_t|$ sea menor que $2\alpha^{-2}$ (de lo contrario sería posible repetir el argumento). Por lo tanto tenemos $3^{n-8\alpha^{-1}} \leqslant 3^{n-t} = |V_t| \leqslant 2\alpha^{-2}$, de donde sigue la cota (3-1). □

Encontrar la cota superior óptima para $\alpha$ en el teorema 3.1 es un problema central en esta área. Por un lado, se intenta reducir la cota superior, y en esta dirección el mejor resultado actualmente conocido que usa el análisis de Fourier es el de (Bateman y Katz 2012), quienes obtuvieron la cota $\alpha = O\big(1/(\log|G|)^{1+\epsilon}\big)$, para algún $\epsilon > 0$. Por otro lado, se intenta encontrar ejemplos de grandes conjuntos en $\mathbb{F}_3^n$ sin progresiones. El conjunto $\{0, 1\}^n$ nos da un primer ejemplo, de cardinalidad $2^n$. El ejemplo de mayor cardinalidad actualmente conocido es el que dio (Edel 2004), demostrando lo siguiente.

**Teorema 3.2** (Edel). *Existe un conjunto $A \subset \mathbb{F}_3^n$ sin progresiones de longitud 3 y de cardinalidad al menos $2.217^n$.*

La diferencia entre la cota superior de Bateman y Katz y la inferior de Edel es enorme. Hasta hace poco tiempo, no se sabía si era posible reducir esta diferencia demostrando que existe $\epsilon > 0$ tal que si $A \subset \mathbb{F}_3^n$ no contiene progresiones de longitud 3 entonces $|A| \leqslant (3 - \epsilon)^n$. Recientemente tuvo lugar un avance muy importante relativo a este problema, cuando Ellenberg y Gijswijt establecieron una cota superior de tipo $|A| = o(2.756^n)$; véase (Ellenberg y Gijswijt 2017). Ellenberg y Gijswijt no usaron el análisis de Fourier, sino que extendieron una innovación crucial en el método polinomial que fue introducida por (Croot, Lev y Pach 2016). Parece natural preguntarse si una

[9]Recuerden la notación de Landau: dos funciones $f(N), g(N)$ satisfacen $f = O(g)$ si existe una constante absoluta $C > 0$ y un $N_0$ tal que $|f(N)| \leqslant C|g(N)|$ para todo $N \geqslant N_0$; las funciones satisfacen $f = o(g)$ si $f(N)/g(N) \to 0$ cuando $N \to \infty$.

cota superior para $|A|$ de este tipo se puede demostrar usando análisis de Fourier.

Pasamos ahora al teorema de Roth.

**Teorema 3.3** (Roth)**.** *Sea $A \subset [N]$, sea $\alpha = |A| / N$, y supongamos que $A$ no contiene ninguna progresión aritmética notrivial de longitud 3. Entonces*

(3-2) $$\alpha = O(1/\log\log N).$$

La prueba usa ideas similares a las de la prueba anterior. En particular, querremos pasar a una larga progresión aritmética donde la densidad de $A$ tenga un incremento de al menos una fracción fija de $\alpha^2$. Técnicamente, obtener este incremento de densidad será un poco más complicado que en la prueba anterior. Uno de los problemas aquí es que el intervalo de enteros $[N]$ no tiene estructura de grupo con la cual hacer análisis de Fourier, de modo que primero tendremos que trasladar el problema a un grupo cíclico. Además, tendremos que hacer esto sin perder de vista que estamos buscando progresiones en los enteros y no solo en grupos cíclicos.

*Prueba.* Denotemos por $\alpha$ la densidad de $A$ en $[N]$, y supongamos que $N > 50/\alpha^2$. Sea $p$ un número primo entre $N/3$ y $2N/3$ (cuya existencia está garantizada por el postulado de Bertrand).

Nótese primero que, denotando por $B$ el conjunto $A \cap [p]$ y por $\beta$ la densidad de $B$ en $[p]$, podemos suponer que $\beta \geqslant \alpha(1-\alpha/200)$, pues de lo contrario tendríamos $|A \cap [p+1, N]| \geqslant \alpha(N-p)+\alpha^2 p/200 \geqslant \alpha(1+\alpha/400)(N-p)$, lo cual implicaría que ya tendríamos un incremento de densidad apropiado en la progresión aritmética $[p+1, N]$, a saber un incremento de al menos $\alpha^2/400$. Trabajemos por lo tanto con este conjunto $B$, usando la transformada de Fourier en $\mathbb{Z}_p$.

De nuestra suposición combinatoria sobre $A$ deduciremos primero una cota inferior notrivial para $\|f_B\|_u$. Sea $B' = B \cap [p/3, 2p/3]$, con densidad en $[p]$ denotada por $\beta'$, y nótese que una progresión aritmética en $B \times B' \times B'$ modulo $p$ es siempre una progresión en los enteros (y contenida en $A$, por supuesto). Por lo tanto, tenemos

$$\Big|\beta\beta'^2 - \mathbb{E}_{\substack{x_1,x_2,x_3\in\mathbb{Z}_p\\ x_1-2x_2+x_3=0}} 1_B(x_1)1_{B'}(x_2)1_{B'}(x_3)\Big| = |\beta\beta'^2 - \beta'/p\,|.$$

La proposición 2.7 nos dice que esto es como mucho $\|f_B\|_u\,\beta'$. Por otro lado, podemos suponer que $\beta' \geqslant \beta/5$. En efecto, de lo contrario tendríamos $|B \cap [p/3]| \geqslant 2\beta p/5$ o bien $|B \cap [2p/3, p]| \geqslant 2\beta p/5$, y entonces $B$ tendría una densidad al menos $(2\beta p/5)/(p/3) = 6\beta/5$ en una de las progresiones $[p/3]$, $[2p/3, p]$. Concluimos que $\|f_B\|_u \geqslant \beta^2/5 - 1/p \geqslant \beta^2/10$. Podemos ahora deducir el incremento de densidad deseado. El lema 2.19 nos dice que existe una progresión $P \subset [p]$ de longitud al menos $(\beta/600)\sqrt{p} \geqslant \alpha\sqrt{N}/2^{10}$ en la cual la densidad de $B$ (y por lo tanto de $A$) se incrementa de $\beta^2/40 \geqslant \alpha^2/80$. Por lo tanto, en todos los casos podemos garantizar un incremento de densidad de al menos $\alpha^2/400$.

A partir de aquí, la prueba consiste en una iteración muy similar a la que ya vimos para el teorema de Meshulam. Repetimos el argumento como mucho $400/\alpha$ veces para que la densidad suba de $\alpha$ a $2\alpha$, y siguiendo así se llega a una densidad 1 después de $t \leqslant 400\alpha^{-1}(1+\frac{1}{2}+\frac{1}{4}+\cdots) \leqslant 800\alpha^{-1}$ repeticiones. En cada repetición pasamos de una progresión de longitud $m$ a una de longitud $(\alpha/2^{10})\sqrt{m}$, con lo cual la longitud de la progresión después de $t$ repeticiones es $(\alpha/2^{10})^{1+1/2+1/4+\cdots} N^{1/2^t} \geqslant (\alpha^2/2^{20}) N^{1/2^{800\alpha^{-1}}}$. Para que no haya contradicción (es decir que no se pueda repetir más el incremento) es necesario que $(\alpha^2/2^{20}) N^{1/2^{800\alpha^{-1}}} \leqslant 50\alpha^{-2}$, de donde se deduce la cota (3-2). □

El problema de mejorar (reducir) la cota (3-2) es tan famoso como en el caso del teorema de Meshulam. De hecho lo es quizás aun más, por sus consecuencias directas en teoría de números. Por ejemplo, una cota de tipo $\alpha \leqslant 1 / (\log N)^{1+c}$, para alguna constante absoluta $c > 0$, implicaría el *teorema de Roth en los números primos*, a saber que todo subconjunto de los números primos de densidad positiva contiene una progresión de longitud 3. Este resultado ya fue demostrado por (Green 2005), pero la idea aquí es que con la cota mencionada el resultado se deduciría mucho más fácilmente, usando sólo el teorema de los números primos. La mejor cota actualmente conocida para el teorema de Roth se debe principalmente a Sanders, quien demostró en (Sanders 2011), combinando varias herramientas recientes en un argumento bastante técnico de incremento de densidad sobre conjuntos de Bohr, que la cota superior en (3-2) se puede reducir a $(\log\log N)^6/\log N$; véase también la mejora de (Bloom 2016) que reduce la potencia 6 a 4.

En cuanto a ejemplos de grandes subconjuntos de $[N]$ sin progresiones, la construcción principal conocida se debe a (Behrend 1946). ¡Esta construcción no se ha mejorado significativamente desde 1946! (Para una exposición clara de una cota algo más fuerte, siguiendo la idea original de Behrend, véase (Green y Wolf 2010).)

**Teorema 3.4** (Behrend). *Existe un conjunto $A \subset [N]$ sin progresiones aritméticas de longitud 3 que verifica*

$$\alpha \geqslant c\ \exp\big(-c\sqrt{\log N}\big),$$

*donde $c > 0$ es una constante absoluta.*

## 3.2 Lema de Bogolyubov.

Recordemos que, dados dos subconjuntos $A$, $B$ de un grupo abeliano, su *conjunto suma* es $A+B = \{a+b : a \in A, b \in B\}$. Como ya mencionamos, este conjunto es precisamente el soporte de la convolución $1_A * 1_B$.

El resultado del que trataremos aquí es un ejemplo central de cómo la transformada de Fourier permite estudiar la estructura de conjuntos suma, gracias a su interacción con la operación de convolución. Hay un principio general en análisis relativo a la operación de convolución, según el cual esta operación "suaviza" las funciones[10]. En particular cuanto más se repite esta operación (tomando $f * f$, después $f * f * f$, etc.) más regular se hace la función. Este principio tiene una versión análoga de gran utilidad en combinatoria aritmética, a saber que tomar sumas (o diferencias) iteradas de un conjunto produce conjuntos con estructura aritmética cada vez más rica. El resultado principal de esta sección ilustra esto con el conjunto de suma y diferencia siguiente:

$$2A - 2A = \{a_1 + a_2 - a_3 - a_4 : a_1, a_2, a_3, a_4 \in A\}.$$

[10] Una instancia precisa de este principio es que para funciones integrables $f, g : \mathbb{R} \to \mathbb{R}$, si $f$ y $g$ se pueden derivar $m$ y $n$ veces respectivamente, entonces la convolución $f * g$ se puede derivar $m + n$ veces.

Usando la transformada de Fourier, se puede demostrar que un tal conjunto siempre contiene un conjunto de Bohr de gran densidad (densidad que depende sólo de $\alpha = |A|/|G|$). Este resultado, llamado comúnmente *lema de Bogolyubov*, se origina en (Bogolioùboff 1939) y juega un papel importante en el área (como indicamos más adelante).

**Lema 3.5** (Bogolyubov)**.** *Sea $A$ un subconjunto de densidad $\alpha$ en un grupo abeliano finito $G$. Entonces existe $S \subset G$ con $|S| \leqslant 2\alpha^{-2}$ y tal que $2A - 2A \supset B(S, \frac{1}{4})$.*

Antes de dar la prueba, motivaremos este resultado describiendo unas propiedades básicas de los conjuntos de Bohr que los hacen semejantes a subgrupos; de modo que, efectivamente, el lema nos dice que $2A - 2A$ tiene una rica estructura aditiva. Nótese también que si tomamos menos sumandos, considerando por ejemplo $A - A$ en vez de $2A - 2A$, entonces el lema falla (véase el ejercicio 7).

**Definición 3.6** (Grupo aproximado)**.** Sea $K$ un entero positivo. Un subconjunto finito $S$ de un grupo abeliano es un *grupo $K$-aproximado* si contiene el elemento 0, es simétrico respecto a 0 (es decir que $S = -S$) y existe un cubrimiento de $S + S$ por $K$ traslaciones de $S$.

Esta noción aparece en Tao y Vu (2006, Definición 2.25).

**Lema 3.7.** *Un conjunto de Bohr $B(S, \delta)$ es un grupo $11^{|S|}$-aproximado.*

Para demostrar esto usamos un resultado muy útil llamado *lema de cubrimiento*, introducido por Ruzsa.

**Lema 3.8** (Ruzsa)**.** *Sean $A, B$ subconjuntos finitos de un grupo abeliano tales que $|A+B| \leqslant K|A|$. Entonces existe un subconjunto $X \subset B$ de cardinalidad $|X| \leqslant K$ tal que $B \subset A - A + X$.*

*Prueba.* Sea $\mathfrak{X}$ la familia de conjuntos $Y \subset B$ tales que los conjuntos $y + A$, $y \in Y$ son disjuntos dos a dos. Sea $X$ un miembro de $\mathfrak{X}$ de cardinalidad máxima. Entonces, por un lado, tenemos $|X|\,|A| = |X + A| \leqslant |B + A| \leqslant K\,|A|$, de donde sigue que $|X| \leqslant K$, y por otro lado, siendo $X$ máximo, para todo $b \in B$ existe $x \in X$ tal que $(b + A) \cap (x + A) \neq \varnothing$, luego $b \in x + A - A$. □

*Prueba del lema 3.7.* Arguyendo como en la prueba de la proposición 2.14 tenemos, para todo $z \in \mathbb{T}^S$,

$$\sum_{x \in G} \prod_{r \in S} 1_{\|r \cdot x - z_r\| < \frac{\delta}{2}}(x, z) \leqslant |B(S, \delta)|.$$

Sea $\lambda > 0$. Definimos el conjunto $Z = \{z \in \mathbb{T}^S : \sup_{r \in S} \|z_r\| \leqslant (\lambda + \frac{1}{2})\,\delta\}$. Tenemos entonces

$$\begin{aligned}
\delta^{|S|}\,|B(S, \lambda\delta)| &= \sum_{x \in B(S,\lambda\delta)} \int_{\mathbb{T}^S} \prod_{r \in S} 1_{\|r \cdot x - z_r\| < \frac{\delta}{2}}(x, z)\, \mathrm{d}\mu_{\mathbb{T}^S}(z) \\
&= \int_{\mathbb{T}^S} \sum_{x \in B(S,\lambda\delta)} \prod_{r \in S} 1_{\|r \cdot x - z_r\| < \frac{\delta}{2}}(x, z)\;\; 1_Z(z)\, \mathrm{d}\mu_{\mathbb{T}^S}(z) \\
&\leqslant \int_Z \sum_{x \in G} \prod_{r \in S} 1_{\|r \cdot x - z_r\| < \frac{\delta}{2}}(x, z)\, \mathrm{d}\mu_{\mathbb{T}^S}(z) \;\leqslant\; |B(S, \delta)|\; \mu_{\mathbb{T}^S}(Z).
\end{aligned}$$

Como $\mu_{\mathbb{T}^S}(Z) = ((2\lambda + 1)\delta)^{|S|}$, deducimos que

$$\forall\, \lambda > 0, \qquad |B(S, \lambda\delta)| \leqslant (2\lambda + 1)^{|S|}|B(S, \delta)|. \tag{3-3}$$

Denotemos $B(S, \lambda\delta)$ por $B_\lambda$. Aplicando (3-3) tenemos $|B_{1/2} + B_2| \leqslant |B_{5/2}| \leqslant 11^{|S|} B_{1/2}$. Por lo tanto el lema de cubrimiento de Ruzsa nos da un conjunto $X \subset B_2$ de cardinalidad como mucho $|B_{5/2}|/|B_{1/2}| \leqslant 11^{|S|}$, tal que $B_2 \subset X + B_{1/2} - B_{1/2} \subset X + B_1$. Por lo tanto, el conjunto $B(S, \delta) + B(S, \delta) \subset B_2$ se puede cubrir con $11^{|S|}$ traslaciones de $B(S, \delta)$, a saber las traslaciones $x + B_1, x \in X$. □

Para demostrar el lema 3.5, introducimos la noción Fourier-analítica siguiente.

**Definición 3.9.** Sea $\epsilon > 0$ y sea $f : G \to \mathbb{C}$ con $\|f\|_{L^\infty} \leqslant 1$. El *$\epsilon$-espectro* de $f$ es el conjunto de frecuencias

$$\mathrm{Spec}_\epsilon(f) = \{r \in \widehat{G} : |\widehat{f}(r)| \geqslant \epsilon\}.$$

La identidad de Parseval nos da una primera cota superior general para la cardinalidad de este conjunto. En efecto, tenemos

$$\|f\|_{L^2}^2 = \|\widehat{f}\|_{\ell^2}^2 \geqslant \sum_{r \in \mathrm{Spec}_\epsilon(f)} |\widehat{f}(r)|^2 \geqslant \epsilon^2\, |\mathrm{Spec}_\epsilon(f)|,$$

luego $|\mathrm{Spec}_\epsilon(f)| \leqslant \epsilon^{-2}\|f\|_{L^2}^2 \leqslant \epsilon^{-2}\|f\|_{L^\infty}^2$.

Podemos ahora demostrar el lema de Bogolyubov.

*Prueba del lema 3.5.*
La primera observación clave es que $2A - 2A = \mathrm{supp}(1_A * 1_A * 1_{-A} * 1_{-A})$. Por otra parte, tenemos para todo $x \in G$

$$1_A * 1_A * 1_{-A} * 1_{-A}(x) = \sum_{r \in \widehat{G}} |\widehat{1}_A(r)|^4\, e_r(x). \tag{3-4}$$

Fijemos $S = \mathrm{Spec}_\lambda(1_A)$, donde $\lambda > 0$ es un parámetro que escogeremos más tarde. La idea es que para probar que $B(S, \frac{1}{4}) \subset 2A - 2A$, basta con demostrar que para todo $x \in B(S, \frac{1}{4})$ la convolución en (3-4) es positiva en $x$. Esta convolución es igual a

$$\mathrm{Re} \sum_r |\widehat{1}_A(r)|^4 e_r(x) = \sum_{r \in S} |\widehat{1_A}(r)|^4 \cos(2\pi\, r \cdot x) + \sum_{r \notin S} |\widehat{1_A}(r)|^4 \cos(2\pi\, r \cdot x).$$

Como $\cos(2\pi\, r \cdot x) > 0$ para todo $x \in B(S, \frac{1}{4})$ y todo $r \in S$, deducimos que

$$1_A * 1_A * 1_{-A} * 1_{-A}(x) \ \geqslant\ |\widehat{1}_A(0)|^4 - \sum_{r \notin S} |\widehat{1}_A(r)|^4 \ = \ \alpha^4 - \sum_{r \notin S} |\widehat{1}_A(r)|^4.$$

La identidad de Parseval nos da que $\sum_{r \notin S} |\widehat{1_A}(r)|^4 \leqslant \lambda^2 \alpha$. Fijando $\lambda = \alpha^{3/2}/\sqrt{2}$, deducimos que la convolución es positiva como deseado. □

Mejorar la cota $2\alpha^{-2}$ en el lema de Bogolyubov es un problema de gran importancia actualmente en combinatoria aritmética, en particular por su relación con la *conjetura polinomial de Freiman–Ruzsa* (para más información sobre esta conjetura, véase los apuntes de Julia (Wolf 2020)). La mejor cota actualmente conocida fue obtenida por (Sanders 2012), y esencialmente reduce $2\alpha^{-2}$ a $O(\log^4(1/\alpha))$. Para una prueba alternativa (y bastante clara) de esta mejora, véase Croot, Łaba y Sisask (2013, Teorema 7.5).

## Referencias

## Ejercicios

**Ejercicio 1.** Demostrar el lema 1.3, a saber que si $G$ es un grupo abeliano finito de exponente $n$, entonces

1. Todo carácter $\chi \in \widehat{G}$ toma sus valores en $\{z \in \mathbb{C} : z^n = 1\}$.
2. Si $G$ es un grupo cíclico $\mathbb{Z}_N$, entonces $\widehat{G} = \{x \mapsto \exp(2\pi i\, rx/N) : r \in \mathbb{Z}_N\}$.

(pista: para la parte (2), usar (1) y el hecho que cada $\chi \in \widehat{\mathbb{Z}_N}$ está determinado por su valor en $x = 1$.)

**Ejercicio 2.** Demostrar que para todo grupo abeliano finito $G$ y $H$, tenemos $G \cong H \Rightarrow \widehat{G} \cong \widehat{H}$ y $\widehat{G \oplus H} \cong \widehat{G} \oplus \widehat{H}$.

**Ejercicio 3.** Demostrar la fórmula de producto para la convolución:

$$\text{para cualquier } r \in G \text{ y } f, g \in \mathbb{C}^G, \quad \widehat{f * g}(r) = \widehat{f}(r)\, \widehat{g}(r).$$

**Ejercicio 4.** Sea $Z$ un grupo abeliano finito, sea $F : Z \to \mathbb{C}$, sea $H$ un subgrupo de $Z$, y sea $g(x) = F(x) 1_H(x)$. Demostrar que $\widehat{1_H} = \frac{|H|}{|G|} 1_{H^\perp}$. Usando esto, demostrar que $\widehat{g}(s) =$ $\frac{|H|}{|G|}\sum_{r\in H^\perp}\widehat{F}(s+r)$ para todo $s\in\widehat{G}$. Deducir la *fórmula de Poisson*:

$$\mathbb{E}_{x\in H}F(x)=\sum_{r\in H^\perp}\widehat{F}(r). \tag{1-5}$$

Usando (1-5), demostrar la fórmula (2-4):

$$\mathbb{E}_{\substack{x_1,\dots,x_t\in G:\\ c_1x_1+\cdots+c_tx_t=0}} f_1(x_1)\,f_2(x_2)\cdots f_t(x_t) = \sum_{r\in G}\widehat{f_1}(c_1r)\,\widehat{f_2}(c_2r)\cdots\widehat{f_t}(c_tr).$$

(pista: aplicar (1-5) con $Z=G^t$ y $H=\{x\in G^t : c_1x_1+\cdots+c_tx_t=0\}$.)

**Ejercicio 5.** Demostrar que la norma $U^2$ de Gowers,

$$\|f\|_{U^2(G)}=\big(\mathbb{E}_{x,h_1,h_2\in G}\ f(x)\ \overline{f(x+h_1)}\ \overline{f(x+h_2)}\ f(x+h_1+h_2)\big)^{1/4},$$

satisface la fórmula $\|f\|_{U^2(G)}=\|\widehat{f}\|_{\ell^4(\widehat{G})}$. Demostrar que $\|f\|_u\leqslant\|f\|_{U^2(G)}\leqslant\|f\|_u^{1/2}\,\|f\|_{L^2(G)}^{1/2}$.

**Ejercicio 6.** El teorema de aproximación de Dirichlet nos dice que dados unos elementos $\theta_1,\dots,\theta_d\in\mathbb{T}$ cualesquiera y un entero $N$, existe $n\in[N-1]$ tal que $\|n\theta_i\|_{\mathbb{T}}\leqslant N^{-1/d}$ para todo $i\in[d]$.

Demostrar que un conjunto de Bohr $B(S,\delta)$ en $\mathbb{Z}_N$ contiene una progresión aritmética centrada en 0 y de cardinalidad al menos $\delta N^{1/|S|}$.

Sea $A$ un subconjunto de $\mathbb{Z}_N$ de cardinalidad $|A|\leqslant\frac{1}{10}\log N$. Demostrar que existe $r\neq 0$ tal que $|\widehat{1_A}(r)|\geqslant\alpha/2$. (Este ejemplo indica que el análisis de Fourier pierde su utilidad combinatoria cuando la densidad del conjunto es demasiado pequeña.)

**Ejercicio 7.** Demostrar que existe un conjunto $A\subset\mathbb{F}_2^n$ de densidad 1/4 tal que $A+A$ no contiene ninguna clase lateral de subespacio de codimensión $\sqrt{n}$.
(Pista: considérese $A=\{x\in\mathbb{F}_2^n : x$ tiene al menos $n/2+\sqrt{n}/2$ coordenadas igual a $1\}$.)

**Ejercicio 8** (Construcción explícita de un conjunto extremadamente quasi-aleatorio).
1) Sea $f_s:\mathbb{Z}_p\to\mathbb{C}$, $x\mapsto e(sx^2/p)$, donde $s\in\mathbb{Z}_p\setminus\{0\}$. Usando el hecho que $|\widehat{f_s}(r)|^2=\mathbb{E}_h\mathbb{E}_x f_s(x+h)\overline{f_s(x)}e(-rh/p)$, demostrar que $\|f_s\|_u\leqslant p^{-1/2}$.
2) Usando el lema 2.16, demostrar que si $P$ es un intervalo de densidad $1/2$ en $\mathbb{Z}_p$ entonces $\sum_{r\in\mathbb{Z}_p}|\widehat{1_P}(r)|\leqslant 10\log p$.
3) Sea $A=\{x\in\mathbb{Z}_p : x^2\in P\}$. Demostrar que $\widehat{1_A}(r)=\sum_s\widehat{1_P}(s)\widehat{f_s}(r)$, y combinar esto con 1) y 2) para deducir que $\sup_{r\neq 0}|\widehat{1_A}(r)|\leqslant 10p^{-1/2}\log p$.
4) La norma $U^3$ de Gowers se define sobre $\mathbb{C}^G$ por

$$\|f\|_{U^3(G)}=\Big(\mathbb{E}_{x,h_1,h_2,h_3\in G}\ f(x)\,\overline{f(x+h_1)}\,\overline{f(x+h_2)}\,f(x+h_1+h_2)\,\overline{f(x+h_3)}f(x+h_1+h_3)f(x+h_2+h_3)\overline{f(x+h_1+h_2+h_3)}\Big)^{1/8}.$$

Demostrar que $\|f_s\|_{U^3(\mathbb{Z}_p)} = 1$ y $\|1_A\|_{U^3(\mathbb{Z}_p)} \gg 1$.

**Ejercicio 9.** † A partir del teorema de Roth, deducir la versión siguiente demostrada por (Varnavides 1959): para todo $\alpha > 0$ existe $c > 0$ tal que todo conjunto $A \subset [N]$ de densidad $\alpha$ verifica $\mathbb{E}_{x,d\in[N]} 1_A(x)1_A(x+d)1_A(x+2d) \geqslant c$. (Un resultado análogo se da para el teorema de Meshulam.)

**Ejercicio 10.** † Una ecuación lineal $c_1x_1 + \cdots + c_tx_t = 0$ es *invariante por traslación* (o simplemente *invariante*) si $c_1 + \cdots + c_t = 0$. Verificar que la prueba del teorema 3.3 se puede adaptar para obtener un resultado análogo para cualquier tal ecuación con $t \geqslant 3$. Se dice de un sistema de ecuaciones $Mx = 0$ con $M \in \mathbb{Z}^{r\times m}$ que es *invariante* si $M(1,\ldots,1) = 0$, y que tiene *complejidad 1* si, para funciones $f : \mathbb{Z}_p \to \mathbb{C}$ con $\|f\|_\infty \leqslant 1$, el promedio $\mathbb{E}_{x\in\mathbb{Z}_p^m : Mx=0} f(x_1)\cdots f(x_m)$ tiende a 0 cuando $\|f\|_u \to 0$. (Por ejemplo, los sistemas de una sola ecuación lineal con al menos tres variables son de complejidad 1; esto se ve siguiendo la prueba de la proposición 2.7.) Verificar que la prueba del teorema 3.3 se puede generalizar a todo sistema invariante de complejidad 1.

(Los ejercicios marcados con † son más difíciles.)

Pablo Candela
pablo.candela@uam.es
Departamento de Matemáticas
Universidad Autónoma de Madrid
Ciudad Universitaria de Cantoblanco
28049 - Madrid
España

# 3 Introducción a la teoría de las curvas elípticas

# INTRODUCCIÓN A LA TEORÍA DE LAS CURVAS ELÍPTICAS

Marusia Rebolledo y Marc Hindry

**Resumen**

La primera parte de este curso expone las bases de la teoría de las curvas elípticas: definiciones, ecuaciones, ley de grupo y luego desarrolla los aspectos más aritméticos: teorema de Hasse (sobre un cuerpo finito), alturas y teorema de Mordell–Weil (sobre un cuerpo de números). la segunda parte introduce la altura canónica de Néron–Tate y el regulador de una curva elíptica. Con el fin de limitar el regulador, se introducen las funciones $L$ y zeta asociadas a un cuerpo de números, a una variedad algebraica, a una representación de Galois o a una forma modular. En el caso de una curva elíptica, las tres últimas son unificadas a través del teorema de Wiles.

## Índice general

## 1 Prefacio

Este curso es una introducción a la teoría geométrica y aritmética de curvas elípticas y a las funciones $L$ que les pueden ser asociadas. Como referencias generales sobre curvas elípticas citamos a Knapp (1992), Milne (2006), Silverman (1986, 1994) y Silverman y Tate (1992).

Tales curvas aparecen naturalmente en el estudio de ecuaciones diofánticas; corresponden al primer ejemplo donde no se puede aplicar sistemáticamente, como sí se hace para las cónicas, el principio de Hasse. La riqueza de estas curvas viene en parte del hecho que el método de la cuerda y tangente proporciona una ley de grupo. Son también las variedades abelianas más simples (dimensión 1). Las estructuras diversas de estas curvas y sus nexos, vía las funciones $L$, con objetos de naturaleza algebraica (representaciones de Galois) o analítica (formas modulares), están en el corazón de numerosos resultados y preguntas actuales en la geometría aritmética. Entre estos resultados, el más conocido es ciertamente el Último Teorema de Fermat. La función $L$ que permite estos nexos, es una serie de Dirichlet del mismo tipo que la función zeta de Riemann. Las series de Dirichlet y la función zeta de Riemann fueron introducidas para demostrar los principales teoremas acerca de la distribución de números primos.

El éxito de este método ha llevado a introducir funciones análogas llamadas funciones $L$ de Hasse–Weil asociadas a las curvas elípticas. Presentaremos en este curso estas series, sus principales propiedades – algunas apenas conjeturadas, tales como la "hipótesis de Riemann" – y sus relaciones con la aritmética de las curvas elípticas.

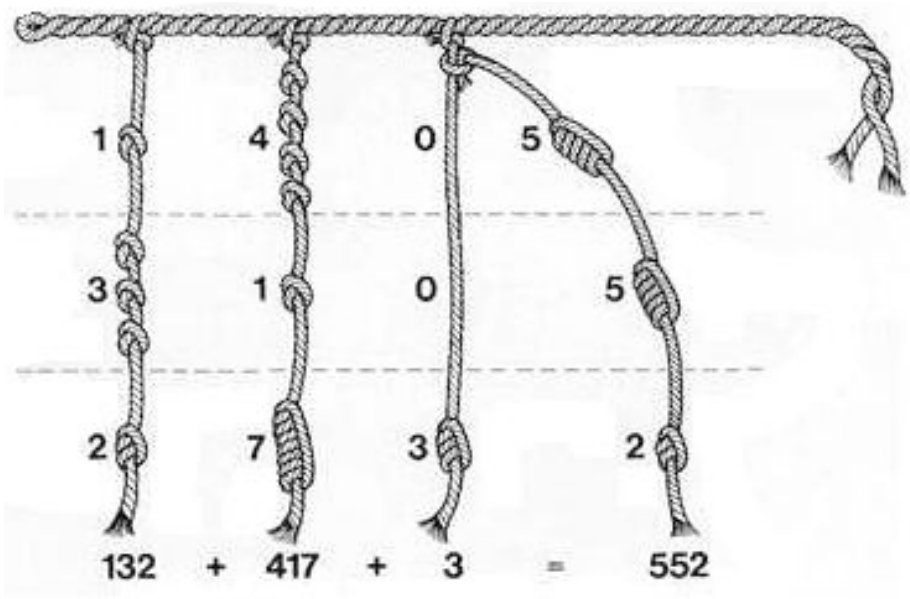

*Completamos este prefacio a través de un ejercicio con dos preguntas que esconden el estudio y la teoría de curvas elípticas. Se encontrarán indicaciones para la solución al final de estas notas.*

**Ejercicio 1.** (En el estilo de Fermat, número se refiere a número natural)

1. ¿Cuáles cubos pueden escribirse como un cuadrado aumentado de dos unidades?

2. ¿Cuáles números pueden escribirse como un producto de dos números consecutivos *y* de tres números consecutivos (ejemplo: $6 = 2 \times 3 = 1 \times 2 \times 3$)?

Advertencia: *En este capítulo, los pequeños* (*Ejer.*) *señalan que hay algo que demostrar en la frase anterior que se deja como ejercicio al lector.*

# 2 Curvas elípticas

## 2.1 Preliminares.

**2.1.1 Espacios afines y proyectivos .** Sea $K$ un cuerpo y $\bar{K}$ una clausura algebraica de $K$. En este curso, lo usual será $K = \mathbb{Q}, \mathbb{Q}_p, \mathbb{C}$, o $\mathbb{F}_q$.

Llamamos *espacio afín de dimensión n sobre* $K$ al espacio vectorial: $\mathbb{A}^n = \mathbb{A}^n_K = \mathbb{A}^n(\bar{K}) = \bar{K}^n$. Llamamos *espacio proyectivo de dimensión n sobre* $K$ al conjunto de las rectas vectoriales de $\bar{K}^{n+1}$ o sea al cociente

$$\mathbb{P}^n = \mathbb{P}^n(\bar{K}) = (\mathbb{A}^{n+1} - \{0\})/\sim$$

por la relación de equivalencia, donde $(x_0, \ldots, x_n) \sim (y_0, \ldots, y_n)$ si y sólo si existe $\lambda \in \bar{K}^*$ tal que $x_i = \lambda y_i$ para todo $i$.

Denotamos por $(x_0 : \cdots : x_n)$ la clase en $\mathbb{P}^n$ de un elemento $(x_0, \ldots, x_n) \in \mathbb{A}^{n+1}$.

**Ejercicio 2.** Verificar que $\sim$ es una relación de equivalencia.

El grupo de Galois $G_K = \mathrm{Gal}(\bar{K}/K)$ actúa sobre $\mathbb{A}^n$: para $\sigma \in G_K$, definimos $(x_1, \ldots, x_n)^\sigma = (x_1^\sigma, \ldots, x_n^\sigma)$ donde $x^\sigma = \sigma(x)$. De manera análoga, $G_K = \mathrm{Gal}(\bar{K}/K)$ actúa sobre $\mathbb{P}^n$. Denotamos $\mathbb{A}^n(K) = \{(x_1, \ldots, x_n) \in \mathbb{A}^n; x_i \in K, i = 1, \ldots n\}$ y $\mathbb{P}^n(K) = \{(x_0 : \cdots : x_n) \in \mathbb{A}^{n+1}; x_i \in K, i = 0, \ldots n\}$ a los elementos invariantes por la acción de $G_K$. Cuidado: $P = (y_0, \ldots, y_n) \in \mathbb{P}^n(K)$ no significa que $y_i \in K$ para todo $i$, sino que existe $\lambda \in \bar{K}^*$ tal que $\lambda y_i \in K$ para todo $i$.

Consideremos $H = \{(x_0 : \cdots : x_n) \in \mathbb{P}^n; x_0 = 0\}$ y $U = \mathbb{P}^n \setminus H$ (definido por una resta de conjuntos). La aplicación

$$\phi : U \longrightarrow \mathbb{A}^n; (x_0 : \cdots : x_n) \mapsto \left(\frac{x_1}{x_0}, \ldots, \frac{x_n}{x_0}\right)$$

es bien definida, biyectiva y de inversa $(y_1, \ldots, y_n) \mapsto (1 : y_1 : \cdots : y_n)$. (Ejer.)

Así, hemos inyectado una copia de $\mathbb{A}^n$ en $\mathbb{P}^n$ de manera que $\mathbb{P}^n$ es la unión disjunta de $U \cong \mathbb{A}^n$ y $H$ que llamamos *hiperplano en el infinito*. Con un procedimiento análogo para $x_i \neq 0$ en lugar de $x_0 \neq 0$, obtenemos otras $n$ inmersiones de $\mathbb{A}^n$ en $\mathbb{P}^n$ y por cada una, un hiperplano en el infinito.

**Ejemplo 1.** La recta proyectiva es la unión de $U_0$ y del punto en el infinito $(0 : 1)$, donde $U_0 \cong \mathbb{A}^1$ via $(x_0 : x_1), x_0 \neq 0 \mapsto (1 : \frac{x_1}{x_0})$. El plano proyectivo $\mathbb{P}^2$ es la unión de $U_0$ y de la recta en el infinito $\mathbb{P}^2 \setminus U_0 = \{(0 : x : y); x, y \in \bar{K}\} \cong \mathbb{P}^1$.

**2.1.2 Curvas planas.** Una *curva algebraica afín plana* es el lugar $C$ en $\mathbb{A}^2$ de los ceros de un polinomio no constante $F \in \bar{K}[x, y]$, o sea, $C = V(F) = \{(x, y) \in \mathbb{A}^2; F(x, y) = 0\}$. Una *curva algebraica proyectiva plana* es el lugar $C \subset \mathbb{P}^2$ de los ceros de un polinomio *homogéneo* no constante $H \in \bar{K}[X, Y, Z]$.

Recordamos que $H \in \bar{K}[X, Y, Z]$ de grado total $d$ es *homogéneo* si para todo $\lambda \in \bar{K}$, $H(\lambda X, \lambda Y, \lambda Z) = \lambda^d H(X, Y, Z)$ (de manera que tiene sentido hablar del lugar de los ceros de $H$ en $\mathbb{P}^2$).

**Ejemplo 2.** Decimos que la curva es una *recta* afín o proyectiva (respectivamente una cónica, una cúbica) si un polinomio que la define es de grado 1 (respectivamente 2, 3) y es homogéneo en el caso proyectivo.

**Paso de lo afín a lo proyectivo:** el plano proyectivo $\mathbb{P}^2$ puede ser cubierto por tres abiertos afines: $U_i = \{(x_0 : x_1 : x_2) \in \mathbb{P}^2; x_i \neq 0\} \cong \mathbb{A}^2, i = 0, 1, 2$, como ha sido ya explicado. Entonces, una curva proyectiva $C = V(H)$ para $H \in \bar{K}[X, Y, Z]$ homogéneo y no constante, es la unión de tres curvas $C_i = C \cap U_i$, $i = 0, 1, 2$ que se identifican con curvas afines, una vez que hemos identificado $U_i$ con $\mathbb{A}^2$. Más precisamente, $C_0 \cong \{(y, z) \in \mathbb{A}^2; H(1, y, z) = 0\}$, $C_1 \cong \{(x, z) \in \mathbb{A}^2; H(x, 1, z) = 0\}$, $C_2 \cong \{(x, y) \in \mathbb{A}^2; H(x, y, 1) = 0\}$. Las curvas afines $C_0, C_1, C_2$ son llamadas *cartas afines* de $C$.

Recíprocamente, si $C'$ es una curva afín plana dada por un polinomio $F \in \overline{K}[x, y]$ de grado $d$, existe una curva proyectiva $C$ que contiene $C'$, a saber, la curva dada por el *polinomio homogeneizado* de $F$: $H(X, Y, Z) := Z^d F(X/Z, Y/Z)$. (Ejer.) La curva $C$ es llamada *completación proyectiva de* $C'$ y la curva $C'$ es una de las cartas afines de $C$; la recubrimos haciendo $C \cap U_2$. En esa carta afín $C'$, los puntos de $C \setminus C'$ son llamados *puntos en el infinito*.

**Ejemplo 3.** Sea $L$ la recta proyectiva definida por $aX + bY + cZ = 0$. Entonces la recta $L$ es la unión de la recta afín $L_2 = L \cap U_2$ dada por la ecuación $ax + by + c = 0$ y de un punto en el infinito (el conjunto de los puntos en el infinito para esta carta afín): $L \cap \{(X : Y : Z) \in L; Z = 0\} = \{(X : Y : Z); aX = -bY\} = \{(-b : a : 0)\}$.

**Ejemplo 4.** Sea $C$ la cúbica proyectiva definida por $Y^2Z = X^3 + Z^3$. En la carta afín $C_2 = C \cap U_2$ definida por $y^2 = x^3 + 1$, hay un solo punto en el infinito pues $C \backslash C_2 = \{(X : Y : 0); X^3 = 0\} = \{(0 : 1 : 0)\}$.

**Puntos racionales.** Si un polinomio definiendo $C$ tiene sus coeficientes en el cuerpo $K$, se dice que la curva está *definida sobre* $K$. En este caso, $G_K$ actúa sobre los puntos de $C$. En efecto, si $P = (x_0, y_0) \in C$ y $\sigma \in G_K$, se define $P^\sigma = (x_0^\sigma, y_0^\sigma)$. Éste es un punto bien definido que anula también al polinomio que define $C$. Ejer. Se puede definir lo mismo con curvas proyectivas. Denotamos $C(K)$ al conjunto de los puntos invariantes por esa acción, llamados *puntos $K$-racionales*.

Cuidado: $C$ no es determinada por sus puntos racionales sobre $K$ (de acuerdo al siguiente ejemplo).

**Ejemplo 5.** Sea $C = V(X^2 + 1)$. Entonces, $C(\mathbb{R}) = \varnothing$.

Denotamos $C_F$ la curva (afín o proyectiva) definida por un polinomio $F$ sin factor cuadrado.

**Singularidades, criterio del Jacobiano –** Decimos que un punto $P$ de una curva afín plana $C_F$ es *singular* si el Jacobiano $\left(\frac{\partial F}{\partial x}, \frac{\partial F}{\partial y}\right)$ se anula en $P$. Decimos que un punto de una curva proyectiva $C$ es singular si lo es para una carta afín que lo contiene, lo cual no depende de la carta afín elegida (como lo muestra el siguiente ejercicio).

**Ejercicio 3.** Sea $P$ un punto de una curva proyectiva $C_H$. Demostrar que $P$ es singular si y sólo si $[H(P) = 0, \left(\frac{\partial H}{\partial X}, \frac{\partial H}{\partial Y}, \frac{\partial H}{\partial Z}\right)(P) = (0,0,0)]$ si y sólo si $P$ es singular en todas las cartas afines que lo contienen.

Para determinar la multiplicidad de un punto singular y las tangentes, se puede usar el desarrollo de Taylor de $F$ en $P = (x_0, y_0)$. En efecto, $F$ se puede escribir como $F(X,Y) = F_1(X - x_0, Y - y_0) + F_2(X - x_0, Y - y_0) + \cdots + F_d(X - x_0, Y - y_0)$ donde, para todo $i$, $F_i$ es un polinomio homogéneo de grado $i$ en $X - x_0$ y en $Y - y_0$. (Dichos polinomios son dados por el desarrollo de Taylor, por ejemplo $F_1(X,Y) = \frac{\partial F}{\partial x}(P)X + \frac{\partial F}{\partial y}(P)Y$.). Así, el punto $P$ es singular si y sólo si $F_1 = 0$. La *multiplicidad* de esa singularidad corresponde al entero más pequeño $m \geqslant 1$ tal que $F_m \neq 0$. En un punto no singular la tangente es la recta de ecuación $F_1(X - x_0, Y - y_0) = 0$. En un punto singular doble, las tangentes son dadas por los factores de la forma cuadrática $F_2(X - x_0, Y - y_0) = 0$.

Una curva es dicha *lisa* si no tiene ningún punto singular.

**Ejemplo 6.** Consideremos la curva $C = C_F \subset \mathbb{A}^2$, con $F(x,y) = y^2 - x^3$. Tenemos $\frac{\partial F}{\partial x}(x,y) = -3x^2$ y $\frac{\partial F}{\partial y}(x,y) = 2y$, de manera que el único punto singular de $C$ es $P = (0,0)$. Del hecho que $F(x,y) = y^2 - x^3 = F_2(x,y) + F_3(x,y)$ con $F_2(x,y) = y^2$, concluimos que $P$ es un punto singular doble con una única tangente de ecuación $y = 0$. En este caso, decimos que $P$ es una *punta*.

**Ejemplo 7.** Sea $C = C_F$ con $F(x,y) = y^2 - x^3 - x^2 \subset \mathbb{A}^2$. Aquí tenemos $\frac{\partial F}{\partial x}(x,y) = -3x^2 - 2x$ y $\frac{\partial F}{\partial y}(x,y) = 2y$, y entonces el único punto singular es $P = (0,0)$. Esta vez $F_2(x,y) = y^2 - x^2$, de tal forma que $P$ es un punto singular doble pero con dos tangentes distintas definidas respectivamente por las ecuaciones $(y - x = 0)$ y $(y + x = 0)$. En este caso, decimos que $P$ es un *nodo*.

**Ejercicio 4.** Sea $C = V(y^2 - x^3 - ax - b) \subset \mathbb{A}^2$, donde $a, b \in K$, $\mathrm{char}(K) \neq 2, 3$. Demostrar que $C$ es lisa si y sólo si $x^3 + ax + b$ no tiene raíz doble, es decir, $4a^3 + 27b^2 \neq 0$.

**2.2 Curvas elípticas: propiedades geométricas.** En esa sección, consideramos de nuevo $K$ un cuerpo y $\bar{K}$ una clausura algebraica de $K$.

**2.2.1 Curvas definidas por una ecuación de Weierstrass.** Sea $C$ una curva proyectiva plana definida por una *ecuación de Weierstrass*:

(1) $$y^2 + a_1xy + a_3y = x^3 + a_2x^2 + a_4x + a_6 \qquad (a_i \in \overline{K}, i = 1, \ldots, 6)$$

o sea, dicho más correctamente, por la ecuación homogénea asociada $Y^2Z + a_1XYZ + a_3YZ^2 = X^3 + a_2X^2Z + a_4XZ^2 + a_6Z^3$. Esto quiere decir que la ecuación (1) corresponde a la ecuación de una carta afín de $C$, a saber la carta $C \cap U_2$. En esa carta, los puntos en el infinito son dados por el conjunto $C \cap \{(X : Y : Z); Z = 0\} = \{(0 : 1 : 0\}$, de manera que solo hay un punto en el infinito.

Si para todo $i \in 1, \ldots, 6$ $a_i \in K$, $C$ es definida sobre $K$.

Cuando $K$ es de característica $\neq 2$, con el cambio de variables $y \mapsto \frac{1}{2}(y - a_1x - a_3)$, obtenemos una ecuación más simple para $C$ (Ejer.) (ver Silverman (1986) Cap. III.1):

$$y^2 = 4x^3 + b_2x^2 + 2b_4x + b_6 \quad \text{donde} \quad b_2 = a_1^2 + 4a_2,\ b_4 = 2a_4 + a_1a_3,\ b_6 = a_3^2 + 4a_6.$$

Cuando $\mathrm{char}(K) \neq 2, 3$, podemos también eliminar el término $x^2$ mediante el cambio de variables $(x,y) \mapsto \left(\frac{x - 3b_2}{36}, \frac{y}{108}\right)$, obteniendo una ecuación llamada *ecuación de Weierstrass reducida*, de la forma siguiente:

$$y^2 = x^3 - 27c_4x - 54c_6, \quad \text{donde} \quad c_4 = b_2^2 - 24b_4 \quad \text{y} \quad c_6 = -b_2^3 + 36b_2b_4 - 216b_6.$$

Los cambios de variables que respetan la forma de la ecuación y dejan invariante el punto en el infinito son exactamente los cambios de la forma

$$(x,y) \mapsto (u^2x + r, u^3y + u^2sx + t)$$

donde $u, r, s, t \in \overline{K}, u \neq 0$. Tales cambios son llamados *admisibles* (ver Silverman (ibíd.) Cap. III.1 para más detalles).

Definimos el *discriminante de la ecuación* (1) por

(2) $$\Delta = -b_2^2b_8 - 8b_4^3 - 27b_6^2 + 9b_2b_4b_6,$$

donde $b_i, i = 2, 4, 6$ son definidos como antes y donde $b_8 = a_1^2a_6 + 4a_2a_6 - a_1a_3a_4 + a_2a_3^2 - a_4^2$.

Un cambio de variables admisible $(u, r, s, t)$ trasforma $\Delta$ en $u^{-12}\Delta$ (de acuerdo a Silverman (ibíd.) Cap. III.1 Tabla 1.2).

**Ejemplo 8.** El discriminante de una ecuación de Weierstrass reducida $y^2 = x^3 + Ax + B$ es $\Delta = -16(4A^3 + 27B^2)$.

Se pueden estudiar las singularidades de la curva $C$ por medio de $\Delta$ y $c_4$:

**Proposición 1.** *Sea $C$ dada por una ecuación de Weierstrass* (1). *La curva $C$ es no singular si y sólo si $\Delta \neq 0$. Si $\Delta = 0$, entonces $C$ tiene solo un punto singular que es un nodo si (además) $c_4 \neq 0$ o una punta si (además) $c_4 = 0$.*

*Prueba.* Consideremos la ecuación homogénea de $C$:

$$F(X, Y, Z) = Y^2Z + a_1XYZ + a_3YZ^2 - X^3 - a_2X^2Z - a_4XZ^2 - a_6Z^3 = 0.$$

En el punto en el infinito $O = [0, 1, 0]$, tenemos $\frac{\partial F}{\partial Z}(O) = 1 \neq 0$, así que el punto $O$ nunca es singular.

Para simplificar las notaciones, supongamos que $K$ es de característica $\neq 2, 3$. Entonces $C$ tiene una ecuación reducida $f(x, y) = y^2 - x^3 + 27c_4x + 54c_6 = 0$ (ver Silverman (ibíd., Cap. III. Proposición 1.4) para el caso general). Tenemos $\frac{\partial f}{\partial y}(x, y) = 2y$. Así, un punto $(x, y)$ es singular si y sólo si $y = 0$ y $x$ es un cero múltiple de $Q(X) = X^3 - 27c_4X - 54c_6$. Tal punto existe si y sólo si el discriminante del polinomio $Q$ de grado 3 es cero, o sea $\Delta = 0$ y en este caso, esta singularidad es única, porque $Q$ tiene a lo más un cero múltiple.

Ahora, supongamos que $\Delta = 0$. Sea $P = (x_0, y_0)$ el único punto singular de $C$. Desde una ecuación de Weierstrass general hacemos los cambios de variables $x = x' + x_0$, $y = y' + y_0$ para enviar al punto $P$ en $(0, 0)$. Esos cambios dejan invariantes $\Delta$ y $c_4$, y entonces podemos suponer que $P = (0, 0)$ y que $C$ tiene una ecuación de la forma $y^2 + a_1xy - x^3 - a_2x^2 = 0$, porque $a_6 = f(0, 0) = 0, a_4 = \frac{\partial f}{\partial x}(0, 0) = 0, a_3 = \frac{\partial f}{\partial y}(0, 0) = 0$. Los factores de la forma cuadrática $F_2(x, y) = y^2 + a_1xy - a_2x^2$ son distintos si y sólo si el discriminante $a_1^2 + 4a_2$ de $F_2$ no es cero. Como $c_4 = (a_1^2 + 4a_2)^2$, eso demuestra que $P = (0, 0)$ es un nodo si y sólo $c_4 \neq 0$ o una punta en el caso contrario.

□

### 2.2.2 Curvas elípticas: definiciones.

**Definición 1.** Una *curva elíptica* $(E, O)$ *sobre* $K$ es el dato de una curva $E$, proyectiva, plana y no singular, definida sobre $K$ y de *género* 1, dotada de un punto racional $O \in E(K)$. Tales curvas pueden ser definidas por una ecuación de Weierstrass (1) con discriminante no nulo.

**Género geométrico**. No vamos a definirlo exactamente pero se trata de un invariante geométrico útil para clasificar las curvas planas. Por ejemplo, para una curva proyectiva plana no singular de grado $d$, el género geométrico es igual al género aritmético, o sea:

$$g = \frac{(d-1)(d-2)}{2}.$$

(Para una curva singular, se deben agregar a esta ecuación cantidades dependientes de las singularidades).

Así, se ve que una curva definida por una ecuación de Weierstrass (1) con discriminante $\Delta \neq 0$ (entonces lisa) y dotada del punto en el infinito $O = (0 : 1 : 0)$, es una curva elíptica definida sobre todo cuerpo que contiene los coeficientes de la ecuación. Ejer.

Recíprocamente, se puede demostrar que toda curva elíptica puede ser definida por una ecuación de Weierstrass (con discriminante $\Delta \neq 0$). Esto requiere el teorema de Riemann–Roch y lo admitimos (ver Silverman (1986, Cap. III.3)).

**Ejemplo 9.** Los puntos $\mathbb{R}$-racionales de la curva elíptica $E : y^2 = x^3 - 3x + 3$ son dados por el dibujo siguiente:

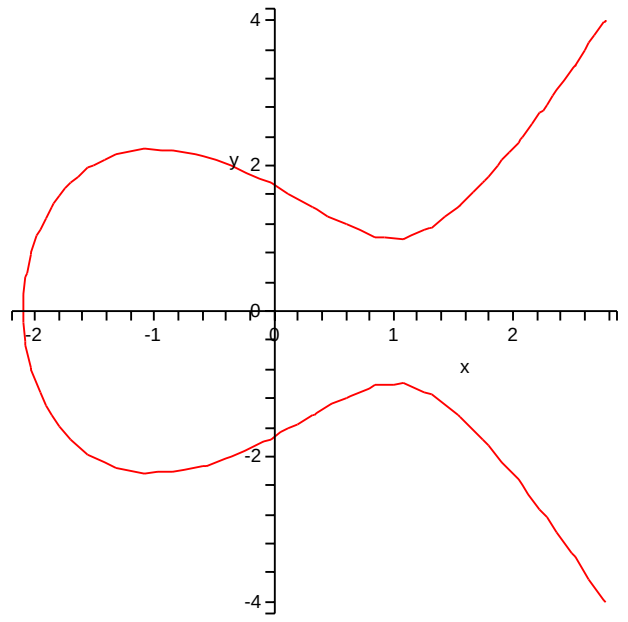

donde hay que pensar en $O$ como un punto en "el infinito".

Los puntos $\mathbb{R}$-racionales de $E : y^2 = x^3 + x$ son dados por el siguiente gráfico:

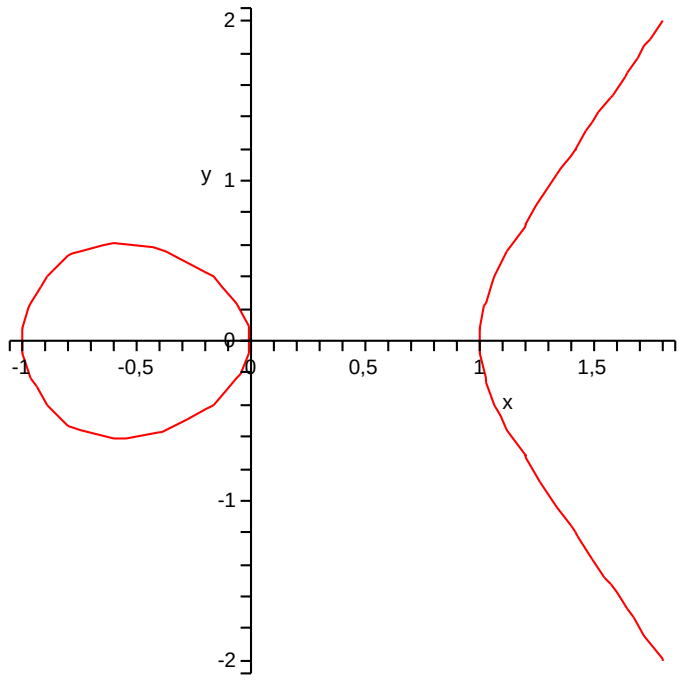

**2.2.3 Ley de grupo. Multiplicidad de intersección con una recta**: Sea $L : aX + bY + cZ = 0$ una recta de $\mathbb{P}^2$, $C_F$ una curva proyectiva plana y $P \in C_F \cap L \neq \varnothing$. Podemos suponer sin pérdida de generalidad que $P = (x_0 : y_0 : z_0)$ está en la carta afín $U_2$, o sea que $z_0 \neq 0$. La recta afín $L \cap U_2$ se puede parametrizar por $x(t) = bt + x_0,\ y(t) = -at + y_0$. Sea $f(x, y) = F(x : y : 1)$ definiendo la carta afín $C \cap U_2$. La multiplicidad del cero $t = 0$ en el polinomio $f(x(t), y(t))$ en $t$ se llama *multiplicidad de intersección de* $L$ *con* $C_F$ *en* $P$ con la convención de que la multiplicidad es infinita si el dicho polinomio es nulo (i. e. $L \subset C_F$).

El siguiente teorema es un caso particular del teorema de Bézout geométrico:

**Teorema 1** (Bézout)**.** *Una recta proyectiva intersecta una curva proyectiva plana de grado* $m$ *en* $m$ *puntos (*$\overline{K}$*-racionales) contados con multiplicidades. Si la curva proyectiva es de grado* $3$ *y está definida sobre* $K$ *y si dos de los puntos de intersección son* $K$*-racionales, entonces el tercero es también* $K$*-racional.*

**Ejercicio 5.**

1. Demostrar que la multiplicidad de intersección de una recta $L$ con una curva $C$ en un punto $P$ es superior a la multiplicidad de $P$ sobre $C$.
2. Demostrar que una cúbica proyectiva plana tiene a los más un punto singular y que es un punto doble.
3. Demostrar, usando el teorema de Bézout, que toda curva proyectiva plana no singular es irreducible.

Sea $E$ una curva elíptica sobre un cuerpo $K$. Definimos sobre $E$ una ley de composición interna del siguiente modo: sean $P, Q \in E$ y $L \subset \mathbb{P}^2$ la recta pasando por $P$ y $Q$ (la tangente a $E$ si $P = Q$). Como $E$ es una cúbica, el teorema de Bézout demuestra que $L$ intersecta $E$ en un tercer punto $R$ (donde contamos las multiplicidades). Sea $L'$ la recta pasando por $R$ y $O$. Denotamos $P + Q$ el tercer punto donde $L'$ corta a $E$.

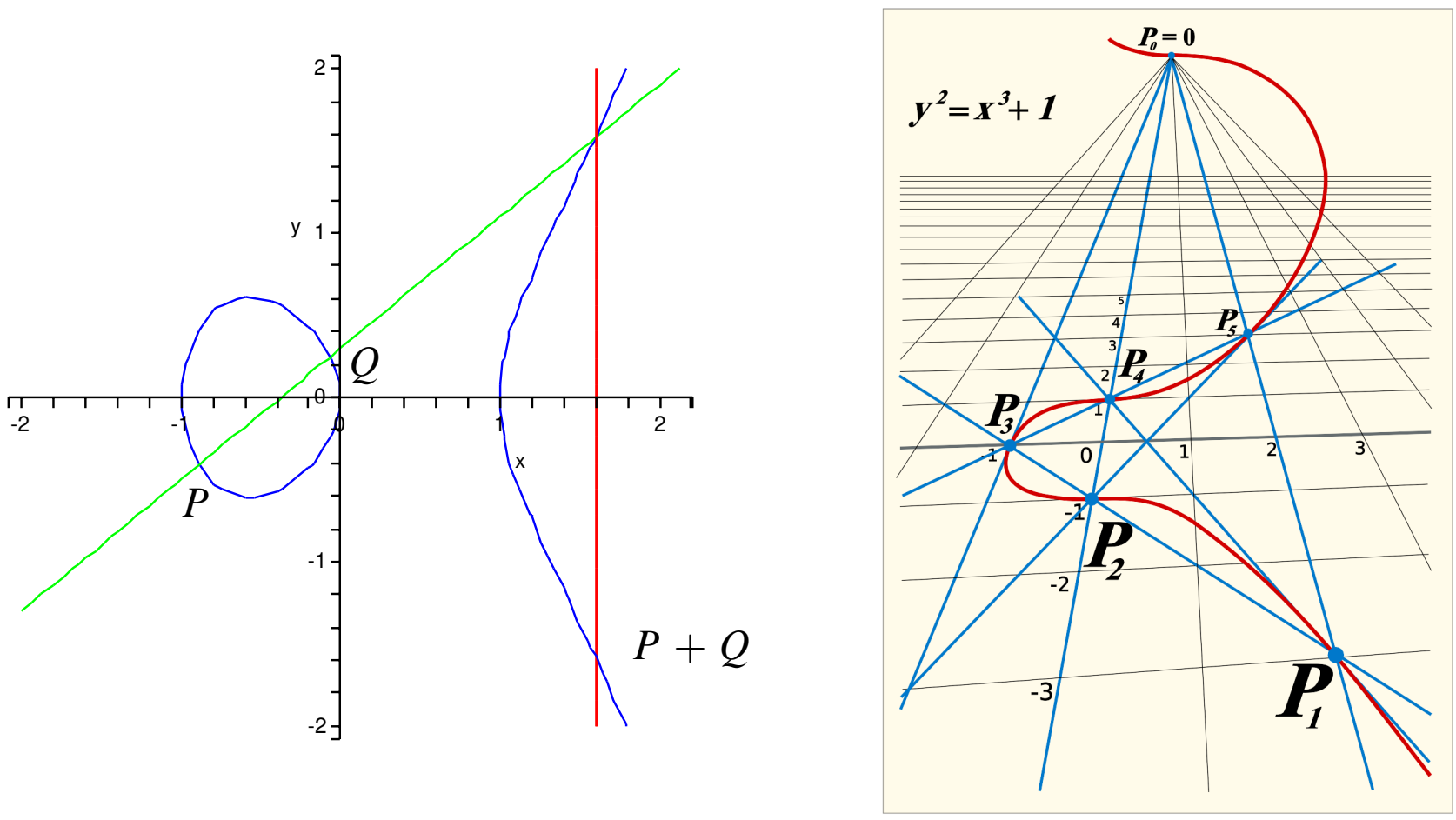

Figura 1: Torsion on cubic curve by Jean Brette for Wikipedia. Licensed under CC BY-SA

**Proposición 2.** *La ley de composición precedente define sobre* $E$ *una estructura de grupo abeliano de elemento neutro* $O$*.*

*Prueba.* La ley es interna por definición. Denotamos por $(PQ)$ al tercer punto (contado con multiplicidad) de intersección de $E$ con la recta pasando por dos puntos $P$ y $Q$ de $E$, de tal modo que $P + Q = (O(PQ))$. Como $(PQ) = (QP)$ obtenemos la conmutatividad. Por otro lado, $O + P = (O(OP)) = P$ para todo punto $P$. Finalmente, para todo punto $P \in E$, tenemos que $P + P' = O$ donde $P' = ((OO)P)$: por definición, los puntos $P$, $P'$, $(OO)$ están sobre una recta $L$, entonces $P + P' = (O(PP')) = (O(OO))$. Además, $(OO)$ es el tercer punto de intersección de $E$ con la tangente a $E$ en $O$ (pues $O$ es un punto de multiplicidad $\geqslant 2$ con su tangente). Entonces $(O(OO)) = O$.

La parte difícil de la demostración, que admitimos aquí, es la asociatividad. El lector se puede referir a Knapp (1992, Cap. III.3) para una prueba elemental o a Silverman (1986, Cap. III.3.4) para una prueba usando el teorema de Riemann–Roch. □

Si elegimos otro punto base $O'$, la ley $+'$ obtenida verifica que $P +' Q = P + Q - O'$ y las estructuras de grupos $(E, +')$ y $(E, +)$ son isomorfas (via $P \mapsto P - O'$). (Ejer.)

Hay fórmulas explícitas para la ley de grupo (ver Silverman (ibíd., Cap. III.2.3) o ejercicio 10). Por ejemplo, si $E$ tiene ecuación (1) y si $P = (x : y : 1)$, entonces $-P = (x, -y - a_1x - a_3)$ y (*fórmula de duplicación*)

$$x(P + P) = x([2]P) = \frac{x^4 - b_4x^2 - 2b_6x - b_8}{4x^3 + b_2x^2 + 2b_4x + b_6}. \tag{3}$$

**2.2.4 Morfismos e isogenias. Funciones regulares.–** Sea $C$ una curva afín plana. Una aplicación $f : C \longrightarrow \mathbb{A}^1$ es una *función regular* si proviene de un polinomio, o sea $f = ((a, b) \in C \mapsto F(a, b))$ donde $F \in \overline{K}[X, Y]$. El conjunto de las funciones regulares sobre $C$ es un anillo que denotamos $\overline{K}[C]$. Si $C$ es además irreducible, es decir no contiene una curva estricta, entonces $\overline{K}[C]$ es un anillo integral. Denotamos $\overline{K}(C) = \mathrm{Frac}(\overline{K}[C])$ al *cuerpo de las funciones racionales de* $C$. Una función $f \in \overline{K}(C)$ es llamada *regular en un punto* $P \in C$ si existen $g, h \in \overline{K}[C]$ tales que $h(P) \neq 0$ y $f = g/h$. Si $f \in \overline{K}(C)$ es regular en todos los puntos de $C$ entonces $f \in \overline{K}[C]$ (entonces la terminología es consistente).

Para una curva proyectiva plana $C$, definimos el *cuerpo de las funciones racionales* como el cuerpo de las funciones de una carta afín de $C$. Eso no depende de la carta elegida y lo denotamos $\overline{K}(C)$. (Ejer.) Una función es *regular en* $P$ si lo es en una carta afín que contiene a $P$. Podemos demostrar que si $C$ es una curva proyectiva plana irreducible, entonces toda función $f \in \overline{K}(C)$ regular sobre $C$ es constante.

**Morfismos de curvas proyectivas** Sean $C$ y $C' = C_H$ dos curvas proyectivas planas donde $C$ irreducible. Una *aplicación racional de* $C$ *hacia* $C'$ está dada por $\phi = (f_0 : f_1 : f_2)$ con $f_i \in \overline{K}(C)$ no todas cero y tal que $H(f_0, f_1, f_2) = 0$. Decimos que esta es una aplicación *regular en* $P$ si existe $g \in \overline{K}(C)^\times$ tal que $gf_i$ $(i = 0, 1, 2)$ son regulares en $P$ y no se anulan simultáneamente en $P$. Un *morfismo* es una aplicación racional que es regular en todos puntos. Si $C$ es lisa, toda aplicación racional desde $C$ es un morfismo.

***j*-invariante de una curva elíptica.–** Sea $E$ una curva elíptica sobre $K$ dada por una ecuación de Weierstrass

$$y^2 + a_1xy + a_3y = x^3 + a_2x^2 + a_4x + a_6 \qquad (a_i \in \overline{K}, i = 1, \dots, 6)$$

denotamos $j = j(E) = \frac{c_4^3}{\Delta}$ al *$j$-invariante de $E$* (con la notación ya introducida para $c_4$). Con esta definición tenemos en particular que si $E$ está dada por una ecuación reducida $y^2 = x^3 + Ax + B$, entonces $j(E) = -1728(4A)^3/\Delta$.

**Proposición 3.** *Dos curvas elípticas sobre $K$ son isomorfas sobre $\bar{K}$ si y sólo si tienen el mismo $j$-invariante. Para $j_0 \in \bar{K}$ dado, existe a menos de un $\bar{K}$-isomorfismo, una única curva elíptica (definida sobre $K(j_0)$) con j-invariante igual a $j_0$.*

*Prueba.*

1. Se puede verificar que los cambios admisibles dejan el $j$-invariante estable. Si $E \cong E'$, existe un cambio admisible pasando de una ecuación de Weierstrass de $E$ a la de $E'$ y entonces $j(E) = j(E')$. Recíprocamente, supongamos que $j(E) = j(E')$ y, para simplificar la demostración, que $\operatorname{car}(K) \neq 2, 3$. Entonces $E$ es dada por una ecuación de la forma $y^2 = x^3 + Ax + B$, $E'$ por otra de la forma $y^2 = x^3 + A'x + B'$. Como $j(E) = j(E')$, tenemos que $\frac{(4A)^3}{4A^3+27B^2} = \frac{(4A')^3}{4A'^3+27B'^2}$ lo cual implica que $A^3B'^2 = A'^3B^2$. Así, se puede verificar sin dificultad que el cambio de variables admisible $(x, y) = (u^2x', u^3y')$ con el siguiente $u$, produce un isomorfismo $E \longrightarrow E'$:

$$u = \begin{cases} (A/A')^{1/4} = (B/B')^{1/6} & \text{si } A \neq 0 \text{ y } B \neq 0 \\ (B/B')^{1/6} & \text{si } A = 0 \\ (A/A')^{1/4} & \text{si } B = 0. \end{cases}$$

   Esas expresiones son bien definidas porque $A = 0$ si y sólo si $j(E) = 0 = j(E')$ y eso si y sólo si $A' = 0$. Además, si $A = 0$ $(= A')$, entonces $B \neq 0$ y $B' \neq 0$ ya que $\Delta \neq 0$ y $\Delta' \neq 0$.

2. Sea $j_0 \in \bar{K}$. Si $j_0 \neq 0, 1728$, entonces la curva definida por la ecuación $y^2 + xy = x^3 - \frac{36}{j_0-1728}x - \frac{1}{j_0-1728}$ tiene discriminante $\Delta = \frac{j_0^2}{(j_0-1728)^3}$ y $j$-invariante $j_0$, como deseado.

   Y si $j_0 = 0$ (respectivamente $j_0 = 1728$), basta considerar la curva definida por $y^2 + y = x^3$ (respectivamente $y^2 = x^3 + x$), que tiene discriminante $-27$ (respectivamente $-64$) y $j$-invariante 0 (respectivamente 1728).

□

**Isogenias.–** Una aplicación $E_1 \longrightarrow E_2$ es una *isogenia* si es un morfismo tal que $\phi(O) = O$. Denotamos $\operatorname{Hom}(E_1, E_2)$ al conjunto de las isogenias que van desde $E_1$ hacia $E_2$ y $\operatorname{End}(E)$ al conjunto de los endomorfismos de $E$, es decir, de las isogenias $E \longrightarrow E$.

Admitimos los dos hechos siguientes (ver Silverman (ibíd., Cap. III.4)). Una isogenia no constante es sobreyectiva. Una isogenia $\phi : E_1 \longrightarrow E_2$ induce una inyección de cuerpos de funciones $\phi^* : \bar{K}(E_2) \longrightarrow \bar{K}(E_1)$. Llamamos *grado* de $\phi$ al grado de la extensión de cuerpos $\bar{K}(E_1)/\phi^*(\bar{K}(E_2))$ y decimos que $\phi$ es *separable* (respectivamente *inseparable*, respectivamente *puramente inseparable*) si la extensión de cuerpos de funciones correspondiente lo es.

**Ejemplo 10** (Endomorfismo de Frobenius)**.** Sea $E$ una curva elíptica definida sobre un cuerpo finito $\mathbb{F}_q$. La aplicación definida por $\phi(O) = O$ y $\varphi : (x, y) \mapsto (x^q, y^q)$ para todo punto $P = (x : y : 1) \neq O$ es un endomorfismo de $E$, llamado *endomorfismo de Frobenius*. Esta es una isogenia inseparable de grado $q$. Usando las fórmulas explícitas para la ley de grupo, se puede ver también que se trata de un endomorfismo de grupo. (Ejer.) Más generalmente, tenemos el teorema siguiente.

**Teorema 2.** *Sea* $\varphi \in \mathrm{Hom}(E_1, E_2)$ *una isogenia de curvas elípticas. Entonces para todos* $P, Q \in E_1$*, tenemos* $\varphi(P + Q) = \varphi(P) + \varphi(Q)$*. Además, el núcleo de* $\varphi$ *es un subgrupo de* $E_1$ *de orden igual al grado separable de* $\varphi$*. Recíprocamente, para todo subgrupo* $C$ *de* $E_1$ *existe una única curva elíptica* $E_2$ *y una isogenia separable* $\phi : E_1 \longrightarrow E_2$ *de núcleo* $C$*. En ese caso, denotamos* $E_2 := E_1/C$*.*

*Prueba.* Ver los teoremas 4.8, 4.10 y la proposición 4.12 de Silverman (1986, Cap. III.). □

**Ejemplo 11** (multiplicación por un entero)**.** Sea $E$ una curva elíptica. Se define de manera natural la *multiplicación por un entero* $m \in \mathbb{Z}$ denotada $[m] : E \longrightarrow E$, según la ley de grupo de $E$. Denotamos por $E[m]$ al núcleo de $[m]$, el cual es llamado *conjunto de los puntos de* $m$*-torsión* de $E$. Denotamos $E_{tors}$ *el conjunto de los puntos de torsión* $E$, es decir $E_{tors} = \bigcup_{m \geqslant 1} E[m]$. Si $E$ está definida sobre $K$, entonces $E_{tors}(K)$ es el subgrupo de los puntos de orden finito en $E(K)$.

**Teorema 3.** *Sea* $E$ *una curva elíptica sobre un cuerpo* $K$ *de característica* $\ell$*.*

1. *Si* $m \neq 0$*,* $[m]$ *es una isogenia no constante de grado* $m^2$*.*
2. *Si* $\ell = 0$ *o* $(\ell, m) = 1$*,* $[m]$ *es separable y tenemos* $E[m] \cong (\mathbb{Z}/m\mathbb{Z})^2$*;*
3. *Si* $\ell \neq 0$*, tenemos* $E[\ell^e] = \{O\}$ *para todo entero* $e \geqslant 1$*, o* $E[\ell^e] \cong \mathbb{Z}/\ell^e\mathbb{Z}$ *para todo entero* $e \geqslant 1$*.*

**Teorema 4** (isogenia dual)**.** *Para toda isogenia no constante* $\phi : E_1 \longrightarrow E_2$ *de grado* $m$*, existe una única isogenia* $\widehat{\phi} : E_2 \longrightarrow E_1$ *llamada* isogenia dual de $\phi$*, tal que* $\widehat{\phi} \circ \phi = [m]$*.*

No detallamos más aquí esta noción, ni la demostración del teorema 3. El lector puede referirse a Silverman (ibíd., Cap. III.6), en particular a los teoremas 6.1, 6.2 y al corolario 6.4.

**2.3 Puntos racionales de una curva elíptica (resultados).** Sea $E$ una curva elíptica sobre un cuerpo $K$. La estructura del conjunto de los puntos $K$-racionales $E(K)$ depende mucho de la naturaleza de $K$. Nos interesa particularmente el caso donde $K$ es un cuerpo de números (o sea una extensión finita de $\mathbb{Q}$) o un cuerpo finito.

**2.3.1 Curvas elípticas sobre $\mathbb{C}$.** Cuando $K = \mathbb{C}$, $E(\mathbb{C})$ tiene una estructura de variedad analítica isomorfa a un toro.

**Teorema 5** (teorema de uniformización)**.** *Sea* $E$ *una curva elíptica sobre* $\mathbb{C}$*. Existe un retículo* $\Lambda \subset \mathbb{C}$ *único salvo por una homotecia y un isomorfismo analítico complejo* $\alpha : \mathbb{C}/\Lambda \longrightarrow E(\mathbb{C})$ *de grupos de Lie complejos.*

Podemos dar explícitamente el isomorfismo $\alpha$ usando las *funciones de Weierstrass*: para más detalles, ver en Silverman (ibíd., pág. VI.5). El lector puede también consultar la presentación de R. Miatello en este volumen Miatello (2020).

**2.3.2 Curvas elípticas sobre un cuerpo finito.** Sea $K = \mathbb{F}_q$ "el" cuerpo finito con $q$ elementos donde $q = p^k$ y $p$ es un número primo.

**Teorema 6** (Hasse). *Para todo $m \geqslant 1$, el grupo $E(\mathbb{F}_{q^m})$ es un grupo abeliano finito y*

$$\sharp E(\mathbb{F}_{q^m}) = q^m + 1 - \alpha^m - \beta^m \quad \textit{con} \quad |\alpha| = |\beta| = \sqrt{q}.$$

*En particular,*

$$|\sharp E(\mathbb{F}_q) - q - 1| \leqslant 2\sqrt{q}.$$

*Observación* 1. Para $q = p$ primo, Deuring demostró que el teorema de Hasse es óptimo, es decir que para todo $a \in \mathbb{Z}$ tal que $|a| \leqslant 2\sqrt{p}$, existe $E/\mathbb{F}_p$ tal que $\sharp E(\mathbb{F}_p) = 1 + p - a$.

*Prueba.* Sea $\varphi$ el endomorfismo de Frobenius. Entonces $E(\mathbb{F}_q) = \{P \in E; \varphi(P) = P\}$ y $E(\mathbb{F}_{q^m}) = \{P \in E; \varphi^m(P) = P\} = \ker(1 - \varphi^m)$. Como $1 - \varphi^m$ es separable, obtenemos $\sharp E(\mathbb{F}_{q^m}) = \deg(1 - \varphi^m)$.

Sea $\ell \neq p$ un número primo. Un endomorfismo $\psi \in \mathrm{End}(E)$ induce para todo entero $n \geqslant 1$ un endomorfismo de $E[\ell^n]$ y entonces induce un endomorfismo $\psi_\ell \in \mathrm{End}(T_\ell(E))$ donde $T_\ell(E)$ es el limite proyectivo de los grupos de $\ell^n$-torsión $E[\ell^n]$. Como $\ell \neq p$, $E[\ell^n]$ es un $(\mathbb{Z}/\ell\mathbb{Z})$-módulo libre de rango 2, y entonces $T_\ell(E)$ es un $\mathbb{Z}_\ell$-módulo libre de rango 2, llamado *módulo de Tate* (ver Silverman (1986, Cap. III.7)). Así, después de la elección de una base de $T_\ell(E)$, podemos considerar $\psi_\ell$ como un elemento de $\mathrm{GL}_2(\mathbb{Z}_\ell)$. Su determinante y su traza verifican (conforme a Silverman (ibíd., V. Proposition 2.3)):

$$\det(\psi_\ell) = \deg(\psi) \quad \text{y} \quad \mathrm{Tr}\,(\psi_\ell) = 1 + \deg(\psi) - \deg(1 - \psi).$$

En particular, tales cantidades están en $\mathbb{Z}$.

Dado que $\sharp E(\mathbb{F}_{q^m}) = \deg(1 - \varphi^m) = \det(1 - \varphi_\ell^m)$, consideremos el polinomio característico de $\varphi_\ell^m$. Tenemos $\det(T - \varphi_\ell^m) = (T - \alpha^m)(T - \beta^m)$ donde $\alpha, \beta$ son las raíces de $\det(T - \varphi_\ell)$ en $\mathbb{C}$. Por lo que precede, el polinomio característico $\det(T - \varphi_\ell) = T^2 - \mathrm{Tr}\,(\varphi_\ell)T + \det(\varphi_\ell)$, de grado 2, tiene coeficientes en $\mathbb{Z}$ y valores positivos sobre los racionales. En efecto, tenemos $\det((m/n) - \varphi_\ell) = \det(m - n\varphi_\ell)/n^2 = \deg(m - n\varphi)/n^2 \geqslant 0$ $(m/n \in \mathbb{Q})$. Así, tiene una raíz doble o dos raíces complejas conjugadas, y de ahí concluimos que $|\alpha| = |\beta|$. Finalmente, $\alpha\beta = \det(\varphi_\ell) = \deg(\varphi) = q$, por lo que $|\alpha| = |\beta| = \sqrt{q}$ y $\sharp E(\mathbb{F}_{q^m}) = (1 - \alpha^m)(1 - \beta^m) = 1 + q^m - \alpha^m - \beta^m$.

□

Recordamos que $E[p^e]$ es cero para todo $e \geqslant 1$, o cíclico de orden $p^e$ para todo $e \geqslant 1$ (en vista del teorema 3). En el primer caso, decimos que la curva elíptica $E$ es *supersingular* sobre $\mathbb{F}_q$. En el segundo caso, decimos que es *ordinaria*.

**2.3.3 Curvas elípticas sobre un cuerpo local.** Si $E$ es una curva elíptica sobre un cuerpo $p$-ádico $K$ (por ejemplo $K = \mathbb{Q}_p$), entonces veremos que $E(K)$ tiene una estructura de grupo de Lie $p$-ádico compacto: es una extensión de un grupo finito por un pro-$p$-grupo $E_1(K)$ (ver Sección 2.4).

**2.3.4 Curvas elípticas sobre un cuerpo de números.**

**Teorema 7** (Mordell–Weil)**.** *Sea $E$ una curva elíptica sobre un cuerpo de números $K$. El grupo $E(K)$, llamado* grupo de Mordell–Weil*, es un grupo abeliano de tipo finito, es decir, existe un isomorfismo de grupos*

$$E(K) \cong \mathbb{Z}^r \oplus E(K)_{tors}$$

*con $r \in \mathbb{Z}_{\geqslant 0}$ y donde $E(K)_{tors}$ es el grupo formado de los puntos de torsión de $E(K)$. En particular $E(K)_{tors}$ es un grupo finito.*

La prueba de este teorema es basada sobre el *método de descenso* de Fermat. El lector puede referirse a Darmon (2004) 1.2 para una breve prueba de este hecho o a Silverman (1986) Cap. VIII por una prueba detallada y completa. Ver también la sección 2.6 de estas notas.

Uno puede determinar sencillamente el grupo de torsión $E(K)_{tors}$ de una curva elíptica dada (ver Sección 3). En cambio, el *rango $r$* es mucho más misterioso y es el objeto de varias conjeturas como por aquella de Birch y Swinnerton-Dyer (ver conjetura 3, sección 3.6.2 ) que lo une al comportamiento de una cierta función de una variedad compleja asociada a la curva elíptica llamada *función $L$ de Hasse–Weil* (según, sección 3.5). Ver la introcción de la sección 3.1.

**2.4 Curvas elípticas sobre un cuerpo local y reducción.** En esta sección, $K$ es un cuerpo local completo para una valuación discreta $v$, por ejemplo $K$ es $\mathbb{Q}_p$ o una extensión finita de $\mathbb{Q}_p$. Denotamos $R$ a su anillo de valuación discreta, $k$ al cuerpo residual, $p$ a su característica y $\pi \in \mathcal{O}_K$ a una uniformizante. (Por ejemplo, $K = \mathbb{Q}_p, R = \mathbb{Z}_p, k = \mathbb{F}_p, \pi = p$).

**2.4.1 Reducción de una curva elíptica. Ecuación minimal.–** Sea $E$ una curva elíptica dada por una ecuación de Weierstrass con coeficientes en $K$:

$$y^2 + a_1xy + a_3y = x^3 + a_2x^2 + a_4x + a_6 \qquad (a_i \in K, i = 1, \ldots, 6). \tag{4}$$

Un cambio de variables admisible $(x', y') = (u^{-2}x, u^{-3}y)$ cambia a $a_i$ en $a'_i = u^i a_i$ y a $\Delta$ en $\Delta' = u^{-12}\Delta$. Si tomamos $u$ como una buena potencia de $\pi$ que elimina los denominadores de los $a_i$, obtenemos una nueva ecuación con coeficientes $a_i \in R$ para todo $i$. (Basta elegir $v(u) \geqslant -\min(v(a_i)/i)$ (Ejer.) ). En este caso, el nuevo discriminante está también en $R$, o sea $v(\Delta) \geqslant 0$. Siendo un conjunto no vacío de $\mathbb{Z}_{\geqslant 0}$, el conjunto de las valuaciones $v(\Delta)$ de discriminante de ecuaciones de Weiestrass de $E$ con coeficientes en $R$, admite un mínimo. Así, una *ecuación de Weierstrass minimal de $E$ en $\pi$* es una ecuación (4) con $a_i \in R, i = 2, 4, 6$ y $v(\Delta)$ minimal para todas ecuaciones con coeficientes en $R$.

**Ejercicio 6.** 1. Demostrar que si la ecuación tiene coeficientes en $R$ y $v(\Delta) < 12$ entonces la ecuación es minimal. (Lo mismo sucede si $v(c_4) < 4$ o si $v(c_6) < 6$.)

2. Demostrar que si $p \neq 2, 3$ y la ecuación es minimal entonces $v(c_4) < 4$ o $v(\Delta) < 12$.

**Proposición 4.** *Una ecuación de Weiestrass minimal para E sobre K es única salvo un cambio admisible de la forma* $(x, y) \mapsto (u^2x + r, u^3y + u^2sx + t)$ *con* $u \in R^\times$ *y* $r, s, t \in R$.

*Prueba.* Consultar Silverman (ibíd., VII. Prop. 1.3). □

**Reducción.–** Sea $E$ una curva elíptica sobre $K$ con ecuación minimal

$$y^2 + a_1xy + a_3y = x^3 + a_2x^2 + a_4x + a_6 \qquad (a_i \in R, i = 1, \dots, 6).$$

Reduciendo los coeficientes $a_i$ módulo $\pi$, obtenemos una ecuación de Weiestrass sobre $k$. La curva definida por esta ecuación es una cúbica (posiblemente singular) llamada *reducción módulo* $\pi$ *de* $E$ y denotada $\widetilde{E}$. La proposición 4 muestra que la ecuación obtenida por reducción es única salvo un cambio de variables admisible para las ecuaciones de Weierstrass sobre $k$.

La curva $\widetilde{E}$ sobre $k$ es de uno de los tipos siguientes:

1. $\widetilde{E}$ es una curva elíptica sobre $k$: eso ocurre cuando $\Delta \in R^\times$. Decimos que $E$ tiene *buena reducción módulo* $\pi$;

2. $\widetilde{E}$ tiene una singularidad que es un nodo: decimos que $E$ *tiene mala reducción multiplicativa* módulo $\pi$. Eso ocurre cuando $\Delta \equiv 0 \pmod{\pi}$ y $c_4 \in R^\times$.

3. $\widetilde{E}$ tiene una singularidad que es una punta: decimos que $E$ *tiene mala reducción aditiva* módulo $\pi$. Eso ocurre cuando $\Delta \equiv c_4 \equiv 0 \pmod{\pi}$.

En el caso de reducción multiplicativa, decimos además que la reducción es *desplegada* si las pendientes de las tangentes en el punto singular están en $k$ y *no desplegada* si no.

**2.4.2 Aplicación de reducción.** Sean $E$ una curva elíptica definida sobre $K$ y $\widetilde{E}$ su reducción módulo $\pi$, quizás singular.

Si $P = (x : y : z) \in \mathbb{P}^2(K)$, existe $\lambda \in K^\times$, tal que $\lambda x, \lambda y, \lambda z$ están en $R$ y al menos uno de ellos está en $R^\times$. Así, podemos definir $\widetilde{P} := (\lambda x \pmod{\pi} : \lambda y \pmod{\pi} : \lambda z \pmod{\pi}) \in \mathbb{P}^2(k)$. Esto define una aplicación llamada *aplicación de reducción* $E(K) \longrightarrow \widetilde{E}(k)$.

La ley de composición cuerdas/tangente definida como en la sección 2.2.3 define un grupo abeliano sobre el subconjunto $\widetilde{E}_{ns}$ de los puntos no singulares de $\widetilde{E}$ (ver Silverman (ibíd., Cap. III.2.5)).

**Proposición 5.** *La aplicación de reducción define una secuencia exacta de grupos abelianos*

$$0 \longrightarrow E_1(K) \longrightarrow E_0(K) \longrightarrow \widetilde{E}_{ns}(k) \longrightarrow 0$$

*donde* $E_0(K) = \{P \in E(K); \widetilde{P} \in \widetilde{E}_{ns}(k)\}$ *y donde* $E_1(K)$ *es el núcleo de la aplicación de reducción.*

*Prueba.* El subconjunto $E_0(K)$ es un subgrupo de $E(K)$ porque una recta entre dos puntos no singulares intersecta a una cúbica en un tercer punto no singular. La aplicación de reducción es un morfismo de grupos porque la ley de grupo es definida por el método cuerda/tangente y la imagen de una recta de $\mathbb{P}^2(K)$ es una recta de $\mathbb{P}^2(k)$. Por definición de $E_1(K)$, la secuencia es un complejo y tenemos una inyección a derecha. Basta demostrar entonces que la aplicación de reducción es sobreyectiva.

La prueba es una aplicación del *lema de Hensel* siguiente: sea $f \in R[X]$, no cero, y $\alpha \in R$ tal que $v(f(\alpha)) > 0$ y $v(f'(\alpha)) = 0$, entonces existe un único $\beta \in R$ tal que $f(\beta) = 0$ y $v(\beta - \alpha) > 0$.

El punto $\widetilde{O} \in \widetilde{E}_{ns}(k)$ es la imagen de $O \in E_0(K)$. Sea $Q = (a, b) \in \widetilde{E}_{ns}(k)$. Queremos levantar $a$ y $b$ en $K$ manteniéndonos siempre sobre la curva. Sea $F(x, y) = 0$ una ecuación minimal de $E$. Denotamos $\bar{F}(x, y) = 0$ a su correspondiente reducción. Como $Q$ es no singular, al menos uno de $\frac{\partial \bar{F}}{\partial X}(a, b)$ o $\frac{\partial \bar{F}}{\partial Y}(a, b)$ no es cero, por ejemplo $\frac{\partial \bar{F}}{\partial X}(a, b) \neq 0$. Sea $y \in R$ tal que $y \equiv b \bmod \pi$. Consideremos el polinomio $P(X) = F(X, y) \in R[X]$. Su reducción $\bar{P} \in k[X]$ verifica $\bar{P}(a) = 0$ y $\bar{P}'(a) \neq 0$. Entonces para $\alpha \in R$ tal que $\alpha \equiv a \bmod \pi$, tenemos $v(P(\alpha)) > 0$ y $v(P'(\alpha)) = 0$. Por el lema de Hensel, deducimos que existe $x \in R$ tal que $P(x) = 0$ y $v(x - \alpha) > 0$. El punto $(x, y)$ está entonces en $E_0(K)$ y tiene reducción $Q$. □

Así, la aplicación de reducción define un isomorfismo $E_0(K)/E_1(K) \cong \widetilde{E}_{ns}(k)$. En particular, este grupo es finito. Tenemos también la proposición siguiente (ver Silverman (1986, VII. Corolario 6.2), o Milne (2006, II. Teorema 4.1) para una prueba sobre $\mathbb{Q}_p$):

**Proposición 6.** *El subgrupo $E_0(K)$ es de índice finito en $E(K)$.*

Un punto $P = (x : y : z)$ está en $E_1(K)$ si $\widetilde{P} = \widetilde{O} = (0 : 1 : 0)$ o sea si existe $\lambda \in K^\times$ tal que $\lambda y \in R^\times$ y $\lambda x, \lambda z \in (\pi)$ es decir, si $v(x) - v(y) \geqslant 1$ y $v(z) - v(y) \geqslant 1$. Eso motiva la definición siguiente: $E_n(K) := \{P = (x : y : z) \in E(K); v(x) - v(y) \geqslant n\}$, para todo entero $n > 1$.

**Proposición 7.** *La filtración $E(K) \supset E_0(K) \supset E_1(K) \supset E_2(K) \cdots \supset E_n(K) \supset \ldots$ tiene las propiedades siguientes:*

1. *para todo $n \geqslant 1$, $E_n(K)$ es un subgrupo de $E_1(K)$ y la aplicación*

$$\begin{array}{rccl} \varpi_n : & E_n(K) & \longrightarrow & k \\ & P & \mapsto & \pi^{-n} x(P)/y(P) \pmod{\pi} \end{array}$$

   *define un morfismo de grupos sobreyectivo de núcleo $E_{n+1}(K)$.*
   *(Por lo tanto, $E_n(K)/E_{n+1}(K) \cong k$).*

2. *la filtración es exhaustiva, es decir, $\cap_{n \geqslant 0} E_n(K) = \{O\}$.*

*Prueba.* 1. Procedemos por recurrencia sobre $n$. Ya vimos que $E_1(K)$ es un subgrupo de $E_0(K)$. Supongamos que $E_n(K)$ es un subgrupo de $E(K)$ y demostremos que $\varpi_n : E_n(K) \longrightarrow k$ es un morfismo de grupos sobreyectivo de núcleo $E_{n+1}(K)$.

Un punto $P \in E_n(K)$ puede ser representado por $P = (x : y : z)$ con $y \in R^\times, x \in \pi^n R$. Escribamos $x = \pi^n x_0$ con $x_0 \in R$. La ecuación de $E$ es verificada por $P$:

$$(5) \qquad y^2 z + a_1 \pi^n x_0 yz + a_3 yz^2 = \pi^{3n} x_0^3 + a_2 \pi^{2n} x_0^2 z + a_4 \pi^n x_0 z^2 + a_6 z^3,$$

cual nos da $v(z) = v(y^2 z + a_1 \pi^n x_0 yz + a_3 yz^2) = v(\pi^{3n} x_0^3 + a_2 \pi^{2n} x_0^2 z + a_4 \pi^n x_0 z^2 + a_6 z^3) \geqslant 3n$. Escribamos $z = \pi^{3n} z_0,\ z_0 \in R$. La ecuación (5) se vuelve

$$y^2 z_0 + a_1 \pi^n x_0 y z_0 + a_3 \pi^{3n} y z_0^2 = x_0^3 + a_2 \pi^{2n} x_0^2 z_0 + a_4 \pi^{4n} x_0 z_0^2 + a_6 \pi^{6n} z_0^3,$$

es decir, el punto $P_0 := (x_0 : y : z_0)$ verifica una ecuación de Weierstrass. Se puede verificar además que esta ecuación tiene discriminante no cero y minimal para $v$. Sea $C$ la curva elíptica definida por esta ecuación. Su reducción $\widetilde{C}$ tiene ecuación $Y^2 Z = X^3$, así que es singular con una punta. En este caso, la aplicación $(x, y) \mapsto x/y$ define un isomorfismo $\widetilde{C}_{ns}(k) \cong k$ (Ejer., conforme a Silverman (1986, ejercicio 3.5)). Además, como $y \in R^\times$, tenemos $\widetilde{P}_0 \in \widetilde{C}_{ns}(k)$. Así la composición de las aplicaciones

$$\begin{array}{ccccccc} E_n(K) & \longrightarrow & C_0(K) & \longrightarrow & \widetilde{C}_{ns}(k) & \xrightarrow{\sim} & k \\ P = (x : y : z) & \mapsto & P_0 = (x_0 : y : z_0) & \mapsto & \widetilde{P}_0 & \mapsto & x_0/y = \pi^{-n} x/y \end{array}$$

define un morfismo de grupos de $E_n(K)$ en $k$ que es precisamente la aplicación $\varpi_n$. Este morfismo es sobreyectivo porque $P \mapsto P_0$ y la aplicación de reducción en $C$ lo son. Para terminar, $P$ está en el núcleo de $\varpi_n$ si y sólo si $\widetilde{P}_0 = (0 : 1 : 0)$, es decir, $x_0 \in (\pi), z_0 \in (\pi)$ (y $y \in R^\times$) y luego $P \in E_{n+1}(K)$.

2. Si $P = (x : y : z) \in \cap_{n \geqslant 1} E_n(K)$ entonces $x = 0$ (porque $v(x) \geqslant n$ para todo $n \geqslant 1$). Por lo tanto, la ecuación de $E$ da $y^2 z + a_3 yz^2 = a_6 z^3$. Si $z = 0$ entonces $P = O$. Si $z \neq 0$, entonces $y^2 + a_3 yz = a_6 z^2$ y $v(z) = 0$ lo que contradice el hecho de que $P \in E_1(K)$. □

De las proposiciones 5 y 7 se deduce el resultado siguiente que será útil para determinar los puntos de torsión y para demostrar el teorema de Mordell–Weil:

**Proposición 8.** *Sea $E$ una curva elíptica sobre $K$ y $m$ un entero primo a $p = \mathrm{char}(k)$. Tenemos:*

1. $E_1(K)[m] = \{O\}$;
2. *Si $E$ tiene buena reducción entonces la aplicación de reducción $E(K)[m] \longrightarrow \widetilde{E}(k)$ es inyectiva.*

*Prueba.* Observamos que 2. se deduce de 1. y de la secuencia exacta de la proposición 5.

Sea $P \in E_1(K)[m]$ con $m$ primo a $p$. Sea $n \geqslant 1$ un entero tal que $P \in E_n(K)$. Como $[m]P = O$ entonces $m\varpi_n(P) = \varpi_n([m]P) = 0$. Como $m$ primo a $p$, tenemos entonces $\varpi_n(P) = 0$ y luego $P \in E_{n+1}(K)$. Así, por recurrencia, $P \in \cap_{n \geqslant 1} E_n(K)$, lo que termina la prueba, usando la proposición 7. □

El teorema siguiente da una cota para los denominadores de los puntos de torsión (ver Silverman (ibíd., VII. Teorema 3.4) para una demostración).

**Teorema 8** (Cassels). *Sean $E$ una curva elíptica sobre $K$ de ecuación de Weierstrass $y^2 + a_1xy + a_3y = x^3 + a_2x^2 + a_4x + a_6$ $(a_i \in R)$ y $P \in E(K)$ un punto de orden exactamente $m \geqslant 2$.*

1. *Si $m$ no es una potencia de $p$, entonces $x(P)$, $y(P)$ están en $R$.*
2. *Si $m = p^n$, entonces $\pi^{2r}x(P)$, $\pi^{3r}y(P)$ están en $R$, donde $r = \left[\frac{v(p)}{p^n - p^{n-1}}\right]$ (aquí $[x]$ es la parte entera de $x$).*

**2.5 Torsión de curvas elípticas sobre un cuerpo de números.** Ahora consideremos una curva elíptica $E$ sobre un cuerpo de números $K$ con anillo de enteros $R$.

Se puede usar la proposición 8 de la manera siguiente para determinar los puntos de torsión de $E$. Denotamos $v$ una valuación discreta de $K$ y $K_v$ el cuerpo completo asociado. Entonces $E(K) \subset E(K_v)$. Si encontramos una ecuación minimal de $E$ sobre $K_v$ y la reducimos, podemos aplicar la proposición a $E/K_v$. Haciendo el mismo procedimiento de reducción, pero para distintas valuaciones, podemos encontrar la torsión.

**Ejemplo 12.** Sea $E/\mathbb{Q} : y^2 + y = x^3 - x + 1$ de discriminante $\Delta = -13.47$. Como $v_2(\Delta) = 0 < 12$, la ecuación es minimal en 2 (y en todo número primo) y $\widetilde{E}$ mod 2 es no singular de ecuación $y^2 + y = x^3 + x + 1$. Tenemos $\widetilde{E}(\mathbb{F}_2) = \{0\}$ así que $E(\mathbb{Q})_{tors} = \{0\}$.

De manera análoga, pueden obtener las condiciones locales para los denominadores dadas por el teorema 8 en varias plazas del cuerpo de números $K$.

**Teorema 9.** *Sea $E$ una curva elíptica sobre un cuerpo de números $K$ dada por una ecuación de Weierstrass con coeficientes en $R$. Sea $P \in E(K)$ de orden exactamente $m \geqslant 2$.*

1. *Si $m$ no es una potencia de un primo, entonces $x(P), y(P)$ están en $R$.*
2. *Si $m = p^n$ es una potencia de un primo $p$, entonces para toda $v \in M_K$ finita,*

$$ord_v(x(P)) \geqslant -2r_v \quad y \quad ord_v(y(P)) \geqslant -3r_v$$

*donde $r_v = \left[\frac{ord_v(p)}{p^n - p^{n-1}}\right]$. En particular, cuando $ord_v(p) = 0$, $x(P), y(P)$ son $v$-enteros.*

Sobre $\mathbb{Q}$, el siguiente corolario ha sido demostrado anteriormente e independientemente por Lutz y Nagell.

**Corolario 1** (Lutz–Nagell). *Sea $E/\mathbb{Q}$ una curva elíptica con una ecuación de Weierstrass $y^2 = x^3 + Ax + B$ con $A, B \in \mathbb{Z}$ y sea $P \in E(\mathbb{Q})$ de orden finito, $P \neq O$. Entonces*

$$i)\ x(P), y(P) \in \mathbb{Z} \quad ; \quad ii)\ y(P) = 0 \ o \ y(P)^2 \mid 4A^3 + 27B^2.$$

*Prueba.* Sea $m$ el orden de $P$. Observemos que, cuando $K = \mathbb{Q}$, $\left[\frac{ord_v(p)}{p^n - p^{n-1}}\right] = \left[\frac{1}{p^n - p^{n-1}}\right]$ vale 0 a menos que $p = 2$ y $n = 1$. Entonces por teorema 9, si $m \neq 2$, $x(P), y(P)$ son $v-$enteros para todo $v$, luego $x(P), y(P) \in \mathbb{Z}$. Si $m = 2$, entonces $y(P) = 0$ y $x(P)$ es un entero, al ser raíz de un polinomio unitario con coeficientes enteros.

Si denotamos $\varphi(X) = X^4-2AX^2-8BX+A^2$ y $\psi(X) = X^3+AX+B$, tenemos $\varphi(x(P)) = 4\psi(x(P))x([2]P)$ por la fórmula de duplicación y $f(X)\varphi(X)-g(X)\psi(X) = 4A^3+27B^2$ donde $f(X) = 3X^2+4A$ y $g(X) = 3X^3-5AX-27B$. Como $y(P)^2 = \psi(x(P))$, obtenemos

$$y(P)^2(4f(x(P))x([2]P) - g(x(P))) = 4A^3 + 27B^2.$$

Aplicando $i)$ a los puntos $P$ y $[2]P$ que son de torsión, tenemos que $x(P), x([2]P)$ están en $\mathbb{Z}$, lo que demuestra que el factor de $y(P)^2$ es entero.

□

Lo anterior demuestra en particular que $E(\mathbb{Q})_{tors}$ es un grupo finito y da un algoritmo para determinarlo. En efecto, hay una lista finita de $y(P)$ verificando las condiciones 1,2, y por lo tanto una lista finita $\mathfrak{L}$ de puntos posibles para la torsión. Para $P \in \mathfrak{L}$, es fácil de ver cuándo $P$ no es de torsión porque en este caso, existe $n \in \mathbb{Z}$ tal que $nP$ no está en $\mathfrak{L}$. Entonces, elegimos $P \in \mathfrak{L}$, calculamos $P, [2]P, [3]P, \dots [n]P$ hasta obtener un punto no perteneciente a $\mathfrak{L}$ u obtener $[n]P = O$. A propósito de esto, el teorema de Mazur (1977) siguiente muestra que $n \leqslant 12$.

Este teorema describe todos los grupos de torsión posibles para las curvas elípticas sobre $\mathbb{Q}$. Estos son un número finito, así que pocos enteros pueden ser el orden de un punto de torsión de una curva elíptica sobre $\mathbb{Q}$.

**Teorema 10** (Mazur). *Sea $E$ una curva elíptica sobre $\mathbb{Q}$. Entonces $E(\mathbb{Q})_{tors}$ es uno de los grupos siguientes*

$$\mathbb{Z}/N\mathbb{Z},\ 1 \leqslant N \leqslant 10,\ N = 12 \qquad \mathbb{Z}/2\mathbb{Z} \times \mathbb{Z}/2N\mathbb{Z}\ ,\ 1 \leqslant N \leqslant 4.$$

*Además, cada uno de esos grupos puede ser realizado efectivamente como grupo de torsión de una curva elíptica sobre $\mathbb{Q}$.*

*Observación* 2. Por cada uno de los grupos dados por el teorema anterior, hay una infinidad de curvas elípticas sobre $\mathbb{Q}$ con este grupo para grupos de torsión sobre $\mathbb{Q}$ (Kubert) (esas curvas forman una familia de un parámetro).

Una pregunta natural es: ¿cómo se generaliza este resultado para $[K : \mathbb{Q}] > 1$? Se dispone de un resultado análogo para los cuerpos cuadráticos (Kamienny, Kenku, Momose). Pero cuando $[K : \mathbb{Q}] = 3$, la lista de los grupos posibles no es completa. El teorema de Kamienny para los cuerpos cuadráticos es uniforme: existe una cota para el tamaño de la torsión de todas las curvas elípticas sobre un cuerpo cuadrático. Eso se generaliza para un cuerpo de números:

**Teorema 11** (Merel). *Para todo $d \geqslant 1$, existe $B(d)$ tal que para toda curva elíptica sobre un cuerpo de números $K$ de grado $d$ sobre $\mathbb{Q}$, $|E(K)_{tors}| \leqslant B(d)$.*

**2.6 “Demostración” del teorema de Mordell (ideas esenciales).** Aquí vamos a dar los pasos claves de la demostración del teorema 7 sobre $\mathbb{Q}$ (teorema de Mordell). Para más detalles, el lector puede referirse a Silverman (1986, Cap. VIII).

Sea $E$ una curva elíptica sobre $\mathbb{Q}$. La demostración consta esencialmente de dos pasos:

1. demostrar que $E(\mathbb{Q})/2E(\mathbb{Q})$ es finito. Vamos a demostrar sin embargo un poco más: $E(\mathbb{Q})/mE(\mathbb{Q})$ es finito para todo entero $m$ (*teorema débil de Mordell*).

2. descenso: la idea de este paso es que no se puede "dividir " indefinidamente un punto racional por 2, ya que, como será visto, la división por 2 va disminuyendo la *altura* del punto y porque hay un número finito de puntos de altura pequeña.

### 2.6.1 Teorema débil de Mordell.

**Teorema 12.** *Sea $E$ una curva elíptica sobre un cuerpo de números $K$. Para todo entero $m \geqslant 2$, $E(K)/mE(K)$ es finito.*

Estamos interesados en enunciar ese teorema para todo cuerpo de números porque el primer paso de la demostración es de agrandar el cuerpo hasta que contenga las coordenadas de los puntos de $m$-torsión. El esquema de la demostración es el siguiente:

**Paso 1.–** El siguiente lema permite agrandar $K$ hasta que $E[m] \subset E(K)$.

**Lema 1.** *Si $L/K$ es galoisiana finita y si $E(L)/mE(L)$ es finito, entonces $E(K)/mE(K)$ es finito.*

**Paso 2.–** Supongamos ahora que $K$ es tal que $E[m] \subset E(K)$. Se define un emparejamiento $\lambda : E(K) \times G_K \longrightarrow E[m]$ llamado *emparejamiento de Kummer* (ver (6) y ejercicio 7).

**Lema 2.** *Sea $L = K([m]^{-1}E(K))$ el compositum de los cuerpos donde los puntos $Q$ tales que $mQ \in E(K)$ son definidos. El emparejamiento de Kummer induce un emparejamiento perfecto*

$$\lambda' : E(K)/mE(K) \times \mathrm{Gal}(L/K) \longrightarrow E[m].$$

De este lema, deducimos que $E(K)/mE(K)$ es finito si y sólo si $L/K$ es una extensión finita.

**Paso 3.–** Demostramos que las extensiones $K(Q)$ de $K$ generadas por las coordenadas de $Q \in [m]^{-1}E(K)$ son no ramificadas afuera de un conjunto finito de plazas de $K$. Un teorema de Minkowski muestra entonces que hay un número finito de tales extensiones $K(Q)$ de $K$ y deducimos de esto que $L/K$ es finita.

Vamos a demostrar ahora los lemas precedentes y el teorema débil de Mordell.

Sea $P \in E(K)$ y $Q \in E(\bar{K})$ tal que $mQ = P$. El punto $Q$ es $K$-racional si y solamente si $Q^\sigma - Q = O$ para todo $\sigma \in G_K = \mathrm{Gal}(\bar{K}/K)$. Queremos estudiar y entender la falta de $K$-racionalidad de $Q$. Por eso, para cada pareja $(P, Q)$ de puntos definidos como antes, consideramos la aplicación

$$\begin{array}{rccl} \lambda_{P,Q} : & G_K & \longrightarrow & E \\ & \sigma & \mapsto & Q^\sigma - Q. \end{array}$$

Observemos que

$$i)\ \lambda_{P,Q}(G_K) \subset E[m] \qquad ii)\ \text{si } \lambda_{P,Q} = \lambda_{P',Q'} \text{ entonces } P - P' \in mE(K).$$

El punto $i)$ emana del hecho que la acción de $G_K$ conmuta con la multiplicación por un entero. Para el punto $ii)$: si para todo $\sigma \in G_K$, $Q^\sigma - Q = Q'^\sigma - Q'$, entonces $Q - Q' \in E(K)$ y luego $P - P' = [m](Q - Q') \in mE(K)$.

Gracias a $i)$, $\lambda_{P,Q}$ define una aplicación desde $G_K$ en $E[m]$, hecho que denotamos $\lambda_{P,Q} \in Ap(G_K, E[m])$.

Vamos a utilizar esos dos hechos $i)$ y $ii)$ para demostrar el lema 1:

*Demostración del lema 1.* Supongamos que $L/K$ es una extensión de cuerpos tal que $E(L)/mE(L)$ es finito. La inyección $E(K) \hookrightarrow E(L)$ define un morfismo de grupos $\varphi : E(K)/mE(K) \longrightarrow E(L)/mE(L)$. Para demostrar que $E(K)/mE(K)$ es finito, basta demostrar que el núcleo $\Phi = (E(K) \cap mE(L))/mE(K)$ de $\varphi$ es finito. Sea $\mathfrak{R}$ un sistema de representantes del cociente $\Phi$ y sea $P \in \mathfrak{R}$. Como $P \in mE(L)$, existe $Q \in E(L)$ tal que $P = [m]Q$. Del hecho que $Q \in E(L)$, vemos que la restricción de $\lambda_{P,Q}$ a $G_L$ es la aplicación cero, por lo tanto $\lambda_{P,Q}$ induce una aplicación $\widetilde{\lambda}_{P,Q} : \mathrm{Gal}(L/K) \longrightarrow E[m]$. Así, el proceso que a la clase de $P$ en $\Phi$ asocia $\widetilde{\lambda}_{P,Q}$ define una aplicación desde $\Phi$ hacia el conjunto finito $Ap(\mathrm{Gal}(L/K), E[m])$. Por la observación $ii)$, esta aplicación es inyectiva. Luego $\Phi$ es finito. □

Gracias al lema 1, podemos suponer que $E[m] \subset E(K)$. De esta manera, $\lambda_{P,Q}$ no depende de $Q \in E(\bar{K})$ tal que $mQ = P$, porque si $Q'$ es tal que $mQ' = P = mQ$ entonces $Q - Q' \in E[m] \subset E(K)$ y $Q^\sigma = Q'^\sigma$ para todo $\sigma \in G_K$. Así, podemos considerar la aplicación

$$\begin{array}{rrcl} \lambda : & E(K) \times G_K & \longrightarrow & E[m] \\ & (P, \sigma) & \mapsto & \lambda(P, \sigma) = \lambda_P(\sigma) = Q^\sigma - Q \end{array} \tag{6}$$

donde $Q \in E(\bar{K})$ es tal que $mQ = P$. Consideremos ahora $L = K([m]^{-1}E(K))$.

**Ejercicio 7.** Demostrar los hechos siguientes:

1. la aplicación $\lambda$ es un emparejamiento bilineal, o sea que para todo par de puntos $P, P'$ en $E(K)$ y todo par de aplicaciones $\sigma, \tau$ en $G_K$,

$$\lambda(P + P', \sigma) = \lambda(P, \sigma) + \lambda(P', \sigma) \quad \text{y} \quad \lambda(P, \sigma\tau) = \lambda(P, \sigma) + \lambda(P, \tau).$$

   Este emparejamiento es llamado *emparejamiento de Kummer*;

2. para todo $P \in mE(K)$ y todo $\sigma \in G_K$, $\lambda(P, \sigma) = O$;

3. para todo $P \in E(K)$ y todo $\sigma \in G_L$, $\lambda(P, \sigma) = O$.

Así, en vista de los resultados del ejercicio 7, el emparejamiento $\lambda$ induce un emparejamiento bilineal

$$\lambda' : E(K)/mE(K) \times \mathrm{Gal}(L/K) \longrightarrow E[m].$$

*Demostración del lema 2.* Por la observación $ii)$, el morfismo $E(K)/mE(K) \longrightarrow \mathrm{Hom}(\mathrm{Gal}(L/K),$ inducido por $\lambda$, es inyectivo. Además, $L$ es generado sobre $K$ por las coordenadas de los puntos $Q$ tal que $mQ = P$, $P \in E(K)$. Por lo tanto, si $\sigma \in \mathrm{Gal}(L/K)$ es tal que para todo $P \in E(K)$, $\lambda(P, \sigma) = Q^\sigma - Q = O$, entonces $\sigma = \mathrm{id}$. Eso demuestra que $\mathrm{Gal}(L/K) \longrightarrow \mathrm{Hom}(E(K)/mE(K), E[m])$ también es inyectivo. Entonces, $\lambda'$ es un emparejamiento perfecto. □

*Paso 3 y fin de la demostración del teorema 12.* Por el lema 2, $E(K)/mE(K)$ es finito si y solamente si $L/K$ es una extensión finita. Esto es lo que vamos a demostrar ahora.

Sea $P \in E(K)$ y $K(Q)/K$ la extensión generada por las coordenadas de $Q \in E(\bar{K})$ tal que $mQ = P$. Vamos a demostrar que $K(Q)$ es no ramificada afuera del conjunto $S$ definido como la unión de las plazas de mala reducción de $E$ sobre $K$, de las plazas $v$ de $K$ tales que $v(m) \neq 0$, y de las plazas en el infinito de $K$. (Para una definición de lo que significa "no ramificada" ver la subsección 3.4.1, en particular la definición 11).

Sea $\mathfrak{p} \notin S$ un ideal de $K$ y $\mathfrak{P}$ un ideal primo de $K(Q)$ sobre $\mathfrak{p}$. Sea $\sigma$ en el grupo de inercia $I_{\mathfrak{P}/\mathfrak{p}} = \{\sigma \in \mathrm{Gal}(K(Q)/K)\ ;\ \sigma(x) - x \in \mathfrak{P}, (x \in \mathcal{O}_{K(Q)})\}$. Sea $P' \in E(K)$ y $Q' \in E(\bar{K})$ tal que $mQ' = P'$. Como $\widetilde{Q'^{\sigma}} = \widetilde{Q}'$ (mod $\mathfrak{P}$), la reducción de $\lambda(P', \sigma) = Q'^{\sigma} - Q$ módulo $\mathfrak{P}$ es cero. Como $E$ tiene buena reducción en $\mathfrak{p} \notin S$, tenemos que $\lambda(P', \sigma) = O$ (por la proposición 8). Y eso sucede para todo $P' \in E(K)$. Como el emparejamiento $\lambda'$ es perfecto, deducimos que $\sigma = \mathrm{id}$. Obtenemos que $I_{\mathfrak{P}} = \{\mathrm{id}\}$ para todo $\mathfrak{P} \mid \mathfrak{p}$, ($\mathfrak{p} \notin S$). Así $K(Q)/K$ es no ramificada afuera de $S$.

Para terminar, usamos el teorema siguiente de Minkowski: Sea $S \subset M_K$ un conjunto finito de plazas de un cuerpo de números $K$. Entonces, hay un número finito de extensiones de grado acotado sobre $K$ y no ramificadas afuera de $S$.

Aplicándolo, obtenemos que, como la extensión $K(Q)$ de $K$ tiene grado $\leqslant m^2$ para $Q \in [m]^{-1}E(K)$ y es no ramificada afuera de $S$, hay un número finito de dichas extensiones finitas $K(Q)$, $(Q \in [m]^{-1}E(K))$. Así, $L/K$ es finita.

□

**2.6.2 Descenso.** Vamos a definir una noción de altura *ingenua* que depende de la ecuación de Weierstrass. Eso bastará para demostrar el teorema de Mordell. Para demostrar sin embargo el teorema de Mordell–Weil (sobre un cuerpo de números de grado $> 1$), se necesita una noción de altura independiente de la ecuación: la altura de Neron–Tate que será definida en el la sección 3.1 definición 4.

**Definición 2** (Altura de un racional)**.** Sea $x = u/v$ con $u, v \in \mathbb{Z}$ primos entre sí, se define $H(x) = \mathrm{Max}(|u|, |v|)$ y $h(x) = \log(H(x)) \geqslant 0$.

Sea $E/\mathbb{Q}$ una curva elíptica dada por una ecuación reducida

$$y^2 = x^3 + Ax + B \quad \text{con} \quad A, B \in \mathbb{Z}. \tag{7}$$

**Definición 3** (Altura ingenua)**.** La *altura sobre $E$ relativa a la ecuación* (7) es la función

$$\begin{array}{rccc} h: & E(\mathbb{Q}) & \longrightarrow & \mathbb{R} \\ & O \neq P = (x, y) & \mapsto & h(x) \\ & O & \mapsto & 0. \end{array}$$

La altura de una curva elíptica $E/\mathbb{Q}$ tiene las siguientes propiedades:

**Proposición 9.** *1. Para toda constante $C$, el conjunto $\{P \in E(\mathbb{Q}); h(P) \leqslant C\}$ es finito.*

*2. Sea $P_0 \in E(\mathbb{Q})$. Existe $C$ tal que para todo $P \in E(\mathbb{Q})$,*

$$h(P + P_0) \leqslant 2h(P) + C.$$

*La constante $C$ depende de $P_0$ como lo precisa el ejercicio 11 en la subsección 3.1.*

3. *Existe $C'$ tal que para todo $P \in E(\mathbb{Q})$, $h([2]P) \geqslant 4h(P) - C'$.*

*Prueba.* Ver Silverman (1986, Cap. VIII. Lemma 4.2) □

**Descenso – Fin de la demostración del Teorema de Mordell.** Por el teorema 12, el grupo abeliano $E(\mathbb{Q})/2E(\mathbb{Q})$ es finito. Sea $\{Q_1, \ldots, Q_s\}$ un sistema finito de representantes de $E(\mathbb{Q})/2E(\mathbb{Q})$ en $E(\mathbb{Q})$.

Consideremos $P \in E(\mathbb{Q})$. Existen $i_1 \in \{1, \ldots, s\}$ y $P_1 \in E(\mathbb{Q})$ tales que $P = Q_{i_1} + 2P_1$. Cómo $P_1 \in E(\mathbb{Q})$ podemos repetir el proceso, construyendo por recurrencia dos sucesiones $(i_j)_{j\geqslant 1} \subset \{1, \ldots, s\}^{\mathbb{Z}_{\geqslant 0}}$ y $(P_j)_{j\geqslant 1} \subset E(\mathbb{Q})^{\mathbb{Z}_{\geqslant 0}}$ tales que,

$$P_j = Q_{i_j} + 2P_{j+1} \quad (j \geqslant 1).$$

De la proposición 9 (9.3. y 9.2.) con $P_0 = -Q_{j+1}$, deducimos que para todo $j \geqslant 1$, existen constantes $C_j$ y $C'$ tales que

$$\begin{aligned} h(P_{j+1}) &\leqslant \frac{1}{4}\,(2h(P_j) + C_j + C') \\ &\leqslant \frac{1}{4}\,(2h(P_j) + C + C') \end{aligned}$$

con $C = \operatorname{Max}(C_j)$.

Deducimos, por inducción, que para todo $j \geqslant 1$,

$$h(P_j) \leqslant \frac{1}{2^j} h(P) + \frac{C + C'}{2}.$$

De esta forma, existen $C''$ y $N$ tales que para todo $j > N$, se tiene $h(P_j) \leqslant C''$. Por lo tanto,

$$P = Q_{i_1} + 2Q_{i_2} + \cdots + 2^N Q_{i_N} + 2^{N+1} P_{N+1}$$

está en el grupo generado por $Q_1, \ldots, Q_s$ y en el conjunto $\{Q \in E(\mathbb{Q}); h(Q) \leqslant C''\}$, que es finito por la Proposición 9, 1. Lo anterior permite terminar la demostración. □

# 3 Funciones zeta y $L$ clásicas

**3.1 Generadores del grupo de Mordell–Weil.** Sabemos que el grupo de Mordell–Weil de una curva elíptica $E/\mathbb{Q}$ es de la forma

$$E(\mathbb{Q}) = E(\mathbb{Q})_{\text{tor}} \oplus \mathbb{Z}\, P_1 \oplus \cdots \oplus \mathbb{Z}\, P_r.$$

Hemos visto que el subgrupo de los puntos de torsión es fácil de calcular y además, es bien comprendido desde el punto de vista teórico. La parte de orden infinito de $E(\mathbb{Q})$ es mucho más extraña. Antes de hablar de funciones zeta y $L$, queremos motivar la introducción de estas funciones por la búsqueda de respuestas a las preguntas siguientes:

1. ¿Existe algún método para calcular el rango $r = \text{rango}\, E(\mathbb{Q})$?
2. ¿Hay alguna interpretación del rango $r = \text{rango}\, E(\mathbb{Q})$?
3. ¿Existe alguna cota superior del tamaño de generadores del grupo de Mordell–Weil?

Confesamos humildemente que: ¡nadie sabe responder de manera matemáticamente completa! Es decir, si bien sabemos calcular el grupo de Mordell–Weil para muchos ejemplos, siempre es necesario un poquito de suerte …Existe una vía – todavía conjetural y quizás no requiere de tanta fortuna – que involucra sin embargo la introducción de objetos analíticos y que vamos a desarrollar en esta segunda parte del curso.

**Ejercicio 8.** Consideramos la curva elíptica de ecuación de Weierstrass

$$y^2 + y = x^3 - x.$$

1. Demostrar que la curva tiene buena reducción en todo $p$, salvo en $p = 37$.
2. Demostrar que $E(\mathbb{F}_2) \cong \mathbb{Z}/5\mathbb{Z}$ y $E(\mathbb{F}_3) \cong \mathbb{Z}/7\mathbb{Z}$.
3. Concluir que $E(\mathbb{Q})_{\text{tor}} = \{0_E\}$ y que el punto $P := (0, 0)$ es de orden infinito.

(Nota: de hecho, el rango de esta curva es 1, aunque eso es más difícil de demostrar. Además el punto $P$ es un generador del grupo $E(\mathbb{Q})$, es decir, $E(\mathbb{Q}) = \mathbb{Z}.P$ (ver las indicaciones para el ejercicio 1 al final de estas notas.)

**Ejercicio 9.** La cúbica de Fermat, dada en coordenadas proyectivas por $X^3 + Y^3 + Z^3 = 0$ es una curva elíptica una vez que se escoge el origen, por ejemplo $0_E = (0, 1, -1)$.

1. Mostrar que los tres puntos $0_E$, $P = (1, 0, -1)$ y $Q = (1, -1, 0)$ forman un subgrupo de $E(\mathbb{Q})$ isomorfo a $\mathbb{Z}/3\mathbb{Z}$.
2. (Teorema de Fermat para $n = 3$) Mostrar que $E(\mathbb{Q}) = \{0_E, P, Q\}$.
3. Escribir la curva de Fermat en forma de Weierstrass minimal.

   [Indicación: poniendo $X = 3x/y$, con $Y = (y - 9)/y$, verificar que $y^2 - 9y = x^3 - 27$.]
4. Verificar en ambos modelos que la curva tiene buena reducción fuera de $p = 3$.

Comencemos por dar una formulación más intrínseca de las preguntas. La altura de puntos racionales (vistas en la subsección 2.6.2) verificam dos propriedades importantes:

$$-c_1 \leqslant h(2P) - 4h(P) \leqslant c_1 \tag{8}$$

$$-c_2 \leqslant h(P + Q) + h(P - Q) - 2h(P) - 2h(Q) \leqslant c_2 \tag{9}$$

Estas propriedades, mas fuertes que en la proposición 9, se pueden demonstrar como sugerido en el ejercicio siguiente.

**Ejercicio 10.** Demostrar las fórmulas precedentes utilizando las fórmulas geométricas para la curva elíptica de ecuación $y^2 = x^3 + ax + b$ que también podrían ser demostradas como ejercicio. (Ver la formula (3) en la subsección 2.23 y Silverman (1986, Cap. III.2.3)):

$$x([2]P) = \frac{x(P)^4 - 2ax(P)^2 - 8bx(P) + a^2}{4x^3 + b_2x^2 + 2b_4x + b_6}$$

$$x(P+Q) + x(P-Q) = \frac{2(x(P) + x(Q))(a + x(P)x(Q)) + 4b}{(x(P) - x(Q)))^2}$$

$$x(P+Q)x(P-Q) = \frac{(x(P)x(Q) - a)^2 - 4b(x(P) + x(Q))}{(x(P) - x(Q)))^2}$$

y las desigualdades aritméticas, para $\alpha, \beta \in \mathbb{Q}$:

$$\frac{1}{2}H(\alpha)H(\beta) \leqslant H(1, \alpha+\beta, \alpha\beta) \leqslant 2H(\alpha)H(\beta)$$

[Para $r, s \in \mathbb{Q}$, se escribe $(1, r, s) \sim (a, b, c)$ en $\mathbb{P}^2$ con $a, b, c$ enteros coprimos y se define $H(1, r, s) = \max(|a|, |b|, |c|)$.]

Utilizaremos el lema elemental siguiente:

**Lema 3.** *(Tate) Sea $S$ un conjunto, $\alpha > 1$ y dos aplicaciones $h : S \to \mathbb{R}$ y $\phi : S \to S$ tales que, para todo $x \in S$, tenemos $|h(\phi(x)) - \alpha h(x)| \leqslant c_1$, entonces, la sucesión $\alpha^{-n}h(\phi^n(x))$ es convergente y la función definida por*

$$\hat{h}(x) := \lim_{n\to\infty} \frac{h(\phi^n(x))}{\alpha^n},$$

*cumple las propiedades*

1. $|\hat{h}(x) - h(x)| \leqslant c_1/(\alpha - 1)$;
2. $\hat{h}(\phi(x)) = \alpha\hat{h}(x)$.

*Demostración.* Empezamos por verificar que $u_n := \alpha^{-n}h(\phi^n(x))$ es una sucesión de Cauchy. En efecto, como $-c_1 \leqslant h(\phi^n(x)) - \alpha h(\phi^{n-1}(x)) \leqslant c_1$, multiplicando por $\alpha^{-n}$ y sumando las desigualdades, obtenemos

$$-c_1\left(\frac{1}{\alpha^n} + \cdots \frac{1}{\alpha^{m+1}}\right) \leqslant u_n - u_m \leqslant c_1\left(\frac{1}{\alpha^n} + \cdots \frac{1}{\alpha^{m+1}}\right)$$

Lo cual prueba que $u_n$ es una sucesión de Cauchy. Haciendo tender $n$ hacia el infinito obtenemos

$$-\frac{c_1}{\alpha^m(\alpha-1)} \leqslant \hat{h}(x) - \alpha^{-m}h(\phi^m(x)) \leqslant \frac{c_1}{\alpha^m(\alpha-1)}$$

y en particular que $|\hat{h}(x) - h(x)|$ es acotada por $c_1/(\alpha - 1)$. Finalmente

$$\hat{h}(\phi(x)) = \lim_{n\to\infty} \frac{h(\phi^n(\phi(x)))}{\alpha^n} = \alpha \lim_{n\to\infty} \frac{h(\phi^{n+1}(x))}{\alpha^{n+1}} = \alpha\hat{h}(x).$$

□

Juntando las fórmulas (8) y (9) podemos introducir la altura de Néron–Tate y el regulador.

**Definición 4.** La altura de Néron–Tate (o altura canónica) de un punto $P \in E(\mathbb{Q})$ es la cantidad

$$\hat{h}(P) := \lim_{n\to\infty} \frac{h([2^n](P))}{4^n}.$$

**Teorema 13.** *La función "altura de Néron–Tate" $\hat{h} : E(\mathbb{Q}) \to \mathbb{R}$ es una forma cuadrática definida positiva; la diferencia entre la altura de Néron–Tate y la altura "ingenua" es una función acotada.*

*Demostración.* Aplicando el lema 3 de Tate a la función altura $h : E(\mathbb{Q}) \to \mathbb{R}$ y a la aplicación $\phi = [2]$, escogiendo $\alpha = 4$ obtenemos que $\hat{h} - h$ es acotada y que $\hat{h}(2P) = 4\hat{h}(P)$ y más generalmente que $\hat{h}(2^m P) = 4^m \hat{h}(P)$. Aplicando la fórmula (9) a los puntos $2^m P$ y $2^m Q$ obtenemos

$$-\frac{c_2}{4^m} \leqslant \frac{h(2^m(P+Q)) + h(2^m(P-Q)) - 2h(2^m P) - 2h(2^m Q)}{4^m} \leqslant \frac{c_2}{4^m}.$$

Tomando el límite se obtiene la ley del paralelogramo siguiente

$$\hat{h}(P+Q) + \hat{h}(P-Q) = 2\hat{h}(P) + 2\hat{h}(Q)$$

la cual sabemos que caracteriza a una forma cuadrática. Decir que $\hat{h}$ es definida positiva significa que, efectuando el producto tensorial con $\mathbb{R}$, la forma cuadrática $\hat{h}_{\mathbb{R}} : E(\mathbb{Q}) \otimes \mathbb{R} \to \mathbb{R}$ es definida positiva en el sentido clásico. Eso es verdad porque $\hat{h}_{\mathbb{R}}$ es positiva y verifica que el conjunto de los puntos $P \in E(\mathbb{Q})$ tales que $\hat{h}(P)$ esta acotada. Ver el ejercicio seguiente. □

**Ejercicio 11.** Estudiamos algunas propiedades de las alturas.

1. Demostrar que un punto $P \in E(\mathbb{Q})$ verifica $\hat{h}(P) = 0$ si y sólo si $P$ es de torsión.
2. Demostrar la versión siguiente del descenso: sean $Q_1, \ldots, Q_s$ representantes de $E(\mathbb{Q})/2E(\mathbb{Q})$ y $c_E := \max \hat{h}(Q_i)$ entonces el conjunto finito

   $$\{P \in E(\mathbb{Q}) \mid \hat{h}(P) \leqslant c_E\}$$

   es un conjunto de generadores del grupo de Mordell–Weil $E(\mathbb{Q})$.
3. Observando que $\hat{h}(P + P_0) \leqslant 2\hat{h}(P) + 2\hat{h}(P_0)$, concluir que existe una constante $c_0$ (independiente de $P$ y $P_0$) tal que

   $$h(P + P_0) \leqslant 2h(P) + 2h(P_0) + c_0.$$

Denotamos $< P, Q >:= \frac{1}{2}\left(\hat{h}(P+Q) - \hat{h}(P) - \hat{h}(Q)\right)$ al producto escalar asociado.

**Definición 5.** Sea $\{P_1, \ldots, P_r\}$ una base de $E(\mathbb{Q})$ sobre $\mathbb{Z}$, módulo la torsión. El regulador de $E/\mathbb{Q}$ es definido por

$$\mathrm{Reg}\,(E/\mathbb{Q}) := \det\,(< P_i, P_j >)_{1\leqslant i,j\leqslant r}\,.$$

El interés de la definición anterior viene del teorema siguiente debido a Minkowski y Hermite.

**Teorema 14.** *Sea $F \cong \mathbb{R}^n$ un espacio euclidiano con norma $||\cdot||$ y $\Lambda \cong \mathbb{Z}^n$ un retículo en $F$. Sea $P_1, \ldots, P_n$ una base de $\Lambda$ y $\mathrm{Reg}\,(\Lambda) := \det(< P_i, P_j >)$, entonces existe $Q_1, \ldots, Q_n$ una base de $\Lambda$ tal que para una constatnte $c_n$,*

$$||Q_1|| \ldots ||Q_n|| \leqslant c_n \left(\mathrm{Reg}\,(\Lambda)\right)^{1/2} \tag{10}$$

Observamos que la desigualdad de Hadamard nos proporciona, para cada base de $\Lambda$, que $\mathrm{Reg}\,(\Lambda)^{1/2} \leqslant ||Q_1|| \ldots ||Q_n||$.

Aplicando el teorema de Hermite–Minkowski al retículo $E(\mathbb{Q})/E(\mathbb{Q})_{\mathrm{tor}}$ en el espacio euclidiano $F = E(\mathbb{Q}) \otimes \mathbb{R}$, provisto de la forma cuadrática $\hat{h}$, obtenemos el corolario siguiente:

**Corolario 2.** *Sea $E/\mathbb{Q}$ una curva elíptica. Existe una base $Q_1, \ldots, Q_r$ (de la parte infinita) del grupo de Mordell–Weil tal que*

$$\hat{h}(Q_1) \ldots \hat{h}(Q_r) \leqslant c'_r \mathrm{Reg}\,(E/\mathbb{Q}). \tag{11}$$

Observamos que, para obtener una cota para las alturas de generadores (minimales), basta tener una cota *superior* para el regulador, y también una cota *inferior* para la altura minimal de un punto de orden infinito. O sea, si denotamos

$$m_E := \min_{Q \notin E(\mathbb{Q})_{\mathrm{tor}}} \hat{h}(Q),$$

y ordenamos en el orden creciente $\hat{h}(Q_1) \leqslant \cdots \leqslant \hat{h}(Q_r)$ obtenemos una cota

$$\hat{h}(Q_i) \leqslant \left( \frac{c'_r \mathrm{Reg}\,(E/\mathbb{Q})}{m_E^{i-1}} \right)^{\frac{1}{r-i+1}}.$$

Es útil observar que existen resultados teóricos sobre $m_E$ (por ejemplo Hindry y Silverman (1988)) pero también que, para un ejemplo dado, es muy fácil escribir una cota inferior para $m_E$. El problema fundamental es entonces el de encontrar una cota superior para $\mathrm{Reg}\,(E/\mathbb{Q})$.

Terminamos esta sección con una breve discusión sobre los parámetros con respeto a los cuales se quiere acotar $\mathrm{Reg}\,(E/\mathbb{Q})$. Escogemos la noción "ingenua" de altura de una curva elíptica siguiente.

**Definición 6.** Sea $E/\mathbb{Q}$ una curva elíptica, dada por una ecuación de Weierstrass minimal $y^2 = x^3 + Ax + B$ (es decir, $A, B \in \mathbb{Z}$ y para cada primo $p$, o bien $p^4$ no divide a $A$ o bien $p^6$ no divide a $B$). Definimos la altura de $E$ como

$$H(E/\mathbb{Q}) := \max \left\{ |A|^{1/4}, |B|^{1/6} \right\} \tag{12}$$

Denotamos también $h(E/\mathbb{Q}) := \log H(E/\mathbb{Q})$. Por ejemplo tenemos la desigualdad $\Delta_E \leqslant c_1 H(E)^{12}$ (Ejer.)

Existe una noción más sofisticada e intrínseca de altura llamada altura de Faltings, pero las dos son comparables en el sentido que hay una desigualdad de la forma $|h_{\mathrm{Faltings}}(E) - h(E)| \leqslant c_1 \log h(E) + c_2$.

Vamos a ver que varias conjeturas de naturaleza analítica implican una cota del tipo siguiente (ver Hindry (2007) y Lang (1983)).

**Conjetura 1.** Para todo $\epsilon > 0$, existe una constante $C_\epsilon > 0$ tal que para toda curva elíptica $E/\mathbb{Q}$ tenemos

$$\operatorname{Reg}(E/\mathbb{Q}) \leqslant C_\epsilon H(E/\mathbb{Q})^{1+\epsilon} \tag{13}$$

**3.2 La función zeta de Riemann.** La función zeta de Riemann, que tal vez debería llamarse función zeta de Euler–Riemann[1], es definida por una serie (del tipo serie de Dirichlet) y producto (del tipo producto de Euler), ambos convergentes para $\Re(s) > 1$, como sigue:

$$\zeta(s) := \sum_{n=1}^{\infty} \frac{1}{n^s} = \prod_{p\in\mathcal{P}} \left(1 - \frac{1}{p^s}\right)^{-1}. \tag{14}$$

La fórmula de producto anterior (debida a Euler) es una versión analítica de la unicidad de la descomposición en factores primos de los enteros.

**Continuación analítica y ecuación funcional.**

No es difícil extender la función al semiplano $\Re(s) > 0$, por ejemplo vía la fórmula (Ejer.)

$$\zeta(s) = s\int_0^\infty [t]t^{-s-1}dt = \sum_{n=1}^{\infty} n^{-s} = \frac{1}{s-1} + 1 + s\int_1^\infty ([t]-t)t^{-s-1}dt, \tag{15}$$

que muestra también que la función zeta tiene un polo simple en $s = 1$ con residuo igual a 1.

Disponemos de más. Recordemos que la función Gamma de Euler es definida para $\Re(s) > 0$ por la integral

$$\Gamma(s) := \int_0^\infty e^{-t}t^{s-1}dt,$$

y luego extendida al plano complejo (con polos simples en los enteros negativos) gracias a la ecuación $\Gamma(s+1) = s\Gamma(s)$ (Ejer.)

**Teorema 15.** *(Ecuación funcional) La función $\zeta(s)$ puede ser extendida en una función sobre todo el plano complejo, holomorfa, salvo un polo simple en $s = 1$ con residuo 1. Además satisface la ecuación funcional siguiente. Sea $\xi(s) := \pi^{-s/2}\Gamma(s/2)\zeta(s)$, entonces, fuera de 0 y 1, la función $\xi(s)$ es acotada en toda banda vertical del plano complejo y verifica:*

$$\xi(s) = \xi(1-s). \tag{16}$$

**Corolario 3.** *La función $\zeta(s)$ no se anula en el semiplano $\Re(s) > 1$; en el semiplano $\Re(s) < 0$ solamente se anula en los enteros negativos pares no nulos $-2, -4, -6, \ldots$. Todos los otros ceros están en la banda crítica $0 \leqslant \Re(s) \leqslant 1$.*

[1] Los valores reales de $\zeta(s)$ han sido estudiados por Euler (1707–1783), luego Dirichlet (1805–1859) empezó el uso de la variable compleja en ese contexto para demostrar el teorema de la progresión aritmética y Riemann (1826–1866) publicó en 1859 su famoso artículo sobre la repartición de los números primos.

*Demostración.* Esbozo de demostración. Utilizamos análisis armónico o análisis de Fourier. Sea $\hat{f}(y) := \int_{\mathbb{R}} f(x)\exp(2\pi i x y)dx$, la transformada de Fourier de una función integrable $f$. Tenemos la fórmula de Poisson:

$$\sum_{n\in\mathbb{Z}} f(n) = \sum_{n\in\mathbb{Z}} \hat{f}(n). \tag{17}$$

Definimos $\theta(u) := \sum_{n\in\mathbb{Z}} \exp\left(-\pi u n^2\right)$. Después, aplicando la fórmula de Poisson a $f(x) = \exp(-\pi u x^2)$ cuya transformada de Fourier es $\hat{f}(y) = \exp(-\pi y^2/u)/\sqrt{u}$ (Ejer.), obtenemos una ecuación funcional para la función theta (cuando veamos la noción de formas modulares, podremos traducirlo como el hecho de que $\theta$ es modular de peso $1/2$):

$$\theta(1/u) = \sqrt{u}\,\theta(u). \tag{18}$$

Podemos hacer el cálculo siguiente, cambiando la variable $u$ por $t = \pi n^2 u$ (este cálculo siendo válido, inicialmente, para $\Re(s) > 1$).

$$\begin{aligned}\xi(s) = \pi^{-s/2}\Gamma(s/2)\zeta(s) &= \sum_{n\geqslant 1}\int_0^\infty e^{-t}t^{s/2}\pi^{-s/2}n^{-s}\frac{dt}{t}\\ &= \int_0^\infty\left(\sum_{n\geqslant 1}\exp(-\pi u n^2)\right)u^{s/2}\frac{du}{u} = \int_0^\infty \tilde{\theta}(u)\frac{u^{s/2}du}{u}\end{aligned}$$

donde

$$\tilde{\theta}(u) := \sum_{n\geqslant 1}\exp\left(-\pi u n^2\right) = \frac{\theta(u)-1}{2}.$$

Observemos que $\tilde{\theta}(u) = O(\exp(-\pi u))$, cuando $u$ crece al infinito. Insertando (18) obtenemos que

$$\tilde{\theta}\left(\frac{1}{u}\right) = \sqrt{u}\,\tilde{\theta}(u) + \frac{1}{2}\left(\sqrt{u}-1\right). \tag{19}$$

Cómo $\int_1^\infty t^{-s}dt = 1/(s-1)$, y utilizando (19), obtenemos que

$$\begin{aligned}\xi(s) &= \int_0^1 \tilde{\theta}(u)\frac{u^{s/2}du}{u} + \int_1^\infty \tilde{\theta}(u)\frac{u^{s/2}du}{u}\\ &= \int_1^\infty \tilde{\theta}(1/u)\frac{u^{-s/2}du}{u} + \int_1^\infty \tilde{\theta}(u)\frac{u^{s/2}du}{u}\\ &= \int_1^\infty\left\{\sqrt{u}\tilde{\theta}(u) + \frac{1}{2}\left(\sqrt{u}-1\right)\right\}\frac{u^{-s/2}du}{u} + \int_1^\infty \tilde{\theta}(u)\frac{u^{s/2}du}{u}\\ &= \int_1^\infty \tilde{\theta}(u)\left\{u^{\frac{s}{2}} + u^{\frac{1-s}{2}}\right\}\frac{du}{u} + \frac{1}{s-1} - \frac{1}{s}.\end{aligned} \tag{20}$$

La última expresión es, a priori, válida para $\Re(s) > 1$, pero es fácil ver que, como $\tilde{\theta}(u) = O(\exp(-\pi u))$, la función definida por la integral es holomorfa sobre $\mathbb{C}$, que además es claramente simétrica con respeto a la transformación $s \mapsto 1-s$. Tenemos finalmente que la función definida por la integral es acotada en cada banda vertical. □

Los ceros en la banda crítica tienen dos simetrías, a saber $s \mapsto \bar{s}$ y $s \mapsto 1-s$, y esto tal vez sugiere la conjetura siguiente.

**Conjetura 2.** (Hipótesis de Riemann) Los ceros en la banda crítica de la función $\zeta(s)$ verifican $\Re(s) = \frac{1}{2}$.

De manera equivalente, esta conjetura afirma que la función $\zeta(s)$ no se anula para $\Re(s) > \frac{1}{2}$.

Esencialmente, solamente conocemos una versión más débil de esta conjetura – pero que es suficiente para demostrar el teorema de los números primos – demostrado por Hadamard y de la Vallée Poussin (ver por ejemplo Iwaniec y Kowalski (2004)):

**Teorema 16.** *(Hadamard – de la Vallée Poussin) La función $\zeta(s)$ no se anula en la recta $\Re(s) = 1$.*

### 3.3 Generalizaciones de la función de Riemann.

**3.3.1 La función zeta de Dedekind.** La función zeta de de Dedekind de un cuerpo de números $K$ es también definida por una serie y un producto, ambos convergentes para $\Re(s) > 1$:

$$\zeta_K(s) := \sum_{I}^{\infty} \frac{1}{N(I)^s} = \prod_{\mathfrak{p}} \left(1 - \frac{1}{N(\mathfrak{p})^s}\right)^{-1} \tag{21}$$

donde esta vez, $I$ recorre los *ideales* no nulos del anillo de los enteros algebraicos $\mathcal{O}_K$, mientras que $\mathfrak{p}$ recorre los ideales *primos* no nulos de dicho anillo y donde $N(I)$ denota el cardinal del cociente $\mathcal{O}_K/I$; la fórmula del producto para la función zeta de Dedekind es una expresión analítica del teorema sobre la unicidad de la descomposición de un ideal en producto de ideales primos en un anillo de Dedekind.

Para enunciar la ecuación funcional, recordamos que un cuerpo de números tiene un *discriminante*, cuyo valor absoluto es denotado $\Delta_K$, y tiene $r_1$ inmersiones reales $K \hookrightarrow \mathbb{R}$ y $r_2$ pares de inmersiones complejas conjugadas $K \hookrightarrow \mathbb{C}$, de manera que $r_1 + 2r_2 = [K : \mathbb{Q}]$.

**Definición 7.** Definimos los factores Gamma real y complejo como:

$$\Gamma_{\mathbb{R}}(s) := \pi^{-s/2}\Gamma(s/2) \qquad \text{y} \qquad \Gamma_{\mathbb{C}}(s) := (2\pi)^{-s}\Gamma(s).$$

**Teorema 17.** *(Ecuación funcional) La función $\zeta_K(s)$ puede ser extendida en una función definida sobre todo el plano complejo y holomorfa salvo en un polo simple en $s = 1$ con residuo $\lambda(K)$. Además satisface la ecuación funcional siguiente. Sea*

$$\xi_K(s) := \Delta_K^{s/2}\Gamma_{\mathbb{R}}(s)^{r_1}\Gamma_{\mathbb{C}}(s)^{r_2}\zeta_K(s)$$

*entonces, fuera de los puntos* 0 *y* 1 *del plano complejo, la función* $\xi(s)$ *es acotada en toda banda vertical y verifica:*

(22) $$\xi_K(s) = \xi_K(1-s).$$

El residuo $\lambda(K)$ es dado por una fórmula magnífica que contiene los invariantes más importantes del cuerpo $K$.

1. El número de clases $h_K$;
2. El regulador $R_K$ de las unidades (ver por ejemplo el libro de Lang (1970), Capitulo V, §1);
3. El número $w_K$ de raíces de la unidad.

El número de clases $h_K$ es el número de clases de ideales módulo los ideales principales de $\mathcal{O}_K$ o sea $h_K = |\text{Pic}(\mathcal{O}_K)|$. Por el teorema de las unidades de Dirichlet, el grupo de las unidades, es decir, el grupo $\mathcal{O}_K^\times$ de los elementos invertibles, es isomorfo a $\mathbb{Z}/w_K\mathbb{Z} \times \mathbb{Z}^r$ con $r = r_1 + r_2 - 1$; **Fórmula para el residuo de $\zeta_K(s)$ en $s = 1$.** (Ver Lang (ibíd.) Capitulo XIII, §3, Teorema 1.)

(23) $$\lambda(K) = \frac{h_K R_K}{\sqrt{\Delta_K}} \cdot \frac{2^{r_1}(2\pi)^{r_2}}{w_K}$$

Interpretando un ideal primo no nulo de $\mathcal{O}_K$ como un ideal maximal, no es difícil generalizar esta definición para la función zeta de un anillo $R$ de tipo finito sobre $\mathbb{Z}$; esta condición con el fin de que para todo ideal maximal $\mathfrak{m}$, el cociente $R/\mathfrak{m}$ sea finito en cardinal, el cual denotamos $\kappa(\mathfrak{m})$. Así, definimos:

(24) $$\zeta_R(s) := \prod_{\mathfrak{m}} \left(1 - \frac{1}{\kappa(\mathfrak{m})^s}\right)^{-1}.$$

Observamos que el producto sobre todos los ideales maximales es convergente para $\Re(s) \gg 1$.

**3.3.2 La función zeta de un esquema de tipo finito sobre $\mathbb{Z}$.** Sea $\mathcal{X}$ un esquema sobre $\mathbb{Z}$ (ver el "apéndice" al fin de este texto para una breve introducción al lenguaje de esquemas y, por ejemplo, el libro de Hartshorne (1977), para un curso completo). Denotamos $|\mathcal{X}|$ el conjunto de los puntos cerrados de $\mathcal{X}$. Por cada punto *cerrado* $x$ de $\mathcal{X}$, tenemos un anillo local $\mathcal{O}_{x,\mathcal{X}}$ con ideal maximal $\mathfrak{M}_{x,\mathcal{X}}$ y cuerpo residual $\kappa(x) := \mathcal{O}_{x,\mathcal{X}}/\mathfrak{M}_{x,\mathcal{X}}$. Denotamos $N(x)$ al cardinal de $\kappa(x)$. La función zeta de $\mathcal{X}$ es dada por el producto:

(25) $$\zeta_{\mathcal{X}}(s) := \prod_{x \in |\mathcal{X}|} (1 - N(x)^{-s})^{-1}$$

Si $R$ es un anillo de tipo finito sobre $\mathbb{Z}$ y $\mathcal{X}$ es el *espectro* de $R$ entonces $\zeta_R(s) = \zeta_{\mathcal{X}}(s)$

Formalmente $\zeta(\mathcal{X}_1 \sqcup \mathcal{X}_2, s) = \zeta(\mathcal{X}_1, s)\zeta(\mathcal{X}_2, s)$. De esta manera, si se descompone el conjunto de puntos cerrados según la característica residual, se obtiene el producto de Euler siguiente, donde $\mathcal{X}_p$ denota la fibra en $p$:

(26) $$\zeta(\mathcal{X}, s) = \prod_p \zeta(\mathcal{X}_p, s).$$

Observemos que $\mathcal{X}_p$ es una variedad (no necesariamente irreducible) sobre el cuerpo finito $\mathbb{F}_p$. Es claro ahora que tenemos que estudiar la función zeta de tales variedades. Será el contenido del próximo párrafo.

**Ejemplo 13.** i) Si $\mathcal{X} = \operatorname{spec}(\mathbb{Z})$ (respectivamente $\mathcal{X} = \operatorname{spec}(\mathcal{O}_K)$), entonces $\zeta(\mathcal{X}, s)$ es simplemente la función zeta de Riemann (respectivamente la función zeta de Dedekind para el cuerpo $K$).

ii) Si $\mathcal{X} = \mathbb{A}^1_{\mathbb{Z}} = \operatorname{spec}(\mathbb{Z}[T])$, podemos identificar los puntos cerrados (ideales maximales de $\mathbb{Z}[T]$) de característica $p$ con los polinomios unitarios irreducibles de $\mathbb{F}_p[T]$, cuyo conjunto denotamos $\mathrm{Irr}_p$; escribimos también $M_p$ para el conjunto de los polinomios unitarios con coeficientes en $\mathbb{F}_p$. Con la notación anterior, vemos que el cálculo siguiente es bastante simple:

$$\begin{aligned}
\zeta(\mathbb{A}^1_{\mathbb{Z}}, s) &= \prod_p \prod_{Q \in \mathrm{Irr}_p} (1 - p^{-s \deg Q})^{-1} \\
&= \prod_p \prod_{Q \in \mathrm{Irr}_p} \sum_{m=0}^{\infty} p^{-sm \deg Q} \\
&= \prod_p \sum_{P \in M_p} p^{-s \deg P} \\
&= \prod_p \sum_{d=0}^{\infty} p^{d-ds} \\
&= \prod_p (1 - p^{1-s})^{-1} \\
&= \zeta(s-1).
\end{aligned}$$

iii) Dividiendo (mediante una partición) los puntos cerrados de $\mathbb{P}^1_{\mathbb{Z}}$ en $\mathbb{A}^1_{\mathbb{Z}} \sqcup \mathbb{A}^0_{\mathbb{Z}}$ obtenemos la fórmula

$$\zeta(\mathbb{P}^1_{\mathbb{Z}}, s) = \zeta(s)\zeta(s-1).$$

**Ejercicio 12.** Utilizando la ecuación funcional de la función zeta de Riemann encontrar una relación de la forma siguiente (con $G(s)$ un producto de valores de la función Gamma):

$$G(2-s)\zeta(\mathbb{P}^1_{\mathbb{Z}}, 2-s) = G(s)\zeta(\mathbb{P}^1_{\mathbb{Z}}, s).$$

### 3.3.3 La función zeta de Weil de una variedad sobre un cuerpo finito.

Sea $X$ una variedad sobre $\mathbb{F}_p$. A fin de tener propiedades más elegantes, vamos a suponer que $X$ es además lisa y proyectiva (aunque esta suposición no es necesaria inmediatamente). Para un punto cerrado $x$ escribimos $d_x = [\mathbb{F}_p(x) : \mathbb{F}_p]$. Podemos calcular

$$\log \zeta_X(s) = \sum_{x,m} \frac{N(x)^{-ms}}{m} = \sum_{x,m} \frac{p^{-d_x ms}}{m} = \sum_{n \geqslant 1} p^{-ns} \left( \sum_{m d_x = n} \frac{1}{m} \right).$$

Escribimos a la última suma como $\sum_{n \geqslant 1} u_n p^{-ns}$ y continuamos el cálculo observando que, para un punto cerrado $x \in |X|$, tener grado residual $d_x$ que divide a $n$ es equivalente a corresponder a una clase de conjugación de puntos en $X(\mathbb{F}_{p^n})$, luego:

$$u_n = \sum_{x\in|X|,\, d_x=\frac{n}{m}} \frac{1}{m} = \frac{1}{n}\sum_{x\in|X|,\, d_x|n} d_x = \frac{1}{n}\#X(\mathbb{F}_{p^n}).$$

Lo anterior sugiere la siguiente definición:

**Definición 8.** Sea $X/\mathbb{F}_p$ una variedad sobre un cuerpo finito, la *función zeta de Weil* es la serie formal

$$Z(X/\mathbb{F}_p, T) = \exp\left(\sum_{m=1}^{\infty} \frac{|X(\mathbb{F}_{p^m})|}{m} T^m\right). \tag{27}$$

El vínculo con las funciones zeta precedentes se exhibe en los siguientes resultados Ejer.

**Teorema 18.** *Sea $X/\mathbb{F}_p$ una variedad sobre un cuerpo finito, entonces:*

$$\zeta(X/\mathbb{F}_p, s) = Z(X/\mathbb{F}_p, p^{-s}). \tag{28}$$

**Ejemplo 14.** Calculamos un ejemplo simple, cuando $X = \mathbb{P}^n$. En este caso

$$\#X(\mathbb{F}_{p^m}) = \frac{p^{m(n+1)}-1}{p^m-1} = p^{mn} + p^{m(n-1)} + \cdots + p^m + 1,$$

luego

$$\sum_{m=1}^{\infty} \frac{\#X(\mathbb{F}_{p^m})}{m} T^m = \sum_{j=0}^{n}\sum_{m=1}^{\infty} \frac{p^{mj}}{m} T^m = -\sum_{j=0}^{n} \log(1-p^j T),$$

y entonces:

$$Z(\mathbb{P}^n, T) = \frac{1}{(1-T)(1-pT)\dots(1-p^{n-1}T)(1-p^nT)}. \tag{29}$$

Por otro lado una verificación simple muestra que $Z$ cumple la ecuación funcional:

$$Z(\mathbb{P}^n, T) = (-1)^{n+1} p^{\frac{n(n+1)}{2}} T^{n+1} Z\left(\mathbb{P}^n, \frac{1}{p^n T}\right).$$

Las propiedades vistas en el ejemplo anterior se generalizan ampliamente. Las propiedades de la función zeta de Weil fueron conjeturadas en su forma general por Weil y demostradas por Grothendieck y Deligne.

**Teorema 19.** *(Conjeturas de Weil) Sea $X$ una variedad lisa y proyectiva $\mathbb{F}_p$ de dimensión $n$.*

1. *(Racionalidad) Existen polinomios $P_j(X,T) = \prod_{i=1}^{b_j}(1-\alpha_{j,i}T) \in \mathbb{Z}[T]$ para $j = 0, \dots, 2n$ tales que*

$$Z(X,T) = \frac{P_1(X,T)\dots P_{2n-1}(X,T)}{P_0(X,T)\dots P_{2n}(X,T)} = \prod_{j=0}^{2n} P_j(X,T)^{(-1)^{j+1}}. \tag{30}$$

*Además $b_0 = b_{2n} = 1$ y, $P_0(X,T) = 1-T$ y $P_{2n}(X,T) = 1-p^nT$.*

2. *(Ecuación funcional) Sea* $\chi(X) = \sum_{j=0}^{2n}(-1)^j b_j$ *la característica de Euler–Poincaré. Entonces*

$$Z(X,T) = \pm p^{\frac{n\chi(X)}{2}} T^{\chi(X)} Z\left(X, \frac{1}{p^n T}\right). \tag{31}$$

3. *(Hipótesis de Riemann) Los enteros algebraicos* $\alpha_{j,i}$ *satisfacen* $|\alpha_{j,i}| = p^{j/2}$.

4. *Los números* $b_j = b_j(X)$ *satisfacen una propiedad de continuidad en familias lisas; en particular, si* $X$ *es la reducción módulo un ideal primo de una variedad* $Y$ *definida sobre un cuerpo de números* $K$*, entonces* $b_j(X)$ *es igual al número de Betti de la variedad compleja* $Y \otimes_K \mathbb{C}$.

El siguiente es el ejemplo clave para este curso.

**Ejemplo 15.** Sea $E$ una curva elíptica sobre $\mathbb{F}_p$. tenemos que $b_0 = b_2 = 1$, $b_1 = 2$ y que los polinomios se escriben $P_0(T) = 1 - T$, $P_2(T) = 1 - pT$ y $P_1(T) = 1 - aT + pT^2$. Entonces, por el teorema 19 hay un entero algebraico $\alpha$ con $\alpha\bar{\alpha} = p$ tal que

$$Z(E,T) = \frac{(1-\alpha T)(1-\bar{\alpha}_p T)}{(1-T)(1-pT)}.$$

Esto es equivalente al teorema de Hasse (Teorema 6) que dice que

$$\#E(\mathbb{F}_{p^m}) = p + 1 - \alpha^m - \bar{\alpha}^m.$$

Las consideraciones precedentes permiten definir la función $\zeta$ (y después la función $L$) asociada a una curva elíptica $E/\mathbb{Q}$ como un producto de Euler, donde el factor para un $p$ de buena reducción será:

$$\zeta(E_p, s) = Z(E_p, p^{-s}) = \frac{(1-\alpha_p p^{-s})(1-\bar{\alpha} p^{-s})}{(1-p^{-s})(1-p^{1-s})} = \frac{1 - a_p p^{-s} + p^{1-2s}}{(1-p^{-s})(1-p^{1-s})}.$$

Observamos además que el teorema de Hasse (Teorema 6) nos dice que los ceros de $\zeta(E_p, s)$ verifican $\Re(s) = \frac{1}{2}$. Esta es la razón por la cual se llama a este resultado (y respectivamente a las generalizaciones posteriores por Weil y Deligne) *hipótesis de Riemann* para curvas elípticas (respectivamente para variedades lisas y proyectivas).

Sin embargo no es claro cómo deben ser los factores para $p$ con mala reducción. Siguiendo a Serre (1969/70), vamos a introducir las representaciones de Galois, para aclarar este punto. En efecto, las representaciones de Galois son una herramienta fundamental en la aritmética moderna; por ejemplo, juegan un papel fundamental en el trabajo de Wiles.

**3.4 La función $L$ asociada a una representación de Galois.** Uno de los objetos centrales de la geometría aritmética es el grupo de Galois absoluto:

$$G_{\mathbb{Q}} := \mathrm{Gal}\,(\bar{\mathbb{Q}}/\mathbb{Q}). \tag{32}$$

Es natural estudiar este enorme grupo profinito a través de sus representaciones. Vamos a considerar dos tipos de representaciones, primero las de coeficientes complejos (representaciones de Artin) y después las representaciones asociadas a curvas elípticas, con coeficientes en $\mathbb{Q}_\ell$, que son el prototipo de representaciones asociadas a la cohomología $\ell$-ádica de variedades algebraicas.

**Definición 9.** Sea $V$ un espacio vectorial de dimensión $n$ sobre un cuerpo $K$, una representación de Galois es un homomorfismo continuo (en estas notas, $K = \mathbb{C}$ o $\mathbb{Q}_\ell$ y la topología es natural) de grupos:

$$\rho : G_{\mathbb{Q}} \to \mathrm{GL}_n(K) = \mathrm{GL}\,(V).$$

**3.4.1 Representaciones de Artin y elementos de Frobenius.** Sea $K/\mathbb{Q}$ una extensión finita. Cada primo racional $p$ se descompone en $K$ como un producto de ideales primos $p\mathcal{O}_K = \mathfrak{P}_1^{e_1}\dots\mathfrak{P}_g^{e_g}$; el primo $p$ siendo *ramificado* en $K/\mathbb{Q}$ si para algún exponente $e_i \geqslant 2$. Si denotamos $f_i := [\mathcal{O}_K/\mathfrak{P}_i : \mathbb{F}_p]$ entonces tenemos que $\sum_{i=1}^{g} e_i f_i = [K : \mathbb{Q}]$.

Supongamos ahora que $K/\mathbb{Q}$ es una extensión de Galois con grupo de Galois $G$.

**Definición 10.** El *grupo de descomposición* de $\mathfrak{P}/p$ es el subgrupo

$$D(\mathfrak{P}/p) := \{\sigma \in G \mid \sigma(\mathfrak{P}) = \mathfrak{P}\}.$$

Cuando $\sigma$ está en el grupo de descomposición, se puede definir $\tilde{\sigma} : \mathcal{O}_K/\mathfrak{P} \to \mathcal{O}_K/\mathfrak{P}$ mediante el diagrama siguiente

$$\begin{array}{ccc} \mathcal{O}_K & \xrightarrow{\sigma} & \mathcal{O}_K \\ \downarrow & & \downarrow \\ \mathcal{O}_K/\mathfrak{P} & \xrightarrow{\tilde{\sigma}} & \mathcal{O}_K/\mathfrak{P}. \end{array}$$

**Definición 11.** El núcleo del mapa $\sigma \to \tilde{\sigma}$ del grupo de descomposición $D(\mathfrak{P}/p)$ hacia $\mathrm{Gal}\,((\mathcal{O}_K/\mathfrak{P})/\mathbb{F}_p)$ es llamado *grupo de inercia* de $\mathfrak{P}/p$ y es denotado $I(\mathfrak{P}/p)$.

Observemos que, si $\#D(\mathfrak{P}/p) = f$ y $\#I(\mathfrak{P}/p) = e$, entonces $\#G = [K : \mathbb{Q}] = efg$. Notamos que el grupo de inercia $I(\mathfrak{P}/p)$ es trivial cuando $\mathfrak{P}/p$ no es ramificado. Por otro lado, la aplicación $\sigma \to \tilde{\sigma}$ de $D(\mathfrak{P}/p)$ hacia $\mathrm{Gal}\,((\mathcal{O}_K/\mathfrak{P})/\mathbb{F}_p)$ es sobreyectiva y se sabe que su imagen es un grupo cíclico con un generador canónico dado por $x \mapsto x^p$.

**Definición 12.** El Frobenius en $\mathfrak{P}$ es el elemento $\mathrm{Frob}_{\mathfrak{P}}$ tal que para todo $x \in \mathcal{O}_K$ tenemos

$$\mathrm{Frob}_{\mathfrak{P}}(x) \equiv x^p \bmod \mathfrak{P}.$$

Observemos que, en el caso ramificado, el Frobenius es en verdad una clase lateral módulo el grupo de inercia. Además, si escogemos un otro primo $\mathfrak{P}'$, encima de $p$, digamos tal que $\mathfrak{P}' = \sigma(\mathfrak{P})$, entonces $\mathrm{Frob}_{\mathfrak{P}'} = \sigma \mathrm{Frob}_{\mathfrak{P}} \sigma^{-1}$, de modo que $\mathrm{Frob}_{\mathfrak{P}}$ depende solamente de $p$ salvo por una conjugación. Notaremos en este caso $\mathrm{Frob}_p$ a la correspondiente clase de conjugación.

*Observación* 3. Las observaciones siguientes son simples pero importantes y pueden aplicarse a los elementos de Frobenius y el subgrupo de inercia. Sea $\rho : G \to \mathrm{GL}\,(V)$ una representación. Si $f = hgh^{-1} \in G$, entonces los polinomios característicos de $\rho(f)$ y $\rho(g)$ actuando sobre $V$ son iguales. Si $H$ es un subgrupo de $G$, notamos

$$V^H := \{v \in V \mid \forall h \in H,\ \rho(h)(v) = v\},$$

al subespacio de los vectores fijos. Si $f \in gH$ y $g$ normaliza $H$ (es decir, tenemos $gHg^{-1} = H$) entonces $g$ y $f$ estabilizan $V^H$ y los polinomios característicos de $\rho(f)$ y $\rho(g)$ actuando sobre $V^H$ son también iguales.

Lo anterior permite, por ejemplo, la definir la función $L$ asociada a una representación de Artin $\rho$ como:

$$L(\rho, s) = \prod_p \det\left(1 - \rho(\mathrm{Frob}_p) p^{-s} \mid V^{I_p}\right)^{-1}. \tag{33}$$

*No disponemos del espacio ni del tiempo para desarrollar la teoría de tales funciones L (para una iniciación, ver Lang (1970)), solamente indicaremos que son uno de los objetos clave para el programa de Langlands y para entender los vínculos entre objetos geométricos – como variedades algebraicas – y analíticos – como las formas modulares y automorfas.*

**3.4.2 Representaciones de Galois asociadas a una curva elíptica.** Observamos que, por razones topológicas, una representación (que siempre se supone que es continua) $\rho : \mathrm{Gal}(\bar{\mathbb{Q}}/\mathbb{Q}) \to \mathrm{GL}_n(\mathbb{C})$ se factoriza a través de un grupo finito $G = \mathrm{Gal}(K/\mathbb{Q})$, pero que esta propiedad no es verdadera para las representaciones $\ell$-ádicas asociadas a una curva elíptica. Son este tipo de representaciones que consideraremos ahora.

Recordemos que el cuerpo de los números $\ell$-ádicos $\mathbb{Q}_\ell$ puede ser construido como la completación del cuerpo $\mathbb{Q}$ con respecto al valor absoluto $|x|_\ell := \ell^{-\mathrm{ord}_\ell x}$. El anillo de los enteros $\mathbb{Z}_\ell = \{x \in \mathbb{Q}_\ell \mid |x|_\ell \leqslant 1\}$ es entonces la completación de $\mathbb{Z}$. Alternativamente, se puede definir $\mathbb{Z}_\ell$ como el límite inverso de los grupos finitos $\mathbb{Z}/\ell^n\mathbb{Z}$ (bajo los homomorfismos evidentes):

$$\mathbb{Z}_\ell = \varprojlim_n \mathbb{Z}/\ell^n\mathbb{Z}.$$

Construimos también el cuerpo $\mathbb{Q}_\ell$ como el cuerpo de fracciones de $\mathbb{Z}_\ell$.

**Definición 13.** Sea $E$ una curva elíptica sobre un cuerpo $K$ de característica 0 o $p$ diferente de $\ell$. El núcleo de la multiplicación por $\ell^n$ es denotado por

$$E[\ell^n] := \mathrm{Ker}\left\{[\ell^n] : E(\bar{K}) \to E(\bar{K})\right\},$$

y es isomorfo (como grupo) a $(\mathbb{Z}/\ell^n\mathbb{Z})^2$. El *módulo de Tate* es el límite proyectivo

$$T_\ell(E) := \varprojlim E[\ell^n] \cong \left(\varprojlim \mathbb{Z}/\ell^n\mathbb{Z}\right)^2 \cong \mathbb{Z}_\ell^2.$$

El grupo de Galois $G_K := \mathrm{Gal}(\bar{K}/K)$ actúa sobre cada $E[\ell^n]$ y entonces sobre $T_\ell(E)$, definiendo así una representación:

$$\rho_{E,\ell} : G_K \to \mathrm{GL}(T_\ell(E)) \cong \mathrm{GL}_2(\mathbb{Z}_\ell) \subset \mathrm{GL}_2(\mathbb{Q}_\ell).$$

Sea ahora $E/\mathbb{Q}$ una curva elíptica. Tenemos, para todo $\ell$ una representación $\rho_\ell$ con coeficientes en $\mathbb{Q}_\ell$. Observemos que, para $p$ con buena reducción y $\ell \neq p$, la reducción módulo $p$ define un isomorfismo $E[\ell^m] \to \tilde{E}_p[\ell^m]$ y en consecuencia un isomorfismo $T_\ell(E) \to T_\ell(\tilde{E}_p)$. Así el polinomio característico de $\rho(\mathrm{Frob}_p)$ es igual al polinomio característico del endomorfismo de Frobenius módulo $p$. A priori, el polinomio característico de $\mathrm{Frob}_p$ es un polinomio con coeficientes en $\mathbb{Z}_\ell$, pero se puede deducir (lo cual es esencialmente la prueba del teorema de Hasse) que en realidad los coeficientes son enteros y además independientes de $\ell$. De hecho $\det \rho_{E,\ell}(\mathrm{Frob}_p) = p$ y $\mathrm{Tr}\, \rho_{E,\ell}(\mathrm{Frob}_p) = a_p = p + 1 - \#E(\mathbb{F}_p)$. Es decir, el polinomio característico mencionado es $1 - a_p T + pT^2$. Este resultado nos permite la definición siguiente, donde denotamos $V_\ell = V_\ell(E)$ al $\mathbb{Q}_\ell$-espacio vectorial $T_\ell(E) \otimes_{\mathbb{Z}_\ell} \mathbb{Q}_\ell$.

**Definición 14.** La función $L$ asociada a las representaciones $\rho_{E,\ell}$ es[2]:

$$L(\rho_E, s) = \prod_{\mathfrak{p}} \det \left(1 - \rho_{E,\ell}(\mathrm{Frob}_{\mathfrak{p}}) N(\mathfrak{p})^{-s} \mid V_\ell^{I_{\mathfrak{p}}}\right)^{-1}. \tag{34}$$

Cuando $p$ divide al discriminante, la curva $E$ módulo $p$ tiene o bien una punta, en cuyo caso $\dim V^I = 0$, y hacemos $a_p = 0$, o bien un nodo, en cuyo caso $\dim V^I = 1$. Cuando tenemos un nodo, las dos tangentes pueden ser racionales sobre $\mathbb{F}_p$ y $\rho(\mathrm{Frob}_p)$ actúa trivialmente y definamos $a_p = 1$, o bien las dos tangentes no son racionales sobre $\mathbb{F}_p$ y $\rho(\mathrm{Frob}_p)$ actúa como la multiplicación por $-1$, por lo que ponemos $a_p = -1$. Así tenemos la expresión explícita:

$$L(\rho_E, s) = \prod_{p \mid \Delta} (1 - a_p p^{-s})^{-1} \prod_{p \nmid \Delta} (1 - a_p p^{-s} + p^{1-2s})^{-1}. \tag{35}$$

## 3.5 La función $L$ de Hasse–Weil de una curva elíptica.

**3.5.1 Definición como producto de Euler.** Uno puede dar la definición más concreta posible de la función $L(E, s)$, obtenida a través de la consideración de representaciones de Galois, del siguiente modo:

**Definición 15.** Sea $E/\mathbb{Q}$ una curva elíptica con discriminante minimal $\Delta_E$. Para un primo $p$ que no divide a $\Delta$, definimos $a_p = a_p(E) = p + 1 - |E(\mathbb{F}_p)|$ y para un primo $p$ que divide a $\Delta_E$ ponemos[3]

$$a_p = \begin{cases} +1 & \text{si la reducción es multiplicativa y } \textit{desplegada} \\ -1 & \text{si la reducción es multiplicativa y } \textit{es no desplegada} \\ 0 & \text{si la reducción es aditiva} \end{cases}$$

Así podemos definir

$$L(E, s) := \sum_{n=1}^{\infty} \frac{a_n}{n^s} = \prod_{p \mid \Delta_E} (1 - a_p p^{-s})^{-1} \prod_{p \nmid \Delta_E} (1 - a_p p^{-s} + p^{1-2s})^{-1}. \tag{36}$$

La serie de Dirichlet y el producto de Euler son ambos convergentes para $\Re(s) > 3/2$. Esto se demuestra utilizando el teorema de Hasse que proporciona $|a_p| \leqslant 2\sqrt{p}$ y más generalmente que $|a_n| \leqslant \sigma(n)\sqrt{n}$, donde $\sigma(n)$ corresponde al número de divisores de $n$ Ejer.

[2] Serre habla de *sistema de representaciones compatibles*.

[3] Desplegada = *split* en inglés = *déployée* en francés.

La función $L(E, s)$ tiene propiedades similares a las otras funciones $\zeta$ y $L$ (definidas en secciones anteriores), pero estas son mucho más difíciles de establecer. En efecto el método de demostración requerido consiste en establecer un vínculo con otros objetos analíticos: las formas modulares.

**3.5.2 Función $L$ asociada a una forma modular.** *Esta parte es breve y recomendamos el curso de Harris–Miatello–Moreno–Pacetti–Tornaria para más detalles.*

Las formas modulares son funciones holomorfas sobre el semiplano de Poincaré $\mathcal{H} := \{\tau \in \mathbb{C} \mid \Im(\tau) > 0\}$. Ellas satisfacen la condición principal de que, para un número $k$, llamado *peso* de la forma modular y un subgrupo $\Gamma$ de

$$\mathrm{SL}\,(2,\mathbb{Z}) = \left\{\gamma = \begin{pmatrix} a & b \\ c & d \end{pmatrix}; a, b, c, d \in \mathbb{Z},\ ad - bc = 1\right\}$$

tenemos, para todo $\gamma \in \Gamma$ y para todo $\tau \in \mathcal{H}$ que

$$f\left(\frac{a\tau + b}{c\tau + d}\right) = (c\tau + d)^k f(\tau)\cdot \tag{37}$$

La segunda condición que ellas satisfacen corresponde a una condición de “holomorfía en el infinito”, la cual vamos a explicar solamente para el único subgrupo que vamos a considerar en este curso:

**Definición 16.** El grupo de congruencia de nivel $N$ es el subgrupo:

$$\Gamma_0(N) := \left\{\begin{pmatrix} a & b \\ c & d \end{pmatrix} \in \mathrm{SL}\,(2,\mathbb{Z}) \mid N \text{ divide } c\right\}\cdot \tag{38}$$

Observamos que la matriz $\begin{pmatrix} 1 & 1 \\ 0 & 1 \end{pmatrix}$ pertenece a $\Gamma_0(N)$, de manera que toda función holomorfa $f : \mathcal{H} \to \mathbb{C}$ que cumple (37) es periódica, es decir es tal que $f(\tau + 1) = f(\tau)$, ey en particular se puede escribir como una serie de Fourier:

$$f = \sum_{n\in\mathbb{Z}} a_n q^n \quad (\text{donde } q = \exp(2\pi i \tau))$$

**Definición 17.** Una función holomorfa de la forma $f = \sum_{n\in\mathbb{Z}} a_n q^n$ es holomorfa en el infinito si $a_n = 0$ para $n < 0$. Tal función es holomorfa y se anula en el infinito si $a_n = 0$ para $n \leqslant 0$.

**Definición 18.** Una forma modular de peso $k$, con respeto al grupo $\Gamma_0(N)$ es una función holomorfa $f : \mathcal{H} \to \mathbb{C}$ que cumple (37) y además, para cada $\gamma \in \Gamma_0(N)$, la función $(c\tau + d)^{-k} f(\gamma(\tau))$ es holomorfa en el infinito; Tal forma es dicha parabólica si además se anula en el infinito.

Se nota $M_k(\Gamma_0(N))$ al espacio vectorial de las formas modulares de peso $k$, con respeto al grupo $\Gamma_0(N)$ (respectivamente $S_k(\Gamma_0(N))$ al espacio de tales formas modulares que además son parabólicas).

Vamos a considerar solamente formas de peso $k = 2$.

Definimos la función $L$ asociada a una forma modular parabólica por la fórmula (donde los $a_n$ corresponden a los coeficientes de $f$):

(39) $$L(f, s) := \sum_{n\geqslant 1} a_n n^s.$$

En este caso, tenemos la relación Ejer.

(40) $$\Gamma_{\mathbb{C}}(s)L(f, s) = (2\pi)^{-s}\Gamma(s)L(f, s) = \int_0^\infty f(it)t^{s-1}dt.$$

**Definición 19.** Sea $f = \sum_n a_n(f)q^n \in M_2(\Gamma_0(N))$. Los *operadores de Hecke* son definidos como sigue.

1. Si $p$ no divide $N$, consisten de los operadores $f \mapsto T_p f$ definidos por:

$$a_n(T_p f) := a_{np}(f) + pa_{n/p}(f),$$

donde, por convención, $a_{n/p} = 0$ si $p$ no divide $n$.

2. Si $p$ divide $N$, consisten de los operadores $f \mapsto U_p f$ definidos por:

$$a_n(U_p f) := a_{np}(f).$$

**Ejercicio 13.** Verificar que las operaciones $T_p$ y $U_p$ dejan estables respectivamente a $M_2(\Gamma_0(N))$ y a $S_2(\Gamma_0(N))$.

**Teorema 20.** (Hecke, ver Diamond y Shurman (2005)) *Los operadores de Hecke conmutan. Además, si $f = \sum_n a_n(f)q^n \in S_2(\Gamma_0(N))$ es un vector propio simultáneamente para cada operador, esto es $T_p f = \lambda_p f$ y $U_p f = \lambda_p f$, entonces $a_p(f) = \lambda_p a_1(f)$. Y si $f$ es normalizada de manera que $a_1(f) = 1$, la función $L(s, f)$ se descompone en producto de Euler del modo siguiente:*

(41) $$L(s, f) = \sum_{n=1}^{\infty} a_n(f)n^{-s} = \prod_{p\,|\,N}(1 - a_p(f)p^{-s})^{-1}\prod_{p\,\nmid N}(1 - a_p(f)p^{-s} + p^{1-2s})^{-1}.$$

□

Esta función $L$ se parece mucho a la función $L$ de una curva elíptica. De hecho veremos que, con una condición suplementaria, ella satisface la ecuación funcional esperada para una curva elíptica.

Observemos ahora que la matriz $W_N := \begin{pmatrix} 0 & 1 \\ -N & 0 \end{pmatrix}$, que no pertenece a $\mathrm{SL}(2, \mathbb{Z})$, normaliza sin embargo al subgrupo $\Gamma_0(N)$, pues

$$W_N \begin{pmatrix} a & b \\ c & d \end{pmatrix} W_N^{-1} = \begin{pmatrix} d & -c/N \\ -bN & a \end{pmatrix}.$$

Así $W_N$ actúa sobre $M_2(\Gamma_0(N))$ (respectivamente $S_2(\Gamma_0(N))$), la acción estando dada por $f_{|w_N}(z) = N^{-1}z^{-2}f\left(-\frac{1}{Nz}\right)$. Notando que $W_N^2 = -NId$, vemos que $w_N$ actúa como una involución. Así

los espacios $M_2(\Gamma_0(N))$ y $S_2(\Gamma_0(N))$ se descomponen en la suma de dos subespacios propios tales que:

$$f\left(-\frac{1}{Nz}\right) = f(W_N \cdot z) = \pm N z^2 f(z)\cdot \tag{42}$$

**Teorema 21.** (Hecke) *Sea* $\epsilon = \pm 1$ *y* $f(z) = \sum_{n\geqslant 1} a_n \exp(2\pi i n z)$ *una forma modular parabólica para* $\Gamma_0(N)$ *(de peso 2) tal que*

$$f\left(-\frac{1}{Nz}\right) = \epsilon N z^2 f(z). \tag{43}$$

*Definiendo* $\Lambda(s, f) := N^{s/2}(2\pi)^{-s}\Gamma(s)L(s, f)$*, donde* $L(s, f) := \sum_{n=1}^{\infty} a_n n^{-s}$*, la función* $\Lambda(s, f)$ *tiene continuación analítica en todo el plano complejo y satisface la ecuación funcional:*

$$\Lambda(s, f) = -\epsilon\Lambda(2 - s, f). \tag{44}$$

*Además,* $\Lambda(s, f)$ *es acotada en toda banda vertical del plano complejo.*

*Demostración.* Observemos que para $z = it$ (con $t \in \mathbb{R}_+$) la ecuación (43) nos proporciona:

$$f\left(\frac{i}{Nt}\right) = -\epsilon N t^2 f(it).$$

La analogía con la ecuación (18) usada para probar la ecuación funcional de la función zeta de Riemann es clara. Así, podemos calcular, utilizando (40) y el cambio de variables $t \mapsto 1/Nt$, lo siguiente:

$$\begin{aligned}
\Lambda(s, f) &= N^{s/2}\int_0^{\infty} f(it)t^{s-1}dt \\
&= N^{s/2}\int_0^{\frac{1}{\sqrt{N}}} f(it)t^{s-1}dt + N^{s/2}\int_{\frac{1}{\sqrt{N}}}^{\infty} f(it)t^{s-1}dt \\
&= N^{-s/2}\int_{\frac{1}{\sqrt{N}}}^{\infty} f(i/Nt)t^{-s-1}dt + N^{s/2}\int_{\frac{1}{\sqrt{N}}}^{\infty} f(it)t^{s-1}dt \\
&= -\epsilon N^{\frac{1}{2}(2-s)}\int_{\frac{1}{\sqrt{N}}}^{\infty} f(it)t^{1-s}dt + N^{s/2}\int_{\frac{1}{\sqrt{N}}}^{\infty} f(it)t^{s-1}dt \\
&= \int_{\frac{1}{\sqrt{N}}}^{\infty} f(it)\left[-\epsilon N^{\frac{1}{2}(2-s)}t^{2-s}dt + N^{s/2}t^s\right]\frac{dt}{t}.
\end{aligned} \tag{45}$$

La última expresión define una función entera (utilizando una cota del tipo $|a_n| = O(n^c)$, y la consecuencia de que $|f(it)| = O(\exp(-2\pi t))$ cuando $t$ crece hasta el infinito). La $(-\epsilon)$-simetría cuando $s$ es reemplazado por $2 - s$ es clara ahora, así como la propiedad de ser acotada en toda banda vertical. □

**3.5.3 Continuación analítica y ecuación funcional de** $L(E,s)$**.** recordemos que en el capítulo, se ha definido el *discriminante* minimal $\Delta_E$ de una curva elíptica $E/\mathbb{Q}$. Existe otro invariante definido para una curva elíptica, parecido (pero distinto) llamado el *conductor* $N_E$; para el cual daremos la definición solamente a menos de una potencia de 2 y 3.

**Definición 20.** El conductor $N_E$ de $E/\mathbb{Q}$ es definido por su descomposición en factores primos como:

$$N_E := \prod_{p \mid \Delta_E} p^{n(E,p)}$$

donde $n(E,p)$ vale $+1$ en el caso de run $p$ con educción multiplicativa, vale $+2$ en el caso de reducción aditiva con $p > 3$, y donde además $2 \leqslant n(E,2) \leqslant 8$ y $2 \leqslant n(E,3) \leqslant 5$ en el caso de reducción aditiva.

El teorema siguiente muestra que la función $L(E,s)$ tiene propiedades similares a su predecesora $\zeta(s)$.

**Teorema 22.** *La función* $L(E,s)$ *se extiende en una función analítica en todo el plano complejo y satisface la ecuación funcional siguiente, denotando* $\Lambda(E,s) = N_E^{s/2}(2\pi)^{-s}\Gamma(s)L(E,s)$

$$\Lambda(E, 2-s) = \pm\Lambda(E,s) \tag{46}$$

¡Lo anterior es esencialmente equivalente al Teorema de Wiles! Más precisamente, lo que demuestra Wiles es la conjetura de Shimura–Taniyama que afirma que los coeficientes $a_n$ (obtenidos a partir de los $a_p$) de la curva elíptica son los coeficientes de una forma modular de peso 2 con respeto al grupo $\Gamma_0(N_E)$. De manera que la ecuación funcional resultante es entonces consecuencia de la ecuación funcional de la función $L$ de una forma modular.

**Comentario.** El teorema de Wiles no utiliza explícitamente las funciones $L$ pero unifica tres objetos, a priori, de origen muy distinto: una curva elíptica $E/\mathbb{Q}$, las representaciones de Galois asociadas y las formas modulares (de peso 2, con respeto a un subgrupo $\Gamma_0(N)$ y que son vectores propios para operadores de Hecke y $W_N$). Sin embargo el hilo conductor entre estos tres objetos son las funciones $L$ asociadas. En este aspecto, el teorema de Wiles puede ser visto como solo una parte (muy importante) del vasto programa de Langlands.

### 3.5.4 El signo de la ecuación funcional.

**Definición 21.** Denotamos $W(E)$ al signo $\pm 1$ que aparece en la ecuación funcional (46).

La importancia del signo $W(E)$ ("*root number*" en inglés) viene del hecho de que este determina la paridad del orden de anulación de $L(E,s)$. O sea, si denotamos $r_{\text{an}} = r_{\text{an}}(E/\mathbb{Q})$ al *rango analítico*, es decir, al orden de anulación de $L(E,s)$ en $s=1$, tenemos

$$W(E) = (-1)^{r_{\text{an}}(E/\mathbb{Q})} \tag{47}$$

En el siguiente teorema, veremos que este signo se puede calcular como producto de signos locales.

**Teorema 23.** *Sea $E$ curva elíptica definida sobre $\mathbb{Q}_p$. Existen signos locales $W_p(E)$ que se pueden calcular por fórmulas explícitas (dadas debajo) y cuyo producto es igual al signo global*[4]. *Es decir, dada una curva elíptica $E$ definida sobre $\mathbb{Q}$, tenemos:debe*

$$W(E/\mathbb{Q}) = \prod_p W_p(E/\mathbb{Q}_p).$$

Las reglas de cálculo son las siguientes:

1. $W_\infty(E) = -1$;
2. Cuando $E$ tiene buena reducción en $p$, entonces $W_p(E) = +1$;
3. Cuando $E$ tiene reducción multiplicativa en $p$, entonces $W_p(E) = -1$ (respectivamente $W_p(E) = +1$) si las dos tangentes son racionales sobre $\mathbb{F}_p$, caso desplegado, "*split, déployée*" (respectivamente si no lo son, caso no desplegado).
4. (Para $p > 2$). Cuando la reducción en $p$ es aditiva y potencialmente multiplicativa, $W_p(E) = \left(\frac{-1}{p}\right)$ (símbolo de Legendre);
5. (Para $p > 3$). Cuando la reducción en $p$ es aditiva y potencialmente buena, denotamos $e := \frac{12}{\text{mcd}\,(\text{ord}\,(\Delta_E),12)}$, entonces

$$W_p(E) = \begin{cases} \left(\frac{-1}{p}\right) & \text{si } e = 2, 6 \\ \left(\frac{2}{p}\right) & \text{si } e = 4 \\ \left(\frac{3}{p}\right) & \text{si } e = 3. \end{cases}$$

Existen también fórmulas en el caso de buena reducción potencial y $p = 2$ o $3$, pero son bastante más complicadas (ver por ejemplo Rizzo (2003)).

**Ejercicio 14.** Reconsideramos la curva elíptica $E$ del ejercicio 8 dada por $y^2 + y = x^3 - x$.

1. Mostrar que $W(E) = W_\infty(E)W_{37}(E) = -W_{37}(E)$.
2. Mostrar que $E$ tiene reducción multiplicativa en $p = 37$ y que la reducción es no desplegada.
3. Concluir que $W(E) = -1$ y $L(E, 1) = 0$.

   [En efecto $L'(E, 1) \neq 0$, de manera que $r_{\text{an}}(E/\mathbb{Q}) = r(E/\mathbb{Q}) = 1$].

**Ejercicio 15.** Reconsideramos la curva elíptica $E$ del ejercicio 9 dada por $y^2 + 9y = x^3 - 27$ (cúbica de Fermat).

1. Mostrar que $W(E) = W_\infty(E)W_3(E) = -W_3(E)$.
2. Admitiendo (o verificando con la ayuda de Rizzo (ibíd.)) que $W_3(E) = -1$, concluir que $W(E) = +1$.

   [En efecto $L(E, 1) \neq 0$, de manera que $r_{\text{an}}(E/\mathbb{Q}) = r(E/\mathbb{Q}) = 0$].

[4]Para que la fórmula sea exacta, se incluir al "primo arquimediano" $\infty$, cuyo signo es $W_\infty(E) = -1$.

**3.6 Valor en $s = 1$.** La primera parte de la conjetura de Birch y Swinnerton-Dyer dice que, para toda curva elíptica $E/\mathbb{Q}$ tenemos:

$$r_{\text{an}}(E/\mathbb{Q}) = r(E/\mathbb{Q}).$$

Esto es, que el orden de anulación de la función $L(E,s)$ en $s = 1$ es igual al rango del grupo de Mordell–Weil.

Es importante remarcar que la conjetura de Birch y Swinnerton-Dyer da muchas más precisiones con respecto a la teoría de curvas elípticas, las cuales vamos ahora a mencionar y desarrollar. Como la fórmula describiendo el residuo de la función $\zeta_K(s)$ en $s = 1$ involucra muchas informaciones aritméticas, el valor de la función $L(E,s)$ en $s = 1$ contiene bastante de la aritmética de la curva $E/\mathbb{Q}$, o debería contener, según esta conjetura.

**3.6.1 El grupo de Shafarevich–Tate (esbozo).** El grupo de *Shafarevich–Tate*[5] define una medida de las obstrucciones cohomológicas al principio de Hasse para cúbicas. De manera más técnica, su definición es:

$$Ш(E/\mathbb{Q}) := \ker \left\{ H^1(\mathrm{Gal}\,(\bar{\mathbb{Q}}/\mathbb{Q}), E) \to \prod_p H^1(\mathrm{Gal}\,(\bar{\mathbb{Q}}_p/\mathbb{Q}_p), E_{\mathbb{Q}_p}) \right\}.$$

No disponemos del tiempo para desarrollar las herramientas (como la cohomología de Galois) para entender mejor a este grupo, así que solamente explicaremos cómo el grupo Ш interviene naturalmente en el cálculo del grupo de Mordell–Weil. El proceso de la demostración del teorema de Mordell–Weil débil conduce naturalmente a establecer una sucesión exacta de grupos finitos o también de módulos de Galois (ver Silverman (1986), Capitulo X, §4).

$$0 \to E(\mathbb{Q})/nE(\mathbb{Q}) \to S^{(n)}(E/\mathbb{Q}) \to Ш(E/\mathbb{Q})[n] \to 0.$$

Donde el grupo en el medio, llamado $n$-grupo de Selmer, es finito y calculable.

La parte anulada por $[n]$ del grupo $Ш(E/\mathbb{Q})$ es finita, aunque no se sabe en general calcular este grupo. Es una parte de la conjetura de Birch y Swinnerton-Dyer que el grupo total sea finito, pero esto solamente se le ha demostrado en casos particulares. El grupo de Selmer en cambio, es mucho más conocido, por ejemplo, resultados excelentes de Bhargava (junto con Shankar) han permitido demostrar resultados involucrando promedios, que son muy precisos.

**Teorema 24.** *(Bhargava–Shankar) Sea $n \in \{2,3,4,5\}$ y $\sigma(n)$ la suma de los divisores de $n$, entonces*

$$\lim_{T\to\infty} \frac{\sum_{H(E)\leqslant T} |S^{(n)}(E/\mathbb{Q})|}{\sum_{H(E)\leqslant T} 1} = \sigma(n).$$

*Como corolario se obtiene por ejemplo que:*

$$\limsup_{T\to\infty} \frac{\sum_{H(E)\leqslant T} 5^{r(E/\mathbb{Q})}}{\sum_{H(E)\leqslant T} 1} \leqslant 6.$$

[5] La notación clásica con la letra rusa o cirílica Ш se debe a un homenaje a Igor Shafarevich (Шафаре́вич).

**3.6.2 Conjetura de Birch y Swinnerton-Dyer.** Hemos definido el regulador $\mathrm{Reg}(E/\mathbb{Q})$ y (brevemente) el grupo de Shafarevich. Las otras cantidades incluidas en la fórmula de Birch y Swinnerton-Dyer son las siguientes.

**Definición 22.** Sea $y^2+a_1xy+a_3y = x^3+a_2x^2+a_4x+a_6$ un modelo de ecuación de Weierstass minimal de una curva elíptica $E/\mathbb{Q}$. Entonces la forma diferencial de Néron está dada por

$$\omega_E := \frac{dx}{2y + a_1x + a_3}.$$

Y el período real de $E$ corresponde al número real positivo:

$$\Omega_E := \int_{E(\mathbb{R})} |\omega_E|.$$

Se puede comparar el tamaño de $\Omega_E$ así (ver Hindry (2007)) con constantes que se podrían explicitar:

$$c_1 H(E) \leqslant \Omega_E^{-1} \leqslant c_2 H(E) \log H(E) \tag{48}$$

**Definición 23.** Sea $E/\mathbb{Q}$ una curva elíptica. El subgrupo $E^0(\mathbb{Q}_p)$ corresponde al subgrupo de puntos que se reducen módulo $p$ en un punto no singular. El número de Tamagawa local es definido como el índice $c_p = (E(\mathbb{Q}_p) : E^0(\mathbb{Q}_p)$, mientras que el número de Tamagawa global es definido como el producto $\prod_p c_p$.

Las informaciones simples sobre $c_p$ son las siguientes: $c_p = 1$ en el caso de buena reducción, $c_p = \mathrm{ord}_p(\Delta_E)$ en el caso de reducción multiplicativa desplegada y $c_p \leqslant 4$ en los otros casos.

Gracias a esto último, podemos ahora expresar la conjetura completa (ver Tate (1995)):

**Conjetura 3.** (Conjetura de Birch y Swinnerton-Dyer) Sea $E/\mathbb{Q}$ una curva elíptica y $L(E,s)$ su función $L$ asociada. Entonces:

1. $r_{\mathrm{an}}(E/\mathbb{Q}) = r(E/\mathbb{Q})$.
2. El grupo $Ш(E/\mathbb{Q})$ es finito.
3. El coeficiente principal de $L(E,s)$ en $s = 1$ está dado por:

$$\lim_{s\to 1} \frac{L(E,s)}{(s-1)^r} = |Ш(E/\mathbb{Q})| \mathrm{Reg}(E/\mathbb{Q}) \Omega_E \frac{\prod_p c_p}{|E(\mathbb{Q})_{\mathrm{tor}}|^2}. \tag{49}$$

Notamos $L^*(E,1)$ a la primera expresión de acá arriba dada por el límite cuando $s$ tiende hasta 1. Damos ahora la aplicación prometida a la estimación de $\mathrm{Reg}(E/\mathbb{Q})$. Aceptando la fórmula conjetural (49), recordando (que $Ш(E/\mathbb{Q})| \prod_p c_p|$ es un entero (ver la definición 23) y que $|E(\mathbb{Q})_{\mathrm{tor}}| \leqslant 16$ obtenemos:

$$\mathrm{Reg}(E/\mathbb{Q}) \leqslant 16^2 L^*(E,1)\Omega_E^{-1} \leqslant c_2 H(E) \log H(E) L^*(E,1).$$

Utilizando la ecuación funcional, las técnicas de teoría analítica de números y en particular, el principio de Phragmén–Lindelöf, nos proporcionan (ver Iwaniec y Kowalski (2004) y Hindry (2010)):

**Lema 4.** *Sea $L(E,s)$ la función $L$ asociada a una curva elíptica $E/\mathbb{Q}$ de conductor $N_E$. Tenemos las estimaciones:*

1. $L^*(E,1) \leqslant 2^r N_E^{1/4}(\log N_E)^2$.
2. *Y si suponemos que la función $L(E,s)$ verifica la hipótesis de Riemann (es decir, $L(E,s) \neq 0$ cuando $\Re(s) > 1$), entonces:*
$$L^*(E,1) \leqslant C_\epsilon N_E^\epsilon \leqslant C'_\epsilon H(E)^\epsilon.$$

Juntando toda la información obtenida hasta ahora, llegamos naturalmente a la conjetura expuesta acá abajo, la cual es implicada por las conjeturas de Birch y Swinnerton-Dyer y de Riemann (aplicada a $L(E,s)$). Observamos que la conjetura es puramente diofántica y no hace referencia a las funciones $L$ (comparar con Lang (1983)).

**Conjetura 4.** Existe, para cada $\epsilon > 0$, una constante $C_\epsilon$, tal que para toda curva elíptica $E/\mathbb{Q}$ tenemos:

$$\operatorname{Reg}(E/\mathbb{Q}) \leqslant C_\epsilon H(E)^{1+\epsilon}. \tag{50}$$

*Observación* 4. Si no queremos utilizar la conjetura de Riemann, recordando que $N_E \leqslant \Delta_E \leqslant cH(E)^{12}$, obtenemos una cota más débil de la forma $\operatorname{Reg}(E/\mathbb{Q}) \leqslant c_\epsilon H(E)^{4+\epsilon}$.

*Observación* 5. El mismo argumento condicional proporciona cotas del mismo orden $H(E)^{1+\epsilon}$ para el tamaño del grupo de Shafarevich–Tate cuando $L(E,1) \neq 0$, es decir, cuando $r = 0$ y $\operatorname{Reg}(E/\mathbb{Q}) = 1$ por convención (comparar con Goldfeld y Szpiro (1995)). Demostrando luego que $\operatorname{Reg}(E/\mathbb{Q}) \geqslant C_\epsilon H(E)^{-\epsilon}$, con tal suposición, este resultado se puede extender la cota a todas las curvas.

Es el momento de aclarar las respuestas (parciales) a las preguntas iniciales (sección 3.1) que las consideraciones de funciones $L$ y $\zeta$ nos hicieron plantear, *condicionalmente* a la conjetura de Birch y Swinnerton-Dyer.

1. La demostración del teorema de Mordell–Weil proporciona solamente una cota superior para el rango $r = r(E/\mathbb{Q})$,

2. El rango analítico es en general fácil de calcular, la paridad se calcula todavía más rápidamente, calculando $W(E)$: calculamos $L(E,1)$, $L'(E,1)$ y así hasta encontrar un valor no nulo [6].

3. Existen generadores $P_i$ del grupo de Mordell–Weil tales que, si se les ordena en orden creciente $\hat{h}(P_1) \leqslant \hat{h}(P_2) \leqslant \ldots$, entonces
$$\hat{h}(P_i) \leqslant C_\epsilon H(E)^{\frac{1}{r-i+1}+\epsilon}$$
Una vez que se tiene una cota, se dispone de un algoritmo obvio para buscar de manera exhaustiva puntos racionales, aunque el algoritmo es pésimo desde el punto de vista de la complejidad computacional.

[6] La única dificultad algorítmica es la de comprobar – o aceptar – que una derivada es nula cuando su valor aproximado es, digamos, menor que $10^{-100}$.

**Indicaciones.** (para el ejercicio 1).
1) La única solución es $3^3 = 5^2 + 2$. Se trata de encontrar los puntos enteros sobre $y^2 = x^3 - 2$. Mostrar que $x, y$ son impares; considerar el anillo, que se resultará ser principal, $A = \mathbb{Z}[i\sqrt{2}]$ y mostrar que $y + i\sqrt{2}$ es un cubo en $A$. Concluir.
2) Las soluciones son $6 = 2\times 3 = 1\times 2\times 3$ y $210 = 14\times 15 = 5\times 6\times 7$. La pregunta es esencialmente equivalente a determinar los puntos enteros de la curva $y^2+y = x^3-x$ considerada en los ejercicios 8 y 14. Sea $Q \in E(\mathbb{Q})$, mostrar que si $mQ$ es un punto entero, entonces $Q$ también. Mostrar que $E(\mathbb{R})$ tiene dos componentes conexas: $C_0$ conteniendo el punto al infinito y $C_1$; determinar en seguida los puntos enteros en el componente $C_1$ (que es compacto en el plano afín). Una vez demostrado que el rango es igual a 1 (es la parte más difícil, se recomienda admitirla o ver Silverman (1986), *ejercicio* 10.9), mostrar que $E(\mathbb{Q})$ es generado por $P = (0,0)$, que $\pm P$, $\pm 2P$, $\pm 3P$, $\pm 4P$, $\pm 6P$ son puntos enteros, pero que $8P$ y $12P$ no lo son (notar que $6P = (6,14)$). Concluir argumentando que si $mP$ es un punto entero, escribiendo $m = 2^u m_1$ (con $m_1$ impar), entonces debemos tener que $m_1 P$ es entero y que pertenece a $C_1$.

## Apéndice: Esquemas explicados para niños

*Ver por ejemplo Hartshorne (1977) para más detalles sobre la noción de esquema.*

Se puede describir una variedad proyectiva sobre un cuerpo $K$, como $\mathbb{P}^n$, por ejemplo una curva elíptica, como recubierta por abiertos afines y pegados (de acuerdo a la noción topológica de pegamento).

**Ejemplo 16.** 1) La línea proyectiva $\mathbb{P}^1$ se puede describir como $U_0 \cup U_1$, donde el abierto afín $U_0 = \{(x_0, x_1) \in \mathbb{P}^n \mid x_0 \neq 0\}$ es isomorfo a la línea $\mathbb{A}^1$ mediante la aplicación $x \mapsto (1, x)$ y el abierto afín $U_1 = \{(x_0, x_1) \in \mathbb{P}^n \mid x_1 \neq 0\}$ es isomorfo a la línea $\mathbb{A}^1$ mediante la aplicación $t \mapsto (t, 1)$, los cuales son pegados posteriormente identificando $t = x^{-1}$.
2) Una curva elíptica $E$ puede ser descrita como la unión de dos abiertos afines $V_1 \cup V_2$ donde $V_1 = \{(x, y) \in \mathbb{A}^2 \mid y^2 = x^3 + ax + b\}$ y $V_2 = \{(u, v) \in \mathbb{A}^2 \mid v^2 = u(bu^3 + au^2 + 1)\}$, pegados a través de la identificación $(u, v) = (1/x, y/x^2)$.

Las variedades afines pueden ser identificadas con su anillo de funciones, es decir, si $V$ es una variedad afín en $\mathbb{A}^n$, se la puede describir como el conjunto de ceros de un ideal $I \subset K[X_1, \ldots, X_n]$; el anillo de funciones (o $K$-álgebra) de $V$ es $\mathcal{O}(V) = K[X_1, \ldots, X_n]/I$; un punto $a = (a_1, \ldots, a_n)$ en $V$ corresponde al ideal generado por (las imagenes de) los $X_i - a_i$. Existe una correspondencia entre aplicaciones $\phi : V \to W$ de variedades afines y homomorfismos de anillos (o $K$-álgebra) $\phi^* : \mathcal{O}(W) \to \mathcal{O}(V)$ (donde $\phi^*(f) = f \circ \phi$).

Los puntos de $\mathbb{A}^1$ (o $U_0$ o $U_1$) se identifican con los ideales maximales de $K[X]$ y las aplicaciones de $\mathbb{A}^1$ hacia $\mathbb{A}^1$ corresponden a endomorfismos de anillo de $K[X]$. Así, $V_1$ se identifica con el anillo $B := K[x, y]/(-y^2 + x^3 + ax + b)$ por ejemplo, y una aplicación de $V_1$ hacia $\mathbb{A}^1$ corresponde a un homomorfismo de anillos $K[X] \to B$.

Grothendieck refinó este concepto introduciendo lo que es el espectro de un anillo $A$ (para un anillo unitario cualquiera) cuyos puntos son los ideales primos de A. El espectro de un anillo $A$, denotado spec $(A)$, es el conjunto de sus ideales primos junto con la topología de Zariski; la cual es definida con la propiedad $P \subset Q$ (de inclusión de ideales primos) significando que $Q$ está en la cerradura de $\{P\}$, como punto de spec $(A)$ (en una descripción más completa, tendríamos que definir también lo que es el *haz* estructural (ver Hartshorne (ibíd.))).

Un esquema afín es simplemente el espectro de un anillo. Los morfismos correspondientes spec $(A) \to$ spec $(B)$ consideran homomorfismos de anillos $f : B \to A$ y son definidos por $P \to f^{-1}(P)$, donde $P$ es un ideal primo. Un esquema general es definido de manera similar, pegando (topológicamente) esquemas afines.

**Ejemplo 17.** 1) Cuando $K$ es un cuerpo, el esquema spec $(K)$ es reducido a un punto; en cambio el esquema spec $(\mathbb{Z})$ tiene una infinidad de puntos cerrados (en correspondencia con números primos o ideales $p\mathbb{Z}$) y un punto denso (que corresponde al ideal nulo o a la imagen de spec $(\mathbb{Q}) \to$ spec $(\mathbb{Z})$ correspondiendo a la inclusión de anillos $\mathbb{Z} \hookrightarrow \mathbb{Q}$).
2) El esquema spec $(\mathbb{Z}[X])$ tiene varios tipos de puntos: los punto cerrados, que corresponden a ideales $I = (p, P)$ generados por un número primo $p$ y un polinomio $P$ irreducible módulo $p$, los puntos de dimensión (o codimensión) uno, correspondiendo a ideales principales generados por un número primo $p$ o un polinomio irreducible $P$, y finalmente el punto genérico, que corresponde al ideal nulo Ejer.
3) El esquema $\mathbb{P}^1_{\mathbb{Z}}$ es recubierto por dos copias de spec $(\mathbb{Z}[X])$ de manera análoga a $\mathbb{P}^1$ sobre $K$.
4) Se puede definir un esquema $E/\mathbb{Z}$, donde $E$ es una curva elíptica, pegando dos esquemas afines $V_1$ y $V_2$ definidos así:
el esquema $V_1$ es el espectro de $\mathbb{Z}[x, y]/(-y^2 + x^3 + ax + b)$
el esquema $V_2$ el espectro de $\mathbb{Z}[u, v]/(-v^2 + u(bu^3 + au^2 + 1))$.

En este lenguaje, una variedad afín "usual" sobre un cuerpo es el espectro de una $K$-álgebra integral (un dominio) de tipo finito. De hecho, el lenguaje de esquemas es mucho más rico, por ejemplo podemos hablar de una línea triple (con multiplicidad tres) $X + Y = 0$ como el esquema $L^{(3)} =$ spec $(K[X, Y]/(X+Y)^3)$ o de una variedad sobre $\mathbb{Z}$, como un esquema $f : \mathfrak{X} \to$ spec $(\mathbb{Z})$. En este caso, para el punto genérico $\eta$ de spec $(\mathbb{Z})$, es decir, el punto que corresponde al ideal nulo, se obtiene la fibra genérica $X = f^{-1}(\eta)$, que es una variedad sobre $\mathbb{Q}$. Además, para cada ideal maximal $p$ de $\mathbb{Z}$, la fibra nos proporciona una variedad $\mathfrak{X}_p = f^{-1}(P)$ que es la reducción de $X$ módulo $p$.

## Referencias

Marusia Rebolledo
marusia.rebolledo@uca.fr
Université Clermont Auvergne
Laboratoire de mathématiques Blaise Pascal

Marc Hindry
marc.hindry@imj-prg.fr
Université Paris Diderot
Institut de mathématiques de Jussieu – Paris rive gauche

# 4 | Crecimiento en grupos y expansores

# CRECIMIENTO Y EXPANSIÓN EN $SL_2$

Harald Andrés Helfgott

## 1 Prefacio

Ésta es una breve introducción al estudio del crecimiento en los grupos finitos, con $SL_2$ como ejemplo. Su énfasis cae sobre los desarrollos de la última década, provenientes en parte de la combinatoria.

El texto – basado en parte en Helfgott (2015) – consiste, en esencia, en notas de clase para un curso en la escuela de invierno AGRA II en la Universidad San Antonio Abad, Cusco, Perú (10–21 agosto 2015), incluyendo algunos ejercicios. El curso fue la primera mitad de una unidad; la segunda mitad, sobre expansores en conexión al espacio hiperbólico, corrió a cargo de M. Belolipetsky.

El tópico tiene una intersección apreciable con varios otros textos, incluyendo el libro de T. Tao (2015) y las notas de E. Kowalski (2013). El tratamiento en el Capítulo 4 difiere un tanto de Helfgott (2015), en la medida que sigue un tratamiento más global (grupos algebraicos) y menos local (álgebras de Lie); en ésto puede detectarse una influencia de Tao (2015) y Kowalski (2013) (y, en última instancia, Larsen y Pink (2011)). No parecen haber desventajas o ventajas decisivas en ésto; simplemente he tomado la oportunidad de explorar un formalismo distinto.

Sin duda, Helfgott (2015) da más detalles que el texto presente, tanto de tipo histórico como de tipo puramente matemático; su tema también es más amplio. La meta principal aquí es dar una introducción concisa, accesible y en cierto sentido participativa al tópico.

Las brevísimas introducciones al grupo $SL_2$ y a la geometría algebraica en general tienen como intención hacer que el texto comprensible sea comprensible para estudiantes de distintas áreas, aparte de ser partes necesarias de la cultura general. Se pide la paciencia del los lectores para los cuales tales introducciones son innecesarias.

## 2 Introducción

Nuestra tema es el crecimiento en los grupos; nuestro ejemplo serán los grupos $SL_2(K)$, $K$ un cuerpo finito.

Qué se quiere decir aquí por *crecimiento*? Hay diferentes puntos de vista, dependiendo del área. La manera más concreta de expresar la cuestión es quizás la siguiente: tenemos un subconjunto finito $A$ de un grupo $G$. Consideremos los conjuntos

$$\begin{aligned} &A, \\ A \cdot A &= \{x \cdot y : x, y \in A\}, \\ A \cdot A \cdot A &= \{x \cdot y \cdot z : x, y, z \in A\}, \\ &\ldots \\ A^k &= \{x_1 x_2 \ldots x_k : x_i \in A\}. \end{aligned}$$

Escribamos $|S|$ por el número de elementos de un conjunto finito $S$. La pregunta es: qué tan rápido crece $|A^k|$ a medida que $k$ se incrementa?

Tal cuestión ha sido estudiada desde la perspectiva de la combinatoria aditiva (caso de $G$ abeliano) y de la teoría de grupos geométrica ($G$ infinito, $k \to \infty$). También hay varios conceptos relacionados, de suma importancia, provenientes de la teoría de grafos y de analogías con la geometría: *diámetros*, *expansores*, etc.

Ahora bien, porqué elegir a los grupos $\mathrm{SL}_2(K)$ como primer caso a estudiar, más allá de la necesidad, en una exposición, de comenzar por un caso concreto?

**2.1 Los grupos $\mathrm{SL}_2(R)$.** Sea $R$ un anillo; por ejemplo, podemos tomar $R = \mathbb{Z}$, o $R = \mathbb{Z}/p\mathbb{Z}$. Definimos

$$\mathrm{SL}_2(R) = \left\{ \begin{pmatrix} a & b \\ c & d \end{pmatrix} : a, b, c, d \in R, \;\; ad - bc = 1. \right\}.$$

La letra “S” en SL viene de e**s**pecial (lo que aquí quiere decir: de determinante igual a 1), mientras que “L” viene de **l**ineal (por tratarse de un grupo de matrices). El número 2 viene del hecho que éstas son matrices 2x2.

El grupo $\mathrm{SL}_2(R)$ puede verse por lo menos de dos formas: como un grupo abstracto, y como un grupo de transformaciones geométricas. Visto de una manera o la otra, se trata de un buen caso a estudiar, pues es, por así decirlo, el objeto más sencillo en demostrar toda una gama de comportamientos complejos. Veamos cómo.

**2.1.1 La estructura del grupo $\mathrm{SL}_2(R)$.** Dado un grupo $G$, nos interesamos en sus subgrupos $H < G$, y, en particular, en sus subgrupos *normales* $H \lhd G$. Dado $H \lhd G$, podemos decir que $G$ se descompone en $H$ y $G/H$. Un grupo $G$ sin subgrupos normales (aparte de $\{e\}$ y $G$) se llama *simple*.

Los grupos simples juegan un rol similar al de los primos en los enteros. Es fácil ver que, para todo grupo finito $G$, existen

$$\{e\} = H_0 \lhd H_1 \lhd H_2 \lhd \cdots \lhd H_k = G \tag{1}$$

tales que $H_i/H_{i-1}$ es simple y notrivial para $1 \leqslant i \leqslant k$. El teorema de Jordan–Hölder nos dice que tal decomposición es en esencia única: los factores $H_{i+1}/H_i$ en (1) están determinados por $G$, y a lo más su orden puede cambiar.

Un grupo resoluble es un grupo que tiene una decomposición tal que $H_{i+1}/H_i$ es abeliano para todo $i$. Como hemos dicho, la combinatoria aditiva ha estudiado tradicionalmente el crecimiento en los grupos abelianos. El estudio del crecimiento en los grupos resolubles esta lejos de ser trivial, o de reducirse por completo al crecimiento en los grupos abelianos. Empero, una decomposición $H \triangleleft G$ reduce los problemas de crecimiento (como muchos otros) al estudio de (a) los grupos $H$ y $G/H$, (b) la *acción* de $G/H$ sobre $H$. Por ello, en últimas cuentas, tiene sentido concentrarse en el estudio de los grupos simples, y, en particular, en el estudio del crecimiento en los grupos simples no abelianos.

Sea $K$ un cuerpo finito. El grupo $\mathrm{SL}_2(K)$ no es ni abeliano ni resoluble. El *centro*

$$Z(G) = \{g \in G : \forall h \in G \quad hg = gh\}$$

de un grupo $G$ es siempre un subgrupo normal de $G$. Ahora bien, para $G = \mathrm{SL}_2(K)$, $Z(G)$ es igual a $\{I, -I\}$; así, a menos que $K$ sea el cuerpo $\mathbb{F}_2$ con dos elementos, $Z(G)$ no es el grupo trivial $\{e\} \neq \{I\}$, y por lo tanto $G = \mathrm{SL}_2(K)$ no es simple. Empero, el cociente

$$\mathrm{PSL}_2(K) := \mathrm{SL}_2(K)/Z(\mathrm{SL}_2(K)) = \mathrm{SL}_2(K)/\{I, -I\}$$

sí es simple, para $|K|$ finito y mayor que 3.[1]

*Comentario de índole cultural.* Así, al considerar $\mathrm{SL}_2(K)$ para $K$ variando sobre todos los cuerpos finitos, obtenemos toda un conjunto infinito de grupos finitos simples $\mathrm{PSL}_2(K)$. Se trata, por así decirlo, de la familia más sencilla de grupos finitos simples, junto con aquella dada por los grupos *alternantes* $A_n$. (El grupo $A_n$ es el único subgrupo de índice 2 del *grupo simétrico* $S_n$, el cual, a su vez, consiste en las $n!$ permutaciones de $n$ elementos, con la composición como operación del grupo.) En verdad, el famoso Teorema de la Clasificación de grupos simples nos dice que hay dos tipos de familias infinitas de grupos simples: las familias de grupos de matrices, como $\mathrm{PSL}_2(K)$, y la familia $A_n$, aparte de un número finito de grupos especiales (como el así llamado "monstruo").

Veamos la estructura de $G = \mathrm{SL}_2(K)$ en más detalle. Si bien $\mathrm{SL}_2(K)$ no tiene subgrupos normales más allá de $\{I\}$, $\{I, -I\}$ y $\mathrm{SL}_2(K)$, tiene subgrupos de varios otros tipos. Los más interesantes para nosotros son los *toros*; en $\mathrm{SL}_2(K)$, aparte del grupo trivial, todos son *toros máximos*. Recordamos que el centralizador de un elemento $g \in G$ es el grupo

(2) $$C(g) = \{h \in G : hg = gh\}.$$

Un *toro máximo* (denotado por $T(K)$) es un grupo $C(g)$ donde $g$ es *regular semisimple*; un elemento $g \in G$ es *regular semisimple* si tiene dos valores propios distintos. Esto es lo mismo que decir que $T(K) = \sigma D \sigma^{-1} \cap \mathrm{SL}_2(K)$, donde $D$ es el grupo de matrices diagonales

$$D = \left\{ \begin{pmatrix} r & 0 \\ 0 & r^{-1} \end{pmatrix} : r \in \overline{K} \right\}$$

y $\sigma \in \mathrm{SL}_2(\overline{K})$. Aquí $\overline{K}$ es la compleción (clausura) algebraica de $K$. Nótese que $\sigma$ puede o puede no estar en $\mathrm{SL}_2(K)$.

*Ejemplo.* Sea $K = \mathbb{R}$. Si $\sigma \in \mathrm{SL}_2(\mathbb{R})$, entonces $T(K) = \sigma D \sigma^{-1} \cap \mathrm{SL}_2(K)$ es de la forma

$$\sigma \left\{ \begin{pmatrix} r & 0 \\ 0 & r^{-1} \end{pmatrix} : r \in \mathbb{R}^* \right\} \sigma^{-1},$$

[1] Este es un hecho no trivial. Para $|K|$ un número primo $> 3$, fue mencionado, pero no probado, por Galois (1831). Fue probado para tales $K$ por Moore (1896), y para $K$ finito general por Jordan (1870). La nota Conrad (s.f.) contiene tanto estos apuntes históricos como una prueba completa para $\mathrm{SL}_n(K)$, donde $n > 2$ o $n = 2$ y $|K| > 3$.

y es más bien de la forma

$$\tau \left\{ \begin{pmatrix} \cos\theta & -\sin\theta \\ \sin\theta & \cos\theta \end{pmatrix} : \theta \in \mathbb{R}/\mathbb{Z} \right\} \tau^{-1} \tag{3}$$

(para algún $\tau \in \mathrm{SL}_2(\mathbb{R})$) si $\sigma \notin \mathrm{SL}_2(\mathbb{R})$. Está claro que, en el segundo caso, $T(K)$ es un círculo, es decir, un toro 1-dimensional (en el sentido tradicional de "toro").

**2.1.2 El grupo de transformaciones $\mathrm{SL}_2(\mathbb{R})$.** Uno de varios modelos equivalentes para la geometría hiperbólica en dos dimensiones es el semiplano de Poincaré, también llamado simplemente semiplano superior:

$$\mathbb{H} = \{(x, y) \in \mathbb{R}^2 : y > 0\}$$

con la métrica dada por

$$ds = \frac{\sqrt{dx^2 + dy^2}}{y}.$$

Las isometrías de $\mathbb{H}$ que preservan la orientación son las transformaciones lineares fraccionales

$$z \mapsto \frac{az + b}{cz + d}, \tag{4}$$

donde $a, b, c, d \in \mathbb{R}$ y $ad - bc \neq 0$. Es fácil ver que esto induce una biyección $\mathrm{PSL}_2(\mathbb{R})$ al conjunto de transformaciones lineares fraccionales. Escribamos $gz$ para la imagen de $z$ bajo una transformación (4) inducida por una matriz correspondiente a un elemento $g \in \mathrm{PSL}_2(\mathbb{R})$. No es nada difícil verificar que, para $g_1, g_2 \in \mathrm{PSL}_2(\mathbb{R})$,

$$g_1(g_2 z) = (g_1 g_2) z.$$

En otras palabras, tenemos un isomorfismo de $\mathrm{PSL}_2(\mathbb{R})$ al grupo (con la composición como operación) de las transformaciones lineares fraccionales.

Podemos considerar subacciones; por ejemplo, el *grupo modular (completo)* $\mathrm{SL}_2(\mathbb{Z})$ actúa sobre $\mathbb{H}$. El cociente $\mathrm{SL}_2(\mathbb{Z})\backslash\mathbb{H}$ es de volúmen finito sin ser compacto. Tambien podemos considerar los *grupos modulares de congruencia*

$$\Gamma(N) = \left\{ \begin{pmatrix} a & b \\ c & d \end{pmatrix} \in \mathrm{SL}_2(\mathbb{Z}) : a \equiv d \equiv 1 \bmod N, \quad b \equiv c \equiv 0 \bmod N \right\}$$

para $N \geqslant 1$. Claro está, $\Gamma(N)$ es el núcleo ("kernel") de la reducción $\mathrm{SL}_2(\mathbb{Z}) \to \mathrm{SL}_2(\mathbb{Z}/N\mathbb{Z})$, y, por lo tanto, $\Gamma(N)\backslash\mathbb{H}$ consiste de $|\mathrm{SL}_2(\mathbb{Z}/N\mathbb{Z})|$ copias de $\mathrm{SL}_2(\mathbb{Z})\backslash\mathbb{H}$.

**2.2 Perspectivas sobre el crecimiento en los grupos.** El crecimiento en los grupos ha sido estudiado de varias perspectivas distintas. Nos concentraremos después en desarrollos relativamente recientes que toman sus herramientas en parte de algunas de estas áreas (clasificación de subgrupos, combinatoria aditiva) y su relevancia de otras (el estudio de los diámetros y la expansión). Hay aún otras areas de suma importancia cuya relación con nuestro tema recién comienza a elucidarse (teoría de modelos, teoría de grupos geométrica). Demos una mirada a vuelo de pájaro.

*Combinatoria aditiva.* Éste es en verdad un nombre reciente para un campo de estudios más antiguo, con una cierta intersección con la *teoría aditiva de los números*. Se puede decir que la combinatoria aditiva se diferencia de esta última cuando considera el crecimiento de conjuntos bastante arbitrarios, y no sólamente el de conjuntos como los primos o los cuadrados. Uno de los resultados claves es el teorema de Freĭman (1973), el cual clasifica los subconjuntos finitos $A \subset \mathbb{Z}$ tales que $A + A$ no es mucho más grande que $A$. Ruzsa dio una segunda prueba Ruzsa (1991), más general y más simple, e introdujo muchos conceptos ahora claves en el área.

El uso del signo $+$ y de la palabra *aditiva* muestran que, hasta recientemente, la combinatoria aditiva estudiaba grupos abelianos, si bien algunas de sus técnicas se generalizan a los grupos no abelianos de manera natural.

*Clasificación de subgrupos.* Sea $A$ un subconjunto de un grupo $G$. Asumamos que $A$ contiene la identidad $e \in G$. Entonces $|A \cdot A| = |A|$ sí y sólo sí $A$ es un subgrupo de $G$; en otras palabras, clasificar los subgrupos de $G$ equivale a clasificar los subconjuntos de $G$ que no crecen.

Clasificar los subgrupos de un grupo es a menudo algo lejos de trivial – más aún si se desea emprender tal tarea sin utilizar la Clasificacíon de los grupos simples (una herramienta muy fuerte, cuya prueba fue inicialmente juzgada incompleta o poco satisfactoria por muchos). Resulta ser que los trabajos en este área basados sobre argumentos elementales, antes que sobre la Clasificación, son a veces robustos: pueden ser adaptados para darnos información, no sólo sobre los subgrupos de $G$, sino sobre los subconjuntos $A$ de $G$ que crecen poco.

*Diámetros y tiempos de mezcla.* Sea $A$ un conjunto de generadores de un grupo $G$; en otras palabras, $A \subset G$ es tal que todo elemento $g$ de $G$ puede escribirse como un producto $g = x_1 x_2 \dots x_r$ para alguna elección de $x_i \in A \cup A^{-1}$. El *diámetro* de $G$ con respecto a $A$ es el $k$ mínimo tal que todo elemento $g$ de $G$ puede escribirse como $g = x_1 x_2 \dots x_r$ con $x_i \in A$ y $r \leqslant k$. Si $G$ es finito, el diámetro es necesariamente finito.

Por qué hablamos de "diámetro"? El *grafo de Cayley* $\Gamma(G, A)$ es el grafo que tiene $G$ como su conjunto de vértices y $\{(g, ag) : g \in G, a \in A\}$ como su conjunto de aristas. Podemos definir la distancia $d(g_1, g_2)$ entre $g_1, g_2 \in G$ como la longitud del camino más corto de $g_1$ a $g_2$, donde se define que la longitud de cada arista es 1. Definimos el diámetro de un grafo como definimos el de cualquier figura: será el máximo de la distancia $d(g_1, g_2)$ para toda elección posible de vértices $g_1$, $g_2$. Es fácil verificar que $\mathrm{diam}(\Gamma(G, A))$ es igual al diámetro de $G$ con respecto a $A$ que acababamos de definir.

La conjetura de Babai y Seress (1988, p. 176) postula que, si $G$ es simple y no abeliano, entonces, para cualquier conjunto de generadores $A$ de $G$,

$$\mathrm{diam}(\Gamma(G, A)) \ll (\log |G|)^{O(1)},$$

donde las constantes implícitas son absolutas.

(Un poco de notación. Sean $f, g$ funciones de un conjunto $X$ a $\mathbb{C}$. Como es habitual en la teoría analítica de números, para nosotros, $f(x) \ll g(x)$, $g(x) \gg f(x)$ y $f(x) = O(g(x))$ quieren decir la misma cosa: hay $C > 0$ y $X_0 \subset X$ finito ("constantes implícitas") tales que $|f(x)| \leqslant C \cdot g(x)$ para todo $x \in X$ fuera de $X_0$. (En verdad, necesitamos que $f(x)$ y $g(x)$ estén bien definidas sólo para $x$ fuera de $X_0$.) Escribimos $\ll_a$, $\gg_a$, $O_a$ si $X_0$ y $C$ dependen de una cantidad $a$ (digamos). Si $X_0$ y $C$ no dependen de nada, las llamamos constantes *absolutas*.)

El *tiempo de mezcla* es el $k$ mínimo tal que, si $x_1, x_2, \ldots, x_k$ son tomados al azar en $A$ con la distribución uniforme en $A$, la distribución del producto $x_1 \cdots x_k$ (o, lo que es lo mismo, la distribución del resultado de una caminata aleatoria de longitud $k$ en $\Gamma(G, A)$) esta cerca de la distribución uniforme en $G$. Hablamos de distintos tiempos de mezcla dependiendo de lo que se quiera decir por "cerca". El estudio de los tiempos de mezcla ha tenido no solo un fuerte color probabilístico (véanse las referencias Diaconis y Saloff-Coste (1993), Levin, Peres y Wilmer (2009)) sino a menudo también algorítmico (e.g. en Babai, Beals y Seress (2004)).

*Expansores y huecos espectrales.* Comenzemos dando una definición elemental de lo que es un expansor. Sea $A$ un conjunto de generadores de un grupo finito $G$. Decimos que el grafo $\Gamma(G, A)$ es un $\epsilon$-expansor (para $\epsilon > 0$ dado) si todo subconjunto $S \subset G$ con $|S| \leqslant |G|/2$ satisface $|S \cup AS| \geqslant (1+\epsilon)|S|$. Es muy simple de ver que todo $\epsilon$-expansor tiene diámetro $O((\log |G|)/\epsilon)$, es decir, muy pequeño; está claro que el diámetro de $\Gamma(G, A)$ es siempre por lo menos $O((\log |G|)/(\log |A|))$.

La alternativa (al final equivalente) es definir los grafos expansores en términos del primer valor propio no trivial $\lambda_1$ del Laplaciano discreto de un grafo de Cayley. La *matriz de adyacencia* (normalizada) $\mathfrak{A}$ de un grafo es un operador linear en el espacio de funciones $f : G \to \mathbb{C}$; envía tal función a la funcíon cuyo valor en $v$ es el promedio de $f(w)$ en los vecinos $w$ de $v$. Para ser explícitos, en el caso del grafo de Cayley $\Gamma(G, A)$,

$$(\mathcal{A} f)(g) = \frac{1}{|A|} \sum_{a \in A} f(ag). \tag{5}$$

El *Laplaciano discreto* es simplemente $\triangle = I - \mathcal{A}$. (Muchos lectores lo reconocerán como el análogo de un Laplaciano sobre una superficie.)

Asumamos $A = A^{-1}$. Entonces $\triangle$ es un operador simétrico, así que todos sus valores propios son reales. Está claro que el valor propio más pequeño es $\lambda_0 = 0$, correspondiente a las funciones propias constantes. Podemos ordenar los valores propios:

$$0 = \lambda_0 \leqslant \lambda_1 \leqslant \lambda_2 \leqslant \ldots .$$

A la cantidad $|\lambda_1 - \lambda_0| = \lambda_1$ se le da el nombre de *hueco espectral*.

Decimos que $\Gamma(G, A)$ es un $\epsilon$-expansor si $\lambda_1 \geqslant \epsilon$. Para $|A|$ acotado, esta definición es equivalente a la primera que dimos (si bien la constante $\epsilon$ difiere en las dos definiciones). Decimos que una familia (conjunto infinito) de grafos $\Gamma(G, A)$ es una *familia de expansores*, o que es una familia con un hueco espectral, si todo grafo en la familia es un $\epsilon$-expansor para algún $\epsilon > 0$ fijo.

Uno de los problemas centrales del área es probar que ciertas familias (si no todas las familias) del tipo

$$\{\Gamma(\mathrm{SL}_2(\mathbb{Z}/p\mathbb{Z}), A_p)\}_{p \text{ primo}}, \qquad A_p \text{ genera } \mathrm{SL}_2(\mathbb{Z}/p\mathbb{Z})$$

son familias de expansores. Los primeros resultados atacaban el problema mediante el estudio del Laplaciano sobre las superficies $\Gamma(p)\backslash\mathbb{H}$. Un resultado clásico de Selberg (1965) nos dice que el Laplaciano sobre $\Gamma(p)\backslash\mathbb{H}$ tiene un hueco espectral independiente de $\epsilon$.

*Teoría de grupos geométrica. Teoría de modelos.* La teoría de grupos geométrica se centra en el estudio del crecimiento de $|A^k|$ para $k \to \infty$, donde $A$ es un subconjunto de un grupo infinito $G$. Por ejemplo, un teorema de Gromov (1981) muestra que, si $A$ genera a $G$ y $|A^k| \ll k^{O(1)}$, entonces $G$ tiene que ser un grupo de un tipo muy particular (virtualmente nilpotente, para ser precisos). Los argumentos de la teoría de grupos geométrica a menudo muestran que, aún si un grupo no está dado a priori de una manera geométrica, el crecimiento de un subconjunto puede darle de manera natural una geometría que puede ser utilizada.

Si bien los problemas tratados por la teoría de grupos geométrica son muy cercanos a los nuestros, tales argumentos aún no son moneda corriente en el subárea que discutiremos en estas notas, lo cual puede decir simplemente que la manera de aplicarlos aún está por descubrirse. Por otra parte, la *teoría de modelos* – en esencia, una rama de la lógica con aplicaciones a las estructuras algebraicas – ha jugado un rol directo en el subárea. Por ejemplo, Hrushovski (2012) dio una nueva prueba del teorema de Gromov, expresando cuestiones del crecimiento en grupos en un lenguaje proveniente de la teoría de modelos; más allá en esta dirección, se debe mencionar a Breuillard, Green y Tao (2012). Se trata de temas que parecen estar lejos de estar agotados.

**2.3 Resultados.** Uno de los propósitos es dar una prueba del siguiente resultado, debido al autor Helfgott (2008) para $K = \mathbb{Z}/p\mathbb{Z}$. No será idéntica a la primera prueba que diera, sino que incluirá las ideas de varios autores posteriores, incluídas algunas que han hecho que el enunciado sea más general que el original, y otras que han hecho que la prueba sea más clara y fácil de generalizar. En cualquier forma, el enunciado deriva claramente su inspiración de la combinatoria aditiva.

**Teorema 2.1.** *Sea $K$ un cuerpo. Sea $A \subset \mathrm{SL}_2(K)$ un conjunto que genera $\mathrm{SL}_2(K)$. Entonces, ya sea*

$$|A^3| \geqslant |A|^{1+\delta}$$

*o $(A \cup A^{-1} \cup \{e\})^k = \mathrm{SL}_2(K)$, donde $\delta > 0$ y $k > 0$ son constantes absolutas.*

Por cierto, gracias a Gowers (2008) y Nikolov y Pyber (2011), $(A \cup A^{-1} \cup \{e\})^k = \mathrm{SL}_2(K)$ puede reemplazarse por $A^3$. Mostraremos por lo menos que podemos tomar $k = 3$. Las primeras generalizaciones a $K$ de orden finito no primo se deben a Dinai (2011) y Varjú (2012); hoy en día, se obtiene la forma general sin mayores complicaciones.

Veamos una consecuencia sencilla.

**Ejercicio 2.1.** *Sea $G = \mathrm{SL}_2(K)$, $K$ un cuerpo finito. Sea $A \subset G$ un conjunto que genera $G$. El diámetro de $\Gamma(G, A)$ es $\ll (\log |G|)^{O(1)}$, donde las constantes implícitas son absolutas.*

Este enunciado es exactamente la conjetura de Babai para $G = \mathrm{SL}_2(K)$.

Cuándo es que $\Gamma(G, A)$ tiene diámetro $\ll \log |G|$? Yendo más lejos – cuándo es un $\epsilon$-expansor?

Se sabía desde los años 80 (ver las referencias en Helfgott (2015)) que la existencia de un agujero espectral para $\Gamma(p)\backslash\mathbb{H}$ (probada en Selberg (1965)) implica que, para

$$A_0 = \left\{ \begin{pmatrix} 1 & 1 \\ 0 & 1 \end{pmatrix}, \begin{pmatrix} 1 & 0 \\ 1 & 1 \end{pmatrix} \right\}, \tag{6}$$

los grafos $\Gamma(\mathrm{SL}_2(\mathbb{Z}/p\mathbb{Z}), A_0 \bmod p)$ forman una familia de expansores (i.e., son todos $\epsilon$-expansores para algún $\epsilon$ fijo). Empero, para, digamos,

$$A_0 = \left\{ \begin{pmatrix} 1 & 3 \\ 0 & 1 \end{pmatrix}, \begin{pmatrix} 1 & 0 \\ 3 & 1 \end{pmatrix} \right\}, \tag{7}$$

no se tenía tal resultado, ni aún una cota razonable para el diámetro. (Ésta diferencia fue resaltada varias veces por Lubotzky.)

**Ejercicio 2.2.** *Sea* $G = \mathrm{SL}_2(\mathbb{Z}/p\mathbb{Z})$*; sea* $A_0$ *como en (7). Pruebe que el diámetro de* $\Gamma(G, A_0 \bmod p)$ *es* $\ll \log |G|$.

Para resolver este ejercicio, es útil saber que $A_0$ genera un subgrupo libre de $\mathrm{SL}_2(\mathbb{Z})$. Decimos que un conjunto $A_0$ en un grupo genera un *subgrupo libre* si no existen $x_i \in A_0$, $1 \leqslant i \leqslant k$, $x_{i+1} \notin \{x_i, x_i^{-1}\}$ para $1 \leqslant i \leqslant k-1$, $x_i \neq e$ para $1 \leqslant i \leqslant k$, y $r_i \in \mathbb{Z}$, $r_i \neq 0$, tales que

$$x_1^{r_1} \cdots x_k^{r_k} = e.$$

Saber que $A_0$ genera un subgrupo libre es particularmente útil en los primeros pasos de la iteración; en los ultimos pasos, podemos utilisar el Teorema 2.1.

En verdad, la aseveración del ejercicio 2.2 sigue siendo válida sin la suposición que el grupo $\langle A_0 \rangle$ generado por $A_0$ sea libre; es suficiente (y fácil) mostrar que $\langle A_0 \rangle$ siempre tiene un subgrupo libre grande (ver el apéndice A). Sí se debe asumir que $\langle A_0 \rangle$ genera un subgrupo *Zariski-denso* de $\mathrm{SL}_2$, para así asegurar que $A_0 \bmod p$ en verdad genere $\mathrm{SL}_2(\mathbb{Z}/p\mathbb{Z})$, para $p$ mayor que una constante $C$.

Bourgain y Gamburd (2008) fueron netamente más lejos: probaron que, si $A_0$ genera un grupo Zariski-denso de $\mathrm{SL}_2$, entonces

$$\{\Gamma(\mathrm{SL}_2(\mathbb{Z}/p\mathbb{Z}), A_0 \bmod p)\}_{p > C,\, p \text{ primo}}$$

es una familia de expansores. Nos concentraremos en dar una prueba del Teorema 2.1. Al final, esbozaremos el procedimiento de Bourgain y Gamburd, basado en parte sobre dicho teorema.

Queda aún mucho por hacer; por ejemplo, no sabemos si la familia de todos los grafos

$$\{\Gamma(\mathrm{SL}_2(\mathbb{Z}/p\mathbb{Z}), A)\}_{p \text{ primo}, A \text{ genera } \mathrm{SL}_2(\mathbb{Z}/p\mathbb{Z})}$$

es una familia de expansores. Por otra parte, si bien hay generalizaciones del teorema 2.1 a otros grupos lineares (Helfgott (2011), Gill y Helfgott (2011), y, de manera más general, Breuillard, Green y Tao (2011) y Pyber y Szabó (s.f.)), aún no tenemos una prueba de la conjetura de Babai para los grupos alternantes $A_n$; la mejor cota conocida para el diámetro de $A_n$ con respecto a un conjunto arbitrario de generadores es la cota dada en Helfgott y Seress (2014), la cual no es tan buena como $\ll (\log |G|)^{O(1)}$.

## 3 Herramientas elementales

**3.1 Productos triples.** La combinatoria aditiva, al estudiar el crecimiento, estudia los conjuntos que crecen lentamente. En los grupos abelianos, sus resultados son a menudo enunciados de tal manera que clasifican los conjuntos $A$ tales que $|A^2|$ no es mucho más grande que $|A|$; en los grupos noabelianos, generalmente se clasifica los conjuntos $A$ tales que $|A^3|$ no es mucho más grande que $|A|$. Por qué?

En un grupo abeliano, si $|A^2| < K|A|$, entonces $|A^k| < K^{O(k)}|A|$ – i.e., si un conjunto no crece después de ser multiplicado por si mismo una vez, no crecera después de ser multiplicado por sí mismo muchas veces. Éste es un resultado de Plünnecke (1970) y Ruzsa (1989); Petridis (2012) dio recientemente una prueba particularmente elegante.

En un grupo no abeliano, puede haber conjuntos $A$ que rompen esta regla.

**Ejercicio 3.1.** *Sea $G$ un grupo. Sean $H < G$, $g \in G \setminus H$ y $A = H \cup \{g\}$. Entonces $|A^2| < 3|A|$, pero $A^3 \supset HgH$, y $HgH$ puede ser mucho más grande que $A$. Dé un ejemplo con $G = \mathrm{SL}_2(\mathbb{Z}/p\mathbb{Z})$.*

Empero, las ideas de Ruzsa sí se aplican al caso no abeliano, como fue indicado en Helfgott (2008) y Tao (2008); en verdad, no hay que cambiar nada en Ruzsa y Turjányi (1985), pues nunca utiliza la condición que $G$ sea abeliano. Lo que obtendremos es que basta con que $|A^3|$ (en vez de $|A^2|$) no sea mucho más grande que $|A|$ para que $|A^k|$ crezca lentamente. Veamos como se hacen las cosas.

**Lema 3.1** (Desigualdad triangular de Ruzsa). *Sean $A$, $B$ y $C$ subconjuntos finitos de un grupo $G$. Entonces*

$$|AC^{-1}||B| \leqslant |AB^{-1}||BC^{-1}|. \tag{8}$$

*Prueba.* Contruiremos una inyección $\iota : AC^{-1} \times B \hookrightarrow AB^{-1} \times BC^{-1}$. Para cada $d \in AC^{-1}$, escojamos $(f_1(d), f_2(d)) = (a, c) \in A \times C$ tal que $d = ac^{-1}$. Sea $\iota(d, b) = (f_1(d)b^{-1}, b(f_2(d))^{-1})$. Podemos recuperar $d = f_1(d)(f_2(d))^{-1}$ de $\iota(d, b)$; por lo tanto, podemos recuperar $(f_1, f_2)(d) = (a, c)$, y así también $b$. Por lo tanto, $\iota$ es una inyección. □

**Ejercicio 3.2.** *Sea $G$ un grupo. Pruebe que*

$$\frac{|(A \cup A^{-1} \cup \{e\})^3|}{|A|} \leqslant \left(3\frac{|A^3|}{|A|}\right)^3 \tag{9}$$

*para todo conjunto finito $A$ de $G$. Muestre también que, si $A = A^{-1}$ (i.e., si $g^{-1} \in A$ para todo $g \in A$), entonces*

$$\frac{|A^k|}{|A|} \leqslant \left(\frac{|A^3|}{|A|}\right)^{k-2}. \tag{10}$$

*para todo $k \geqslant 3$. Concluya que*

$$\frac{|A^k|}{|A|} \leqslant \frac{|(A \cup A^{-1} \cup \{e\})^k|}{|A|} \leqslant 3^{k-2}\left(\frac{|A^3|}{|A|}\right)^{3(k-2)} \tag{11}$$

*para todo $A \subset G$ y todo $k \geqslant 3$.*

Esto quiere decir que, de ahora en adelante, si obtenemos que $|A^k|$ no es mucho más grande que $|A|$, podemos concluir que $|A^3|$ no es mucho más grande que $|A|$. Por cierto, gracias a (9), podremos suponer en varios contextos que $e \in A$ y $A = A^{-1}$ sin pérdida de generalidad.

**3.2 El teorema de órbita-estabilizador para los conjuntos.** Una de las ideas recurrentes en la investigación del crecimiento en los grupos es la siguiente: muchos enunciados acerca de los subgrupos – así como sus métodos de prueba – pueden generalizarse a los subconjuntos. Si el método de prueba es constructivo, cuantitativo o probabilístico, esto es un indicio que la prueba podría generalizarse de tal manera.

El *teorema de órbita-estabilizador* es un buen ejemplo, tanto por su simplicidad (realmente debería llamarse "lema") como por subyacer a un número sorpredente de resultados sobre el crecimiento.

Primero, un poco de lenguaje. Una *acción* $G \curvearrowright X$ es un homomorfismo de un grupo $G$ al grupo de automorfismos de un objeto $X$. Estudiaremos el caso en el que $X$ es simplemente un conjunto; su "grupo de automorfismos" es simplemente el grupo de biyecciones de $X$ a $X$ (con la composición como operación de grupo.) Para $A \subset G$ y $x \in X$, la *órbita* $Ax$ ("órbita de $x$ bajo la acción de $A$") es el conjunto $Ax = \{g \cdot x : g \in A\}$. El *estabilizador* $\mathrm{Estab}(x) \subset G$ está dado por $\mathrm{Estab}(x) = \{g \in G : g \cdot x = x\}$.

El enunciado que daremos es como en Helfgott y Seress (2014, §3.1).

**Lema 3.2** (Teorema de órbita-estabilizador para conjuntos)**.** *Sea $G$ un grupo actuando sobre un conjunto $X$. Sea $x \in X$, y sea $A \subseteq G$ no vacío. Entonces*

$$|(A^{-1}A) \cap \mathrm{Estab}(x)| \geqslant \frac{|A|}{|Ax|} \tag{12}$$

*y, para $B \subseteq G$,*

$$|BA| \geqslant |A \cap \mathrm{Estab}(x)||Bx|. \tag{13}$$

El teorema de órbita-estabilizador usual, que se enseña usualmente en un primer curso de teoría de grupos dice que, para $H$ un subgrupo de $G$,

$$|H \cap \mathrm{Estab}(x)| = \frac{|H|}{|Hx|}.$$

Éste es un caso especial del lema que estamos por probar – el caso $A = B = H$.

**Ejercicio 3.3.** *Pruebe el Lema 3.2. Sugerencia: para (12), use el principio de los palomares.*

El grupo $G$ tiene la acción evidente "por la izquierda" sobre sí mismo: $g \in G$ actúa sobre los elementos $h \in H$ por multiplicación por la izquierda, i.e.,

$$g \mapsto (h \mapsto g \cdot h).$$

Está tambien, claro está, la acción por la derecha

$$g \mapsto (h \mapsto h \cdot g^{-1}).$$

(Por qué es que $g \mapsto (h \mapsto hg)$ no es una acción?) Ninguna de estas dos acciones son interesantes cuando se trata de aplicar directamente el Lema 3.2, pues los estabilizadores son triviales. Empero, tenemos también la acción *por conjugación*

$$g \mapsto (h \mapsto ghg^{-1}).$$

El estabilizador de un punto $h \in G$ no es sino su *centralizador* $C(h)$, definido en (2); la órbita de un punto $h \in G$ bajo la acción de todo el grupo $G$ es la *clase de conjugación*

$$\mathrm{Cl}(h) = \{ghg^{-1} : g \in G\}.$$

Así, tenemos el siguiente resultado, crucial en lo que sigue. Su importancia consiste en hacer que las cotas superiores (como las que derivaremos más tarde) sobre intersecciones con $\mathrm{Cl}(g)$ impliquen cotas inferiores sobre intersecciones con $C(g)$. La importancia de esto último es que siempre es útil saber que disponemos de muchos elementos dentro de un *variedad* (tal como un toro).

**Lema 3.3.** *Sea $A \subset G$ un conjunto no vacío. Entonces, para todo $g \in A^l$, $l \geqslant 1$,*

$$|A^{-1}A \cap C(g)| \geqslant \frac{|A|}{|A^{l+1}A^{-1} \cap \mathrm{Cl}(g)|}.$$

*Prueba.* Sea $G \curvearrowright G$ la acción de $G$ sobre sí mismo por conjugación. Aplique (12) con $x = g$; la órbita de $g$ bajo la acción de $A$ es un subconjunto de $A^{l+1}A^{-1} \cap \mathrm{Cl}(g)$. □

Es instructivo ver otras consecuencias de (12). La siguiente nos muestra, por así decirlo, que si obtenemos que la intersección de $A$ con un subgrupo $H$ de $G$ crezca, entonces hemos mostrado que $A$ mismo crece.

**Ejercicio 3.4.** *Sea $G$ un grupo y $H$ un subgrupo de $G$. Sea $A \subset G$ un conjunto no vacío tal que $A = A^{-1}$. Pruebe que, para todo $k > 0$,*

$$|A^{k+1}| \geqslant \frac{|A^k \cap H|}{|A^2 \cap H|}|A|.$$

*(Sugerencia: considere la acción $G \curvearrowright G/H$ por multiplicación por la izquierda, es decir, $g \mapsto (aH \mapsto gaH)$.)*

# 4 Intersecciones con variedades

## 4.1 Geometría algebraica extremadamente básica.

**4.1.1 Variedades.** Una *variedad* (algebraica y afín) en un espacio vectorial de $n$ dimensiones sobre un cuerpo $K$ consiste en todos los puntos $(x_1, x_2, \ldots, x_n) \in \overline{K}^n$ que satisfacen un sistema de ecuaciones

$$P_i(x_1, \ldots, x_n) = 0, \quad 1 \leqslant i \leqslant k, \tag{14}$$

donde $P_i$ son polinomios con coeficientes en $K$.

(Hay distintas maneras alternativas de formalizar el mismo concepto. Podríamos definir formalmente la variedad no exactamente como el sistema de ecuaciones en sí, sino como el conjunto de todas las ecuaciones polinomiales implicadas por el sistema en (14). También hay definiciones de apariencia mucho más abstractas, basada en la teoría de *esquemas* (Grothendieck), pero no necesitaremos entrar allí.)

Se dice generalmente que los puntos $(x_1, \ldots, x_n)$ que satisfacen (14) *yacen* sobre la variedad, así como hablamos de puntos que yacen sobre una curva o superficie algebraica; claro está, las curvas y las superficies son casos especiales de variedades. Dada una variedad $V$ definida sobre $K$ y un cuerpo $L$ tal que $K \subset L \subset \overline{K}$, escribimos $V(L)$ por el conjunto de todos los puntos $(x_1, x_2, \ldots, x_n) \in L^n$ que yacen sobre $V$.

El caso trivial es el de la variedad $\mathbb{A}^n$ (*espacio afín*) definida por el sistema vacío de ecuaciones (o por la ecuación $P(x_1, \ldots, x_n) = 0$, donde $P$ es el polinomio 0). Claramente, $\mathbb{A}^n(L) = L^n$.

Dadas dos variedades $V_1$, $V_2$, tanto $V_1 \cap V_2$ como $V_1 \cup V_2$ son variedades: la variedad $V_1 \cap V_2$ está dada por la unión de las ecuaciones que definen $V_1$ y aquellas que definen $V_2$, mientras que, si $V_1$ está definida por (14) y $V_2$ está definida por

$$Q_j(x_1, \ldots, x_n) = 0, \quad 1 \leqslant j \leqslant k', \tag{15}$$

donde $Q_i$ son polinomios con coeficientes en $K$, entonces la variedad $V_1 \cup V_2$ está dada por las ecuaciones

$$(P_i \cdot Q_j)(x_1, \ldots, x_n) = 0, \quad 1 \leqslant i \leqslant k, \quad 1 \leqslant i' \leqslant k'.$$

Consideremos ahora los grupos lineares, como $\mathrm{SL}_2(K)$. Está claro que $\mathrm{SL}_2(K)$ está contenido en

$$M_2(\overline{K}) = \left\{ \begin{pmatrix} x_1 & x_2 \\ x_3 & x_4 \end{pmatrix} : x_1, x_2, x_3, x_4 \in \overline{K} \right\},$$

el cual es un espacio vectorial (de dimensión 4) sobre $\overline{K}$. Por lo tanto, tiene sentido hablar de variedades $V$ en $M_2$. Por ejemplo, tenemos la variedad $V$ de elementos que tienen una traza dada:

$$\operatorname{tr} \begin{pmatrix} x_1 & x_2 \\ x_3 & x_4 \end{pmatrix} = C, \quad \text{i.e.,} \quad x_1 + x_4 - C = 0. \tag{16}$$

Nuestro grupo $\mathrm{SL}_2$ es también una variedad, dada por la ecuación $x_1x_4 - x_2x_3 = 1$. Es así un ejemplo de un *grupo algebraico*. Estrictamente hablando, son los puntos $\mathrm{SL}_2(K)$ (o $\mathrm{SL}_2(L)$) del grupo algebraico $\mathrm{SL}_2$ los que forman un grupo en el sentido usual del término. (Claro está, la operación de grupo $\cdot : \mathrm{SL}_2 \times \mathrm{SL}_2 \to \mathrm{SL}_2$ está bien definida como un *morfismo de variedades afines* – una aplicación de una variedad a otra dado por polinomios.)

Es fácil ver que un toro máximo $T = C(g)$ también es un grupo algebraico, puesto que la ecuación

$$hg = gh$$

es un sistema de ecuaciones polinomiales (lineares, en verdad) sobre los coeficientes de $h \in \mathrm{SL}_2$. Así, $T$ es un subgrupo algebraico de $G$.

**4.1.2 Dimensión, grado, intersecciones.** Dada una variedad $V$, una *subvariedad* $W$ es una variedad contenida en él. Una variedad $V$ se dice irreducible si no es la unión de dos subvariedades no vacías $W, W' \subsetneq V$. Por ejemplo, $\mathrm{SL}_2$ es *irreducible*; ésto es básicamente una consecuencia del hecho que el polinomio $x_1x_4 - x_2x_3 - 1$ es irreducible.

Podemos definir la *dimensión* de un variedad irreducible $V$ como el entero $d$ máximo tal que haya una cadena de variedades irreducibles no vacías

$$V_0 \subsetneq V_1 \subsetneq V_2 \subsetneq \ldots \subsetneq V_d = V.$$

Ésto coincide con el concepto intuitivo de *dimensión*: un plano contiene una línea, que contiene a un punto, y así un plano es de dimensión por lo menos 2; en verdad, es de dimensión exactamente 2.

*Hecho.* La dimensión de $\mathbb{A}^n$ es $n$.

En particular, la dimensión de una variedad en $\mathbb{A}^n$ es siempre finita ($\leqslant n$).

**Ejercicio 4.1.** *Probemos que la intersección $\bigcap_{i\in I} V_i$ de una colección finita o infinita de variedades $V_i \in \mathbb{A}^n$, $i \in I$, es una variedad.*

*Para una colección finita, ésto es evidente. Por ende, para el caso infinito, bastará si mostramos que hay un $S \subset I$ finito tal que $\bigcap_{i\in I} V_i = \bigcap_{i\in S} V_i$. Muestre que esto se reduce a mostrar que una cadena de variedades*

$$\mathbb{A}^n \supsetneq W_1 \supsetneq W_2 \supsetneq W_3 \supsetneq \dots$$

*debe ser finita.* (Esta propiedad se llama *propiedad Noetheriana*.)

*Reduzca esto a su vez al caso de $W_1$ irreducible. Concluya por inducción en la dimensión de $W_1$.*

**Ejercicio 4.2.** *(a) Sea $V$ una variedad. Pruebe que se puede expresar $V$ como una unión finita de variedades irreducibles $W_1, W_2, \dots, W_k$. (Pista: use la propiedad Noetheriana.)*

*(b) Sea $V$ irreducible. Muestre que, si $V$ es una unión finita de variedades*

$$V_1, V_2, \dots, V_k,$$

*entonces existe un $1 \leqslant i \leqslant k$ tal que $V_i = V$.*

*(c) Sea $V$ una variedad. Muestre que, si imponemos la condición que $W_i \not\subset W_j$ para $i \neq j$, la descomposición $V = W_1 \cup W_2 \cup \dots \cup W_k$ en variedades irreducibles $W_i$ es única (excepto que, claro, los $W_i$ pueden permutarse). Las variedades $W_i$ se llaman* componentes irreducibles *de $V$.*

Si una variedad es una unión de variedades irreducibles todas de dimensión $d$, decimos que es “de dimensión pura”, y podemos decir que es de dimensión $d$.

Dadas dos variedades irreducibles $W \subset V$, la *codimensión* $\operatorname{codim}(W)$ de $W$ en $V$ es simplemente $\dim(V) - \dim(W)$. Es fácil ver que, para $V$ irreducible, o bien $W = V$, o bien $\operatorname{codim}(W) > 0$. Si $\operatorname{codim}(W) > 0$ (o si $W$ es una unión de variedades irreducibles de codimensión positiva), podemos pensar en los elementos de $W(\overline{K})$ como *especiales*, y en los elementos de $V(\overline{K})$ que no están en $W(\overline{K})$ como *genéricos*.

Es posible (y muy recomendable) considerar variedades más generales que las variedades afines. Por ejemplo, podemos considerar las *variedades proyectivas*, definidas por sistemas de ecuaciones $P(x_0, x_1, \dots, x_n) = 0$ donde cada $P$ es un polinomio homogéneo en $n+1$ variables. Los puntos en un variedad proyectiva viven en el *espacio proyectivo* $\mathbb{P}^n$. Los puntos del espacio proyectivo sobre un campo $L$ son elementos de $L^{n+1}$ (excepto $(0, 0, \dots, 0)$), donde se identifica dos elementos $x, x' \in L^{n+1}$ si $x$ es un múltiplo escalar de $x'$, i.e., si $x = \lambda x'$ para algún $\lambda \in L$. En otras palabras,

$$\mathbb{P}^n(L) = (L^{n+1} \setminus \{(0, 0, \dots, 0)\})/\sim, \quad \text{donde } x \sim x' \text{ si } \exists \lambda \in L \text{ t.q. } x = \lambda x'.$$

La opción de trabajar con variedades proyectivas nos da mucha libertad: en particular, es posible mostrar que podemos hablar de la variedad proyectiva de todas las líneas (o todos los planos) en el espacio $n$-dimensional (proyectivo). El ejemplo más simple es la variedad de todas las líneas en el plano: como una línea en el plano proyectivo ($n = 2$) está dada por una ecuación linear homogénea

$$c_0x_0 + c_1x_1 + c_2x_2 = 0,$$

y como dos tales ecuaciones dan la misma línea si sus triples $(c_1, c_2, c_3)$ son múltiples el uno del otro, tenemos que las líneas en el plano proyectivo están en correspondencia uno-a-uno con $\mathbb{P}^2$ mismo. En $\mathbb{P}^n$, como en $\mathbb{A}^n$, podemos hablar de subvariedades, codimensión, elementos genéricos. Podemos hacer una inmersíon de $\mathbb{A}^n$ en $\mathbb{P}^n$:

$$(x_1, x_2, \ldots, x_n) \mapsto (1, x_1, x_2, \ldots, x_n).$$

El complemento es la subvariedad de $\mathbb{P}^n$ dada por $x_0 = 0$; para $n = 2$, se le llama *recta en el infinito*.

Contemplemos ahora una curva irreducible $C$ en el plano, es decir, una subvariedad de $\mathbb{A}^2$ de dimensión 1. (En verdad, no es necesaria la irreducibilidad; la supondremos realmente sólo al trabajar con análogos en dimensiones superiores.) Consideremos también una línea $\ell$ en $\mathbb{A}^2$. Podría ser que $\ell$ fuera tangente a $C$, pero no es difícil demostrar que tal es el caso sólo cuando $\ell$ yace en una subvariedad de $\mathbb{P}^2$ (la *curva dual a* $C$). En otras palabras, una línea genérica (se dice también: *en posición general*) no es tangente a $C$. El número de puntos de intersección en $\mathbb{P}^2(\overline{K})$ de la curva $C$ con una línea genérica $\ell$ resulta ser independiente de $\ell$; llamamos a ese número el *grado* de $C$.

Resulta ser que el grado de una curva irreducible en el plano dada por una ecuación

$$P(x_1, x_2) = 0$$

(o por una ecuación $P(x_0, x_1, x_2) = 0$, $P$ homogéneo) es simplemente el grado de $P$. La ventaja de la definición que dimos del grado de una curva es que es más conceptual y se generaliza de manera natural. En dimensiones superiores, la misma variedad puede ser definida por distintos sistemas de ecuaciones de grados distintos; deseamos una definición de *grado* que dependa sólo de la variedad, y no del sistema que la define.

Si $V$ es una variedad de dimensión 2 en $\mathbb{A}^3$, definimos su grado como su número de intersecciones con una recta genérica; si $V$ es de dimensión 1 en $\mathbb{A}^3$, definimos su grado como su número de intersecciones con un *plano* genérico. Así como hablamos de líneas y planos, podemos definir, en general, una *variedad linear* mediante ecuaciones lineares. Para $V$ una variedad irreducible en $\mathbb{A}^n$ de dimensión $d$, definimos el *grado* $\deg(V)$ de $V$ como el número de intersecciones de $V$ con una variedad linear genérica de codimensión $d$ en $\mathbb{A}^n$. La misma definición es válida cuando $V$ es no necesariamente irreducible pero de dimensión pura e igual a $d$.

Si bien este grado, como decíamos, no tiene porque corresponder al grado de ninguna de las ecuaciones en un sistema de ecuaciones que defina a $V$, puede ser acotado por una constante que depende sólo de los grados de tales ecuaciones y su número. Ese es un caso especial de lo que estamos por discutir.

El *teorema de Bézout*, en el plano, nos dice que, para dos curvas irreducibles distintas $C_1$, $C_2$ en $\mathbb{A}^2$, el número de puntos de la intersección $(C_1 \cap C_2)(\overline{K})$ es a lo más $d_1d_2$. (En verdad, si consideramos $C_1$ y $C_2$ genéricos, o trabajamos en $\mathbb{P}^2$ y contamos "multiplicidades", el número de puntos de intersección es *exactamente* $d_1 d_2$.)

En general, si $V_1$ y $V_2$ son variedades irreducibles, y escribimos $V_1 \cap V_2$ como una unión de variedades irreducibles $W_1, W_2, \ldots, W_k$, con $W_i \not\subset W_j$ para $i \neq j$, una generalización del teorema de Bézout nos dice que

$$\sum_{i=1}^{k} \deg(W_k) \leqslant \deg(V_1)\deg(V_2). \tag{17}$$

(Véase, por ejemplo, Danilov y Shokurov (1998, p.251), donde se menciona a Fulton y MacPherson en conexión a (17) y enunciados más generales.) Toda variedad irreducible de dimensión 0 consiste en un único punto; por ello (17) implica el teorema de Bézout habitual.

**4.2 Escape de subvariedades.** Sea $G$ un grupo que actúa por transformaciones lineares sobre el espacio $n$-dimensional $K^n$, $K$ un cuerpo. (En otras palabras, se nos es dado un homomorfismo $\phi : G \to \mathrm{GL}_n(K)$ de $G$ al grupo de matrices invertibles $\mathrm{GL}_n(K)$.) Sea $W$ una variedad de codimensión positiva en $\mathbb{A}^n$. Estábamos llamando a los elementos de $W(\overline{K})$ *especiales*, y a los otros elementos de $\mathbb{A}^n(\overline{K})$ *genéricos*.

Sea $A$ un conjunto de generadores de $G$ y $x$ un punto de $W$. Muy bien podría ser que la órbita $A \cdot x$ esté contenida por entero en $W$. Empero, como veremos ahora, si $Gx$ no está contenida en $W$, entonces siempre es posible *escapar* de $W$ en un número acotado de pasos: no sólo que habrá (por definición) algún producto $g$ de un número finito de elementos de $A$ y $A^{-1}$ tal que $g \cdot x$ está fuera de $W$, sino que habrá un producto (a decir verdad, muchos productos) $g \in (A \cup A^{-1})^k$, *$k$ acotado*, tal que $g \cdot x$ está fuera de $W$. En otras palabras, si escapamos por lo menos una vez, eventualmente, de $W$, escapamos de muchas maneras de $W$, después de un número acotado $k$ de pasos.

La prueba[2] procede por inducción en la dimensión, controlando el grado.

**Proposición 4.1.** *Sean dados:*

- *$G$, un grupo actuando por transformaciones lineares sobre $K^n$, $K$ un cuerpo;*
- *$W \subsetneq \mathbb{A}^n$, una variedad,*
- *un conjunto de generadores $A \subset G$;*
- *un elemento $x \in \mathbb{A}^n(K)$ tal que $G \cdot x$ no está contenido en $W$.*

*Entonces hay constantes $k$, $c$ que dependen sólo del número, dimensión y grado de los componentes irreducibles de $W$, tales que hay por lo menos* $\max(1, c|A|)$ *elementos $g \in (A \cup A^{-1} \cup \{e\})^k$ tales que $gx \notin W(K)$.*

Para aclarar el proceso de inducción, daremos primero la prueba en un caso particular. Una variedad *linear* es simplemente una línea, un espacio, etc.; en otras palabras, es una variedad definida por ecuaciones lineares.

[2]El enunciado de la proposición es cómo en Helfgott (2011), basado en Eskin, Mozes y Oh (2005), pero la idea es probablemente más antigua.

*Prueba para $W$ linear e irreducible.* Sea $W$ linear e irreducible. Podemos asumir sin pérdida de generalidad que $A = A^{-1}$ y $e \in A$.

Procederemos por inducción en la dimensión de $W$. Si $\dim(W) = 0$, entonces $W$ consiste en un sólo punto $x_0$, y el enunciado que queremos probar es cierto: existe un $g \in A$ tal que $gx \neq x_0$ (por qué?); si hay menos de $|A|/2$ tales elementos, escogemos un $g_0 \in A$ tal que $g_0 x_0 \neq x_0$ (por qué existe?), y entonces, para cada uno de los más de $|A|/2$ elementos $g \in A$ tales que $gx = x_0$, tenemos que $g_0 g x = g_0 x_0 \neq x_0$.

Asumamos, entonces, que $\dim(W) > 0$, y que el enunciado ha sido probado para todas las variedades lineares irreducibles $W'$ con $\dim(W') < \dim(W)$. Si $gW = W$ para todo $g \in A$, entonces ya sea (a) $gx \notin W(K)$ para todo $g \in A$, y el enunciado es inmediato, o (b) $gx \in W(K)$ para todo $g \in G$ (puesto que $G$ está generado por $A$), lo cual está en contradicción con nuestras suposiciones. Podemos asumir, entonces, que $gW \neq W$ para algun $g \in A$.

Entonces $W' = gW \cap W$ es una variedad linear irreducible de dimensión $\dim(W') < \dim(W)$. Por lo tanto, por la hipótesis inductiva, hay $\geqslant \max(1, c'|A|)$ elementos $g'$ de $A^{k'}$ (donde $c'$ y $k'$ dependen sólo de $\dim(W)$) tales que $g'x$ no yace en la variedad $W' = gW \cap W$. Entonces, para cada tal $g'$, ya sea $g^{-1}g'x$ o $g'x$ no yace en $W$. Así, hemos probado la proposicion con $c = c'/2$, $k = k' + 1$. □

**Ejercicio 4.3.** *Generalize la prueba que acabamos de dar de tal manera que dé Prop. 4.1 para $W$ arbitrario.* Sugerencia: *como un primer paso, generalice la prueba de tal manera que funcione para toda unión $W$ de variedades lineares irreducibles. (Ésto ya exigirá adaptar el proceso de inducción de tal manera que se controle de alguna manera el número de componentes en cada paso. Claro está, la intersección de dos uniones $W$ de $d$ variedades lineares irreducibles tiene a lo más $d^2$ componentes.) Luego muestre que la prueba es válida para toda union $W$ de variedades irreducibles, no necesariamente lineares, utilizando la generalización (17) del teorema de Bézout.*

**4.3 Estimaciones dimensionales.** Dado un conjunto de generadores $A \subset \mathrm{SL}_2(K)$ (o $A \subset \mathrm{SL}_n(K)$, o lo que se desee) y una subvariedad $V$ de codimensión positiva en $\mathrm{SL}_2$, sabemos que una proporción positiva de los elementos de $A^k$, $k$ acotado, yacen fuera de $V$: éste es un caso especial de la Proposición 4.1 (con $x$ igual a la identidad $e$).

Aunque esto desde ya implica una cota superior para el número de elementos de $A^k$ en $V(K)$, podemos dar una cota mucho mejor. Los estimados de este tipo pueden trazarse en parte a Larsen y Pink (2011) (caso de $A$ un subgrupo, $V$ general) y en parte a Helfgott (2008) y Helfgott (2011) ($A$ un conjunto en general, pero $V$ especial). Tales cotas tienen en general la forma

$$|A \cap V(K)| \ll |(A \cup A^{-1} \cup \{e\})^k|^{\frac{\dim V}{\dim G}}. \tag{18}$$

Se lograron cotas completamente generales del tipo (18) en Breuillard, Green y Tao (2011) y Pyber y Szabó (s.f.) ($A$ y $V$ arbitrarios, $G$ un grupo linear algebraico simple, como en Larsen y Pink (2011)).

Como primero paso hacia la estrategia general, veamos un caso particular de manera muy concreta (aunque no lo usemos al final). La prueba es básicamente la misma que en Helfgott (2008, §4).

**Lema 4.1.** *Sea $G = \mathrm{SL}_2$, $K$ un cuerpo, y $T$ un toro máximo. Sea $A \subset G(K)$ un conjunto de generadores de $G(K)$. Entonces*

$$|A \cap T(K)| \ll |(A \cup A^{-1} \cup \{e\})^k|^{1/3} \tag{19}$$

*donde k y la constante implícita son constantes absolutas.*

*Prueba.* Podemos suponer sin pérdida de generalidad que $|K|$ es mayor que una constante, pues, de lo contrario, la conclusión es trivial. También podemos suponer sin pérdida de generalidad que $A = A^{-1}$, $e \in A$, y que $|A|$ es mayor que una constante, reemplazando $A$ por $(A \cup A^{-1} \cup \{e\})^c$, $c$ constante, de ser necesario. Podemos también escribir los elementos de $T$ como matrices diagonales, conjugando por un elemento de $\mathrm{SL}_2(\overline{K})$.

Sea

$$g = \begin{pmatrix} a & b \\ c & d \end{pmatrix} \tag{20}$$

un elemento cualquiera de $\mathrm{SL}_2(\overline{K})$ con $abcd \neq 0$. Consideremos la aplicación $\phi : T(K) \times T(K) \times T(K) \to G(K)$ dada por

$$\phi(x, y, z) = x \cdot gyg^{-1} \cdot z.$$

Queremos mostrar que esta aplicación es en algún sentido cercana a ser inyectiva. (La razón de tal estrategia? Si la aplicación fuera inyectiva, y tuvieramos $g \in A^\ell$, $\ell$ una constante, entonces tendríamos

$$|A \cap T(K)|^3 = |\phi(A \cap T(K), A \cap T(K), A \cap T(K))| \leqslant |AA^\ell AA^{-\ell}A| = |A^{2\ell+3}|,$$

lo cual implicaría inmediatamente el resultado que queremos. Simplemente estamos usando el hecho que el tamaño de la imagen $\phi(D)$ de una inyección $\phi$ tiene el mismo número de elementos que su dominio $D$.)

Multiplicando matrices, vemos que, para

$$x = \begin{pmatrix} r & 0 \\ 0 & r^{-1} \end{pmatrix}, \quad y = \begin{pmatrix} s & 0 \\ 0 & s^{-1} \end{pmatrix}, \quad z = \begin{pmatrix} t & 0 \\ 0 & t^{-1} \end{pmatrix},$$

$\phi((x, y, z))$ es igual a

$$\begin{pmatrix} rt(sad - s^{-1}bc) & rt^{-1}(s^{-1} - s)ab \\ r^{-1}t(s - s^{-1})cd & r^{-1}t^{-1}(s^{-1}ad - sbc) \end{pmatrix}. \tag{21}$$

Sea $s \in \overline{K}$ tal que $s^{-1} - s \neq 0$ y $sad - s^{-1}bc \neq 0$. Un breve cálculo muestra que entonces $\phi^{-1}(\{\phi((x, y, z))\})$ tiene a lo más 16 elementos: tenemos que

$$rt^{-1}(s^{-1} - s)ab \cdot r^{-1}t(s - s^{-1})cd = -(s - s^{-1})^2 abcd,$$

y, como $abcd \neq 0$, hay a lo más 4 valores de $s$ dado un valor de $-(s - s^{-1})^2 abcd$ (el producto de las esquinas superior derecha e inferior izquierda de (21)); para cada tal valor de $s$, el producto y el cociente de las esquinas superior izquierda y superior derecha de (21) determinan $r^2$ y $t^2$, respectivamente, y obviamente hay sólo 2 valores de $r$ y 2 valores de $t$ para $r^2$ y $t^2$ dados.

Ahora bien, hay a lo más 4 valores de $s$ tales que $s^{-1} - s = 0$ o $sad - s^{-1}bc = 0$. Por lo tanto, tenemos que

$$|\phi(A \cap T(K), A \cap T(K), A \cap T(K))| \geqslant \frac{1}{16}|A \cap T(K)|(|A \cap T(K)| - 4)|A \cap T(K)|,$$

y, como antes, $\phi(A \cap T(K), A \cap T(K), A \cap T(K)) \subset AA^{\ell}AA^{-\ell}A = A^{2\ell+3}$. Si $|A \cap T(K)|$ es menor que 8 (o cualquier otra constante) entonces la conclusión (19) es trivial. Por lo tanto, concluimos que

$$|A \cap T(K)|^3 \leqslant 2|A \cap T(K)|(|A \cap T(K)| - 4)|A \cap T(K)| \leqslant 32|A^{2\ell+3}|,$$

i.e., (19) es cierta.

Sólo queda verificar que existe un elemento (20) de $A^{\ell}$ con $abcd \neq 0$. Ahora bien, $abcd = 0$ define una subvariedad $W$ de $\mathbb{A}^4 \sim M_2$; más aún, para $|K| > 2$, existen elementos de $G(K)$ fuera de tal variedad. Por lo tanto, las condiciones de Prop. 4.1 se cumplen (con $x$ igual a la identidad $e$). Así, obtenemos que existe $g \in A^{\ell}$ ($\ell$ una constante) tal que $g \notin W(K)$, lo cual era lo que necesitábamos. □

Hagamos abstracción de lo que acabamos de hacer, para así poder generalizar el resultado a una variedad arbitraria $V$ en vez de $T$. Trataremos el caso de $V$ de dimensión 1, por conveniencia. La estrategia de la prueba del Lema 4.1 consiste en construir un morfismo $\phi : V \times V \times \cdots \times V \to G$ ($r$ copias de $V$, donde $r = \dim(G)$) de la forma

$$\phi(v_1, \ldots, v_r) = v_1 g_1 v_2 g_2 \cdots v_{r-1} g_{r-1} v_r, \tag{22}$$

donde $g_1, g_2, \ldots, g_{r-1} \in A^{\ell}$, y mostrar que, para $v = (v_1, \ldots, v_r)$ genérico (es decir, fuera de una subvariedad de $V \times \cdots \times V$ de codimensión positiva), la preimagen $\phi^{-1}(\phi(v))$ tiene dimensión 0. En verdad, como acabamos de ver, es suficiente mostrar que esto es cierto para $(g_1, g_2, \ldots, g_{r-1})$ un elemento genérico de $G^{r-1}$; el argumento de escape (Prop. 4.1) se encarga del resto.

Para hacer que el argumento marche para $V$ general (y $G$ general), es necesario asumir algunos fundamentos. Esencialmente, tenemos la elección de ya sea trabajar sobre el álgebra de tipo Lie o introducir un poco más de geometría algebraica. La primera elección (tomada en Helfgott (2015), siguiendo a Helfgott (2011)) asume que el lector tiene cierta familiaridad con los grupos y álgebras de Lie, y que sabe, o está dispuesto a creer, que la relación entre grupos y álgebras de Lie sigue la misma si trabajamos sobre un cuerpo finito en vez de $\mathbb{R}$ o $\mathbb{C}$. La segunda elección – que tomaremos aquí – requiere saber, o estar listo a aceptar, un par de hechos básicos sobre morfismos, válidos sobre cuerpos arbitrarios.

No importa gran cosa si se sigue el uno u el otro formalismo. Los fundamentos, en uno y otro caso, se sentaron sólidamente en la primera mitad del siglo XX (Zariski, Chevalley, etc.) y son relativamente accesibles. Los lectores que sientan interés en estudiar las bases del camino que seguiremos están invitados a leer Mumford (1999, Ch. 1) (de por sí una excelente idea) o cualquier texto similar.

Está claro que, si $\phi : \mathbb{A}^n \to \mathbb{A}^m$ es un morfismo y $V \subset \mathbb{A}^m$ es una variedad, entonces la preimagen $\phi^{-1}(V)$ es una variedad (por qué?). Algo nada evidente que utilizaremos es el hecho que, si $\phi$ es como dijimos y $V \subset \mathbb{A}^n$ es una variedad, entonces $\phi(V)$ es un *conjunto construíble*, lo cual quiere decir una unión finita de términos de la forma $W \setminus W'$, donde $W$ y $W' \subset W$ son variedades. (Por ejemplo, si $V \subset \mathbb{A}^2$ es la variedad dada por $x_1 x_2 = 1$ (una hipérbola), entonces su imagen bajo el morfismo $\phi(x_1, x_2) = x_1$ es el conjunto construíble $\mathbb{A}^1 \setminus \{0\}$.) Éste es un teorema de Chevalley Mumford (ibíd., §I.8, Cor. 2); encapsula parte del campo clásico llamado *teoría de la eliminación*. Es fácil deducir, que, en general, para $V$ construíble, $\phi(V)$ es construíble.

Siempre podemos expresar un conjunto construíble $S$ como una unión $\cup_i (W_i \setminus W_i')$ con $\dim(W_i') < \dim(W_i)$. (Por qué?) La *clausura de Zariski* $\overline{S}$ del conjunto construíble $S$ es entonces $\cup_i W_i$.

El siguiente lema es Tao (2015, Prop. 1.5.30), lo cual es a su vez en esencia Larsen y Pink (2011, Lemma 4.5). (Murmullo de un mundo paralelo: en el formalismo que no seguimos, esto corresponde al hecho, básico pero no trivial, que el álgebra de un grupo de tipo Lie simple es simple.)

Decimos que un grupo algebraico $G$ es *casi simple* si no tiene ningún subgrupo algebraico normal $H$ de dimensión positiva y menor que $\dim(G)$. (Por ejemplo, $SL_n$ es casi simple para todo $n \geqslant 2$.)

**Lema 4.2.** *Sea $G \subset SL_n$ un grupo algebraico irreducible y casi simple definido sobre un cuerpo $K$. Sean $V', V \subsetneq G$ subvariedades con* $\dim(V') > 0$. *Entonces, para todo $g \in G(\overline{K})$ fuera de una subvariedad $W \subsetneq G$, algún componente de la clausura de Zariski $\overline{V'gV}$ tiene dimensión* $> \dim(V)$.

*Más aún, el número de componentes de $W$ y sus grados están acotados por una constante que depende sólo de $n$ y del número y grados de los componentes de $V'$ y $V$.*

*Prueba.* Podemos suponer sin pérdida de generalidad que $V$ y $V'$ son irreducibles, y que $e \in V'(\overline{K})$.

Sea $g \in G(\overline{K})$. Supongamos que $\overline{V'gV}$ tiene dimensión $\leqslant \dim(V)$. Para todo $v' \in V'(\overline{K})$, $v'gV$ es una variedad de dimensión $\dim(V)$ – uno del número finito de componentes $W_i$ de $\overline{V'gV}$ de dimensión $\dim(V)$. Para cada tal $W_i$, los puntos $v'$ tales que $v'gV = W_i$ forman una variedad $V_i$, al ser la intersección de variedades

$$\bigcap_{v \in V(\overline{K})} W_i v^{-1} g^{-1}.$$

Así, $V'$ es la unión de un número finito de variedades $V_i$; como $V'$ es irreducible, esto implica que $V' = V_i$ para algún $V$. En particular, $gV = e \cdot gV = W_i$. Por lo tanto, $V'gV = gV$.

Ahora bien, los elementos $g \in G(\overline{K})$ tales que $V'gV = gV$ son una intersección

$$\bigcap_{\substack{v \in V \\ v' \in V'}} \phi_{v',v}^{-1}(V)$$

de variedades $\phi_{v',v}^{-1}(V)$, donde $\phi_{v',v}(g) = g^{-1}v'gv$. Por lo tanto, tales elementos constituyen una subvariedad $W$ de $G$; más aún, gracias a Bézout (17), su número de componentes así como el grado de estos están acotados por una constante que depende sólo de $n$ y del número y grados de los componentes de $V'$ y $V$.

Falta sólo mostrar que $W \neq G$. Supongamos que $W = G$. Entonces $V'gV = gV$ para todo $g \in G(\overline{K})$. Ahora bien, el estabilizador $\{g \in G(\overline{K}) : gV = V\}$ de $V$ no sólo es un grupo, sino que es (el conjunto de puntos de) una variedad (nuevamente: para ver esto, expréselo como una intersección de variedades). Llamemos a tal variedad $\mathrm{Estab}(V)$. Tenemos que $g^{-1}V'g \subset \mathrm{Estab}(V)$ para todo $g \in G(\overline{K})$, y por lo tanto

$$V' \subset \bigcap_{g \in G(\overline{K})} g\,\mathrm{Estab}(V) g^{-1}.$$

Esto muestra que la variedad $\cap_{g\in G(\overline{K})} g\,\mathrm{Estab}(V)g^{-1}$ es de dimensión $\geqslant \dim(V') > 0$. Al mismo tiempo, dicha variedad es un subgrupo algebraico normal de $G$, contenido en $\mathrm{Estab}(V)$. Como $\mathrm{Estab}(V) \subsetneq G$, tenemos un subgrupo algebraico normal de $G$, de dimensión positiva y estrictamente contenido en $G$. En otras palabras, $G$ no es un grupo algebraico casi simple. Contradicción. □

Sean $G$, $V$ y $V'$ tales que satisfagan las hipótesis del Lema 4.2, y sea $g \in G(\overline{K})$ como en la conclusión del Lema, i.e., fuera de la subvariedad $W \subsetneq G$. Asumamos que $\dim(V') = 1$, y consideremos el morfismo $\phi : V' \times V \to \overline{V'gV}$ dado por

$$\phi(v', v) = v'gv.$$

La dimensión de la imagen de un morfismo no es mayor que la dimension de su dominio (ejercicio), así que

$$\dim(\phi(V' \times V)) \leqslant \dim(V' \times V) = \dim(V') + \dim(V) = \dim(V) + 1.$$

Al mismo tiempo, por el Lema 4.2, $\dim(\phi(V', V)) > \dim(V)$. Por lo tanto,

$$\phi(V' \times V) = \dim(V) + 1 = \dim(V' \times V).$$

Si un morfismo $\phi : X \to X'$ es tal que $\overline{\phi(X)} = X'$, decimos que $\phi$ es *dominante*. (Por ejemplo, el morfismo $\phi$ que acabamos de considerar es dominante.) Aceptemos el hecho que, si $\phi$ es dominante y $\dim(X') = \dim(X)$, entonces hay una subvariedad $Y \subsetneq X$ tal que, para todo $x \in X(\overline{K})$ que no yazca en $Y(\overline{K})$, la variedad $\phi^{-1}(\phi(\{x\}))$ tiene dimensión 0. (Ésta es una consecuencia inmediata del Mumford (1999, §1.8, Thm. 3).) Más aún, el número de componentes de $Y$ y su grado están acotados en términos del grado de $\phi$ y del número, dimensión y grado de los componentes de $X$ y $X'$.

(Para que lo que acabamos de llamar una consecuencia inmediata de algo en otra parte se vuelva intuitivamente claro, considere el caso $K = \mathbb{R}$. Entonces $Y$ es la subvariedad de $X$ definida por la condición “la determinante $D\phi(x)$ de $\phi$ en el punto $x$ tiene determinante 0”. Hay varias maneras de ver que $Y$ es una subvariedad $\subsetneq X$ para $K$ arbitrario: como dijimos, es posible definir derivadas sobre cuerpos arbitrarios, o, alternativamente, proceder como en Mumford (ibíd., §1.8, Thm. 3).)

Aplicando esto a la aplicación $\phi$ que teníamos, obtenemos que hay una subvariedad $Y \subsetneq V' \times V$ tal que, para todo $x \in (V' \times V)(\overline{K})$ que no yace en $Y$, $\phi^{-1}(\phi(x))$ tiene dimensión 0.

Dado ésto, podemos probar una generalización del Lema 4.1. Se trata realmente de (18) para toda variedad de dimensión 1.

**Proposición 4.2.** *Sea $K$ un cuerpo y $G \subset \mathrm{SL}_n$ un grupo algebraico casi simple tal que $|G(K)| \geqslant c|K|^{\dim(G)}$, $c > 0$. Sea $Z \subset G$ una variedad de dimensión 1. Sea $A \subset G(K)$ un conjunto de generadores de $G(K)$. Entonces*

$$|A \cap Z(K)| \ll |(A \cup A^{-1} \cup \{e\})^k|^{1/\dim(G)}, \tag{23}$$

*donde $k$ y la constante implícita dependen solamente de $n$, de $c$ y del número y grado de los componentes irreducibles de $G$ y $Z$.*

Obviamente, $G = \mathrm{SL}_n$ es una elección válida, pues es casi simple y $|\mathrm{SL}_n(K)| \gg |K|^{n^2-1} = |K|^{\dim(G)}$.

**Ejercicio 4.4.** *Pruebe la Proposición 4.2. He aquí un esbozo:*

*(a) Muestre el siguiente lema básico: si $W \subset \mathbb{A}^N$ es una variedad de dimensión $d$, entonces el número de puntos $(x_1, \ldots, x_N) \in \mathbb{A}^n(K)$ que yacen en $W$ es $\ll |K|^d$, donde la constante implícita depende sólo de $N$ y del número y grado de componentes irreducibles de $W$. (Sugerencia: para $d = 0$, ésto está claro. Para $d > 0$, considere la proyección $\pi : \mathbb{A}^N \to \mathbb{A}^{N-1}$ a las primeras $N-1$ coordenadas, o más bien dicho la restricción $\pi|W$ de $\pi$ a $W$. Reduzca al caso de dimensión $d-1$ – la manera de hacerlo depende de si $\pi|W$ es o no es dominante.)*

*(b) Utilizando el escape de subvariedades (Prop. 4.1) y el Lema 4.2, muestre que, dadas las condiciones del Lema 4.2, existe un elemento $g \in (A \cup A^{-1} \cup \{e\})^{\ell}$, $\ell$ una constante (dependiendo de esto y aquello), tal que algún componente de la clausura de Zariski $\overline{V'gV}$ tiene dimensión $> \dim(V)$. Esto es rutina, pero no olvide mostrar que hay algun punto de $G(K)$ fuera de $W$ (usando el lema básico que acaba de probar).*

*(c) Aplicando esto (y las consecuencias discutidas inmediatamente despues del Lema 4.2) de manera iterada, muestre que existen $g_1, \ldots, g_{r-1} \in (A \cup A^{-1} \cup \{e\})^{\ell'}$ y una subvariedad $Y \subsetneq Z \times Z \times \ldots \times Z$ ($r = \dim(G)$ veces) tales que, para todo $x \in (Z \times Z \times \ldots \times Z)(\overline{K})$ que no yace en $Y$, $\phi^{-1}(\phi(x))$ es de dimensión $0$, donde $\phi$ es como en (22).*

*(d) Usando nuevamente un argumento que distingue si una projección (esta vez de $Z \times \ldots \times Z$ ($r$ veces) a $Z \times \ldots \times Z$ ($r-1$ veces)) es dominante, e iterando, muestre que hay a lo más $O(|A \cap Z(K)|^{r-1})$ elementos de $(A \cup Z(K)) \times \ldots \times (A \cup Z(K))$ ($r$ veces) en $Y$.*

*(e) Concluya que la proposición Prop. 4.2 es cierta.*

En general, se puede probar (18) para $\dim(V)$ arbitrario siguiendo argumentos muy similares, mezclados con una inducción sobre la dimensión de la variedad $V$ en (18). Ilustraremos el proceso básico haciendo las cosas en detalle para $G = \mathrm{SL}_2$ y para el tipo de variedad $V$ que realmente necesitamos.

Se trata de la variedad $V_t$ definida por

$$\det(g) = 1, \mathrm{tr}(g) = t \tag{24}$$

para $t \neq \pm 2$. Tales variedades nos interesan por el hecho que, para cualquier $g \in \mathrm{SL}_2(K)$ regular semisimple (lo cual en $\mathrm{SL}_2$ quiere decir: con dos valores propios distintos), la clase de conjugacion $\mathrm{Cl}(g)$ está contenida en $V_{\mathrm{tr}(g)}$.

**Proposición 4.3.** *Sea $K$ un cuerpo; sea $A \subset \mathrm{SL}_2(K)$ un conjunto de generadores de $\mathrm{SL}_2(K)$. Sea $V_t$ dada por (24). Entonces, para todo $t \in K$ aparte de $\pm 2$,*

$$|A \cap V_t(K)| \ll |(A \cup A^{-1} \cup \{e\})^k|^{\frac{2}{3}}, \tag{25}$$

*donde $k$ y la constante implícita son constantes absolutas.*

Claro está, $\dim(\mathrm{SL}_2) = 3$ y $\dim(V_t) = 2$, así que este es un caso particular de (18).

*Prueba.* Consideremos la aplicación $\phi : V_t(K) \times V_t(K) \to \mathrm{SL}_2(K)$ definida por

$$\phi(y_1, y_2) = y_1 y_2^{-1}.$$

Está claro que

$$\phi(A \cap V_t(K), A \cap V_t(K)) \subset A^2.$$

Así, si $\phi$ fuera inyectiva, tendríamos inmediatamente que $|A \cap V_t(K)|^2 \leqslant |A^2|$. Ahora bien, $\phi$ no es inyectiva. La preimagen de $\{h\}$, $h \in \mathrm{SL}_2(K)$, es

$$\phi^{-1}(\{h\}) = \{(w, h^{-1}w) : \mathrm{tr}(w) = t, tr(h^{-1}w) = t\}.$$

Debemos preguntarnos, entonces, cuántos elementos de $A$ yacen en la subvariedad $Z_{t,h}$ de $G$ definida por

$$Z_{t,h} = \{(w, hw) : \mathrm{tr}(w) = t, tr(h^{-1}w) = t\}.$$

Para $h \neq \pm e$, $\dim(Z_{t,h}) = 1$ (verificar), y el número y grado de componentes de $Z_{t,h}$ esta acotado por una constante absoluta. Así, aplicando la Proposición 4.2, obtenemos que, para $h \neq \pm e$,

$$|A \cap Z_{t,h}(K)| \ll |A^{k'}|^{1/3},$$

donde $k'$ y la constante implícita son absolutas.

Ahora bien, para cada $y_1 \in V_t(K)$, hay por lo menos $|V_t(K)| - 2$ elementos $y_2 \in V_t(K)$ tales que $y_1 y_2^{-1} \neq \pm e$. Concluímos que

$$|A \cap V(K)|(|A \cap V(K)| - 2) \leqslant |A^2| \cdot \max_{g \neq \pm e} |A \cap Z_{t,h}(K)| \ll |A^2||A^{k'}|^{1/3}.$$

Podemos asumir que $|A \cap V(K)| \geqslant 3$, pues de lo contrario la conclusión deseada es trivial. Obtenemos, entonces, que

$$|A \cap V(K)| \ll |A^k|^{2/3}$$

para $k = \max(2, k')$, como queríamos. □

Pasemos a la consecuencia que nos interesa.

**Corolario 4.1.** *Sea $K$ un cuerpo y $G = \mathrm{SL}_2$. Sea $A$ un conjunto de generadores de $G(K)$; sea $g \in A^\ell$ ($\ell \geqslant 1$) regular semisimple. Entonces*

$$|A^{-1}A \cap C(g)| \gg \frac{|A|}{|(A \cup A^{-1} \cup \{e\})^{k\ell}|^{2/3}}, \tag{26}$$

*donde $k$ y la constante implícita son absolutas.*

*En particular, si $|A^3| \leqslant |A|^{1+\delta}$, entonces*

$$|A^{-1}A \cap C(g)| \gg_\ell |A|^{1/3 - O(\delta\ell)}. \tag{27}$$

*Prueba.* La Proposición 4.3 y el Lema 3.3 implican (26) inmediatamente, puesto que $\mathrm{Cl}(g) \subset V_{\mathrm{tr}(g)}$, donde $V_t$ se define como en (24). De manera también muy sencilla, la conclusión (27) se deduce de (26) a través de (11). □

Veamos ahora dos problemas cuyos resultados no utilizaremos; son esenciales, empero, si se quiere trabajar en $\mathrm{SL}_n$ para $n$ arbitrario. El primer problema es relativamente ambicioso, pero ya hemos visto todos los elementos esenciales para su solución. En esencia, sólo se trata de saber organizar la recursión.

**Ejercicio 4.5.** *Generalice la Proposicion 4.2 a $Z$ de dimensión arbitraria.*

En general, un elemento $g \in \mathrm{SL}_n(K)$ es *regular semisimple* si tiene $n$ valores propios distintos. Claro está, todo elemento de $C(g)$ tiene los mismos vectores propios que $g$. Cuando $G = \mathrm{SL}_n$, como para $\mathrm{SL}_2$, los elementos de $C(g)$ son los puntos $T(K)$ de un subgrupo algebraico abeliano $T$ de $G$, llamado un toro máximo. Tenemos que $\dim(T) = n-1$ y $\dim(\overline{\mathrm{Cl}(g)}) = \dim(G) - \dim(T)$.

**Ejercicio 4.6.** *Generalice 4.1 a $G = \mathrm{SL}_n$, para $g$ semisimple. En vez de (27), la conclusión reza como sigue:*

$$|A^{-1}A \cap C(g)| \gg |A|^{\frac{\dim(T)}{\dim(G)} - O(\delta)}, \tag{28}$$

*donde las constantes implícitas dependen solo de $n$.*

Terminemos por una breve nota con un lado anecdótico. Una versión del Corolario 4.1 fue probada en Helfgott (2008), donde jugó un rol central. Luego fue generalizada a $\mathrm{SL}_n$ en Helfgott (2011), dando, en esencia, (28).

Empero, estas versiones tenian una debilidad: daban (27) y (28) para la mayoría de los $g \in A^{\ell}$, y no para *todo* $g \in A^{\ell}$. Esto hacía que el resto del argumento – la parte que estamos por ver – fuera más complicado y difícil de generalizar que lo es hoy en día.

La moraleja es, por supuesto, que no hay que asumir que las técnicas y argumentos que a uno le son familiares son óptimos – y que para simplificar una prueba vale la pena tratar de probar resultados intermedios más fuertes.

# 5 El crecimiento en $\mathrm{SL}_2(K)$

**5.1 El caso de los subconjuntos grandes.** Veamos primero que pasa con $A \cdot A \cdot A$ cuando $A \subset \mathrm{SL}_2(\mathbb{F}_q)$ es grande con respecto a $G = \mathrm{SL}_2(\mathbb{F}_q)$. En verdad no es difícil mostrar que, si $|A| \geqslant |G|^{1-\delta}$, $\delta > 0$ suficientemente pequeño, entonces $(A \cup A^{-1} \cup \{e\})^k = G$, donde $k$ es una constante absoluta. Probaremos algo más fuerte: $A^3 = G$. La prueba se debe a Nikolov y Pyber (2011); está basada sobre una idea clásica, desarrollada en este contexto por Gowers (2008). Nos dará la oportunidad de revisitar el tema de los valores propios de la matriz de adyacencia $\mathcal{A}$ de $\Gamma(G, A)$. (Los comenzamos a discutir en §2.2.)

Primero, recordemos que una *representación compleja* de un grupo $G$ es un homomorphismo $\phi : G \to \mathrm{GL}_d(\mathbb{C})$; decimos, naturalmente, que $d \geqslant 1$ es la *dimensión* de la representación. Una representación $\phi$ es *trivial* si $\phi(g) = e$ para todo $g \in G$.

El siguiente resultado se debe a Frobenius (1896), por lo menos para $q$ primo. Se puede mostrar simplemente examinando una tabla de caracteres, como en Shalom (1999) (que da también análogos de esta proposición para otros grupos de tipo Lie). Para $q$ primo, hay una prueba breve y elegante; véase, e.g., Tao (2015, Lemma 1.3.3).

**Proposición 5.1.** *Sea* $G = \mathrm{SL}_2(\mathbb{F}_q)$, $q = p^\alpha$. *Entonces toda representación compleja no trivial de* $G$ *tiene dimensión* $\geqslant (q-1)/2$.

Ahora bien, para cada valor propio $\nu$ de $\mathscr{A}$, podemos considerar su *espacio propio* – el espacio vectorial que consiste en todas las funciones propias $f : G \to \mathbb{C}$ con valor propio $\nu$. Como puede verse de la definición de $\mathscr{A}$ (inmediatamente después de (5)), tal espacio es invariante bajo la acción de $G$ por multiplicación por la derecha. En otras palabras, es una representación de $G$ - y puede ser trivial sólo si se trata del espacio (uni-dimensional) que consiste de las funciones constantes, i.e., el espacio propio que corresponde al valor propio $\nu_0 = 1$. Por lo tanto, todo los otros valores propios tienen multiplicidad $\geqslant (q-1)/2$. Asumamos, como es nuestra costumbre, que $A = A^{-1}$, lo cual implica que todos los valores propios son reales:

$$\ldots \leqslant \nu_2 \leqslant \nu_1 \leqslant \nu_0 = 1.$$

La idea es ahora es obtener un hueco espectral, i.e., una cota superior para $\nu_j$, $j > 0$. Es muy común usar el hecho que la traza de una potencia $\mathscr{A}^r$ de una matriz de adyacencia puede expresarse de dos maneras: como el número (normalizado por el factor $1/|A|^r$, en nuestro caso) de ciclos de longitud $r$ en el grafo $\Gamma(G, A)$, por una parte, y como la suma de potencias $r$-ésimas de los valores propios de $\mathscr{A}$, por otra. En nuestro caso, para $r = 2$, esto nos da

$$\frac{|G||A|}{|A|^2} = \sum_j \nu_j^2 \geqslant \frac{q-1}{2}\nu_j^2, \tag{29}$$

para cualquier $j \geqslant 1$, y, por lo tanto,

$$|\nu_j| \leqslant \sqrt{\frac{|G|/|A|}{(q-1)/2}}. \tag{30}$$

Ésta es una cota superior muy pequeña para $|A|$ grande. Esto quiere decir que unas cuantas aplicaciones de $\mathscr{A}$ bastan para hacer que una función se "uniformice", i.e., se vuelva casi constante, pues cualquier componente ortogonal al espacio propio de funciones constantes es multiplicado por algún $\nu_j$, $j \geqslant 1$, en cada paso. La prueba siguiente simplemente aplica esta observación.

**Proposición 5.2** (Nikolov y Pyber (2011)). *Sea* $G = \mathrm{SL}_2(\mathbb{F}_q)$, $q = p^\alpha$. *Sea* $A \subset G$, $A = A^{-1}$. *Asumamos* $|A| \geqslant 2|G|^{8/9}$. *Entonces*

$$A^3 = G.$$

La suposicion $A = A^{-1}$ es en verdad innecesaria, gracias al trabajo adicional puesto en Gowers (2008) para el caso no simétrico.

*Prueba.* Supongamos $g \in G$ such that $g \notin A^3$. Entonces el producto escalar

$$\langle \mathscr{A}1_A, 1_{gA}\rangle = \sum_{x\in G} (\mathscr{A}1_A)(x) \cdot 1_{gA^{-1}}(x)$$

es igual a 0. Podemos asumir que los vectores propios $v_j$ satisfacen $\langle v_j, v_j\rangle = 1$. Entonces

$$\begin{aligned}\langle \mathscr{A}1_A, 1_{gA}\rangle &= \langle \sum_{j\geqslant 0} \nu_j \langle 1_A, v_j\rangle v_j, 1_{gA}\rangle \\ &= \nu_0\langle 1_A, v_0\rangle\langle v_0, 1_{gA^{-1}}\rangle + \sum_{j>0} \nu_j \langle 1_A, v_j\rangle\langle v_j, 1_{gA^{-1}}\rangle.\end{aligned}$$

Ahora bien, $v_0$ es la función constante, y, al satisfacer $\langle v_0, v_0\rangle = 1$, es igual a $1/\sqrt{|G|}$. Luego

$$\nu_0\langle 1_A, v_0\rangle\langle v_0, 1_{gA^{-1}}\rangle = 1 \cdot \frac{|A|}{\sqrt{|G|}} \cdot \frac{|g^{-1}A|}{\sqrt{|G|}} = \frac{|A|^2}{|G|}.$$

Al mismo tiempo, gracias a (30) y Cauchy–Schwarz,

$$\left|\sum_{j>0} \nu_j \langle 1_A, v_j\rangle\langle v_j, 1_{gA^{-1}}\rangle\right| \leqslant \sqrt{\frac{2|G|/|A|}{q-1}} \sqrt{\sum_{j\geqslant 1} |\langle 1_A, v_j\rangle|^2} \sqrt{\sum_{j\geqslant 1} |\langle v_j, 1_{gA^{-1}}\rangle|^2}$$
$$\leqslant \sqrt{\frac{2|G|/|A|}{q-1}} |1_A|_2 |1_{gA^{-1}}|_2 = \sqrt{\frac{2|G||A|}{q-1}}.$$

Como $|G| = (q^2-q)q$, tenemos que $|A| \geqslant 2|G|^{8/9}$ implica que

$$\frac{|A|^2}{|G|} > \sqrt{\frac{2|G||A|}{q-1}},$$

y por lo tanto $\langle \mathcal{A}1_A, 1_{gA^{-1}}\rangle$ es mayor que 0. Contradicción. □

**5.2 El crecimiento en $\mathrm{SL}_2(K)$, $K$ arbitrario.** Probemos finalmente el teorema 2.1. En esta parte nos acercaremos más a tratamientos nuevos (en particular, Pyber y Szabó (s.f.)) que al tratamiento original en Helfgott (2008); estos tratamientos nuevos se generalizan más fácilmente. Si bien sólo deseamos presentar una prueba para $\mathrm{SL}_2$, notaremos el punto o dos en la prueba dónde hay que trabajar un poco a la hora de generalizarla para $\mathrm{SL}_n$.

La primera prueba de este teorema en la literatura utilizaba el *teorema de la suma y producto*, un resultado no trivial de combinatoria aditiva. La prueba que daremos no lo utiliza, pero sí tiene algo en común con su prueba: la inducción, usada de una manera particular. En esencia, si algo es cierto para el paso $n$, pero no para el paso $n+1$, se trata de usar ese mismo hecho para obtener la conclusión que deseamos de otra manera (lo que se llama un "fulcro" (*pivot*) en la prueba que estamos por ver). El hecho que estemos en un grupo sin un orden natural ($n, n+1$, etc.) resulta ser irrelevante.

*Prueba del Teorema 2.1.* Gracias a (9), podemos asumir que $A = A^{-1}$ y $e \in A$. También podemos asumir que $|A|$ es mayor que una constante absoluta, pues de lo contrario la conclusión es trivial. Escribamos $G = \mathrm{SL}_2$.

Supongamos que $|A^3| < |A|^{1+\delta}$, donde $\delta > 0$ es una pequeña constante a ser determinada más tarde. Por escape (Prop. 4.1), existe un elemento $g_0 \in A^c$ regular semisimple (esto es, $\mathrm{tr}(g_0) \neq \pm 2$), donde $c$ es una constante absoluta. (A decir verdad, $c = 2$; ejercicio opcional.) Su centralizador en $G(K)$ es $C(g) = T(\overline{K}) \cap G(K)$ para algún toro maximal $T$.

Llamemos a $\xi \in G(K)$ un *fulcro* si la función $\phi_g : A \times C(g) \to G(K)$ definida por

$$(a, t) \mapsto a\xi t\xi^{-1} \tag{31}$$

es inyectiva en tanto que función de $\pm e \cdot A/\{\pm e\} \times C(g)/\{\pm e\}$ a $G(K)/\{\pm e\}$.

*Caso (a): Hay un fulcro $\xi$ en $A$.* Por el Corolario 4.1, existen $\gg |A|^{1/3-O(c\delta)}$ elementos de $C(g)$ en $A^2$. Por lo tanto, por la inyectividad de $\phi_\xi$,

$$\left|\phi_\xi(A, A^2 \cap C(g))\right| \geqslant \frac{1}{4}|A||A^2 \cap C(g)| \gg |A|^{\frac{4}{3}-O(c\delta)}.$$

Al mismo tiempo, $\phi_\xi(A, A^2 \cap C(g)) \subset A^5$, y por lo tanto

$$|A^5| \gg |A|^{4/3-O(c\delta)}.$$

Para $|A|$ mayor que una constante y $\delta > 0$ menor que una constante, esto nos da una contradicción con $|A^3| < |A|^{1+\delta}$ (por Ruzsa (10)).

*Caso (b): No hay fulcros $\xi$ en $G(K)$.* Entonces, para todo $\xi \in G(K)$, hay $a_1, a_2 \in A$, $t_1, t_2 \in T(K)$, $(a_1, t_1) \neq (\pm a_2, \pm t_2)$ tales que $a_1 \xi t_1 \xi^{-1} = \pm e \cdot a_2 \xi t_2 \xi^{-1}$, lo cual da

$$a_2^{-1} a_1 = \pm e \cdot \xi t_2 t_1^{-1} \xi^{-1}.$$

En otras palabras, para cada $\xi \in G(K)$, $A^{-1}A$ tiene una intersección no trivial con el toro $\xi T \xi^{-1}$:

$$A^{-1}A \cap \xi T(K) \xi^{-1} \neq \{\pm e\}. \tag{32}$$

(Por cierto, esto sólo es posible si $K$ es un cuerpo finito $\mathbb{F}_q$. Por qué?)

Escoja cualquier $g \in A^{-1}A \cap \xi T(K)\xi^{-1}$ con $g \neq \pm e$. Entonces $g$ es regular semisimple (nota: esto es peculiar a $SL_2$) y su centralizador $C(g)$ es igual a $\xi T(K) \xi^{-1}$ (por qué?). Por lo tanto, por el corolario 4.1, obtenemos que hay $\geqslant c'|A|^{1/3-O(\delta)}$ elementos de $\xi T(K)\xi^{-1}$ en $A^2$, donde $c'$ y la constante implícita son absolutas.

Por lo menos $(1/2)|G(K)|/|T(K)|$ toros máximos de $G$ son de la forma $\xi T \xi^{-1}$, $\xi \in G(K)$ (demostrar!). También tenemos que todo elemento de $G$ que no sea $\pm e$ puede estar en a lo más un toro máximo (de nuevo algo peculiar a $SL_2$). Por lo tanto,

$$|A^2| \geqslant \frac{1}{2}\frac{|G(K)|}{|T(K)|}(c'|A|^{1/3-O(\delta)} - 2) \gg q^2 |A|^{1/3-O(\delta)}.$$

Por lo tanto, ya sea $|A^2| > |A|^{1+\delta}$ (en contradicción con la suposición que $|A^3| \leqslant |A|^{1+\delta}$) o $|A| \geqslant |G|^{1-O(\delta)}$. En el segundo caso, la proposición 5.2 implica que $A^3 = G$.

*Caso (c): Hay elementos de $G(K)$ que son fulcros y otros que no lo son.* Como $\langle A \rangle = G(K)$, esto implica que existe un $\xi \in G$ que no es un fulcro y un $a \in A$ tal que $a\xi \in G$ sí es un fulcro. Como $\xi$ no es un fulcro, (32) es cierto, y por lo tanto hay $|A|^{1/3-O(\delta)}$ elementos de $\xi T \xi^{-1}$ en $A^k$.

Al mismo tiempo, $a\xi$ es un fulcro, i.e., la aplicación $\phi_{a\xi}$ definida en (31) es inyectiva (considerada como una aplicación de $A/\{\pm e\} \times C(g)/\{\pm e\}$ a $G(K)/\{\pm e\}$). Por lo tanto,

$$\left|\phi_{a\xi}(A, \xi^{-1}(A^k \cap \xi T \xi^{-1})\xi)\right| \geqslant \frac{1}{4}|A||A^k \cap \xi T \xi^{-1}| \geqslant \frac{1}{4}|A|^{\frac{4}{3}-O(\delta)}.$$

Como $\phi_{a\xi}(A, \xi^{-1}(A^k \cap \xi T \xi^{-1})\xi) \subset A^{k+3}$, obtenemos que

$$|A^{k+3}| \geqslant \frac{1}{4}|A|^{4/3-O(\delta)}. \tag{33}$$

Gracias otra vez a Ruzsa (10), esto contradice $|A^3| \leqslant |A|^{1+\delta}$ para $\delta$ suficientemente pequeño. □

## A Expansión en $\mathrm{SL}_2(\mathbb{Z}/p\mathbb{Z})$

Daremos aquí un esbozo de cómo Bourgain y Gamburd probaron que, para $A_0 \subset \mathrm{SL}_2(\mathbb{Z})$ tal que $\langle A_0 \rangle$ es Zariski-denso, entonces

$$\{\Gamma(\mathrm{SL}_2(\mathbb{Z}/p\mathbb{Z}), A_0 \bmod p)\}_{p > C,\, p \text{ primo}}$$

es una familia de expansores, i.e., tiene un hueco espectral constante.

Primero, clarifiquemos que quiere decir "Zariski-denso". Esto quiere decir simplemente que no existe ninguna subvariedad $V \subsetneq \mathrm{SL}_2(\mathbb{C})$ que contenga a $\langle A_0 \rangle$. Como dijimos en (2.3), es un hecho conocido que esto implica que $A_0 \mod p$ genera $\mathrm{SL}_2(\mathbb{Z}/p\mathbb{Z})$ para $p$ mayor que una constante $C$ (Weisfeiler (1984), Matthews, Vaserstein y Weisfeiler (1984), Nori (1987) y Hrushovski y Pillay (1995) lo prueban para $\mathrm{SL}_2$ y para muchos grupos más).

Es sencillo pasar a un subgrupo libre:

$$\Gamma(2) = \{g \in \mathrm{SL}_2(\mathbb{Z}) : g \equiv I \mod 2\}$$

es libre, y es un resultado estándar (Nielsen–Schreier) que todo subgrupo de un grupo libre es libre. Por lo tanto $\langle A_0 \rangle \cap \Gamma(2)$ es libre. El índice de $\langle A_0 \rangle \cap \Gamma(2)$ en $\langle A_0 \rangle$ es finito (por qué?), y podemos encontrar un conjunto finito que genera $\langle A_0 \rangle \cap \Gamma(2)$ (generadores de Schreier, por ejemplo). Esto es suficiente para que podamos asumir, sin pérdida de generalidad, que $\langle A_0 \rangle$ es libre.

Lo que ahora haremos es considerar la función

$$\mu(x) = \begin{cases} \frac{1}{|A_0 \bmod p|} & \text{si } x \in A_0 \bmod p, \\ 0 & \text{si } x \notin A_0 \bmod p \end{cases}$$

y sus convoluciones. La convolución $f \cdot g$ de dos funciones $f, g : G \to \mathbb{C}$ se define por

$$(f \cdot g)(x) = \sum_{y \in G} f(xy^{-1}) g(y).$$

La norma $\ell_p$ de una función $f : G \to \mathbb{C}$ es

$$|f|_p = \left( \sum_{y \in G} |f(y)|^p \right)^{1/p}.$$

Es fácil ver que la convolución $\mu^{(\ell)} := \mu * \mu * \ldots * \mu$ ($\ell$ veces) tiene norma $\ell_1$ igual a 1. Empero, la norma $\ell_2$ varía. Para todo $f$, $|f * \mu|_2 \leqslant |f|_2$, por Cauchy–Schwarz, con igualdad sólo si $f$ es uniforme, es decir, constante (ejercicio). Por lo tanto, $|\mu^{(\ell)}|_2$ decrece cuando $\ell$ aumenta.

Nos interesa saber que tan rápido decrece, pues ésto nos da información sobre los valores propios de $\mathcal{A}$. Veamos por qué. El operador $\mathcal{A}$ no es sino la convolución por $\mu$. Podemos comparar, como en (29), dos expresiones para la traza. Por una parte, la traza de $\mathcal{A}^{2\ell}$ es igual a la suma, para todo $g$, del número de maneras de ir de $g$ a $g$ tomando productos por $A$ exactamente $2\ell$ veces, dividido por $|A_0|^{2\ell}$; esto es

$$|G| \mu^{(2\ell)}(e) = |G| \sum_{x \in G} \mu^{(\ell)}(x^{-1}) \mu^{(\ell)}(x) = |G| |\mu^{(\ell)}|_2^2,$$

donde $G = \mathrm{SL}_2(\mathbb{Z}/p\mathbb{Z})$. (Como de costumbre, asumimos que $A_0 = A_0^{-1}$.) Por otra parte, la traza de $\mathcal{A}^{2\ell}$ es igual a $\sum_i \nu_i^{2\ell}$, donde $1 = \nu_0 > \nu_1 \geqslant \ldots$ son los valores propios de $\mathcal{A}$.

Como ya discutimos en §5.1, todo valor propio $\nu_j$, $j \geqslant 1$, tiene multiplicidad $\geqslant (p-1)/2$. Por lo tanto, tenemos que, para todo $j \geqslant 1$,

$$\frac{p-1}{2}\nu_j^{2\ell} \leqslant \sum_{j\geqslant 0} \nu_j^{2\ell} = |G||\mu^{(\ell)}|_2^2.$$

Nuestra meta será mostrar que, para algún $\ell \leqslant C \log p$, $C$ una constante suficientemente grande, la función $\mu^{(\ell)}$ es razonablemente uniforme, o "llana", por lo menos del punto de vista de su norma $\ell_2$: $|\mu^{(\ell)}|_2^2 \ll 1/|G|^{1-\epsilon}$. (La distribución uniforme tiene norma $\ell_2$ igual a $1/|G|$, naturalmente.) Entonces tendremos

$$\nu_j^{2\ell} \ll \frac{|G|^\epsilon}{p} \ll \frac{1}{p^{1-3\epsilon}}$$

(puesto que $|G| \ll p^3$) y por lo tanto

$$\nu_j \leqslant e^{-\frac{(1-3\epsilon)\log p}{C \log p}} \leqslant 1 - \delta,$$

donde $\delta > 0$ es una constante (que, como $C$, puede depender de $A_0$). Esto es lo que deseamos.

(El uso de la multiplicidad de $\nu_j$ en este contexto particular remonta a Sarnak y Xue (1991).)

Lo que queda es, como decíamos, mostrar que $|\mu^{(\ell)}|$ decrece rápidamente cuando $\ell$ aumenta. Ésto esta estrechamente ligado a mostrar que $|A_0^\ell|$ decrece (en particular, lo implica), pero no es trivialmente equivalente.

La prueba tiene dos pasos. Primero, igual que para $|A_0^\ell|$ (ver el ejercicio 2.2 y los comentarios que lo siguen), está el caso de lo que pasa para $\ell \leqslant \epsilon' \log p$, donde $\epsilon'$ es lo suficientemente pequeño como para que, para elementos $g_1, g_2, \ldots, g_{2\ell} \in A_0 = A_0 \cup A_0^{-1}$ cualesquiera, tengamos que ninguno de los coeficientes de la matriz $g_1 g_2 \ldots g_{2\ell} \in \mathrm{SL}_2(\mathbb{Z})$ tenga valor absoluto $\geqslant p-1$. Entonces, tenemos que no existen $x_i \in A_0 \bmod p$, $1 \leqslant i \leqslant k$, $x_{i+1} \notin \{x_i, x_i^{-1}\}$ para $1 \leqslant i \leqslant k-1$, $x_i \neq e$ para $1 \leqslant i \leqslant k$, y $r_i \in \mathbb{Z}$, $r_i \neq 0$, $\sum_{1\leqslant i\leqslant k} |r_i| \leqslant 2\ell$, tales que

$$x_1^{r_1} \cdots x_k^{r_k} = e.$$

(Idea: si un elemento de $\mathrm{SL}_2(\mathbb{Z})$ es congruente $\bmod p$ a la identidad sin ser la identidad: entonces por lo menos uno de sus coeficientes de matriz tiene valor absoluto por lo menos $p-1$.)

Esto implica inmediatamente que los productos de elementos de $A_0 \bmod p$ de longitud $\ell$ son todos diferentes (excepto por las igualdades obvias del tipo $x \cdot e = x$ y $x \cdot x^{-1} = e$). Por lo tanto, $|(A \bmod p)^\ell|$ crece exponencialmente: $|(A_0 \bmod p)^\ell| \geqslant (|A_0| - 2)^\ell$. Ésta era la parte crucial a solución del ejercicio 2.2. Mostrar que $|\mu^\ell|_2$ decrece también exponencialmente no es mucho más difícil, sobre todo porque podemos asumir que $A_0$ es más grande que una constante. (Para $A_0$ más pequeño que una constante, sería un asunto más delicado: se trata de un resultado clásico de Kesten (1959) sobre los grupos libres.)

Queda por ver como decrece $|\mu^\ell|_2$ para $\epsilon' \log p \leqslant \ell \leqslant C \log p$. Aquí que Bourgain y Gamburd muestran que, si tuvieramos

$$|\mu^{2\ell}|_2 > |\mu^\ell|_2^{1+\delta'},$$

$\delta' > 0$, entonces existe un conjunto $A' \subset \mathrm{SL}_2(\mathbb{Z}/p\mathbb{Z})$ tal que $|A'^3| < |A'|^{1+O(\delta')}$. (La herramienta principal es el teorema de Balog y Szemerédi (1994), fortalecido por Gowers (2001) y generalizado por Tao (2008) al caso no conmutativo.) Muestran también que $\mu(A')$ es grande, por lo cual $A'$ es menor que $|G|^{1-O(\delta')}$ a menos que $\mu$ ya sea tan uniforme como deseamos. Un argumento auxiliar muestra que $A'$ genera $\mathrm{SL}_2(\mathbb{Z}/p\mathbb{Z})$. Por lo tanto, $|A'^3| < |A'|^{1+O(\delta')}$ entra en contradicción con el Teorema 2.1.

Esto muestra que $|\mu^{2\ell}|_2 \leqslant |\mu^{\ell}|_2^{1+\delta'}$ para $\epsilon' \log p \leqslant \ell \leqslant C \log p$, y termina la prueba. Concluimos que $\nu_1 \leqslant 1 - \delta$, que era lo que queríamos demostrar.

## Referencias

Harald Andrés Helfgott
helfgott@math.univ-paris-diderot.fr
Universität Göttingen / CNRS / Université de Paris VI/VII

# EXPANSORES GEOMÉTRICOS

Mikhail Belolipetsky

**Resumen**

En esta parte vamos a estudiar las propiedades geométricas y topológicas de los espacios asociados a las secuencias de los subgrupos de congruencia. En particular, consideraremos algunos resultados de P. Buser, P. Sarnak, R. Brooks, M. Lackenby y M. Gromov en el contexto de expansores geométricos.

## 1 El plano y el espacio hiperbólico

Vamos a usar principalmente el modelo del semiespacio superior para el *espacio hiperbólico*. En este modelo tenemos:

$$\mathbb{H}^3 = \{(x, y, t) \in \mathbb{R}^3 \mid t > 0\}$$

con la métrica inducida por

$$ds^2 = \frac{dx^2 + dy^2 + dt^2}{t^2}.$$

Con esta métrica, $\mathbb{H}^3$ se convierte en una variedad riemanniana completa con curvatura constante $-1$. Tal variedad, que es además conexa y simplemente conexa, es única salvo isometrías.

Del mismo modo, el *plano hiperbólico* puede ser representado por el semiplano superior

$$\mathbb{H}^2 = \{(x, t) \in \mathbb{R}^2 \mid t > 0\}$$

con la métrica inducida por $ds^2 = (dx^2 + dt^2)t^{-2}$.

Observar que esta métrica en $\mathbb{H}^2$ es la restricción de la métrica hiperbólica en $\mathbb{H}^3$ al plano $y = 0$.

El grupo $\mathrm{PSL}_2(\mathbb{C})$ es el cociente del grupo $\mathrm{SL}_2(\mathbb{C})$ de todas las matrices $2 \times 2$ con entradas complejas y determinante 1 por el centro $\{\pm\mathrm{Id}\}$:

$$\mathrm{PSL}_2(\mathbb{C}) = \mathrm{SL}_2(\mathbb{C})/\{\pm\mathrm{Id}\}.$$

Nosotros vamos a considerar frecuentemente elementos de un subgrupo $\Gamma$ de $\mathrm{PSL}_2(\mathbb{C})$ como matrices, e ignorar la diferencia entre $\mathrm{PSL}_2(\mathbb{C})$ y $\mathrm{SL}_2(\mathbb{C})$.

Con la acción inducida por la transformación de Möbius

$$z \to \frac{az + b}{cz + d},$$

los elementos $\gamma = \begin{pmatrix} a & b \\ c & d \end{pmatrix}$ de $\mathrm{PSL}_2(\mathbb{C})$ actúan sobre la esfera de Riemann $\widehat{\mathbb{C}} = \mathbb{C} \cup \{\infty\}$ y esta acción es biholomorfa.

La acción de cada $\gamma \in \mathrm{PSL}_2(\mathbb{C})$ se extiende a $\mathbb{H}^3$ con la extensión de Poincaré:

> $\gamma$ = un producto del número par de las inversiones en los círculos y líneas en $\mathbb{C}$. (**Ejercicio:** Verificar este hecho.)
> Ahora considera $\widehat{\mathbb{C}}$ como la frontera $t = 0$ de $\mathbb{H}^3$ — la esfera en el infinito. Podemos extender cada línea/círculo en $\mathbb{C}$ a un plano/hemisferio en $\mathbb{H}^3$ ortogonal a $\mathbb{C}$. El producto de las correspondientes reflexiones en los planos y las inversiones en hemisferios da lugar a la extensión de $\gamma$ para $\mathbb{H}^3$.

Podemos verificar que las reflexiones e inversiones son isometrías de la métrica hiperbólica en $\mathbb{H}^3$ y que ellas generan el grupo de isometrías $\mathrm{Isom}(\mathbb{H}^3)$. Por eso con la extensión de Poincaré tenemos

$$\mathrm{PSL}_2(\mathbb{C}) \cong \mathrm{Isom}^+(\mathbb{H}^3),$$

el grupo de isometrías de $\mathbb{H}^3$ que preservan la orientación.

Para isometrías del plano hiperbólico $\mathbb{H}^2$ tenemos también el isomorfismo

$$\mathrm{PSL}_2(\mathbb{R}) \cong \mathrm{Isom}^+(\mathbb{H}^2).$$

Las isometrías de $\mathbb{H}^3$ actúan transitivamente en los hemisferios y planos ortogonales a $\mathbb{C}$. Así que podemos enviar $\mathbb{H}^2 = \{y = 0 \text{ en } \mathbb{H}^3\}$ para cualquier otro plano de este tipo o hemisferio. Por tanto todos los hemisferios y planos en $\mathbb{H}^3$, los cuales son ortogonales a $\mathbb{C}$ con la restricción de la métrica hiperbólica de $\mathbb{H}^3$, son modelos de $\mathbb{H}^2$. Ellos dan el conjunto de *planos geodésicos* en $\mathbb{H}^3$. Si dos de estos planos se cruzan en $\mathbb{H}^3$, podemos definir el ángulo diedral entre ellos. Este ángulo degenera a cero si los planos son tangentes en la esfera $\widehat{\mathbb{C}}$ en el infinito. En otros casos, los planos tienen un único perpendicular común, cuya longitud define la *distancia* entre ellos.

Las líneas geodésicas en $\mathbb{H}^3$ son semicírculos y líneas rectas ortogonales a $\mathbb{C}$. La geometría de puntos, líneas, planos y sus relaciones de incidencia en el espacio hiperbólico es bien conocida. En particular, tenemos formulas trigonométricas para los triángulos hiperbólicos que son muy útiles. También tenemos fórmulas para calcular longitudes, áreas y volúmenes que luego usaremos. Los volúmenes son calculados con respecto al elemento de volumen hiperbólico $dV$ inducido por la métrica. Para nuestro modelo de $\mathbb{H}^3$, tenemos $dV = \frac{1}{t^3}dxdydt$.

**Ejercicios 1.1.** *1. Pruebe que* $\mathrm{PSL}_2(\mathbb{C})$ *actúa transitivamente sobre geodésicas en* $\mathbb{H}^3$.
*2. Obtener la fórmula de la distancia entre un punto y un plano hiperbólico en* $\mathbb{H}^3$.
*3. Obtener las fórmulas trigonométricas para un triángulo hiperbólico.*
*4. Calcular el área de un triángulo en* $\mathbb{H}^2$ *y el volumen de un símplex ideal (i.e. con todos los vértices en infinito) en* $\mathbb{H}^3$.

Ahora vamos a considerar los *subgrupos* de $\mathrm{PSL}_2(\mathbb{C})$ y la geometría relacionada a ellos.

Para comenzar, para los elementos $\gamma \neq \mathrm{Id}$ tenemos la siguiente clasificación:

- $\gamma$ es *elíptico* si $\mathrm{tr}(\gamma) \in \mathbb{R}$ y $|\mathrm{tr}(\gamma)| < 2$;
- $\gamma$ es *parabólico* si $\mathrm{tr}(\gamma) = \pm 2$;

- $\gamma$ es *loxodrómico* (o *hiperbólico*) en otro caso.

Tenemos las siguientes propiedades geométricas correspondientes a esta clasificación:

$\gamma$ es parabólico $\iff \gamma$ tiene un único punto fijo en $\widehat{\mathbb{C}} = \partial\mathbb{H}^3$.

**Ejemplo.** $\gamma : z \to z+1$ tiene matriz $\begin{pmatrix} 1 & 1 \\ 0 & 1 \end{pmatrix}$ y fija únicamente el punto $\infty \in \widehat{\mathbb{C}}$.

Podemos verificar que todos los otros elementos parabólicos son conjugados en $\mathrm{PSL}_2(\mathbb{C})$ a este $\gamma$.

$\gamma$ es elíptico $\iff \gamma$ tiene dos puntos fijos $p, q \in \widehat{\mathbb{C}}$ y es una rotación alrededor de la geodésica $A_\gamma = [p,q] \subset \mathbb{H}^3$.

$\gamma$ es hiperbólico $\iff \gamma$ tiene dos puntos fijos $p, q \in \widehat{\mathbb{C}}$ y es un movimiento de tornillo a lo largo de la geodésica $A_\gamma = [p,q] \subset \mathbb{H}^3$.

La geodésica $A_\gamma$ es llamada *eje* de un elemento elíptico o hiperbólico $\gamma$. Para elementos hiperbólicos tenemos la fórmula importante para el *desplazamiento* $\ell$ a lo largo del eje:

$$\cosh(\ell(\gamma)/2) = |\operatorname{tr}(\gamma)|/2. \tag{1}$$

Observe que sólo los elementos elípticos poseen puntos fijos en $\mathbb{H}^3$.

El grupo $\mathrm{PSL}_2(\mathbb{C})$ actúa transitivamente en los puntos de $\mathbb{H}^3$, de modo que el estabilizador de cualquier punto en $\mathbb{H}^3$ es conjugado al estabilizador de $(0,0,1)$, que es $\mathrm{SO}_3(\mathbb{R})$ — un subgrupo compacto máximo de $\mathrm{PSL}_2(\mathbb{C})$. Tenemos que

$$\mathbb{H}^3 \cong \mathrm{PSL}_2(\mathbb{C})/\,\mathrm{SO}_3(\mathbb{R}),$$

como un *espacio simétrico*.

Igualmente, tenemos una acción transitiva de $\mathrm{PSL}_2(\mathbb{C})$ en $\widehat{\mathbb{C}}$, entonces los estabilizadores de los puntos en infinito son conjugados a

$$B = \mathrm{St}(\infty) = \left\{ \begin{pmatrix} a & b \\ 0 & a^{-1} \end{pmatrix} \middle|\, a \in \mathbb{C}^*,\ b \in \mathbb{C} \right\},$$

el *subgrupo de Borel* de $\mathrm{PSL}_2(\mathbb{C})$.

Cualquier subgrupo finito de $\mathrm{PSL}_2(\mathbb{C})$ debe tener un punto fijo en $\mathbb{H}^3$ y así ser conjugado a un subgrupo de $\mathrm{SO}_3(\mathbb{R})$. Los subgrupos finitos de $\mathrm{SO}_3(\mathbb{R})$ son bien conocidos:

> cíclicos ($C_n$), diedrales ($D_n$) y los grupos de simetrías de poliedros regulares $\mathrm{A}_4$, $\mathrm{S}_4$ o $\mathrm{A}_5$.

**Definición 1.2.** *Sea* $\Gamma$ *un subgrupo de* $\mathrm{PSL}_2(\mathbb{C})$*:*

- $\Gamma$ *es* reducible *si todos los elementos* $\gamma \in \Gamma$ *tienen un punto fijo común en* $\widehat{\mathbb{C}}$. *Si* $\Gamma$ *no es reducible, este es llamado* irreducible.
- $\Gamma$ *es* elemental *si posee una órbita finita en* $\mathbb{H}^3 \cup \widehat{\mathbb{C}}$.

Observe que reducible implica elemental, pero el recíproco no es cierto. Por ejemplo, un subgrupo elemental puede tener elementos parabólicos y elípticos.

**Teorema 1.3.** *Cada subgrupo no elemental de* $\mathrm{PSL}_2(\mathbb{C})$ *contiene una infinidad de elementos loxodrómicos tal que ningún par de estos elementos posee un punto fijo en común.*

**Definición 1.4.** *Un* grupo kleiniano *es un subgrupo discreto de* $\mathrm{PSL}_2(\mathbb{C})$.

Esta condición es equivalente a exigir que $\Gamma$ actúa discontinuamente en $\mathbb{H}^3$, que significa que para cada compacto $K \subset \mathbb{H}^3$ el conjunto $\{\gamma \in \Gamma | \gamma K \cap K\}$ es finito.

El estabilizador en $\Gamma$ de un punto $\xi$ en infinito puede conjugarse a un subgrupo de Borel $B$. Los subgrupos discretos de $B$ son clasificados (*Ejercicio:* Descríbelos), el caso que es particularmente importante para la geometría da la siguiente definición:

**Definición 1.5.** *El punto* $\xi \in \widehat{\mathbb{C}}$ *es una* cúspide *de un grupo kleiniano* $\Gamma$ *si el estabilizador* $\Gamma_\xi$ *contiene un grupo abeliano libre de rango* 2.

(En este caso $\Gamma_\xi \cong (\mathbb{Z} \oplus \mathbb{Z}) \rtimes F$, $F$ — un grupo finito.)

Puesto que un grupo kleiniano actúa discontinuamente en $\mathbb{H}^3$, podemos construir un *dominio fundamental* para esta acción. Por definición, un dominio fundamental es un subconjunto cerrado $\mathfrak{F} \subset \mathbb{H}^3$ tal que:

- $\bigcup_{\gamma \in \Gamma} \gamma \mathfrak{F} = \mathbb{H}^3$;
- $\mathfrak{F}^\circ \bigcap \gamma \mathfrak{F}^\circ = \varnothing$, para todo $\gamma \neq \mathrm{Id}$, donde $\mathfrak{F}^\circ$ denota el interior de $\mathfrak{F}$;
- la frontera de $\mathfrak{F}$ tiene medida cero.

Podemos obtener un dominio fundamental para $\Gamma$ usando la *construcción de Dirichlet:*

Para $P \in \mathbb{H}^3$ un punto tal que $\gamma(P) \neq P$ para todo $\gamma \in \Gamma \setminus \{\mathrm{Id}\}$, define

$$\mathfrak{F}_P(\Gamma) := \{Q \in \mathbb{H}^3 \mid d(Q,P) \leqslant d(\gamma(Q),P)\ \forall \gamma \in \Gamma\}.$$

**Ejercicio 1.6.** *Muestre que* $\mathfrak{F}_P(\Gamma)$ *es un dominio fundamental del grupo kleiniano* $\Gamma$.

**Definición 1.7.** *Un grupo kleiniano es llamado geométricamente finito si él admite un dominio de Dirichlet con un número finito de lados.*

Los grupos geométricamente finitos son finitamente generados.

Un grupo kleiniano tiene *covolumen finito* si posee dominio fundamental $\mathfrak{F}$ de volumen hiperbólico finito. En este caso, el covolumen de $\Gamma$ es

$$\mathrm{covol}(\Gamma) = \mathrm{vol}(\mathfrak{F}) = \int_{\mathfrak{F}} dV.$$

Un grupo es llamado *cocompacto* si tiene un dominio fundamental compacto.

**Proposición 1.8.** *Sean* $\mathfrak{F}_1$, $\mathfrak{F}_2$ *dos dominios fundamentales de un grupo kleiniano* $\Gamma$. *Entonces, si* $\int_{\mathfrak{F}_1} dV$ *es finito,* $\int_{\mathfrak{F}_2} dV$ *es también finito y son iguales.*

Si $\Gamma$ es cocompacto, su covolumen es finito. En la otra dirección tenemos:

**Teorema 1.9.** *Si* $\Gamma$ *tiene covolumen finito, entonces existe* $P \in \mathbb{H}^3$ *tal que* $\mathfrak{F}_P(\Gamma)$ *tiene un número finito de lados. En particular,* $\Gamma$ *es geométricamente finito y finitamente generado.*

Sea $\Gamma$ conteniendo un elemento parabólico $\gamma$. Podemos suponer que el punto fijo de $\gamma$ es $\infty \in \widehat{\mathbb{C}}$. En este caso existe una *horobola*

$$H_\infty(t_0) = \{(x, y, t) \in \mathbb{H}^3 \mid t > t_0\}$$

tal que la acción de $\Gamma$ en $H_\infty(t_0)$ es la misma que la acción de $\Gamma_\infty$:

$$\text{para } x, y \in H_\infty(t_0) \text{ existe } \delta \in \Gamma \mid \delta(x) = y \iff \delta \in \Gamma_\infty.$$

Entonces $\Gamma_\infty$ actúa sobre $H_\infty(t_0)$, y como $\Gamma_\infty \subset B$ él actúa en la *horoesfera*, el borde de la horobola (i.e. $\{(x, y, t_0)\}$), como un grupo de transformaciones euclidianas. Así tenemos una descripción precisa de la acción de un grupo kleiniano en las cercanías de una cúspide.

Si $\Gamma$ contiene un elemento parabólico podemos ver que $\Gamma$ no es cocompacto. Además con la condición de covolumen finito podemos obtener un mejor resultado:

**Teorema 1.10.** *Sea $\Gamma$ un grupo kleiniano de covolumen finito. Si $\Gamma$ no es cocompacto, entonces $\Gamma$ tiene que contener un elemento parabólico $\gamma$. Si $\xi$ es el punto fijo de $\gamma$, entonces $\xi$ es una cúspide. Además, existe sólo un número finito de las clases de $\Gamma$-equivalencia de cúspides y para esto sus vecindades horobolas pueden ser disjuntas.*

Si empezamos con $\mathrm{PSL}_2(\mathbb{R})$ en lugar de $\mathrm{PSL}_2(\mathbb{C})$ buena parte de la discusión anterior permanece válida. Nos limitaremos a indicar la definición.

**Definición 1.11.** *Un subgrupo discreto de* $\mathrm{PSL}_2(\mathbb{R})$ *es llamado un* grupo fuchsiano.

Para dar un ejemplo consideremos el grupo $\Gamma = \mathrm{PSL}_2(\mathbb{Z}[\sqrt{-3}]) < \mathrm{PSL}_2(\mathbb{C})$. Este es un ejemplo de un grupo kleiniano no cocompacto de covolumen finito. El grupo $\Gamma$ tiene un dominio fundamental

$$\mathfrak{F} = \left\{(x, y, t) \in \mathbb{H}^3 \mid x^2 + y^2 + t^2 \geqslant 1,\ -\frac{1}{2} \leqslant x \leqslant \frac{1}{2},\ 0 \leqslant y \leqslant \frac{\sqrt{3}}{2}\right\}.$$

El covolumen de $\Gamma$ puede ser calculado, y es igual a $\frac{|\Delta_k|^{3/2}\zeta_k(2)}{4\pi^2}$ donde $k = \mathbb{Q}(\sqrt{-3})$, el discriminante $\Delta_k = -3$, y $\zeta_k(2) = 1.28519096\ldots$ que es la función zeta de Dedekind del cuerpo $k$.

El espacio cociente $\Gamma\backslash\mathbb{H}^3$ tiene el espacio de cubrimiento de grado 12 que es homeomorfo al complemento del nudo "figura ocho" en la esfera. Este fue uno de los primeros ejemplos de motivación en la geometría hiperbólica de dimensión 3 estudiados por Riley (1975) y Thurston (1979). Para otros ejemplos de dominios fundamentales de grupos del tipo $\mathrm{PSL}_2(\mathbb{Z}[\sqrt{-d}])$ podemos referir a Elstrodt, Grunewald y Mennicke (1998).

Una *variedad hiperbólica (sin borde) de dimensión n* es una variedad riemanniana $M$ tal que cada punto en $M$ tiene una vecindad isométrica a un abierto en $\mathbb{H}^n$. Si $\Gamma$ es un grupo kleiniano sin torsión, entonces $\Gamma$ actúa libremente en $\mathbb{H}^3$ y el cociente $\Gamma\backslash\mathbb{H}^3$ es una 3-variedad hiperbólica orientable. Recíprocamente, cada 3-variedad hiperbólica orientable $M$ puede ser obtenida como

$$M = \Gamma\backslash\mathbb{H}^3,$$

con $\Gamma$ un grupo kleiniano isomorfo al grupo fundamental $\pi_1(M)$.

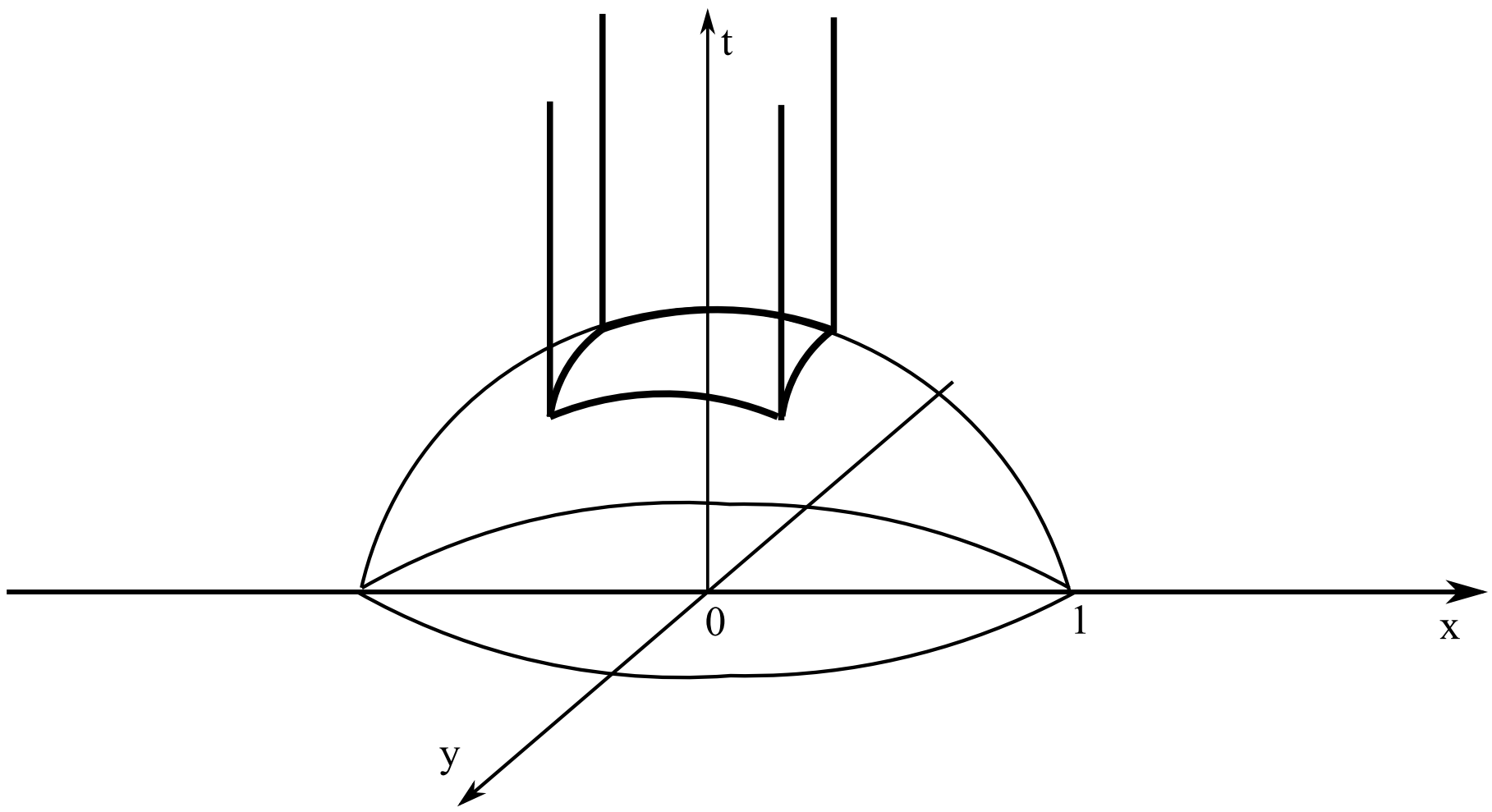

Figura 1: Un dominio fundamental de $\Gamma = \mathrm{PSL}_2(\mathbb{Z}[\sqrt{-3}])$.

Si $\Gamma$ es un grupo kleiniano que puede tener elementos no triviales de orden finito, entonces el espacio $O = \Gamma\backslash\mathbb{H}^3$ es una variedad singular denominada *orbifold*.

Otra vez, la misma terminología aplica para los grupos fuchsianos y variedades hiperbólicas de dimensión 2 y también para las dimensiones $n > 3$.

Por el lema de Selberg cada orbifold hiperbólico tiene un cubrimiento finito que es una variedad hiperbólica:

**Lema 1.12** (Lema de Selberg)**.** *Si* $\Gamma$ *es un subgrupo finitamente generado de* $\mathrm{GL}_n(\mathbb{C})$*, entonces* $\Gamma$ *tiene un subgrupo sin torsión de índice finito.*

**Ejercicio 1.13.** *Dar una prueba del lema.*

Referimos a los libros Ratcliffe (2006) y Maclachlan y Reid (2003) para más información sobre la geometría hiperbólica y los grupos kleinianos.

## 2 Sístole y volumen

Sea $\mathbb{H}^n$ el espacio hiperbólico de dimensión $n$ (para nosotros principalmente $n = 2$ o $3$), y $\mathbb{B}^n_r$ una bola de radio $r$ en $\mathbb{H}^n$. Como antes, $\mathrm{vol}(\cdot)$ denota el volumen hiperbólico. La *sístole*, denotada $\mathrm{sys}_1(\cdot)$, es la longitud de la geodésica cerrada no contráctil más corta (ver Figura 2).

Debido a la curvatura negativa

$$\mathrm{vol}(\mathbb{B}^n_r) \geqslant C_n e^{k_n r}, \text{ para constantes } C_n,\ k_n > 0 \text{ y } r \gg 0.$$

Consideramos una variedad hiperbólica compacta $M = \Gamma\backslash\mathbb{H}^n$. Si $\ell = \mathrm{sys}_1(M)$, entonces $M$ tiene que contener una bola $\mathbb{B}^n_r$ de radio $r = \frac{\ell}{2}$. (En general, el máximo $r$ tal que para cualquier $x \in M$ la variedad riemanniana $M$ contiene una bola $\mathbb{B}_r(x)$ es llamado el *radio de inyectividad* de $M$.)

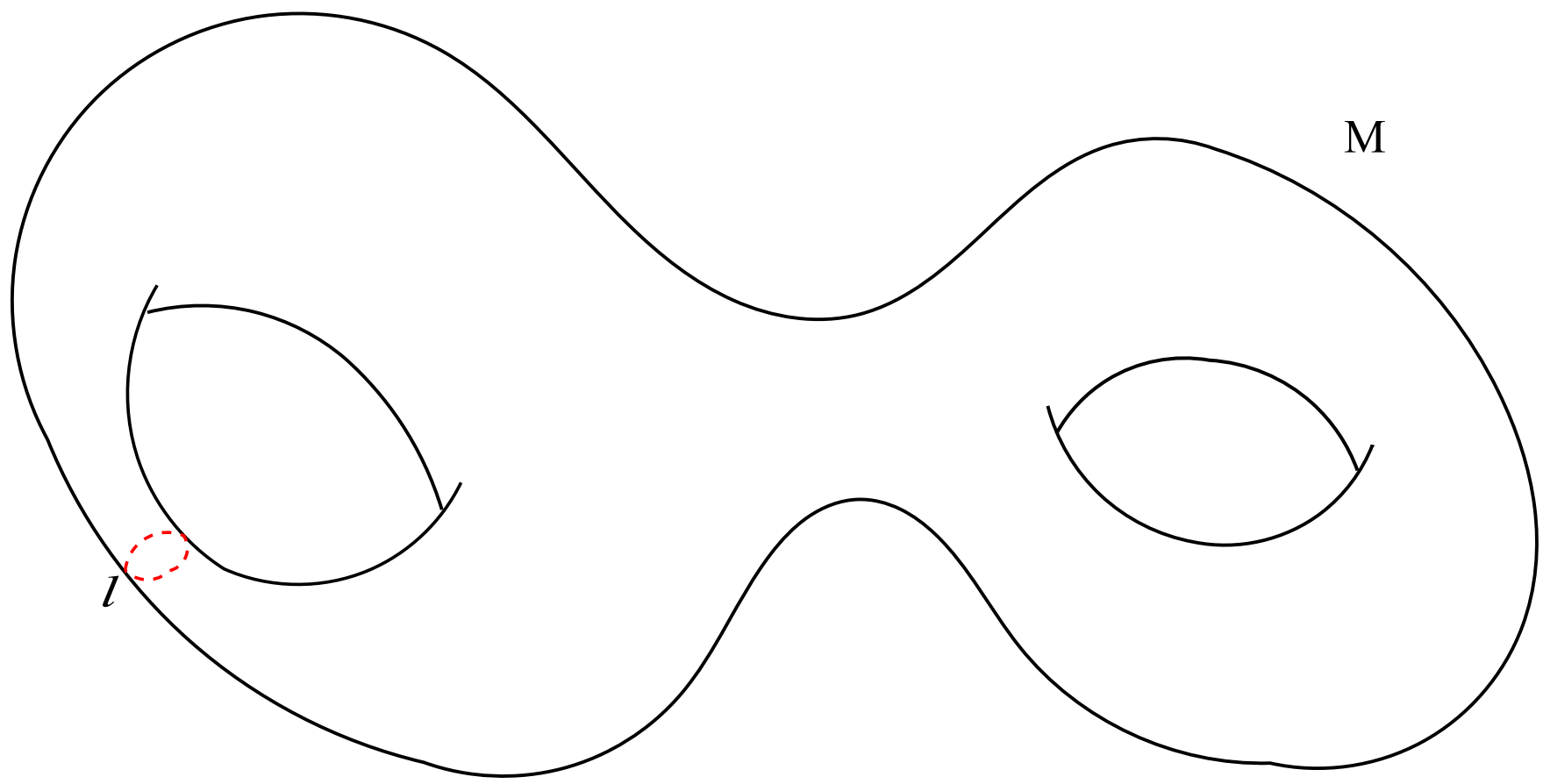

Figura 2: La sístole de la superficie $M$.

Por tanto

$$\operatorname{vol}(M) \geqslant \operatorname{vol}(\mathbb{B}_r^n) \geqslant C_n e^{k_n \ell/2};$$
$$\operatorname{sys}_1(M) \leqslant C \log\big(\operatorname{vol}(M)\big).$$

Así probamos

**Lema 2.1.** *La sístole de una variedad hiperbólica compacta $M$ satisface*

$$\operatorname{sys}_1(M) \leqslant C \log\big(\operatorname{vol}(M)\big),$$

*donde $C$ es una constante positiva que sólo depende de la dimensión de $M$.*

**Ejercicio 2.2.** *Calcula el valor de constante $C$ en el lema.*

La pregunta importante es si esta desigualdad es asintóticamente óptima, i.e. si existen una constante $C' > 0$ y una secuencia de variedades $M_i$ tal que

$$\operatorname{vol}(M_i) \to \infty \quad \text{y} \quad \operatorname{sys}_1(M_i) \geqslant C' \log\big(\operatorname{vol}(M_i)\big).$$

La respuesta es "sí", pero sólo conocemos una manera de obtener tales secuencias. Para esto tenemos que utilizar los *expansores geométricos*.

Un primer ejemplo de una secuencia de expansores geométricos es

$$M_p = \Gamma(p)\backslash\mathbb{H}^2 \text{ con } \Gamma(p) = \ker\big(\mathrm{PSL}_2(\mathbb{Z}) \to \mathrm{PSL}_2(\mathbb{Z}/p\mathbb{Z})\big),\ p \text{ primo.}$$

Tenemos que $M_p \to M$ es un cubrimiento de grado

$$|\,\mathrm{PSL}_2(\mathbb{Z}/p\mathbb{Z})| = \frac{1}{2}p(p-1)(p+1).$$

**Ejercicio 2.3.** *Calcula el género y el número de cúspides de $M_p$.*

Vimos que, en particular, $M$ y todos $M_p$ son no compactos. Para tener la imagen geométrica más clara es mejor tener una secuencia de variedades compactas. Tales secuencias fueron construidas por Buser y Sarnak (1994).

Sean $a, b \in \mathbb{Z}_{>0}$ y $\mathbb{D} = \left(\frac{a,b}{\mathbb{Q}}\right)$ el *álgebra de cuaterniones* sobre $\mathbb{Q}$ generada por

$$1, i, j, k \text{ tales que } i^2 = a,\ j^2 = b,\ ij = -ji = k.$$

Elegimos $a$, $b$ tal que la forma cuadrática

$$N(X) = \mathrm{Norm}(x_0 + x_1 i + x_2 j + x_3 k) = x_0^2 - a x_1^2 - b x_2^2 + ab x_3^2$$

no representa cero para $(x_0, x_1, x_2, x_3) \in \mathbb{Q}^4$, $X \neq 0$. Por ejemplo, podemos tomar $a = 2$, $b = 3$. En este caso $\mathbb{D}$ es un álgebra de división (i.e. cada $X \in \mathbb{D} \setminus \{0\}$ tiene un inverso $X^{-1} \in \mathbb{D}$).

Consideramos subgrupos

$$\widetilde{\Gamma} = \{X \in \mathbb{D}(\mathbb{Z}) \mid N(X) = 1\} \quad \text{y} \quad \widetilde{\Gamma}(p) = \{X \in \widetilde{\Gamma} \mid X \equiv 1(p)\}.$$

Tenemos un isomorfismo de $\widetilde{\Gamma}$ con un subgrupo de $\mathrm{PSL}_2(\mathbb{R})$ dado por

$$X = x_0 + x_1 i + x_2 j + x_3 k \to \begin{pmatrix} x_0 + x_1\sqrt{a} & x_2 + x_3\sqrt{a} \\ b(x_2 - x_3\sqrt{a}) & x_0 - x_1\sqrt{a} \end{pmatrix}.$$

Vamos a denotar las imagenes de $\widetilde{\Gamma}$, $\widetilde{\Gamma}(p)$ por $\Gamma$ y $\Gamma(p)$, respectivamente. Entonces $\{\Gamma(p) < \Gamma\}_{p \text{ primo}}$ es una secuencia de *subgrupos de congruencia* de $\Gamma$, como antes, pero ahora $\Gamma(p)\backslash\mathbb{H}^2$ es compacto para todo $p$. La última afirmación se deduce del hecho de que $\mathbb{D}$ es un álgebra de división (ver Gel'fand, Graev y Pyatetskii-Shapiro (1990, el apéndice para el capítulo 1) para una exposición buena, corta y clara de este material). Además, para $p > 2$ tenemos que $\Gamma(p)$ no tiene torsión (i.e. elementos elípticos de orden finito), entonces $S_p = \Gamma(p)\backslash\mathbb{H}^2$ son superficies riemannianas compactas.

Por la fórmula de Riemann–Hurwitz, el género es igual a

$$\begin{aligned} g(S_p) &= \frac{1}{2}p(p-1)(p+1)\frac{|\chi(\Gamma\backslash\mathbb{H}^2)|}{2} + 1 \\ &= p(p-1)(p+1)\nu + 1, \end{aligned}$$

donde $\nu = \nu(a, b)$ es una constante positiva definida por la característica de Euler de $\Gamma\backslash\mathbb{H}^2$.

Ahora vamos a estimar la sístole de $S_p$. Sea $p > 0$ un número primo y

$$\alpha = x_0 + x_1 i + x_2 j + x_3 k \in \widetilde{\Gamma}(p).$$

Tenemos que $p | x_j$ para $j = 1, 2, 3$, y

$$1 = N(\alpha) = x_0^2 - a x_1^2 - b x_2^2 + ab x_3^2.$$

Por lo tanto

$$\begin{aligned} x_0^2 &\equiv 1(p^2); \\ x_0 &\equiv \pm 1(p^2). \end{aligned}$$

Si $\alpha \neq \pm 1$, entonces $x_0 \neq \pm 1$ (pues $\Gamma$ no tiene elementos parabólicos), y por eso

$$|x_0| \geqslant p^2 - 1.$$

Entonces, si $\gamma \in \Gamma(p)$ y $\gamma \neq 1$, tenemos que

$$|\operatorname{tr}(\gamma)| \geqslant 2p^2 - 2.$$

En particular, esto confirma que $\gamma \in \Gamma(p)$ son hiperbólicos y entonces $S_p = \Gamma(p)\backslash\mathbb{H}^2$ es de hecho una superficie suave.

La fórmula (1) implica que

$$\ell(\gamma) > 2\log\big(|\operatorname{tr}(\gamma)| - 1\big),$$

y entonces para $\gamma \in \Gamma(p)$,

$$\ell(\gamma) > 2\log(2p^2 - 3).$$

Esto implica:

**Teorema 2.4** (Buser y Sarnak (1994))**.** *Para la secuencia de variedades $S_p$ tenemos*

$$\operatorname{sys}_1(S_p) \geqslant \frac{4}{3}\log\big(g(S_p)\big) - C,$$

*donde $C > 0$ no depende de $p$.*

Este resultado fue posteriormente generalizado para otras variedades aritméticas de dimensión 2 y 3 por Katz, Schaps y Vishne (2007).

**Ejercicios 2.5.** *1. Calcula la constante $C$.*

*2. Generaliza el resultado para variedades no compactas como $M_p$.*

Otra mirada a las sístoles de los espacios de recubrimiento proviene de la teoría geométrica de grupos. Según Gromov, podemos definir la noción de una *sístole de un grupo* finitamente generado por un conjunto $X$, de modo que si $\Gamma = \pi_1(M)$ tenemos que

$$\operatorname{sys}_1(\Gamma, X) \simeq \operatorname{sys}_1(M).$$

Más precisamente, sea $M$ una variedad riemanniana compacta aesférica (i.e. el recubrimiento universal de $M$ es contráctil) y $\Gamma = \pi_1(M)$. Consideramos una secuencia de recubrimientos:

$$M_i \to M \text{ de grado } d_i.$$

Podemos considerar la secuencia correspondiente de subgrupos $\Gamma_i = \pi_1(M_i) < \Gamma$ de índice $d_i$. Ahora fijemos un conjunto de generadores $X \subset \Gamma$ y definimos $\operatorname{sys}_1(\Gamma_i, X)$ como la mínima longitud de $X$-palabras de un elemento no trivial de $\Gamma_i$.

*Propiedades:*

- Si $\Gamma_i \lhd \Gamma$ es un subgrupo normal, $X$ genera un subgrupo libre de $\Gamma$ y es simétrico (i.e. $X = X^{-1}$), entonces tenemos que $\operatorname{sys}_1(\Gamma_i, X)$ es igual a la sístole del grafo de Cayley del grupo $\Gamma/\Gamma_i$ con respecto a la proyección de los generadores de $X$.

- Para cualquier otro conjunto de generadores $X'$ de $\Gamma$ tenemos

$$C_1 \leqslant \frac{\mathrm{sys}_1(\Gamma_i, X)}{\mathrm{sys}_1(\Gamma_i, X')} \leqslant C_2,$$

donde $C_1, C_2 > 0$ no dependen de $i$.
Por esta propiedad podemos suprimir $X$ de la notación de sístole cuando se estudia su comportamiento asintótico.

- Existen también dos constantes positivas $D_1$, $D_2$ que no dependen de $i$ tales que

$$D_1 \leqslant \frac{\mathrm{sys}_1(M_i)}{\mathrm{sys}_1(\Gamma_i, X)} \leqslant D_2.$$

*Ejemplos:*

(1) Si $\Gamma$ es un grupo abeliano libre de rango $n$, entonces $\mathrm{sys}_1(\Gamma_i) \lesssim d_i^{1/n}$ y esta desigualdad es óptima, i.e. existe una secuencia $\{\Gamma_i\}$ con $\mathrm{sys}_1(\Gamma_i) \simeq d_i^{1/n}$.

(2) Si $\Gamma$ es un grupo nilpotente sin torsión de crecimiento polinomial de grado $m$, entonces $\mathrm{sys}_1(\Gamma_i) \lesssim d_i^{1/m}$ y esta desigualdad es óptima.

(3) Si $\Gamma$ tiene crecimiento exponencial, entonces $\mathrm{sys}_1(\Gamma_i) \lesssim \log(d_i)$. Un ejemplo de $\Gamma$ de este tipo es $\Gamma = \pi_1(M)$ para una variedad riemanniana $M$ de curvatura estrictamente negativa.

La desigualdad en el ejemplo (3) es asintóticamente óptima gracias al lema siguiente:

**Lema 2.6.** *Sea $\Gamma < \mathrm{SL}_n(\mathbb{Z})$ con conjunto de generadores $X$ y sean $\{\Gamma(p) < \Gamma\}$ sus subgrupos de congruencia módulo primos $p$. Entonces*

$$\mathrm{sys}_1(\Gamma(p), X) \geqslant C \log p, \text{ con } C = C(X) \text{ una constante positiva.}$$

(En Gromov (1996), Gromov supone que $\Gamma$ no tiene los elementos unipotentes pero esto no es necesario.)

**Ejercicio 2.7.** *Dar una prueba del lema.*

Los resultados de Buser–Sarnak y Katz–Schaps–Vishne son casos particulares de esta lema, pero notar que en estos casos ellos saben también el *valor de la constante $C$*. El lema también se aplica para subgrupos de congruencia de grupos *no aritméticos* finitamente generados. Infelizmente, la prueba del lema no da información para calcular $C$, sólo estimaciones muy aproximadas.

El famoso teorema (lema) de Milnor y Schwartz muestra que si un grupo $\Gamma$ actúa geométricamente en un espacio riemanniano $X$, entonces $\Gamma$ es finitamente generado y quasi-isométrico con $X$, en particular, $\Gamma$ y $X$ tienen el mismo tipo de crecimiento. Este resultado explica la relación entre sístoles de espacios y de grupos.

Después de la clase, Andrzej Żuk mencionó el siguiente problema interesante:

**Problema abierto 2.8.** *¿Existen un espacio métrico $M$ y una secuencia infinita de cubrimientos $\{M_i \to M\}$ tal que $\mathrm{sys}_1(M_i) \geqslant C \log\big(\mathrm{vol}(M_i)\big) - C'$ con $C > \frac{4}{3}$?*

## 3 Propiedades geométricas de las superficies de congruencia

En esta sección vamos a investigar principalmente superficies riemannianas (sin frontera) de curvatura $-1$. Sea $S$ una superficie hiperbólica dada por $S = \Gamma\backslash\mathbb{H}^2$, con $\Gamma < \mathrm{PSL}_2(\mathbb{R})$ un subgrupo fuchsiano.

Sea $\Delta$ el operador de Laplace en $S$ (entonces $\Delta(f) = -\mathrm{div}\,\mathrm{grad}(f)$), y sean $0 = \lambda_0 < \lambda_1 \leqslant \ldots$ los autovalores de $\Delta$. El primer autovalor positivo $\lambda_1$ es llamado el *gap espectral* de $S$. Podemos calcular este autovalor por la fórmula de Rayleigh:

$$\lambda_1(S) = \inf_f \frac{\int_S \|\mathrm{grad}\, f\|^2 dV}{\int_S f^2 dV}, \tag{2}$$

donde el ínfimo es tomado sobre las funciones $f$ con soporte compacto en $S$ y tal que $\int_S f dV = 0$. Aquí $dV$, grad, etc. son definidos usando la métrica hiperbólica que desciende de $\mathbb{H}^2$.

Es fácil ver que $\lambda_1$ se puede hacer arbitrariamente pequeño, incluso para superficies $S$ de un género fijo.

**Ejemplo 3.1.** Consideremos una secuencia de superficies $\{S_i\}$ de género 2 con $\mathrm{sys}_1(S_i) = \ell(\gamma_i)$ y tal que $\gamma_i$ separa $S_i$ en dos partes. En este caso

$$\text{si } \ell(\gamma_i) \to 0 \text{ cuando } i \to \infty, \text{ tenemos } \lambda_1(S_i) \to 0.$$

Ciertamente, podemos tomar $f = f_i$ en (2) tal que $f_i = const$ fuera de una vecindad del collar de $\gamma_i$. En este caso $\int_{S_i} \|\mathrm{grad} f_i\|^2 dV$ es proporcional a $\ell(\gamma_i)$ y al mismo tiempo podemos tener $\int_{S_i} f_i^2 dV = 1$ (y $\int_{S_i} f_i dV = 0$). Entonces $\lambda_1(S_i) \leqslant c\ell(\gamma_i) \to 0$ por (2).

Ahora sea $M$ una variedad riemanniana compacta de dimensión $n$ y con curvatura $-1 \leqslant \kappa \leqslant -a^2$ ("pinched negative"). Tenemos los siguientes resultados:

**Teorema 3.2** (Cheeger (1970))**.**

$$\lambda_1(M) \geqslant \frac{1}{4} h(M)^2,$$

*donde $h$ es la* constante isoperimétrica de Cheeger *definida por*

$$h(M) = \inf_N \frac{\mathrm{area}(N)}{\min(\mathrm{vol}(M_1), \mathrm{vol}(M_2))},$$

*donde la hipersuperficie $N$ tiene dimensión $n-1$ y divide $M$ en dos partes, $M_1$ y $M_2$.*

**Teorema 3.3** (Buser (1982))**.**

$$\lambda_1(M) \leqslant c_1 h(M) + c_2 h(M)^2,$$

*donde $c_1$ y $c_2$ son constantes positivas explícitas.*

**Observación 3.4.** El teorema de Cheeger es verdad para cualquier variedad riemanniana compacta y tiene una generalización para las variedades no compactas; la constante $\frac{1}{4}$ en el teorema es óptima. En el teorema de Buser para dimensión 2 y curvatura $-1$ podemos tomar $c_1 = 2$, $c_2 = 10$, pero estos valores pueden no ser óptimos.

**Teorema 3.5** (Brooks (1992)).

$$\operatorname{diam}(M) \leqslant C_1(h,r)\log\operatorname{vol}(M) + C_2(h,r),$$

*donde las constantes $C_1(h,r)$, $C_2(h,r) > 0$ dependen sólo de la constante isoperimétrica $h(M)$ y el radio de inyectividad $r(M)$.*

La prueba del teorema de Brooks es corta y bonita:

*Prueba.* Sea $x \in M$, denotemos por $V(t,x)$ el volumen de una bola de radio $t$ y centro $x$ en $M$.

Si la curvatura de $M$ es $\leqslant -a^2$ y $t \leqslant r$, el radio de inyectividad, podemos acotar $V(t,x)$ por debajo por el volumen de una bola de radio $t$ en el espacio hiperbólico de curvatura constante $= -a^2$.

Ahora suponga que $t > r$. En este caso tenemos

$$\frac{V'(t,x)}{V(t,x)} \geqslant h,$$

mientras que $V(t,x) < \frac{1}{2}\operatorname{vol}(M)$. (Esto es porque el cambio infinitesimal del volumen es el área de la frontera.)

Mediante la integración obtenemos

$$V(t,x) \geqslant e^{h(t-r)}V(r,x) \text{ hasta } V(t,x) = \frac{1}{2}\operatorname{vol}(M),$$

que ocurre cuando

$$t = t_0 = \frac{1}{h}\Big(\log\big(\operatorname{vol}(M)\big) - \log(2) - \log\big(V(r)\big)\Big) + r.$$

Entonces para cualquier $x, y \in M$ tenemos $\mathbb{B}_{t_0}(x) \cap \mathbb{B}_{t_0}(y) \neq 0$, y esto implica $\operatorname{diam}(M) \leqslant 2t_0$. □

**Ejercicio 3.6.** *Investigar $h(S)$ y $\operatorname{diam}(S)$ para $S$ en el Ejemplo 3.1.*

Ahora volvamos a las superficies de congruencia $\{M_p \to M\}$ con

$$M_p = \Gamma(p)\backslash\mathbb{H}^2 \text{ con } \Gamma(p) = \ker\big(\mathrm{PSL}_2(\mathbb{Z}) \to \mathrm{PSL}_2(\mathbb{Z}/p\mathbb{Z})\big).$$

Un problema fundamental abierto es

**Conjetura 3.7** (Selberg (1965)). $\lambda_1(M_p) \geqslant \frac{1}{4}$.

Si la conjetura es cierta, la estimación es óptima. Esta conjetura tiene muchas aplicaciones importantes, ver Sarnak (1995) y las referencias allí.

Selberg mostró en el mismo artículo el siguiente resultado:

**Teorema 3.8** (Selberg). *Tenemos $\lambda_1(M_p) \geqslant \frac{3}{16}$.*

Nosotros vamos a probar un resultado un poco más débil siguiendo a Sarnak y Xue (1991) (ver también Sarnak (1995)).

El grupo de transformaciones de recubrimiento $M_p \to M$ es $G = \mathrm{PSL}_2(\mathbb{Z}/p\mathbb{Z})$. Suponga que existe un autovalor $\lambda$ *excepcional*, i.e. $0 \leqslant \lambda \leqslant 1/4$. Primero podemos ver que $\lambda$ tiene que tener una *multiplicidad alta*. Sea $V_\lambda$ el espacio propio que corresponde a $\lambda$. El laplaciano conmuta con las transformaciones de recubrimiento, entonces $\Delta$ actúa en $V_\lambda$, y por esto $V_\lambda$ tiene que contener una representación irreducible no trivial de $G$. Por un teorema de Frobenius, cualquier representación irreducible no trivial de $\mathrm{PSL}_2(\mathbb{Z}/p\mathbb{Z})$ tiene dimensión $\geqslant (p-1)/2$. Llegamos a la conclusión que $\dim(V_\lambda) \geqslant (p-1)/2$. La idea es mostrar que un pequeño autovalor no puede tener una multiplicidad $m(\lambda, M_p)$ tan grande.

Sarnak y Xue mostraron que para $\epsilon > 0$ existe $C_\epsilon > 0$ tal que

$$m(\lambda, M_p) \leqslant C_\epsilon \, |\mathrm{PSL}_2(\mathbb{Z}/p\mathbb{Z})|^{1-2\nu+\epsilon}, \tag{3}$$

donde $\lambda = 1/4 - \nu^2$ $(0 \leqslant \lambda \leqslant 1/4)$.

**Ejercicio 3.9** ($\star$). *Dar una prueba de* (3).

La combinación de (3) con la desigualdad para $\dim(V_\lambda)$ y la información sobre $|G|$ implica:

$$p - 1 \leqslant K\, p^{3(1-2\nu+\epsilon)}.$$

Para $p$ suficientemente grande, eso es posible sólo si $3(1-2\nu+\epsilon) \geqslant 1$, lo cual implica que

$$\lambda_1(M_p) \geqslant \frac{5}{36} - \epsilon.$$

Para el último resultado de este tipo que da $\lambda_1(M_p) \geqslant \frac{975}{4096}$ ver Kim y Sarnak (2003). Vignéras (1983) probó que $\lambda_1(S_p) \geqslant \frac{3}{16}$ para cualquier superficie aritmética de congruencia, incluidas las superficies compactas de Buser y Sarnak. En Burger y Sarnak (1991) los autores probaron un resultado de este tipo para variedades hiperbólicas de dimensión $n \geqslant 3$, mostrando que $\lambda_1(M_p^n) \geqslant \frac{2n-3}{4}$. El caso general de la conjetura de Selberg es conocido como la *conjetura de Ramanujan generalizada*.

**Ejercicio 3.10.** *Investigar* $h(S_p)$ *y* $\operatorname{diam}(S_p)$ *para las superficies de Buser–Sarnak.*

Una idea de Kazhdan fue que podemos obtener estimaciones no triviales para $\lambda_1(M_p)$ usando solamente las informaciones sobre $\mathrm{sys}_1(M_p)$ y $\dim(V_\lambda)$, pero Brooks (1988) mostró que $\mathrm{sys}_1(M_p) \sim \frac{4}{3}\log(|G|)$ no es suficiente para eso. Con cualquier constante $> \frac{4}{3}$ se podría dar una prueba completamente geométrica de $\lambda_1(M_p) \geqslant \delta > 0$.

En artículo Brooks (1992), superficies de congruencia fueron caracterizadas como:

*(a) cortas y gordas;*

*(b) con simetrías interesantes.*

La primera propiedad está relacionada con el diámetro y el sístole, mientras que la segunda es sobre el grupo $\mathrm{PSL}_2(\mathbb{Z}/p\mathbb{Z})$. Lo que es interesante aquí no es sólo el número de simetrías, sino también las representaciones del grupo.

**Problema 3.11.** *Investigar cómo estas propiedades dependen del espacio simétrico* $\mathfrak{X}$, *que en nuestro caso fue el plano hiperbólico* $\mathbb{H}^2$.

# 4 La descomposición de Heegaard de la 3-variedad

En esta sección vamos a revisar algunos resultados de topología en dimensión 3. El objetivo es entender mejor la topología de 3-variedades de congruencia y por eso es necesario conocer algunos invariantes topológicos. Más información sobre este tema se puede encontrar, por ejemplo, en las notas de Scharlemann (2003) y Scharlemann, Schultens y Saito (2016). Un recurso imprescindible para la topología de 3-variedades es las notas de Thurston (1979).

Vamos asumir en toda la sección que $M$ es una variedad topológica compacta orientable de dimensión 3 posiblemente con borde.

Un teorema fundamental de Moise dice que toda 3-variedad $M$ puede ser triangulada, ver **Moise1977**:

**Teorema 4.1** (Moise, 1952)**.** *(1) Cualquier* $3$*-variedad compacta* $M$ *es homeomorfa a un complejo simplicial finito.*

*(2) Si* $M$ *es homeomorfa a dos complejos* $K$ *y* $L$*, entonces el homeomorfismo* $K \to L$ *es isotópico a un homeomorfismo PL (lineal por partes).*

Supongamos que $(M, \partial M)$ y $(N, \partial N)$ son variedades compactas. Diremos que la inclusión $(N, \partial N) \hookrightarrow (M, \partial M)$ es *propia* si $\partial M \cap N = \partial N$.

La topología PL implica las siguientes propiedades:

1) Un punto en el interior de $M$ posee una vecindad (cerrada) homeomorfa a un 3-bola.

2) Un arco $\alpha$ propiamente encajado en $(M, \partial M)$ tiene una vecindad homeomorfa a $\alpha \times \mathbb{D}^2$, donde $\mathbb{D}^2$ es un disco.

3) Si $M$ es orientable y $c$ es un círculo en el interior de $M$, entonces $c$ tiene una vecindad homeomorfa a $c \times \mathbb{D}^2$.

4) Si $M$ es orientable y $S$ es una superficie orientable propiamente encajada en $M$, entonces $S$ tiene una vecindad homeomorfa a $S \times I$, con $I = [0, 1]$.

**Adjuntando una asa:** Sea $M$ una 3-variedad con borde.

- Adjuntando una 1-asa a $M$ significa la construcción que define una variedad $M \cup_h (I \times \mathbb{D}^2)$, donde $\mathbb{D}^2$ es un disco y $h : (\partial I) \times \mathbb{D}^2 \to \partial M$ es un encaje.
- Adjuntando una 2-asa a $M$ significa $M \cup_h (\mathbb{D}^2 \times I)$, donde $h : (\partial \mathbb{D}^2) \times I \to \partial M$ es un encaje.

Si $h$ preserva la orientación, entonces $M$ con un asa es orientable. Un *cubo con asas* de género $g$ es una bola $\mathbb{B}^3$ con $g$ 1-asas adjuntadas.

**Proposición 4.2.** *Si* $M$ *es una* $3$*-variedad compacta sin borde y* $K$ *es una triangulación de* $M$*, entonces una vecindad cerrada* $\eta(K^1)$ *y su complemento* $M \setminus \mathrm{int}(\eta(K^1))$ *son cubos con asas.*

*Prueba.* Sea $\Gamma = K^1$, un grafo finito conectado, y $\Gamma'$ un árbol máximo de $\Gamma$. Entonces $\eta(\Gamma)$ puede ser obtenida adjuntando $e(K^1) - e(\Gamma')$ 1-asas a $\eta(\Gamma')$, que es una bola (aquí $e(\cdot)$ denota el número de aristas de un grafo).

Para $M \setminus \mathrm{int}(\eta(K^1))$ tenemos que usar el mismo argumento para una triangulación de $M$ dual a $K$. □

**Ejercicio 4.3.** *Calcular el género de $\eta(K^1)$ en función de los números de vértices y aristas de $K^1$.*

**Corolario 4.4.** *Sea $M$ una 3-variedad compacta sin borde. Entonces $M = H_1 \cup_h H_2$, donde $H_1$ y $H_2$ son cubos con asas y $h : \partial H_1 \to \partial H_2$ es un homeomorfismo.*

**Definición 4.5.** *La descomposición $M = H_1 \cup_h H_2$ (o $M = H_1 \cup_S H_2$ con $S = \partial H_1 = \partial H_2$) es llamada la* descomposición de Heegaard *y el mínimo posible género de $S$ es llamado el* género de Heegaard de $M$.

La única 3-variedad compacta sin borde de género de Heegaard igual a 0 es la 3-esfera, y las variedades de género 1 son los espacios lente. Los espacios hiperbólicos tienen el género de Heegaard $\geqslant 2$.

**Definición 4.6.** *Sea $F$ una superficie cerrada no homeomorfa a la esfera $\mathbb{S}^2$. Un* cuerpo de compresión *$H$ es una 3-variedad obtenida de $F \times I$ adjuntando 2-asas a $F \times \{1\}$. Vamos a denotar $F \times \{0\}$ por $\partial_- H$ y $\partial H - F \times \{0\}$ por $\partial_+ H$. Sin asas adjuntadas, $F \times I$ es un cuerpo de compresión* trivial.

**Observaciones 4.7.**

1) Sea $S$ una superficie cerrada y $M$ la variedad obtenida de $S \times I$ adjuntando 2-asas y llenando 2-esferas por 3-bolas. Entonces $M$ es un cubo con asas o un cuerpo de compresión.
2) Una descomposición de Heegaard $H_1 \cup_h H_2$ es llamada *trivial* si $H_1$ o $H_2$ es un cuerpo de compresión trivial, i.e. homeomorfo a $S \times I$.
3) También puede ser definida la descomposición de Heegaard de variedades con borde.

Vamos a considerar las superficies en 3-variedades.

**Definición 4.8.** *Sea $F$ una superficie.*

- *Una curva $\gamma$ simple cerrada en $F$ es llamada* no esencial *si ella divide a $F$ en dos partes, una de las cuales es un disco; caso contrario $\gamma$ es* esencial.
- *Un arco $\alpha$ propiamente encajado en $F$ es llamado* no esencial *si junto con un arco en $\partial F$ delimita un disco, y en el caso contrario $\alpha$ es* esencial.

Sea $\Gamma$ una 1-variedad propiamente encajada en la superficie $F$. Si una componente de $\Gamma$ es una curva simple no esencial, entonces una componente $\alpha$ de $\Gamma$ encierra un disco en $F$ que es disjunto de $\Gamma$. Esta curva $\alpha$ es llamada un *círculo interior* de $\Gamma$.

**Definición 4.9.** *Sea $T$ una superficie encajada en una 3-variedad $M$. La superficie $T$ es llamada* compresible *si $T$ encierra una 3-bola en $M$ o existe una curva esencial simple cerrada $\gamma \subset T$ que encierra un disco $D$ en $M$ tal que $\mathrm{int}(D) \cap T = \varnothing$. En caso contrario $T$ es llamada* no compresible.

El siguiente resultado es el famoso *teorema del lazo:*

**Teorema 4.10** (Papakyriakopoulos, 1956). *Sean $M$ una 3-variedad, $T$ una superficie propiamente encajada en $M$ e $i : T \hookrightarrow M$ la inclusión. Si $i_* : \pi_1(T) \to \pi_1(M)$ no es inyectiva, entonces $T$ es compresible.*

**Corolario 4.11.** *Cualquier superficie propiamente encajada $T$ en $\mathbb{S}^3$ es compresible.*

*Prueba.* Si $T$ no es homeomorfa a la esfera $\mathbb{S}^2$ entonces $i_* : \pi_1(T) \to \pi_1(\mathbb{S}^3)$ no es inyectiva, y por el teorema $T$ es compresible. Si $T$ es homeomorfa a $\mathbb{S}^2$, entonces $T$ encierra una 3-bola en $\mathbb{S}^3$ y $T$ es compresible por la definición. $\square$

Sea $T$ una superficie propiamente encajada en una 3-variedad $M$. Si existe un disco $D$ en $M$ tal que $D \cap T = \partial D$, entonces existe una 3-bola $B = D \times I$ tal que $B \cap T = \partial D \times I$. Sea

$$T' := (T - \partial D \times I) \cup D \times \{0, 1\}.$$

Esta construcción es llamada *2-cirugía en $T$ a lo largo de $D$*. Tenemos

$$\chi(T') = \chi(T) + 2.$$

**Definición 4.12.** *Una 3-variedad $M$ es llamada* reducible *si $M$ contiene una 2-esfera no compresible; caso contrario $M$ es* irreducible.

Sea $M$ una 3-variedad reducible con esfera $P$ que la reduce y con un disco $D$ tal que $D \cap P = \partial D$. Podemos considerar superficie $P'$ obtenida por 2-cirugía en $P$ a lo largo de $D$. En este caso $P'$ tiene dos componentes y una de éstas es también una esfera que reduce $M$.

Suponga que $M_1$ y $M_2$ son dos 3-variedades, y que $S_i$ son 2-esferas que encierran 3-bolas $B_i$ en $M_i$, $i = 1, 2$. La variedad $(M_1 - \mathrm{int}(B_1)) \cup (M_2 - \mathrm{int}(B_2))$, denotada por $M_1 \# M_2$, es llamada la *suma conexa* de $M_1$ y $M_2$. El siguiente teorema muestra que cada 3-variedad es una suma conexa de variedades irreducibles.

**Teorema 4.13** (Kneser, 1929; Milnor, 1957). *Sea $M$ una 3-variedad compacta, orientable. Entonces $M = M_1 \# \ldots \# M_n \# (\mathbb{S}^1 \times \mathbb{S}^2) \ldots \# (\mathbb{S}^1 \times \mathbb{S}^2)$, donde cada $M_i$ es irreducible. Además, la descomposición es única salvo cambios de orden.*

Ahora volvemos a las descomposiciones de Heegaard.

**Definición 4.14.** *Sea $M$ una 3-variedad con una descomposición de Heegaard $M = H_1 \cup_S H_2$.*

- *La descomposición $H_1 \cup_S H_2$ es llamada* estabilizada *si tenemos dos discos propiamente encajados $D_i \subset H_i$, $i = 1, 2$ tales que $\partial D_1$ intersecta $\partial D_2$ transversalmente en un punto.*
- *La descomposición es llamada* reducible *si tenemos dos discos propiamente encajados $D_i \subset H_i$ tales que $\partial D_1 = \partial D_2$; en el caso contrario la descomposición es* irreducible.

**Observación 4.15.** Una definición equivalente de la descomposición de Heegaard reducible dice que existe una 2-esfera $P$ en $M$ tal que $P \cap S$ es esencial en $S$ y tiene sólo una componente.

**Ejercicio 4.16.** *Prueba la equivalencia de las definiciones.*

**Proposición 4.17.** *Sea $H_1 \cup_S H_2$ una descomposición de Heegaard estabilizada. Entonces ella es reducible o ella es una descomposición estándar de género uno de $\mathbb{S}^3$.*

*Prueba.* Sean $D_i \subset H_i$ discos propiamente encajados tales que $|\partial D_1 \cap \partial D_2| = 1$. Sea $B$ la unión de las vecindades bicollares de $D_1$ en $H_1$ y $D_2$ en $H_2$. En este caso $B$ es una 3-bola con el borde $P = \partial B$. Podemos mover un poco $P$ de modo que ella cruce cada $H_i$ en un hemisferio y la curva $\gamma = P \cap S$ separe de $S$ un toro pinchado. Si $\gamma$ es esencial en $S$ entonces $P$ es una esfera que reduce la descomposición. Si $\gamma$ no es esencial, entonces $S$ es un toro que divide $M$ en dos toros sólidos cuyos meridianos se cruzan en un solo punto. Ésta es una descomposición de Heegaard de género 1 de $\mathbb{S}^3$. □

Suponga que $M$ es una 3-variedad con una descomposición de Heegaard $H_1 \cup_S H_2$ y $M'$ es otra 3-variedad con una descomposición $H_1' \cup_{S'} H_2'$. Entonces la variedad $M\#M'$ tiene una descomposición de Heegaard natural definida como sigue: Sea $D$ un disco en $S$ y $D'$ un disco en $S'$. Definimos

$$H_1'' = H_1 \cup_{D=D'} H_1' \text{ y } H_2'' = H_2 \cup_{D=D'} H_2'.$$

Podemos verificar que $H_1''$ y $H_2''$ son los cuerpos de compresión, y que

$$H_1'' \cup_{\partial_+ H_1'' = \partial_+ H_2''} H_2''$$

es una descomposición de Heegaard reducible de la variedad $M\#M'$. La pregunta natural es si es verdad que cualquier descomposición de Heegaard de una variedad reducible es reducible. La respuesta está dada por el conocido teorema de Haken:

**Teorema 4.18** (Haken, 1968)**.** *Sea $M$ una 3-variedad reducible y $M = H_1 \cup_S H_2$ una descomposición de Heegaard. Entonces existe una esfera $P$ que reduce $M$ y tal que $P \cap S$ es un círculo, i.e. la descomposición de Heegaard es reducible.*

*Idea de la prueba.* Podemos suponer que $P$ cruza uno de los cuerpos de compresión $H_1$ o $H_2$ sólo en discos. La idea es minimizar el número de los discos inductivamente para reducir al caso cuando tenemos sólo uno de ellos.

Sea $P$ tal que cruza $H_2$ sólo en discos y consideremos la superficie $P_2 := P \cap H_1$. Apliquemos todas las posibles compresiones para $P_1$ y veamos qué ocurre con $P$: ella se convierte en una unión de 2-esferas; por lo menos una de las cuales es una esfera de compresión. Vamos a denotarla por $P_2$. Ahora podemos repetir la misma operación con $P_2$, etc. Es posible verificar que en cada paso el número de discos disminuye. □

**Observación 4.19.** Otro método para probar el mismo teorema usa grafos ("spine") de cuerpos de compresión.

**Ejercicio 4.20** ($\star$)**.** *Dar los detalles de la prueba.*

## 5 Los cubrimientos de congruencia de una 3-variedad aritmética

En esta sección vamos a investigar el comportamiento del género de Heegaard en las secuencias de expansores geométricos.

Primero tenemos que definir estos tipos de variedades en dimensión 3. Sea $M = \Gamma\backslash\mathbb{H}^3$, con $\Gamma = \mathrm{PSL}_2(\mathcal{O}_d) < \mathrm{PSL}_2(\mathbb{C})$, donde $\mathcal{O}_d$ denota el anillo de los enteros del cuerpo $\mathbb{Q}(\sqrt{-d})$, $d \in \mathbb{Z}_{>0}$. Estos grupos $\Gamma$ son llamados los *grupos de Bianchi*. Por ejemplo, para $d = 1$ tenemos $\mathcal{O}_d = \mathbb{Z}[i]$, los enteros Gaussianos, y $\Gamma = \mathrm{PSL}_2(\mathbb{Z}[i]) < \mathrm{PSL}_2(\mathbb{C})$ es el *grupo de Picard*. Uno de estos grupos fue considerado en el ejemplo en Sección 1.

Ahora, similarmente a la Sección 2, tomemos una secuencia

$$M_{\mathcal{P}} = \Gamma(\mathcal{P})\backslash\mathbb{H}^3 \text{ con}$$
$$\Gamma(\mathcal{P}) = \ker\big(\mathrm{PSL}_2(\mathcal{O}_d) \to \mathrm{PSL}_2(\mathcal{O}_d/\mathcal{P})\big),\ \mathcal{P} \text{ es un ideal primo en } \mathcal{O}_d.$$

Esto da una secuencia de cubrimientos $M_{\mathcal{P}} \to M$ de grados finitos. Para los ideales $\mathcal{P}$ con la norma suficientemente grande las variedades $M_{\mathcal{P}}$ son suaves. Como en el caso de $\mathrm{PSL}_2(\mathbb{Z})$, todas ellas son no compactas pero completas y con volúmenes finitos. Para definir variedades compactas de este tipo necesitamos otra vez usar álgebras de cuaterniones (ver Maclachlan y Reid (2003)).

Lema 2.6 implica que para $||\mathcal{P}|| \gg 1$,

$$\mathrm{sys}_1(M_{\mathcal{P}}) \geqslant C \log(||\mathcal{P}||).$$

En Bachman, Cooper y White (2004) los autores probaron un resultado interesante que da una estimación del género de Heegaard $g$ de una variedad hiperbólica compacta dependiendo de su radio de inyectividad $r$:

$$g \geqslant \frac{\cosh(r)}{2}.$$

Ya que en el en caso compacto $r(M) = \mathrm{sys}_1(M)/2$, la combinación de estos resultados implica que existe una constante $\gamma > 0$ tal que para $\mathcal{P}$ suficiente grande tenemos

$$g(M_{\mathcal{P}}) \geqslant \mathrm{vol}(M_{\mathcal{P}})^{\gamma}.$$

El siguiente resultado muestra que podemos tomar $\gamma = 1$, que es asintóticamente óptima. Este teorema fue probado por Lackenby (2006) e independientemente por Gromov (2009). Las pruebas usan algunos resultados importantes de geometría diferencial.

**Teorema 5.1** (Lackenby (2006) y Gromov (2009))**.** *El género de Heegaard de variedades compactas congruencia $M_{\mathcal{P}} \to M$ satisface*

$$C_1\mathrm{vol}(M_{\mathcal{P}}) \leqslant g(M_{\mathcal{P}}) \leqslant C_2\mathrm{vol}(M_{\mathcal{P}}),$$

*donde $C_1$, $C_2$ son constantes positivas que dependen sólo de $M$.*

*Prueba.* Sea $M = M_1 \cup_S M_2$ una descomposición de Heegaard de la 3-variedad $M$ con género $g = g(S)$ que es igual al género de Heegaard de $M$. Podemos asumir que $g \geqslant 2$.

La superficie de Heegaard $S$ define un *sweepout* $S_t, t \in [0,1]$, de $M$ (i.e. una aplicación continua $f : [0,1] \times S \to M$ tal que $f(t) = S_t$, $S_0$ y $S_1$ son grafos, $S_t$ es isotópico a $S$ por todos $0 < t < 1$ y $\mathrm{Im}(f) = M$). Supongamos que $S_{t_0}$ tiene área máxima entre las superficies $S_t$, $t \in (0,1)$. El ínfimo de estas áreas entre todos los posibles sweepoutes $S_t$ es llamado *área mínmax* de $S$. Por

Pitts y Rubinstein (1986), en este caso existe una *superficie mínima* $F$ en $M$ que tiene área igual al área mínmax y el género $g(F) \leqslant g(S)$.

Sean $h_{ij}$ las componentes de la segunda forma fundamental de $F$ y $R_{ij}$ las componentes de la curvatura seccional de $M$. Por el teorema de Gauss–Bonnet y la fórmula de Gauss tenemos

$$\int_F R_{12} + h_{11}h_{22} - h_{12}^2 ds = 2\pi\chi(F).$$

Como $F$ es una superficie mínima, $h_{11} + h_{22} = 0$, y al ser $M$ hiperbólica, $R_{12} = -1$, entonces:

$$\int_F -1 - h_{11}^2 - h_{12}^2 ds = 2\pi\chi(F);$$

$$\int_F 1 ds = \int_F -h_{11}^2 - h_{12}^2 ds - 2\pi\chi(F) \leqslant -2\pi\chi(F) = 4\pi(g(F) - 1);$$

$$\text{área}(F) \leqslant 4\pi(g(F) - 1) \leqslant 4\pi(g(S) - 1).$$

Entonces para cualquier $\epsilon > 0$ existe un sweepout $S_t^\epsilon$ tal que el área máxima de superficies en $S_t^\epsilon$ es $\leqslant 4\pi(g-1) + \epsilon$.

Sea $G$ la superficie de $S_t^\epsilon$ que divide $M$ en dos partes iguales (tal superficie siempre existe por continuidad). Por la definición de la constante de Cheeger tenemos:

$$h(M) \leqslant \frac{\text{área}(G)}{\text{vol}(M)/2} \leqslant \frac{8\pi(g-1)}{\text{vol}(M)} + \frac{2\epsilon}{\text{vol}(M)};$$

$$g \geqslant \frac{h(M)\text{vol}(M)}{8\pi} - \frac{\epsilon}{4\pi} + 1.$$

Ahora recordamos que para las variedades de congruencia $\lambda_1(M_\wp) \geqslant c_1$ (por Buser y Sarnak (1994)) y entonces por la desigualdad de Buser (ver Teorema 3.3), $h(M_\wp) \geqslant c_2 > 0$. Esto implica $g \geqslant C_1 \text{vol}(M_\wp)$.

La desigualdad opuesta es fácil: Es suficiente levantar una triangulación de $M$ para $M_\wp \to M$ y ver que la descomposición de Heegaard de $M_\wp$ definida por esta triangulación tiene el género $\leqslant c_3 \text{grado}(M_\wp \to M) = C_2 \text{vol}(M_\wp)$ (con $C_2 = c_3 \text{vol}(M)$). □

**Ejercicio 5.2.** *Generalice este resultado para las* 3*-variedades de congruencia no compactas.*

**Ejercicio 5.3.** *Investigue las propiedades de descomposiciones de Heegaard obtenidas en la prueba del teorema.*

Con las ideas similares Gromov mostró:

**Teorema 5.4** (Gromov (2009))**.** *Para cualquier aplicación suave genérica* $F : M_\wp \to \mathbb{R}$ *existe una fibra de* $F$ *con suma de números de Betti* $\geqslant C\text{vol}(M_\wp)$.

Podemos comparar esta *desigualdad de fibra complicada* con la definición de los grafos expansores. La conclusión es que las secuencias de 3-variedades de congruencia pueden ser consideradas como *expansores geométricos* (o topológicos) de dimensión tres.

**Problema abierto 5.5.** *¿Hay un análogo de esta propiedad para las variedades aritméticas en dimensiones superiores?*

**Agradecimientos.** Me gustaría agradecer a Maria Campana Ramia, a Emilio Lauret, a Plinio Pino Murillo y al referí por haber revisado este texto.

## Referencias

MIKHAIL BELOLIPETSKY
mbel@impa.br
INSTITUTO NACIONAL DE MATEMÁTICA PURA E APLICADA (IMPA)
ESTRADA DONA CASTORINA, 110
22460-320 RIO DE JANEIRO RJ
BRAZIL

# 5 Análisis y geometría en grupos

# ELEMENTOS DE ANÁLISIS, GEOMETRÍA Y PROBABILIDAD DE GRUPOS Y GRAFOS

Andrzej Żuk

## 1 Amenabilidad

**Definición** ***von Neumann (1929)***: Un grupo $\Gamma$ es *amenable* si existe $\mu : 2^{\Gamma} \longrightarrow [0;1]$ tal que :

1. $\mu(A \cup B) = \mu(A) + \mu(B)$ para cualquier par $(A, B)$ de subconjuntos disjuntos de $\Gamma$;
2. $\mu(\Gamma) = 1$;
3. $\mu(gA) = \mu(A)$ para todo subconjunto $A$ de $\Gamma$ y todo $g$ en $\Gamma$.

*Notas* :

- Si $\Gamma$ es finito, necesariamente $\mu(A) = \frac{|A|}{|\Gamma|}$.
- Como se demostrará más tarde, un grupo puede ser amenable si y sólo si todos sus subgrupos finitamente generados son amenables. Podemos, por lo tanto, limitarnos a estudiar el caso de un grupo numerable.

Esta noción de "amenable" fue estudiada originalmente para entender la descomposición paradójica de la esfera $S^2$. De hecho, hay 8 subconjuntos $A_1, \cdots, A_4$ y $B_1, \cdots, B_4$ de la esfera $S^2$ y elementos $g_i, h_i$ de $SO(3)$ tales que :

$$S^2 = A_1 \dot{\cup} \cdots \dot{\cup} A_4 \dot{\cup} B_1 \dot{\cup} \cdots \dot{\cup} B_4$$

$$S^2 = g_1(A_1) \dot{\cup} \cdots \dot{\cup} g_4(A_4) \qquad\qquad S^2 = h_1(B_1) \dot{\cup} \cdots \dot{\cup} h_4(B_4)$$

Este teorema se mantiene en gran generalidad y se populariza a menudo como una prueba de que se puede cortar una manzana en piezas finito que se pueden juntar para obtener dos copias idénticas de la manzana inicial. Fue considerado como uno de los logros más espectaculares de la matemática pura en la primera mitad del siglo XX. Por ejemplo, Feynman escribe en su autobiografía que este resultado fue presentado como un argumento para estudiar las matemáticas en lugar de la física.

Banach (1923) probó que $S^1$ no puede tener una descomposición paradójica como $S^2$.

**Definiciones** :

- Una palabra en las letras $a$, $a^{-1}$, $b$, y $b^{-1}$ se dice irreducible si $a\,a^{-1}$, $a^{-1}\,a$, $b\,b^{-1}$ y $b^{-1}\,b$ nunca aparecen.

- El grupo libre en dos elementos $F_2 = < a, b >$ es el conjunto de palabras irreducibles en $a$, $a^{-1}$, $b$, y $b^{-1}$, equipado con la operación que consiste en juntar y reducir: por ejemplo,

$$(ab^3a^{-1})(ab) = ab^3a^{-1}ab = ab^4$$

  y reducción. Su elemento neutro es la palabra vacía, la cual denotamos por $e$.

**Proposicion 1** (von Neumann (1929)). *$F_2$, el grupo libre de rango 2 no es amenable.*

**Demostración**
Para demostrar la proposición, se exhibirá la descomposición paradójica similar a la de la esfera $S^2$. Probaremos que no puede existir tal medida. Sean

$$A_+ = \{\text{palabras que empiezan con } a\} \qquad A_- = \{\text{palabras que empiezan con } a^{-1}\}$$
$$B_+ = \{\text{palabras que empiezan con } b\} \qquad B_- = \{\text{palabras que empiezan con } b^{-1}\}$$

A continuación
$$F_2 = A_+ \dot{\cup} A_- \dot{\cup} B_+ \dot{\cup} B_- \dot{\cup} \{e\}$$

y

$$a^{-1}A_+ \dot{\cup} A_- = F_2,$$
$$b^{-1}B_+ \dot{\cup} B_- = F_2.$$

Aplicando la medida $\mu$ :

$$1 = \mu(F_2) = \mu(A_+ \dot{\cup} A_- \dot{\cup} B_+ \dot{\cup} B_- \dot{\cup} \{e\})$$

Dado que la medida es invariante por translación, la medida de cualquier elemento es nula; en particular, $\mu(\{e\}) = 0$.

Como $\mu$ es una medida invariante por translación, $\mu(A_+) = \mu(a^{-1}A_-)$. Por lo tanto,

$$1 = \mu(a^{-1}A_+) + \mu(A_-) + \mu(b^{-1}B_+) + \mu(B_-) = \mu(F_2) + \mu(F_2) = 2$$

Contradicción. □

*Nota* : Podemos demostrar que un grupo no es amenable si y sólo si admite una descomposición paradójica.

**Proposicion 2** (Banach (1923)). *El grupo $\mathbb{Z}$ es amenable.*

**Demostración**
La idea es definir $\mu(A)$ como el límite $\lim_{n\to+\infty}\frac{|A\cap[-n;n]|}{2n+1}$. El problema es que este límite no tiene sentido en principio. Podemos darle un uso de los ultrafiltros, pero es esencialmente lo mismo que mostrar que $\mathbb{Z}$ es amenable.

Definiremos $\mu$ utilizando una funcional $l^\infty(\mathbb{Z})\to\mathbb{R}$, también denotada $\mu$. Vamos a obtener el valor de $\mu(A)$ definiendo : $\mu(A)=\mu(\mathbb{1}_A)$. Aplicaremos la normalización

$$\mu(c.\mathbb{1}_{\mathbb{Z}})=c.$$

Sea $g\in l^\infty(\mathbb{Z})$, y $n\cdot g$ la función $g$ trasladada en $n$. Más precisamente, $n\cdot g$ se define por:

$$n\cdot g(k)=g(n+k)$$

Entonces podemos definir

$$H=\left\langle\sum_{\text{finito}}\big(g_i-n_i\cdot(g_i)\big),n_i\in\mathbb{Z},g_i\in l^\infty\right\rangle\subseteq l^\infty(\mathbb{Z})$$

y pedimos que

$$\mu(h)=0$$

para todo $h$ en $H$.

Gracias al teorema de Hahn–Banach, entonces podemos extender $\mu$ a todo el espacio. Para aplicar este teorema, tenemos que demostrar que

$$\|\mu_{|\mathbb{C}\mathbb{1}\oplus H}\|\leqslant 1 \tag{1-1}$$

En primer lugar, probamos que

$$\inf_{n\in\mathbb{Z}}h(n)\leqslant 0$$

para cada $h\in\mathcal{H}$. Es fácil demostrar y que es equivalente a

$$\sup_{n\in\mathbb{Z}}h(n)\geqslant 0$$

para cada $h\in\mathcal{H}$.

De esto sigue para $h\in\mathcal{H}$

$$\|c\mathbb{1}+h\|_{l^\infty}\geqslant|c|$$

lo que demuestra (1-1). □

*Nota* : De hecho, podemos demostrar que todo grupo abeliano es amenable.

En general, si un grupo $\Gamma$ contiene un grupo $F_2$, entonces no es amenable. De hecho, uno puede encontrar un descomposición paradójica de $\Gamma$ a partir de la de $F_2$.

Una pregunta es si a la recíproca ($\Gamma$ no amenable $\Longrightarrow F_2\subseteq\Gamma$) es verdadera. Esta pregunta natural fue un problema abierto. Finalmente ha sido demostrado en los años 80 que no es cierta.

Un contraejemplo es el siguiente:

$$B(2,665) = < a, b | w^{665}(a,b) = 1 >$$

donde la notación $w^{665}(a,b) = 1$ significa que cada palabra del grupo, elevado a la potencia 665 es el elemento neutro.

Un teorema de Adyan y Novikov demuestra que este grupo es infinito. Volviendo a la prueba del teorema, Adyan probó que este grupo no es amenable.
*Nota* : Sabemos que los grupos $B(2,n)$ para $n = 2, 3, 4$ y 6 son finitos. Para 5 y 7 a 664, no lo sabemos...

**Teorema 3** (Følner (1955)). *Sea* $\Gamma$ *un grupo numerable. Es amenable si y sólo si existe una sucesión de subconjuntos finitos* $A_n$ *de* $\Gamma$ *tal que:*

$$\lim_{n\to\infty} \frac{|A_n \Delta g A_n|}{|A_n|} = 0 \textit{ para cualquier } g \textit{ en } \Gamma$$

*donde* $gA_n$ *es el trasladado del* $A_n$ *por* $g$.

**Definición** : Un grupo $\Gamma$ se dice de tipo finito (o *finitamente generado*), si tiene un subconjunto $S \subseteq \Gamma$ finito que lo genera.

Se puede demostrar que un grupo puede ser amenable si y sólo si todos sus subgrupos de tipo finito son amenables. Por lo tanto, podemos reducir nuestro estudio a estos grupos.

**Definición** : Un grupo $G$ de tipo finito en $S$ se dice de crecimiento subexponencial si $b(n) = |S^n|$ crece más lentamente que cualquier función exponencial.

**Proposicion 4.** *Si* $G$ *es un grupo generado por un subconjunto* $S$ *finito que es un grupo de crecimiento subexponencial, entonces existe una sucesión estrictamente creciente de enteros* $(n_k)$ *tales que* $S^{n_k}$ *es un conjunto de Følner.*

**Demostración**
Supongamos que no. En este caso, existe una constante $C > 0$ tal que para todo $n$

$$|\partial S^n| \geqslant C|S^n|.$$

Pero, $|S^{n+1}| - |S^n| = |S^{n+1} \setminus S^n| = |\partial S^n| \geqslant C|S^n|$, entonces

$$|S^{n+1}| \geqslant (1+C)|S^n| \geqslant (1+C)^{n+1}$$

La última desigualdad se obtiene por inducción en $n$. □

**Problema abierto :** $S^n$ es la sucesión de Følner, en toda su generalidad (sin necesidad de quitarle una subsucesión) ? Se sabe que esto se da por ejemplo para grupos nilpotentes.

## 2 Grupos de autómatas

Los trabajos de von Neumann (1929) y Day (1957) establecieron que – como lo hemos demostrado –los grupos de crecimiento subexponencial son amenables y que esta clase es cerrada con respecto a las operaciones básicas: extensiones, cocientes, subgrupos y los límites directos. Antes de la construcción del grupo generado por el autómata de la Figura 1, todos los grupos amenables conocidos podrían ser obtenidos a partir de grupos de crecimiento subexponencial utilizando las operaciones básicas descritas anteriormente.

Sea $SG_0$ la clase de grupos tal que todo subgrupo finitamente generado es de crecimiento subexponencial. Supongamos que $\alpha > 0$ es un ordinal y que hemos definido $SG_\beta$ para cada ordinal $\beta < \alpha$. Entonces si $\alpha$ es un ordinal límite,

$$SG_\alpha = \bigcup_{\beta<\alpha} SG_\beta.$$

Si $\alpha$ no es un ordinal límite, sea $SG_\alpha$ la clase de grupos que se puede obtener de los grupos en $SG_{\alpha-1}$ utilizando extensiones y los límites directos. Sea

$$SG = \bigcup_{\alpha} SG_\alpha.$$

Grupos en esta clase se llaman subexponencial amenables.

$SG$ es la clase más pequeña de grupos que contienen los grupos de crecimiento subexponencial, y es cerrada con respecto a las operaciones básicas. Las clases $SG_\alpha$ son cerradas con respecto a cocientes y a subgrupos.

**Definición de grupo generado para un autómata.** Estudiaremos autómatas finitos, reversibles, con el mismo alfabeto en la entrada y salida, por ejemplo $D = \{0, 1, \ldots, d-1\}$ para algún entero $d > 1$. En un tal autómata $A$ se asocia un conjunto finito de estados $Q$, una función de transición $\varphi : Q \times D \to Q$ y una función de salida $\psi : Q \times D \to D$; el autómata $A$ se caracteriza por la terna $(D, Q, \varphi, \psi)$.

El autómata $A$ se llama inversible si para cada $q \in Q$, la función $\psi(q, \cdot) : D \to D$ es una biyección. En este caso, $\psi(q, \cdot)$ se pueden identificar con el elemento correspondiente $\sigma_q$ del grupo simétrico $S_d$ con $d = |D|$ símbolos.

Hay una manera conveniente de representar un autómata finito por un grafo orientado y etiquetado $\Gamma(A)$ cuyos vértices corresponden a los elementos de $Q$. Dos estados $q, s \in Q$ son unidos por una flecha etiquetada con $i \in D$ si $\varphi(q, i) = s$; cada vértice $q \in Q$ está etiquetado por el correspondiente elemento $\sigma_q$ del grupo simétrico.

Los autómatas que acabamos de definir son los autómatas no iniciales. Para hacer de este autómata un autómata inicial debemos elegir un estado $q \in Q$, el estado inicial. El autómata inicial $A_q = (D, Q, \varphi, \psi, q)$ actúa en una sucesión finita o infinita en $D$ de la siguiente manera. Para cada símbolo $x \in D$ inmediatamente da la salida $y = \psi(q, x)$ y cambia su estado inicial por $\varphi(q, x)$.

Al unirse a la salida de $A_q$ con la entrada de otro autómata

$$B_s = (D, S, \alpha, \beta, s),$$

se obtiene un mapa que corresponde a un autómata llamado la composición de $A_q$ y $B_s$ denotado por $A_q \star B_s$.

Este autómata se describe formalmente como el autómata cuyo conjunto de estados es $Q \times S$ y cuyas funciones de transición $\Phi$ y salida $\Psi$ se definen por

$$\Phi((x, y), i) = (\varphi(x, i), \alpha(y, \psi(x, i))),$$

$$\Psi((x, y), i) = \beta(y, \psi(x, i)),$$

con el estado inicial $(q, s)$.

La composición $A \star B$ de dos autómatas no iniciales se define por las mismas fórmulas para las funciones de entrada y de salida, pero sin especificar el estado inicial.

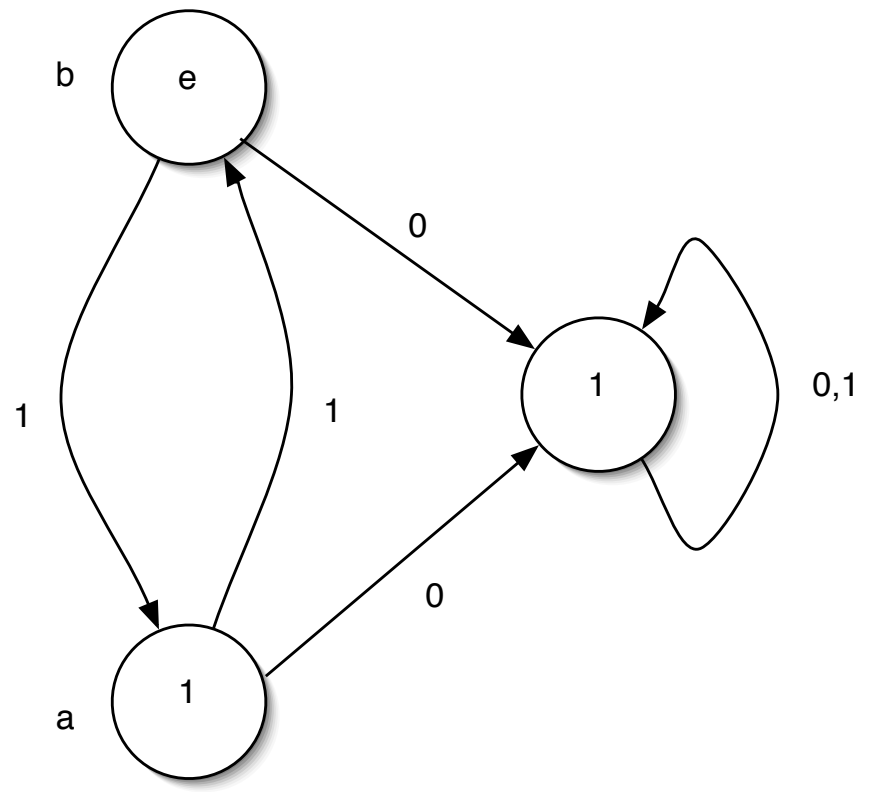

Figura 1: Autómata asociado a un grupo amenable.

Se dice que dos autómatas iniciales son equivalentes si determinan la misma aplicación por todos sucesiones finitas y infinitas en $D$. Existe un algoritmo para minimizar el número de estados.

El autómata que produce el mapa de identidad en el conjunto de sucesiones se llama trivial. Si $A$ es inversible entonces para cada estado $q$ el autómata $A_q$ admite un autómata inverso $A_q^{-1}$ tal que $A_q \star A_q^{-1}$, $A_q^{-1} \star A_q$ son equivalentes a los autómatas triviales. El autómata inverso formalmente puede ser descripto como el autómata $(D, Q, \widetilde{\varphi}, \widetilde{\psi}, q)$ donde $\widetilde{\varphi}(s, i) = \varphi(s, \sigma_s(i))$ y $\widetilde{\psi}(s, i) = \sigma_s^{-1}(i)$ para $s \in Q$. Las clases de equivalencia de autómatas finitos inversibles sobre un alfabeto $D$ es un grupo que se llama el grupo de los autómatas finitos, que depende de $D$. Cada conjunto de autómatas iniciales genera un subgrupo de este grupo.

Ahora $A$ es un autómata inversible no inicial. Sea $Q = \{q_1, \ldots, q_t\}$ el conjunto de estados de $A$ y sea $A_{q_1}, \ldots, A_{q_t}$ todos los autómatas iniciales que pueden ser obtenidos a partir de $A$. El grupo $G(A) = \langle A_{q_1}, \ldots, A_{q_t} \rangle$ se llama el grupo generado por $A$.

Para más información sobre grupos generados para un autómata, ver Żuk (2008).

**Teorema 5** (Grigorchuk y Żuk (2002))**.** *El grupo $G$ generado por el autómata de la Figura 1 no es subexponencial amenable, i.e. $G \notin SG$.*

**Demostración (Esbozo. Para más detalles ver Żuk (2008))**
Supongamos que $G \in SG_\alpha$ para $\alpha$ minimal. Entonces $\alpha$ no puede ser 0 ya que $G$ es de crecimiento exponencial porque contiene un semigrupo libre. Además, $\alpha$ no es un ordinal límite, porque si $G \in SG_\alpha$ para un ordinal límite $\alpha$ entonces $G \in SG_\beta$ para un ordinal $\beta < \alpha$. Además, $G$ no es límite directo (de una sucesión creciente de grupos) ya que es de tipo finito. Así existe $N$, $H \in SG_{\alpha-1}$ tal que la siguiente sucesión es exacta:

$$1 \to N \to G \to H \to 1.$$

Para el grupo $G$ y cada subgrupo normal $N \lhd G$ que no es trivial, se cumple la siguiente propiedad: existe un subgrupo de $N$ con $G$ como cociente. Para demostrarlo se utilizan las acciones de grupos de autómatas en árboles arraigados. No es difícil demostrar que cada clase $SG_\alpha$ es cerrada con respecto a cocientes y a subgrupos. Deducimos que $G \in SG_{\alpha-1}$. Contradicción. □

Para demostrar la amenabilidad de $G$ utilizamos un criterio de Kesten de paseos al azar en $G$.

## 3 Paseos al azar

Ahora definiremos un operador $M : l^2(\Gamma) \longrightarrow l^2(\Gamma)$ por la fórmula:

$$Mf(g) = \frac{1}{|S|} \sum_{s \in S} f(sg)$$

Vamos a considerar el caso especial donde $S = S^{-1}$. $M$ es autoadjunto, y tenemos el siguiente teorema:

**Teorema 6** (Kesten (1959))**.** *$\Gamma$ es amenable si y sólo si $\|M\| = 1$.*

*Notas* :

- Es fácil ver que $\|M\| \leqslant 1$.
- Es fácil demostrar que amenabilidad implica que $\|M\| = 1$.

  De hecho, para una secuencia de Følner $A_n$, tenemos $\mathbb{1}_{A_n} \in l^2(\Gamma)$ y $\frac{\|M\mathbb{1}_{A_n}\|}{\|\mathbb{1}_{A_n}\|} \to 1$, por lo tanto, el resultado!

Hay muchas caracterizaciones diferentes de la noción de amenabilidad (lo que demuestra que el concepto es rico), pero esta caracterización dada por Kesten es particularmente importante.

**Definición** : Sea un grupo $\Gamma$ de tipo finito y un subconjunto $S$. Decimos que $\mathrm{Cay}(\Gamma, S)$ es un grafo de Cayley si es un grafo cuyos vértices son los elementos $\gamma \in \Gamma$ y cuyas aristas son los pares $(\gamma, s\gamma)$ para $s \in S$ y $\gamma \in \Gamma$.
*Nota* : Un grafo de Cayley depende de los generadores elegidos. Por ejemplo, en el caso de $\mathbb{Z}$, el grafo es evidente si se toma como subconjunto generador a $S = \{\pm 1\}$ :

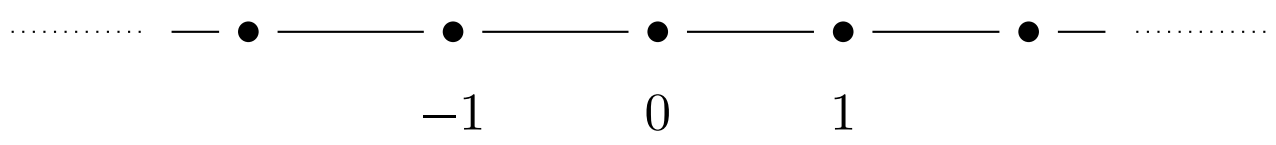

pero también puede ser más sorprendente si se tiene $S = \{\pm 2, \pm 3\}$ :

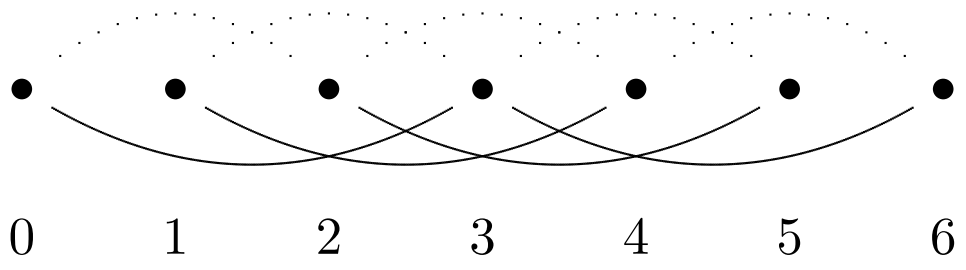

*Ejemplos* :

- En el caso de grupo $\mathbb{Z}^2$, podemos tomar $S = \{(0, \pm 1), (\pm 1, 0)\}$. El grafo de Cayley es entonces una retícula regular en el plano.
- En el caso del grupo $F_2 =< a, b >$, uno puede elegir $S = \{a^{\pm 1}, b^{\pm 1}\}$. El grafo tiene la siguiente forma :

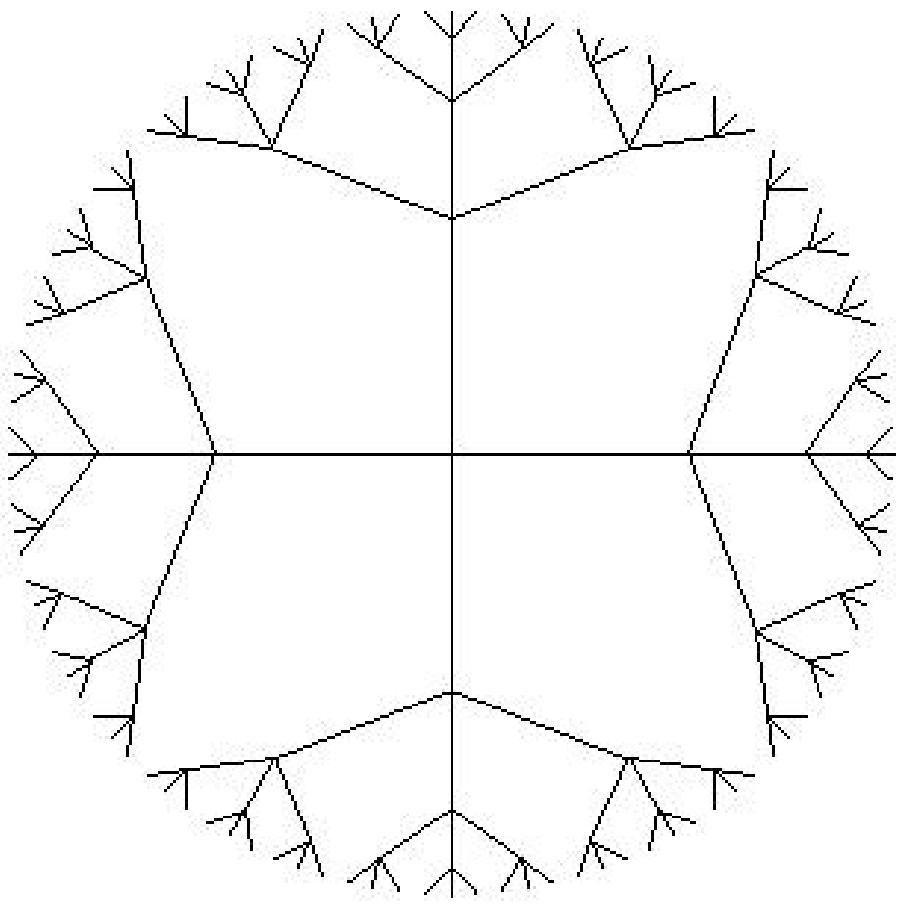

Figura 2: Grafo de Cayley de $F_2$, con $S = \{a^{\pm 1}, b^{\pm 1}\}$.

Se trata de un árbol, es decir, no tiene ciclos.

**Definiciones** :

- Si $G$ es un grupo de tipo finito, un subconjunto finito $S$ que genera $G$ es simétrico si $s \in S \implies s^{-1} \in S$.
- Por un grupo $G$ de tipo finito y un subconjunto simétrico $S$ que genera $G$, un paseo al azar en $G$ es *simple* si todos los elementos de $S$ son equiprobables.

**Teorema 7** (Kesten (1959))**.** *Sea $G$ un grupo de tipo finito y un subconjunto simétrico $S$ que genera $G$. Para un paseo al azar simple, le notemos $p_n(id, id)$ la probabilidad de volver al origen después de $n$ pasos. El grupo $G$ es amenable si y sólo si:*

$$\lim_{n\to\infty} \sqrt[2n]{p_{2n}(id, id)} = 1.$$

## 4 Funciones de tipo positivo

**Definición** : Una función $\varphi : G \to \mathbb{C}$ se llama *de tipo positivo* si para todo $g_1, \cdots, g_n \in G$ y todo $\lambda_1, \cdots, \lambda_n \in \mathbb{C}$, tenemos

$$\sum_{i=1}^{n}\sum_{j=1}^{n} \overline{\lambda_i}\lambda_j \varphi(g_i^{-1}g_j) \geqslant 0$$

*Ejemplo* : Sea $\pi : G \to B(\mathscr{H})$ una representación del grupo $G$ en un espacio de Hilbert $\mathscr{H}$. Si fijamos un vector $\xi \in \mathscr{H}$, entonces

$$\varphi(g) = < \pi(g)\xi, \xi >$$

es una función de tipo positivo.

**Teorema 8.** *Un grupo $G$ es amenable si y sólo si existe una sucesión de funciones de tipo positivo $\varphi_i$ con soporte finito tal que*

$$\lim_{i\to\infty} \varphi_i(g) = 1$$

*para cualquier $g \in G$.*

**Demostración**

…

($\Rightarrow$) Sea $\lambda$ la representación regular de $G$. Entonces, para $A_n$ una sucesión de Følner y $\xi_n = \frac{1}{|A_n|}\chi_{A_n}$,

$$\varphi_n(g) =< \lambda(g)\xi_n, \xi_n >$$

tiene las propiedades requeridas. □

## 5 Grafos expansores y aplicaciones

**Definiciones** :

- Sea $X$ un grafo finito, y $A$ un subconjunto de $X$. Definiremos el *borde* $\partial A$ *de* $A$ como el conjunto de aristas tales que un extremo está en $A$ y el otro en $A^c$.
- Definiremos la *constante isoperimétrica* de un grafo $X$ como

$$h(X) = \min\left\{\frac{|\partial A|}{|A|} : A \subseteq X, 1 \leqslant |A| \leqslant \frac{|X|}{2}\right\}$$

*Nota* : En el caso de los grafos infinitos, tomamos subconjuntos finitos $A$ de $X$, y eliminamos la condición $|A| \leqslant |X|/2$.

*Ejemplos* :

- Si $X = \mathrm{Cay}(\mathbb{Z}, \{\pm 1\})$, es claro que $h(X) = 0$. Basta considerar los intervalos de números enteros $\{-n, -n+1, \cdots, n-1, n\}$ (cuyo borde es de cardinal 2) y $n$ tiende a infinito.

- De manera más general, para cualquier grupo amenable se tiene que $h(X) = 0$. En efecto, si $A_n$ es una sucesión de Følner, entonces $\frac{|\partial A_n|}{|A_n|} \longrightarrow 0$, pues $|\partial A_n| \leqslant \sum_{s \in S} |A_n \Delta s A_n|$ y $S$ es finito.

- Si $X = \mathrm{Cay}(F_2, \{a^{\pm 1}, b^{\pm 1}\})$, entonces $h(X) = 2$.

**Definición** : Una sucesión de grafos finitos $X_n$ ($|X_n| < \infty$) de grado $k$ (donde $k \geqslant 3$) es una *sucesión de grafos expansores* si $|X_n| \to \infty$ y existe $c > 0$ tal que $h(X_n) \geqslant c > 0$.

No podemos producir grafos expansores de los grupos amenables como demuestra la siguiente proposición:

**Proposicion 9.** *Sea* $\Gamma$ *un grupo amenable generado por un subconjunto finito* $S$, *y* $\Gamma_n$ *una sucesión de cocientes finitos de* $\Gamma$.

*Luego Cay*$(\Gamma_n, S)$ *no es una sucesión de grafos expansores.*

*Ejemplo* : Veamos un ejemplo donde $\Gamma$ y $\Gamma_n$ verifican las condiciones : tomamos $\Gamma = \mathbb{Z}$ y $\Gamma_n = \mathbb{Z}/n\mathbb{Z}$.

Para que la noción que se nos acaba de presentar tenga sentido debemos preguntarnos si hay grafos expansores.

Demostraremos que casi todos los grafos son grafos expansores. Podemos comparar este enfoque con la demostración de la existencia de los números trascendentes en $\mathbb{R}$. En esta demostración, se demuestra que casi todos los reales son trascendentes! Tengamos en cuenta que es bastante fácil de probar la existencia de los números trascendentes, pero es mucho más difícil probar que $\pi$ o $e$ son transcendentes.

A los efectos de nuestra demostración, introducimos conjuntos de grafos.

**Definición** : Los conjuntos de grafos $X(n,k)$ se definen de la siguiente manera:

- Sus vértices son el conjunto $\{0,1\} \times \{1, \cdots, n\}$, que puede ser representado como dos filas de $n$ puntos cada uno, uno frente al otro. Así, cada grafo en $X(n,k)$ tiene $2n$ vértices.

- Las aristas se definen mediante $k$ permutaciones $\pi_1, \cdots, \pi_k \in S_n$: existe una arista entre $(0,i)$ y $(1,j)$ si existe un $k$ tal que $\pi_k(i) = j$.

Es fácil comprobar que estos grafos son de grado $k$. También tengamos en cuenta que no identificamos grafos isomorfos.

**Teorema 10.** *Supongamos que $k \geqslant 5$. Entonces*

$$\lim_{n\to\infty} \frac{\#\{X \in X(n,k) : h(X) \geqslant 1/2\}}{\#X(n,k)} = 1.$$

Para demostrar el teorema, vamos a utilizar una "aproximación" de la noción de constante isoperimétrica, válida para los grafos de los tipos de $X(n,k)$, que es mucho más fácil de usar que el concepto general.

**Definiciones** :

- Definiremos $\partial' A' = \{x \in X \setminus A' : \text{existe una arista } (x,y) : y \in A'\}$
- Sea $X \in X(n,k)$. Definiremos $I$ la primera fila de vértices de $X$, i.e. los puntos $\{(0,1), \cdots, (0,n)\}$, y $O$ la segunda fila. Definiremos :

$$h'(X) = \min\left\{\frac{|\partial' A'|}{|A'|} : A' \subseteq I, |A'| \leqslant \frac{|I|}{2} = \frac{n}{2}\right\}$$

**Proposicion 11.** *Sea $X \in X(n,k)$, entonces $h(X) \geqslant h'(X) - 1$.*

**Demostración del teorema**
En lugar de demostrar que $h(X) \geqslant \frac{1}{2}$, vamos a demostrar que $h'(X) \geqslant \frac{3}{2}$.

El número de todas las permutaciones $(\pi_1, \cdots, \pi_n)$ es $(n!)^k$.

Vamos a acotar el número de $(\pi_1, \cdots, \pi_n)$ "inválidos", es decir, los que dan $X$ con $h'(X) < \frac{3}{2}$. Si $h'(X) \leqslant \frac{3}{2}$, existe $A \subseteq I$ con $|A| < \frac{1}{2}n$ y existe $B \subseteq O$ con $|B| = \frac{3}{2}|A|$ tal que :

$$\partial' A \subseteq B$$

De este modo podemos acotar el número de $(\pi_1, \cdots, \pi_n)$ "inválidos" por

$$\sum_{\substack{A\subseteq I\\ |A|\leqslant \frac{n}{2}}} \sum_{\substack{B\subseteq O\\ |B|=\frac{3}{2}|A|}} \left(\frac{|B|!(n-|A|)!}{(|B|-|A|)!}\right)^k$$

Podemos reescribir la última suma sobre el cardinal $i$ de $A$ que se convierte

$$\sum_{i=1}^{n/2} \binom{n}{i}\binom{n}{3i/2}\left(\frac{(3i/2)!(n-i)!}{(i/2)!}\right)^k = \sum_{i=1}^{n/2}\alpha(i)$$

La función $\alpha(i)$ está disminuyendo para $i = 1, \ldots, n/3$. Por lo tanto, la suma anterior de $i$ de 1 a $n/3$ puede ser limitada por

$$n^3((n-1)!)^k. \tag{5-1}$$

Por $n/3 \leqslant i \leqslant n/2$ uno tiene

$$\binom{n}{i}\binom{n}{3i/2} \leqslant 2^{2n}$$

and $\beta(i) = \left(\frac{(3i/2)!(n-i)!}{(i/2)!}\right)^k$ tiene su máximo para $i = n/3$ o $n/2$. Así, la suma en este rango puede ser limitada por

$$n2^{2n}(\beta(n/3) + \beta(n/2)). \tag{5-2}$$

No es difícil demostrar que expresión (5-1) y expresión (5-2) dividido por $(n!)^k$ tiende a 0, lo que demuestra el teorema. □

## 6 Propiedad (T)

**Definiciones** :

- Sea $\pi : G \longrightarrow B(\mathcal{H})$ una representación de un grupo $G$ sobre $\mathcal{H}$. Esta representación es *unitaria* si $\pi(g) \in \mathcal{U}(\mathcal{H})$ para cualquier $g$ en $G$.
- Diremos que una representación $\pi$ *casi tiene un vector invariante* si para cada $\varepsilon > 0$ y para cada subconjunto compacto $K$ de $G$, existe un vector $\xi$ de $\mathcal{H}$ :

$$\|\pi(k)\xi - \xi\| < \varepsilon\|\xi\| \text{ para cada } k \text{ en } K.$$

*Nota* : Si $G$ es un grupo discreto, substituimos la condición $K$ compacto por $K$ finito.

*Ejemplo* : Sea $G = \mathbb{Z}$, $\mathcal{H} = l^2(\mathbb{Z})$ y consideremos la representación regular $\lambda$ :

$$\lambda(n)f(m) = f(m+n) \text{ para todo } f \in l^2(\mathbb{Z}) \text{ y } n \in \mathbb{Z}.$$

$\lambda$ casi tiene un vector invariante. De hecho, sea $\xi_n = \mathbb{1}_{\{-n,\cdots,n\}}$. Consideramos un subconjunto finito $K$ de $\mathbb{Z}$. Entonces

$$\lim_{n\to\infty} \frac{\|\xi_n - \lambda(k)\xi_n\|}{\|\xi_n\|} = 0 \text{ para todo } k \text{ en } K.$$

Sin embargo, la representación $\lambda$ no tiene un vector invariante. (Tal vector tendría que ser una función constante y de norma $l^2$ finita, lo cual es absurdo.)

**Proposicion 12.** *Las siguientes afirmaciones son equivalentes :*

*i. El grupo G es amenable;*

*ii. La representación regular de G casi tiene un vector invariante.*

**Demostración**
La dirección $i \implies ii$ es muy fácil utilizando las sucesiones de Følner. En concreto, tomamos $\xi_n = \mathbb{1}_{A_n}$ y entonces se utiliza la definición de sucesiones de Følner.

Para demostrar el converso, vamos a construir una sucesión de Følner. Sea $K$ un subconjunto compacto o finito de $G$. Consideremos un vector casi invariante (de norma 1) para la representación regular, i.e. $f \in l^2(G)$ tal que $\sum_{g\in G} f(g)^2 = 1$ y $\sum_{g\in G} |f_k(g) - f(g)|^2 \leqslant \varepsilon$ para todo $k$ en $K$. (Recordemos que $f_k$ representa $f$ trasladado por $k$.)

Sea $F = f^2$. Entonces $F \in l^1(G)$, y $\|F\| = 1$. Fijamos $k$ en $K$. Podemos escribir :

$$\begin{aligned}\|F_k - F\|_1 &= \sum_{g\in G} |f_k(g)^2 - f(g)^2| = \sum_{g\in G} |f_k(g) - f(g)|\,|f_k(g) + f(g)| \\ &\leqslant \left(\sum_{g\in G}(f_k(g) - f(g))^2\right)^{1/2}\left(\sum_{g\in G}|f_k(g) + f(g)|^2\right)^{1/2} \\ &\leqslant \varepsilon^{1/2} 2\|f\|_2 = 2\varepsilon^{1/2} = \varepsilon'\end{aligned}$$

y podemos hacer que la cantidad sea arbitrariamente pequeña.

Para todo $a$ positivo, sea $\mathfrak{U}_a = \{g \in G : F(g) \geqslant a\}$. Entonces

$$1 = \sum_{g\in G} F(g) = \int_0^\infty |\mathfrak{U}_a| da$$

y

$$||F_k - F||_1 = \sum_{g\in G} |F_k(g) - F(g)| = \int_0^\infty |\mathfrak{U}_a \Delta k \mathfrak{U}_a| da \leqslant \varepsilon'$$

Deducimos que $\int_0^\infty |\mathfrak{U}_a \Delta k \mathfrak{U}_a| da \leqslant \varepsilon' \int_0^\infty |\mathfrak{U}_a| da$ y así

$$\int_0^\infty \sum_{k\in K} |\mathfrak{U}_a \Delta k \mathfrak{U}_a| da \leqslant \varepsilon' |K| \int_0^\infty |\mathfrak{U}_a| da.$$

Por lo tanto existe $a$ tal que $\sum_{k\in K} |\mathfrak{U}_a \Delta k \mathfrak{U}_a| \leqslant \varepsilon' |K|\,|\mathfrak{U}_a|$.

Finalmente sea $\varepsilon'' = \varepsilon'|K|$. Entonces

$$|\mathfrak{U}_a \Delta k \mathfrak{U}_a| \leqslant \varepsilon'' |\mathfrak{U}_a| \text{ para todo } k \text{ en } K.$$

□

**Definición *(Každan (1967))***: $G$ tiene la *propiedad (T)* si cualquier representación unitaria de $G$ que casi tiene un vector invariante, tiene un vector invariante.

*Ejemplo* : Como muestra el ejemplo al principio de esta sección $\mathbb{Z}$ no tiene la propiedad (T). De manera más general, según la proposición 12 un grupo amenable casi tiene un vector invariante pero tiene vectores invariantes si y sólo si es finito.

**Proposicion 13** (Margulis (1973))**.** *Sea $G$ un grupo de tipo finito, y $S$ un subconjunto finito que genera $G$. Supongamos que $G$ tiene la propiedad (T). Si $G_n$ es una sucesión de cocientes finitos de $G$, entonces existe una constante $c > 0$ tal que*

$$h\left(Cay(G_n, S)\right) \geqslant c,$$

*i.e. $Cay(G_n, S)$ es una sucesión de grafos expansores.*

**Demostración**
Supongamos que esta constante $c$ no existe. Vamos a construir una representación de $G$ que casi tiene un vector invariante, sin que haya vectores invariantes.

Si esta constante $c$ no existe, para todo $\varepsilon > 0$ existe un grupo $G_n$ tal que $h(\mathrm{Cay}(G_n, S)) \leqslant \varepsilon$, i.e. existe un subconjunto $A_n \subseteq G_n$ con $|A_n| \leqslant \frac{|G_n|}{2}$ tal que para todo $s$ de $S$,

$$\frac{|A_n \Delta s A_n|}{|A_n|} \leqslant \varepsilon$$

Consideramos $\mathbb{1}_{A_n} \in l^2(G_n)$. Entonces

$$\|\mathbb{1}_{A_n} - s\mathbb{1}_{A_n}\|_2^2 = \|\mathbb{1}_{A_n} - \mathbb{1}_{s\,A_n}\|_2^2 \leqslant \varepsilon \|\mathbb{1}_{A_n}\|_2^2 \tag{6-1}$$

i.e. el vector $\mathbb{1}_{A_n}$ es casi invariante. El problema es que $\mathbb{1}_{G_n}$ es un vector invariante.

Para evitar este problema, consideramos el espacio vectorial

$$l_0^2(G_n) = \{f \in l^2(G_n), \sum_{g \in G} f(g) = 0\}$$

Sea $a$ una constante positiva tal que $|A_n| = a|A_n^c|$. Entonces $\mathbb{1}_{A_n} - a\mathbb{1}_{A_n^c} \in l_0^2(G)$. Podemos acotar $a : |A_n| \leqslant \frac{|G_n|}{2}$ implica que $a \leqslant 1$.

Por otra parte, podemos escribir $\mathbb{1}_{A_n^c} = \mathbb{1}_{G_n} - \mathbb{1}_{A_n}$. Esto nos lleva a $\mathbb{1}_{A_n} - a\mathbb{1}_{A_n^c} = (1+a)\mathbb{1}_{A_n} - a\mathbb{1}_{G_n}$.

Para probar el análogo de (6-1) de este vector, se escribe

$$\begin{aligned}\|(1+a)\mathbb{1}_{A_n} - a\mathbb{1}_{G_n} - ((1+a)\mathbb{1}_{s\,A_n} - a\mathbb{1}_{s\,G_n})\| &= (1+a)\,\|\mathbb{1}_{A_n} - \mathbb{1}_{s\,A_n}\| \\ &\leqslant \varepsilon(1+a)\|\mathbb{1}_{A_n}\| \\ &\leqslant 2\varepsilon\|\mathbb{1}_{A_n} - a\mathbb{1}_{A_n^c}\|\end{aligned}$$

La primera desigualdad proviene de (6-1). La segunda de $a \leqslant 1$ y de $\|\mathbb{1}_{A_n}\| \leqslant \|\mathbb{1}_{A_n} - a\mathbb{1}_{A_n^c}\|$. Este último aumento se demuestra usando que $\mathbb{1}_{A_n}$ y $\mathbb{1}_{A_n^c}$ son ortogonales en $l^2(G_n)$. Consideremos $\oplus_{n\geqslant 1} l_0^2(G_n)$. Así, hemos construido una representación de $G$ que casi tiene un vector invariante y no tiene un vector invariante lo que contradice el hecho de que $G$ tiene la propiedad (T). □

Cuando Kazhdan introdujo la propiedad (T) tenía la idea de usarla para estudiar los subgrupos discretos de los grupos de Lie. Podemos plantear el siguiente problema:

> Si $\Gamma$ es un subgrupo discreto de $\mathrm{SL}(n, \mathbb{R})$ con $\mathrm{Vol}(\mathrm{SL}(n, \mathbb{R})/\Gamma) < \infty$, nos preguntamos si $\Gamma$ es de tipo finito.

Si $n = 2$, la respuesta es positiva y se conoce desde los trabajos de Poincaré. Para que la misma afirmación sea cierta para $n \geqslant 3$, Kazhdan introdujo un nuevo enfoque.

**Proposicion 14** (Každan (1967)). *Si $G$ es numerable y tiene la propiedad (T), entonces $G$ es de tipo finito.*

**Demostración**
Numeremos los elementos de $G : G = \{g_1, g_2, \cdots\}$. Ahora definiremos $G_n$ como el subgrupo de $G$ generado por $\{g_1, g_2, \cdots, g_n\}$. Si $G$ no es de tipo finito, para todo $n$, $|G : G_n| = \infty$. Notamos que $G/G_n$ es un cociente pero no es siempre un grupo.

Consideramos $\oplus_{n\in\mathbb{N}}\, l^2(G/G_n)$. Para todo $n$, $G$ actúa en $G/G_n$. Para todo $n_0 \in \mathbb{N}$, podemos considerar $\delta_n \in \oplus l^2(G/G_n)$ que se define por :

$$\delta_{n_0} = (0, \cdots, 0, \underbrace{\delta_{G_{n_0}}}_{\in l^2(G/G_{n_0})}, 0, \cdots)$$

donde $\delta_{G_{n_0}}$ es una función que es 1 en $G_{n_0}$ y 0 de otro modo.

Claramente, $\delta_n$ es invariante por $G_n$. Pero por cada subconjunto finito $K$ de $G$, existe $n$ tal que $K \subseteq G_n$. Así que esta representación casi tiene vectores invariantes, pero no hay un vector invariante lo cual contradice la propiedad (T). □

Todavía tenemos que preguntarnos cómo probar la propiedad (T) ! La respuesta está dada por el teorema de Kazhdan:

**Teorema 15.** *Si $G$ es un grupo de Lie, y si $\Gamma \subseteq G$ es una retícula (i.e. $\Gamma$ discreto y $Vol(G/\Gamma) < \infty$) entonces $G$ tiene la propiedad (T) si y sólo si $\Gamma$ tiene la propiedad (T).*

**Demostración** $(\Rightarrow)$
Para demostrar la implicación que nos interesa $(\Rightarrow)$, sólo se tiene que inducir a la representación de $\Gamma$ en $G$. □

Se puede demostrar que todos los grupos de Lie simples de rango por lo menos dos tienen la propiedad (T) Margulis (1991).

**Proposicion 16** (Každan (1967))**.** *Si $G$ tiene la propiedad (T) entonces $|G/G'| < \infty$, donde $G'$ es la abelianización de $G$.*

**Demostración**
Si $G/G'$ es infinito, tenemos un epimorfismo $\varphi : G \longrightarrow \mathbb{Z}$. Dado que $\mathbb{Z}$ no tiene la propiedad (T), podríamos construir una representación unitaria de $G$ que casi tiene vectores invariantes, y que no tiene vectores invariantes distintos de cero. □

Como los subgrupos de índices finitos de grupos con propiedad (T) también tienen esta propiedad, se obtiene como consecuencia que todos los subgrupos de índices finitos de los grupos de Kazhdan tienen abelianización finita. Esta propiedad puede ser muy difícil de probar por otros métodos.

Recuerde el siguiente problema, que data de los años 20:

> La medida de Lebesgue es la única medida finitamente aditiva en $S^n$, invariante bajo $SO(n+1)$, definida en los conjuntos medibles de Lebesgue de medida total 1?

Para $n = 1$, Banach (1923) en los años 20 dio una respuesta negativa a la pregunta.

Para $n \geqslant 4$, Margulis (1982) y Sullivan (1981) demostraron que la respuesta es afirmativa. Para demostrar esta propiedad, que utiliza la propiedad (T).

Drinfel'd (1984) finalmente demostró que la respuesta es afirmativa para los casos más difíciles, $n = 3$ y 2.

La unicidad se demostró con la siguiente propiedad: hay $k$ elementos $\rho_i$ de $SO(n+1)$ y $\varepsilon > 0$ tal que para todo $f \in L^2_0(S^n)$

$$\max_i \|f - \rho_i f\|_2 \geqslant \varepsilon \|f\|_2,$$

lo que implica que no hay vectores casi invariantes y que explica la conexión con la propiedad (T)

No obstante, sigue siendo un problema abierto, si la propiedad vale si elegimos las rotaciones $\rho_i$ (por ejemplo dos) al azar. Si esto se mantiene, sería otro ejemplo de la situación cuando para una propiedad dada casi todas las opciones, pero es bastante difícil exponer un ejemplo explícito.

## 7 Funciones condicionalmente de tipo negativo

**Definición** : Una función $\varphi : G \to \mathbb{R}$ es condicionalmente *de tipo negativo* si

1. $\varphi(\mathrm{Id}) = 0$
2. para $g_1, \cdots, g_n \in G$ y $\lambda_1, \cdots, \lambda_n \in \mathbb{R}$ cualesquiera tales que

$$\sum_{i=1}^{n} \lambda_i = 0,$$

tenemos que

$$\sum_{i=1}^{n} \sum_{j=1}^{n} \lambda_i \lambda_j \varphi(g_i^{-1} g_j) \leqslant 0$$

Existe una relación directa entre las funciones de tipo positivo y las funciones condicionalmente de tipo negativo. Es un resultado fácil y bien conocido:

**Teorema 17.** *Una función $\varphi : G \to \mathbb{R}$ tal que $\varphi(\mathrm{Id}) = 0$ es condicionalmente de tipo negativo si y sólo si para todo $t > 0$,*

$$e^{-t\varphi}$$

*es de tipo positivo.*

No es obvio, pero podemos mostrar la siguiente caracterización muy útil de la propiedad $(T)$.

**Teorema 18.** *Un grupo $G$ tiene la propiedad $(T)$ si y sólo si cualquier función condicionalmente de tipo negativo en $G$ es acotada.*

En los últimos años esta condición se utilizó para probar la negación de la propiedad (T) para muchas clases de grupos. Hay muchas maneras de definir una función condicionalmente de tipo negativo en $G$. Por ejemplo para un grupo $G$ generado por un subconjunto finito $S$, podemos definir la norma $\|g\|_S$ para $g \in G$ como el mínimo número de generadores necesarios para escribir $g$.

Se puede demostrar que los grupos libres con respecto a un conjunto de generadores estándar son condicionalmente de tipo negativo.

## 8 Propiedad (T) para grupos discretos

La siguiente condición simple permite probar la propiedad (T) para muchos grupos discretos.

Sea $\Gamma$ un grupo generado por un conjunto finito $S$ tal que $S$ sea simétrico, es decir, $S = S^{-1}$, y que el elemento de identidad $e$ no pertenezca a $S$.

**Definición** :

Definimos un grafo finito $L(S)$, de la siguiente manera:

1. vértices de $L(S) = \{s; s \in S\}$,
2. aristas de $L(S) = \{(s, s'); s, s', s^{-1}s' \in S\}$.

Supongamos que el grafo $L(S)$ es conexo. Esta condición no es restrictiva, porque para un grupo finitamente generado $\Gamma$ siempre se puede encontrar un conjunto generador finito y simétrico $S$, que no contiene $e$, tal que $L(S)$ sea conexo. (Por ejemplo, $S \cup S^2 \setminus e$ satisface esta condición.) Esto puede ser visto en el caso simple de $\Gamma = \mathbb{Z}$; si $S = \{-1, 1\}$ entonces el grafo $L(S)$ no es conexo, pero si añadimos al conjunto de generadores $\{-2, 2\}$ el grafo vuelve conexo.

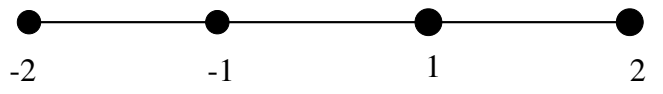

Figura 3: El grafo $L(S)$ por $\Gamma = \mathbb{Z}$ y $S = \{-2, -1, 1, 2\}$

Para un vértice $s \in L(S)$ sea $deg(s)$ su grado, es decir, el número de aristas adyacentes a $s$. Sea $\Delta$ un operador de Laplace discreto actuando sobre funciones definidas en vértices de $L(S)$, es decir, para $f \in l^2(L(S), deg)$

$$\Delta f(s) = f(s) - \frac{1}{deg(s)} \sum_{s' \sim s} f(s'),$$

donde $s' \sim s$ significa que el vértice $s'$ es adyacente al vértice $s$.

El operador $\Delta$ es un operador nonegativo, autoadjunto en $l^2(L(S), deg)$. Si $L(S)$ está conectado entonces el cero es un valor propio simple de $\Delta$. Sea $\lambda_1(L(S))$ el menor valor propio $\neq 0$ de $\Delta$ actuando sobre $l^2(L(S), deg)$.

**Teorema 19** (Żuk (2003))**.** *Sea $\Gamma$ un grupo generado por un subconjunto finito $S$, tal que $S$ es simétrico y $e \notin S$. Si el grafo $L(S)$ está conectado y*

$$\lambda_1(L(S)) > \frac{1}{2} \tag{8-1}$$

*entonces $\Gamma$ tiene la propiedad (T) de Kazhdan.*

**Observación** La condición indicada en el teorema 1 es óptima. Para ver esto, consideremos el grupo $\Gamma = \mathbb{Z}$ con el conjunto de generadores $S = \{\pm 1, \pm 2\}$. Entonces el grafo $L(S)$ son cuatro vértices y tres bordes (véase la figura 3). Para este grafo $\lambda_1(L(S)) = \frac{1}{2}$ y el grupo $\mathbb{Z}$ no tiene la propiedad (T).

La condición anterior se aplica a varios grupos. Por ejemplo, se puede dar una nueva prueba de propiedad (T) para algunos retículos. Veamos esto para las retículos en $SL(3, \mathbb{Q}_p)$.

Sea una familia de grupos que actúan transitivamente sobre los vértices de alguna clase de complejos simpliciales que se llaman edificios del tipo $\widetilde{A}_2$. Estos grupos están parametrizados por un entero $q$ que es una potencia de un número primo. Admiten una presentación tal que $L(S)$ es el gráfico de incidencia del plano proyectivo $\mathbb{P}^2(\mathbb{F}_q)$ sobre el cuerpo finito $\mathbb{F}_q$, i.e.

$$\begin{aligned} \text{vértices de } L(S) &= \{\text{puntos } p \text{ y lineas } l \text{ tales que } p, l \in \mathbb{P}^2(\mathbb{F}_q)\}, \\ \text{aristas de } L(S) &= \{(p, l); p \in l\}. \end{aligned}$$

Figura 4 muestra el grafo $L(S)$ para $q = 2$.

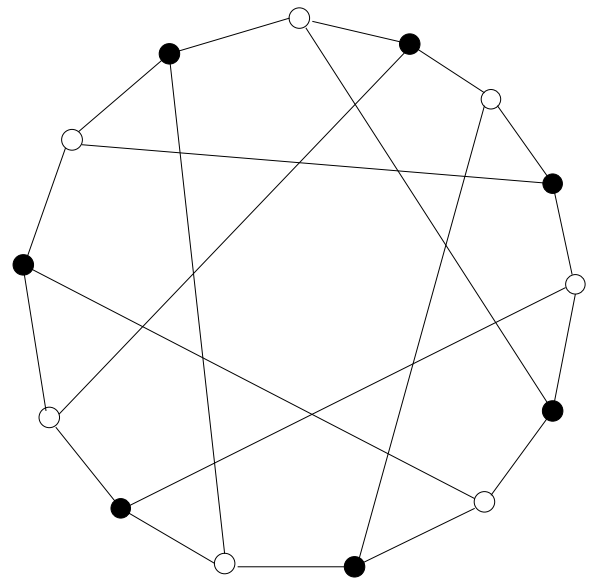

Figura 4: El grafo $L(S)$ para un retículo en $SL(3, \mathbb{Q}_2)$

Feit y Higman han calculado el espectro del operador de Laplace en grafos de incidencia de planos proyectivos finitos.

**Proposicion 20** (Feit–Higman)**.** *Sea $L$ el grafo de incidencia $\mathbb{P}^2(\mathbb{F}_q)$. Entonces*

$$\lambda_1(L) = 1 - \frac{\sqrt{q}}{q+1}.$$

Para cualquier $q \geqslant 2$ tenemos $\lambda_1(L) = 1 - \frac{\sqrt{q}}{q+1} > \frac{1}{2}$ y por lo tanto estos grupos tienen la propiedad (T).

Usando el Teorema 19 se puede demostrar que casi todo grupo discreto tiene la propiedad (T).

Para más información sobre la amenabilidad y los grupos con la propiedad $(T)$, ver Colin de Verdière (1998), Lubotzky (1994) y Sarnak (1990).

## Referencias

ANDRZEJ ŻUK
andrzej.zuk@imj-prg.fr
INSTITUT DE MATHÉMATIQUES
UNIVERSITÉ PARIS 7
13 RUE ALBERT EINSTEIN
75013 PARIS
FRANCE

# 6 | Equidistribución y análisis diofántico

# ALTURAS Y DISTRIBUCIÓN DE PUNTOS ALGEBRAICOS

Ricardo Menares

**Resumen**

Las raíces del polinomio $x^n - 1$ se sitúan sobre el círculo unitario formando los vértices de un polígono regular de n lados. Cuando n crece, los polígonos aproximan al círculo cada vez mejor. En otras palabras, las raíces de la unidad se reparten de manera uniforme sobre el círculo cuando el orden tiende a infinito. Por otro lado, la familia de polinomios $(x-1)^n$ tiene sólo una raíz. El contraste entre la distribución límite de las raíces (uniforme en el primer ejemplo, concentrada en un punto en el segundo) se explica por la manera en que crecen los coeficientes del polinomio (nada en el primer caso y exponencialmente en el segundo). El ejemplo de las raíces de la unidad se extiende, más generalmente, al caso de sucesiones de puntos algebraicos "pequeños" en una curva (Teoremas de equidistribución de Bilu y de Szpiro–Ullmo–Zhang). La manera correcta de cuantificar el tamaño de los puntos algebraicos es a través de la teoría de alturas, que presentaremos en la primera parte del curso. La demostración de los teoremas de equidistribución mencionados requiere herramientas de geometría algebraica, teoría de números y análisis armónico. Trataremos de presentar, en líneas generales, la interacción de estas técnicas en este contexto.

## Índice general

Se contó con apoyo de Fondecyt 11110225 y Proyecto Investigación Asociativa PUCV 037.469/2015.

# 1 Introducción

Estas notas corresponden a dos charlas a dictarse el contexto de la escuela Aritmética, Grupos y Análisis (AGRA 2), que tendrá lugar en el Cusco, Perú, en Agosto de 2015. Aprovecho este espacio para agradecer a los organizadores la posibilidad de exponer el presente tópico.

Se ha tratado de mantener los requisitos necesarios para seguir estas charlas al mínimo. Para poder exponer el tema en el tiempo que se nos ha dado, hemos decidido concentrarnos en un caso especial del Teorema de Bilu, que data de 1997 (Teorema 5.7 en el texto).

# 2 Números algebraicos

**Definición 2.1.** Un elemento $\alpha \in \mathbb{C}$ se dice *número algebraico* si existe un polinomio no constante, con coeficientes racionales, $f(x) \in \mathbb{Q}[x]$, tal que $f(\alpha) = 0$.

Notar que todo número racional es algebraico. En efecto, si $\alpha = \frac{a}{b}$, con $a, b \in \mathbb{Z}$ y $b \neq 0$, entonces podemos tomar $f(x) = x - \frac{a}{b}$. Más ejemplos:

**Ejemplos 2.2.**

1. $\alpha = i, \quad f(x) = x^2 + 1$
2. $\alpha = \sqrt{2}, \quad f(x) = x^2 - 2$
3. $\alpha = \sqrt[3]{2}, \quad f(x) = x^3 - 2$
4. $\alpha = \zeta_n := e^{2\pi i/n}, \quad f(x) = x^n - 1, \quad n \in \mathbb{Z}_{>0}$
5. $\alpha = \frac{1+\sqrt{5}}{2}, \quad f(x) = x^2 - x - 1$

Denotamos $\overline{\mathbb{Q}}$ al conjunto de los números algebraicos. Notar que $\overline{\mathbb{Q}}$ es un conjunto numerable. En efecto, $\mathbb{Q}[x]$ es numerable y cada elemento no constante de $\mathbb{Q}[x]$ tiene a lo más un número finito de raíces. Dado que $\mathbb{C}$ no es numerable, se desprende que existen números que no son algebraicos (y de hecho son mayoría). Sin embargo, no es fácil identificar un número no algebraico. Por ejemplo, se sabe que $e$ (Hermite 1873) y $\pi$ (Lindemann 1882) no son algebraicos, pero a la redacción de estas líneas no se sabe decidir si $e + \pi$ es algebraico o no.

El polinomio $f(x)$ que figura en la Definición 2.1 no es único. Cualquier polinomio de la forma $h(x) = f(x)g(x)$, con $g(x) \in \mathbb{Q}[x]$, sirve también. Sin embargo, nos será útil contar con un polinomio asociado de manera canónica a un número algebraico.

**Proposición 2.3.** *Sea $\alpha \in \overline{\mathbb{Q}}$. Entonces existe un único polinomio no constante $f_\alpha(x)$ que cumple*

1. *$f_\alpha(x)$ tiene coeficientes racionales*
2. *$f_\alpha(\alpha) = 0$*
3. *$f_\alpha(x)$ es mónico (es decir, el coeficiente del término dominante es 1)*
4. *el grado de $f_\alpha(x)$ es mínimo entre los polinomios que satisfacen (1), (2) y (3)*

**Demostración**: el principio del buen orden nos asegura que existe un polinomio $f(x) \in \mathbb{Q}[x]$ que satisface las cuatro propiedades del enunciado. Si $g(x)$ es otro polinomio que cumple las mismas propiedades, entonces podemos aplicar división de polinomios

$$f(x) = q(x)g(x) + r(x), \quad q(x), r(x) \in \mathbb{Q}[x], \quad \deg r(x) < \deg g(x).$$

Como $r(\alpha) = f(\alpha) - q(\alpha)g(\alpha) = 0$, de la minimalidad del grado de $f(x)$ concluimos que $r(x)$ es el polinomio nulo, es decir $f(x) = q(x)g(x)$. Como $f$ y $g$ tienen el mismo grado y son mónicos, entonces $f = g$ □

**Definición 2.4.**
- Decimos que el polinomio $f_\alpha(x)$ dado por la Proposición 2.3 es el *polinomio mínimo* de $\alpha$.
- Definimos el *grado de* $\alpha$ por $\deg \alpha := \deg f_\alpha$.

**Definición 2.5.** Un polinomio $f(x) \in \mathbb{Q}[x]$ se dice $\mathbb{Q}$-reductible si se puede escribir como producto de polinomios, con coeficientes racionales, de menor grado. Un polinomio en $\mathbb{Q}[x]$ que no es $\mathbb{Q}$-reductible se dice $\mathbb{Q}$-irreductible.

**Observación**: la minimalidad del grado de $f_\alpha(x)$ asegura que éste es un polinomio $\mathbb{Q}$-irreductible. Recíprocamente, se tiene el siguiente resultado.

**Proposición 2.6.** *Sea* $\alpha \in \overline{\mathbb{Q}}$ *y sea* $f(x) \in \mathbb{Q}[x]$ *un polinomio* $\mathbb{Q}$*-irreductible y mónico tal que* $f(\alpha) = 0$. *Entonces* $f = f_\alpha$.

**Demostración**: aplicar el algoritmo de la división para polinomios □

**Definición 2.7.** Dado un polinomio

$$f(x) = a_n x^n + a_{n-1}x^{n-1} + \cdots + a_1 x + a_0 \in \mathbb{C}[x],$$

definimos el *polinomio reverso*

$$f^*(x) = a_0 x^n + a_1 x^{n-1} + \cdots + a_{n-1}x + a_n \in \mathbb{C}[x].$$

**Lema 2.8.** *Si* $f(x) \in \mathbb{Q}[x]$ *es* $\mathbb{Q}$*-irreductible, entonces* $f^*(x)$ *también lo es.*

**Demostración**: basta notar que $f^*(x) = x^{\deg f} f(1/x)$ y usar la Definición 2.5 □

De la Proposición 2.6 y el Lema 2.8, se deduce

**Corolario 2.9.** *Si* $\alpha \in \overline{\mathbb{Q}}$ *y* $\alpha \neq 0$, *entonces* $1/\alpha \in \overline{\mathbb{Q}}$. *Más aún,* $f_{1/\alpha} = f_\alpha^*$.

No siempre es fácil determinar el polinomio mínimo de un número algebraico. En los ejemplos 2.2, los polinomios indicados son todos irreductibles (luego coinciden con el polinomio mínimo) excepto en el ejemplo 2.2, (4). En efecto,

$$x^n - 1 = (x-1)(x^{n-1} + x^{n-2} + \cdots + x + 1).$$

Dos herramientas básicas para decidir la irreductibilidad de un polinomio con coeficientes racionales son el Lema de Gauss y el Criterio de Eisenstein, que procedemos a explicar.

**Definición 2.10.** Un polinomio con coeficientes enteros $p(x) \in \mathbb{Z}[x]$ se dice $\mathbb{Z}$-reductible si existen polinomios $f(x), g(x) \in \mathbb{Z}[x]$ tales que

- $\deg f, \deg g < \deg p$
- $p(x) = f(x)g(x)$

Diremos que $p(x) \in \mathbb{Z}[x]$ es $\mathbb{Z}$-irreductible si no es $\mathbb{Z}$-reductible.

**Teorema 2.11.** *(Lema de Gauss) Un polinomio $p(x) \in \mathbb{Z}[x]$ que es $\mathbb{Z}$-irreductible es también $\mathbb{Q}$-irreductible.*

**Teorema 2.12.** *(Criterio de Eisentein) Sea $f(x) = a_n x^n + a_{n-1}x^{n-1} + \cdots + a_1 x + a_0$ un polinomio con coeficientes enteros. Suponga que existe un primo $p$ tal que*

- $p | a_i, \quad i = 0, 1, \ldots, n-1$
- $p \nmid a_n$
- $p^2 \nmid a_0$

*Entonces $f(x)$ es $\mathbb{Q}$-irreductible.*

La demostración de estos teoremas está esbozada en los Ejercicios 2.15, (4) y (5).

Denotamos por $\Phi_n(x)$ al polinomio mínimo de $\zeta_n = e^{2\pi i/n}$.

**Observación 2.13.** Se tiene que $\Phi_n(x)$ es un divisor (en $\mathbb{Q}[x]$) del polinomio $x^n - 1$ (cf. Ejercicio 2.15 (1)). Del lema de Gauss, deducimos que $\Phi_n(x) \in \mathbb{Z}[x]$.

Como aplicación de los teoremas anteriores, demostraremos el siguiente

**Lema 2.14.** *Sea $p$ un número primo. Entonces $\Phi_p(x) = x^{p-1} + x^{p-2} + \cdots + x + 1$. En particular,*

$$x^p - 1 = (x-1)\Phi_p(x).$$

**Demostración**: Sea $w(x) = x^{p-1} + x^{p-2} + \cdots + x + 1$. Veamos que $w(x)$ es un polinomio $\mathbb{Q}$-irreductible. Esto último equivale a demostrar que $h(x) := w(x+1)$ es $\mathbb{Q}$-irreductible. Usando la identidad $x^p - 1 = (x-1)w(x)$, se tiene

$$h(x) = w(x+1) = \frac{(x+1)^p - 1}{x} = x^{p-1} + \left( \sum_{k=2}^{p-1} \binom{p}{k} x^{k-1} \right) + p.$$

Dado que $p \left| \binom{p}{k} \right.$, para todo $0 < k < p$, el criterio de Eisenstein permite concluir que $h(x)$ es $\mathbb{Q}$-irreductible. Tenemos entonces que $w(x)$ es $\mathbb{Q}$-irreductible y mónico, luego $w(x) = \Phi_p(x)$ por la Proposición 2.6 □

**Ejercicios 2.15.** 1. Sea $f(x) \in \mathbb{Q}[x]$ un polinomio $\mathbb{Q}$-irreductible y sea $\alpha \in \mathbb{C}$ tal que $f(\alpha) = 0$. Sea $g(x) \in \mathbb{Q}[x]$ un polinomio no constante tal que $g(\alpha) = 0$. Demuestre que existe un polinomio $h(x) \in \mathbb{Q}[x]$ tal que $g(x) = f(x)h(x)$.

2. Un polinomio $p(x)$ se dice que tiene raíces repetidas si se puede factorizar sobre $\mathbb{C}$ de la forma
$$p(x) = (x - z)^2 h(x), \quad z \in \mathbb{C}, \quad h(x) \in \mathbb{C}[x]$$
(de manera equivalente, $p(z) = p'(z) = 0$). Demuestre que un polinomio $\mathbb{Q}$-irreductible $p(x) \in \mathbb{Q}[x]$ no puede tener raíces repetidas.

3. Sea $p$ un primo. Denotamos $\mathbb{F}_p := \mathbb{Z}/p\mathbb{Z}$, el cuerpo finito de $p$ elementos. Sea
$$\nu : \mathbb{Z} \longrightarrow \mathbb{Z}/p\mathbb{Z}$$
el morfismo canónico. Para un polinomio con coeficientes enteros $f(x) = a_n x^n + a_{n-1}x^{n-1} + \cdots + a_1 x + a_0$, definimos
$$\bar{f}(x) = \nu(a_n)x^n + \nu(a_{n-1})x^{n-1} + \cdots + \nu(a_1)x + \nu(a_0) \in \mathbb{F}_p[x].$$
Muestre que la operación $f \mapsto \bar{f}$ define un morfismo de anillos $\mathbb{Z}[x] \longrightarrow \mathbb{F}_p[x]$ (es decir, verifique $\overline{f+g} = \bar{f} + \bar{g}$ y $\overline{fg} = \bar{f}\bar{g}$).

4. *Lema de Gauss*. Decimos que un polinomio $f(x) \in \mathbb{Z}[x]$ es *reducido* si el máximo común divisor de sus coeficientes es 1.

   (a) Sean $f(x), g(x) \in \mathbb{Z}[x]$ dos polinomios y sea $h(x) = f(x)g(x)$. Sea $N$ el máximo común divisor de los coeficientes de $h(x)$. Suponga $N \neq 1$ y tome un primo $p$ tal que $p|N$. Usando el morfismo de anillos del ejercicio anterior, demuestre que $p$ o bien divide a todos los coeficientes de $f(x)$ o bien divide a todos los coeficientes de $g(x)$

   (b) Deduzca que el producto de dos polinomios reducidos es reducido

   (c) Demuestre el Teorema 2.11

5. *Criterio de Eisenstein*. Sean $f(x) \in \mathbb{Z}[x]$ un polinomio y $p$ un primo que satisfacen las hipótesis del Teorema 2.12.

   (a) Suponga que se puede factorizar $f(x) = g(x)h(x)$ con $g(x), h(x) \in \mathbb{Z}[x]$. Muestre que entonces se tiene
   $$\nu(a_n)x^n = \bar{g}(x)\bar{h}(x), \quad \text{en } \mathbb{F}_p[x].$$

   (b) Justifique que $\bar{g}(x) = ux^a, \bar{h} = vx^b$, con $u, v \in \mathbb{F}_p^*$ y $a + b = n$.

   (c) Demuestre el Teorema 2.12.

6. Pruebe que $\overline{\mathbb{Q}}$ es un cuerpo. Es decir,
$$\alpha, \beta \in \overline{\mathbb{Q}}, \alpha \neq 0 \Rightarrow \frac{1}{\alpha}, \alpha\beta, \alpha + \beta \in \overline{\mathbb{Q}}.$$

### 2.1 Órbita galoisiana.

**Definición 2.16.**
- Sean $\alpha \in \overline{\mathbb{Q}}$ y $f_\alpha(x)$ su polinomio mínimo. Decimos que $\beta \in \mathbb{C}$ *es un conjugado de* $\alpha$ si $f_\alpha(\beta) = 0$.
- Definimos la *órbita galoisiana* de $\alpha$ por

$$\begin{aligned} G(\alpha) &:= \{\beta \in \mathbb{C} : \beta \text{ es un conjugado de } \alpha\} \\ &= \{\beta \in \mathbb{C} : f_\alpha(\beta) = 0\} \end{aligned}$$

**Observación 2.17.** El número de conjugados de $\alpha$ es exactamente el grado de $f_\alpha(x)$ (ver ejercicio 2.15, (2) ).

**Ejemplos 2.18.**
- $G(i) = \{i, -i\}$
- $G(\sqrt{2}) = \{\sqrt{2}, -\sqrt{2}\}$
- $G(\sqrt[3]{2}) = \{\sqrt[3]{2}, \zeta_3 \sqrt[3]{2}, \zeta_3^2 \sqrt[3]{2}\}$
- $G(\frac{1+\sqrt{5}}{2}) = \{\frac{1+\sqrt{5}}{2}, \frac{1-\sqrt{5}}{2}\}$

Es de interés para nosotros calcular la órbita galoisiana de $\zeta_n$. Definimos

$$\mu_n := \{z \in \mathbb{C} : z^n = 1\} = \{\text{raíces de la unidad de orden } n\}$$
$$\tilde{\mu}_n := \{z \in \mu_n : z^k \neq 1, \forall 1 \leqslant k < n\} = \{\text{raíces primitivas de la unidad de orden } n\}$$

Notar que $(\mu_n, \cdot)$ es un grupo cíclico de orden $n$. Más precisamente, si $\zeta_n = e^{2\pi i/n}$, entonces $\mu_n = \{\zeta_n^j : 0 \leqslant j \leqslant n-1\}$. Más aún, se tiene

$$\tilde{\mu}_n = \{\zeta_n^j : mcd(j,n) = 1\}. \tag{2.1}$$

En particular, el número de elementos de $\tilde{\mu}_n$ es

$$\varphi(n) := |(\mathbb{Z}/n\mathbb{Z})^*| = \#\{1 \leqslant j \leqslant n-1 : mcd(j,n) = 1\}$$

(función de Euler).

**Proposición 2.19.** *Sea* $\zeta_n = e^{2\pi i n}$. *Entonces*

$$G(\zeta_n) = \tilde{\mu}_n.$$

La demostración de este hecho está esbozada en los Ejercicios 2.21.

**Corolario 2.20.** *Se tiene* $\deg \Phi_n = \varphi(n)$*, para todo entero positivo* $n$.

**Ejercicios 2.21.** 1. Demuestre (2.1).

2. (a) Demuestre que para todo $\alpha \in \overline{\mathbb{Q}}$ y todo entero positivo $k$, se tiene
$$G(\alpha^k) = \{\beta^k : \beta \in G(\alpha)\}.$$
Deduzca $\deg(\alpha^k) \leqslant \deg(\alpha)$.
   (b) Muestre que para $\alpha \in \overline{\mathbb{Q}}^*$, se tiene $G(1/\alpha) = \{1/\beta : \beta \in G(\alpha)\}$. Deduzca $\deg(\alpha) = \deg(1/\alpha)$.
3. Sean $n$ un entero y $p$ un número primo que no divide $n$. Sea $\alpha \in \mathbb{C}$ tal que $\Phi_n(\alpha) = 0$. El objetivo de este problema es demostrar que $\Phi_n(\alpha^p) = 0$.
   (a) Muestre que $\alpha^p \in \overline{\mathbb{Q}}$ y que su polinomio mínimo tiene coeficientes enteros.
   (b) Suponga que $\Phi_n(\alpha^p) \neq 0$. Sea $g(x) \in \mathbb{Z}[x]$ el polinomio mínimo de $\alpha^p$. Sea $h(x) = g(x^p)$. Muestre que $\Phi_n(x) | h(x)$ en $\mathbb{Z}[x]$.
   (c) Sea $j(x) \in \mathbb{F}_p[x]$ un factor irreductible de $\overline{\Phi}_n(x) \in \mathbb{F}_p[x]$ (cf. la operación
$$f \in \mathbb{Z}[x] \mapsto \bar{f} \in \mathbb{F}_p[x]$$
definida en el Ejercicio 2.15 (3)). Muestre que $j(x) | \bar{g}(x)$.
   (d) Use lo anterior para mostrar que $j(x)^2 | x^n - 1$. Concluya.
4. Use el ejercicio anterior para demostrar la Proposición 2.19.

## 3 Medida de Mahler de un polinomio en una variable

Sea $f(x) \in \mathbb{C}[x]$ un polinomio no nulo con coeficientes complejos. Definimos su *medida de Mahler* por
$$M(f) = \exp\Big(\frac{1}{2\pi}\int_0^{2\pi} \log|f(e^{i\theta})|d\theta\Big).$$

No es difícil ver que la integral es absolutamente convergente, incluso si el polinomio tiene raíces sobre el círculo unitario.

**Lema 3.1.** *Para todo polinomio no nulo $f \in \mathbb{C}[x]$, se tiene $M(f) = M(f^*)$.*

**Demostración**: si $n = \deg f$, tenemos $f^*(x) = x^n f(1/x)$, luego
$$M(f^*) = \exp\Big(\frac{1}{2\pi}\int_0^{2\pi} \log|f(e^{-i\theta})|d\theta\Big).$$
El enunciado resulta entonces de un cambio de variable apropiado □

Nos será útil en lo que sigue usar la función $\log^+$, dada por
$$\log^+ z = \begin{cases} \log z & \text{si } z \geqslant 1 \\ 0 & \text{si } 0 < z < 1 \end{cases}$$

La medida de Mahler satisface las propiedades siguientes:

**Teorema 3.2.** *1. $M(fg) = M(f)M(g)$, para todo $f, g \in \mathbb{C}[x]$, $fg \neq 0$.*

*2. $M(f) > 0$, para todo $f \in \mathbb{C}[x]$, no nulo.*

*3. (Fórmula de Jensen) Sea $f(x) \in \mathbb{C}[x]$, un polinomio no nulo de grado $n$, que escribimos de la forma*

$$f(x) = a(x-\alpha_1)(x-\alpha_2)\cdots(x-\alpha_n), \quad a, \alpha_1, \ldots, \alpha_n \in \mathbb{C}. \tag{3.2}$$

*Entonces*

$$\log M(f) = \log|a| + \sum_{i=1}^{n} \log^+ |\alpha_i|. \tag{3.3}$$

**Demostración**: Las primeras dos propiedades se deducen directamente de la definición. Usando (1), nos reducimos a mostrar que para todo $\alpha \in \mathbb{C}$, se tiene

$$\frac{1}{2\pi}\int_0^{2\pi} \log|e^{i\theta} - \alpha| d\theta = \log^+ |\alpha|. \tag{3.4}$$

Para demostrar esta identidad, usaremos que el promedio de una función armónica $h$ sobre el borde del círculo unitario es $h(0)$.

Si $|\alpha| > 1$, entonces la función $h(x) = \log|x-\alpha|$ es armónica en el disco unitario, de manera que el término izquierdo de (3.4) es $h(0) = \log|\alpha|$.

Si $|\alpha| < 1$, entonces la función $w(x) = \log|1-\alpha\bar{x}|$ es armónica en el disco unitario y coincide con $h(x)$ cuando $|x| = 1$. Por lo tanto el término izquierdo de (3.4) es $w(0) = \log|1| = 0$.

El caso $|\alpha| = 1$ se deduce por la continuidad de la función $\alpha \mapsto \int_0^{2\pi} \log|e^{i\theta} - \alpha| d\theta$ □

**3.1 Discriminante y Resultante.** Representamos al polinomio $f(x) \in \mathbb{C}[x]$ de la forma (3.2). Notar que las raíces $\alpha_1, \alpha_2, \ldots, \alpha_n$ no son necesariamente distintas. Se define el *discriminante* de $f$ por

$$D(f) = a^{2n-2}\prod_{i<j}(\alpha_i - \alpha_j)^2.$$

Es inmediato de la definición que $D(f) = 0$ si y sólo si $f$ tiene raíces repetidas.

**Definición 3.3.** Sean

$$f(x) = a_n x^n + a_{n-1}x^{n-1} + \cdots + a_1 x + a_0 \in \mathbb{C}[x]$$
$$g(x) = b_m x^m + b_{m-1}x^{m-1} + \cdots + b_1 x + b_0 \in \mathbb{C}[x].$$

Formamos la siguiente matriz de tamaño $(n+m) \times (n+m)$

$$M(f,g) = \begin{pmatrix} a_n & a_{n-1} & \ldots & a_0 & 0 & \ldots & 0 \\ 0 & a_n & a_{n-1} & \ldots & a_0 & \ldots & 0 \\ & \ldots & & & & \ldots & \\ 0 & 0 & \ldots & a_n & a_{n-1} & \ldots & a_0 \\ b_m & b_{m-1} & \ldots & b_0 & & & \\ 0 & b_m & b_{m-1} & \ldots & b_0 & & \\ & \ldots & & & & & \\ 0 & 0 & \ldots & b_m & b_{m-1} & \ldots & b_0 \end{pmatrix}$$

Definimos la *resultante* $R(f,g)$ por

$$R(f,g) := \det M(f,g).$$

**Teorema 3.4.** *(Lang (2002), Proposición 8.3, p. 202) Escribimos $g(x) = b_m(x-\beta_1)\cdots(x-\beta_m)$. Entonces se tiene*

$$R(f,g) = a_n^m b_m^n \prod_{i=1}^{n}\prod_{j=1}^{m}(\alpha_i - \beta_j). \tag{3.5}$$

**Observación 3.5.** De este Teorema vemos que $R(f,g) = 0$ si y sólo si $f$ y $g$ tienen una raíz en común.

Central en lo que sigue es la siguiente

**Proposición 3.6.** *Si $f(x)$ tiene coeficientes enteros y $a_n$ es el coeficiente del término dominante, entonces $a_n D(f) \in \mathbb{Z}$.*

**Demostración**: tomemos $g = f'$ en (3.5). Entonces $m = n-1$ y $b_m = na_n$. Si escribimos

$$f'(x) = na_n(x-\beta_1)\cdots(x-\beta_{n-1}),$$

entonces

$$R(f,f') = n^n a_n^{2n-1}\prod_{i=1}^{n}\frac{f'(\alpha_i)}{na_n} = a_n^{n-1}\prod_{i=1}^{n} f'(\alpha_i).$$

Pero $f'(\alpha_i) = a_n \prod_{\substack{j=1\\ j\neq i}}^{n}(\alpha_i - \alpha_j)$, luego obtenemos

$$\begin{aligned}R(f,f') &= a_n^{2n-1}\prod_{i=1}^{n}\prod_{\substack{j=1\\ j\neq i}}^{n}(\alpha_i-\alpha_j)\\ &= (-1)^{n(n-1)/2}a_n^{2n-1}\prod_{i<j}(\alpha_i-\alpha_j)^2\\ &= (-1)^{n(n-1)/2}a_n D(f).\end{aligned}$$

Como $f$ y $f'$ tienen coeficientes enteros, de la definición de resultante tenemos $R(f,f') \in \mathbb{Z}$□

**Ejercicio 3.7.** Demostrar

$$|D(f)| \leqslant n^n M(f)^{2n-2}.$$

Indicación: exprese el discriminante como un determinante de Vandermonde.

## 4 Altura de un número algebraico

Sea $\alpha \in \overline{\mathbb{Q}}^*$ y sea $f_\alpha(x) \in \mathbb{Q}[x]$ su polinomio mínimo. Existe un único entero positivo $m$ tal que el polinomio

$$g_\alpha(x) := mf_\alpha(x)$$

satisface

- $g_\alpha(x)$ tiene coeficientes enteros
- el máximo común divisor de los coeficientes de $g_\alpha(x)$ es 1.

Decimos que $g_\alpha(x)$ es el *polinomio mínimo sobre $\mathbb{Z}$ de $\alpha$*.

El claro que la órbita galoisiana de $\alpha$ se calcula por

$$G(\alpha) = \{\beta \in \mathbb{C} : g_\alpha(\beta) = 0\}$$

(cf. Definición 2.16).

**Definición 4.1.** Sea $\alpha \in \overline{\mathbb{Q}}^*$. Definimos su *altura* (de Weil) por

$$h(\alpha) = \frac{1}{\deg \alpha} \log M(g_\alpha).$$

De los Lemas 2.8 y 3.1 deducimos

**Lema 4.2.** *Para todo $\alpha \in \overline{\mathbb{Q}}^*$, se tiene $h(\alpha) = h(1/\alpha)$.*

Escribiendo

$$g_\alpha(x) = a_n x^n + a_{n-1} x^{n-1} + \cdots + a_1 x + a_0,$$

de la fórmula de Jensen (3.3) se tiene

$$h(\alpha) = \frac{1}{\deg \alpha}\Big(\log a_n + \sum_{\beta \in G(\alpha)} \log^+ |\beta|\Big). \tag{4.6}$$

**Proposición 4.3.** *Se tiene que $h(\alpha) \geqslant 0$, para todo $\alpha \in \overline{\mathbb{Q}}^*$. Más aún, $h(\alpha) = 0$ si y sólo si $\alpha$ es una raíz de la unidad.*

**Demostración**: en la expresión (4.6), $a_n$ es un entero positivo, de donde es claro que $h(\alpha) \geqslant 0$. Si $\alpha$ es una raíz de la unidad, la Proposición 2.19 asegura que $|\beta| = 1$, para todo $\beta \in G(\alpha)$. Además, de la Observación 2.13 tenemos que su polinomio mínimo sobre $\mathbb{Z}$ es mónico. Usando nuevamente (4.6), tenemos $h(\alpha) = 0$.

Supongamos ahora $h(\alpha) = 0$. Sea $n = \deg(\alpha)$. Escribamos

$$g_\alpha(x) = a_n x^n + a_{n-1} x^{n-1} + \cdots + a_1 x + a_0, \quad a_i \in \mathbb{Z}. \tag{4.7}$$

Entonces la expresión (4.6) garantiza que $a_n = 1$ y

$$|\beta| \leqslant 1, \tag{4.8}$$

para todo $\beta \in G(\alpha)$.

La desigualdad (4.8) también vale para productos de elementos en $G(\alpha)$. Entonces, de la fórmula $g_\alpha(x) = \prod_{\beta \in G(\alpha)}(x - \beta)$, vemos que los coeficientes de $g_\alpha(x)$ satisfacen

$$|a_j| \leqslant n, \quad \forall j = 0, 1, 2, \ldots, n. \tag{4.9}$$

Se $S_n$ el conjunto formado por todos los polinomios con coeficientes enteros, de grado a lo más $n$, de la forma (4.7) que además satisfacen (4.9). Tenemos que $S_n$ es un conjunto finito.

Ahora bien, si $k$ es un entero positivo, usando el Ejercicio 2.21, (2) vemos que el razonamiento anterior también se aplica a $\alpha^k$, es decir, $\{g_{\alpha^k}(x) : k \in \mathbb{Z}_{>0}\} \subseteq S_n$. Como $S_n$ es finito, el conjunto de raíces de polinomios en $S_n$ también lo es, luego el conjunto

$$\{\alpha^k : k \in \mathbb{Z}_{>0}\}$$

es finito. Por lo tanto, existen enteros $k_1 \neq k_2$ tales que $\alpha^{k_1} = \alpha^{k_2}$ (equivalentemente, $\alpha^{k_1-k_2} = 1$), lo que prueba que $\alpha$ es una raíz de la unidad □

Terminamos esta sección con un problema abierto.

**Pregunta 4.4.** (Lehmer (1933)). Decidir si la siguiente afirmación es cierta: existe una constante $c > 0$ tal que

$$h(\alpha) \geqslant \frac{c}{\deg \alpha}, \quad \forall \alpha \in \overline{\mathbb{Q}}^* \text{ que no es una raíz de la unidad.}$$

**Ejercicio 4.5.** Adapte la demostración de la Proposición 4.3 para demostrar el teorema de Northcott: Sean $A, B > 0$. El conjunto

$$\{\alpha \in \overline{\mathbb{Q}}^* : \deg(\alpha) \leqslant A, \quad h(\alpha) \leqslant B\}$$

es finito.

# 5 Teorema de Bilu

## 5.1 Convergencia de medidas y enunciado del Teorema de Bilu.

**Definición 5.1.** Sea $M$ el conjunto de medidas de probabilidad sobre $\mathbb{C}$. Definimos

$$C_0(\mathbb{C}) = \{f : \mathbb{C} \to \mathbb{C} \text{ continua y de soporte compacto}\}.$$

Decimos que una sucesión $(\nu_n) \subset M$ converge a $\nu \in M$ si para toda función $f \in C_0(\mathbb{C})$, se tiene

$$\lim_{n\to\infty} \int_{\mathbb{C}} f \nu_n = \int_{\mathbb{C}} f \nu.$$

En este caso, escribimos

$$\lim_{n\to\infty} \nu_n = \nu.$$

**Ejemplos 5.2.** Algunas medidas de probabilidad:

- Denotamos por $\nu_S \in M$ a la medida uniforme sobre el círculo. Está caracterizada por

$$\int_{\mathbb{C}} f \nu_S = \frac{1}{2\pi} \int_0^{2\pi} f(e^{i\theta}) d\theta, \quad \forall f \in C_0(\mathbb{C}).$$

- Dado un punto $z \in \mathbb{C}$, denotamos por $\delta_z \in M$ a la medida de Dirac soportada en $z$. Está caracterizada por

$$\int_{\mathbb{C}} f\delta_z = f(z), \quad \forall f \in C_0(\mathbb{C}).$$

- Más generalmente, sea $E \subset \mathbb{C}$ un conjunto finito. Definimos

$$\delta_E := \frac{1}{\#E} \sum_{z \in E} \delta_z \in M. \tag{5.10}$$

Esta medida está caracterizada por

$$\int_{\mathbb{C}} f\delta_E = \frac{1}{\#E} \sum_{z \in E} f(z), \quad \forall f \in C_0(\mathbb{C}).$$

Las medidas de esta forma son ejemplos de *medidas atómicas*.

**Definición 5.3.** Sean $C \subset \mathbb{C}$ y $\nu \in M$. Decimos que $\nu$ *está soportada en* $C$ si para toda $f \in C_0(\mathbb{C})$ cuyo soporte es disjunto a $C$, se tiene

$$\int_{\mathbb{C}} f\mu = 0.$$

Si $\nu$ está soportada en un conjunto compacto, decimos que $\nu$ tiene soporte compacto.

Denotamos por $M_C \subseteq M$ al conjunto de medidas soportadas en $C$.

Es claro que la medida $\nu_S$ está soportada en el círculo unitario y una medida de la forma $\delta_E$ está soportada en $E$. Todas estas medidas tienen soporte compacto.

**Observación 5.4.** Cuando $C$ es compacto, el espacio $M_C$ es secuencialmente compacto (Billingsley (1999)). Es decir, toda sucesión $(\nu_n) \subset M_C$ admite una subsucesión convergente.

**Definición 5.5.** Consideremos una sucesión $E_n \subset \mathbb{C}$, donde cada $E_n$ es un conjunto finito. Decimos que la familia $(E_n)$ *se equidistribuye con respecto a* $\nu \in M$ si

$$\lim_{n \to \infty} \delta_{E_n} = \nu.$$

En otras palabras, se pide que para toda $f \in C_0(\mathbb{C})$,

$$\lim_{n \to \infty} \frac{1}{\#E_n} \sum_{z \in E_n} f(z) = \int_{\mathbb{C}} f\nu.$$

**Definición 5.6.** Decimos que una sucesión $(x_n) \subset \mathbb{C}$ es *genérica* si para todo $m \in \mathbb{N}$, el conjunto

$$\{k \in \mathbb{N} : x_k = x_m\}$$

es finito.

**Teorema 5.7.** *(Bilu (1997)) Sea $\nu_S$ la medida de probabilidad uniforme sobre el círculo. Sea $(\alpha_n) \subset \overline{\mathbb{Q}}^*$ una sucesión genérica de números algebraicos tal que*

$$\lim_{n\to\infty} h(\alpha_n) = 0.$$

*Entonces la familia de conjuntos $G(\alpha_n)$ se equidistribuye con respecto a $\nu_S$.*

En particular, se tiene el siguiente resultado notable:

**Corolario 5.8.** *La sucesión de conjuntos $(\tilde{\mu}_n)$ se equidistribuye con respecto a la medida $\nu_S$.*

**Demostración**: tomando en cuenta la Proposición 2.19 y la Proposición 4.3, el enunciado se deduce directamente del Teorema de Bilu □

**Observación 5.9.** El propósito de los ejercicios de esta sección es de indicar una demostración del Corolario 5.8 que no utiliza el Teorema de Bilu.

**Ejercicios 5.10.** 1. Sea $X$ un espacio métrico compacto. El espacio $M(X)$ (resp. $C_0(X)$) se define análogamente a $M$ (resp. a $C_0(\mathbb{C})$) en la Definición 5.1 reemplazando $\mathbb{C}$ por $X$ (resp. reemplazando $f : \mathbb{C} \to \mathbb{C}$ por $f : X \to \mathbb{C}$). La convergencia de medidas en $M(X)$ se enuncia de manera análoga la Definición 5.1, reemplazando $\int_{\mathbb{C}}$ por $\int_X$ y $C_0(\mathbb{C})$ por $C_0(X)$. El espacio $C_0(X) = C(X)$ se encuentra dotado de la topología dada por la norma supremo:

$$\|f\| = \sup_{z\in X} |f(z)|.$$

Sea $H \subset C(X)$ una subálgebra (es decir, $f, g \in H \to fg, f + g \in H$). Sea $V \subset H$ un conjunto generador (es decir, todo elemento de $H$ es una combinación $\mathbb{C}$-lineal finita de elementos de $V$). Suponga que $H$ es conjunto denso. Muestre que la sucesión de medidas $(\nu_n) \in M(X)$ converge a $\nu \in M(X)$ si y sólo si

$$\int_X f\nu_n = \int_X f\nu, \quad \forall f \in V.$$

2. (a) En la notación del ejercicio anterior, tome $X = \{z \in \mathbb{C} : |z| = 1\}$. Para cada $k \in \mathbb{Z}$, definimos $h_k : X \to \mathbb{C}$ por $h_k(x) = x^k$. Sea $H \subset C(X)$ la subálgebra generada por $\{h_k : k \in \mathbb{Z}\}$. Muestre que $H$ es densa en $C(X)$. Indicación: use el teorema de Stone–Weierstrass.

   (b) (Criterio de Weyl) Sea $E_n \subset X$ una sucesión de conjuntos finitos. Deduzca que la familia $(E_n)$ se equidistribuye con respecto a $\nu_S$ si y sólo si para todo $k \in \mathbb{Z}$, con $k \neq 0$, se tiene

$$\lim_{n\to\infty} \frac{1}{\#E_n} \sum_{z\in E_n} z^k = 0.$$

   (c) Utilizando el punto anterior, demuestre que la sucesión de conjuntos $\mu_n$ se equidistribuye con respecto a $\nu_S$.

3. El propósito de este ejercicio es demostrar el Corolario 5.8, usando los ejercicios anteriores.

(a) Sea $\mu : \mathbb{N} \to \{-1, 0, 1\}$ la función de Möbius, dada por

$$\mu(n) = \begin{cases} 0 & \text{si existe un entero } m \geqslant 2 \text{ tal que } m^2|n \\ (-1)^r & \text{si } n \text{ es el producto de } r \text{ primos distintos.} \end{cases}$$

Sea $\delta_1 : \mathbb{N} \to \{0, 1\}$ dada por

$$\delta_1(n) = \begin{cases} 1 & \text{si } n = 1 \\ 0 & \text{si } n > 1 \end{cases}$$

Demuestre

$$\sum_{d|n} \mu(d) = \delta_1(n).$$

(b) *Sumas de Ramanujan.* Denotamos por $(a, b)$ el máximo común divisor del par $a, b$. Para $n \in \mathbb{N}$ y $k \in \mathbb{Z} - \{0\}$, definimos

$$S(n,k) = \sum_{\substack{0 \leqslant j \leqslant n-1 \\ (n,j)=1}} \zeta_n^{jk}, \quad \zeta_n = e^{2\pi i/n}.$$

Muestre que

$$S(n,k) = \sum_{d|(n,k)} d\mu\Big(\frac{n}{d}\Big).$$

Indicación: escriba $S(n,k) = \sum_{j=0}^{n-1} \delta_1\big((n,k)\big)\zeta_n^{jk}$ y aplique el punto anterior.

(c) Demuestre el Corolario 5.8.

**5.2 Energía.** Sea $\nu \in M$ una medida con soporte compacto. Entonces se define su energía por

$$E(\nu) := -\int_{\mathbb{C}} \int_{\mathbb{C}} \log|z - w|\nu(z)\nu(w) \in \mathbb{R} \cup \{\infty\}.$$

**Ejemplos 5.11.**

- Usando (3.4), vemos que $E(\nu_S) = 0$.
- Para todo conjunto finito $G \subset \mathbb{C}$, se tiene

$$E(\delta_G) = \infty.$$

En vista del último ejemplo, es útil introducir, para todo conjunto finito $G \subset \mathbb{C}$, la cantidad

$$E'(\delta_G) := -\frac{1}{(\#G)^2} \sum_{\substack{z,w \in G \\ z \neq w}} \log|z - w| \in \mathbb{R}.$$

El nexo entre la teoría de alturas y la teoría del potencial está dado por el siguiente

**Lema 5.12.** *(Comparación Energía–Altura) Sea $\alpha \in \overline{\mathbb{Q}}^*$. Entonces se tiene*

$$E'(\delta_{G(\alpha)}) \leqslant 2h(\alpha).$$

**Demostración**: sea $n = \deg \alpha$, $a_n$ el coeficiente dominante de $g_\alpha$ y escribimos $G(\alpha) = \{\alpha_1, \ldots, \alpha_n\}$. Como $g_\alpha$ es irreductible, tenemos $D(g_\alpha) \neq 0$. Además, $a_n D(g_\alpha) \in \mathbb{Z}$ (Proposición 3.6), de donde

$$\log(a_n |D(g_\alpha)|) \geqslant \log 1 = 0. \tag{5.11}$$

Se tiene

$$\begin{aligned} 0 \leqslant \log a_n + \log |D(g_\alpha)| &= (2n-1)\log a_n + 2\sum_{i<j} \log|\alpha_i - \alpha_j| \\ &= (2n-1)\log a_n + \sum_{i\neq j} \log|\alpha_i - \alpha_j| \\ &= (2n-1)\log a_n - n^2 E'(\delta_{G(\alpha)}). \end{aligned}$$

Por otro lado, de (4.6) tenemos

$$\log a_n = nh(\alpha) - \sum_i \log^+ |\alpha_i| \leqslant nh(\alpha).$$

Combinando ambas desigualdades se obtiene el resultado □

En lo que sigue utilizaremos dos resultados de la teoría de capacidades. Las demostraciones pueden consultarse en Rumely (1989).

**Proposición 5.13.** *Sea $G_n \subset \mathbb{C}$ una sucesión de conjuntos finitos tales que $\#G_n$ tiende a infinito con $n$. Supongamos que existe $\nu \in M$ tal que la familia $(G_n)$ se equidistribuye con respecto a $\nu$. Entonces se tiene*

$$E(\nu) \leqslant \liminf_{n\to\infty} E'(\delta_{G_n}).$$

**Proposición 5.14.** *Sea $M_S$ el conjunto de todas las medidas en $M$ soportadas en $\{z \in \mathbb{C} : |z| = 1\}$. Entonces*

$$E(\nu) \geqslant 0, \quad \forall \nu \in M_S.$$

*Más aún,*

$$\nu \in M_S, \quad E(\nu) = 0 \Leftrightarrow \nu = \nu_S.$$

**5.3 Esbozo de la demostración del Teorema de Bilu.** En toda esta sección asumimos que la sucesión $(\alpha_n) \subset \overline{\mathbb{Q}}^*$ es genérica y cumple

$$h(\alpha_n) \to 0, \quad n \to \infty.$$

Unidas al Ejercicio 4.5, estas hipótesis aseguran que $\deg(\alpha_n)$ tiende a infinito con $n$.

La demostración del enunciado siguiente se encuentra al final de esta sección.

**Lema 5.15.** *(Equidistribución en radio). Para cada $n \in \mathbb{N}$, se puede escoger un conjunto $E_n \subseteq G(\alpha_n)$ de manera que*

1. *se tiene*

$$\lim_{n\to\infty} \frac{\#E_n}{\#G(\alpha_n)} = 1.$$

2. *Para todo $\varepsilon > 0$, existe $n_0$ tal que para todo $n \geqslant n_0$ se tiene*

$$z \in E_n \Rightarrow 1-\varepsilon \leqslant |z| \leqslant 1+\varepsilon.$$

Supongamos que podemos establecer que $\delta_{E_n} \to \nu_S$. Entonces la condición (1) en el Lema 5.15 permite mostrar que también se tiene $\delta_{G(\alpha_n)} \to \nu_S$. Luego en lo que sigue, trabajaremos con la familia $(E_n)$.

La condición (2) en el Lema 5.15, permite mostrar que existe un compacto $K \subset \mathbb{C}$ que contiene al círculo unitario tal que $\delta_{E_n} \in M_K$, para todo $n$. Como $M_K$ es secuencialmente compacto (cf. Observación 5.4), tenemos que la sucesión $\delta_{E_n}$ admite una subsucesión $\delta_{E_{n_j}}$ convergente a una medida $\nu \in M_K$.

Usando de nuevo la condición (2) en el Lema 5.15, podemos concluir que $\nu$ está soportada en el círculo unitario. De la Proposición 5.14, deducimos que $E(\nu) \geqslant 0$. Ahora bien,

$$\begin{aligned}
E(\nu) &\leqslant \liminf E'(\delta_{E_{n_j}}) && \text{Proposición 5.13}\\
&\leqslant \liminf E'(\delta_{G(\alpha_{n_j})}) && \text{condición (1) en el Lema 5.15}\\
&\leqslant 2\liminf h(\alpha_{n_j}) && \text{Lema 5.12}\\
&= 0.
\end{aligned}$$

Por lo tanto, $E(\nu) = 0$. Aplicando nuevamente la Proposición 5.14, deducimos $\nu = \nu_S$.

Hemos establecido que toda subsucesión convergente de $\delta_{E_n}$ tiene como límite a la medida $\nu_S$. Esto muestra que $\delta_{E_n}$ converge a $\nu_S$, terminando la demostración del Teorema de Bilu.

**Demostración del Lema 5.15**: Notar que por el Lema 4.2, también tenemos

$$h(1/\alpha_n) \to 0, \quad n \to \infty.$$

Procederemos por contradicción. Pasando a una subsucesión si fuese necesario, podemos suponer que existen subconjuntos $V_n \subseteq G(\alpha_n)$ tales que

- $\lim_{n\to\infty} \frac{\#V_n}{\#G(\alpha_n)} = c > 0.$
- existe $a > 1$ tal que para todo $n$, se tiene

$$z \in V_n \Rightarrow |z| > a \text{ ó } \left|\frac{1}{z}\right| > a.$$

Descomponemos $V_n = A_n \cup B_n$, donde

$$z \in A_n \Leftrightarrow |z| > a, \qquad z \in B_n \Leftrightarrow \left|\frac{1}{z}\right| > a.$$

Usando la fórmula de Jensen tenemos

$$h(\alpha_n) \geqslant \frac{1}{\deg(\alpha_n)}\Big(\sum_{\beta\in A_n} \log|\beta|\Big) \geqslant \frac{\#A_n}{\deg(\alpha_n)} \log a.$$

Por otro lado, usando el Ejercicio 2.21 (2) también se tiene

$$h(1/\alpha_n) \geqslant \frac{1}{\deg(\alpha_n)}\Big(\sum_{\beta\in B_n} \log|1/\beta|\Big) \geqslant \frac{\#B_n}{\deg(\alpha_n)}\log a.$$

Luego

$$h(\alpha_n) + h(1/\alpha_n) \geqslant \frac{\#V_n}{\deg\alpha_n}\log a.$$

Tomando el límite cuando $n \to \infty$, obtenemos

$$0 \geqslant c\log a > 0,$$

lo que es absurdo □

# 6 Bibliografía sugerida

Hemos presentado el estudio más somero posible de los números algebraicos. Sin embargo, se trata de una teoría profunda e importante dentro de la matemática pura contemporánea. Existen excelentes libros que exponen este material, por ejemplo Borevich y Shafarevich (1966), Neukirch (1999), Marcus (1977).

La noción de altura es una herramienta muy útil en teoría de números, cuyos orígenes se pueden trazar hasta Fermat. En muchos problemas modernos del área la dificultad se concentra en establecer una variante adecuada de esta noción. Recomendamos Bombieri y Gubler (2006) y Silverman (1986).

El nexo entre la teoría del potencial y la teoría de alturas está bien explicado en Rumely (1989).

# Referencias

Ricardo Menares
rmenares@mat.uc.cl
Departamento de Matemáticas, PUC
Vicuña Mackenna 4860
Santiago
Chile

# Índice de Autores

## Títulos Publicados

**Introdução à Análise Funcional** – *César R. de Oliveira*

**Introdução à Topologia Diferencial** – *Elon Lages Lima*

**Les Équations Différentielles Algébriques et les Singularités Mobiles** – *Ivan Pan y Marcos Sebastiani*

**Criptografia, Números Primos e Algoritmos** – *Manoel Lemos*

**Introdução à Economia Dinâmica e Mercados Incompletos** – *Aloísio Araújo*

**Conjuntos de Cantor, Dinâmica e Aritmética** – *Carlos Gustavo T. de A. Moreira*

**Geometria Hiperbólica** – *João Lucas Marques Barbosa*

**Introdução à Economia Matemática** – *Aloísio Araújo*

**Superfícies Mínimas** – *Manfredo Perdigão do Carmo*

**The Index Formula for Dirac Operators: an Introduction** – *Levi Lopes de Lima*

**Introduction to Symplectic and Hamiltonian Geometry** – *Ana Cannas da Silva*

**Primos de Mersenne (e outros primos muito grandes)** – *Carlos Gustavo T. de A. Moreira y Nicolau Saldanha*

**The Contact Process on Graphs** – *Márcia Salzano*

**Canonical Metrics on Compact almost Complex Manifolds** – *Santiago R. Simanca*

**Introduction to Toric Varieties** – *Jean-Paul Brasselet*

**Birational Geometry of Foliations** – *Marco Brunella*

**Introduction to Nonlinear Dispersive Equations** – *Felipe Linares y Gustavo Ponce*

**Introdução à Teoria das Probabilidades** – *Pedro J. Fernandez*

**Teoria dos Corpos** – *Otto Endler*

**Introdução à Dinâmica de Aplicações do Tipo Twist** – *Clodoaldo G. Ragazzo, Mário J. Dias Carneiro y Salvador Addas Zanata*

**Elementos de Estatística Computacional usando Plataformas de Software Livre/Gratuito** – *Alejandro C. Frery y Francisco Cribari-Neto*

**Uma Introdução a Soluções de Viscosidade para Equações de Hamilton-Jacobi** – *Helena J. Nussenzveig Lopes y Milton C. Lopes Filho*

**Elements of Analytic Hypoellipticity** – *Nicholas Hanges*

**Métodos Clássicos em Teoria do Potencial** – *Augusto Ponce*

**Variedades Diferenciáveis** – *Elon Lages Lima*

**O Método do Referencial Móvel** – *Manfredo do Carmo*

**A Student's Guide to Symplectic Spaces, Grassmannians and Maslov Index** – *Paolo Piccione y Daniel Victor Tausk*

**Métodos Topológicos en el Análisis no Lineal** – *Pablo Amster*

**Tópicos em Combinatória Contemporânea** – *Carlos Gustavo T. de A. Moreira y Yoshiharu Kohayakawa*

**Uma Iniciação aos Sistemas Dinâmicos Estocásticos** – *Paulo Ruffino*

**Compressive Sensing** – *Adriana Schulz, Eduardo A. B. da Silva y Luiz Velho*

**O Teorema de Poncelet** – *Marcos Sebastiani*

**Cálculo Tensorial** – *Elon Lages Lima*

**Aspectos Ergódicos da Teoria dos Números** – *Alexander Arbieto, Carlos Matheus y Carlos Gustavo T. de A. Moreira*

**A Survey on Hyperbolicity of Projective Hypersurfaces** – *Simone Diverio y Erwan Rousseau*

**Algebraic Stacks and Moduli of Vector Bundles** – *Frank Neumann*

**O Teorema de Sard e suas Aplicações** – *Edson Durão Júdice*

**Tópicos de Mecânica Clássica** – *Artur Lopes*

**Holonomy groups in Riemannian Geometry** – *Andrew Clarke y Bianca Santoro*

**Linear ODEs: an Algebraic Perspective** – *Letterio Gatto*

**Introdução à Análise Exploratória de Dados Multivariados** – *Pedro J. Fernandez y Victor Yohai*

**Aritmética, Grupos y Análisis, AGRA II** – *Mikhail Belolipetsky, Harald Andrés Helfgott y Carlos Gustavo Moreira (Editores)*

# Aritmética, Grupos y Análisis

He aquí las notas de la escuela AGRA II - Aritmética, Grupos y Análisis - realizada entre el 8 y el 22 de agosto del 2015, en la Universidad de San Antonio Abad del Cusco (UNSAAC), Perú.
La meta primaria de la serie de escuelas AGRA es la formación de estudiantes graduados y jóvenes investigadores en el proceso de especializarse en teoría de números, teoría de grupos, formas modulares y áreas afines.
El programa del AGRA II consistió de las unidades: Formas modulares y curvas de Shimura, Combinatoria aditiva, Introducción a la teoría de las curvas elípticas, Crecimiento en grupos y expansores y Análisis y geometría en grupos com charlas de Roberto Miatello, Gonzalo Tornaría, Ariel Pacetti, Michael Harris, Juanjo Rué, Julia Wolf, Javier Cilleruelo, Pablo Candela, Marusia Rebolledo, Marc Hindry, Harald Helfgott, Mikhail Belolipetsky y Andrzej Zuk.
Aparte de las unidades, hubo dos charlas de Ricardo Menares, cuyas notas se encuentran en el presente volúmen. Hubo también una charla por Fernando Rodríguez Villegas y una serie de dos charlas de Carlos Moreno sobre valores de funciones zeta y funciones $L$.
Las notas se han visto suplementadas por videos de algunas charlas, tomados por personal local y editados por ICTP. Los videos están disponibles en YouTube.